\documentclass{article}
\usepackage{amsmath}
\usepackage{amssymb}
\usepackage{amsthm}
\usepackage{setspace}
\usepackage{comment}
\usepackage[all,cmtip,2cell]{xy}
\usepackage{amscd}
\usepackage[margin=1in]{geometry}
\usepackage{mathrsfs}
\usepackage[bottom]{footmisc}
\usepackage{enumitem}
\usepackage{MnSymbol}
\usepackage{tablefootnote}
\usepackage{framed}
\usepackage{xparse}
\usepackage{url}

\usepackage{breakurl}
\usepackage[breaklinks]{hyperref}

\UseTwocells
\xyoption{2cell}
\everymath{\displaystyle}
\theoremstyle{plain}
\newtheorem{theorem}[subsubsection]{Theorem}
\newtheorem{lemma}[subsubsection]{Lemma}
\newtheorem{proposition}[subsubsection]{Proposition}
\newtheorem{proposition-construction}[subsubsection]{Proposition-Construction}
\newtheorem{corollary}[subsubsection]{Corollary}

\newtheorem{claim}[subsubsection]{Claim}

\newtheorem{conjecture}[subsubsection]{Conjecture}

\newtheorem{expectation}[subsubsection]{Expectation}
\newtheorem{variant}[subsubsection]{Variant}
\newtheorem{warning}[subsubsection]{Warning}

\theoremstyle{definition}
\newtheorem{definition}[subsubsection]{Definition}
\newtheorem{notation}[subsubsection]{Notation}
\newtheorem{convention}[subsubsection]{Convention}
\newtheorem{condition}[subsubsection]{Condition}

\theoremstyle{remark}
\newtheorem{construction}[subsubsection]{Construction}
\newtheorem{remark}[subsubsection]{Remark}
\newtheorem{example}[subsubsection]{Example}

\newcommand{\mb}[1]{\mathbb{#1}}
\newcommand{\mf}[1]{\mathfrak{#1}}
\newcommand{\mc}[1]{\mathcal{#1}}

\newcommand{\ms}[1]{\mathsf{#1}}
\newcommand{\msc}[1]{\mathscr{#1}}
\newcommand{\tx}{\textnormal}
\newcommand{\tb}{\textbf}

\newcommand{\alert}{\emph}
\DeclareMathOperator{\Hom}{Hom}

\DeclareMathOperator{\im}{im}

\newcommand{\mCF}{\mathcal{F}}

\DeclareMathOperator{\limit}{lim}
\DeclareMathOperator{\colimit}{colim}

\NewDocumentCommand{\colim}{e{_^}}{
  \mathbin{\mathop{\on{colim}}\displaylimits
    \IfValueT{#1}{_{#1}\,}
    \IfValueT{#2}{^{#2}\,}
  }
}

\DeclareMathOperator{\inj}{\hookrightarrow}
\DeclareMathOperator{\surj}{\twoheadrightarrow}
\DeclareMathOperator{\tensor}{\otimes}
\newcommand{\xrightleftarrows}[1]{\mathrel{\substack{\xrightarrow{#1} \\[-.9ex] \xleftarrow{#1}}}}
\DeclareMathOperator{\adjoint}{\xrightleftarrows{\rule{0.5cm}{0cm}} }

\newcommand{\on}{\operatorname}

\DeclareMathOperator{\oblv}{\tx{oblv}}

\DeclareMathOperator{\res}{\tx{res}}

\DeclareMathOperator{\Rep}{{Rep}}

\renewcommand{\mod}{\operatorname{-mod}}
\newcommand{\DMod}{\tx{DMod}}
\newcommand{\QCoh}{\tx{QCoh}}
\newcommand{\IndCoh}{\tx{IndCoh}}
\newcommand{\Coh}{\tx{Coh}}

\newcommand{\Ran}{\tx{Ran}}

\newcommand{\Vect}{\tx{Vect}}
\newcommand{\KL}{\tx{KL}}

\newcommand{\dR}{\tx{dR}}

\newcommand{\KLT}{\KL_{\kappa'}(T)}
\newcommand{\semiinf}{\frac{\infty}{2}}

\newcommand{\baseOp}{\mc{O}^{\tensor}}
\newcommand{\CM}{\ms{CM}}
\newcommand{\fSetnu}{{\ms{Comm}_\tx{nu}}}
\newcommand{\CMnu}{{\ms{CM}_\tx{nu}}}

\newcommand{\CommOperad}{\ms{Comm}}

\newcommand{\Setli}{\ms{Comm}}
\newcommand{\op}{\tx{op}}
\newcommand{\un}{\tx{un}}
\newcommand{\nonun}{\tx{nu}}

\newcommand{\weightLat}{\check{\Lambda}}
\newcommand{\coweightLat}{\Lambda}
\newcommand{\corootLat}{\Delta}
\newcommand{\rootLat}{\check{\Delta}}

\newcommand{\clambda}{{\check{\lambda}}}

\newcommand{\cmu}{\check{\mu}}
\newcommand{\crho}{\check{\rho}}

\newcommand{\ren}{\tx{ren}}
\newcommand{\SI}{\frac{\infty}{2}}

\newcommand{\DModRH}{\DMod^\tx{rh}}
\newcommand{\DModrhs}{\DMod^{\tx{hol}, \tx{RS}}}

\newcommand{\DGCat}{\tx{DGCat}}
\newcommand{\DGCatMarked}{\tx{DGCat}_*}
\newcommand{\VVLCatMarked}{\mathcal{CAT}_*}
\newcommand{\VVLCat}{\mathcal{CAT}}
\newcommand{\VVLSpc}{\mathcal{SPC}}

\newcommand{\ECOp}{\ms{Dbl}}
\newcommand{\ETOp}{\ms{Tpl}}

\newcommand{\cT}{\check{T}}
\newcommand{\PartitionMark}{\tx{Part}}
\newcommand{\ConstructibleMark}{\tx{constr}}

\newcommand{\deRhamGerbe}{{\mathscr{G}_\dR}}
\newcommand{\BettiGerbe}{\mathscr{G}_{\tx{Betti}}}

\newcommand{\LWeakMark}{\tx{L-weak}}
\newcommand{\WeakMark}{\tx{weak}}
\newcommand{\CoMark}{\tx{dual}}

\newcommand{\BettiU}{{U_\tx{Betti}}}
\newcommand{\TheGerbe}{\mathscr{G}}
\newcommand{\doubleLoop}{\mathbf{\weightLat}}

\newcommand{\Loop}{\mc{L}}
\newcommand{\Arc}{\mc{L}^+}
\newcommand{\TheFactAlg}{\Omega^\tx{Lus}_{\tx{KM}, \kappa}}

\newcommand{\StObj}{\Delta}
\newcommand{\CoStObj}{\nabla}

\newcommand{\modulePoints}{\overrightarrow{x}}

\newcommand{\kmgmod}{\hat{\mf{g}}\mod_{\kappa}}
\newcommand{\nkmgmod}{\hat{\mf{g}}\mod_{-\kappa}}

\newcommand{\AnomalyTerm}{E_\tx{anom}}
\newcommand{\DualTwisting}{[\check{\kappa}, 0]_{\leftarrow \omega^{\crho}}}

\newcommand{\Enalg}{ { \textup{Alg}^{\mathbb{E}_n} } } 
\newcommand{\Enminusonealg}{ { \textup{Alg}^{\mathbb{E}_{n-1}} } } 
\newcommand{\Eonealg}{ { \textup{Alg}^{\mathbb{E}_1} } } 
\newcommand{\lmod}{ { \textup{-mod} } }
\newcommand{\rmod}{ { \textup{-mod}^{r} } }
\newcommand{\Enminusonemod}{ { -\textup{mod}^{\mathbb{E}_{n-1}} } }
\newcommand{\Enmod}{ { -\textup{mod}^{\mathbb{E}_{n}} } }

\newcommand{\alglmodpair}{ { (\tx{Alg},\tx{Mod}^r) } }

\newcommand{\ratLat}{{\coweightLat^\sharp}}
\newcommand{\Waffextrat}{W^{\tx{aff}, \tx{ext}, \sharp}}
\newcommand{\Wfin}{W^\tx{fin}}
\newcommand{\RanBetti}{\Ran_{\tx{Betti}}}

\usepackage[style=alphabetic]{biblatex}
\DeclareFieldFormat{postnote}{#1}
\DeclareFieldFormat{multipostnote}{#1} 
\makeatletter
\newcommand*{\da@rightarrow}{\mathchar"0\hexnumber@\symAMSa 4B }
\newcommand*{\da@leftarrow}{\mathchar"0\hexnumber@\symAMSa 4C }
\newcommand*{\xdashedrightarrow}[2][]{%
  \mathrel{%
    \mathpalette{\da@xarrow{#1}{#2}{}\da@rightarrow{\,}{}}{}%
  }%
}
\newcommand{\xdashedleftarrow}[2][]{%
  \mathrel{%
    \mathpalette{\da@xarrow{#1}{#2}\da@leftarrow{}{}{\,}}{}%
  }%
}
\newcommand*{\da@xarrow}[7]{%
  \sbox0{$\ifx#7\scriptstyle\scriptscriptstyle\else\scriptstyle\fi#5#1#6\m@th$}%
  \sbox2{$\ifx#7\scriptstyle\scriptscriptstyle\else\scriptstyle\fi#5#2#6\m@th$}%
  \sbox4{$#7\dabar@\m@th$}%
  \dimen@=\wd0 %
  \ifdim\wd2 >\dimen@
    \dimen@=\wd2 %
  \fi
  \count@=2 %
  \def\da@bars{\dabar@\dabar@}%
  \@whiledim\count@\wd4<\dimen@\do{%
    \advance\count@\@ne
    \expandafter\def\expandafter\da@bars\expandafter{%
      \da@bars
      \dabar@ 
    }%
  }%
  \mathrel{#3}%
  \mathrel{%
    \mathop{\da@bars}\limits
    \ifx\\#1\\%
    \else
      _{\copy0}%
    \fi
    \ifx\\#2\\%
    \else
      ^{\copy2}%
    \fi
  }%
  \mathrel{#4}%
}
\makeatother

\makeatletter
\gdef\tshortstack{\@ifnextchar[\@tshortstack{\@tshortstack[c]}}
\gdef\@tshortstack[#1]{%
  \leavevmode
  \vtop\bgroup
    \baselineskip-\p@\lineskip 3\p@
    \let\mb@l\hss\let\mb@r\hss
    \expandafter\let\csname mb@#1\endcsname\relax
    \let\\\@stackcr
    \@ishortstack}
\makeatother

\title{An Extension of the Kazhdan--Lusztig Equivalence}
\author{Lin Chen and Yuchen Fu}
\date{}

\begin{document}

\maketitle

\begin{abstract}
We prove a tamely ramified version of the Kazhdan--Lusztig equivalence using factorization algebras.
More precisely, we establish an equivalence between the DG category of Iwahori--integrable affine Lie algebra representations and the DG category of representations of the ``mixed'' quantum group.
This confirms a conjecture by D. Gaitsgory in \cite{gaitsgory2021conjectural}.
\end{abstract}

\tableofcontents

\newpage

\section{Introduction}

\subsection{What is Done in This Paper}

One of the most beautiful and perhaps unexpected discoveries of representation theory in the 1980s was the link between affine Lie algebras and quantum groups. Let us briefly recall the story.

Let $G$ be a reductive group over $\mb{C}$, and $\mf{g}$ its Lie algebra. There are two ways one can \emph{deform} the category $\Rep(G)$:
\begin{itemize}
\item On one hand, we have the abelian category $(\hat{\mf{g}}_{-\kappa}\mod^{\Arc G})^{\heartsuit}$ of  smooth representations of $\hat{\mf{g}}_{-\kappa}$ which are integrable with respect to the arc group $\Arc G$. Here $-\kappa$ is an invariant symmetric bilinear form on the weight lattice, and $\hat{\mf{g}}_{-\kappa}$ is the corresponding affine Lie algebra. We shall denote this category by $\tx{KL}_{-\kappa}(G)^{\heartsuit}$;
\item On the other hand, we have $\Rep_{q^{-1}}(G)^{\heartsuit}$, the abelian category of representations of the Lusztig quantum group $U^\tx{Lus}_{q^{-1}}(G)$ at level\footnote{Our convention in this paper is that $\kappa$ will always correspond to $q$, and $-\kappa$ always to $q^{-1}$.} $q^{-1}$, which is a $\mb{C}^\times$-valued quadratic form on the weight lattice.
\end{itemize}

In a sequence of celebrated papers, Kazhdan and Lusztig \cite{kazhdanlusztig1, kazhdanlusztig2, kazhdanlusztig3, kazhdanlusztig4, lusztig1994monodromic} proved:
\begin{theorem}
\label{thm:KL}
Let $-\kappa$ be either an \emph{irrational} or \emph{negative rational} level, avoiding small roots of unity. Then there exists a \emph{braided monoidal} equivalence:
\[\KL_{-\kappa}(G)^{\heartsuit} \simeq \tx{Rep}_{q^{-1}}(G)^\heartsuit.\]
\end{theorem}

The situation is most interesting when the level is negative rational, at which the category captures the behavior of the modular representation theory of $G$ up to the first Frobenius kernel (c.f. \cite{ajs1994quantum}).

In \cite{gaitsgory2021conjectural}, Gaitsgory proposed an extension of the above equivalence to a class with weaker integrability conditions:
\begin{itemize}
\item Let $I \subseteq \Arc G$ denote the Iwahori subgroup of $\Loop G$. We also have $\kmgmod^I$, the \emph{renormalized} derived category of Iwahori-integrable smooth representations of $\hat{\mf{g}}$. This category usually goes under the name ``integral block of the affine BGG category $\mc{O}$'';
\item There is a ``mixed'' quantum group $U_q^\tx{mxd}(G)$, which a graded Hopf algebra that is Lusztig's quantum group on the $\mf{n}^+$ part, and Kac-De Concini quantum group on the $\mf{n}^-$ part. We let $\tx{Rep}_q^\tx{mxd}(G)$ denote the \emph{renormalized} derived category of graded representations of $U_q^\tx{mxd}(G)$ (see Section \ref{sect:prep-q} for precise definitions).
\end{itemize}
The proposed extension, which we prove in this paper, is that
\begin{theorem}[Main Theorem]
\label{thm:gI-mxd-equiv}
	Let $\kappa$ be a non-critical level avoiding small roots of unity, then there exists an equivalence of DG categories
	\[\mathtt{F}_\kappa: \kmgmod^I \simeq \tx{Rep}_q^\tx{mxd}(G).\]
\end{theorem}

The precise definition of ``avoiding small roots of unity'' is given in as follows:

\begin{condition}
\label{cond:avoid-small-roots-of-unity}
	For every simple factor $\mf{g}_i$ of $\mf{g}$, we have $\kappa|_{\mf{g}_i} = \frac{c_i - h_i^\vee}{2 h_i^\vee} \tx{Kil}_{\mf{g}_i}$ for some $c_i \in \mb{C}$, where $\tx{Kil}_{\mf{g}_i}$ is the Killing form on $\mf{g}_i$, $h_i^\vee$ is the dual Coxeter number. Let $h_i$ denote the Coxeter number of this simple factor, and $d_i$ its lacing number. We say $\kappa$ \emph{avoids small roots of unity} if each $c_i$ satisfies the following condition, which is the combination of Definition \ref{def:kappa-almost-admissible} and Convention \ref{conv:q-small-roots}:
	\begin{center}
	Either $c \not\in \mb{Q}$, or $c = \frac{p}{q}$ for $(p, q) = 1$, such that
	\[\begin{cases}
	p \ge \max(d_i + 1, h_i - 1) & \tx{if}~(q, d_i) = 1 \\
	p \ge \max(d_i + 1, h_i^\vee - 1) & \tx{else}
	\end{cases}\]
	\end{center}
\end{condition}
More explicitly, we have the following table of conditions on the minimal value of $p$ required for the equivalence to hold, in comparison with that required by Kazhdan and Lusztig:
\begin{table}[h!]
\centering
\begin{tabular}{ccc}
Simple Factor & Present Paper & Kazhdan--Lusztig\tablefootnote{For type $\tb{E}_8$ \cite[Page 413, Corollary 2]{kazhdanlusztig4} originally listed 26, but we were unable to verify this claim; to be safe we are following \cite{tanisaki2004character}, who listed 32 instead. We also note that values for non-simply-laced types for the original equivalence do not seem to exist in literature.} \\ \hline
$A_n$ & $\max(2, n)$  & 1 \\ \hline
$B_n$ & $\begin{cases}
\max(3, 2n - 2) & q~\tx{odd} \\
\max(3, 2n - 1) & q~\tx{even}
\end{cases}
$ & - \\ \hline
$C_n$ & $\begin{cases}
\max(3, n) & q~\tx{odd} \\
\max(3, 2n - 1) & q~\tx{even}
\end{cases}$ & - \\ \hline
$D_{2n}$ & $\max(2, 4n - 3)$ & 1 \\ \hline
$D_{2n+1}$ & $\max(2, 4n - 1)$ & 3 \\ \hline
$E_6$ & 11 & 14 \\ \hline
$E_7$ & 17 & 20 \\ \hline
$E_8$ & 29 & 32 \\ \hline
$F_4$ & $\begin{cases}
8 & q~\tx{odd} \\
11 & q~\tx{even}
\end{cases}
$ & - \\ \hline
$G_2$ & $\begin{cases}
4 & (3, q) = 1 \\
5 & \tx{else}
\end{cases}
$ & -
\end{tabular}
\label{tbl:minimal-value}
\caption{Minimal value for $p$}
\end{table}

Several remarks are in order:
\begin{itemize}
\item Our proof is independent of the original proof by Kazhdan and Lusztig (in fact, we do not yet know that it agrees with the original equivalence: see Expectation \ref{ex:backward-compatible} below);
\item The equivalence in Theorem \ref{thm:gI-mxd-equiv} is not $t$-exact on any rational level;
\item The RHS of the equivalence possesses a braided monoidal structure, hence the category $\kmgmod^I$ also admits a braided monoidal structure. To the best of our knowledge, this is not previously known. Because of the non-$t$-exactness of our equivalence, this structure is not easy to describe, and it is highly desirable to have a more concrete understanding of this structure. One potential way of doing this is by considering the \emph{semi-infinite} version of the LHS instead (c.f. Expectation \ref{ex:NKTO-equivalence} below);
\item The equivalence in Theorem~\ref{thm:gI-mxd-equiv} is expected to satisfy a number of compatibilities and extensions. However, some technical necessities are still undergoing construction. We hope to realize some of the expectations listed below in subsequent articles.
\end{itemize}

\subsection{Proof via Factorization Algebras}
Now we give an overview of our proof strategy.

The basic idea lies in the following categorical pattern. Suppose we have a lax monoidal functor $F: \mc{C} \to \mc{D}$ between two monoidal categories, then $F$ automatically factors as
\[\mc{C} \simeq \tb{1}_\mc{C}\mod(\mc{C}) \xrightarrow{F^\tx{enh}} F(\tb{1})\mod(\mc{D}) \xrightarrow{\oblv} \mc{D}\]
where $\tb{1}_{\mc{C}}$ is the monoidal unit. Often, $F^\tx{enh}$ has a better chance of being an equivalence than $F$ does. Similarly, if $\mc{M}$ is a $\mc{C}$-module category, and $\mc{N}$ a $\mc{D}$-module category, then a functor $F_\tx{Mod}$ compatible with these module structures would upgrade to a functor
\[\mc{M} \simeq \tb{1}_\mc{C}\mod(\mc{M}) \xrightarrow{F_\tx{Mod}^\tx{enh}} F(\tb{1}_\mc{C})\mod(\mc{N}) \xrightarrow{\oblv} \mc{N}\]

In fact, looking at Theorem \ref{thm:KL}, we have even more data: both sides are \emph{braided monoidal} categories. Since we work with DG categories, we should use the homotopical analogue of ``braided monoidal'', i.e. the little 2-disk operad (a.k.a.\ $\mb{E}_2$ operad). So, one could hope that there are some \emph{lax $\mb{E}_2$} functors that translate our categories into categories of $\mb{E}_2$-modules over the image of $\tb{1}$, then hope that these images are amenable to explicit identification.

First problem arises: $\mb{E}_2$, being a purely topological notion, interacts poorly with algebraic geometry and Lie theory. But not all hope is lost: it turns out we do have an algebro-geometric analogue of $\mb{E}_2$-structures, known as \emph{factorization structure}.

Intuitively, given a smooth complex curve $X$, a factorization structure on an object $\mc{C}$ (an algebra, a space, a category, etc.) is a sheaf of such objects over the moduli space of finite subsets of $X$, such that its fiber over any $k$-tuple of disjoint points $(x_1, \ldots, x_k)$ is $\mc{C}^{\tensor k}$. The behavior of fibers as we move from two disjoint points to a single point then encodes a structure similar to the monoidal structure of an $\mb{E}_2$-category.

A slightly more rigorous definition is as follows.
Let $\fSetnu$ denote the category of non-empty finite sets and surjective morphisms.
When we describe an object or structure (a prestack, a sheaf, etc) $A$ as \emph{factorizable}, such an object or structure will live over the \emph{Ran space} (the space of all finite subsets of our curve $X$) and will satisfy the \emph{factorization property} of being multiplicative on the disjoint locus. For instance, a factorizable D-module is the following data:
\begin{itemize}
\item For each $I \in \fSetnu$, a D-module $\mc{F}_I \in \DMod(X^I)$;
\item For any map $J \surj I$ and the resulting diagonal embedding $\Delta_{IJ}: X^I \to X^J$, an isomorphism
\[\Delta_{IJ}^!(\mc{F}_J) \simeq \mc{F}_I;\]
\item Define $(X^I \times X^J)_{\tx{disj}}$ as the open subscheme of $X^I \times X^J$ whose $S$-point is the subset of $(X^I \times X^J)(S)$ that are \emph{disjoint} pairs of points on $X(S)$. For any $I$ and $J$, the data also includes an isomorphism
\[j^!(\mc{F}_I \boxtimes \mc{F}_J) \simeq u^!(\mc{F}_{I \sqcup J}),\]
where $j: (X^I \times X^J)_{\tx{disj}} \inj X^I \times X^J$ is the embedding of the disjoint locus, and $u: (X^I \times X^J)_{\tx{disj}} \to X^{I \sqcup J}$ is the union map;
\item As well as all higher homotopy data.
\end{itemize}

The basic observation is that this description also makes sense for constructible sheaves. And, in the topological world, a fundamental result of Lurie tells us that $\mb{E}_2$-structures are equivalent to \emph{topological} factorization structures.
So this brings us to our basic plan of attack:

\[
\xymatrix{
\hat{\mf{g}}\mod^I_\kappa \ar@{-->}[rr] \ar_{}[dd] & & \Rep_q^{\tx{mxd}}(G) \ar^{}[d] \\
 & & \mb{E}_2\tx{-Module} \ar^{\tx{Higher Algebra}}[d] \\
\tx{Algebraic Factorization Module} \ar^{}_{\tx{Riemann-Hilbert}}[rr] & & \tx{Topological Factorization Module}
}
\]

\paragraph{The Betti Side}
Let us first look at RHS, where the above intuition holds almost literally. Let $\Rep_q(T)$ denote the representation category for the quantum torus. There exists an $\mb{E}_2$-algebra $\Omega_{q}^{\tx{Lus}}$ \emph{internal to} the braided monoidal category $\tx{Rep}_q(T)$, which can be obtained by taking Koszul dual of the Kac-De Concini quantum Hopf algebra $U_q^\tx{KD}(N^-)$. We will establish a non-commutative version of result from \cite{ben2010integral} and \cite{francis2013tangent} which says that
\begin{theorem}[Proposition \ref{prop:HH-is-E_n-mod}]
\label{thm:quantum-factMod-is-mxd}
We have an equivalence
\[\tx{Rep}_q^\tx{mxd}(G) \simeq \Omega_{q}^{\tx{Lus}}\mod^{\mb{E}_2}(\tx{Rep}_q(T))\]
for all values of $q$.
\end{theorem}

\begin{remark}
In some sense, we skipped a step: we directly matched mixed quantum groups with $\mb{E}_2$-modules without saying what the lax $\mb{E}_2$ functor involved is. Though we do not need this for the present paper, let us note that it is given by cohomology with respect to the Lusztig algebra $U_q^\tx{Lus}(N)$.
\end{remark}

An important caveat arises: the category that holds $\mb{E}_2$-modules is now $\Rep_q(T)$, which is itself only an \emph{$\mb{E}_2$-category}. Lurie established in \cite{HA} that to each $\mb{E}_2$-category one can attach a factorizable \emph{cosheaf} of categories; as we mentioned before, however, our desired notion of factorization category is a \emph{sheaf}. Since Verdier duality does not exist in the unstable setting, we will need to establish a \emph{categorical Verdier duality} to the following effect:
\begin{theorem}[Proposition \ref{prop:E2-to-fact-cat}, Proposition \ref{prop:E2Alg-Mod-to-FactAlg-Mod}]
\label{thm:categorical-Verdier}
To each $\mb{E}_2$ DG category $\mc{C}$ one can attach a topological factorization category $\tx{Fact}(\mc{C})$, such that the category of (non-unital) $\mb{E}_2$-algebras (resp.\ modules) in $\mc{C}$ correspond precisely to the category of factorization algebras (resp.\ modules).
\end{theorem}

We further show (Corollary \ref{cor:RepqT-is-GrcT}) that
\[\tx{Fact}(\Rep_q(T)) \simeq \tx{Shv}_{\check{\kappa}}(\tx{Gr}_{\check{T}}(\mb{C}))\]
here $\tx{Shv}_{\check{\kappa}}$ denote the category of $\Vect$-valued topological sheaves (twisted by a certain gerbe), and $\tx{Gr}_{\check{T}}$ is the \emph{affine Grassmannian}.

\begin{corollary}
\label{cor:E2-mod-is-factmod}
	There exists an equivalence of categories
\[\tx{FAlg}: \mb{E}_2\tx{-Alg}_{\nonun}(\Rep_q(T)) \to \tx{FactAlg}(\tx{Shv}_{\check{\kappa}}(\tx{Gr}_{\check{T}}(\mb{C})))\]
(here $\nonun$ indicates non-unitality), and, for every non-unital $\mb{E}_2$-algebra $A$ in $\Rep_q(T)$, an equivalence of categories
\[\tx{FMod}: A\mod^{\mb{E}_2}(\Rep_q(T)) \simeq \tx{FAlg}(A)\tx{-FactMod}(\tx{Shv}_{\check{\kappa}}(\tx{Gr}_{\check{T}}(\mb{C}))).\]
\end{corollary}

\begin{remark}
Since the pioneering work of \cite{finkelberg1994localization}, \cite{bezrukavnikov2006factorizable} it has been known that quantum groups (and more generally, sufficiently small Hopf algebras) and their representations can be encoded with topological factorization data. Because we work with derived objects, we opt to use the Koszul dual description of $\mb{E}_2$-algebras. In particular, our proof method is independent of \cite{bezrukavnikov2006factorizable}.
\end{remark}

\paragraph{The de Rham Side}

Now we look at Kac-Moody representations. First we have to answer the question: who is our acting category $\mc{C}$? The initial guess would be $\kmgmod^I$ itself, but this turns out to be not quite right, because $I$ as a group scheme \emph{does not factorize}; intuitively, this means this category doesn't know how to move around on the curve (see, however, Expectation \ref{ex:NKTO-equivalence}).

On the other hand, $\tx{KL}_\kappa(G)$, the derived version of $\tx{KL}_\kappa(G)^\heartsuit$, does. Indeed, we will show that $\tx{KL}_\kappa(G)$ admits an enhancement as a \emph{factorization category}, and $\kmgmod^I$ is a \emph{factorization module category} over $\KL_\kappa(G)$. This should be considered as encoding the fact that $\kmgmod^I$ is a module over $\KL_\kappa(G)$ under the Kazhdan--Lusztig fusion tensor product, an observation tracing back to \cite{finkelberg1993fusion}.

Let us look at our intuition again. Note that, in the very first step we used the (almost tautological) equivalence of $\mc{M} \simeq \tb{1}_\mc{C}\mod(\mc{M})$. It is, however, less obvious how to make sense of this step in algebraic geometry.

The additional input is that of a \emph{unital factorization structure}. In the same D-module example, this is following \emph{additional data}:
\begin{itemize}
\item For every \emph{injective map} $p: I \inj J$ and the corresponding projection $\pi_p: X^J \to X^I$, a map
\[\pi_p^!(A_I) \to A_J\]
in a way compatible with all of the data above;
\item Additionally, we also allow $I$ to be $\emptyset$ now. This gives, in particular, a map of factorization algebras $\omega_{\Ran} \to \mc{F}$. The object coming from $I = \emptyset$ will be called the \emph{vacuum};
\end{itemize}

\begin{example}
	In the case of the Kazhdan--Lusztig category, the vacuum is precisely the Kac-Moody vacuum vertex algebra $\mc{A}_{\mf{g}, \kappa}$, considered as an element of $\tx{KL}_\kappa(G)$.
\end{example}

\begin{remark}
We emphasize that the unital structure  \emph{fundamentally} alters the behavior of the categories involved. One may intuitively think of this as due to the fact that the unital Ran space has very different topology from the non-unital one.
\end{remark}

Using \emph{factorizable renormalized ind-coherent sheaves} (this is a new technical tool we introduce), we will show that category $\KL_\kappa(G)$ is indeed an \emph{unital} factorization category, carrying $\mc{A}_{\mf{g}, \kappa}$ as its vacuum. The unital structure, intuitively speaking, encodes how these vacuum modules can appear at all places and how they can interact with our given Kac-Moody representation (considered as fixed at a point on the curve).

Next we want to find a candidate for $F_\tx{Mod}$. A crucial observation, made at the abelian level by Beilinson and Drinfeld in \cite[Chapter 3.8]{beilinson2004chiral}, is that $C_*^{\SI}(\Loop \mf{n}, -)$, the \emph{semi-infinite cohomology} functor (introduced by \cite{feigin1984semi}) with respect to the nilpotent subalgebra $\Loop \mf{n}$, is a \emph{factorizable} functor. Let $\KLT$ denote the Kazhdan--Lusztig category for the torus, the general theory of \emph{unital factorization categories} then upgrades the functor
\[C_*^{\SI}(\Loop \mf{n}, -): \kmgmod^I \to \KLT\]
to a functor
\[C_*^{\SI}(\Loop \mf{n}, -)^\tx{enh}: \kmgmod^I \to C_*^{\SI}(\Loop \mf{n}, \mb{V}^0_\kappa)\tx{-FactMod}(\KLT)\]
where RHS is the category of factorization modules over the \emph{factorization algebra} $C_*^{\SI}(\Loop \mf{n}, \mb{V}^0_\kappa)$. Here $\kappa'$ denotes the \emph{shifted level} that semi-infinite cohomology maps into.

To produce actual sheaves, we will need to prove the following folklore:
\begin{proposition}[Proposition~\ref{prop:toric-FLE}, ``Toric FLE'']
There exists an equivalence of unital factorization categories
\[\tx{FLE}_T: \KLT \simeq \DMod_{\check{\kappa}'}(\tx{Gr}_{\check{T}}).\]
\end{proposition}
Let us set $\Omega_{\tx{KM}, \kappa} := \tx{FLE}_T \circ C_*^{\SI}(\Loop \mf{n}, \mb{V}^0_\kappa)$, so we can define
\[J_{\tx{KM}, \kappa} := \tx{FLE}_T \circ C_*^{\SI}(\Loop \mf{n}, -)^\tx{enh}: \kmgmod^I \to \Omega_{\tx{KM}, \kappa}\tx{-FactMod}(\DMod_{\check{\kappa}'}(\tx{Gr}_{\check{T}}))\]

Our main result is that
\begin{theorem}[Theorem \ref{thm:jacquet-is-equivalence}]
\label{thm:Jacquet-encodes-km}
At positive (or irrational) levels, $J_{\tx{KM}, \kappa}$ is an equivalence of categories.
\end{theorem}
We emphasize that this equivalence \emph{only} holds for positive levels.

\begin{remark}
The appearance of semi-infinite cohomology deserves some explanation. The category $\kmgmod$ is a category with a strong $\Loop G$ action, and our category $\kmgmod^I$ is its $I$-invariant category. Now, \cite[Section 6.2.1]{raskin2016chiral} showed that for any category $\mc{C}$ with a strong $\Loop G$-action, we have equivalence of categories
	\[\mc{C}_{\Loop N \Arc T} \simeq \mc{C}_I \simeq \mc{C}^{I} \simeq \mc{C}^{\Loop N \Arc T},\]
	here superscript (resp.\ subscript) denotes invariant (resp.\ coinvariant). Therefore, the \emph{affine} category $\mc{O}$ coincides with the coinvariant\footnote{The semi-infinite invariance and coinvariance categories both admit factorization structure, but they \emph{do not agree} as factorization categories. For reasons related to local Langlands, we think of the \emph{coinvariance} category as the right choice at positive levels.} \emph{semi-infinite} category $\mc{O}$ $(\kmgmod)_{\Loop N \Arc T}$, and the functor $J_{\tx{KM}, \kappa}$ should really be thought of as mapping out of the latter.
\end{remark}

In order to prove Theorem \ref{thm:Jacquet-encodes-km}, we observe that both sides have the structure of highest-weight categories, and proceed to match standard and costandard objects from both sides. For the same reason as in the previous remark, our choice of standard objects on the LHS will be the (dual objects of) \emph{semi-infinite} Verma modules, also known as Wakimoto modules.

It is not easy, however, to directly show that standard objects go to standard objects: in the setting of factorization modules, the standard objects are characterized by having simple $*$-fibers, which we have very little control of in general. To handle this $*$-fiber, we will establish a non-trivial \emph{local-global compatibility}, which is the following diagram (we use a slight variant in the actual proof):
\[\xymatrix{
\tx{KL}_\kappa(G)_x \ar^(0.3){C_*^{\semiinf}(\Loop \mf{n}, -)}[rr] \ar^{\tx{Loc}_x}[d] & & C_*^{\semiinf}(\Loop \mf{n}, \mb{V}_\kappa^0)\tx{-FactMod}(\KLT_x) \ar^{\tx{Loc}^\tx{ch}_x}[d] \\
\DMod_\kappa(\tx{Bun}_G(\mb{P}^1)) \ar^{\tx{CT}_*[-d]}[rr] & & \DMod_{\kappa'}(\tx{Bun}_T(\mb{P}^1))
}\]
here we consider $x$ as a fixed point on $\mb{P}^1$, so that LHS can be seen as the usual localization procedure from Kac-Moody representation to D-modules on $\tx{Bun}_G$. The RHS is \emph{chiral localization}, which roughly means doing localization at every point \emph{at the same time}, with the vacuum $C_*^{\semiinf}(\Loop \mf{n}, \mb{V}_\kappa^0)$ filling in every other point on the curve. Finally, $\tx{CT}_*$ is the \emph{geometric constant term functor}, and $d$ is a cohomological shift. This allows us to relate the $*$-fibers to \emph{conformal blocks}, from which the result follows.

To recap, we have:

\[
\xymatrix{
\hat{\mf{g}}\mod^I_\kappa \ar@{-->}[rr] \ar_{\tx{Geometry}}[d] & & \Rep_q^{\tx{mxd}}(G) \ar^{\tx{Koszul Duality}}[d] \\
\tx{Algebraic}~\mc{A}_{\mf{g}, \kappa}\tx{-Factorization Module} \ar_{\tx{Semi-infinite Cohomology}}[d] & & \tx{Twisted}~\mb{E}_2\tx{-Module} \ar^{\tx{Higher Algebra}}[d] \\
\tx{Algebraic}~\Omega_{\tx{KM}, \kappa}\tx{-Factorization Module} \ar^{}_{\tx{Riemann-Hilbert}}[rr] & & \tx{Topological}~\tx{FAlg}(\Omega^\tx{Lus}_{q})\tx{-Factorization Module}
}
\]

To link the two worlds, we show that
\begin{proposition}[Corollary \ref{cor:riemann-hilbert}, Corollary \ref{cor:RH-for-fact-mod}]
\label{prop:fact-mods-match-up}
There exists a factorizable system of equivalence of categories
\[\tx{RH}: \DMod_{\check{\kappa}'}(\tx{Gr}_{\check{T}}) \simeq \tx{Shv}_{\check{\kappa}'}(\tx{Gr}_{\check{T}}(\mb{C})).\]
The functor $\tx{RH}$ preserves factorization algebras, and, for every factorization algebra $A$, induces an equivalence of categories
\[A\tx{-FactMod}(\KLT) \simeq \tx{RH}(A)\tx{-FactMod}(\tx{Shv}_{\check{\kappa}'}(\tx{Gr}_{\check{T}}(\mb{C}))).\]
\end{proposition}

The last ingredient of the proof is that
\begin{proposition}[Proposition \ref{prop:KM-E2-Lus-match}]
	\label{prop:algebras-match}
	For $\kappa$ positive and avoiding small roots of unity, we have
	\[\tx{RH}(\Omega_{\tx{KM}, \kappa}) \simeq \tx{FAlg}(\Omega^\tx{Lus}_{q}).\]
\end{proposition}

The proof of this statement relies on the fact that our factorization algebra is in fact a \emph{perverse sheave} that admit a combinatorial description (obtained in \cite{gaitsgory2021factorization}), so that to match the two sides it suffices to compute their $!$- and $*$- fibers up to some low cohomological degrees.

Then, Theorem \ref{thm:Jacquet-encodes-km}, Theorem \ref{thm:quantum-factMod-is-mxd}, Corollary \ref{cor:E2-mod-is-factmod}, Proposition \ref{prop:fact-mods-match-up} and Proposition \ref{prop:algebras-match} together imply Theorem \ref{thm:gI-mxd-equiv} when $\kappa$ is positive (or irrational).

\paragraph{Negative Levels}
For every DG category $\mc{C}$, we have its dual category $\mc{C}^\vee := \tx{Fun}_{\tx{cont}}(\mc{C}, \Vect)$ consisting of colimit-preserving functors from $\mc{C}$ to $\Vect$. On the Kac-Moody side, we have a \emph{categorical} (a.k.a.\ cohomological) duality (c.f. Definition \ref{def:km-categorical-duality}):
\[\hat{\mf{g}}_{-\kappa}\mod^I \simeq (\kmgmod^I)^\vee\]
On the quantum side, we have the \emph{contragredient duality} (c.f. Definition \ref{def:mixed-contragredient-duality}):
\[\Rep^\tx{mxd}_{q^{-1}}(G) \simeq (\Rep^\tx{mxd}_q(G))^\vee.\]
Thus for $-\kappa$ negative, there also exists an equivalence of categories
\[\mathtt{F}_{-\kappa}: \nkmgmod^I \xrightarrow{\simeq} \tx{Rep}_{q^{-1}}^{\tx{mxd}}(G)\]
defined as the \emph{dual functor} of $\mathtt{F}_{\kappa}$.

As mentioned in the beginning, it is an important --- and non-trivial --- question whether our equivalence agrees with that of Kazhdan and Lusztig.
\begin{expectation}
\label{ex:backward-compatible}
We expect that $\mathtt{F}_{-\kappa}$ induces a $t$-exact equivalence on the $\Arc G$-invariant category at the negative level, which is the functor in Theorem \ref{thm:KL}.
\end{expectation}

\begin{remark}
The Jacquet functor $J_{\tx{KM}, \kappa}$ described in this paper exists for negative levels, but is \emph{not} an equivalence there.
\end{remark}

\subsection{Technical Aspects}

The bulk of the present paper's content is \emph{theory-building}: indeed, once all the tools have been set up, it takes only a few pages' calculation to get to the desired result. So here we make a list of main technical aspects that constitute the backbone of this work, and offer some perspectives on their necessity:
\begin{enumerate}
\item We further extended Raskin's theory of \emph{renormalized ind-coherent sheaves} (\cite{raskin2020homological}) to provide it with extra functoriality and more base-change. Other than being of independent algebro-geometric interest, we emphasize that this is what makes everything work (factorizability of semi-infinite cohomology and local-global comparison, among other things);
\item As far as we are aware of, the present work produces the first instance of a \emph{factorizable system of ind-coherent sheaves}. Let us emphasize that even the plausibility of such an object is not a priori obvious, because ind-coherent sheaves do not form a sheaf of categories: indeed, only after renormalizing\footnote{We note this is a different renormalization from the one on Kac-Moody representations.} the presheaf of categories do we obtain such a system. We expect this technique to be more widely applicable, for instance in the study of \emph{factorizable weak group actions};
\item We fully implemented the \emph{categorical Verdier duality} necessary for converting twisted Hopf algebras into twisted factorizable sheaves. This has been a folklore in the community for over 10 years. We also note that this method is fully homotopic in nature and completely bypasses\footnote{On the other hand, one advantage of directly working with Hopf algebras is the ability to stay with the heart of the $t$-structure, which reduces homotopy coherence to a finite level.} the hyperplane arrangement argument in \cite{bezrukavnikov2006factorizable}.
\end{enumerate}

\subsection{Additional Remarks}

\paragraph{Hecke Equivariance}
Recall we have the Iwahori Hecke category $\DMod_\kappa(I \backslash \Loop G / I)$ which naturally acts on $\kmgmod^I$ via convolution. Using the equivalence we establish, it also acts on $\Rep_q^\tx{mxd}(G)$.

Let $H$ be the metaplectic dual of $G$ (c.f. Definition \ref{def:metaplectic-dual-group}). We have the following metaplectic extension of R. Bezrukavnikov's result in \cite{bezrukavnikov2016two} to rational levels (work-in-progress of Dhillon, Yun and Zhu):
\[\DMod_\kappa(I \backslash \Loop G / I) \simeq \IndCoh(\left(\widetilde{\mc{N}}_{H} 
\times_{\mf{h}} \widetilde{\mc{N}}_{H}\right)/H)\]
(Note that the RHS is not quite right: we should have a \emph{singular support condition} on it. Here and below we will momentarily ignore this problem.)
\begin{expectation}
	We expect $\IndCoh(\left(\widetilde{\mc{N}}_{H} 
\times_{\mf{h}} \widetilde{\mc{N}}_{H}\right)/H)$ to act naturally on $\Rep_q^\tx{mxd}(G)$, and that $\ms{F}_{\kappa}$ is equivariant with respect to this action via the equivalence above.
\end{expectation}

\begin{remark}
In \cite[Section 7.2.3]{gaitsgory2021conjectural} Gaitsgory conjectured an action of $\QCoh(\mf{n}_H / B_H)$ on $\Rep_q^\tx{mxd}(G)_\ren$. Within the paradigm above it would be provided by the map
\[\QCoh(\mf{n}_H / B_H) \xrightarrow{\Upsilon} \IndCoh(\mf{n}_H / B_H) \simeq \IndCoh(\widetilde{\mc{N}}_H / H) \xrightarrow{\Delta_*} \IndCoh(\left(\widetilde{\mc{N}}_{H} 
\times_{\mf{h}} \widetilde{\mc{N}}_{H}\right)/H).\]
A description of this action is now known and will be included in a subsequent article.
\end{remark}

\paragraph{Factorizable Equivalence}
By \cite{HA}, the category $\Rep_q^\tx{mxd}(G)$ is an $\mb{E}_2$ category itself, so Theorem \ref{thm:categorical-Verdier} allows us to attach to it a topological factorization category $\tx{Fact}(\Rep_q^\tx{mxd}(G))$. On the other hand, recall that we have $\kmgmod^I \simeq (\kmgmod)_{\Loop N \Arc T}$ and the RHS has an algebraic factorization category structure.
One is thus naturally led to wonder about the existence of a factorizable system of equivalences
\[\mathtt{F}_\kappa^\tx{Fact}: (\kmgmod)_{\Loop N \Arc T} \simeq \tx{Fact}(\Rep_q^\tx{mxd}(G))\]
that extends the functor $\mathtt{F}_\kappa$ at a point.

\begin{expectation}
\label{ex:NKTO-equivalence}
We expect some version of this equivalence to exist.
\end{expectation}


\subsection{Fundamental Local Equivalence}
The result of this paper fits naturally into the local geometric Langlands program, initiated in Frenkel and Gaitsgory in \cite{frenkel2006local}. We briefly summarize the ideas here.

Let $\kappa$ denote the quantum parameter (for $G$ a simple group it is just an element of $\mb{P}^1(\mb{C})$), and let $\check{\kappa}$ denote the Langlands dual parameter. Let $\check{G}$ denote the Langlands dual group of $G$. An expected output of the 2-categorical formalism of local geometric Langlands is:
\begin{conjecture}[Spherical Fundamental Local Equivalence (FLE)]
\label{conj:spherical-FLE}
For any level $\kappa$, We have an equivalence of (unital) factorization categories
\[\KL_{\kappa}(G) \simeq \tx{Whit}(\DMod_{\check{\kappa}}(\tx{Gr}_{\check{G}}));\]
here $\tx{KL}_\kappa(G)$ is as defined above, $\tx{Gr}_{\check{G}}$ is the affine Grassmannian, and the RHS (known as the spherical Whittaker category) is the strong $\Loop N$-invariant of $\DMod_{\check{\kappa}}(\tx{Gr}_{\check{G}})$ with respect to a non-dgenerate character $\chi$, $N$ being the unipotent subgroup of a fixed Borel of $G$.
\end{conjecture}
We emphasize that the word ``factorization'' is what makes the statement most interesting, as being a factorization category is akin to (though inequivalent to) having a braided monoidal structure on the category.

\begin{remark}
	The global geometric Langlands program, which aims to prove the following equivalence of categories
	\[\DMod(\tx{Bun}_G) \simeq^? \IndCoh_{\tx{Nilp}}(\tx{LocSys}_{\check{G}}),\]
	has the following guiding philosophy: the LHS and RHS are respectively fairly close to the factorization homology of $\tx{Whit}(\DMod_{\kappa_\tx{crit}}(\tx{Gr}_G))$ and $\tx{KL}_{\infty}(\check{G}) \simeq \Rep(\check{G})$ over the curve (here $\kappa_\tx{crit}$ is the critical level). Because of this, FLE is considered a central component towards establishing (classical and quantum) global geometric Langlands. For a more precise formulation of the role Conjecture \ref{conj:spherical-FLE} plays in the global setting, see \cite{gaitsgory2015outline}.
\end{remark}

\begin{remark}
While the factorizable version of this conjecture remains quite open (except for when $G = T$ is a torus, which we establish here in Proposition \ref{prop:toric-FLE}), its fiberwise statement (i.e. with the word ``factorization'' removed) is more tractable. For instance, the cases of $\kappa = \kappa_\tx{crit}$ and $\kappa = \infty$ were known by \cite{bezrukavnikov2019iwahori} and \cite{frenkel2009local} respectively, and at almost all levels this is now known by the work of \cite{campbell2019fundamental} via affine Soergel bimodules.
\end{remark}

One strategy proposed by Gaitsgory for proving the FLE (factorizably) is by linking both sides to quantum groups:
\[\xymatrix{
\KL_{\kappa}(G) \ar@{-->}^{\tx{FLE}}[rr] \ar@{-->}^{\simeq}_{\tx{KL}}[rd] & &  \tx{Whit}(\DMod_{\check{\kappa}}(\tx{Gr}_{\check{G}})) \\
 & \tx{Rep}_{q}(G) \ar@{-->}^{\simeq}_{\tx{CS}}[ru]
}\]
where on the left we have (a conjectural factorizable version of) the Kazhdan--Lusztig equivalence, and on the right we have the (conjectural) \emph{quantum Casselman-Shalika formula}. It is expected that both (KL) and (CS) can be proven factorizably.

\begin{remark}
	We warn that it is a non-trivial matter to convert the braided monoidal structure of the classical Kazhdan--Lusztig equivalence to a form that is usable for geometric Langlands' purpose, as the braided monoidal structure there is not local (or even algebraic) by nature.
	
	Instead the logic goes the reverse order: we expect the methods established in the present paper will eventually yield a proof of (KL) as an factorizable equivalence, from which we hope to extract a \emph{new proof} of the Kazhdan--Lusztig braided monoidal equivalence.
\end{remark}

For the Iwahori subgroup, the counterpart of FLE is
\begin{proposition}[Iwahori / Tamely Ramified FLE (IFLE)]
For any level $\kappa$ we have an equivalence of categories
\[\hat{\mf{g}}_{\kappa}\mod^{I} \simeq \tx{Whit}(\DMod_{\check{\kappa}}(\tx{Fl}_{\check{G}}));\]
here $\tx{Fl}_{\check{G}}$ is the affine Flag variety.
\end{proposition}
Note that the word ``factorization'' has disappeared from the formulation, since the Iwahori group does not have a factorizable analogue (see, however, Expectation \ref{ex:NKTO-equivalence} above). The corresponding picture involving quantum groups becomes
\[\xymatrix{
\hat{\mf{g}}_{\kappa}\mod^{I} \ar^{\tx{IFLE}}[rr] \ar^{\simeq}_{\tx{Theorem~\ref{thm:gI-mxd-equiv}}}[rd] & &  \tx{Whit}(\DMod_{\check{\kappa}}(\tx{Fl}_{\check{G}})) \\
 & \tx{Rep}^\tx{mxd}_{q}(G) \ar^{\simeq}_{\tx{CS}_\tx{Fl}}[ru]
}\]
where now all arrows are actual theorems: the LHS is the main result of this paper, and RHS can be obtained by combining \cite{yang2021twisted} with Theorem \ref{thm:quantum-factMod-is-mxd}, Corollary \ref{cor:E2-mod-is-factmod} and Proposition \ref{prop:fact-mods-match-up}.
We note that the top arrow is also independently established in \cite{campbell2019fundamental}; however, if we interpret the top arrow as the equivalence in \emph{loc.cit.}, then the commutativity of the triangle is not yet known.

\subsection{Acknowledgment}
Our deepest gratitude is to Dennis Gaitsgory, whose visions and ideas can be found shining everywhere in this project. We are also grateful for his suggesting this project to us in 2018.

Our profound admiration is to Jacob Lurie, whose bravery and dedication provided the technical bedrock on which all aspects of this paper stand.

We thank Gurbir Dhillon, John Francis, Ben Knudsen, Sam Raskin and Nick Rozenblyum for valuable communications and suggestions regarding various parts of the project. We thank Kevin Lin, James Tao, David Yang, Ruotao Yang, Yifei Zhao for helpful conversations related to this paper.

The work of the first author was supported by a grant from the Simons Foundation (816048, LC) during his stay at the Institute for Advanced Study.

\section{Proof Outline}

For reader's convenience, we will first present the main proof in this section, after stating all of the definitions and background results necessary to process it. Proofs of individual statements can be found in later sections.

\subsection{Conventions}

\paragraph{Base Field}
We will work over $k = \mb{C}$. Most results unrelated to Riemann-Hilbert hold true for any algebraically closed field of characteristic $0$, but we do not need this generality.

\paragraph{Category Theory}
We use the theory of $\infty$-categories set up in \cite{HTT}. At several points we will make use of the theory of $(\infty, 2)$-categories; for some foundational setup, we refer our readers to the appendix of \cite{GR-DAG1}. By a \emph{category} we always mean an $(\infty, 1)$-category unless specified otherwise; similarly by a \emph{2-category} we mean an $(\infty, 2)$-category. We let
\begin{itemize}
\item $\tx{Spc}$ denote the category of $\infty$-groupoids;
\item $1\tx{-Cat}$ denote the category of all (small) $(\infty, 1)$-categories, and $1\tx{-Cat}^{2\tx{-Cat}}$ denote the enhancement of it into a 2-category;
\item $\tx{Pr}^\tx{L, St}$ denote the (symmetric monoidal) category of stable, presentable categories and continuous (colimit-preserving) functors;
\item $\tx{Vect}$ denote the category of chain complexes of vector spaces over $k$;
\item $\tx{DGCat}$ denote the (symmetric monoidal) category $\tx{Vect}\mod(\tx{Pr}^\tx{L, St})$, and $\DGCat_*$ denote the (symmetric monoidal) category of pairs $(c, \mc{C})$, where $\mc{C} \in \DGCat$ and $c \in \mc{C}$.
\end{itemize}

\paragraph{Higher Algebra}
We make frequent use of terminologies and facts about homotopical algebra as set up in \cite{HA}.

\paragraph{Derived Categories}
The derived category of an Grothendieck abelian category $\mc{A}$ for us will refer to the $\infty$-category of \emph{unbounded} chain complexes in $\mc{A}$ localized with respect to quasi-isomorphisms. We refer readers to \cite[Chapter 1]{HA} for a detailed discussion. We adopt the \emph{cohomological} grading convention.

When we write $\Hom(X, Y)$ for objects in a DG category, we always mean the entire mapping chain complex; we write $\tx{Ext}^n$ when we want to refer to its $n$th cohomology.

\paragraph{Algebraic Geometry}
For foundation of derived algebraic geometry, we refer readers to \cite{GR-DAG1} and \cite{lurie2018spectral}. Our terminologies and conventions will follow that of \cite{GR-DAG1}.

Let $\tx{CAlg}^{\le 0}$ be the category of connective commutative algebra objects in $\tx{Vect}$. We let
\begin{itemize}
\item $\tx{Sch}^\tx{aff} := (\tx{CAlg}^{\le 0})^\op$ denote the category of affine schemes;
\item $\tx{PreStk}$ denote the $(\infty, 1)$-category of $\tx{Spc}$-valued presheaves over $\tx{Sch}^\tx{aff}$;
\item $\tx{LaxPreStk}$ denote the $(\infty, 2)$-category of $1\tx{-Cat}$-valued presheaves over $\tx{Sch}^{\tx{aff}}$.
\end{itemize}

Given a prestack $\mc{Y}$, we use $\mc{Y}_\dR$ to denote the de Rham prestack defined by $\mc{Y}_\dR(S) := \mc{Y}(S^\tx{red, cl})$, where $S^\tx{red, cl}$ is the corresponding classical reduced affine scheme.

\paragraph{Factorization Structure}
In the introduction we gave an informal definition of what a factorization structure is; an equivalent (see Section \ref{sect:factspc} for proof) way of formulating it is as follows.

We consider $\Ran_{\nonun, \dR}$, the non-unital \emph{Ran space} of $X_\dR$, which is the prestack whose space of $S$-points is the space of finite subsets of $X_\dR(S)$.
A non-unital factorization space is then an object $\mc{Y}_{\Ran_\dR} \in \tx{PreStk}$ over $\Ran_{\nonun, \dR}$ along with the following extra data: for every point $\{x_1, \ldots, x_n\} \in \Ran_{\nonun, \dR}(S)$ that lies in the disjoint locus\footnote{For $S = \mb{C}$ the meaning of ``disjoint'' is clear; for general $S$ see Section \ref{sect:factspc}. For the sake of intuition it does little harm to only consider the case of $\mb{C}$-points.}, we require an isomorphism
\[\mc{Y}_{\Ran_\dR} |_{\{x_1, \ldots, x_n\}} \simeq \mc{Y}_{\Ran_\dR}|_{x_1} \times_S \mc{Y}_{\Ran_\dR}|_{x_1} \times_S \ldots \times_S \mc{Y}_{\Ran_\dR}|_{x_n}.\]

If one replaces ``prestack over $\Ran_{\nonun, \dR}$'' with ``D-modules over $\Ran_{\nonun, \dR}$'', one obtains the notion of non-unital factorization \emph{algebras}. If one instead uses ``sheaves of categories over $\Ran_{\nonun, \dR}$,'' one obtains the notion of non-unital factorization \emph{categories}. By considering sheaves of \emph{pointed} categories, one obtains the notion of non-unital factorization algebras \emph{internal} to a factorization category. For a fully rigorous definition, see Section \ref{sect:factspc}.

Let $x_0 \in X_\dR(\mb{C})$ be a fixed $\mb{C}$-point. Consider $\Ran_{\nonun, x_0, \dR}$, the prestack whose space of $S$-points is the space of finite subsets of $X_\dR(S)$ containing $S \to \tx{Spec}(\mb{C}) \xrightarrow{x_0} X_\dR$. Fix a non-unital factorization space $\mc{A}$; a non-unital factorization module space is a prestack $\mc{Z}$ over $\Ran_{\nonun, x_0, \dR}$ along with the following extra data: for every $\{x_0, x_1, \ldots, x_n\} \in \Ran_{\nonun, x_0, \dR}(S)$ that lies in the disjoint locus, we require an isomorphism
\[\mc{Z}|_{\{x_0, \ldots, x_n\}} \simeq \mc{Z}|_{x_0} \times_S \mc{A}|_{x_1} \times_S \ldots \times_S \mc{A}|_{x_n};\]
The notion of factorization module categories and (internal) factorization module algebras can be defined analogously; we can also replace $x_0$ by a finite set of fixed points.

A word on unitality: as we mentioned in the introduction, the notion of unital structure in the factorization setting plays many crucial technical roles; on the other hand, it does render some of the definitions more cumbersome. We thus suggest readers to ignore all things related to unital structure on the first pass by considering it as an analogue of the unital structure in the $\mb{E}_2$-world, and refer to the appendix for a rigorous treatment if needs be.

\paragraph{Reductive Group}
For $G$ a connected reductive group, we write $\mc{L}G$ (resp. $\mc{L}^+ G$) for the prestacks given by
\[\mc{L}G(A) := G(A (\!(t)\!)) \hspace{1em} \mc{L}^+G(A) := G(A[\![t]\!]).\]
We refer to them as the loop (resp. arc) group of $G$. It is well-known that $\mc{L}G$ is an indscheme and $\mc{L}^+ G$ a scheme, albeit both of infinite type. We let $\tx{Gr}_G$ denote the fppf sheaf quotient $\mc{L} G / \mc{L}^+ G$.

There are well-known factorizable spaces over $\Ran_\dR$, which we'll denote by
\[[\mc{L} G]_{\Ran_\dR}, [\mc{L}^+ G]_{\Ran_\dR}, [\tx{Gr}_G]_{\Ran_\dR}\]
respectively, whose fiber at any $\mb{C}$-point on $X$ recovers the prestacks defined above.

We fix a standard Borel $B$ and a maximal torus $T$ of $G$. We let $N$ denote the unipotent part of $B$. We let $B^-$ and $N^-$ denote the opposite Borel (resp. unipotent). We let $\mf{g}$ (resp. $\mf{b}, \mf{n}, \mf{t}$) refer to the Lie algebras.

We let $(\weightLat, \rootLat, \coweightLat, \corootLat)$ denote the root datum corresponding to $(G, T)$. Here $\coweightLat$ denote the coweight lattice of $G$, $\weightLat$ its weight lattice, $\corootLat$ its coroot lattice and $\rootLat$ its root lattice.

The superscript $<0$ (resp.\ $\le 0, >0, \ge 0$) on the (co)weight lattice will always refer to the submonoid generated by negative (resp.\ non-positive, etc.) (co)\emph{roots}; in other words, the inequality signs refer to the standard Bruhat order. The superscript $+$ and $-$ will mean \emph{dominant} and \emph{anti-dominant} respectively.

We let $W_\tx{fin}$ denote the finite Weyl group of $(G, T)$, $W_\tx{aff}$ the affine Weyl group, and $W_{\tx{aff,ext}}$ the extended affine Weyl group of $G$.

We let $I$ be the Iwahori subgroup, defined as the preimage of $B$ under $\mc{L}^+(G) \xrightarrow{t = 0} G$.

\paragraph{Global Structure}
Let $X$ be a smooth complete algebraic curve over $\mb{C}$. For each reductive group $G$, we let $\tx{Bun}_G(X)$ denote the moduli stack of $G$-bundles on $X$; this stack is smooth and classical, of dimension $\dim(\mf{g})(g(X) - 1)$, where $g(X)$ is the genus of $X$. When $X$ is clear from the context, we will denote it simply by $\tx{Bun}_G$.
We fix $\omega_X^{1/2}$, a square root of the canonical bundle on $X$.

\paragraph{Quantum Parameters}
A \emph{quantum parameter} (introduced in \cite{zhao2020tame}) is, by definition, a pair $[\kappa, E]$ where $\kappa$ is a $W$-invariant symmetric bilinear form on $\mf{t}$ and $E$ is an extension of sheaf of abelian groups of the form
\[0 \to \omega_X \to E \to \mf{z} \tensor \mc{O}_X \to 0.\]

\begin{remark}
The parameter $E$ plays a technically significant but conceptually minor role in the story. For reading of the main proof it is safe to ignore it on the first pass.
\end{remark}

Let $\mf{g}_1, \ldots, \mf{g}_r$ be the simple factors of $\mf{g}$, then $\kappa$ can always be written uniquely as $\sum_{i = 1}^{r} [c_i] \cdot \tx{Kil}_i + \mu_\mf{z}$, where $\mu_\mf{z} \in \tx{Sym}^2(\mf{z}^*)$, $\tx{Kil}_i$ is the Killing form on $\mf{g}_i$ and we set $[c_i] := \frac{c_i - h_i^\vee}{2 h_i^\vee}$ for $c_i \in \mb{C}$, $h_i^\vee$ being the dual Coxeter number of $\mf{g}_i$.
We write $\kappa_{\tx{crit}, G} := -\frac{1}{2}\tx{Kil}_\mf{g}$ and refer to $[\kappa_{\tx{crit}, G}, 0]$ as the critical level; note that it corresponds to $c_i = 0$ for each $i$. We refer to levels where each $c_i \in \mb{Q}^{< 0}$ (resp. $\mb{Q}^{> 0}, \mb{C} \setminus \mb{Q}, \mb{Z} \setminus \{0\}$) as \emph{negative} (resp. \emph{positive}, \emph{irrational}, \emph{integral}).

\begin{convention}[Critical Shift]
\label{conv:critical-shift}
The superscript $\kappa$ (without the square bracket) will always denote the \emph{critically shifted} level $[\kappa + \kappa_{\tx{crit}, G}, 0]$. For a bilinear form $\kappa$, we let $-\kappa$ denote the bilinear form reflected along $\kappa_{\tx{crit}, G}$. In this article, $\kappa$ will always be a \emph{positive} level, and $-\kappa$ always a negative level.
\end{convention}

\paragraph{Langlands Duality}
Langlands duality attaches to every $(G, T)$ a dual group $(\check{G}, \check{T})$ whose root datum is $(\coweightLat, \corootLat, \weightLat, \rootLat)$. We let $[\check{\kappa}, \check{E}]$ denote the \emph{Langlands dual} parameter: on each simple factor it sends $c_i$ to $-\frac{1}{d_i c_i}$ where $d_i$ is the lacing number of $\mf{g}_i$, and $\check{E}$ is the extension induced by $\mf{z} \simeq \check{\mf{z}}$.

\subsection{Factorizable Kazhdan--Lusztig Category}

At this point we will start to freely use the theory of factorizable renormalized ind-coherent sheaves developed in Section \ref{sect:appendix-indcoh}. The reader is suggested to take on faith the following claim, which summarizes the main output of said theory, and refer to the appendix when necessary.

\begin{claim}
For a class of sufficiently nice prestacks $\mc{Y}$, there exists a compactly generated DG category $\IndCoh_{*, \tx{ren}}(\mc{Y})$ equipped with a (not necessarily left complete) $t$-structure, such that the bounded below category coincides with that of $\QCoh(\mc{Y})$. We let $\IndCoh^{!, \tx{ren}}(\mc{Y})$ denote its dual category.

Let $B_\dR$ be the de Rham prestack associated with a laft prestack $B$ (we shall mostly consider $B = \Ran$). Suppose $\mc{Y} \to B_\dR$ is a prestack whose fiber at every $x: \tx{pt} \to B_\dR$ is reasonably nice. Then there exists a \emph{crystal of categories} $\mc{I}nd\mc{C}oh_{*, \tx{ren}}(\mc{Y})$ (and a dual version $\mc{I}nd\mc{C}oh^{!, \tx{ren}}(\mc{Y})$) over $B_\dR$, whose fiber at every $x$ is given by the renormalized ind-coherent sheaf category $\IndCoh_{*, \tx{ren}}(Y_x)$ (resp.\ $\IndCoh^{!, \tx{ren}}(Y_x)$).
\end{claim}

To each quantum parameter $[\kappa, E]$ attached to a reductive group $H$ (in our use cases we only consider $H = G$ and $H = T$) we can associate an (unital) factorization line bundle $\mathscr{L}_{[\kappa, E]}$ on 
\[[\widehat{\tx{Hecke}}_H]_{\Ran_\dR} := [\Arc H]_{\Ran_\dR} \backslash ([\Loop H]_{\Ran_\dR})_{[\Arc H]_{\Ran_\dR}}^\wedge / [\Arc H]_{\Ran_\dR}\]
which admits a multiplicative structure under convolution (c.f. Proposition \ref{prop:yifei-twisting-conversion}). By restricting, we can also consider $\mathscr{L}_{[\kappa, E]}$ as a factorization line bundle for $[\widehat{\tx{Hecke}}_{H'}]_{\Ran_\dR}$, $H' \subseteq H$ any affine subgroup.

We have the following key construction, whose details are deferred to Section \ref{sect:fact-kl-category}.

Let $H$ be a finite type affine algebraic group, and let $\mathscr{L}$ be a unital factorization line bundle on $[\widehat{\tx{Hecke}}_H]_{\Ran_\dR}$.
Let $\mc{R}ep(\Arc H)$ denote $\mc{I}nd\mc{C}oh^{!, \tx{ren}}(\mb{B} \Arc H)$; this category admits a unital factorization category structure. The category $\mc{I}nd\mc{C}oh^{!, \tx{ren}}([\widehat{\tx{Hecke}}_H]_{\Ran_\dR})$ acts on it in a way compatible with the factorization structure. We set
\[\mc{KL}(H)_{\mathscr{L}} := \mathscr{L}\mod(\mc{R}ep(\Arc H))\]
and refer to it as the factorizable Kazhdan--Lusztig category; it again possesses the structure of a unital factorization category.

Because the category $\mc{R}ep(\Arc H)$ is canonically self-dual (Remark \ref{rem-Rep^!=Rep^*}), it also admits a left action by $\mc{I}nd\mc{C}oh_{*, \tx{ren}}([\widehat{\tx{Hecke}}_H]_{\Ran_\dR})$. Dualizing, we obtain a \emph{right} action by $\mc{I}nd\mc{C}oh^{!, \tx{ren}}(\mb{B} \Arc H)$ on it. We define\footnote{The meaning of the negative sign in $\mathcal{KL}_*(H)_{\mathscr{L}^{-1}}$ will be clear in Remark \ref{rem-KL-twisted-as-totalization-*}.}
\[\mc{KL}(H)_{\mathscr{L}^{-1}, \tx{co}} := \mathscr{L}\mod^r(\mc{Rep}(\Arc H));\]
One can also think of $\mc{KL}(H)_{\mathscr{L}^{-1}, \tx{co}}$ as the dual factorization category of $\mc{KL}(H)_{\mathscr{L}}$ (Remark \ref{rem-dual-!-version-*-version-KL-twisted}).

A salient feature of this category is the existence of an equivalence of unital factorization categories (Theorem \ref{thm-tate-twist})
\[\mc{KL}(H)_{\mathscr{L}} \simeq \mc{KL}(H)_{\mathscr{L} \tensor \mathscr{L}^\tx{Tate}, \tx{co}}\]
where $\mathscr{L}^\tx{Tate}$ is the Tate twist associated with $H$ (c.f. Examples \ref{example:tate-bundle}). Over each closed point $x$, $\mathscr{L}^\tx{Tate}$ is the unique line bundle whose pullback to $[\Loop H]_{\Arc H}^\wedge$ is the restriction of the determinant line bundle on $\Loop H$; over the entire curve, however, it involves an additional $E$-term.

\begin{example}
Let $\check{Z}_G^\circ$ be the dual torus to $Z_G^\circ$ (the neutral component of the center of $G$). We can consider $2 \crho := \det(\mf{n})$ as a character for $Z_G^\circ$, thus a dual character for $\check{Z}_G^\circ$. We have the Atiyah bundle
\[0 \to \tx{Lie}(\check{Z}_G^\circ) \tensor O_X \to \tx{At}(\omega_X^{\check{\rho}}) \to T_X \to 0\]
associated to the $\check{Z}_G^\circ$-bundle induced from $\omega_X^{1/2}$ via $2\crho$.
The Tate bundle $\mathscr{L}^\tx{Tate}_B$ for $B$ is the pullback from $[\widehat{\tx{Hecke}}_T]_{\Ran_\dR}$ of the line bundle $\mathscr{L}_{[-\kappa_{\tx{crit}, G}, \AnomalyTerm]}$, where $\AnomalyTerm$ (the \emph{anomaly term}) is defined via the following pullback diagram:
\[
\xymatrix{
0 \ar[r] & \omega_X \ar[d] \ar[r] & E_\tx{anom} \ar[d] \ar[r] & \mf{t} \tensor O_X \ar[d] \ar[r] & 0 \\
0 \ar[r] & \omega_X \ar[r] & \tx{At}(\omega_X^{\check{\rho}})^* \ar[r] & \mf{z} \tensor O_X \ar[r] & 0
}
\]
where the second row is the monoidal dual of the Atiyah bundle. (This follows from Remark \ref{rem-compatible-pair-two-gerbe}, Example \ref{exam-tate-det-local-global} and \cite[Section 6]{zhao2017quantum}.)
\end{example}

As explained in Remark \ref{rem-gerbe-vs-central-extension}, the fiber of the factorization category $\mc{KL}(G)_{\mathscr{L}_{[\kappa, 0]}}$ at any $\mb{C}$-point is given by
\[\tx{KL}_\kappa(G) := \mathscr{L}_{[\kappa + \kappa_{\tx{crit}, G}, 0]}\mod(\Rep(\Arc G));\]
here $\Rep(\Arc G) := \IndCoh^{!,\tx{ren}}(\mc{L}^+ G)$ is acted on by $\IndCoh^{!, \tx{ren}}(\Arc G\backslash \Loop G_{\Arc G}^\wedge / \Arc G)$ via convolution. By \emph{loc.cit}, this category coincides with the \emph{renormalized} derived category of $(\hat{\mf{g}}, \mc{L}^+ G)$ Harish-Chandra modules as defined in \cite{gaitsgory2021conjectural}.

\begin{definition}
By restriction, $\mathscr{L}_{[\kappa + \kappa_{\tx{crit}, G}, 0]}$ can be also see as a line bundle over $H \backslash \Loop G_{\Arc G}^\wedge / H$ for any compact open $H \subseteq \Loop^+ G$. We define
\[\kmgmod^H := \mathscr{L}_{[\kappa + \kappa_{\tx{crit}, G}, 0]}\mod(\Rep(H));\]
here $\Rep(H) := \IndCoh^{!,\tx{ren}}(H)$ is acted on by $\IndCoh^{!, \tx{ren}}(H \backslash \Loop G_{\Arc G}^\wedge / H)$ via convolution.
\end{definition}

This is likewise the renormalized derived category of $(\hat{\mf{g}}, H)$ Harish-Chandra modules. We shall be mainly interested in $\kmgmod^I$, known as the affine BGG category $\mc{O}$.

Let $\mathscr{L} := \mathscr{L}_{[\kappa + \kappa_{\tx{crit}, G}, 0]}$.
Fix a $\mb{C}$-point $x$ on the curve $X$. As explained in Section \ref{ssec-fact-mod-KL}, this category extends to a $\mc{KL}(G)_{\mathscr{L}}$-factorization module category $\widetilde{\mc{L} \mf{g}}\mod_{\mathscr{L}}^{\Arc G \times_{G_x} B_x}$, whose fiber at $x$ recovers $\kmgmod^I$.

Letting $\mathscr{L}$ also denote the restriction along $B \subseteq G$, then we have a (c.f. Construction \ref{constr-KL-functorial}) factorizable functor
\[\tx{res}: \widetilde{\mc{L} \mf{g}}\mod_{\mathscr{L}}^{\Arc G \times_{G_x} B_x} \to \mc{KL}(B)_{\mathscr{L}}\]
defined by $!$-pullback of ind-coherent sheaves. On the other hand, $*$-pushforward (which is well-defined for the renormalized $*$-IndCoh category) gives a (lax-unital) factorization functor
\[C^\semiinf: \mc{KL}(B)_{\mathscr{L}} \simeq \mc{KL}(B)_{\mathscr{L} \tensor \mathscr{L}^\tx{Tate}_B, \tx{co}} \to \mc{KL}(T)_{\mathscr{L}_{[\kappa, \AnomalyTerm]}, \tx{co}} \simeq \mc{KL}(T)_{\mathscr{L}_{[\kappa, \AnomalyTerm]}};\]
The meaning of the superscript will become evident in the next subsection.

\subsection{Affine Lie algebra Representations}

\label{sect:prep_km}

In this section, we collect some facts about $\kmgmod^I$, mostly following \cite{gaitsgory2021conjectural} and \cite{frenkel2006local}.

\begin{convention}
In this subsection only, $\kappa$ will denote any non-critical level unless specified otherwise.
\end{convention}

\begin{definition}
For each $\clambda \in \weightLat^+$, we have the Weyl modules
\[\mb{V}_\kappa^\clambda := \tx{Ind}_{\Arc \mf{g} \oplus \mb{C}}^{\hat{\mf{g}}} V^\clambda \in \tx{KL}(G)_\kappa^\heartsuit,\] where $V^\clambda$ is the standard finite-dimensional irreducible representation of $G$.
We shall refer to $\mb{V}^0_\kappa$ as the \emph{vacuum} representation at level $\kappa$.
\end{definition}

\begin{definition}
Similarly, for each $\clambda \in \weightLat$, we have the affine Verma module
\[\mb{M}_\kappa^\clambda := \tx{Ind}_{\tx{Lie}(I)}^{\hat{\mf{g}}} k^\clambda \in \kmgmod^{I, \heartsuit},\]
We let $\mb{M}_\kappa^{\clambda, \vee}$ denote the dual affine Verma modules, which are their right orthogonals, i.e. the unique object in $\kmgmod^{I}$ satisfying
\[\tx{Hom}_{\kmgmod^I}(\mb{M}_\kappa^\clambda, \mb{M}_\kappa^{\clambda', \vee}) = \begin{cases}
k & \tx{if }\clambda = \clambda', \\
0 & \tx{otherwise.}
\end{cases}
\]
\end{definition}

\begin{definition}[Averaging]
	For every $K_1 \subseteq K_2$, we let $\tx{Av}_{*}^{K_2 / K_1}$ denote the right adjoint to the forgetful functor $\kmgmod^{K_2} \to \kmgmod^{K_1}$.
\end{definition}

\begin{definition}[Categorical Duality]
\label{def:km-categorical-duality}
Our general study of ind-coherent sheaves, in particular Theorem \ref{thm-tate-twist}, recovers the classically known fact that for every compact open $H$ there exists a duality pairing of categories
\[\langle -, - \rangle_{\kmgmod^H}: \kmgmod^H \tensor \nkmgmod^H \to \Vect\]
(c.f. \cite[Section 1.3]{gaitsgory2021conjectural}).
\end{definition}
We let $\mb{D}_{\tx{KM}}$ denote the resulting functor on compact objects\footnote{This ought not to cause confusion: the duality for different $H$ are compatible under forgetful functors.}. This functor has the following properties:
\begin{itemize}
\item For any $\kappa$ and $\clambda$ we have $\mb{D}_{\tx{KM}}(\mb{M}_{\kappa}^{\clambda}) \simeq \mb{M}_{-\kappa}^{-\clambda - 2 \crho}[\dim(\mf{n})]$;
\item For any $\kappa$ we have $\mb{D}_{\tx{KM}}(\mb{V}_\kappa^0) \simeq \mb{V}_{-\kappa}^0$.
\end{itemize}

\paragraph{Highest Weight Structure}

The collection of objects $\mb{M}_\kappa^{\clambda}$ and $\mb{M}_\kappa^{\clambda, \vee}$ form a \emph{highest weight structure} for $\kmgmod^I$. However, it turns out that they do not correspond to the \emph{quantum} (dual) Verma modules under our equivalence. Instead, we will have to use the following objects:

\begin{definition}
\label{def:KM-st-cost-def}
At $\kappa$ positive (or irrational), we will refer to
\[\StObj_{\tx{KM}, \kappa}^{\clambda} := \mb{D}_{\tx{KM}}(\mb{W}_{-\kappa}^{1, -\clambda - 2\crho}[\dim(\mf{n})])\]
as the \emph{standard object} of highest weight $\clambda$, and
\[\CoStObj_{\tx{KM}, \kappa}^{\clambda} := \mb{W}_{\kappa}^{w_0, \clambda}\]
as the \emph{costandard object} of highest weight $\clambda$. Here $1$ and $w_0$ are respectively the shortest and longest element in the finite Weyl group, and $\mb{W}_\kappa^{w, \clambda}$ is the \emph{Wakimoto} module of level $\kappa$, type $w$ and highest weight $\clambda$.
\end{definition}

To define Wakimoto modules and check that the definition is valid (i.e. that $\mb{W}_{-\kappa}^{1, -\clambda - 2\crho}$ is compact), we need a bit more representation theory.

\paragraph{Sheaves on Affine Flag Variety}
Let $I^\circ$ denote the unipotent radical of $I$, i.e. the preimage of $N$ under $\Arc G \xrightarrow{t = 0} G$. Let $\tx{Fl}_G := \Loop G / I$ denote the affine flag variety, and $\widetilde{\tx{Fl}_G} := \Loop G / I^\circ$ denote the enhanced affine flag variety.

As in \cite[Section 20.4]{frenkel2006local}, for each $\clambda \in \weightLat$ the pair $(-\kappa, \clambda)$ determines a twisting on $\tx{Fl}_G$, and we let $\DMod_{(-\kappa, \clambda)}(\tx{Fl}_G)$ denote the resulting category of twisted D-modules. When $\clambda = 0$, we will simply use the subscript $-\kappa$.
The category $\DMod_{-\kappa}(\tx{Fl}_G)^I$ admits a monoidal structure known as \emph{convolution}, which we will denote by $\star$. It acts on $\nkmgmod^I$ via this convolution structure as well.

\paragraph{Metaplectic Duality}
Take $\kappa$ non-critical, then we can consider $\kappa$ as a symmetric bilinear form on $\coweightLat$. Assume $\kappa$ is \emph{rational}, which means its values are torsion as elements in $\mb{C} / \mb{Z}$. Let $\coweightLat^\sharp$ be its kernel, and $\weightLat^\sharp$ the dual lattice to $\coweightLat^\sharp$.
To each $\alpha \in \corootLat$, we attach an element $\tx{ord}(q(\alpha)) \cdot \alpha$ where $q(\alpha) := \frac{\kappa(\alpha, \alpha)}{2} \in \mb{C} / \mb{Z}$; these form a set $\corootLat^\sharp$. Similarly, we form $\rootLat^\sharp \subseteq \rootLat \tensor_{\mb{Z}} \mb{Q}$ by attaching to each $\check{\alpha}$ the element $\frac{1}{\tx{ord}(q(\alpha))} \cdot \check{\alpha}$.
$(\weightLat^\sharp, \rootLat^\sharp, \coweightLat^\sharp, \corootLat^\sharp)$ form another root datum, corresponding to a pair $(G^\sharp, T^\sharp)$.
\begin{definition}
\label{def:metaplectic-dual-group}
We let $(H, T_H)$ be the Langlands dual of $(G^\sharp, T^\sharp)$ and call it the \emph{metaplectic dual} of $(G, T)$. (For instance, at integral levels we have $(H, T_H) \simeq (\check{G}, \check{T})$.) For convenience's sake, we extend this definition to all levels by declaring the metaplectic datum to be trivial when $\kappa$ is irrational.
\end{definition}

We let $\Waffextrat := \Wfin \ltimes \ratLat \subseteq W^{\tx{aff}, \tx{ext}}$ denote the extended affine Weyl group corresponding to $\ratLat$.

\begin{definition}
For $\tilde{w} \in \Waffextrat$, we let $j_{\tilde{w}, *}, j_{\tilde{w}, !} \in \DMod_{-\kappa}(\tx{Fl}_G)^I$ denote the $!$ and $*$-extensions of the constant sheaf supported on the orbit $I \tilde{w} I/I$ of $\tx{Fl}_G$; we normalize them such that $!$-fibers at $\tilde{w}$ is placed in degree $\ell(w)$, and refer to them as (co)standard sheaves on the affine flag variety.
\end{definition}
These (co)standard sheaves satisfy the following properties:
\begin{itemize}
\item For $\tilde{w}_1, \tilde{w}_2$ such that $\ell(\tilde{w}_1) + \ell(\tilde{w}_2) = \ell(\tilde{w}_1 \tilde{w_2})$, we have
\[j_{\tilde{w}_1, *} \star j_{\tilde{w}_2, *} \simeq j_{\tilde{w}_1 \tilde{w}_2, *} \hspace{1em} j_{\tilde{w}_1, !} \star j_{\tilde{w}_2, !} \simeq j_{\tilde{w}_1 \tilde{w}_2, !}\]
\item $j_{\tilde{w}, *} \star j_{\tilde{w}^{-1}, !} \simeq \delta_{1} \simeq j_{\tilde{w}^{-1}, !} \star j_{\tilde{w}, *}$;
\end{itemize}
We will write $j_{\mu, !}$ (resp.\ $j_{\mu, *}$) for $\mu \in \ratLat$ by considering it as a subgroup of $\Waffextrat$.

Next we introduce \emph{Wakimoto} sheaves on the affine flag variety, which is a collection of objects $\{J_\lambda\}_{\lambda \in \ratLat} \in \DMod_{-\kappa}(\tx{Fl}_G)^I$ such that\footnote{We follow the convention of \cite{gaitsgory2021conjectural} which is opposite from that of \cite{arkhipov2009perverse}.}
\begin{itemize}
\item For $\mu \in (\ratLat)^{-}$ we have $J_\mu \simeq j_{\mu, *}$;
\item For $\mu \in (\ratLat)^{+}$ we have $J_\mu \simeq j_{\mu, !}$;
\item For any $\mu_1, \mu_2 \in \ratLat$ we have $J_{\mu_1} \star J_{\mu_2} \simeq J_{\mu_1 + \mu_2}$.
\end{itemize}
For convenience's sake we will also extend the definition to all of $\Waffextrat$ by declaring $J_{\tilde{w}} := J_{\mu} \star j_{w_f}$ where $\tilde{w} = \mu w_f, \mu \in \ratLat, w_f \in W_\tx{fin}$. For more details regarding these objects in at rational levels, see \cite[Section 5]{yang2021twisted}.

\begin{warning}
\label{warn:cotangent-fiber}
We warn that we are \emph{not} keeping track of the cotangent fiber as defined in \cite[Section 2.5.3]{gaitsgory2021conjectural}. Such additional bookkeeping would be necessary if one were to consider the \emph{factorizable} version of our main theorem (c.f. Expectation \ref{ex:NKTO-equivalence}).
\end{warning}

\paragraph{Kashiwara-Tanisaki Localization}
We have an $\IndCoh$-global section functor
\[\Gamma_{-\kappa}(\tx{Fl}_G, -)^{I^\circ}: \DMod_{(-\kappa, \clambda)}(\tx{Fl}_G)^{I^\circ} \to \nkmgmod^{I^\circ}\]
as in \cite{kashiwara1996kazhdan} \footnote{Because our Kac-Moody category is renormalized, some care is needed in comparing it with the classical notion; we refer readers to \cite[Appendix A]{raskin2020affine} for more discussion.}. This functor also respects the remaining strong $T$-action, and we have an induced functor on the strong $T$-invariants:
\[\Gamma_{-\kappa}(\tx{Fl}_G, -)^{I}: \DMod_{(-\kappa, \clambda)}(\tx{Fl}_G)^{I} \to \nkmgmod^{I}\]

\begin{definition}
\label{def:admissible-weight}
An element $\clambda \in \weightLat$ is called $\kappa$-admissible if:
\begin{itemize}
\item $\langle \alpha, \clambda + \crho \rangle \not\in \mb{Z}^{>0}$ every positive coroot $\alpha$; and
\item $\langle \alpha, \clambda + \crho \rangle + n \cdot \frac{-\kappa(\alpha, \alpha)}{2} \not\in \mb{Z}^{>0}$ every positive coroot $\alpha$ and every $n \in \mb{Z}^{>0}$.
\end{itemize}
\end{definition}

\begin{theorem}[\cite{kashiwara1996kazhdan}]
\label{thm:kashiwara-tanisaki}
If $\clambda$ is $\kappa$-admissible, then the functor $\Gamma_{-\kappa}(\tx{Fl}_G, -)^{I^\circ}$ and $\Gamma_{-\kappa}(\tx{Fl}_G, -)^I$ are $t$-exact, send each standard sheaf $j_{\tilde{w}, !}$ to $\mb{M}_{-\kappa}^{\tilde{w} \cdot \clambda}$ and each costandard sheaf $j_{\tilde{w}, *}$ to $\mb{M}_{-\kappa}^{\tilde{w} \cdot \clambda, \vee}$.
\end{theorem}

\begin{remark}
	Under the more stringent condition of $\clambda$ being \emph{regular $\kappa$-admissible} (which means replacing $\not\in \mb{Z}^{>0}$ with $\not\in \mb{Z}^{\ge 0}$ in Definition \ref{def:admissible-weight}), it is also known that $\Gamma_{-\kappa}(\tx{Fl}_G, -)^{I^\circ}$ is a fully faithful embedding of a direct summand of $\nkmgmod^{I^\circ}$ (c.f. \cite[Theorem 5.5]{frenkel2009localization}). We do not invoke this version.
\end{remark}

\paragraph{Semi-infinite Cohomology}
We refer readers to \cite[Section 2]{gaitsgory2021conjectural} for a detailed discussion on semi-infinite cohomology.
\begin{definition}[Definition-Construction]
	Write $\Loop N$ as a union of group indschemes $\bigcup_{\alpha} N_\alpha$, each containing $N_0 := \Arc N$; let $\mf{n}_\alpha$ denote their Lie algebras. For $N_{\alpha_1} \subseteq N_{\alpha_2}$, there is a canonically defined connecting map $C^*(\mf{n}_{\alpha_1}, -) \Rightarrow C^*(\mf{n}_{\alpha_2}, - \tensor \det(\mf{n}_{\alpha_2} / \mf{n}_{\alpha_1}))$. There is a \emph{semi-infinite cohomology} functor
	\[C_*^{\semiinf}(\Loop \mf{n}, \Arc \mf{n}, -): \kmgmod^{\Arc T} \to \tx{KL}_{\kappa}(T)\]
	whose composition with the forgetful to $\Vect$ is given by
	\[V \mapsto \colimit_{\alpha} C^*(\mf{n}_\alpha, V \tensor \det(\mf{n}_{\alpha} / \Arc \mf{n}));\]
	For each $\clambda \in \weightLat$ we set
	\[C_*^{\semiinf}(\Loop \mf{n}, \Arc \mf{n}, -)^{\clambda} := \Hom_{\Rep(\Arc T)}(k^\clambda, C_*^{\semiinf}(\Loop \mf{n}, \Arc \mf{n}, -)): \kmgmod^{\Arc T} \to \Vect.\]
\end{definition}

\paragraph{Wakimoto Modules}
For any $w \in \Wfin, \clambda \in \weightLat$, the \emph{Wakimoto module of type $w$}, denoted $\mb{W}_\kappa^{w, \clambda} \in \kmgmod^{I, \heartsuit}$, is defined in \cite[Section 2.4]{gaitsgory2021conjectural}. Here we will give a ``quick and dirty'' definition:
\begin{definition}
	For $\clambda \in \weightLat$, the type 1 Wakimoto module $\mb{W}_\kappa^{1, \clambda}$ is the unique object in $\kmgmod^I$ satisfying
	\[\langle \mb{W}_\kappa^{1, \clambda}, - \rangle_{\hat{\mf{g}}\mod^I} \simeq C_*^{\semiinf}(\Loop \mf{n}, \Arc \mf{n}, - \tensor \det(\mf{n}^-)^{\tensor -1})^{-\clambda} : \nkmgmod^{I} \to \Vect;\]
	For $w \in \Wfin$, we define the type $w$ Wakimoto module $\mb{W}_{\kappa}^{w, \clambda}$ as
	\[\mb{W}_\kappa^{w, \clambda} := j_{w, *} \star \mb{W}_{\kappa}^{1, \clambda}.\]
\end{definition}

\begin{remark}
Wakimoto modules are the ``limiting behavior'' of Verma modules as they are twisted by elements of the affine Weyl group; more precisely we have (\cite[Section 3.3.2]{gaitsgory2021conjectural})
\[\mb{W}_{\kappa}^{1, \clambda} \simeq \colimit_{\mu \in \coweightLat^{\sharp, +}} j_{-\mu, *} \star \mb{M}_{\kappa}^{\clambda + \mu}.\]
\end{remark}

We will mainly be interested in the \emph{type 1} Wakimoto modules $\mb{W}_\kappa^{1, \clambda}$ and the \emph{type $w_0$} version $\mb{W}_\kappa^{w_0, \clambda}$; these should be considered as semi-infinite Verma at the negative level (resp.\ semi-infinite dual Vermas at the positive level).
We will mostly only need their formal properties listed below:

\tb{For all levels}
The following properties are true for Wakimoto modules at all levels:
\begin{itemize}
\item (\cite[Section 2.5.9]{gaitsgory2021conjectural}) For any $\mu \in \ratLat$, $w \in \Wfin, \clambda \in \weightLat$, we have $J_\mu \star \mb{W}_{\kappa}^{w, \clambda} \simeq \mb{W}_\kappa^{w, \clambda + \mu}$;
\item (\cite[Section 2.4.7]{gaitsgory2021conjectural}]) $C_*^{\semiinf}(\Loop \mf{n}^-, \Arc \mf{n}^-, \mb{W}_{\kappa}^{1, \clambda})^{\clambda'} = \begin{cases}
\det(\mf{n}^-)^{\tensor -1} & \tx{if}~\clambda' = \clambda + 2 \crho \\
0 & \tx{else}
\end{cases}$
\item (\cite[Section 1.3.7]{quantum-vs-semiinf}) $C_*^{\semiinf}(\Loop \mf{n}, \Arc \mf{n}, \mb{W}_{\kappa}^{w_0, \clambda})^{\clambda'} = \begin{cases}
\mb{C} & \tx{if}~\clambda' = \clambda \\
0 & \tx{else}
\end{cases}$
\end{itemize}

\tb{For negative levels}
At negative (and critical) levels, Kac-Kazhdan theorem implies that Wakimoto modules will converge to the opposite Verma within finite time. 
For $\clambda$ sufficiently dominant, we have
\[\mb{W}_{-\kappa}^{1, \clambda} \simeq \mb{M}_{-\kappa}^{\clambda} \hspace{1em} \mb{W}_{-\kappa}^{w_0, \clambda} \simeq \mb{M}_{-\kappa}^{\vee, \clambda};\]
whereas for $\cmu$ sufficiently anti-dominant, we have
\[\mb{W}_{-\kappa}^{1, \cmu} \simeq \mb{M}_{-\kappa}^{\vee, \cmu} \hspace{1em} \mb{W}_{-\kappa}^{w_0, \cmu} \simeq \mb{M}_{-\kappa}^{\cmu}.\]

Now we check that Definition \ref{def:KM-st-cost-def} is a valid one:

\begin{lemma}
\label{lemma:wakimoto-generates}
	The type 1 Wakimoto modules $\mb{W}_{-\kappa}^{1, \clambda}$ form a set of compact generators of $\nkmgmod^I$.
\end{lemma}

\begin{proof}
It follows from the properties of Wakimoto sheaves that, for every $\clambda$, we can choose $\mu \in (\ratLat)^+$ such that $\clambda + \mu \in \weightLat^{+}$, and have $\mb{W}_{-\kappa}^{1, \clambda} \simeq j_{-\mu, *} \star \mb{M}_{-\kappa}^{\clambda + \mu}$. Since $j_{-\mu, *}$ is an invertible operation and Vermas are compact generators, we know that $\mb{W}_{-\kappa}^{1, \clambda}$ are compact as well. To show they generate, it suffices to show that every Verma can be written as a \emph{finite} colimit of Wakimoto modules. Fix some $\clambda_0$ admissible and $\tilde{w} \in \Waffextrat$ such that $\clambda = \tilde{w} \cdot \clambda_0$, and invoke Kashiwara-Tanisaki localization. Then we have that in such case,
\[\Gamma_{-\kappa}(\tx{Fl}_G, j_{\tilde{w}, !} \star \delta_1)^I \simeq \mb{M}^{\clambda}_{-\kappa} \hspace{1em} \Gamma_{-\kappa}(\tx{Fl}_G, J_{\tilde{w}} \star \delta_1)^I \simeq \mb{W}^{1, \clambda}_{-\kappa},\]
that the support range (more precisely, images in $K_0$) of $j_{\tilde{w}, !}$ and $J_{\tilde{w}}$ are the same, and the nonzero map $\mb{M}^{\clambda}_{-\kappa} \to \mb{W}^{1, \clambda}_{-\kappa}$ is induced by a map $j_{\widetilde{w}, !} \to J_{\widetilde{w}}$ that is isomorphism on the open strata (\cite[Theorem 5, Lemma 27]{arkhipov2009perverse}). The claim now follows by induction on support dimension.
\end{proof}

\begin{remark}
Our standard objects are designed such that for every $M \in \kmgmod^I$, we have
\[\Hom_{\kmgmod^I}(\StObj_{\tx{KM}, \kappa}^{\clambda}, M) \simeq C_*^\semiinf(\Loop \mf{n}, \Arc N, M)^{\clambda};\]
and the costandards are their right orthogonals.
\end{remark}

\subsection{Jacquet Functor and Local-Global Correspondence}

It follows from definition that the fiber of the functor
\[C^\semiinf \circ \tx{res}: \widetilde{\mc{L} \mf{g}}\mod_{\mathscr{L}}^{\Arc G \times_{G_x} B_x} \to \mc{KL}(T)_{\mathscr{L}_{[\kappa, \AnomalyTerm]}}\]
at the distinguished point $x$ is precisely $C_*^\semiinf(\Arc \mf{n}, \Loop N, -): \kmgmod^I \to \tx{KL}_{\kappa}(T)$. Now, because $\mc{KL}(G)_\kappa$ is an \emph{unital} factorization category with vacuum $\mb{V}_\kappa^0$, we have (Proposition \ref{prop:unital-basic-properties}) an equivalence of categories
\[(\kmgmod^I) \simeq \mb{V}_\kappa^0\tx{-FactMod}_{\un}(\widetilde{\mc{L} \mf{g}}\mod_{\mathscr{L}_{[\kappa, 0]}}^{\Arc G \times_{G_x} B_x});\]
using the fact that $C^\semiinf \circ \tx{res}$ is lax-unital we obtain a functor
\[\tb{J}^{\tx{KM}}_\kappa: \kmgmod^I \simeq \mb{V}_\kappa^0\tx{-FactMod}_{\un}(\widetilde{\mc{L} \mf{g}}\mod_{\mathscr{L}_{[\kappa, 0]}}^{\Arc G \times_{G_x} B_x}) \xrightarrow{C^\semiinf \circ \tx{res}} C^\semiinf(\tx{res}(\mb{V}_\kappa^0))\tx{-FactMod}_\un(\mc{KL}(T)_{\mathscr{L}_{[\kappa, \AnomalyTerm]}}),\]
whose composition with taking the fiber at $x$ recovers $C_*^\semiinf(\Arc \mf{n}, \Loop N, -)$. We will refer to this functor as the \emph{Kac-Moody Jacquet\footnote{The name ``Jacquet functor'' has its origins from the Whittaker side of the fundamental local equivalence.} functor}.

For the rest of this subsection, consider the special case $X = \mb{P}^1$, and let $\modulePoints := (0, \infty)$ be the two distinguished $\mb{C}$-points of $\mb{P}^1$. The natural multi-point generalization of the above construction gives a functor
\[\tb{J}^{\tx{KM}}_\kappa: (\kmgmod^I)_{0} \tensor (\kmgmod^I)_\infty \to C^\semiinf(\tx{res}(\mb{V}_\kappa^0))\tx{-FactMod}_\un([\mc{KL}(T)_{\mathscr{L}_{[\kappa, \AnomalyTerm]}}]_{\Ran_{\modulePoints}(\mb{P}^1)_{\un, \dR}}),\]
where on the LHS we use subscripts to indicate that we identify these two categories as fibers as $0$ resp.\ $\infty$, and on RHS we use subscript to indicate that we work with the multi-marked-point version of the unital Ran space (c.f. Example \ref{example:Ran-module-spaces}).

Set
\[\tx{Bun}_G(\mb{P}^1)^{\tx{level}_B, \modulePoints} := \tx{Bun}_G(\mb{P}^1) \times_{[(\tx{pt} / G)_0 \times (\tx{pt} / G)_{\infty}]} [(\tx{pt} / B)_0 \times (\tx{pt} / B)_\infty]\]
(where subscripts $0$ and $\infty$ are to clarify where these reductions take place), and consider the diagram
\[\xymatrix{
 & \tx{Bun}_B(\mb{P}^1) \ar^{\ms{q}}[rd] \ar_{\ms{p}}[ld] \\
\tx{Bun}_G(\mb{P}^1)^{\tx{level}_B, \modulePoints} & & \tx{Bun}_T(\mb{P}^1)
}\]
Recall from \cite{gaitsgory2017strange} the \emph{geometric constant term functor}
\[\tx{CT}_*:= \ms{q}_* \circ \ms{p}^!: \DMod_{[\kappa + \kappa_{\tx{crit}, G}, 0]}(\tx{Bun}_G(\mb{P}^1)^{\tx{level}_B, \modulePoints}) \to \DMod_{[\kappa + \kappa_{\tx{crit}_G}, 0]}(\tx{Bun}_T(\mb{P}^1));\]
also recall the line bundle $\mc{L}_{\tx{Tate}(\mf{n})}$ on $\tx{Bun}_T(\mb{P}^1)$, which satisfies
\[\mc{L}_{\tx{Tate}(\mf{n})}|_{\mc{P}_T} \simeq \det R \Gamma(\mb{P}^1, \mf{n}_{\mc{P}_T}[1])\]
here $\mf{n}_{\mc{P}_T}$ is the associated vector bundle of the adjoint action, and the (integral) twisting corresponding to this line bundle is the \emph{anomaly term} $[-\kappa_{\tx{crit}_G}, \AnomalyTerm]$. We incorporate this shift into our definition by setting
\[\tx{CT}_*^{\tx{shifted}} := [\mc{L}_{\tx{Tate}(\mf{n})}[-d] \tensor -] \circ \tx{CT}_*: \DMod_{[\kappa + \kappa_{\tx{crit}, G}, 0]}(\tx{Bun}_G(\mb{P}^1)^{\tx{level}_B, \modulePoints}) \to \DMod_{[\kappa , \AnomalyTerm]}(\tx{Bun}_T(\mb{P}^1)),\]
where $d$ is a number that depends on the connected component of $\tx{Bun}_T(\mb{P}^1)$: over the $\lambda$-component it is $-\langle \lambda, 2 \crho \rangle - \dim(\mf{n})$.

In Section \ref{sect:km-localization}, we will construct two \emph{chiral localization} functors $[\tx{Loc}]_{\modulePoints}$ and $[\tx{Loc}]^{\tx{ch}}_{\modulePoints}$, such that we have the following crucial \emph{local-global compatibility} diagram:
\begin{equation}
\label{eq:local-global-diagram}
\xymatrix{
(\kmgmod^I)_{x_0} \tensor (\kmgmod^I)_{x_\infty} \ar^(0.35){\tb{J}_\kappa^{\tx{KM}}}[r] \ar^{[\tx{Loc}]_{\modulePoints}}[d] & C_\tx{co}^\semiinf(\tx{res}(\mb{V}_\kappa^0))\tx{-FactMod}_\un([\mc{KL}(T)_{\mathscr{L}_{[\kappa, \AnomalyTerm]}}]_{(\Ran_{\modulePoints}(\mb{P}^1)_{\un, \dR}}) \ar^{[\tx{Loc}]^{\tx{ch}}_{\modulePoints}}[d] \\
\DMod_{[\kappa + \kappa_{\tx{crit}, G}, 0]}(\tx{Bun}_G(\mb{P}^1)^{\tx{level}_B, \modulePoints}) \ar^{\tx{CT}_*^{\tx{shifted}}}[r] & \DMod_{[\kappa , \AnomalyTerm]}(\tx{Bun}_T(\mb{P}^1))
}
\end{equation}
The role this diagram plays will become evident during the proof of the main theorem.

\subsection{Conversion into Sheaves}
\label{sect:conf-space-in-minimal}
Let $[\kappa, 0]$ be a non-degenerate level and $[\check{\kappa}, 0]$ be its Langlands dual level. As we explained in the introduction, we have a unital factorization category $\DMod_{[\check{\kappa}, 0]}([\tx{Gr}_{\check{T}}]_{\Ran_{\un}})$ whose fiber at any point on the curve is given by $\check{\kappa}$-twisted D-modules over $(\tx{Gr}_{\check{T}})$, and there is an equivalence of unital factorization categories
\[\tx{FLE}_T: \mc{KL}(T)_{[\kappa, 0]} \simeq \DMod_{[\check{\kappa}, 0]}([\tx{Gr}_{\check{T}}]_{\Ran_{\un}});\]
we would like to pass through this equivalence to transform the output of $\tb{J}_\kappa^\tx{KM}$ into D-modules. For technical reasons, however, we will use a slight variant as follows.

A first caveat is that the category at the receiving end of our Jacquet functor has a nonzero $E$-term in its parameter, namely $\AnomalyTerm$. This can be taken into account as follows.
\begin{definition}
\label{def:rho-shifted-gr}
Let $[\tx{Gr}_{\check{T}}^{\omega^{\crho}}]_{\Ran_\dR}$ be the \emph{$\crho$-shifted affine Grassmannian}, whose $S$-points classify:
\begin{itemize}
\item a point $x: S^{\tx{red}, \tx{cl}} \to \Ran$;
\item a $\check{T}$-bundle $\mc{P}$ on the adic disc $\mc{D}'_{\Gamma_x}$ around the graph of $x$;
\item an isomorphism $\alpha$ of $\mc{P}$ on the punctured disc with the restriction of $\omega^{\crho}$;
\end{itemize}
\end{definition}
There is an isomorphism of unital factorization spaces $- \tensor \omega^{\crho}: [\tx{Gr}_{\check{T}}]_{\Ran_\dR} \to [\tx{Gr}_{\check{T}}^{\omega^{\crho}}]_{\Ran_\dR}$. We let $\mathscr{L}_{\DualTwisting}$ denote the unique twisting on $\tx{Gr}_{\check{T}}^{\omega^{\crho}}$ whose pullback along this map gives $\mathscr{L}_{[\check{\kappa}, 0]}$. There is then a unital factorization category $\DMod_{\DualTwisting}([\tx{Gr}_{\check{T}}^{\omega^{\crho}}]_{\Ran_{\un, \dR}})$, and we have (c.f. Section \ref{sect:tori}) another equivalence of unital factorization categories
\[\tx{FLE}^{\omega^{\crho}}_T: \mc{KL}(T)_{[\kappa, \AnomalyTerm]} \simeq \DMod_{\DualTwisting}([\tx{Gr}_{\check{T}}^{\omega^{\crho}}]_{\Ran_{\un}})\]
induced by $\tx{FLE}_T$.

For technical purposes it is better to work with \emph{negatively-weighted} factorization algebras on $\tx{Gr}_{\check{T}}^{\omega^{\crho}}$; geometrically this means working with the \emph{configuration space}. We refer readers to \cite[Section 4]{gaitsgory2019metaplectic} for proofs of statements in this subsection.

Throughout this subsection, let $\modulePoints = (x_1, \ldots, x_k)$ is a fixed set of $\mb{C}$-points of $X$. For each smooth curve $X$, we have
\begin{itemize}
\item $\tx{Conf}(X)$ the scheme of $\weightLat^{<0} \setminus \{0\}$-colored divisors on $X$;
\item $\tx{Conf}^\circ(X)$, the open part where each point is marked by a negative simple root;
\item $\tx{Conf}_{\infty \cdot \modulePoints}$, the ind-scheme where the marking on $x_1, \ldots, x_k$ are now allowed to be any element of $\weightLat$;
\item $\tx{Conf}^{=\clambda}$ and $\tx{Conf}_{\infty \cdot \modulePoints}^{=\clambda}$, the stratification by total weights.
\end{itemize}
The space $\tx{Conf}(X)$ admits an additive semigroup structure, and $\tx{Conf}_{\infty \cdot \modulePoints}$ is a module space of it. Using $\tx{Conf}^\circ(X)$ in replacement of $[\Ran \times \Ran]_{\tx{disj}}$ one still has a notion of factorizable structures over the configuration space; we will use the same notations as before to denote these variants as well.

Let $\mathscr{G}$ be any factorizable twisting on the $\crho$-shifted affine Grassmannian. As we will show in Section \ref{sect:conf-compare-gr}, to each factorization algebra $A$ supported on a certain subspace $[\tx{Gr}_{\check{T}}^{\omega^{\crho}, <0}]_{\Ran}$ (the fibers of which corresponds to $\weightLat^{<0} \subseteq \weightLat$), one can attach a factorization algebra $A_{\tx{red}, \tx{Conf}}$ over the configuration space, and we have a reduction functor
\[\tx{Red}_{\tx{Conf}}: A\tx{-FactMod}_\un(\DMod_{\mathscr{G}}([\tx{Gr}_{\check{T}}^{\omega^{\crho}}]_{\Ran_{\modulePoints}})) \to A_{\tx{red}, \tx{Conf}}\tx{-FactMod}_{\nonun}(\DMod_{\mathscr{G}}(\tx{Conf}_{\infty \cdot \modulePoints})).\]

\begin{remark}
There are a few key differences in using configuration spaces vs. the Ran space: it is an actual scheme (in the module case, ind-scheme), and has very strong finiteness properties (so that operations such as Verdier duality are well-behaved). Also note that for the configuration space, the factorization condition does not involve \emph{diagonal restriction}.
\end{remark}

Let us return to our Jacquet functor. Consider 
\begin{equation}
\label{eq:jacquet-part-i}
\begin{split}
& (\kmgmod^I)_{x_0} \xrightarrow{\tb{J}^\tx{KM}_\kappa} C^\semiinf(\tx{res}(\mb{V}_\kappa^0))\tx{-FactMod}_\un([\mc{KL}(T)_{\mathscr{L}_{[\kappa, \AnomalyTerm]}}]_{(\Ran_{x_0})_{\un, \dR}}) \\
& \xrightarrow[\simeq]{\tx{FLE}^{\crho}_T} \tx{FLE}^{\crho}_T \circ C^\semiinf(\tx{res}(\mb{V}_\kappa^0))\tx{-FactMod}_\un([\DMod_{\DualTwisting}(\tx{Gr}_{\check{T}}^{\omega^{\crho}})]_{(\Ran_{x_0})_{\un, \dR}}).
\end{split}
\end{equation}

We remind that at the fiber level, the functor $C^\semiinf$ is given by $C_*^\semiinf(\Loop \mf{n}, \Arc \mf{n}, -)$. Cohomological estimates from Proposition \ref{prop:KM-vacuum-verification} below therefore guarantees that the process lands in the full subcategory of factorization modules supported on the configuration space, and we can continue the process above as
\begin{equation}
\label{eq:jacquet-part-ii}
\xrightarrow{\tx{Red}_\tx{Conf}} [\tx{FLE}^{\crho}_T \circ C^\semiinf(\tx{res}(\mb{V}_\kappa^0))]_{\tx{red}, \tx{Conf}}\tx{-FactMod}_\nonun(\DMod_{\DualTwisting}(\tx{Conf}_{\infty \cdot x_0})).
\end{equation}

We define
\[\TheFactAlg := \tx{FLE}^{\crho}_T \circ C^\semiinf(\tx{res}(\mb{V}_\kappa^0)) \in \tx{FactAlg}_{\nonun}(\DMod_{\DualTwisting}(\tx{Conf}))\]
and refer to the combination of Equation \eqref{eq:jacquet-part-i} and \eqref{eq:jacquet-part-ii} as
\[\tb{J}_{\kappa}: \kmgmod^I \to \TheFactAlg\tx{-FactMod}_\nonun(\DMod_{\DualTwisting}(\tx{Conf}_{\infty \cdot x_0})).\]
We shall proceed to show that this functor is an equivalence for $X = \mb{A}^1$, under mild constraint on $\kappa$.
\begin{remark}
We note (without proof) that $\tx{Red}_\tx{Conf}$ in this case is in fact an equivalence, so that $\tb{J}_\kappa^\tx{KM}$ is already an equivalence of categories.
\end{remark}

\subsubsection{Highest Weight Structure}
\label{sect:conf-hw}

Another benefit of working with configuration spaces is the following highest weight structure for factorization modules, developed in \cite[5.3]{gaitsgory2019metaplectic}. Assume $A_{\tx{red}, \tx{Conf}}$ is a non-unital factorization algebra on the configuration space as above, and such that it is \emph{ind-holonomic} as a D-module. Let us consider the case $\modulePoints = (x)$ a single fixed point. Then:
\begin{itemize}
\item For every $\clambda \in \weightLat$ there exists a \emph{standard} object $\StObj_{\tx{Conf}}^\clambda$ and a \emph{costandard} object $\CoStObj_{\tx{Conf}}^\clambda$ within the category $A_{\tx{red}, \tx{Conf}}\tx{-FactMod}_{\nonun}(\DMod_{\mathscr{G}}(\tx{Conf}_{\infty \cdot x}))$;
\item The standard objects form a set of compact generators of $ A_{\tx{red}, \tx{Conf}}\tx{-FactMod}_{\nonun}(\DMod_{\mathscr{G}}(\tx{Conf}_{\infty \cdot x}))$, and the costandards are their right orthogonals, in the sense that
\[\Hom_{A_{\tx{red}, \tx{Conf}}\tx{-FactMod}_{\nonun}(\DMod_{\mathscr{G}}(\tx{Conf}_{\infty \cdot x}))}(\StObj_{\tx{Conf}}^{\clambda}, \CoStObj_{\tx{Conf}}^{\cmu}) \simeq \begin{cases}
	\mb{C} & \clambda = \cmu \\
	0 & \tx{else}
\end{cases}\]
Moreover, this property characterizes the standard object $\StObj_{\tx{Conf}}^\clambda$ (resp.\ costandard object $\CoStObj_{\tx{Conf}}^\clambda$) as the unique object that exhibits such property with respect to all costandard (resp.\ standard) objects;
\item The standard object has the following universal property: for every $M$, we have
\[\Hom_{A_{\tx{red}, \tx{Conf}}\tx{-FactMod}_{\nonun}(\DMod_{\mathscr{G}}(\tx{Conf}_{\infty \cdot x}))}(\StObj_{\tx{Conf}}^{\clambda}, M) \simeq \iota_{\mu \cdot x}^!(M)\]
where $\iota_{\mu \cdot x}$ is the inclusion of the point $(\mu \cdot x)$ into $\tx{Conf}_{\infty \cdot x}$;
\item The costandard object has the following universal property: for every $M$, we have
\[\Hom_{A_{\tx{red}, \tx{Conf}}\tx{-FactMod}_{\nonun}(\DMod_{\mathscr{G}}(\tx{Conf}_{\infty \cdot c}))}(M, \CoStObj_{\tx{Conf}}^{\clambda}) \simeq \iota_{\mu \cdot x}^*(M);\]
\item All standard and co-standard objects are regular (resp.\ regular holonomic)  if $A$ is.
\end{itemize}

\subsection{Jacquet Functor is an Equivalence}

We will first need the following fact:

\begin{proposition}
For any smooth curve $X$, $\TheFactAlg$ is the \emph{trivial} local system when restricted to a single copy of the curve $X_{-\clambda}$ for $-\clambda$ a simple negative root.
\end{proposition}

\begin{proof}
We start by noting that $\TheFactAlg$ is an \emph{universal} factorization algebra in the sense of \cite{emily2019universal}: indeed, our construction evidently commutes with \'etale pullbacks along curves. This reduces us to the case of $X = \mb{A}^1$.
It is easy to verify that torus FLE is compatible with pullback along $\mb{A}^1 \to \mb{A}^1_\dR$, which implies that
\[\oblv_{\mb{A}^1_\dR \to \mb{A}^1}(\Omega_{\tx{KM}, \kappa}^{\tx{Lus}}) \in \IndCoh_{\DualTwisting}(([\tx{Gr}_{\check{T}}]_{\mb{A}^1_\dR})_\dR \times_{\mb{A}^1_\dR} \mb{A}^1)\]
is produced from
\[\oblv_{\mb{A}^1_\dR \to \mb{A}^1}(C^\semiinf(\tx{res}(\mb{V}_\kappa^0))) \in [\mc{KL}(T)_{\mathscr{L}_{[\kappa, \AnomalyTerm]}}]_{\mb{A}^1}\]
via FLE. Next note that the action of $\mb{A}^1$ on the curve itself induces a weak $\mb{A}^1$ action on both these two categories above, and this weak action is respected by torus FLE; moreover, the algebra $\oblv_{\mb{A}^1_\dR \to \mb{A}^1}(C^\semiinf(\tx{res}(\mb{V}_\kappa^0)))$ admits an weak $\mb{A}^1$-equivariance structure with respect to this action, as we have $[\mb{B} \Arc G]_{\mb{A}^1} \simeq \mb{B} \Arc G \times \mb{A}^1$ (and similarly for $[\widehat{\Loop G}]_{\mb{A}^1}$, etc). If we let $j_{-\clambda}: \mb{A}^1_{-\clambda} \to (([\tx{Gr}_{\check{T}}]_{\mb{A}^1_\dR})_\dR \times_{\mb{A}^1_\dR} \mb{A}^1)$ denote the inclusion of the $-\clambda$-component, then by this weak equivariance we know that $\oblv(j_{-\clambda}^!(\Omega_{\tx{KM}, \kappa}^{\tx{Lus}})) \in \IndCoh(\mb{A}^1_{-\lambda})$ is a constant coherent sheaf with $!$-fiber being $k$ everywhere (c.f. Proposition \ref{prop:KM-vacuum-verification}). As such it is a local system, but any local system on $\mb{A}^1$ is trivial.
\end{proof}

\begin{corollary}
\label{cor:omega-km-constructible}
The factorization algebra $\TheFactAlg$ is constructible along the diagonal stratification of the Ran space; in particular, it is ind-regular holonomic.
\end{corollary}
\begin{proof}
Follows formally from factorization property.
\end{proof}

Given the Proposition above, Section \ref{sect:conf-hw} tells us that $\StObj_{\tx{Conf}, \kappa}^{\clambda}$ and $\CoStObj_{\tx{Conf}, \kappa}^{\clambda}$ are well-defined objects for $\TheFactAlg$.
From now on, let $X = \mb{A}^1$ be our background curve, and let $x_0 \in X(\mb{C})$ denote the origin. Our main theorem is that
\begin{theorem}
\label{thm:jacquet-is-equivalence}
For $\kappa$ satisfying Condition \ref{cond:avoid-small-roots-of-unity}, $\tb{J}_{\kappa}$ is an equivalence of categories, and satisfies, for every $\clambda \in \weightLat$,
\[\tb{J}_{\kappa}(\StObj_{\tx{KM}, \kappa}^{\clambda}) \simeq \StObj_{\tx{Conf}, \kappa}^{\clambda} \hspace{1em} \tb{J}_{\kappa}(\CoStObj_{\tx{KM}, \kappa}^{\clambda}) \simeq \CoStObj_{\tx{Conf}, \kappa}^{\clambda}.\]
\end{theorem}
The proof of this will be delayed to Section \ref{sect:main-computation}.

\begin{remark}
At this point we \emph{do not yet know} what $\TheFactAlg$ is; in fact, we will only be able to identify it \emph{after} this theorem has been established.
\end{remark}

\subsection{Quantum Torus}
\label{sect:prep-q}

In this section, $q$ will denote a $\mb{C}^\times$-valued (multiplicative) quadratic form on $\weightLat$. In order to utilize the results established in \cite{gaitsgory2021factorization} we need some restrictions on $q$ imposed in \emph{loc.cit.}, namely that:
\begin{itemize}
\item[(non-degenerate)] $q(\check{\alpha}) \neq 1$ for every root $\check{\alpha}$;
\item[\smash{\tshortstack[r]{(avoids \\ small torsion)}}] For every simple factor $\mf{g}_i$ and every long root $\check{\alpha}_l$ in $\mf{g}_i$ we have $\tx{ord}(q(\check{\alpha}_l)) \ge d_i + 1$, where $d_i$ is the lacing number of $\mf{g}_i$;
\item[\smash{\tshortstack[r]{(almost \\ admissible)}}] For every simple factor $\mf{g}_i$, every short root $\check{\alpha}_{s}$ and every long root $\check{\alpha}_{l}$, we have
\[\tx{ord}(q(\check{\alpha}_{s})) \ge h_i - 1 \hspace{1em} \tx{ord}(q(\check{\alpha}_{l})) \ge h_i^\vee - 1,\]
where $h_i$ the Coxeter number of $\mf{g}_i$ and $h_i^\vee$ the dual Coxeter number of $\mf{g}_i$ (this condition appears in Section 2.3.1 of \emph{loc.cit.}).
\end{itemize}

\paragraph{Condition in \cite{kazhdanlusztig4}}
We make a quick detour to compare this restriction with the restriction imposed by \cite{kazhdanlusztig4}. Suppose our group is simple, and we are starting with a \emph{negative rational level}, given by a bilinear form
\[(c + h^\vee) \kappa_{\tx{min}} \hspace{1em} c = \frac{-p}{q} \in \mb{Q}^{< 0} \hspace{1em} (p, q) = 1;\]
here $\kappa_{\tx{min}} := \frac{1}{2 h^\vee} \tx{Kil}_\mf{g}$ is the minimal bilinear form on $\coweightLat$.
Then, tracing through definition, one sees that the quadratic form $q$ produced from our equivalence (extended to the negative level via duality) will be $\exp\left(\frac{- \pi i}{d c} \check{\kappa}_{\tx{min}}\right)$; here $d$ is the lacing number of $\mf{g}$, and $\check{\kappa}_{\tx{min}}$ is the minimal bilinear form for the Langlands dual group. We thus obtain
\[q(\check{\alpha})  = \begin{cases}
\exp\left(2 \pi i \frac{q}{d p}\right) & \check{\alpha}~\tx{short} \\ \\
\exp\left(2 \pi i \frac{q}{p}\right) & \check{\alpha}~\tx{long}
\end{cases}
\]
which results in Table \ref{tbl:minimal-value}.
In particular, one sees that the almost admissibility condition from Definition \ref{def:kappa-almost-admissible} is precisely the almost admissibility condition above.

\begin{convention}
\label{conv:q-small-roots}
For the rest of this subsection, we will assume the restrictions above on $q$ are met.
\end{convention}

\paragraph{Quantum Torus}
\emph{Choose} a bilinear form $b'$ such that $b'(\lambda, \lambda) = q(\lambda)$. We let $\Rep_{q, b'}(T)^{\heartsuit}$ denote the \emph{braided monoidal} Grothendieck abelian (1, 1)-category whose underlying monoidal category is the standard monoidal category $\Rep(T)^\heartsuit$ and whose braided structure is given by
\[k^\lambda \tensor k^\mu \xrightarrow{b'(\lambda, \mu)} k^\mu \tensor k^\lambda.\]
As explained in \emph{loc.cit.}, for any other choice $b''$, we have a \emph{canonical} braided monoidal equivalence
\[\varphi_{b', b''}: \Rep_{q, b'}(T)^\heartsuit \simeq \Rep_{q, b''}(T)^{\heartsuit};\]
moreover, this identification can be made compatible for any triple $(b', b'', b''')$. Since the category of braided monoidal (1, 1)-category is only a (2, 2)-category, the category is in fact independent of the choice of $b'$ and we are allowed to write simply $\Rep_q(T)^{\heartsuit}$.

Let $\tx{Ch}(\Rep_q(T)^{\heartsuit})$ be the category of unbounded chain complexes in $\Rep_q(T)^{\heartsuit}$. Equipped with the injective model structure \cite[1.3.5.3]{HA}), it is a braided monoidal model category. %
Observe that the category $\tx{Ch}(\Rep_q(T)^{\heartsuit})^{\tx{cofib}}$ of cofibrant objects clearly form a braided monoidal subcategory, the unit $k$ is cofibrant, and monoidal tensor with any object preserve weak equivalences.
It follows from \cite[5.1.2.4, 4.1.7.4]{HA} that
\[\mc{D}(\Rep_q(T)^\heartsuit) \simeq N(\tx{Ch}(\Rep_q(T))^\tx{cofib})[\tx{quasi-iso}^{-1}]\]
admits a canonical $\mb{E}_2$ DG category structure.

\begin{definition}
\label{def:quantum-torus}
	The \emph{quantum torus} $\Rep_q(T)$ is the $\mb{E}_2$ DG category constructed above.
\end{definition}

\subsubsection{Koszul Duality}
\label{sect:E2-koszul-duality}
By inverting the braiding (but keeping the monoidal structure intact), there is a canonical equivalence of $\mb{E}_2$ DG categories $\Rep_q(T) \simeq \Rep_{q^{-1}}(T)^{\tx{rev-br}}$, where the superscript on RHS reminds ourselves of the reversed braiding. We have a perfect pairing
\[\tensor^0: \Rep_q(T) \tensor \Rep_{q^{-1}}(T) \to \Vect\]
obtained by ind-extending the pairing
\[(c_1 \in \Rep_q(T)^c, c_2 \in \Rep_{q^{-1}}(T)^c) \mapsto \Hom_{\Vect^{\weightLat}}(k, c_1 \tensor c_2).\]
As explained in \cite[8.1]{gaitsgory2021conjectural}, this gives identifications of $\mb{E}_2$ DG categories
\begin{equation}
\Rep_{q^{-1}}(T) \simeq \Rep_q(T)^{\tx{rev}} \simeq \Rep_q(T)^\vee
\label{eq:Repqinv-is-rev}
\end{equation}
where the second category is obtained from $\Rep_q(T)$ by reversing the multiplication \emph{and} the braiding, and the third one inherits an $\mb{E}_2$ structure from the second one.

\newcommand{\convTensor}{\odot}

\begin{warning}
One should not think of $\Rep_q(T)^\tx{rev}$ as having a forgetful functor to $\Rep_q(T)$: indeed, the functor $\Rep_q(T)^\vee \to \Rep_q(T)$ is given by ind-extending the graded linear dual functor. The usage of ``rev'' here is in accordance with \cite[4.1.1.7]{HA}.
\end{warning}

Let $A$ be an \emph{augmented} $\mb{E}_2$-algebra internal to $\Rep_q(T)$, then \cite[5.2.5.38]{HA} defined the \emph{$\mb{E}_2$-Koszul dual} of $A$, which is the universal augmented $\mb{E}_2$-algebra
\[\tx{KD}_{\mb{E}_2}(A) \in \mb{E}_2\tx{-Alg}_\tx{aug}(\Rep_{q^{-1}}(T))\]
equipped with an $(A, \tx{KD}_{\mb{E}_2}(A)^\tx{rev})$-bimodule\footnote{Here the action is internal to the category of $\mb{E}_1$-algebras in $\Rep_q(T)$.} structure on $k \in \mb{E}_1\tx{-Alg}(\Rep_q(T))$ lifting the $A$-left module and $\tx{KD}_{\mb{E}_2}(A)^\tx{rev}$-right module structure coming from the augmentations.
From \cite[5.2.5]{HA} we know that $\tx{KD}_{\mb{E}_2}(A)$ exists for those $A$ compact within $\Rep_q(T)$ and is computed as $(\tx{Bar}^{(2)}(A)^\vee)^\tx{rev}$, where:
\begin{itemize}
\item $\tx{Bar}^{(2)}(A)$ is the $\mb{E}_2$-coalgebra obtained from applying the bar construction twice;
\item $\vee$ denotes the (right) dual object which obtains the structure of an $\mb{E}_2$-algebra in $\Rep_q(T)$;
\item $\tx{rev}: \mb{E}_2\tx{-Alg}_\tx{aug}(\Rep_q(T)) \simeq \mb{E}_2\tx{-Alg}_\tx{aug}(\Rep_{q^{-1}}(T))$ is the equivalence induced by $\ms{rev}$, the automorphism of $\mb{E}_2$-operad that reverses the last coordinate.
\end{itemize}

\begin{definition}
Recall that for any $\mb{E}_2$-category DG $\mc{C}$ there is an equivalence \cite[5.4.4.10]{HA}
\[\tx{AugIdeal}: \mb{E}_2\tx{-Alg}_{\tx{aug}}(\mc{C}) \simeq \mb{E}_2\tx{-Alg}_\nonun(\mc{C}): \bullet \oplus \tb{1}_{\mc{C}}\]
giving by taking augmentation ideal and formally adjoining the monoidal unit.
For $A \in \mb{E}_2\tx{-Alg}_{\nonun}(\Rep_q(T))$ we define
\[\tx{KD}_{\mb{E}_2}^\nonun(A) := \tx{AugIdeal}(\tx{KD}_{\mb{E}_2}(A \oplus k)).\]
\end{definition}

\paragraph{Bialgebras}
For any $\mb{E}_2$ category $\mc{V}$  we introduce the category of bi-augmented bialgebras
\[\tx{BiAlg}(\mc{V})_{\tx{biaug}} := \mb{E}_1\tx{-CoAlg}_{\tx{coaug}}(\mb{E}_1\tx{-Alg}_{\tx{aug}}(\mc{V}));\]
We have adjunctions
\[\tx{Bar}: \mb{E}_2\tx{-Alg}_{\tx{aug}}(\mc{V}) \adjoint \tx{BiAlg}_{\tx{biaug}}(\mc{V}): \tx{coBar} \hspace{1em} \tx{Bar}: \tx{BiAlg}_{\tx{biaug}}(\mc{V}) \adjoint \mb{E}_2\tx{-CoAlg}_{\tx{coaug}}(\mc{V}): \tx{coBar}\]
computable via the (co)bar construction.

Let $H \in \tx{BiAlg}_{\tx{biaug}}(\Rep_q(T))$ be a bialgebra within $\Rep_q(T)$, and assume $H$ is dualizable as an object of $\Rep_q(T)$. Then the dual object $H^\vee$ is again an element of $\mb{E}_1\tx{-Alg}_{\tx{aug}}(\mb{E}_1\tx{-CoAlg}_{\tx{coaug}}(\Rep_q(T))) \simeq \tx{BiAlg}_{\tx{biaug}}(\Rep_q(T))$. As a graded vector space, its $\clambda$-piece is the linear dual of the $(-\clambda)$-piece of $H$. The automorphism $\ms{rev}$ as above induces a map
\[\tx{rev-comult}: \tx{BiAlg}_{\tx{biaug}}(\Rep_q(T)) \simeq \tx{BiAlg}_{\tx{biaug}}(\Rep_{q^{-1}}(T))\]
which maintains the multiplication structure on a bialgebra but reverses the comultiplication structure.

\begin{lemma}
\label{lemma:KD-vs-coBar}
	Assume that $H$ lies in the heart of $\Rep_q(T)$, is compact, is concentrated in non-positive weights, has 1-dimensional zero weight component and is ind-conilpotent as a coalgebra. Then we have
	\[\tx{KD}_{\mb{E}_2}(\tx{coBar}(H)) \simeq \tx{coBar}(\tx{rev-comult}(H^\vee)) \in \mb{E}_2\tx{-Alg}_\tx{aug}(\Rep_{q^{-1}}(T)).\]
\end{lemma}

\begin{proof}
We apply the standard ``red-shift'' trick: namely, consider the braided monoidal auto-equivalence
\[\tx{Shift}: \Rep_q(T) \to \Rep_q(T) \hspace{1em} \left(A = \bigoplus_{\clambda} A^\clambda\right) \mapsto \bigoplus_{\clambda} A^{\clambda}[2 \tx{ht}(\clambda)]\]
Our assumption on $H$ guarantees that $\tx{Shift}(H)$ is a \emph{connected} conilpotent bialgebra (i.e. its coaugmentation coideal is concentrated in purely positive cohomological degrees). In such case, \cite[Corollary 0.2]{ginot2016deformation} (more precisely, its evident generalization to $\weightLat$-graded vector spaces) guarantees that $\tx{Bar} \circ \tx{coBar}(\tx{Shift}(H)) \simeq \tx{Shift}(H)$, from which we know that $\tx{Bar} \circ \tx{coBar}(H) \simeq H$ as well.
Thus we compute
\[\tx{KD}_{\mb{E}_2}(\tx{coBar}(H)) \simeq (\tx{Bar}^{(2)}(\tx{coBar}(H))^\vee)^\tx{rev} \simeq (\tx{Bar}(H)^\vee)^\tx{rev} \simeq \tx{coBar}(\tx{rev-comult}(H^\vee))\]
where the last step follows from definition.
\end{proof}

\subsection{Quantum Groups}

\begin{definition}[Quantum Groups]
\label{def:lus-and-kd}
Let $U_q^\tx{free}(N^-)$ denote the bialgebra internal to $\Rep_q(T)$ generated by $\{f_{i}\}$, where $i$ ranges over the set of negative simple roots $\{-\alpha_i\}$, and $f_i$ is of degree $-\alpha_i$. It lies in the heart and has a \emph{Hopf algebra} structure, and has comultiplication given by
\[f_i \mapsto f_i \tensor 1 + 1 \tensor f_i.\]
\emph{Assuming the condition on $q$ above}, we define $U_q^\tx{KD}(N)$ be the quotient of $U_q^{\tx{free}}(N^-)$ by the two-sided ideal generated by the quantum Serre relations as defined in \cite[Section 3.2.7]{gaitsgory2021factorization}.

We define $U_q^{\tx{Lus}}(N)$ as the graded dual of $U_q^\tx{KD}(N^-)$; it is again a Hopf algebra in $\Rep_q(T)$.
\end{definition}

\begin{definition}
Given any algebra $A$ in $\Rep_q(T)$, we set
\[A\mod(\Rep_q(T))_{\tx{loc.nilp}}\]
as the ind-completion of the category of finite-dimensional $A$-modules in $\Rep_q(T)$.
\end{definition}

\begin{remark}
If $A$ is a bialgebra, then $A\mod(\Rep_q(T))_{\tx{loc.nilp}}$ gains the structure of an $\mb{E}_1$-algebra \emph{internal to} $\Rep_q(T)\mod^r(\DGCat)$.
\end{remark}

\begin{definition}
\label{def:rel-drinfeld-center}
We define
\[\tx{Rep}_q(B) := U_q^\tx{Lus}(N)\mod(\tx{Rep}_q(T))_{\tx{loc.nilp}}.\]
and its relative Drinfeld center
\[\tx{Rep}_q^\tx{mxd}(G) := Z_{\tx{Dr}, \tx{Rep}_q(T)}(\tx{Rep}_q(B)) := \tx{HH}_{\mb{E}_2, /\Rep_q(T)}(A\mod(\Rep_q(T))_{\tx{loc.nilp}})\]
via Definition \ref{def:En-HH-definition}.
\end{definition}

\begin{remark}
In other words, $\Rep_q^\tx{mxd}(G)$ is the universal $\mb{E}_2$-category acting from the left on $\Rep_q(B)$, in a way that is compatible with the $\Rep_q(T)$-action on the right. For an even more explicit description, see \cite[Section 4.4.2]{gaitsgory2021conjectural}.
\end{remark}

By Theorem \ref{thm-HH-E_n}, $\Rep_q^\tx{mxd}(G)$ coincides with $\Omega_q\mod^{\mb{E}_2}(\Rep_q(T))$, where $\Omega_q := \tx{coBar}(U_q^{\tx{KD}}(N^-))$. There is an obvious forgetful functor
\[\oblv_{\tx{mxd} \to \tx{Lus}}: \Rep_q^\tx{mxd}(G) \to \Rep_q(B);\]
whereas the following is less evident:
\begin{lemma}
\label{lemma:oblv-mxd-to-KD}
The forgetful functor $\Rep_q^\tx{mxd}(G) \xrightarrow{\oblv_{\tx{mxd} \to \tx{Lus}}} \Rep_q(B) \xrightarrow{\oblv} \Rep_q(T)$ factors through another forgetful functor
\[\oblv_{\tx{mxd} \to \tx{KD}^-}: \Rep_q^\tx{mxd}(G) \to U_q^\tx{KD}(N^-)\mod(\Rep_q(T)).\]
\end{lemma}
The proof of this will be delayed until Section \ref{ssec-enveloping-algebras}.

\begin{remark}
As explained in \cite{gaitsgory2021conjectural}, both $\Rep_q(B)$ and $\Rep_q^\tx{mxd}(G)$ are equipped with canonical $t$-structures, such that heart of $\tx{Rep}_q^\tx{mxd}(G)$ can be explicitly described as graded vector spaces equipped with a locally nilpotent action of $U_q^\tx{Lus}(N)$ and another action of $U_q^\tx{KD}(N^-)$. We will not need this description in this paper.
\end{remark}

\begin{definition}[Definition-Construction]
Let $\tx{ind}_{\tx{Lus} \to \tx{mxd}}: \Rep_q(B) \to \Rep_q^\tx{mxd}(G)$ denote the left adjoint of $\oblv_{\tx{mxd} \to \tx{Lus}}$, and let $\tx{coind}_{\tx{KD}^- \to \tx{mxd}}$ denote the right adjoint of $\oblv_{\tx{mxd} \to \tx{KD}^-}$.

We define the quantum \emph{Verma} modules $\mb{M}_{q}^{\clambda, \tx{mxd}}$, for each $\clambda \in \weightLat$, as
\[\mb{M}_{q, \tx{mxd}}^{\clambda} := \tx{ind}_{\tx{Lus} \to \tx{mxd}}(k^{\clambda}).\]
and the quantum dual Verma modules $\mb{M}_{q}^{\clambda, \tx{mxd}, \vee}$ as
\[\mb{M}_{q, \tx{mxd}}^{\clambda, \vee} := \tx{coind}_{\tx{KD}^- \to \tx{mxd}}(k^{\clambda}).\]
\end{definition}

It then follows from definition that these objects satisfy the right orthogonal property, and that $\mb{M}_q^{\clambda, \tx{mxd}, \vee}$ form a set of compact generators of this category.

\begin{definition}[Contragredient Duality]
\label{def:mixed-contragredient-duality}
As shown in \cite[Section 8.3.1]{gaitsgory2021conjectural}, the graded linear duality on $\Rep_q(T)$ induces a duality
\[\Rep_q^\tx{mxd}(G)^\vee \simeq \Rep_{q^{-1}}^\tx{mxd}(G),\]
in a way compatible with the forgetful functor down to the quantum torus. We will let $\mb{D}_{\tx{Quant}}$ denote the resulting contravariant functor on compact objects. As expected, we have
\[\mb{D}_{\tx{Quant}}(\mb{M}_q^{\clambda, \tx{mxd}}) \simeq \mb{M}_{q^{-1}, \tx{mxd}}^{-\clambda - 2\crho}[\dim(\mf{n})].\]
\end{definition}

\subsection{Categorical Verdier Duality}

\begin{convention}
	Let $\mb{E}_2$ be the (rectlinear) little 2-disk operad introduced in \cite[5.1.0.2]{HA}, and let $(\mb{E}_2)_{\tx{BSO}(2)}$ be the variant introduced in \cite[5.4.2.16]{HA} where we allow orientation-preserving rotation of the disk. By considering the canonical tangent $\tx{SO}(2)$-bundle, a $(\mb{E}_2)_{\tx{BSO}(2)}$-algebra structure gives rise to an $\mb{E}_M$-algebra structure (defined in \cite[5.4.5.1]{HA}) for any smooth 2-dimensional manifold $M$.
\end{convention}

\begin{definition}
By an $\mb{E}_2$ DG category $\mc{C}$, we mean an element of the category $\mb{E}_2\tx{-Alg}_\nonun(\DGCat)$. Likewise we have the notion of a $(\mb{E}_2)_{\tx{BSO}(2)}$ DG category and a $\mb{E}_M$ DG category.
\end{definition}

Let $X$ be a smooth algebraic curve over $\mb{C}$. In \cite[Section 29.4]{gaitsgory2019metaplectic}, authors of \emph{loc.cit}.\ appealed to (without proof) the following folklore result:
\begin{itemize}
\item Given a $(\mb{E}_2)_{\tx{BSO}(2)}$ DG category $\mc{C}$, one can (by considering the tangent bundle) obtain a $\mb{E}_{X(\mb{C})}$ DG category then a \emph{topological} factorization category $\tx{Fact}(\mc{C})$ over $X(\mb{C})$;
\item Moreover, to each $(\mb{E}_2)_{\tx{BSO}(2)}$ algebra $A$ internal to $\mc{C}$, one can attach first a $\mb{E}_{X(\mb{C})}$-algebra then a factorization algebra $\tx{Fact}(A)$ within $\tx{Fact}(\mc{C})$, and there is an equivalence of DG categories between the category of $\mb{E}_{X(\mb{C})}$-modules of $A$ inside $\mc{C}$ and (non-unital) factorization modules of $\tx{Fact}(A)$ supported at one point inside $\tx{Fact}(\mc{C})$.
\end{itemize}

In Section \ref{sect:factvse2} we give a full proof of this result; the precise statement is given in Proposition \ref{prop:E2-to-fact-cat}, Corollary \ref{cor:E2Mod-to-FactModCat} and Proposition \ref{prop:E2Alg-Mod-to-FactAlg-Mod}.

\begin{remark}
One categorical level down, the strategy is clear: first apply \cite[5.5.4.10]{HA}, then apply covariant Verdier duality. On the categorical level there is no Verdier duality, so the crux of the proof is to construct a functor from the latter to the former.
\end{remark}

\begin{remark}
	The concept of an (algebraic) factorization category, a topological factorization category, and an $\mb{E}_2$-category are \emph{all inequivalent} concepts. For intuition purpose one can sometimes think of the first two as various kinds of \emph{lax} $\mb{E}_2$-structure; we do not attempt to make this intuition precise.
\end{remark}

Our input $\mb{E}_2$ category is $\Rep_q(T)$. In Section \ref{sect:q-torus}, we further identify the output with the topological factorization category of $q$-twisted analytic sheaves over the $\mb{C}$-points of $\tx{Gr}_{\check{T}}^{\omega^{\crho}}$. We let $\Omega_{\mb{E}_2, q}^{\tx{Lus}}$ denote the topological factorization algebra obtained from (the augmentation ideal of) $\Omega_q \in \mb{E}_2\tx{-Alg}(\Rep_q(T))$. It thus remains to show:

\begin{proposition}
\label{prop:KM-E2-Lus-match}
	$\TheFactAlg$ and $\Omega_{\mb{E}_2, q}^{\tx{Lus}}$ correspond to each other under the Riemann-Hilbert correspondence (defined in Section \ref{sect:rh}).
\end{proposition}

The proof of this statement will again be delayed until Section \ref{sect:main-computation}. As explained in the introduction, this finishes the proof of our main theorem.

\begin{remark}
By this comparison result and the $*$-fiber computation from Proposition \ref{prop:KM-vacuum-verification}, we obtain that
\[C_*(U_{q^{-1}}^{\tx{KD}}(N^-), k) \simeq \begin{cases}
\mb{C}[-\ell(w_f)] & \clambda = w_f(\crho) - \crho, w_f \in W_\tx{fin} \\
0 & \tx{else}
\end{cases}\]
holds true for all levels that satisfy Condition \ref{cond:avoid-small-roots-of-unity}. For $q$ a root of unity of odd order (along with some slightly more stringent conditions) this was established purely algebraically in \cite{kostant-vigre}; we note that our result also provides a proof for the even order case that was not covered by \emph{loc.cit}.
\end{remark}

\begin{remark}
It follows from construction that the equivalence (at positive level) sends $\StObj_{\tx{KM}, \kappa}^{\clambda}$ to $\mb{M}_{q, \tx{mxd}}^{\clambda}$, $\CoStObj_{\tx{KM}, \kappa}^{\clambda}$ to $\mb{M}_{q, \tx{mxd}}^{\clambda, \vee}$ and intertwines $C_*^{\semiinf}(\Loop \mf{n}, \Arc N, -)$ with $C_*(U_q^\tx{Lus}(N), -)$.
\end{remark}

\section{Factorizable Kazhdan--Lusztig Category}

\label{sect:fact-kl-category}

We begin by developing the factorizable Kazhdan--Lusztig category. Throughout this section we will freely use results from Section \ref{sect:appendix-indcoh}.

\subsection{The Players}

\begin{notation} \label{notn-def-KL}
In this section, we use the following notations:
\begin{itemize}
	\item $X$ is a smooth separated algebraic curve, and $\Ran_{\on{un}}:=\Ran_{\on{un}}(X)$ is the corresponding unital Ran space;
	\item $H$ is a finite type affine algebraic group;
	\item $\mathcal{L} H_{\Ran_{\on{un,dR}}}$ is the corr-unital factorization group whose fiber at a closed point $x\in X$ is the loop group $\mathcal{L} H_x$ (see Section \ref{sect:factspc});
	\item  $\mathcal{L}^+ H_{\Ran_{\on{un,dR}}}$ is the co-unital factorization group whose fiber at a closed point $x\in X$ is the arc group $\mathcal{L}^+ H_x$ (see Section \ref{sect:factspc});
	\item We have a strictly \emph{unital} morphism $\mathcal{L}^+ H_{\Ran_{\on{un,dR}}}\to \mathcal{L} H_{\Ran_{\on{un,dR}}}$;
	\item  $[\on{Gr}_H]_{\Ran_{\on{un,dR}}}$ is the unital factorization space whose fiber at a closed point $x\in X$ is the affine Grassmannian $[\on{Gr}_{H}]_x\simeq \mathcal{L} H_x/\mathcal{L}^+ H_x$ (see Section \ref{sect:factspc});
	\item We have a strictly \emph{co-unital} morphism $\mathcal{L} H_{\Ran_{\on{un,dR}}}\to [\on{Gr}_H]_{\Ran_{\on{un,dR}}}$;

	\item $\widehat{\mathcal{L}H}_{\Ran_{\on{un,dR}}}$ is the formal completion of $\mathcal{L} H_{\Ran_{\on{un,dR}}}$ along $\mathcal{L}^+ H_{\Ran_{\on{un,dR}}}$, which inherits a corr-unital factorization structure whose fiber at a closed point $x\in X$ is  $\widehat{\mathcal{L} H}_x \simeq  \mathcal{L}H_x \times_{[\mathcal{L}H_x]_\dR} [\mathcal{L}^+H_x]_\dR$;

	\item $[\widehat{\on{Gr}}_H]_{\Ran_{\on{un,dR}}}$ is the formal completion of $[\on{Gr}_H]_{\Ran_{\on{un,dR}}}$ along its unit section, which inherits a unital factorization structure whose fiber at a closed point $x\in X$ is  $[\widehat{\on{Gr}}_H]_x \simeq \widehat{\mathcal{L} H}_x/ \mathcal{L}^+ H_x $;
	\item $[\mathbb{B}\mathcal{L}^+ H]_{\Ran_{\on{un,dR}}}$ is the classifying space of $\mathcal{L}^+ H_{\Ran_{\on{un,dR}}}$, which inherits a co-unital factorization structure whose fiber at a closed point $x\in X$ is the quotient stack $\mathbb{B}\mathcal{L}^+ H_x$;
	\item $[\widehat{\on{Hecke}}_H^\bullet ]_{\Ran_{ \on{un,dR} }}$ is the Cech nerve of $ [\mathbb{B}\mathcal{L}^+ H]_{\Ran_{\on{un,dR}}} \to [\mathbb{B}\widehat{\mathcal{L}H}]_{\Ran_{\on{un,dR}}}$, which is a cosimplicial object in the category of corr-unital factorization spaces and strictly unital morphisms. Explicitly,
	$$ [\widehat{\on{Hecke}}_H^m ]_{\Ran_{ \on{un,dR} }} := [ \mathcal{L}^+ H\backslash \widehat{\mathcal{L}H} \times^{ \mathcal{L}^+ H } \widehat{\mathcal{L}H} \times^{ \mathcal{L}^+ H } \cdots \times^{  \mathcal{L}^+ H  }  \widehat{\mathcal{L}H}/  \mathcal{L}^+ H]_{ \Ran_{\on{un,dR}} }, $$
	where there are $m$ terms of $\widehat{\mathcal{L}H}$ on the RHS. In the degenerate case $m=0$, 
 	$$ [\widehat{\on{Hecke}}_H^0 ]_{\Ran_{ \on{un,dR} }} := [\mathbb{B}\mathcal{L}^+ H]_{\Ran_{\on{un,dR}}}.$$
\end{itemize}
\end{notation}

\begin{remark} \label{rem-description-unital-structure-hecke}
The following description is helpful for understanding the corr-unital structure on $[\widehat{\on{Hecke}}_H^m ]_{\Ran_{ \on{un,dR} }}$. Let $\phi:I\to J$ be an injective map between finite sets. Recall the above corr-unital structure should provide the following correspondence
$$  [\widehat{\on{Hecke}}_H^m ]_{J} \gets [\widehat{\on{Hecke}}_H^m ]_{\phi}  \to [\widehat{\on{Hecke}}_H^m ]_{I}\times_{X^I} X^J$$
between (convergent) prestacks. Unwinding the definitions, for an affine test scheme $S\in \,^{>-\infty}\on{AffSch}$, the groupoid $[\widehat{\on{Hecke}}_H^m ]_{\phi}(S)$ classifies the following data (i)-(iv):
\begin{itemize}
		\item[(i)] Maps $x_j:S\to X$ labelled by $J$. 
		\item[]	Let $\Gamma_J\to S\times X$ (resp. $\Gamma_I\to S\times X$) be the schema-theoretic sum of the graphs of all the $x_j$ (resp. of $x_{\phi(i)}$ for $i\in I$). Note that we have a closed embedding $\Gamma_I \to \Gamma_J$.
		\item[(ii)] $H$-torsors $\mathcal{P}_0,\cdots,\mathcal{P}_m$ on the adic multi-disk $\mathcal{D}_{\Gamma_J}'$ and isomorphisms 
		$$  [\mathcal{P}_0]|_{ \mathcal{D}_{\Gamma_J}'\setminus \Gamma_I } \simeq \cdots \simeq [\mathcal{P}_m]|_{ \mathcal{D}_{\Gamma_J}'\setminus \Gamma_I } $$
		between their restrictions on $\mathcal{D}_{\Gamma_J}'\setminus \Gamma_I $.

		\item[(iii)] Isomorphisms 
		$$  [\mathcal{P}_0]\times_S S^{\on{cl,red}} \simeq \cdots \simeq [\mathcal{P}_m]\times_S S^{\on{cl,red}}, $$
		where we view $[\mathcal{P}_k]\times_S S^{\on{cl,red}}$ as a $G$-torsor on $\mathcal{D}_{\Gamma_J}'\times_S S^{\on{cl,red}} $.

		\item[] By restriction, (ii) and (iii) provides two chains of isomorphisms
		$$ [\mathcal{P}_0]|_{ \mathcal{D}_{\Gamma_J}'\setminus \Gamma_I } \times_S S^{\on{cl,red}} \simeq \cdots \simeq [\mathcal{P}_m]|_{ \mathcal{D}_{\Gamma_J}'\setminus \Gamma_I } \times_S S^{\on{cl,red}}.$$

		\item[(iv)] An isomorphism between the above two chains.	
\end{itemize}
Recall $[\widehat{\on{Hecke}}_H^m ]_{J}=[\widehat{\on{Hecke}}_H^m ]_{\on{Id}_J}$. Now the correspondence
$$  [\widehat{\on{Hecke}}_H^m ]_{J}(S) \gets [\widehat{\on{Hecke}}_H^m ]_{\phi}(S)  \to [\widehat{\on{Hecke}}_H^m ]_{I}\times_{X^I} X^J(S)$$
is induced by the obvious op-correspondence
$$ \mathcal{D}_{\Gamma_J}'\setminus \Gamma_J \to  \mathcal{D}_{\Gamma_J}'\setminus \Gamma_I  \gets  \mathcal{D}_{\Gamma_I}'\setminus \Gamma_I.$$
\end{remark}

\begin{remark} In the special case when $m=1$ and $I=\emptyset$, the correspondence
$$  [\widehat{\on{Hecke}}_H^m ]_{J} \gets [\widehat{\on{Hecke}}_H^m ]_{\phi}  \to [\widehat{\on{Hecke}}_H^m ]_{I}\times_{X^I} X^J$$ 
is just
$$   [ \mathcal{L}^+ H\backslash \widehat{\mathcal{L}H}/\mathcal{L}^+ H ]_{J} \gets [\mathbb{B}\mathcal{L}^+ H ]_{J} \to X^J.$$
\end{remark}

We need the following techinical lemmas:

\begin{lemma} \label{lem-Hecke-is-nice-enough}
$[\widehat{\on{Hecke}}_H^m ]_{\Ran_{ \on{un,dR} }}$ satisfies $(\spadesuit)$ in Construction \ref{constr-unit-fact-prestack-indcoh-fact-cat}.
\end{lemma}

\proof By Remark \ref{rem-constr-unit-fact-prestack-indcoh-fact-cat} and Remark \ref{rem-description-unital-structure-hecke}, we only need to check 
\begin{itemize}	
	\item[(a)] $[\widehat{\on{Hecke}}_H^m ]_{\phi}$ is of quotient type over $X^J$;

	\item[(b)] $[\widehat{\on{Hecke}}_H^m ]_{\phi} \to [\widehat{\on{Hecke}}_H^m ]_{I}\times_{X^I} X^J$ is $*$-pullable;

	\item[(c)] $[\widehat{\on{Hecke}}_H^m ]_{\phi} \to [\widehat{\on{Hecke}}_H^m ]_{J}$ is an ind-(afp, affine and proper) apaisant indshematic morphism.
\end{itemize}
Consider the prestack $[\widehat{\on{Gr}}_H^{\times m} ]_{\phi}$ that classifies data (i)-(iv) as in Remark \ref{rem-description-unital-structure-hecke} and 
\begin{itemize}
	\item[(v)] A trivialization of $\mathcal{P}_0$. 
\end{itemize}
We have
$$ [\widehat{\on{Hecke}}_H^m ]_{\phi}\simeq  [\mathbb{B}\mathcal{L}^+ H]_I\backslash [\widehat{\on{Gr}}_H^{\times m} ]_{\phi} . $$
Note that when $\phi = \on{Id}_J$, $[\widehat{\on{Gr}}_H^{\times m} ]_{J}$ is just the formal completion of $[\on{Gr}_H]_J \times_{X^J} \cdots \times_{X^J} [\on{Gr}_H]_J$ along its unit section. Moreover, using Beauville-Laszlo descent theorem, one can show that the map
$$ [\widehat{\on{Gr}}_H^{\times m} ]_{\phi} \to [\widehat{\on{Gr}}_H^{\times m} ]_{I}\times_{X^I} X^J $$
is an equivalence. This implies $[\widehat{\on{Gr}}_H^{\times m} ]_{\phi}$ is a classical laft ind-scheme which is ind-(affine and proper) over $X^J$. As in Example \ref{exam-hecke-apaisant-quot-stack}, this implies (a).

(b) follows from the fact that 
$$[\widehat{\on{Gr}}_H^{\times m} ]_{\phi} \simeq [\widehat{\on{Gr}}_H^{\times m} ]_{I}\times_{X^I} X^J$$
is a flat apaisant cover of both $[\widehat{\on{Hecke}}_H^m ]_{\phi}$ and $ [\widehat{\on{Hecke}}_H^m ]_{I}\times_{X^I} X^J $.

(c) follows from the fact that 
$$ [\widehat{\on{Gr}}_H^{\times m} ]_{I}\times_{X^I} X^J\to [\widehat{\on{Gr}}_H^{\times m} ]_{J} $$
is a closed embedding.

\qed

\begin{lemma} \label{lem-Hecke-is-nice-enough-2}
For any $[m]\to [n]$ in $\Delta$ and any $\phi:I\to J$ between finite sets, the morphism
$$ [\widehat{\on{Hecke}}_H^n ]_{\phi }  \to [\widehat{\on{Hecke}}_H^m ]_{\phi }$$
is $(?,\ren)$-pullable.
\end{lemma}

\proof We can assume $[m]\to [n]$ sends $0$ to $0$. Then the proof of Lemma \ref{lem-Hecke-is-nice-enough} implies the desired map can be obtained from
$$ [\widehat{\on{Gr}}_H^n ]_{I }\times_{X^I} X^J  \to [\widehat{\on{Gr}}_H^m ]_{I }\times_{X^I} X^J $$
by taking $[\mathcal{L}^+ H]_J$-quotient. This makes the claim manifest.

\qed

\subsection{The Non-Twisted Version}
\label{ssec-def-KL-non-twisted}

\begin{definition}\label{def-KL-non-twisted}
 By Lemma \ref{lem-Hecke-is-nice-enough}, Lemma \ref{lem-Hecke-is-nice-enough-2}, Construction \ref{constr-unit-fact-prestack-indcoh-fact-cat} and Construction \ref{constr-functorial-unital-fact-indcoh}, we obtain simplicial objects
$$ \mathcal{I}nd\mathcal{C}oh_{*,\ren}( \widehat{\on{Hecke}}_H^\bullet ), \, \mathcal{I}nd\mathcal{C}oh^{!,\ren}( \widehat{\on{Hecke}}_H^\bullet ) \in \on{Fun}(\Delta^\op, \on{FactCat}_{\on{un}}) $$
whose connecting morphisms are \emph{strictly unital} factorization functors, and are given respectively by the $*$-pushforward functors and the left adjoint of the $!$-pullback functors. We define:
\begin{eqnarray*}
 \mathcal{KL}(H) &:=\mathcal{L}\mathfrak{h}\on{-}mod^{\mathcal{L}^+H} &:=\colim_{\Delta^\op, {!\on{-push}} } \mathcal{I}nd\mathcal{C}oh^{!,\ren}( \widehat{\on{Hecke}}_H^\bullet ), \\
 \mathcal{KL}(H)_{\on{co}} &:=  \mathcal{L}\mathfrak{h}\on{-}mod^{\mathcal{L}^+H}_{\on{co}} &:= \colim_{\Delta^\op, *\on{-push} } \mathcal{I}nd\mathcal{C}oh_{*,\ren}( \widehat{\on{Hecke}}_H^\bullet ),
\end{eqnarray*}
where the colimits are taken in $\on{FactCat}_{\on{un}}$. We refer them respectively as the \emph{non-twisted unital factorization Kazhdan--Lusztig category} and the \emph{co-version of the non-twisted unital factorization Kazhdan--Lusztig category}.
\end{definition}

Several explanations are in order.

\begin{remark} \label{rem-Rep^!=Rep^*}
The unital factorization categories $\mathcal{I}nd\mathcal{C}oh^{!,\ren}( \widehat{\on{Hecke}}_H^0 )$ and $\mathcal{I}nd\mathcal{C}oh_{*,\ren}( \widehat{\on{Hecke}}_H^0 )$ are canonically equivalent. In other words:
$$ \mathcal{I}nd\mathcal{C}oh^{!,\ren}( \mathbb{B}\mathcal L^+ H ) \simeq \mathcal{I}nd\mathcal{C}oh_{*,\ren}( \mathbb{B}\mathcal L^+ H ).$$
To construct this equivalence, we construct a self-duality on $\mathcal{I}nd\mathcal{C}oh^{!,\ren}( \mathbb{B}\mathcal L^+ H )$ in $[\on{FactCat}_{\on{un}}]_{\on{lax-unital}}$ as follows. One can check directly that the following morphisms provide a self-duality\footnote{In the definition of the second functor, we used the fact that the diagonal map $ \mathbb{B}\mathcal L^+ H\to\mathbb{B}\mathcal L^+ H\times \mathbb{B}\mathcal L^+ H $ are $?$-pushable (or equivalently, $*$-pullable), which is not true for $\widehat{\on{Hecke}}_H^m$ with $m>0$. Hence we can not obtain a canonical isomorphism between $\mathcal{I}nd\mathcal{C}oh^{!,\ren}( \widehat{\on{Hecke}}_H^m )$ and $\mathcal{I}nd\mathcal{C}oh_{*,\ren}( \widehat{\on{Hecke}}_H^m)$ using \emph{this} method. Instead, in the proof of Theorem \ref{thm-tate-twist}, we will construct $(m+1)$ different isomorphisms between them.}:
\[
\begin{aligned}
	\mathcal{I}nd\mathcal{C}oh^{!,\ren}( \mathbb{B}\mathcal L^+ H ) \otimes \mathcal{I}nd\mathcal{C}oh^{!,\ren}( \mathbb{B}\mathcal L^+ H ) 
	 \xrightarrow{\otimes^!} \mathcal{I}nd\mathcal{C}oh^{!,\ren}( \mathbb{B}\mathcal L^+ H)\xrightarrow{?-\on{push}}  \mathcal{I}nd\mathcal{C}oh^{!,\ren}(\on{pt}) \simeq \mathbf{1},\\
	 \mathbf{1}\simeq  \mathcal{I}nd\mathcal{C}oh^{!,\ren}(\on{pt})  \xrightarrow{!-\on{pull}}  \mathcal{I}nd\mathcal{C}oh^{!,\ren}( \mathbb{B}\mathcal L^+ H) 
	 \xrightarrow{?-\on{push}} \mathcal{I}nd\mathcal{C}oh^{!,\ren}( \mathbb{B}\mathcal L^+ H \times  \mathbb{B}\mathcal L^+ H)\simeq \\ \simeq \mathcal{I}nd\mathcal{C}oh^{!,\ren}( \mathbb{B}\mathcal L^+ H ) \otimes \mathcal{I}nd\mathcal{C}oh^{!,\ren}( \mathbb{B}\mathcal L^+ H ).
\end{aligned}
\]

By construction, via the above equivalence, the chain
$$\mathcal{I}nd\mathcal{C}oh^{!,\ren}( \on{pt} ) \xrightarrow{!\on{-pull}} \mathcal{I}nd\mathcal{C}oh^{!,\ren}( \mathbb{B}\mathcal L^+ H ) \xrightarrow{!\on{-pull}} \mathcal{I}nd\mathcal{C}oh^{!,\ren}( \on{pt} ) $$
goes to
$$\mathcal{I}nd\mathcal{C}oh_{*,\ren}( \on{pt} ) \xrightarrow{*\on{-pull}} \mathcal{I}nd\mathcal{C}oh_{*,\ren}( \mathbb{B}\mathcal L^+ H )  \xrightarrow{*\on{-pull}} \mathcal{I}nd\mathcal{C}oh_{*,\ren}( \on{pt} ) .$$
\end{remark}

\begin{notation} By Remark \ref{rem-Rep^!=Rep^*}, we use
$$ \mathcal{R}ep( \mathcal{L}^+ H )$$
to denote both $\mathcal{I}nd\mathcal{C}oh^{!,\ren}( \mathbb{B}\mathcal L^+ H )$ and $\mathcal{I}nd\mathcal{C}oh_{*,\ren}( \mathbb{B}\mathcal L^+ H )$.
\end{notation}

\begin{remark} By definition, for any closed point $x\in X$, the fiber of $\mathcal{R}ep( \mathcal{L}^+ H )$ at $x$
$$\on{Rep}( \mathcal{L}^+H_x ):= \mathcal{R}ep( \mathcal{L}^+ H )_x \simeq \IndCoh_{*,\ren}( \mathbb{B} \mathcal{L}^+H_x )  $$
is compactly generated, and the heart $\Coh( \mathbb{B} \mathcal{L}^+H_x )^\heartsuit$ of its full subcategory of compact objects can be identified with the category of finite-dimentional representations of $\mathcal{L}^+H_x$. Note that our notation $\on{Rep}( \mathcal{L}^+H_x )$ is compatible with that in \cite{raskin2020homological}.
\end{remark}

\begin{remark} \label{rem-dual-!-version-*-version-KL}
 $\mathcal{KL}(H)$ and $\mathcal{KL}(H)_{\on{co}}$ are indeed canonically dual to each other in the symmetric monoidal category $[\on{FactCat}_{\on{un}}]_{\on{lax-unital}}$ of unital factorization categories and lax-unital factorization functors. See Remark \ref{rem-duality-unital-factorization-category} for what this means. This follows formally from Remark \ref{rem-section-KL-compact-gen} below.
\end{remark}

\begin{remark} \label{rem-section-KL-compact-gen}
The sections of $\mathcal{KL}(H)$ and $\mathcal{KL}(H)_{\on{co}}$ over $f_\dR: T_\dR \to \Ran_{\on{un,dR}}$ are given by
\begin{eqnarray} \label{eqn-KL-as-limit-!}
 \mathbf{\Gamma}(f_\dR, \mathcal{KL}(H)) &\simeq& \colim_{\Delta^\op, !\on{-push} } \mathbf{\Gamma}(f_\dR,  \mathcal{I}nd\mathcal{C}oh^{!,\ren}( \widehat{\on{Hecke}}_H^\bullet )) \simeq \lim_{\Delta, {!\on{-pull}} } \mathbf{\Gamma}(f_\dR,  \mathcal{I}nd\mathcal{C}oh^{!,\ren}( \widehat{\on{Hecke}}_H^\bullet )) \\
 \label{eqn-KL-as-limit-*}
 \mathbf{\Gamma}(f_\dR, \mathcal{KL}(H)_{\on{co}}) &\simeq& \colim_{\Delta^\op, *\on{-push} } \mathbf{\Gamma}(f_\dR,  \mathcal{I}nd\mathcal{C}oh_{*,\ren}( \widehat{\on{Hecke}}_H^\bullet) \simeq \lim_{\Delta, ?\on{-pull} } \mathbf{\Gamma}(f_\dR,  \mathcal{I}nd\mathcal{C}oh_{*,\ren}( \widehat{\on{Hecke}}_H^\bullet).
\end{eqnarray}
Therefore they are both compactly generated and canonically dual to each other. It follows that we also have
\begin{eqnarray*}
 \mathcal{KL}(H) &\simeq& \lim_{\Delta^\op, {!\on{-pull}} } \mathcal{I}nd\mathcal{C}oh^{!,\ren}( \widehat{\on{Hecke}}_H^\bullet ), \\
 \mathcal{KL}(H)_{\on{co}} &\simeq& \lim_{\Delta^\op, {?\on{-pull}} } \mathcal{I}nd\mathcal{C}oh_{*,\ren}( \widehat{\on{Hecke}}_H^\bullet ),
\end{eqnarray*}
where the limit is taken in $[\on{FactCat}_{\on{un}}]_{\on{lax-unital}}$.
\end{remark}

\begin{remark} One should think $\mathcal{KL}(H)$ as the unital factorization category for the $\IndCoh^{!,\ren}$-theory on $[\mathbb{B}\widehat{\mathcal{L}H}]_{\Ran_{\on{un,dR}}}$. However, since $\widehat{\mathcal{L}H}$ is only a group \emph{ind}-scheme, we do not have, or rather have not set up, a well-behaved theory of ind-coherent sheaves on its classifying prestack. 
\end{remark}

\begin{remark} \label{rem-KL=heckemod-repG(O)}
Using the base-change isomorphisms, the totalizations (\ref{eqn-KL-as-limit-!}) satisfies the Beck-Chevalley condition in \cite[Theorem 4.7.5.2]{HA}, hence we have
$$  \mathbf{\Gamma}(f_\dR, \mathcal{KL}(H))  \simeq  \mathbf{\Gamma}(f_\dR, \widehat{\mathcal{H}ecke}_H)\on{-mod}[ \mathbf{\Gamma}(f_\dR,\mathcal{R}ep( \mathcal{L}^+ H ) )   ],$$
where $\widehat{\mathcal{H}ecke}_H$ is a lax unital factorization monad on $ \mathcal{R}ep( \mathcal{L}^+ H )\simeq \mathcal{I}ndCoh^{!,\ren}(\widehat{\on{Hecke}}_H^0)$ whose underlying endo-functor is given by
$$ \mathcal{I}ndCoh^{!,\ren}(\widehat{\on{Hecke}}_H^0) \xrightarrow{ p_1^{!,\ren} }  \mathcal{I}ndCoh^{!,\ren}( \widehat{\on{Hecke}}_H  ) \xrightarrow{ (p_0)^\ren_! } \mathcal{I}ndCoh^{!,\ren}(\widehat{\on{Hecke}}_H^0) ,$$
where $p_0,p_1$ are the connecting morphisms in the cosimplicial diagram $\widehat{\on{Hecke}}_H^\bullet$. In other words, we have:
\begin{equation} \label{eqn-KL-as-monad-module}
  \mathcal{KL}(H) \simeq \widehat{\mathcal{H}ecke}_H\mod( \mathcal{R}ep( \mathcal{L}^+ H )  ).
 \end{equation}

The monad structure on $\widehat{\mathcal{H}ecke}_H$ can be described as follows. The unital factorization category $\mathcal{I}ndCoh^{!,\ren}( \widehat{\on{Hecke}}_H  )$ can be canonically upgraded to an object
$$ [\mathcal{I}ndCoh^{!,\ren}( \widehat{\on{Hecke}}_H  ),\star] \in \on{Alg}( \on{FactCat}_{\on{un}} ),$$
where the convolution operator $\star$ is the following \emph{strictly}\footnote{Note that $(p_{01},p_{12})$ is strictly \emph{co}-unital.} unital factorization functor:
\[
\begin{aligned}
 \mathcal{I}ndCoh^{!,\ren}( \widehat{\on{Hecke}}_H  )\otimes \mathcal{I}ndCoh^{!,\ren}( \widehat{\on{Hecke}}_H  ) \simeq \mathcal{I}ndCoh^{!,\ren}( \widehat{\on{Hecke}}_H \times \widehat{\on{Hecke}}_H  ) \to \\ \xrightarrow{(p_{01},p_{12})^{!,\ren}}  \mathcal{I}ndCoh^{!,\ren}( \widehat{\on{Hecke}}_H^2   )  \xrightarrow{(p_{02})^\ren_!}  \mathcal{I}ndCoh^{!,\ren}( \widehat{\on{Hecke}}_H  ).
\end{aligned}
\]
We have a canonical left action of $[\mathcal{I}ndCoh^{!,\ren}( \widehat{\on{Hecke}}_H  ),\star]$ on $\mathcal{I}ndCoh^{!,\ren}(\widehat{\on{Hecke}}_H^0) $ given by the following strictly unital factorization functor:
\begin{equation} \label{eqn-geometric-left-action-Hecke-on-Rep}
\begin{aligned}
 \mathcal{I}ndCoh^{!,\ren}( \widehat{\on{Hecke}}_H  )\otimes \mathcal{I}ndCoh^{!,\ren}( \widehat{\on{Hecke}}_H^0  ) \simeq \mathcal{I}ndCoh^{!,\ren}( \widehat{\on{Hecke}}_H \times \widehat{\on{Hecke}}_H^0  ) \to \\ \xrightarrow{(\on{Id},p_{1})^{!,\ren}}  \mathcal{I}ndCoh^{!,\ren}( \widehat{\on{Hecke}}_H   ) \xrightarrow{(p_0)^\ren_!}  \mathcal{I}ndCoh^{!,\ren}( \widehat{\on{Hecke}}_H^0  ).
 \end{aligned}
\end{equation}

The unital factorization algebra\footnote{$\omega_{\widehat{\on{Hecke}}_H}$ is unital because the $!$-pullback functor $\mathcal{I}ndCoh(\on{pt}) \to \mathcal{I}ndCoh^{!,\ren}( \widehat{\on{Hecke}}_H  )$ is a lax unital factorization functor.} $\omega_{\widehat{\on{Hecke}}_H}\in \mathcal{I}ndCoh^{!,\ren}( \widehat{\on{Hecke}}_H  )$ can be canonically upgraded to an object in
$$  \omega_{\widehat{\on{Hecke}}_H} \in \on{Alg}([\mathcal{I}ndCoh^{!,\ren}( \widehat{\on{Hecke}}_H  ),\star]) .$$
Explicitly, the multplication morphism is given by
$$ \omega_{\widehat{\on{Hecke}}_H}\star \omega_{\widehat{\on{Hecke}}_H} \simeq (p_0)^\ren_!  \omega_{ \widehat{\on{Hecke}}_H^2  } \simeq ((p_0)^\ren_! \circ p_{0}^{!,\ren} ( \omega_{\widehat{\on{Hecke}}_H} ) \to \omega_{\widehat{\on{Hecke}}_H}. $$
Then the lax unital factorization monad $\widehat{\mathcal{H}ecke}_H$ on $ \mathcal{R}ep( \mathcal{L}^+ H )$ is given by the action of $\omega_{\widehat{\on{Hecke}}_H}$. In other words, we obtain a unital\footnote{See Remark \ref{rem-KL=heckemod-repG(O)-unital} below for the explanation of the unital structure.} factorization equivalence
\begin{equation} \label{eqn-KL-as-hecke-mod}
 \mathcal{KL}(H) \simeq   \omega_{\widehat{\on{Hecke}}_H}\on{-mod}( \mathcal{R}ep( \mathcal{L}^+ H )).
\end{equation}
\end{remark}

\begin{warning} $\omega_{\widehat{\on{Hecke}}_H}\on{-mod}( \mathcal{R}ep( \mathcal{L}^+ H ))$ is \emph{not} the category of factorization modules for $\omega_{\widehat{\on{Hecke}}_H}$.
\end{warning}

\begin{remark} The equivalence (\ref{eqn-KL-as-hecke-mod}) should be view as an analogue of the equivalence
$$ \mathfrak{g}\mod^K \simeq  \omega_{  K\backslash G^{\wedge}_K /K  }\on{-mod}( \on{Rep}(K) ),$$
where $(\IndCoh(K\backslash G^{\wedge}_K /K ),\star)$ acts on $\on{Rep}(K)  \simeq \IndCoh(\on{pt}/K)$.
\end{remark}

\begin{remark} By (\ref{eqn-KL-as-hecke-mod}) (and the canonical duality between $\IndCoh^{!,\ren}$ and $\IndCoh_{*,\ren}$), the fiber $[\mathcal{KL}(H)]_x $ of the factorization category $\mathcal{KL}(H)$ at a closed point $x\in X$ can be identified with the weak invariance category $\Vect^{  \widehat{\mathcal{L}H}_x,w}$ defined in \cite[$\S$ 7.15]{raskin2020homological}, which is then identified with $\mathcal{L}\mathfrak{h}_x\on{-mod}^{\mathcal{L}^+ H_x}$ in \cite[$\S$ 9.12]{raskin2020homological}. Hence (when $H$ is reductive) $\mathcal{KL}(H)$ is indeed a factorization version of the non-twisted Kazhdan--Lusztig category.

 By duality, $[\mathcal{KL}(H)_{\on{co}}]_x $ can be identified with the weak coinvariance category $\Vect_{  \widehat{\mathcal{L}H}_x,w}$.
\end{remark}

\begin{remark} \label{rem-KL*=heckemod-repG(O)}
Let us forget the fact that $\mathcal{R}ep( \mathcal{L}^+ H )$ is self dual. Then (\ref{eqn-KL-as-monad-module}) and (\ref{eqn-KL-as-hecke-mod}) implies
\begin{equation} \label{eqn-KL-*-as-hecke-mod}
 \mathcal{KL}(H)_{\on{co}} \simeq   [\widehat{\mathcal{H}ecke}_H]^\vee \mod( \mathcal{R}ep( \mathcal{L}^+ H )^\vee  ) \simeq \omega_{\widehat{\on{Hecke}}_H}\on{-mod}^r( \mathcal{R}ep( \mathcal{L}^+ H )^\vee),
\end{equation}
where $\mathcal{R}ep( \mathcal{L}^+ H )^\vee$ is equipped with the canonical \emph{right} $[\mathcal{I}ndCoh^{!,\ren}( \widehat{\on{Hecke}}_H  ),\star]$-action induced from the left $[\mathcal{I}ndCoh^{!,\ren}( \widehat{\on{Hecke}}_H  ),\star]$-action on $\mathcal{R}ep( \mathcal{L}^+ H )$. 

Via the canonical equivalence 
$$ \mathcal{R}ep( \mathcal{L}^+ H )^\vee \simeq \mathcal{R}ep( \mathcal{L}^+ H )\simeq \mathcal{I}nd\mathcal{C}oh_\ren^*(  \widehat{\on{Hecke}}_H^0 ) ,$$
the underlying endo-functor of the monad $ [\widehat{\mathcal{H}ecke}_H]^\vee$ is given by
$$ \mathcal{I}ndCoh_{*,\ren}(\widehat{\on{Hecke}}_H^0) \xrightarrow{  (p_0)_\ren^?  }  \mathcal{I}ndCoh_{*,\ren}( \widehat{\on{Hecke}}_H  ) \xrightarrow{ (p_1)_{*,\ren} } \mathcal{I}ndCoh_{*,\ren}(\widehat{\on{Hecke}}_H^0) ,$$
this can also be obtained from the fact that (the rotation of ) the totalization (\ref{eqn-KL-as-limit-*}) satisfies the Beck-Chevalley condition.
\end{remark}

\begin{warning} Consider the monoidal equivalence
$$ \sigma: [\mathcal{I}ndCoh^{!,\ren}( \widehat{\on{Hecke}}_H  ),\star] \simeq [\mathcal{I}ndCoh^{!,\ren}( \widehat{\on{Hecke}}_H  ),\star^{\on{rev}}]$$
induced by the ``taking inverse'' map
$$ \widehat{\on{Hecke}}_H \to \widehat{\on{Hecke}}_H,\; h\mapsto h^{-1}. $$
It induces a \emph{right} action of $[\mathcal{I}ndCoh^{!,\ren}( \widehat{\on{Hecke}}_H  ),\star]$ on $\mathcal{R}ep( \mathcal{L}^+ H )$. Explicitly, this right action is given by the following strictly unital factorization functor:
\begin{equation} \label{eqn-geometric-right-action-Hecke-on-Rep}
\begin{aligned}
 \mathcal{I}ndCoh^{!,\ren}( \widehat{\on{Hecke}}_H^0  )\otimes \mathcal{I}ndCoh^{!,\ren}( \widehat{\on{Hecke}}_H  ) \simeq \mathcal{I}ndCoh^{!,\ren}( \widehat{\on{Hecke}}_H^0 \times \widehat{\on{Hecke}}_H ) \to \\ \xrightarrow{(p_{0},\on{Id})^{!,\ren}}  \mathcal{I}ndCoh^{!,\ren}( \widehat{\on{Hecke}}_H   )  \xrightarrow{(p_1)^\ren_!}  \mathcal{I}ndCoh^{!,\ren}( \widehat{\on{Hecke}}_H^0  ).
\end{aligned}
\end{equation}

However, the canoncal equivalence
$$ \mathcal{R}ep( \mathcal{L}^+ H )\simeq \mathcal{R}ep( \mathcal{L}^+ H )^\vee $$
in Remark \ref{rem-Rep^!=Rep^*} is \emph{not} compatible with the right $[\mathcal{I}ndCoh^{!,\ren}( \widehat{\on{Hecke}}_H  ),\star]$-module structures on both sides. More details will be provided in the proof of Theorem \ref{thm-tate-twist}.

To avoid ambiguity, we refer the above right $[\mathcal{I}ndCoh^{!,\ren}( \widehat{\on{Hecke}}_H  ),\star]$-module on $\mathcal{R}ep( \mathcal{L}^+ H )$ as the \emph{geometric} one.
\end{warning}

\begin{remark} \label{rem-KL=heckemod-repG(O)-unital}
Via the equivalence (\ref{eqn-KL-as-hecke-mod}), the unital factorization structure on the LHS corresponds to the following canonical unital structure on the RHS. Let $\mathcal{C}$ be any unital factorization functor and $L:\mathcal{C}\to \mathcal{C}$ be a lax unital factorization monad such that $L\on{-mod}(\mathcal{C})$ is a factorization category\footnote{This is true if and only if for any $S_\dR\to T_\dR \to \Ran_\dR$, the canonical functor
$$\mathbf{\Gamma}(T_\dR,L)\on{-mod}(\mathbf{\Gamma}(T_\dR, \mathcal{C})) \otimes_{\DMod(T)} \DMod(S) \to  \mathbf{\Gamma}(T_\dR,L)\on{-mod}[\mathbf{\Gamma}(T_\dR, \mathcal{C}) \otimes_{\DMod(T)} \DMod(S) ]$$
is an equivalence.}. Then $L\on{-mod}(\mathcal{C})$ has a canonical unital factorization structure such that $\on{ind}:\mathcal{C} \to L\on{-mod}(\mathcal{C})$ is strictly unital. To see this, let $T$ be an affine test point and $\mathfrak{t}: f\to g$ be a morphism in $\Ran_{\on{un}}(T)$. The desired unital functor
$$   \mathbf{\Gamma}( f_\dR, L  )\on{-mod}( \mathbf{\Gamma}( f_\dR, \mathcal{C}  )  ) \to  \mathbf{\Gamma}( g_\dR, L  )\on{-mod}( \mathbf{\Gamma}( g_\dR, \mathcal{C}  )  ) $$
is induced by the functor
$$\mathbf{\Gamma}( f_\dR, \mathcal{C}  ) \xrightarrow{\mathfrak{t}_!} \mathbf{\Gamma}( g_\dR, \mathcal{C}  ) \xrightarrow{\on{ind}_{ g_\dR }}  \mathbf{\Gamma}( g_\dR, L  )\on{-mod}( \mathbf{\Gamma}( g_\dR, \mathcal{C}  )  ) $$
and the following canonical action of $ \mathbf{\Gamma}( f_\dR, L  )$ on this functor:
$$   \on{ind}_{ g_\dR }\circ  \mathfrak{t}_!    \circ   \mathbf{\Gamma}( f_\dR, L  ) \to \on{ind}_{ g_\dR }\circ  \mathbf{\Gamma}( g_\dR, L  ) \circ  \mathfrak{t}_!  \to \on{ind}_{ g_\dR }\circ  \mathfrak{t}_!, $$
where the first natural transformation is provided by the lax unital structure of $L$, and the second natural transformation is the action of $ \mathbf{\Gamma}( g_\dR, L  )$ on $\on{ind}_{ g_\dR }$.
\end{remark}

\subsection{The Twisted Version}
\label{sect:km-fact-twisted}
\begin{definition} A \emph{multiplicative (unital) factorization line bundle} on $\widehat{\mathcal{L}H}_{\Ran_{\on{un,dR}}}$ is a co-monoidal moprhism 
$$ \mathscr{L}: \mathbf{1}  \to [\mathcal{I}nd\mathcal{C}oh^{!}( \widehat{\mathcal{L}H}), \on{mult}^{!} ]  $$
between co-algebra objects in $[\on{FactCat}_{\on{un}}]_{\on{lax-unital}}$ such that the corresponding unital factorization algebra in $\mathcal{I}nd\mathcal{C}oh^{!}( \widehat{\mathcal{L}H})$ is invertible for the symmetric monoidal structure $[\mathcal{I}nd\mathcal{C}oh^{!}( \widehat{\mathcal{L}H}) ,\overset{!}\otimes]$.

We denote this unital factorization algebra in $\mathcal{I}nd\mathcal{C}oh^{!}( \widehat{\mathcal{L}H})$ by the same symbol $\mathscr{L}$, and treat the other data in the above definition as implicit.
\end{definition}

\begin{example} The $!$-pullback functor along $\widehat{\mathcal{L}H}_{\Ran_{\on{un,dR}}}\to \on{pt}_{\Ran_{\on{un,dR}}}$ provides a multiplicative factorization line bundle, whose underlying unital factorization algbera is $\omega_{ \widehat{\mathcal{L}H}}$.
\end{example}

\begin{remark} By definition, for a multiplicative factorization line bundle $\mathscr{L}$ on $\widehat{\mathcal{L}H}_{\Ran_{\on{un,dR}}}$, we have 
$$  \on{mult}^{!}( \mathscr{L}) \simeq \mathscr{L}\boxtimes \mathscr{L}, \; e^{!}( \mathscr{L}) \simeq \omega_{\on{pt}} $$
as unital factorization algebras in $\mathcal{I}nd\mathcal{C}oh^{!}( \widehat{\mathcal{L}H} \times  \widehat{\mathcal{L}H})$ (resp. in $\mathcal{I}nd\mathcal{C}oh^{!}(\on{pt})$).
\end{remark}

\begin{remark} \label{rem-mul-fact-line-bundle-indcoh-vs-qcoh}
Let $\mathscr{L}$ be a multiplicative factorization line bundle on $\widehat{\mathcal{L}H}_{\Ran_{\on{un,dR}}}$. For any $s: S\to \Ran_{\on{un,dR}}$ with $S\in \,^{>-\infty}\on{Sch}_{ \on{aft}}$, we obtain an invertible object
$$ \mathscr{L}_s \in [ \IndCoh^{!}(     \widehat{\mathcal{L}H}_s ) ,\overset{!}\otimes ]. $$
By the proof of \cite[Proposition 9.4.2]{raskin2020homological}, there exists a unique invertible object 
$$ \mathscr{L}_s^{\on{QCoh}} \in [ \QCoh(     \widehat{\mathcal{L}H}_s ) ,\otimes ] $$
such that $\Upsilon( \mathscr{L}_s^{\on{QCoh}}) \simeq \mathscr{L}_s$. Such an object is classified by a morphism $\widehat{\mathcal{L}H}_s \to \mathbb{BG}_m \times \mathbb{Z} $. The multiplicative structure on $\mathscr{L}$ can upgrade this morphism to a group homomorphism (lying over the projection $S\to \on{pt}$). Since the fibers of $\widehat{\mathcal{L}H}_s\to S$ are connected, this group homomorphism factors through $\widehat{\mathcal{L}H}_s \to \mathbb{BG}_m$. It follows that $\mathscr{L}_s^{\on{QCoh}}$ is indeed a line bundle in $\QCoh(     \widehat{\mathcal{L}H}_s )^\heartsuit$ and our definition of multiplicative factorization line bundle (via $\IndCoh^!$-theory) is essentially equivalent to that in \cite{zhao2017quantum} defined via $\on{QCoh}$-theory.
\end{remark}

\begin{remark} \label{rem-fact-line-bunlde-being-mult-property}
By Remark \ref{rem-mul-fact-line-bundle-indcoh-vs-qcoh}, for a unital factorization line bundle $\mathscr{L}$ on $\widehat{\mathcal{L}H}_{\Ran_{\on{un,dR}}}$ and fixed equivalences $\on{mult}^{!}( \mathscr{L}) \simeq \mathscr{L}\boxtimes \mathscr{L}, \; e^{!}( \mathscr{L}) \simeq \omega_{\on{pt}}$, there is at most one multiplicative structure on $\mathscr{L}$ whose structure equivalences are given by these fixed equivalences. In other words, the higher compatibilities in the definition of multiplicative factorization line bundles are \emph{properties rather than additional data}.
\end{remark}

\begin{construction} For multiplicative factorization line bundles $\mathscr{L}, \mathscr{L}'$ on $\widehat{\mathcal{L}H}_{\Ran_{\on{un,dR}}}$, we define multiplicative factorization line bundles $\mathscr{L}\overset{!}\otimes \mathscr{L}'$ and $\mathscr{L}^{-1}$ in the obvious way.
\end{construction}

\begin{definition} Let $\mathscr{L}$ be a multiplicative factorization line bundle on $\widehat{\mathcal{L}H}_{\Ran_{\on{un,dR}}}$. Let
$$  \mathscr{L}|_{ \mathcal{L}^+ H } \in  \mathcal{I}nd\mathcal{C}oh^{!}( \widehat{\mathcal{L}^+H}) $$
be the unital factorization algebra obtained as the image of $\mathscr{L}$ under the lax-unital factorization functor
$$!\on{-pull}:   \mathcal{I}nd\mathcal{C}oh^{!}( \widehat{\mathcal{L}H}) \to \mathcal{I}nd\mathcal{C}oh^{!}( \widehat{\mathcal{L}^+H}).$$
We say \emph{$\mathscr{L}$ is trivial on $\mathcal{L}^+ H$}, or \emph{$\mathscr{L}|_{ \mathcal{L}^+ H }$ is trivial} if the unital morphism $\omega_{ \mathcal{L}^+ H  } \to \mathscr{L}|_{ \mathcal{L}^+ H }$ is an isomorphism.
\end{definition}

\begin{construction}
\label{constr:hecke-line-bundle-is-algebra}
Let $\mathscr{L}$ be a multiplicative factorization line bundle on $\widehat{\mathcal{L}H}_{\Ran_{\on{un,dR}}}$ that is trivial on $\mathcal{L}^+ H$. Then the restriction of the group homomorphism $\widehat{\mathcal{L}H}_s \to \mathbb{BG}_m$ in Remark \ref{rem-mul-fact-line-bundle-indcoh-vs-qcoh} on $\mathcal{L}^+H_s$ is trivialized. Hence $\widehat{\mathcal{L}H}_s \to \mathbb{BG}_m$ factors through a map $[\widehat{\on{Hecke}}_H]_s \to \mathbb{BG}_m$. Consider the universal line bundle $\mathcal{L}_{\on{univ}} \in \IndCoh(\mathbb{BG}_m )$ and its image $[\mathscr{L}_{\widehat{\on{Hecke}}_H}]_s \in \IndCoh^{!,\ren}(  [\widehat{\on{Hecke}}_H]_s )$ under the $!$-pullback functor. By construction, the image of the $!$-pullback of $ [\mathscr{L}_{\widehat{\on{Hecke}}_H}]_s $ along $\widehat{\mathcal{L}H}_s\to  [\widehat{\on{Hecke}}_H]_s$ can be identified with $\mathscr{L}_s$. The above construction is functorial in $S$ and compatible with the factorization structure. Hence we obtain a unital factorization algebra $\mathscr{L}_{\widehat{\on{Hecke}}_H}$ in $\mathcal{I}nd\mathcal{C}oh^{!,\ren}(  \widehat{\on{Hecke}}_H    )$. We also write $\mathscr{L}_{\widehat{\on{Hecke}}}$ for $\mathscr{L}_{\widehat{\on{Hecke}}_H}$ when there is no danger of ambiguity.

By construction, for any $[m]$, we have an isomorphism\footnote{When $m=0$, this isomorphism says 
$$ !\on{-pull}: \mathcal{I}nd\mathcal{C}oh^{!,\ren}(  \widehat{\on{Hecke}}_H    ) \to \mathcal{I}nd\mathcal{C}oh^{!,\ren}( \mathbb{B} \mathcal{L}^+H    ) $$
sends $\mathscr{L}_{\widehat{\on{Hecke}}}$ to $\omega$.}
$$ 
 (p_{0m})^{!,\ren}\mathscr{L}_{\widehat{\on{Hecke}}} \simeq  (p_{01})^{!,\ren}\mathscr{L}_{\widehat{\on{Hecke}}}\overset{!}\otimes \cdots \overset{!}\otimes (p_{(m-1)m})^{!,\ren}\mathscr{L}_{\widehat{\on{Hecke}}}$$
 between unital factorization algebras in $\mathcal{I}nd\mathcal{C}oh^{!,\ren}(  \widehat{\on{Hecke}}_H^\bullet   )$. We also have higher compatibilities for these isomorphisms. It follows that we can upgrade $\mathscr{L}_{\widehat{\on{Hecke}}}$ to an algebra
 $$ \mathscr{L}_{\widehat{\on{Hecke}}} \in \on{Alg}(  [  \mathcal{I}nd\mathcal{C}oh^{!,\ren}(  \widehat{\on{Hecke}}_H    ) ,\star ]  ).$$
 Explicitly, the multiplication morphism is given by
 $$  \mathscr{L}_{\widehat{\on{Hecke}}}\star \mathscr{L}_{\widehat{\on{Hecke}}} \simeq (p_{0m})^{\ren}_![ (p_{01})^{!,\ren}\mathscr{L}_{\widehat{\on{Hecke}}}\overset{!}\otimes  (p_{12})^{!,\ren}\mathscr{L}_{\widehat{\on{Hecke}}}  ] \simeq (p_{0m})^{\ren}_!\circ (p_{0m})^{!,\ren}(\mathscr{L}_{\widehat{\on{Hecke}}}) \to  \mathscr{L}_{\widehat{\on{Hecke}}}. $$
\end{construction}

\begin{example}
\label{example:tate-bundle}
 We have the following multiplicative factorization line bundle $\mathscr{L}^{\on{Tate}}$ on $\widehat{\mathcal{L}H}_{\Ran_{\on{un,dR}}}$ that is trivial on $\mathcal{L}^+ H$, which is known as the \emph{Tate line bundle}. For any injective $\phi:I\to J$ and $S=\on{Spec}(R) \to [\widehat{\on{Hecke}}_H]_{\phi}$ with $R$ being a \emph{classical} commutative algebra, the restriction of $\mathscr{L}^{\on{Tate}}_{\widehat{\on{Hecke}}}$ on $\on{Spec}(R)$ is given as follows. Let $\mathcal{P}_0$ and $\mathcal{P}_1$ be as in Remark \ref{rem-description-unital-structure-hecke}. Consider the vector bundles $\mathfrak{h}_{\mathcal{P}_0},\mathfrak{h}_{\mathcal{P}_1}$ on $\mathcal{D}_{\Gamma_J}'$ over $\mathcal{D}_{\Gamma_J}'$ associated to the adjoint representation $\mathfrak{h}\in \on{Rep}(H)$. Via the projection $\mathcal{D}_{\Gamma_J}' \to S$, we view them as Tate vector bundles on $S$, and denote them by $(\mathcal{L}^+\mathfrak{h})_{\mathcal{P}_0}, (\mathcal{L}^+\mathfrak{h})_{\mathcal{P}_1} $. Similarly, the restriction of $\mathfrak{h}_{\mathcal{P}_0},\mathfrak{h}_{\mathcal{P}_1}$ on $\mathcal{D}_{\Gamma_J}'\setminus \Gamma_I $ provides Tate vector bundles $(\mathcal{L}\mathfrak{h})_{\mathcal{P}_0}, (\mathcal{L}\mathfrak{h})_{\mathcal{P}_1} $ on $S$ which contain  $(\mathcal{L}^+\mathfrak{h})_{\mathcal{P}_0}, (\mathcal{L}^+\mathfrak{h})_{\mathcal{P}_1} $ as lattice sub-bundles. By definition, we have an identification $(\mathcal{L}\mathfrak{h})_{\mathcal{P}_0}\simeq (\mathcal{L}\mathfrak{h})_{\mathcal{P}_1}$ which allows us to define the relative determinant line bundle $\on{rel.det}( (\mathcal{L}^+\mathfrak{h})_{\mathcal{P}_0},(\mathcal{L}^+\mathfrak{h})_{\mathcal{P}_1}  ) \in \on{QCoh}(S)$. Then the restriction of $\mathscr{L}^{\on{Tate}}_{\widehat{\on{Hecke}}}$ on $S$ is
    $$\Upsilon(  \on{rel.det}( (\mathcal{L}^+\mathfrak{h})_{\mathcal{P}_0},(\mathcal{L}^+\mathfrak{h})_{\mathcal{P}_1}  )   ) \in \on{IndCoh}^!(S).$$
\end{example}

\begin{definition}
\label{def:kl-factorizable}
Let $\mathscr{L}$ be a multiplicative factorization line bundle on $\widehat{\mathcal{L}H}_{\Ran_{\on{un,dR}}}$ that is trivial on $\mathcal{L}^+ H$. We define 
\begin{eqnarray*} 
 \mathcal{KL}(H)_{\mathscr{L}} & :=  \widetilde{\mathcal{L}\mathfrak{h}}\on{-}mod_{\mathscr{L}}^{\mathcal{L}^+H}  &:=   \mathscr{L}_{\widehat{\on{Hecke}}} \on{-mod}( \mathcal{R}ep( \mathcal{L}^+ H )) \\
 \mathcal{KL}(H)_{\mathscr{L}^{-1},\on{co}} & := \widetilde{\mathcal{L}\mathfrak{h}}\on{-}mod_{\mathscr{L}^{-1},\on{co}}^{\mathcal{L}^+H}&:=    \mathscr{L}_{\widehat{\on{Hecke}}}\on{-mod}^r( \mathcal{R}ep( \mathcal{L}^+ H )^\vee),
\end{eqnarray*}
where $\mathcal{R}ep( \mathcal{L}^+ H )$ (resp, $\mathcal{R}ep( \mathcal{L}^+ H )^\vee$) is equipped with the left (resp. right) $ [\mathcal{I}ndCoh^{!,\ren}( \widehat{\on{Hecke}}_H  ),\star]$-module structure described in Remark \ref{rem-KL=heckemod-repG(O)} and Remark \ref{rem-KL*=heckemod-repG(O)}.
\end{definition}

Several explanations are in order.

\begin{remark} \label{rem-KL-twisted-as-totalization}
$\mathcal{KL}(H)_{\mathscr{L}} $ can also be defined as a geometric realization or totalization as follows. There is a functor
$$\mathcal{I}nd\mathcal{C}oh^{!,\ren}( \widehat{\on{Hecke}}_H^\bullet )_{\mathscr{L},\Delta}: \Delta \to [\on{FactCat}_{\on{un}}]_{\on{lax-unital}} $$
that sends $[m]$ to $\mathcal{I}nd\mathcal{C}oh^{!,\ren}( \widehat{\on{Hecke}}_H^m )$ but sends a morphism $f:[m]\to [n]$ in $\Delta$ to the functor 
$$  (p_f)^{!,\ren}(-)\overset{!}\otimes (p_{ f^{\gets}})^{!,\ren} ( \mathscr{L}_{ \widehat{\on{Hecke}} } ):  \mathcal{I}nd\mathcal{C}oh^{!,\ren}( \widehat{\on{Hecke}}_H^m ) \to \mathcal{I}nd\mathcal{C}oh^{!,\ren}( \widehat{\on{Hecke}}_H^n ), $$
where $f^{\gets}: [1]\to [m]$ is the map $f^{\gets}(0)=0, f^{\gets}(1)=f(0)$. The higher compatibilities are provided by the multiplicative structure\footnote{For instance, the composition laws are provided by 
$$ p_{ (g\circ f)^\gets }^{!,\ren}(  \mathscr{L}_{ \widehat{\on{Hecke}} } ) \simeq p_{ g^\gets }^{!,\ren}(  \mathscr{L}_{ \widehat{\on{Hecke}} }) \overset{!}\otimes   p_{ (g\circ  f^\gets) }^{!,\ren}(  \mathscr{L}_{ \widehat{\on{Hecke}} } ) . $$} of $\mathscr{L}$.
By passing to left adjoints, we also have a functor
$$\mathcal{I}nd\mathcal{C}oh^{!,\ren}( \widehat{\on{Hecke}}_H^\bullet )_{\mathscr{L},\Delta^\op} : \Delta^\op \to \on{FactCat}_{\on{un}} $$
whose connecting functors are
$$   (p_f)^{\ren}_! [-  \overset{!}\otimes (p_{ f^{\gets}})^{!,\ren} ( \mathscr{L}_{ \widehat{\on{Hecke}} }^{-1}) ]:  \mathcal{I}nd\mathcal{C}oh^{!,\ren}( \widehat{\on{Hecke}}_H^n ) \to \mathcal{I}nd\mathcal{C}oh^{!,\ren}( \widehat{\on{Hecke}}_H^m ). $$

As in Remark \ref{rem-KL=heckemod-repG(O)}, the totalization of the sections of $\mathcal{I}nd\mathcal{C}oh^{!,\ren}( \widehat{\on{Hecke}}_H^\bullet )_{\mathscr{L}}$ satisfies the Bech-Chevalley condition, and the corresponding monad can be identified with the action of $ \mathscr{L}_{ \widehat{\on{Hecke}} }$ on $\mathcal{R}ep( \mathcal{L}^+ H )$. Therefore
$$ \mathcal{KL}(H)_{\mathscr{L}}  \simeq \lim \mathcal{I}nd\mathcal{C}oh^{!,\ren}( \widehat{\on{Hecke}}_H^\bullet )_{\mathscr{L},\Delta} \simeq \colim   \mathcal{I}nd\mathcal{C}oh^{!,\ren}( \widehat{\on{Hecke}}_H^\bullet )_{\mathscr{L},\Delta^\op}.$$

In particular, the sections of $\mathcal{KL}(H)_{\mathscr{L}}$ are compactly generated.
\end{remark}

\begin{remark} \label{rem-KL-twisted-as-totalization-*}
Passing to duality, the simplicial diagrams in Remark \ref{rem-KL-twisted-as-totalization} provides
$$\mathcal{I}nd\mathcal{C}oh_{*,\ren}( \widehat{\on{Hecke}}_H^\bullet )_{ \mathscr{L}^{-1},\Delta }: \Delta \to [\on{FactCat}_{\on{un}}]_{\on{lax-unital}} $$
whose connecting functors are
$$ (p_{ f^{\gets}})^{!,\ren} (  \mathscr{L}_{ \widehat{\on{Hecke}} }^{-1} )  \overset{\on{act}}\otimes (p_f)_{\ren}^? ( - )   :  \mathcal{I}nd\mathcal{C}oh_{*,\ren}( \widehat{\on{Hecke}}_H^m ) \to \mathcal{I}nd\mathcal{C}oh_{*,\ren}( \widehat{\on{Hecke}}_H^n ),$$
where $\overset{\on{act}}\otimes$ means the $[\mathcal{I}nd\mathcal{C}oh^{!,\ren}(-), \overset{!}\otimes]$-action on $\mathcal{I}nd\mathcal{C}oh_{*,\ren}(-)$. By passing to left adjoints, we also have a functor
$$\mathcal{I}nd\mathcal{C}oh_{*,\ren}( \widehat{\on{Hecke}}_H^\bullet )_{\mathscr{L}^{-1},\Delta^\op}: \Delta^\op \to \on{FactCat}_{\on{un}} $$
whose connecting functors are
$$(p_f)_{*,\ren} [(p_{ f^{\gets}})^{!,\ren} (  \mathscr{L}_{ \widehat{\on{Hecke}} } )  \overset{\on{act}}\otimes - ]   :  \mathcal{I}nd\mathcal{C}oh_{*,\ren}( \widehat{\on{Hecke}}_H^n ) \to \mathcal{I}nd\mathcal{C}oh_{*,\ren}( \widehat{\on{Hecke}}_H^m ).$$
As before, we have
$$\mathcal{KL}(H)_{\mathscr{L}^{-1},\on{co}}  \simeq \lim \mathcal{I}nd\mathcal{C}oh_{*,\ren}( \widehat{\on{Hecke}}_H^\bullet )_{\mathscr{L}^{-1},\Delta} \simeq \colim   \mathcal{I}nd\mathcal{C}oh_{*,\ren}( \widehat{\on{Hecke}}_H^\bullet )_{\mathscr{L}^{-1},\Delta^\op}.$$

\end{remark}

\begin{remark} \label{rem-dual-!-version-*-version-KL-twisted}
 Similar to the non-twisted version, $\mathcal{KL}(H)_{\mathscr{L}}$ and $\mathcal{KL}(H)_{ \mathscr{L}^{-1},\on{co} }$ are dual to each other in the symmetric monoidal category $[\on{FactCat}_{\on{un}}]_{\on{lax-unital}}$.
\end{remark}

\begin{remark} One should think $\mathscr{L}$ as a \emph{factorization gerbe} on $[\mathbb{B}\widehat{\mathcal{L}H}]_{\Ran_{\on{un,dR}}}$ equipped with a trivialization on $[\mathbb{B}\mathcal{L}^+H]_{\Ran_{\on{un,dR}}}$, and $\mathcal{KL}(H)_\mathscr{L}$ as the unital factorization category for the $\IndCoh^{!,\ren}$-theory on $[\mathbb{B}\widehat{\mathcal{L}H}]_{\Ran_{\on{un,dR}}}$ twisted by this gerbe.
\end{remark}

\begin{remark}\label{rem-gerbe-vs-central-extension}
 Let $x\in X$ be a closed point. The group homomorphism $\widehat{\mathcal{L} H}_x \to \mathbb{BG}_m$ in Remark \ref{rem-mul-fact-line-bundle-indcoh-vs-qcoh} induces a central extension of $\widehat{\mathcal{L} H}_x$ that splits over ${\mathcal{L}^+ H}_x$. By \cite[$\S$ 9,4]{raskin2020homological}, this also corresponds to a weak $\widehat{\mathcal{L} H}_x$-action on $\on{Vect}$ whose restriction on ${\mathcal{L}^+ H}_x$ is trivialized. Let $\on{Vect}_{\mathscr{L}_x}$ be the correponding object in $\widehat{\mathcal{L}^+ H}_x\mod^{\on{weak}}$. It follows from Remark \ref{rem-KL-twisted-as-totalization-*} that $[\mathcal{KL}(H)_{\mathscr{L}}]_x$ can be identified with $(\on{Vect}_{\mathscr{L}_x}){^{  \widehat{\mathcal{L}H}_x,w}}$, which is then identified in \cite[Section 9.12]{raskin2020homological} as $\widetilde{\mathcal{L}\mathfrak{h}_x}\on{-mod}^{\mathcal{L}^+ H_x}_{\mathscr{L}}$, where $(\widetilde{\mathcal{L}\mathfrak{h}_x})_{\mathscr{L}}$ is the corresponding central extension of $\mathcal{L}\mathfrak{h}_x$. Hence (when $H$ is reductive) $\mathcal{KL}(H)_{\mathscr{L}}$ is indeed a factorization version of the Kazhdan--Lusztig category.
\end{remark}

\subsection{The Tate Line Bundle}

The goal of this subsection is to prove the following result:

\begin{theorem} \label{thm-tate-twist}
Let $\mathscr{L}$ be a multiplicative factorization line bundle on $\widehat{\mathcal{L}H}_{\Ran_{\on{un,dR}}}$ that is trivial on $\mathcal{L}^+ H$. Then there exists a canonical equivalence
$$ \mathcal{KL}(H)_{ \mathscr{L} }  \simeq \mathcal{KL}(H)_{\mathscr{L}\otimes \mathscr{L}^{\on{Tate}},\on{co} }$$
in $\on{FactCat}_{\on{un}}$.
\end{theorem}

\begin{remark} For any closed point $x\in X$, the combination of \cite[Proposition 7.18.2]{raskin2020homological} and \cite[Theorem 9.16.1]{raskin2020homological} provides \emph{an} equivalence $ [\mathcal{KL}(H)_{ \mathscr{L} }]_x  \simeq [\mathcal{KL}(H)_{\mathscr{L}\otimes \mathscr{L}^{\on{Tate}} ,\on{ co}}]_x$. We will explain in Remark \ref{rem-modular-character-compare-with-sam} how to identify it with the restriction of our equivalence at $x$.
\end{remark}

\proof 

\emph{Step 1}: In this step, for any $[m]$, we will construct $(m+1)$ different $[ \mathcal{I}nd\mathcal{C}oh^{!,\ren}( \widehat{\on{Hecke}}_H^m ),\overset{!}\otimes]$-linear equivalences $\theta_i: \mathcal{I}nd\mathcal{C}oh^{!,\ren}( \widehat{\on{Hecke}}_H^m ) \to \mathcal{I}nd\mathcal{C}oh_{*,\ren}( \widehat{\on{Hecke}}_H^m )$.

For each $i: [0]\to [m]$, consider the corresponding projection $p_i: \widehat{\on{Hecke}}_H^m \to \mathbb{B} \mathcal{L}^+ H $. One can check directly that the following functors provide a self-duality on $\mathcal{I}nd\mathcal{C}oh^{!,\ren}( \widehat{\on{Hecke}}_H^m )$:
\[
\begin{aligned}
\mathcal{I}nd\mathcal{C}oh^{!,\ren}( \widehat{\on{Hecke}}_H^m ) \otimes \mathcal{I}nd\mathcal{C}oh^{!,\ren}( \widehat{\on{Hecke}}_H^m ) \xrightarrow{\otimes^!} \mathcal{I}nd\mathcal{C}oh^{!,\ren}( \widehat{\on{Hecke}}_H^m ) \xrightarrow{  (p_i)^\ren_! } \mathcal{I}nd\mathcal{C}oh^{!,\ren}( \mathbb{B} \mathcal{L}^+ H ) \to \\ \xrightarrow{?-\on{push}}  \mathcal{I}nd\mathcal{C}oh^{!,\ren}(\on{pt}) \simeq \mathbf{1} \\
\mathbf{1} \simeq  \mathcal{I}nd\mathcal{C}oh^{!,\ren}(\on{pt})  \xrightarrow{!-\on{pull}}  \mathcal{I}nd\mathcal{C}oh^{!,\ren}(\widehat{\on{Hecke}}_H^m ) \xrightarrow{ !-\on{push} }  \mathcal{I}nd\mathcal{C}oh^{!,\ren}(\widehat{\on{Hecke}}_H^m \times_{ (\mathbb{B} \mathcal{L}^+ H)  } \widehat{\on{Hecke}}_H^m ) \to\\
\xrightarrow{?\on{-push}} \mathcal{I}nd\mathcal{C}oh^{!,\ren}(\widehat{\on{Hecke}}_H^m \times \widehat{\on{Hecke}}_H^m ) \simeq \mathcal{I}nd\mathcal{C}oh^{!,\ren}( \widehat{\on{Hecke}}_H^m ) \otimes \mathcal{I}nd\mathcal{C}oh^{!,\ren}( \widehat{\on{Hecke}}_H^m ).
\end{aligned}
\]
Therefore for each $i: [0]\to [m]$, we obtain an equivalence
$$ \theta_i: \mathcal{I}nd\mathcal{C}oh^{!,\ren}( \widehat{\on{Hecke}}_H^m ) \to \mathcal{I}nd\mathcal{C}oh_{*,\ren}( \widehat{\on{Hecke}}_H^m ).$$
Note that the pairing functor in the above duality factors through the $!$-tensor product functor. It follows formally that $\theta_i$ has an $[ \mathcal{I}nd\mathcal{C}oh^{!,\ren}( \widehat{\on{Hecke}}_H^m ),\overset{!}\otimes]$-linear structure. 

The above construction is functorial in $[m]$ in the sense that for $f:[m]\to [n]$, we have a canonical commutative diagram
$$
\xymatrix{
    \mathcal{I}nd\mathcal{C}oh^{!,\ren}( \widehat{\on{Hecke}}_H^m ) \ar[r]_-\simeq^-{\theta_i} \ar[d]^-{ (p_f)^{!,\ren}  }
     & \mathcal{I}nd\mathcal{C}oh_{*,\ren}( \widehat{\on{Hecke}}_H^m )
     \ar[d]^-{ (p_f)^{?}_\ren  }  \\
     \mathcal{I}nd\mathcal{C}oh^{!,\ren}( \widehat{\on{Hecke}}_H^n ) \ar[r]_-\simeq^-{\theta_{f\circ i}} & \mathcal{I}nd\mathcal{C}oh_{*,\ren}( \widehat{\on{Hecke}}_H^n ).
}
$$

\emph{Step 2}: In this step, we use the equivalences $\theta_i$ constructed in Step 1 to construct a multiplicative factorization line bundle $\mathscr{L}^{\on{diff}}$ on $\widehat{\mathcal{L}H}_{\Ran_{\on{un,dR}}}$ that is trivial on $\mathcal{L}^+ H$.

Consider the two $[\mathcal{I}nd\mathcal{C}oh^{!,\ren}( \widehat{\on{Hecke}}_H ),\overset{!}\otimes]$-linear equivalences
$$ \theta_0,\theta_1: \mathcal{I}nd\mathcal{C}oh^{!,\ren}( \widehat{\on{Hecke}}_H ) \to \mathcal{I}nd\mathcal{C}oh_{*,\ren}( \widehat{\on{Hecke}}_H )$$
constructed in Step 1. There exists a unique unital factorization line bundle $\mathscr{L}^{\on{diff}}_{ \widehat{\on{Hecke}}_H}$ on $[\widehat{\on{Hecke}}_H ]_{\Ran_{\on{un,dR}}}$ such that
$$ \theta_0(-) \simeq  \theta_1(\mathscr{L}^{\on{diff}}_{ \widehat{\on{Hecke}}_H} \overset{!}\otimes - ) \simeq \mathscr{L}^{\on{diff}}_{ \widehat{\on{Hecke}}_H} \overset{\on{act}}\otimes \theta_1(-) $$
Let $\mathscr{L}^{\on{diff}}$ be its $!$-pullback along $\widehat{\mathcal{L}H}_{\Ran_{\on{un,dR}}}\to [\widehat{\on{Hecke}}_H ]_{\Ran_{\on{un,dR}}}$. Then the functoriality of $\theta_i$ described in Step 1 can upgrade $\mathscr{L}^{\on{diff}}$ to a multiplicative factorization line bundle on $\widehat{\mathcal{L}H}_{\Ran_{\on{un,dR}}}$ that is trivial on $\mathcal{L}^+ H$.

\emph{Step 3}: It follows from construction that we have an equivalence
$$ \mathcal{I}nd\mathcal{C}oh^{!,\ren}( \widehat{\on{Hecke}}_H^\bullet )_{\mathscr{L},\Delta} \simeq \mathcal{I}nd\mathcal{C}oh_{*,\ren}( \widehat{\on{Hecke}}_H^\bullet )_{\mathscr{L} \otimes\mathscr{L}^{\on{diff}} ,\Delta}  $$
between cosimplicial objects in $[\on{FactCat}_{\on{un}}]_{\on{lax-unital}}$, whose value on any term is the equivalence 
$$\theta_0: \mathcal{I}nd\mathcal{C}oh^{!,\ren}( \widehat{\on{Hecke}}_H^m ) \simeq \mathcal{I}nd\mathcal{C}oh_{*,\ren}( \widehat{\on{Hecke}}_H^m ) .$$

Hence we obtain an equivalence 
$$\mathcal{KL}(H)_{ \mathscr{L} }  \simeq \mathcal{KL}^*(H)_{\mathscr{L}\otimes \mathscr{L}^{\on{diff}} }$$
in $[\on{FactCat}_{\on{un}}]_{\on{lax-unital}}$, which is automatically an equivalence in $\on{FactCat}_{\on{un}}$. Hence to finish the proof, we only need to show $\mathscr{L}^{\on{diff}} \simeq \mathscr{L}^{\on{Tate}}$.

\emph{Step 4}: In this step, we construct two different $[\mathcal{I}nd\mathcal{C}oh^{!,\ren}(   \widehat{\mathcal{L}H} ),\overset{!}\otimes ]$-linear equivalences 
$$\mathcal{I}nd\mathcal{C}oh^{!,\ren}(   \widehat{\mathcal{L}H} )\to \mathcal{I}nd\mathcal{C}oh_{*,\ren}(   \widehat{\mathcal{L}H} ),$$
and we show that the difference of them is given by a unital factorization line bundle $\mathscr{L}^{\on{diff}'}$.

We claim the following functor is a perfect pairing, i.e., induces a self-duality on $\mathcal{I}nd\mathcal{C}oh^{!,\ren}(   \widehat{\mathcal{L}H} )$:
\[
\begin{aligned}
\mathcal{I}nd\mathcal{C}oh^{!,\ren}(   \widehat{\mathcal{L}H} ) \otimes \mathcal{I}nd\mathcal{C}oh^{!,\ren}(   \widehat{\mathcal{L}H} ) \xrightarrow{\otimes} \mathcal{I}nd\mathcal{C}oh^{!,\ren}(   \widehat{\mathcal{L}H} ) \xrightarrow{?\on{-push}} \mathcal{I}nd\mathcal{C}oh^{!,\ren}(   \widehat{\mathcal{L}H}/ \mathcal{L}^+H ) \to \\ \xrightarrow{!\on{-push}} \mathcal{I}nd\mathcal{C}oh^{!,\ren}(   \on{pt} ) \simeq \mathbf{1}.
\end{aligned}
\]
To prove this, we only need to verify the factorization functor $\mathcal{I}nd\mathcal{C}oh^{!,\ren}(   \widehat{\mathcal{L}H} )\to \mathcal{I}nd\mathcal{C}oh_{*,\ren}(   \widehat{\mathcal{L}H} )$ induced by it is an equivalence. Consider the following lax-unital factorization functor
$$\mathbf{1} \xrightarrow{?\on{-pull}} \mathcal{I}nd\mathcal{C}oh_{*,\ren}( \widehat{\mathcal{L}H} / \mathcal{L}^+H ) \xrightarrow{*\on{-pull}} \mathcal{I}nd\mathcal{C}oh_{*,\ren}(   \widehat{\mathcal{L}H} )$$
and the corresponding unital factorization algebra $\mathcal{F}_0$. Unwinding the definitions, the aforementioned functor is given by
$$ -\overset{\on{act}}\otimes \mathcal{F}_0: \mathcal{I}nd\mathcal{C}oh^{!,\ren}(   \widehat{\mathcal{L}H} )\to \mathcal{I}nd\mathcal{C}oh_{*,\ren}(   \widehat{\mathcal{L}H} ). $$
To prove its invertibility, we only need to show for any \emph{classical} test scheme $s:S\to \Ran_{\on{un,dR}}$, the functor $-\overset{\on{act}}\otimes [\mathcal{F}_0]_s: \IndCoh^{!}(   \widehat{\mathcal{L}H}_s )\to \IndCoh_{*}(   \widehat{\mathcal{L}H}_s )$ is an equivalence. Let $[ \widehat{\mathcal{L}H}/\mathcal{L}^+H]_s \simeq \colim_\alpha Y_\alpha$ be a presentation of $[ \widehat{\mathcal{L}H}/\mathcal{L}^+H]_s$ as a classical ind-finite type ind-proper indscheme over $S$. Define $Z_\alpha:= Y_\alpha \times_{  \mathbb{B} \mathcal{L}^+H_s } S$ and $\pi_\alpha: Z_\alpha \to Y_\alpha$. Using Proposition \ref{prop-indcoh-on-ind-stack}, we only need to show the action of $\IndCoh^!(Z_\alpha)$ on $(\pi_\alpha)^*( \omega_{Y_\alpha})$ induces an equivalence $\IndCoh^!(Z_\alpha) \to \IndCoh_*(Z_\alpha)$. Write $Z_\alpha = \lim_k Y_\alpha^k$ such that $q_k:Y_\alpha^k \to Y_\alpha$ are smooth. Using Proposition \ref{prop-indcoh-on-quot-stack}, we only need to show the action of $\IndCoh(Y_\alpha^k)$ on $(q_k)^*(\omega_{Y_\alpha})$ induces an automorphism on $\IndCoh(Y_\alpha^k)$, but this follows from the fact that $q_k$ is smooth. This finishes the proof of the claim.

Let $\mathcal{F}_1$ be the unital factorization algebra that corresponds to
$$\mathbf{1} \xrightarrow{?\on{-pull}} \mathcal{I}nd\mathcal{C}oh_{*,\ren}( \mathcal{L}^+H \backslash\widehat{ \mathcal{L}H} ) \xrightarrow{*\on{-pull}} \mathcal{I}nd\mathcal{C}oh_{*,\ren}(   \widehat{\mathcal{L}H} ).$$
By symmetry, the functor
$$ -\overset{\on{act}}\otimes \mathcal{F}_1: \mathcal{I}nd\mathcal{C}oh^{!,\ren}(   \widehat{\mathcal{L}H} )\to \mathcal{I}nd\mathcal{C}oh_{*,\ren}(   \widehat{\mathcal{L}H} ) $$
is also an equivalence. It follows that there exists a unique unital factorization line bundle $\mathscr{L}^{\on{diff}'}$ on $[\widehat{\mathcal{L}H}]_{\on{Ran}_{\on{un,dR}}}$ such that $\mathcal{F}_0 \simeq \mathscr{L}^{\on{diff}'}\overset{\on{act}}\otimes \mathcal{F}_1 $.

\emph{Step 5}: In this step, we show $\mathscr{L}^{\on{diff}'}$ has a canonical multiplicative structure, and $\mathscr{L}^{\on{diff}'}\simeq \mathscr{L}^{\on{diff}}$.

It is clear that $\mathscr{L}^{\on{diff}'}\boxtimes \mathscr{L}^{\on{diff}'}$ is the unique unital factorization line bundle on $[\widehat{\mathcal{L}H}\times \widehat{\mathcal{L}H}]_{\on{Ran}_{\on{un,dR}}}$ such that $\mathcal{F}_0 \boxtimes \mathcal{F}_0 \simeq (\mathscr{L}^{\on{diff}'}\boxtimes \mathscr{L}^{\on{diff}'})\overset{\on{act}}\otimes (\mathcal{F}_1\boxtimes\mathcal{F}_1)$.

On the other hand, consider the composition 
$$\on{mult}: [\widehat{\mathcal{L}H}\times \widehat{\mathcal{L}H}]_{\on{Ran}_{\on{un,dR}}} \xrightarrow{q} [\widehat{\mathcal{L}H}\times^{\mathcal{L}^+H} \widehat{\mathcal{L}H}]_{\on{Ran}_{\on{un,dR}}} \xrightarrow{\on{m}} [\widehat{\mathcal{L}H}]_{\on{Ran}_{\on{un,dR}}}.$$
Using the base-change isomorphisms, we can obtain $q^*_\ren \circ m^{?}_\ren(  \mathcal{F}_i ) \simeq \mathcal{F}_i \boxtimes \mathcal{F}_i$. This implies $\mathcal{F}_0 \boxtimes \mathcal{F}_0 \simeq ( \on{mult}^{!,\ren}(   \mathscr{L}^{\on{diff}'}  )   ) \overset{\on{act}}\otimes (\mathcal{F}_1\boxtimes\mathcal{F}_1)$. 

Therefore we have $\on{mult}^{!,\ren}(   \mathscr{L}^{\on{diff}'})\simeq \mathscr{L}^{\on{diff}'}\boxtimes \mathscr{L}^{\on{diff}'}$. Similarly we have $e^{!,\ren}(   \mathscr{L}^{\on{diff}'})\simeq \omega_{\on{pt}}$. One can work over $ [\widehat{\mathcal{L}H}\times \widehat{\mathcal{L}H}\times \widehat{\mathcal{L}H}]_{\on{Ran}_{\on{un,dR}}}$ to check these equivalences provide a multiplicative structure\footnote{By Remark \ref{rem-fact-line-bunlde-being-mult-property}, this is a property rather than additional structure.} on $\mathscr{L}^{\on{diff}'}$.

Now consider the projection $\pi: [\widehat{\mathcal{L}H}]_{\on{Ran}_{\on{un,dR}}} \to  [\widehat{\on{Hecke}}_H ]_{\Ran_{\on{un,dR}}}$. Using the base-change isomorphisms, we have $\pi^*_\ren \circ \theta_i(\omega) \simeq \mathcal{F}_i$. This implies $\pi^{!,\ren}( \mathscr{L}^{\on{diff}}_{ \widehat{\on{Hecke}}_H}  ) \simeq  \mathscr{L}^{\on{diff}'}$ as unital factorization line bundles. One can check that this equivalence is compatible with the multiplicative structures.

\emph{Step 6}: In this last step, we show $\mathscr{L}^{\on{diff}} \simeq \mathscr{L}^{\on{Tate}}$.

Consider the autmorphism $\alpha=(\on{mult},\on{pr}_2)$ on $[\widehat{\mathcal{L}H}\times \widehat{\mathcal{L}H}]_{\on{Ran}_{\on{un,dR}}}$. Using the base-change isomorphisms, we have $\alpha^*_\ren(  \mathcal{F}_1 \boxtimes  \mathcal{F}_1 )\simeq  \mathcal{F}_1 \boxtimes  \mathcal{F}_1$. On the other hand, we claim 
\begin{equation} \label{eqn-proof-thm-tate-twist-1}
\alpha^*_\ren(  \mathcal{F}_0 \boxtimes  \mathcal{F}_0 )\simeq(\omega \boxtimes \mathscr{L}^{\on{Tate}})\overset{\on{act}}\otimes (\mathcal{F}_0\boxtimes\mathcal{F}_0).
\end{equation}

To simplify the notations, we prove the non-factorization version of the claim. The factorization version can be proved in the same way be replacing $\widehat{\mathcal{L}H}$ by $[\widehat{\mathcal{L}H}]_s$ for any classical test map $s:S\to \Ran_{ \on{un,dR} }$. We view $\widehat{\mathcal{L}H}\times \widehat{\mathcal{L}H}$ as a group indscheme parameterized by its second factor and consider its closed sub-group scheme $\mathcal{K}$ containing those $(u,v)$ with $\on{Ad}_v u\in \mathcal{L}^+H$. Then we have a commutative diagram
$$
\xymatrix{
    \widehat{\mathcal{L}H}\times \widehat{\mathcal{L}H} \ar[r]^-{\beta} \ar[d]^-\alpha &
    (\widehat{\mathcal{L}H}\times \widehat{\mathcal{L}H})/\mathcal{K} \ar[r]^-\gamma \ar[d]^-{\overline{\alpha}} &
    \widehat{\mathcal{L}H} \ar[d] \\
     \widehat{\mathcal{L}H}\times \widehat{\mathcal{L}H} \ar[r] &
    \widehat{\mathcal{L}H}/\mathcal{L}^+H \times \widehat{\mathcal{L}}H/\mathcal{L}^+H \ar[r]^-{\on{pr}_2} &
    \widehat{\mathcal{L}H}/\mathcal{L}^+H 
}
$$
such that the right square is Cartesian. Using the base-change isomorphisms, we obtain 
\begin{equation} \label{eqn-proof-thm-tate-twist-2}
\alpha^*_\ren(  \mathcal{F}_0 \boxtimes  \mathcal{F}_0 )\simeq \beta^{*}_\ren\circ \gamma^{?}_\ren(\mathcal{F}_0).
\end{equation}
On the other hands, we can also view $\mathcal{K}':= \mathcal{L}^+H \times \widehat{\mathcal{L}H} $ as a closed sub-group scheme of $\widehat{\mathcal{L}H}\times \widehat{\mathcal{L}H}$ parameterized by their second factors. Consider the chain $\widehat{\mathcal{L}H}\times \widehat{\mathcal{L}H} \xrightarrow{\beta'}
    (\widehat{\mathcal{L}H} \times \widehat{\mathcal{L}H})/\mathcal{K}' \xrightarrow{\gamma'}
    \widehat{\mathcal{L}H}$. We claim
\begin{equation} \label{eqn-proof-thm-tate-twist-3}
 (\beta')^{*}_\ren\circ (\gamma')^{?}_\ren(  \mathscr{L}^{\on{Tate}}\overset{\on{act}}\otimes- ) \simeq (\beta)^{*}_\ren\circ (\gamma)^{?}_\ren(-). 
 \end{equation}
To prove this claim, consider the (classical) intersection $\mathcal{K}''$ of $\mathcal{K}$ and $\mathcal{K}'$ and the chains
$$ \mathbb{B}\mathcal{K}'' \xrightarrow{\phi} \mathbb{B}\mathcal{K} \xrightarrow{\xi}  \widehat{\mathcal{L}H},\; \mathbb{B}\mathcal{K}'' \xrightarrow{\phi'} \mathbb{B}\mathcal{K}' \xrightarrow{\xi'}   \widehat{\mathcal{L}H}.$$
Note that $\phi$ and $\phi'$ are of finite type and are nil-isomorphisms. Using Proposition \ref{prop-indcoh-on-ind-stack} and Proposition \ref{prop-indcoh-on-quot-stack}, one can deduce the following formula from its analogue in the finite type case:
\begin{equation} \label{eqn-proof-thm-tate-twist-4} 
\phi^?_\ren(-) \simeq \phi^*_\ren( \xi^{!,\ren}( \mathcal{L})\overset{\on{act}}\otimes- ),\; (\phi')^?_\ren(-) \simeq (\phi')^*_\ren( (\xi')^{!,\ren}( \mathcal{L}')\overset{\on{act}}\otimes- )
\end{equation}
where $\mathcal{L}, \mathcal{L}'\in \IndCoh^{!,\ren}(  \widehat{\mathcal{L}H} )$ are (cohomologically shifted) line bundles whose restrictions along a classical test map $s:S\to \widehat{\mathcal{L}H}$ are respectively the determinant line bundles 
$$\on{det}( \on{Lie}(\mathcal{K}|_s)/\on{Lie}(\mathcal{K}''|_s) [1]),\;\on{det}( \on{Lie}(\mathcal{K}'|_s)/\on{Lie}(\mathcal{K}''|_s) [1]).$$
By definition, $\mathcal{L}' \simeq \mathscr{L}^{\on{Tate}} \otimes \mathcal{L}$. Combining this with (\ref{eqn-proof-thm-tate-twist-4}), we obtain the desired equivalence (\ref{eqn-proof-thm-tate-twist-3}).

Now the desired equivalence (\ref{eqn-proof-thm-tate-twist-1}) follows from (\ref{eqn-proof-thm-tate-twist-2}) and (\ref{eqn-proof-thm-tate-twist-3}). Pluging in $\mathcal{F}_0 \simeq \mathscr{L}^{\on{diff}}\overset{\on{act}}\otimes \mathcal{F}_1 $, we obtain
\begin{eqnarray*}
(\omega \boxtimes \mathscr{L}^{\on{Tate}})\overset{\on{act}}\otimes (\mathcal{F}_0\boxtimes\mathcal{F}_0) & \simeq & \alpha^*_\ren(  \mathcal{F}_0 \boxtimes  \mathcal{F}_0 ) \\
& \simeq & \alpha^{!,\ren}(  \mathscr{L}^{\on{diff}}\boxtimes \mathscr{L}^{\on{diff}} )\overset{\on{act}} \otimes \alpha^*_\ren (  \mathcal{F}_1 \boxtimes  \mathcal{F}_1 ) \\
& \simeq & (  \mathscr{L}^{\on{diff}}\boxtimes (\mathscr{L}^{\on{diff}})^{\otimes 2} )\overset{\on{act}} \otimes (  \mathcal{F}_1 \boxtimes  \mathcal{F}_1 )  \\
& \simeq & (\omega \boxtimes \mathscr{L}^{\on{diff}})\overset{\on{act}}\otimes (\mathcal{F}_0\boxtimes\mathcal{F}_0).
\end{eqnarray*}
Since acting on $(\mathcal{F}_0\boxtimes\mathcal{F}_0)$ is an equivalence, we obtain $\omega \boxtimes \mathscr{L}^{\on{Tate}}\simeq \omega \boxtimes \mathscr{L}^{\on{diff}}$ and therefore $\mathscr{L}^{\on{Tate}}\simeq \mathscr{L}^{\on{diff}}$ as unital factorization line bundles. One can work over $ [\widehat{\mathcal{L}H}\times \widehat{\mathcal{L}H}\times \widehat{\mathcal{L}H}]_{\on{Ran}_{\on{un,dR}}}$ to check this equivalence is compatible with the multiplicative structures. This finishes the proof of the theorem.

\qed

\begin{remark} \label{rem-modular-character-compare-with-sam}
In this remark, we sketch how to identify the restriction of our equivalence in Theorem \ref{thm-tate-twist} at a closed point $x\in X$ with the equivalence constructed in \cite{raskin2020homological}.

We need to translate Raskin's construction to our setting. \cite[Lemma 7.17.4]{raskin2020homological} says the monoidal category $A:=[\IndCoh^{!,\ren}( [\widehat{\on{Hecke}}_H ]_x ),\star]$ is rigid\footnote{This is related to the following fact, which we do not prove nor use in this paper. Consider the lax unital factorization functor:
$$ \lambda:  \mathcal{I}ndCoh^{!,\ren}( \widehat{\on{Hecke}}_H^1) \overset{e^{!,\ren}}\to  \mathcal{I}ndCoh^{!,\ren}( \widehat{\on{Hecke}}_H^0) \overset{?\on{-push}} \to \mathcal{I}ndCoh^{!,\ren}( \on{pt} ) \simeq \mathbf{1},$$
where $e: \widehat{\on{Hecke}}_H^0\to \widehat{\on{Hecke}}_H^1$ is the morphism corresponding to the unique map $[1]\to [0]$. Then $[ \mathcal{I}ndCoh^{!,\ren}( \widehat{\on{Hecke}}_H),\star ,\lambda]$ is a Frobenius algebra object in $[\on{FactCat}_{\on{un}}]_{\on{lax-unital}}$ (see \cite[Definition 4.6.5.1]{HA}). A rigid monoidal DG category is in particular a Frobenius algebra object in $\on{DGCat}_{\on{cont}}$, and Raskin's construction actually only used this Frobenius algebra structure.

One can use this fact to obtain an equivalence $\mathcal{KL}(H)_{ \mathscr{L} }  \simeq \mathcal{KL}^*(H)_{\mathscr{L}\otimes \mathscr{L}^{\on{Serre}} }$, where $\mathscr{L}^{\on{Serre}}$ is a multiplicative factorization line bundle such that the \emph{Serre automorphism} for the above Frobenius algebra (see \cite[Remark 4.6.5.4]{HA}) is given by tensoring up with $\mathscr{L}^{\on{Serre}}$. However, it is harder to identify $\mathscr{L}^{\on{Serre}}$ with $\mathscr{L}^{\on{Tate}}$. For instance, Raskin proved (the non-factorization version) of this identification in \cite[Theorem 9.16.1]{raskin2020homological} with the help of semi-infinite cohomologies.}. Consider the \emph{geometric} left and right actions of $A$ on $M:=\on{Rep}( \mathcal{L}^+ H_x )$ defined respectively by (\ref{eqn-geometric-left-action-Hecke-on-Rep}) and (\ref{eqn-geometric-right-action-Hecke-on-Rep}), it is clear that these actions commute hence make $M$ an $(A,A)$-bimodule. One can verify that tensoring up with $[\mathscr{L}_{ \widehat{\on{Hecke}_H} }]_{x}$ is a monoidal automorphism on $A$. Using this automorphism, we can twist the $A$-actions on $\on{Rep}( \mathcal{L}^+ H_x )$ and obtain another $(A,A)$-bimodule $M_{\mathscr{L}}$.

In Raskin's language, the left $A$-module structure on $M_{\mathscr{L}}$ is the same as a weak $\widehat{\mathcal{L}H}_x$-action on $\on{Vect}$. Denote the corresponding weak $\widehat{\mathcal{L}H}_x$-module by $C_{\mathscr{L}}$. The invariance and coinvariance are defined respectively as $(C_{\mathscr{L}})^{\on{inv}}:=\on{Hom}_{A\mod^l}( M, M_{\mathscr{L}} )$ and $(C_{\mathscr{L}})_{\on{coinv}}:=M\otimes_A M_{\mathscr{L}}$. One can verify that the cosimplicial (resp. simplicial) diagram defining them satsifies the Beck-Chevalley conditions, which implies 
$$(C_{\mathscr{L}})^{\on{inv}} \simeq [\mathcal{KL}(H)_{ \mathscr{L} }]_x ,\;(C_{\mathscr{L}})_{\on{coinv}}\simeq [ \mathcal{KL}^*(H)_{\mathscr{L}}]_x.$$
In \cite[Proposition 7.18.2]{raskin2020homological}, Raskin used the general theory of rigid monoidal categories to obtain an equivalence 
$$(C_{\mathscr{L}})^{\on{inv}}\simeq  (C_{\mathscr{L}}\otimes \chi_{\on{Serre}})_{\on{coinv}},$$
where $\chi_{\on{Serre}}$ is the weak $\widehat{\mathcal{L}H}_x$-module twisted from the trivial one by the \emph{Serre automorphism} $\on{Se}_A$ (see \cite[Remark 4.6.5.4]{HA}) on the rigid monoidal category $A$. Hence we only need to identify $\on{Se}_A$ with $[\mathscr{L}^{\on{Tate}}_{ \widehat{\on{Hecke}_H} }]_x \otimes^! (-)$ as symmetric monoidal functors. Recall the underlying endo-functor of $\on{Se}_A$ is the difference of two self-dualities on $A$ whose pairing functors are 
$$A\otimes A \xrightarrow{\star} A \xrightarrow{\on{Maps}( \on{unit}_\star,- )} \on{Vect},\;A\otimes A \xrightarrow{\star^{\on{rev}}} A \xrightarrow{\on{Maps}( \on{unit}_\star,- )} \on{Vect}.$$
Using the base-change isomorphisms, we have $\on{Maps}( \on{unit}_\star, x\star y ) \simeq <\sigma^{!,\ren}(x),y>_1$, where $<-,->_1$ is the pairing in Step 1 of our proof of Theorem \ref{thm-tate-twist} for the case $i=m=1$. Similarly we have $\on{Maps}( \on{unit}_\star, y\star x ) \simeq <\sigma^{!,\ren}(x),y>_0$. Since $\sigma^{!,\ren}$ is an involution, we obtain $\on{Se}_A\simeq \theta_1^{-1}\circ \theta_0$. Hence our proof of Theorem \ref{thm-tate-twist} provides an equivalence $\on{Se}_A\simeq [\mathscr{L}^{\on{Tate}}_{ \widehat{\on{Hecke}_H} }]_x \otimes^! (-)$ as endo-functors. Now the monoidal structure on $\on{Se}_A$ would correspond to \emph{some} multiplicative structure on $\mathscr{L}^{\on{Tate}}$. We only need to show it coincides with the standard one. But this follows from unwinding the definitions.
\end{remark}

\subsection{Variant: Factorization Module Category}
\label{ssec-fact-mod-KL}
We also need a factorization module version of the previous constructions. We introduce the following notations.

\begin{notation} \label{notn-fact-module-km}
We use the following notations:
\begin{itemize}
    \item $K\to H$ is a fixed closed subgroup of $H$.
    \item $x:\on{pt}\to \Ran_{\on{un}}$ is a fixed point of $\Ran_{\on{un}}$, which corresponds to a finite set $I_0\subset \on{Maps}(\on{pt},X)$, or equivalently, to closed points $x_i$ of $X$ labelled by $I_0$.
    \item $\Ran_{x,\on{un}}$ is the laft lax prestack defined in Section \ref{sect:factspc}, and $\Ran_{x,\on{un,dR}}$ is its de-Rham lax prestack.
    \item $[\mathcal{L}^+ H]_{\Ran_{x,\on{un,dR}}}:= [\mathcal{L}^+ H]_{\Ran_{\on{un,dR}}} \times_{  \Ran_{\on{un,dR}} } \Ran_{x,\on{un,dR}}$, and similarly define $[\widehat{\mathcal{L} H}]_{\Ran_{x,\on{un,dR}}}$.
    \item $H_x:= \prod_{j\in J} H_{x_j}:= \prod_{j\in J} H$ and similarly define $K_x$.
    \item $[\mathcal{L}^+ H\times_{H_x} K_x ]_{ \Ran_{x,\on{un,dR}} }:= [\mathcal{L}^+ H]_{\Ran_{x,\on{un,dR}}} \times_{H_x} K_x$, which is a co-unital factorization module space of $\mathcal{L}^+ H_{\Ran_{\on{un,dR}}}$ supported on $x$.
    \item $[\widehat{\mathcal{L}H}\times_{ H_{x,\dR} } K_{x,\dR}  ]_{ \Ran_{x,\on{un,dR}}  }:=[\widehat{\mathcal{L} H}]_{\Ran_{x,\on{un,dR}}} \times_{H_{x,\dR}} K_{x,\dR} $ is the formal completion of $[\mathcal{L} H]_{\Ran_{x,\on{un,dR}}}$ along $[\mathcal{L}^+ H\times_{H_x} K_x ]_{ \Ran_{x,\on{un,dR}} }$, which is a corr-unital factorization module space of $[\mathcal{L} H]_{\Ran_{\on{un,dR}}}$ supported on $x$.

    \item $[ \widehat{\on{Hecke}}_{H}^{K\on{-level}_x,\bullet}  ]_{\Ran_{ x,\on{un,dR} }}$ is the Cech nerve of $ [\mathbb{B}(\mathcal{L}^+ H\times_{H_x} K_x)]_{\Ran_{x,\on{un,dR}}} \to [\mathbb{B}(\widehat{\mathcal{L}H}\times_{ H_{x,\dR} } K_{x,\dR} )]_{\Ran_{x,\on{un,dR}}}$, which is a cosimplicial corr-unital factorization module space of $[\widehat{\on{Hecke}}_H^\bullet ]_{\Ran_{ \on{un,dR} }}$.
\end{itemize}
\end{notation}

\begin{remark} \label{rem-level-structure-loc-hecke}
By definition, we have
$$ [ \widehat{\on{Hecke}}_{H}^{K\on{-level}_x,\bullet}  ]_{\Ran_{ x,\on{un,dR} }} \simeq [ \widehat{\on{Hecke}}_{H}^\bullet  ]_{\Ran_{ x,\on{un,dR} }} \times_{ \mathbb{B}H_x^{\times (\bullet+1) }  }  \mathbb{B} K_x^{\times (\bullet+1) }.$$
\end{remark}

\begin{variant} 
\label{variant:gmodI}
As in the definitions of $\mathcal{KL}(H)_{\mathscr{L}}$ and $ \mathcal{KL}(H)_{\mathscr{L},\on{co}}$, we can define crystal of categories over $\Ran_{x,\on{un,dR}}$:
$$  \widetilde{\mathcal{L}\mathfrak{h}}\on{-}mod^{\mathcal{L}^+H\times_{H_x} K_x }_{\mathscr{L}}  ,\;   \widetilde{\mathcal{L}\mathfrak{h}}\on{-}mod^{\mathcal{L}^+H\times_{H_x} K_x }_{\mathscr{L},\on{co}} ,$$
equipped respectively with unital factorization module structures for the unital factorization categories $\mathcal{KL}(H)_{\mathscr{L}}$ and $\mathcal{KL}(H)_{\mathscr{L},\on{co}}$. 

We can also define a crystal of categories over $\Ran_{x,\on{un,dR}}$:
$$  \mathcal{L}^+\mathfrak{h}\on{-}mod^{\mathcal{L}^+H\times_{H_x} K_x }$$
equipped rwith unital factorization module structures for the unital factorization categories $\mathcal{R}ep( \mathcal{L}^+H )$. 
\end{variant}

\begin{example} Recall a map $y:\on{pt}\to \Ran_{x,\on{un}}$ corresponds to a finite subset $I\subset \on{Maps}(\on{pt},X)$ containing $I_0$. Then 
$$[\widetilde{\mathcal{L}\mathfrak{h}}\on{-}mod^{\mathcal{L}^+H\times_H K}_{\mathscr{L}}  ]_{y} \simeq \bigotimes_{x_i\in I_0} \widetilde{\mathcal{L}\mathfrak{h}}_{x_i}\mod^{\mathcal{L}^+H_{x_i}\times_{H_{x_i}} K_{x_i}}_{\mathscr{L}} \bigotimes \bigotimes_{x_i\in I\setminus I_0} \on{KL}(H)_{\mathscr{L},x_i} .$$
\end{example}

\subsection{Functoriality For the Group}
\begin{construction} \label{constr-KL-functorial}
Let $(K'\to H')$ be another closed embedding between finite type affine algebraic groups, and $(K'\to H')\to (K\to H)$ be group homomorphisms. Note that $\mathscr{L}$ induces a multiplicative factorization line bundle $\mathscr{L}'$ on $\mathcal{L}H'$.

Then the $!$-pullback functors induce lax unital morphisms
$$\on{res}: \mathcal{KL}(H)_{\mathscr{L}} \to \mathcal{KL}(H')_{\mathscr{L}'},  \;\on{res}: \widetilde{\mathcal{L}\mathfrak{h}}\on{-}mod^{\mathcal{L}^+H\times_{H_x} K_x }_{\mathscr{L}} \to  \widetilde{\mathcal{L}\mathfrak{h}'}\on{-}mod^{\mathcal{L}^+H'\times_{H'_x} K'_x }_{\mathscr{L}'}$$
compatible with the factorization structures, which we refer as the \emph{factorization restriction functors}.

Dually, the $*$-pushforward functors lax unital morphisms
$$\on{C}^{\frac{\infty}{2}}_{\on{co}}:  \mathcal{KL}(H')_{\mathscr{L}',\on{co}} \to \mathcal{KL}(H)_{\mathscr{L},\on{co}}, \;\on{C}^{\frac{\infty}{2}}_{\on{co}}:  \widetilde{\mathcal{L}\mathfrak{h}'}\on{-}mod^{\mathcal{L}^+H'\times_{H'_x} K'_x }_{\mathscr{L}',\on{co}} \to \widetilde{\mathcal{L}\mathfrak{h}}\on{-}mod^{\mathcal{L}^+H\times_{H_x} K_x }_{\mathscr{L},\on{co}}$$
compatible with the factorization structures, which we call the \emph{co-version of the factorization semi-infinite cohomology functor}.

Using (the factorization module version of) Theorem \ref{thm-tate-twist}, we obtain a morphism
$$\on{C}^{\frac{\infty}{2}}:  \mathcal{KL}(H')_{\mathscr{L}'\otimes \mathscr{L}^{-\on{Tate}}_{H'}} \to \mathcal{KL}(H)_{\mathscr{L}\otimes \mathscr{L}^{-\on{Tate}}_{H}}, \;\on{C}^{\frac{\infty}{2}}:  \widetilde{\mathcal{L}\mathfrak{h}'}\on{-}mod^{\mathcal{L}^+H'\times_{H'_x} K'_x }_{\mathscr{L}'\otimes \mathscr{L}^{-\on{Tate}}_{H'}} \to \widetilde{\mathcal{L}\mathfrak{h}}\on{-}mod^{\mathcal{L}^+H\times_{H_x} K_x }_{\mathscr{L}\otimes \mathscr{L}^{-\on{Tate}}_{H}}$$
compatible with the factorization structures, which we call the \emph{the factorization semi-infinite cohomology functor}.
\end{construction}

\section{Localization and Local-Global Correspondence}

\label{sect:km-localization}

\begin{notation} \label{notn-def-localization}
We keep Notation \ref{notn-fact-module-km}. Also,
\begin{itemize}
	\item $X$ is assumed to be proper.

	\item The multiplicative factorization line bundle $\mathscr{L}$ is constructed from a central $\on{Lie}^*$-extension $\widetilde{\mathfrak{h}}_{\mathcal{D}}$ of $\mathfrak{h}_{\mathcal{D}}:=\mathfrak{h}\otimes \on{Diff}_X$ using the method in \cite{zhao2017quantum}.

	\item $\on{Bun}_H$ is the moduli stack of $H$-torsors on $X$.
	 
	\item Let $(\on{Bun}_H\times\on{Ran}_{\on{un,dR}})^{H\on{-level}_{\infty}}$ be the prestack such that $(\on{Bun}_H\times\on{Ran}_{\on{un,dR}})^{H\on{-level}_{\infty}}(S)$ classifies
\begin{itemize}
	\item a finite subset $I\subset \on{Maps}(S^{\on{cl,red}},X)$;
	\item[] Let $\Gamma_I\to S\times X$ be the schema-theoretic sum of the graphs of the maps $S^{\on{cl,red}} \to X$;
	\item an $H$-torsor on $X\times S$ equipped with a trivialization on the adic multi-disk $\mathcal{D}_{\Gamma_I}'$.
\end{itemize}

	\item $\on{Bun}_H ^{K\on{-level}_x} \simeq \on{Bun}_H\times_{\mathbb{B}H_x } \mathbb{B}K_x$ is the moduli stack of $H$-torsors on $X$ equipped with $K$-structures at $x_j, j\in J$.

	\item $\mathscr{T}$ is the twisting on $\on{Bun}_H$ constructed from the central $\on{Lie}^*$-extension $\widetilde{\mathfrak{h}}_{\mathcal{D}}$ using the method in \cite{zhao2017quantum}.

	\item $\DMod(\on{Bun}_H ^{K\on{-level}_x})_{\mathscr{T}}$ is the DG category of $\mathscr{T}$-twisted D-modules on $\on{Bun}_H ^{K\on{-level}_x}$.

\end{itemize}
\end{notation}

\begin{remark} \label{rem-compatible-pair-two-gerbe}
Consider 
$$(\on{Bun}_H\times\on{Ran}_{x,\on{un,dR}})^{H\on{-level}_{\infty}}:=(\on{Bun}_H\times\on{Ran}_{\on{un,dR}})^{H\on{-level}_{\infty}}\times_{\on{Ran}_{\on{un,dR}} } \on{Ran}_{x,\on{un,dR}},$$
which is an $[\mathcal{L}^+ H\times_{H_x} K_x]_{\on{Ran}_{x,\on{un,dR}} }$-torsor on $\on{Bun}_H^{K\on{-level}_x}\times \on{Ran}_{x,\on{un,dR}}$. A standard re-gluing construction extend this $[\mathcal{L}^+ H\times_{H_x} K_x]_{\on{Ran}_{x,\on{un,dR}} }$-action to an $[\mathcal{L}H]_{\on{Ran}_{x,\on{un,dR}}}$-action. The restricted $[\widehat{\mathcal{L}H}\times_{H_{x,\dR}} K_{x,\dR}]_{\on{Ran}_{x,\on{un,dR}}}$-action preserves the projection to $(\on{Bun}_{H}^{K\on{-level}_x})_\dR\times\on{Ran}_{\on{un,dR}}$. Hence $\mathscr{T}$ induces a $\mathbb{G}_m$-gerbe on 
$$(\on{Bun}_H\times\on{Ran}_{x,\on{un,dR}})^{H\on{-level}_{\infty}}/[\widehat{\mathcal{L}H}\times_{H_{x,\dR}} K_{x,\dR}]_{\on{Ran}_{x,\on{un,dR}}}$$
equipped with a trivialization on 
$$(\on{Bun}_H\times\on{Ran}_{x,\on{un,dR}})^{H\on{-level}_{\infty}}/[\mathcal{L}^+H\times_{H_{x,\dR}} K_{x,\dR}]_{\on{Ran}_{x,\on{un,dR}}} \simeq  \on{Bun}_H\times_{\mathbb{B}H_{x,\dR} } \mathbb{B}K_{x,\dR} \times \on{Ran}_{x,\on{un,dR}} $$
Note that $\on{Bun}_H\times_{\mathbb{B}H_{x,\dR} } \mathbb{B}K_{x,\dR} $ is just the formal completion of $\on{Bun}_H$ along $\on{Bun}_H ^{K\on{-level}_x}$.

On the other hand, recall that the multiplicative factorization line bundle $\mathscr{L}$ induces a $\mathbb{G}_m$-gerbe on $[\mathbb{B}\widehat{\mathcal{L}H}]_{\on{Ran}_{\on{un,dR}}}$ equipped with a trivialization on $[\mathbb{B}\mathcal{L}^+H]_{\on{Ran}_{\on{un,dR}}}$. By restriction, we obtain another $\mathbb{G}_m$-gerbe on $(\on{Bun}_H\times\on{Ran}_{x,\on{un,dR}})^{H\on{-level}_{\infty}}/[\widehat{\mathcal{L}H}\times_{H_{x,\dR}} K_{x,\dR}]_{\on{Ran}_{x,\on{un,dR}}}$ equipped with a trivialization on $\on{Bun}_H\times_{\mathbb{B}H_{x,\dR} } \mathbb{B}K_{x,\dR} \times \on{Ran}_{x,\on{un,dR}} $.

It follows from the assumption that the above two $\mathbb{G}_m$-gerbes are canonical isomorphic.
\end{remark}

\begin{remark} \label{rem-pullback-compatible-pair}
Let $H'\to H$ be a group homomorphism. Then $\widetilde{\mathfrak{h}}_{\mathcal{D}}$ induces a central $\on{Lie}^*$-extension $\widetilde{\mathfrak{h}'}_{\mathcal{D}}$ of $\mathfrak{h}'_\mathcal{D}$. We often write $(\mathscr{T}', \mathscr{L}')$ for the corresponding twisting and factorization multiplicative line bundle.
\end{remark}

\begin{example}\label{exam-tate-det-local-global}
  Since $\on{Bun}_H $ is smooth, the dualizing complex $K\in \on{QCoh}(\on{Bun}_H )$ is invertible. Let $K[-\on{dim}( \on{Bun}_H )]$ be the corresponding line bundle, where $\on{dim}( \on{Bun}_H )$ is the locally constant function of dimension. Let $\mathscr{K}$ be the integral twisting on $\on{Bun}_H$ given by this line bundle. It follows from definitions that $-\mathscr{K}$ corresponds to the Tate central $\on{Lie}^*$-extension.
\end{example}

The following result summarizes what we need to know about chiral localization.

\begin{theorem} \label{thm-localization-KL}
There are canonical unital morphisms
\begin{eqnarray*}
 \on{Loc}:  &
  \widetilde{\mathcal{L}\mathfrak{h}}\on{-}mod^{\mathcal{L}^+H\times_{H_x} K_x }_{\mathscr{L}} \to \DMod(\on{Bun}_H ^{K\on{-level}_x})_{\mathscr{T}} \boxtimes \mathbf{1},\\
\on{Loc}_{\on{co}}: &
  \widetilde{\mathcal{L}\mathfrak{h}}\on{-}mod^{\mathcal{L}^+H\times_{H_x} K_x }_{\mathscr{L},\on{co}} \to \DMod(\on{Bun}_H ^{K\on{-level}_x})_{\mathscr{T}} \boxtimes \mathbf{1}
\end{eqnarray*}
between crystal of categories on $\on{Ran}_{x,\on{un,dR}}$, where the RHS's are the constant crystals of categories. These morphisms are required to satisfy the following:

\begin{itemize}
	\item[(1)](Compatibility with chiral homology)
	 Let 
	$$[\on{Loc}]_{ x}: [ \widetilde{\mathcal{L}\mathfrak{h}}\on{-}mod^{\mathcal{L}^+H\times_{H_x} K_x }_{\mathscr{L}}  ]_x  \to \DMod(\on{Bun}_H ^{K\on{-level}_x})_{\mathscr{T}}  $$
	be the fiber of $\on{Loc}$ at $x:\on{pt}\to \on{Ran}_{x,\on{un,dR}}$. Then there is a canonical commutative diagram
	$$
	\xymatrix{
	[ \widetilde{\mathcal{L}\mathfrak{h}}\on{-}mod^{\mathcal{L}^+H\times_{H_x} K_x }_{\mathscr{L}}  ]_x  \ar[rr]^-{ [\on{Loc}]_{ x} }  \ar[d]^-{\on{ins}_{\on{unit}}}
	& & \DMod(\on{Bun}_H ^{K\on{-level}_x})_{\mathscr{T}} \\
	[ \widetilde{\mathcal{L}\mathfrak{h}}\on{-}mod^{\mathcal{L}^+H\times_{H_x} K_x }_{\mathscr{L}}  ]_{\on{Ran}_{x,\on{un,dR}}}  \ar[rr]^-{ [\on{Loc}]_{ \on{Ran}_{x,\on{un,dR}}  } } 
	& & \DMod(\on{Bun}_H ^{K\on{-level}_x})_{\mathscr{T}} \otimes \DMod(\on{Ran}_{x,\on{un}}) \ar[u]^-{\mathbf{Id}\otimes \on{C}^{\on{ch}}_{x,\on{un}} }
	},
	$$
	where the left vertical functor is the \emph{inserting unit functor} (see Construction \ref{const-prop-unit=insert-unit} for its definition), and 
	$$\on{C}^{\on{ch}}_{x,\on{un}}:=p_!: \DMod(\on{Ran}_{x,\on{un}})\to \on{Vect}$$
	is \emph{taking chiral homology}, i.e., the $!$-pushforward\footnote{We remind that the well-definedness of this functor relies on the \emph{homological left cofinality} of $\Ran_\dR \to \Ran_{\un, \dR}$; see \cite{gaitsgory2015atiyah} for more details.} functor along $p: \on{Ran}_{x,\on{un}}\to \on{pt}$. There is a similar canonical commutative diagram for $\on{Loc}_{\on{co}}$.

	\item[(2)](Compatibility with Tate twist) When $x$ corresponds to the empty subset $\emptyset\subset \on{Maps}( \on{pt},X )$, there is a canonical commutative diagram
	$$
	\xymatrix{
	 \mathcal{KL}(H)_{\mathscr{L}\otimes \mathscr{L}^{\on{Tate}},\on{co} }
	\ar[r]^-{ \on{Loc}_{\on{co}} }  \ar[d]^-\simeq
	&  \DMod(\on{Bun}_H )_{\mathscr{T}-\mathscr{K}}  \boxtimes \mathbf{1}  \ar[d]^-{\simeq}_-{K\otimes - }
	\\
	\mathcal{KL}(H)_{ \mathscr{L} }
	\ar[r]^-{ \on{Loc} } 
	&  \DMod(\on{Bun}_H )_{\mathscr{T}} \boxtimes \mathbf{1},
	}
	$$
	where the left vertical equivalence is due to Theorem \ref{thm-tate-twist}, and the right equivalence is bacause the twisting $\mathscr{K}$ is integral (i.e., given by a shifted line bundle).

	\item We use 
	$$[\on{Loc}]_x^{\on{ch}}: [ \widetilde{\mathcal{L}\mathfrak{h}}\on{-}mod^{\mathcal{L}^+H\times_{H_x} K_x }_{\mathscr{L}}  ]_{\on{Ran}_{x,\on{un,dR}}} \to \DMod(\on{Bun}_H ^{K\on{-level}_x})_{\mathscr{T}}   $$
	to denote the composition $  (\mathbf{Id}\otimes \on{C}^{\on{ch}}_{x,\on{un}} )\circ   [\on{Loc}]_{ \on{Ran}_{x,\on{un,dR}}  } $.
	For a unital factorization algebra $A$ in $ \mathcal{KL}(H)_{\mathscr{L}}$, we abuse notation and write
	$$ [\on{Loc}]_x^{\on{ch}}: A\on{-FactMod}_{\on{un}}( \widetilde{\mathcal{L}\mathfrak{h}}\on{-}mod^{\mathcal{L}^+H\times_{H_x} K_x }_{\mathscr{L}}  )\to \DMod(\on{Bun}_H ^{K\on{-level}_x})_{\mathscr{T}}   $$
	for the composition $[\on{Loc}]_x^{\on{ch}}\circ \oblv_{A}$.
	
	\item[(3)](Compatibility with restriction) Let $(K'\to H')\to (K\to H)$ be as in Construction \ref{constr-KL-functorial}. Then for any unital factorization algebra $A$ in $ \mathcal{KL}(H)_{\mathscr{L}}$, there is a canonical commutative diagram

$$
	\xymatrix{
	A\on{-FactMod}_{\on{un}}( \widetilde{\mathcal{L}\mathfrak{h}}\on{-}mod^{\mathcal{L}^+H\times_{H_x} K_x }_{\mathscr{L}}  ) \ar[d]^-{\on{res}} 
	\ar[r]^-{[\on{Loc}]_x^{\on{ch}} } 
	& \DMod(\on{Bun}_H ^{K\on{-level}_x})_{\mathscr{T}}  \ar[d]^-{!\on{-pull}} \\
	\on{res}(A)\on{-FactMod}_{\on{un}}( \widetilde{\mathcal{L}\mathfrak{h}'}\on{-}mod^{\mathcal{L}^+H'\times_{H'_x} K'_x }_{\mathscr{L}'}  ) 
	\ar[r]^-{[\on{Loc}]_x^{\on{ch}} } 
	& \DMod(\on{Bun}_{H'}^{K'\on{-level}_x})_{\mathscr{T}'},
	}
	$$
	where the left vertical functor is induced by the lax unital factorization functor $\on{res}$ in Construction \ref{constr-KL-functorial}.

	\item[(4)](Compatibility with semi-infinite cohomology) In the situation of (3), $K'=H',K=H$ and $H'\to H$ is a surjection with unipotent kernel. Then for any unital factorization algebra $A'$ in $ \mathcal{KL}(H')_{\mathscr{L}',\on{co}}$, there is a canonical commutative diagram
	$$
	\xymatrix{
	A'\on{-FactMod}_{\on{Ran}_{x,\on{un,dR}}}    \ar[d]^-{\on{C}^{\frac{\infty}{2}}_{\on{co}}}  \ar[r]^-{[\on{Loc}]_x^{\on{ch}} } 
	& \DMod(\on{Bun}_{H'})_{\mathscr{T}'}  \ar[d]^-{*\on{-push}} \\
	\on{C}^{\frac{\infty}{2}}_{\on{co}}(A')\on{-FactMod}_{\on{Ran}_{x,\on{un,dR}}}
	\ar[r]^-{[\on{Loc}]_x^{\on{ch}} } 
	& \DMod(\on{Bun}_H)_{\mathscr{T}},
	}
	$$
	where the left vertical functor is induced by the lax unital factorization functor $\on{C}^{\frac{\infty}{2}}_{\on{co}}$ in Construction \ref{constr-KL-functorial}.
\end{itemize}
\end{theorem}

\begin{remark} Consider the following correspondence
\[
\begin{aligned}
  (\on{Bun}_H^{K\on{-level}_x})_\dR \times \on{Ran}_{x,\on{un,dR}} \gets  (\on{Bun}_H\times\on{Ran}_{x,\on{un,dR}})^{H\on{-level}_{\infty}}/[\widehat{\mathcal{L}H}\times_{H_{x,\dR}} K_{x,\dR}]_{\on{Ran}_{x,\on{un,dR}}} \to\\
  \to [\mathbb{B}\widehat{\mathcal{L}H}\times_{ H_{x,\dR} } K_{x,\dR}  ]_{ \Ran_{x,\on{un,dR}}  }
 \end{aligned}
 \]
The morphism $\on{Loc}$ should be viewed as the pull-push along this correspondence for a suitable twisted renormalized $\IndCoh$-theory.
\end{remark}

\begin{warning} \label{warn-no-fact-structure-on-Loc}
We did not claim $\on{Loc}$ or $\on{Loc}_{\on{co}}$ have any factorization structure. In fact, $\DMod(\on{Bun}_H )_{\mathscr{T}} \boxtimes \mathbf{1}$ is not even a weak factorization category.
\end{warning}

\subsection{Construction of the Localization Functors}
\label{ssec-constr-localization}

We first review the \emph{global infinitesimal Hecke stacks} constructed in \cite{rozenblyum2011connections} (see also \cite{rozenblyum2021connections}), which are global analogue of $\widehat{\on{Hecke}}_{H,K\on{-level}_x}^\bullet$.

\begin{construction}
Let $[\widehat{\on{Hecke}}_{H ,\on{glob}}^{K\on{-level}_x, \bullet}]_{\on{Ran}_{x,\on{un,dR}}}$ be the Cech nerve of the map 
$$\on{Bun}_H^{K\on{-level}_x}\times \on{Ran}_{x,\on{un,dR}}\to (\on{Bun}_H\times\on{Ran}_{x,\on{un,dR}})^{H\on{-level}_{\infty}}/[\widehat{\mathcal{L}H}\times_{H_{x,\dR}} K_{x,\dR}]_{\on{Ran}_{x,\on{un,dR}}}.$$
We have isomorphisms $[(\widehat{\on{Hecke}}_{H ,\on{glob}}^{K\on{-level}_x, \bullet})_\dR]_{\on{Ran}_{x,\on{un,dR}}} \simeq (\on{Bun}_{H}^{K\on{-level}_x})_\dR\times\on{Ran}_{\on{un,dR}}$.

\end{construction}

\begin{remark}  \label{rem-level-structure-glob-hecke}
In the special case when $x:\on{pt}\to \on{Ran}_{\on{un}}$ corresponds to the empty set $\emptyset\subset \on{Maps}(\on{pt},X)$, we obtain $[\widehat{\on{Hecke}}_{H ,\on{glob}}^{\bullet}]_{\on{Ran}_{\on{un,dR}}}$. By definition, we have
$$ [ \widehat{\on{Hecke}}_{H,\on{glob}}^{K\on{-level}_x,\bullet}  ]_{\Ran_{ x,\on{un,dR} }} \simeq [ \widehat{\on{Hecke}}_{H,\on{glob}}^\bullet  ]_{\Ran_{ x,\on{un,dR} }} \times_{ \mathbb{B}H_x^{\times (\bullet+1) }  }  \mathbb{B} K_x^{\times (\bullet+1) }.$$
\end{remark}

\begin{lemma} \label{lem-global-Hecke-unital}
Each $[\widehat{\on{Hecke}}_{H ,\on{glob}}^{K\on{-level}_x, m}]_{ \on{Ran}_{x,\on{un,dR}}}$ is a \emph{unital} lax prestack over $\on{Ran}_{x,\on{un,dR}}$. In particular, the connecting morphisms are strictly unital.
\end{lemma}

\proof By Remark \ref{rem-level-structure-loc-hecke} and Remark \ref{rem-level-structure-glob-hecke}, we can assume $x:\on{pt} \to \on{Ran}_{\on{un}}$ corresponds to the empty set (i.e., we ignore the $K$-level structures).

Let $\phi:I\to J$ be an injective map between finite sets. To prove the claim, we only need to show the correspondence
$$ [\widehat{\on{Hecke}}_{H ,\on{glob}}^{ m}]_J \gets [\widehat{\on{Hecke}}_{H ,\on{glob}}^{ m}]_{\phi} \to [\widehat{\on{Hecke}}_{H ,\on{glob}}^{ m}]_{I}\times_{X^I} X^J $$
(that provides the corr-unital structure) has a degenerate right arm. Explicitly, $[\widehat{\on{Hecke}}_{H ,\on{glob}}^{ m}]_{\phi }(S)$ classifies
\begin{itemize}
		\item[(i)] Maps $x_j:S\to X$ labelled by $J$. 
		\item[]	Let $\Gamma_J\to S\times X$ (resp. $\Gamma_I\to S\times X$) be the schema-theoretic sum of the graphs of all the $x_j$ (resp. of $x_{\phi(i)}$ for $i\in I$). Note that we have a closed embedding $\Gamma_I \to \Gamma_J$.
		\item[(ii)] $H$-torsors $\mathcal{P}_0,\cdots,\mathcal{P}_m$ on $S\times X$ and isomorphisms 
		$$  [\mathcal{P}_0]|_{ (S\times X)\setminus \Gamma_I } \simeq \cdots \simeq [\mathcal{P}_m]|_{ (S\times X)\setminus \Gamma_I } $$
		between their restrictions on $(S\times X)\setminus \Gamma_I $.

		\item[(iii)] Isomorphisms 
		$$  [\mathcal{P}_0]\times_S S^{\on{cl,red}} \simeq \cdots \simeq [\mathcal{P}_m]\times_S S^{\on{cl,red}}, $$
		where we view $[\mathcal{P}_k]\times_S S^{\on{cl,red}}$ as a $G$-torsor on $ S^{\on{cl,red}}\times X $.

		\item[] By restriction, (ii) and (iii) provides two chains of isomorphisms
		$$ [\mathcal{P}_0]|_{ (S\times X)\setminus \Gamma_I  } \times_S S^{\on{cl,red}} \simeq \cdots \simeq [\mathcal{P}_m]|_{(S\times X)\setminus \Gamma_I  } \times_S S^{\on{cl,red}}.$$

		\item[(iv)] An isomorphism between the above two chains.	
\end{itemize}
This makes the claim manifest because (ii)-(iv) do not depend on $J$.

\qed

\begin{construction} It follows from Remark \ref{rem-description-unital-structure-hecke} and the proof of Lemma \ref{lem-global-Hecke-unital} that there are simplicial strictly \emph{co-unital} maps
$$ [ \widehat{\on{Hecke}}_{H,\on{glob}}^{K\on{-level}_x,\bullet}  ]_{\Ran_{ x,\on{un,dR} }} \to [ \widehat{\on{Hecke}}_{H}^{K\on{-level}_x,\bullet}  ]_{\Ran_{ x,\on{un,dR} }} .$$

Also, for any $[m]\to [n]$, the following commutative diagram is Cartesian:
\begin{equation} \label{eqn-global-local-hecke-cart}
\xymatrix{
	[ \widehat{\on{Hecke}}_{H,\on{glob}}^{K\on{-level}_x,n}  ]_{\Ran_{ x,\on{un,dR} }} \ar[r] \ar[d]
	& [ \widehat{\on{Hecke}}_{H,\on{glob}}^{K\on{-level}_x,m}  ]_{\Ran_{ x,\on{un,dR} }} \ar[d] \\
	[ \widehat{\on{Hecke}}_{H}^{K\on{-level}_x,n}  ]_{\Ran_{ x,\on{un,dR} }} \ar[r] 
	& [ \widehat{\on{Hecke}}_{H}^{K\on{-level}_x,m}  ]_{\Ran_{ x,\on{un,dR} }} .
}
\end{equation}
\end{construction}

\begin{notation} For any \emph{quasi-compact} open substack $U\subset \on{Bun}_H$, we write
$$[ \widehat{\on{Hecke}}_{H,U}^{K\on{-level}_x,\bullet}  ]_{\Ran_{ x,\on{un,dR} }}:=[ \widehat{\on{Hecke}}_{H,\on{glob}}^{K\on{-level}_x,\bullet}  ]_{\Ran_{ x,\on{un,dR} }} \times_{\on{Bun}_{H,\dR}} U_\dR.$$
As before, we also use the notations
$$[ \widehat{\on{Hecke}}_{H,U}^{K\on{-level}_x}  ]_{\Ran_{ x,\on{un,dR} }}:= [ \widehat{\on{Hecke}}_{H,U}^{K\on{-level}_x,1}  ]_{\Ran_{ x,\on{un,dR} }},\; [U^{K\on{-level}_x}]_{\Ran_{ x,\on{un,dR} }}:= [ \widehat{\on{Hecke}}_{H,U}^{K\on{-level}_x,0}  ]_{\Ran_{ x,\on{un,dR} }}.$$
\end{notation}

\begin{lemma} $[ \widehat{\on{Hecke}}_{H,U}^{K\on{-level}_x,m}  ]_{\Ran_{ x,\on{un,dR} }}$ is laft and satisfies (the $\Ran_{ x,\on{un,dR} }$-version of) $(\spadesuit)$ in Construction \ref{constr-unit-fact-prestack-indcoh-fact-cat}.
\end{lemma}

\proof Follows from the Cartesian diagram (\ref{eqn-global-local-hecke-cart}), Remark \ref{rem-level-structure-loc-hecke}, Remark \ref{rem-level-structure-glob-hecke} and the fact that $ U \times X^I \to \mathbb{B}[\mathcal{L}^+ H]_I $ is $*$-pullable.

\qed

\begin{lemma} \label{lem-global-map-to-local-hecke-*-pullable}
$[ \widehat{\on{Hecke}}_{H,U}^{K\on{-level}_x,m}  ]_{\Ran_{ x,\on{un,dR} }}\to [ \widehat{\on{Hecke}}_{H}^{K\on{-level}_x,m}  ]_{\Ran_{ x,\on{un,dR} }}$ satisfies (the $\Ran_{ x,\on{un,dR} }$-version of) the assumptions in Construction \ref{constr-functorial-unital-fact-indcoh}(b).
\end{lemma}

\proof Follows from the Cartesian diagram (\ref{eqn-global-local-hecke-cart}), and the fact that $[U^{K\on{-level}_x}]_{\Ran_{ x,\on{un,dR} }}$ is QCA over $\Ran_{ \on{un,dR} }$.

\qed

Recall 
$$  \widetilde{\mathcal{L}\mathfrak{h}}\on{-}mod^{\mathcal{L}^+H\times_{H_x} K_x }_{\mathscr{L}} \simeq \mathscr{L}_{\widehat{\on{Hecke}}_H^{K\on{-level}_x} } \on{-mod}( \mathcal{R}ep( \mathcal{L}^+ H\times_{H_x} K_x )),$$
where $\mathscr{L}_{\widehat{\on{Hecke}}_H^{K\on{-level}_x} } $ is the line bundle in $\mathcal{I}nd\mathcal{C}oh^{!,\ren}(\widehat{\on{Hecke}}_{H}^{K\on{-level}_x})$ defined as the pullback of $\mathscr{L}_{ \widehat{\on{Hecke}}_H }$. Recall we also have
$$  \widetilde{\mathcal{L}\mathfrak{h}}\on{-}mod^{\mathcal{L}^+H\times_{H_x} K_x }_{\mathscr{L}} \simeq \lim \mathcal{I}nd\mathcal{C}oh^{!,\ren}(   \widehat{\on{Hecke}}_{H}^{K\on{-level}_x,\bullet}  )_{\mathscr{L},\Delta} \simeq \colim \mathcal{I}nd\mathcal{C}oh^{!,\ren}(   \widehat{\on{Hecke}}_{H}^{K\on{-level}_x,\bullet}  )_{\mathscr{L},\Delta^\op}, $$
where the limit is taken in $[\on{CrysCat}(\on{Ran}_{x,\on{un,dR}})]_{\on{lax-unital}}$, while the colimit is taken in $\on{CrysCat}(\on{Ran}_{x,\on{un,dR}})$. Also these diagrams are ``twisted'' by $\mathscr{L}_{\widehat{\on{Hecke}}_H^{K\on{-level}_x} }$ in a way similar to Remark \ref{rem-KL-twisted-as-totalization}.

\begin{construction} \label{const-local-to-global-hecke}
For any quasi-compact open substack $U\subset \on{Bun}_H$, the Cartesian squares (\ref{eqn-global-local-hecke-cart}) provide an $ [\mathcal{I}nd\mathcal{C}oh^{!,\ren}(\widehat{\on{Hecke}}_{H}^{K\on{-level}_x}),\star]$-action on $ \mathcal{I}nd\mathcal{C}oh( U^{K\on{-level}_x})$ (viewed as objects in $[\on{CrysCat}(\on{Ran}_{x,\on{un,dR}})]_{\on{lax-unital}}$). Consider\footnote{Since $\widehat{\on{Hecke}}_{H,U}^{K\on{-level}_x,\bullet}$ is laft, we do not need to use the renormalized $\IndCoh$-theories.}
$$  \mathscr{L}_{\widehat{\on{Hecke}}_H^{K\on{-level}_x} } \on{-mod}( \mathcal{I}nd\mathcal{C}oh( U^{K\on{-level}_x})  ).$$
Similar to Remark \ref{rem-KL-twisted-as-totalization}, this crystal of categories can be written as a geometric realization or totalization as follows. Consider the functor
$$\mathcal{I}nd\mathcal{C}oh(   \widehat{\on{Hecke}}_{H,U}^{K\on{-level}_x,\bullet}  )_{\mathscr{L},\Delta}: \Delta \to [\on{CrysCat}(\on{Ran}_{x,\on{un,dR}})]_{\on{lax-unital}} $$
that sends $[m]$ to $\mathcal{I}nd\mathcal{C}oh(   \widehat{\on{Hecke}}_{H,U}^{K\on{-level}_x,m}  )$, and sends a morphism $f:[m]\to [n]$ in $\Delta$ to the $!$-pullback functor twisted by the line bundle on $ \widehat{\on{Hecke}}_{H,U}^{K\on{-level}_x}$ that is the pullback of $\mathscr{L}_{ \widehat{\on{Hecke}} }$ along $ \widehat{\on{Hecke}}_{H,U}^{K\on{-level}_x} \to \widehat{\on{Hecke}}_{H}$. Passing to left adjoints, we also have a functor
$$\mathcal{I}nd\mathcal{C}oh(   \widehat{\on{Hecke}}_{H,U}^{K\on{-level}_x,\bullet}  )_{\mathscr{L},\Delta^\op}: \Delta^\op \to \on{CrysCat}(\on{Ran}_{x,\on{un,dR}}) $$
connected by twisted pushforward functors (which are simultaneously $!$-pushforward and $*$-pushforward functors). Then we have
$$
\xymatrix{
	\widetilde{\mathcal{L}\mathfrak{h}}\on{-}mod^{\mathcal{L}^+H\times_{H_x} K_x }_{\mathscr{L}} \ar[r] \ar[d]^-\simeq &
	\mathscr{L}_{\widehat{\on{Hecke}}_H^{K\on{-level}_x} } \on{-mod}( \mathcal{I}nd\mathcal{C}oh( U^{K\on{-level}_x})  )\ar[d]^-\simeq \\
	\colim \mathcal{I}nd\mathcal{C}oh^{!,\ren}(   \widehat{\on{Hecke}}_{H}^{K\on{-level}_x,\bullet}  )_{\mathscr{L},\Delta^\op} \ar[r]^-{!\on{-pull}}\ar[d]^-\simeq
	&
	\colim \mathcal{I}nd\mathcal{C}oh(   \widehat{\on{Hecke}}_{H,U}^{K\on{-level}_x,\bullet}  )_{\mathscr{L},\Delta^\op}\ar[d]^-\simeq \\
	\lim \mathcal{I}nd\mathcal{C}oh^{!,\ren}(   \widehat{\on{Hecke}}_{H}^{K\on{-level}_x,\bullet}  )_{\mathscr{L},\Delta} \ar[r]^-{!\on{-pull}}
	&
	\lim \mathcal{I}nd\mathcal{C}oh(   \widehat{\on{Hecke}}_{H,U}^{K\on{-level}_x,\bullet}  )_{\mathscr{L},\Delta},
}
$$
where the limit is taken in $[\on{CrysCat}(\on{Ran}_{x,\on{un,dR}})]_{\on{lax-unital}}$, while the colimit is taken in $\on{CrysCat}(\on{Ran}_{x,\on{un,dR}})$. Note that the second horizontal morphism is well-defined because of the Cartesian diagram (\ref{eqn-global-local-hecke-cart}). Since the horizontal morphisms are induced by $!$-pullback functors along \emph{co-unital} maps, they are unital morphisms between crystals of categories.

Similarly, we have 
\begin{equation}
\label{eq:coloc-construction}
\xymatrix{
	\widetilde{\mathcal{L}\mathfrak{h}}\on{-}mod^{\mathcal{L}^+H\times_{H_x} K_x }_{\mathscr{L},\on{co}}  \ar[d]^-\simeq &
	\mathscr{L}_{\widehat{\on{Hecke}}_H^{K\on{-level}_x} } \on{-mod}( \mathcal{I}nd\mathcal{C}oh( U^{K\on{-level}_x})  )\ar[d]^-\simeq \ar[l] \\
	\colim \mathcal{I}nd\mathcal{C}oh_{*,\ren}(   \widehat{\on{Hecke}}_{H}^{K\on{-level}_x,\bullet}  )_{\mathscr{L},\Delta^\op} \ar[d]^-\simeq
	&
	\colim \mathcal{I}nd\mathcal{C}oh(   \widehat{\on{Hecke}}_{H,U}^{K\on{-level}_x,\bullet}  )_{\mathscr{L},\Delta^\op}\ar[d]^-\simeq  \ar[l]_-{*\on{-push}}\\
	\lim \mathcal{I}nd\mathcal{C}oh_{*,\ren}(   \widehat{\on{Hecke}}_{H}^{K\on{-level}_x,\bullet}  )_{\mathscr{L},\Delta} 
	&
	\lim \mathcal{I}nd\mathcal{C}oh(   \widehat{\on{Hecke}}_{H,U}^{K\on{-level}_x,\bullet}  )_{\mathscr{L},\Delta} \ar[l]_-{*\on{-push}}.
}
\end{equation}
By Lemma \ref{lem-global-map-to-local-hecke-*-pullable}, we can pass to left adjoints along the horizontal direction and obtain a unital morphism
$$  \widetilde{\mathcal{L}\mathfrak{h}}\on{-}mod^{\mathcal{L}^+H\times_{H_x} K_x }_{\mathscr{L},\on{co}} \to \mathscr{L}_{\widehat{\on{Hecke}}_H^{K\on{-level}_x} } \on{-mod}( \mathcal{I}nd\mathcal{C}oh( U^{K\on{-level}_x})  ).$$

\end{construction}

\begin{construction} \label{constr-Hecke-indcoh-vs-dmod}
We have a unital morphism
\begin{equation*}
\begin{aligned}
\on{ind}:\mathscr{L}_{\widehat{\on{Hecke}}_H^{K\on{-level}_x} } \on{-mod}( \mathcal{I}nd\mathcal{C}oh( U^{K\on{-level}_x})  )\simeq  \colim \mathcal{I}nd\mathcal{C}oh(   \widehat{\on{Hecke}}_{H,U}^{K\on{-level}_x,\bullet}  )_{\mathscr{L},\Delta^\op} \to\\
\to  \IndCoh( (U^{K\on{-level}_x})_\dR)_{\mathscr{T}}\boxtimes \mathbf{1} \simeq \DMod( U^{K\on{-level}_x} )_{\mathscr{T}}\boxtimes \mathbf{1} 
\end{aligned}
\end{equation*}
induced by twisted pushforward functors. Here we used the fact that $\mathscr{L}$ and $\mathscr{T}$ are compatible (see Remark \ref{rem-compatible-pair-two-gerbe}). Namely, by restriction along the augementation $ \widehat{\on{Hecke}}_{H,U}^{K\on{-level}_x,\bullet}\to (U^{K\on{-level}_x})_\dR$, the twisting $\mathscr{T}$ induces a cosimplicial $\mathbb{G}_m$-gerbe on $ \widehat{\on{Hecke}}_{H,U}^{K\on{-level}_x,\bullet}$ equipped with a trivialization of its 0-term. Then the compatibility statement implies $ \mathcal{I}nd\mathcal{C}oh(   \widehat{\on{Hecke}}_{H,U}^{K\on{-level}_x,\bullet}  )_{\mathscr{L},\Delta^\op} $ is the simplicial diagram twisted by this gerbe.

Passing to right adjoints, we obtain a lax-unital morphism
$$\on{oblv}: \DMod( U^{K\on{-level}_x} )_{\mathscr{T}}\boxtimes \mathbf{1} \to \mathscr{L}_{\widehat{\on{Hecke}}_H^{K\on{-level}_x} } \on{-mod}( \mathcal{I}nd\mathcal{C}oh( U^{K\on{-level}_x})  )$$
induced by twisted $!$-pullback functors.
\end{construction}

\begin{definition} Recall we have $\DMod(\on{Bun}_H^{K\on{-level}_x })_{\mathscr{T}} \simeq \lim_{U} \DMod( U^{K\on{-level}_x} )_{\mathscr{T}}$, where the connecting functors are given by pullbacks. Consider the following compositions of unital morphisms:
\begin{eqnarray*}
 \on{Loc}_U: \widetilde{\mathcal{L}\mathfrak{h}}\on{-}mod^{\mathcal{L}^+H\times_{H_x} K_x }_{\mathscr{L}} \xrightarrow{!\on{-pull}} \mathscr{L}_{\widehat{\on{Hecke}}_H^{K\on{-level}_x} } \on{-mod}( \mathcal{I}nd\mathcal{C}oh( U^{K\on{-level}_x})  ) \xrightarrow{\on{push}} \DMod( U^{K\on{-level}_x} )_{\mathscr{T}}\boxtimes \mathbf{1}  \\
  \on{Loc}_{U,\on{co}}: \widetilde{\mathcal{L}\mathfrak{h}}\on{-}mod^{\mathcal{L}^+H\times_{H_x} K_x }_{\mathscr{L},\on{co}} \xrightarrow{*\on{-pull}} \mathscr{L}_{\widehat{\on{Hecke}}_H^{K\on{-level}_x} } \on{-mod}( \mathcal{I}nd\mathcal{C}oh( U^{K\on{-level}_x})  ) \xrightarrow{\on{push}} \DMod( U^{K\on{-level}_x} )_{\mathscr{T}}\boxtimes \mathbf{1}.
 \end{eqnarray*}
 They are functorial in $U$, i.e., compatible with the pullback functors $\DMod( U_1^{K\on{-level}_x} )\to \DMod( U_2^{K\on{-level}_x} )$. Passing to limits, we obtain unital morphisms
 \begin{eqnarray*}
 \on{Loc}: \widetilde{\mathcal{L}\mathfrak{h}}\on{-}mod^{\mathcal{L}^+H\times_{H_x} K_x }_{\mathscr{L}} \to \DMod( \on{Bun}_H^{K\on{-level}_x} )_{\mathscr{T}}\boxtimes \mathbf{1}  \\
  \on{Loc}_{\on{co}}: \widetilde{\mathcal{L}\mathfrak{h}}\on{-}mod^{\mathcal{L}^+H\times_{H_x} K_x }_{\mathscr{L},\on{co}} \to \DMod( \on{Bun}_H^{K\on{-level}_x} )_{\mathscr{T}}\boxtimes \mathbf{1}.
 \end{eqnarray*}
\end{definition}

\subsection{Compatibility with Chiral Homology}
Since $\on{Loc}$ is unital, we have a commutative square
$$
	\xymatrix{
	[ \widetilde{\mathcal{L}\mathfrak{h}}\on{-}mod^{\mathcal{L}^+H\times_{H_x} K_x }_{\mathscr{L}}  ]_x  \ar[rr]^-{ [\on{Loc}]_{ x} }  \ar[d]^-{\on{ins}_{\on{unit}}}
	& & \DMod(\on{Bun}_H ^{K\on{-level}_x})_{\mathscr{T}}
	 \ar[d]^-{\on{ins}_{\on{unit}}}\\
	[ \widetilde{\mathcal{L}\mathfrak{h}}\on{-}mod^{\mathcal{L}^+H\times_{H_x} K_x }_{\mathscr{L}}  ]_{\on{Ran}_{x,\on{un,dR}}}  \ar[rr]^-{ [\on{Loc}]_{ \on{Ran}_{x,\on{un,dR}}  } } 
	& & \DMod(\on{Bun}_H ^{K\on{-level}_x})_{\mathscr{T}} \otimes \DMod(\on{Ran}_{x,\on{un}})
	}.
	$$
Note that the right vertical functor sends $\mathcal{F}$ to $\mathcal{F}\boxtimes \omega_{ \on{Ran}_{x,\on{un}} }$. Recall $\on{Ran}_{x,\on{un}} $ is homologically contractible (see \cite[Theorem 4.1.8, 4.6.2]{gaitsgory2015atiyah}). Hence $\on{C}^{\on{ch}}_{x,\on{un}}( \omega_{ \on{Ran}_{x,\on{un}} } ) \simeq \omega_{\on{pt}}$. This proves Theorem \ref{thm-localization-KL}(1).

\subsection{Compatibility with Tate Twist}
Unwinding the definitions, we only need to construct the following commutative diagram for cosimplicial crystals of categories:
$$
\xymatrix{
\mathcal{I}nd\mathcal{C}oh_{*,\ren}(\widehat{\on{Hecke}}_{H}^{\bullet}  )_{\mathscr{L}\otimes \mathscr{L}^{\on{Tate}},\Delta}
	\ar[r]^-{*\on{-pull}}   &
	\mathcal{I}nd\mathcal{C}oh(\widehat{\on{Hecke}}_{H,U}^{\bullet}  )_{\mathscr{L}\otimes \mathscr{L}^{\on{Tate}},\Delta} 
\\
\mathcal{I}nd\mathcal{C}oh^{!,\ren}(\widehat{\on{Hecke}}_{H}^{\bullet}  )_{\mathscr{L},\Delta}  \ar[r]^-{!\on{-pull}} \ar[u]^-\simeq_-{\theta_0}  &
	\mathcal{I}nd\mathcal{C}oh(\widehat{\on{Hecke}}_{H,U}^{\bullet}  )_{\mathscr{L},\Delta}   \ar[u]^-\simeq_-{\vartheta_0}
}
$$
such that
\begin{itemize}
	\item The horizontal morphisms were defined in Construction \ref{const-local-to-global-hecke}.
	\item The left vertical equivalences were defined in Step 3 of the proof of Theorem \ref{thm-tate-twist}. Namely, for each $[m]$, we use the unique $\mathcal{I}nd\mathcal{C}oh^{!,\ren}(\widehat{\on{Hecke}}_{H}^{m}  )$-linear equivalence $\theta_0^{[m]}$ such that
 	$$\mathcal{I}nd\mathcal{C}oh(\on{pt})\xrightarrow{!\on{-pull}}  \mathcal{I}nd\mathcal{C}oh^{!,\ren}(\widehat{\on{Hecke}}_{H}^{0}  ) \xrightarrow{\theta_0^{[m]}}  \mathcal{I}nd\mathcal{C}oh_{*,\ren}(\widehat{\on{Hecke}}_{H}^{m}  ) $$
 	is equivalent to
 	$$ \mathcal{I}nd\mathcal{C}oh(\on{pt})\xrightarrow{*\on{-pull}}  \mathcal{I}nd\mathcal{C}oh_{*,\ren}(\widehat{\on{Hecke}}_{H}^{0}  ) \xrightarrow{(p_0)_\ren^? }  \mathcal{I}nd\mathcal{C}oh_{*,\ren}(\widehat{\on{Hecke}}_{H}^{m}  ) .$$

 	\item For each $[m]$, the right vertical equivalence $\vartheta_0^{[m]}$ is tensoring up with
 	$ p_0^!\circ \Upsilon ( K^{-1}  )$.
\end{itemize}

We only need to show $\vartheta_0^{[m]}$ is the \emph{unique} $\mathcal{I}nd\mathcal{C}oh^{!,\ren}(\widehat{\on{Hecke}}_{H}^{m}  )$-linear equivalence such that
$$
\xymatrix{
\mathcal{I}nd\mathcal{C}oh_{*,\ren}(\widehat{\on{Hecke}}_{H}^{m}  )
	\ar[r]^-{*\on{-pull}}   &
	\mathcal{I}nd\mathcal{C}oh(\widehat{\on{Hecke}}_{H,U}^{m}  )
\\
\mathcal{I}nd\mathcal{C}oh^{!,\ren}(\widehat{\on{Hecke}}_{H}^{m}  )  \ar[r]^-{!\on{-pull}} \ar[u]^-\simeq_-{\theta_0}  &
	\mathcal{I}nd\mathcal{C}oh(\widehat{\on{Hecke}}_{H,U}^{m}  )   \ar[u]^-\simeq_-{\vartheta_0^{[m]}}
}
$$
commutes. Using the base-change isomorphisms, we only need to prove the case $m=0$, which then follows from the fact that $\vartheta_0^{[0]} =  \Upsilon ( K^{-1}  )\otimes -$ is the unique $\mathcal{I}nd\mathcal{C}oh(U )$-linear equivalence such that
$$\mathcal{I}nd\mathcal{C}oh(\on{pt})\xrightarrow{!\on{-pull}}  \mathcal{I}nd\mathcal{C}oh(U   ) \xrightarrow{\vartheta_0^{[0]}} \mathcal{I}nd\mathcal{C}oh(U  )$$
 	is equivalent to the $*$-pullback functor.

\subsection{Compatibility with Propogation}
Our strategy for the proof of Theorem \ref{thm-localization-KL}(3-4) is to reduce to the case when the factorization algebra $A$ (resp. $A'$) is the unit factorization algebra. To achieve this reduction, we need to use the tool of \emph{propogation of factorization modules}, and prove $\on{Loc}$ and $\on{Loc}_{\on{co}}$ are compatible with propogations. We first explain what these propogations mean.

\begin{remark}
The propagation method was initially introduced by Beilinson and Drinfeld in \cite[4.4.3 and 4.4.9]{beilinson2004chiral}. We emphasize that this technique \emph{only exists in the unital setting}.
\end{remark}

\begin{construction}
Consider the laft lax prestack
\[\on{Arrow}( \on{Ran}_{x,\on{un}} ) := S \mapsto \tx{Fun}([2], \Ran_{x, \un}(S)).\]
We have two projections
$$p_s,p_t: \on{Arrow}( \on{Ran}_{x,\on{un}} )\to \on{Ran}_{x,\on{un}} $$
``remembering'' respectively the source and the target. Note that we have a $2$-morphism $\mathfrak{i}:p_s\to p_t$. We equip $\on{Arrow}( \on{Ran}_{x,\on{un}} )$ with the factorization module space structure (for the factorization space $\on{Ran}_{\on{un}}$) such that $p_t$ is compatible with the factorization module structures but $p_s$ is not. For a map $y:\on{pt}\to \on{Ran}_{\on{un}}$ corresponding to a finite subset $J$ with $I_t\cap J = \emptyset$, the map 
$$(\on{Ran}_{\on{un}} \times \on{Arrow}( \on{Ran}_{x,\on{un}} ))_{\on{disj}} \to \on{Arrow}( \on{Ran}_{x,\on{un}} )$$
sends $(y,z)$ to the point corresponding to $I_0\subset I_s \subset J\sqcup I_t$.
\end{construction}

\begin{remark}
Explicitly, a map $z:\on{pt}\to \on{Arrow}( \on{Ran}_{x,\on{un}} )$ corresponds to a chain $I_0\subset I_s\subset I_t$ of finite subsets of $\on{Maps}(\on{pt},X)$ such that $I_0\subset I_s$ (resp. $I_0\subset I_t$) corresponds to the composition $\on{pt}\to \on{Arrow}( \on{Ran}_{x,\on{un}} ) \xrightarrow{p_s} \on{Ran}_{x,\on{un}}$ (resp. $\on{pt}\to \on{Arrow}( \on{Ran}_{x,\on{un}} ) \xrightarrow{p_t} \on{Ran}_{x,\on{un}}$).
\end{remark}

\begin{construction}[Propogation of factorization modules] 
\label{constr-propogation-fact-mod}
Let $V\to W$ be a morphism between laft prestacks. For any unital factorization category $\mathcal{A}$ over $W\times \Ran_{\on{un,dR}}$ (or equivalently, a $W$-family of unital factorization categories over $\Ran_{\on{un,dR}}$), and for any factorization $\mathcal{A}$-module $\mathcal{C}$ over $V\times \on{Ran}_{x,\on{un,dR}}$, the (co)restriction $\mathbf{cores}_{p_{t,\dR}}(\mathcal{C})$, which is a crystal of categories over $V\times  \on{Arrow}( \on{Ran}_{x,\on{un}} )_\dR$, inherits a factorization $\mathcal{A}$-module structure. Moreover, for any unital factorization algebra $A$ in $\mathcal{A}$, there is a canonical functor
$$ \on{prop}_A: A\on{-FactMod}_{\on{un}}( \mathcal{C} ) \to A\on{-FactMod}_{\on{un}}( \mathbf{cores}_{p_{t,\dR}}(\mathcal{C}) ) $$
such that the following diagram commutes
$$
\xymatrix{
	A\on{-FactMod}_{\on{un}}( \mathcal{C} ) \ar[d]^-{\on{prop}_A} \ar[r]^-{\oblv}
	& \mathbf{\Gamma}( \on{Id}_\dR, \mathcal{C} ) \ar[d]^-{\on{res}} \\
	A\on{-FactMod}_{\on{un}}( \mathbf{cores}_{p_{t,\dR}}(\mathcal{C}) ) \ar[r]^-{\oblv}
	&
	\mathbf{\Gamma}( p_{t,\dR}, \mathcal{C} ),
}
$$
where the bottom functor is the restriction functor for sections of $\mathcal{C}$. More generally, for any morphism $f: A \to B$ between unital factorization algebras, we can use restriction of factorization modules (as defined in Proposition \ref{prop:restriction-functor}) to define the composition
$$ \on{prop}_{A}: B\on{-FactMod}_{\on{un}}( \mathcal{C} ) \xrightarrow{\on{prop}_B} B\on{-FactMod}_{\on{un}}( \mathbf{cores}_{p_{t,\dR}}(\mathcal{C}) ) \xrightarrow{\on{res}_f}  A\on{-FactMod}_{\on{un}}( \mathbf{cores}_{p_{t,\dR}}(\mathcal{C}) ) $$
and have a natural transformation
\begin{equation} \label{eqn-nt-prop-res-oblv}
\xymatrix{
	B\on{-FactMod}_{\on{un}}( \mathcal{C} ) \ar[d]^-{\on{prop}_{A}} \ar[r]^-{\oblv}
	& \mathbf{\Gamma}( \on{Id}_\dR, \mathcal{C} ) \ar[d]^-{\on{res}} \\
	A\on{-FactMod}_{\on{un}}( \mathbf{cores}_{p_{t,\dR}}(\mathcal{C}) )\ar@{=>}[ru] \ar[r]^-{\oblv}
	&
	\mathbf{\Gamma}( p_{t,\dR}, \mathcal{C} ),
}
\end{equation}
\end{construction}

\begin{construction} \label{const-prop-unit=insert-unit}
In the special case when $A=\on{unit}$ is the unit factorization algebra of $\mathcal{A}$, we have the following commutative diagram
$$
\xymatrix{
	B\on{-FactMod}_{\on{un}}( \mathcal{C} ) \ar[d]^-{\on{prop}_{\on{unit}}} \ar[r]^-{\oblv}
	& \mathbf{\Gamma}( \on{Id}_\dR, \mathcal{C} ) \ar[d]^-{\on{ins}_{\on{unit}}}\\
	\on{unit}\on{-FactMod}_{\on{un}}( \mathbf{cores}_{p_{t,\dR}}(\mathcal{C}) ) \ar[r]^-{\oblv}
	 &
	\mathbf{\Gamma}( p_{t,\dR}, \mathcal{C} ),
}
$$
where $\on{ins}_{\on{unit}}$ is the functor of \emph{unit insertion}
$$ \on{ins}_{\on{unit}}: \mathbf{\Gamma}( \on{Id}_\dR, \mathcal{C} )\xrightarrow{\on{res}}\mathbf{\Gamma}( p_{s,\dR}, \mathcal{C} ) \xrightarrow{\mathfrak{i}_!}  \mathbf{\Gamma}( p_{t,\dR}, \mathcal{C} ).$$
\end{construction}

\paragraph{Digression: what is propagation?}
Let us explain the name. Suppose $A$ is a unital factorization algebra and $N$ is its, such that its fiber at a $\mb{C}$-point
\[y := [\{x\} \subseteq \{x, x_1, \ldots, x_n\}] \in \on{Arrow}( \on{Ran}_{x,\on{un}} )(\mb{C})\]
is given by some vector space $N_0 \tensor A_0^{\tensor (n - 1)}$. The map $y$ induces a map
\[\left(\Ran_{\un, y} := S \mapsto \{\alpha: S \to \Ran_{\un} \mid (x, x_1, \ldots, x_n) \subseteq \alpha\}\right) \to \on{Arrow}( \on{Ran}_{x,\on{un}} )\]
Then the pullback of $\tx{prop}_A(N)$ along this map produces a unital $A$-factorization module on the former, which should be thought of as consisting of $(N_0, A_0, \ldots, A_0)$ supported at $(x, x_1, \ldots, x_n)$ respectively. Thus, it can be seen as ``propagating'' the algebra $A$ to other points on the curve.

Another perspective is as follows. As explained in \cite{raskin2015chiral}, via \'etale hyperdescent the category of $A$-unital factorization modules in $C$ \emph{itself} form a weak unital factorization category under \emph{chiral fusion}, and the underlying sheaf of $\tx{prop}_A(N)$ can be seen as observing $A$ as an element $A^\tx{enh} \in \tx{FactAlg}_{\un}(A\tx{-FactMod}_\un(\mc{C}))$ and $N$ as an element $N^\tx{enh} \in A^\tx{enh}\tx{-FactMod}_\un(A\tx{-FactMod}_\un(\mc{C}))$.
This is, of course, evocative of the ``enhancement'' equivalences
\[\tx{AlgEnh}: \mc{O}\tx{-Alg}(\mc{C})_{A /} \simeq \mc{O}\tx{-Alg}(A\mod^{\mc{O}}(\mc{C}))\]
\[\tx{ModEnh}: A\mod^{\mc{O}}(\mc{C}) \simeq \tx{AlgEnh}(A)\mod^{\mc{O}}(A\mod^{\mc{O}}(\mc{C}))\]
which holds for any coherent $\infty$-operad $\mc{O}$, any (weakly) $\mc{O}$-monoidal category $\mc{C}$ and any $\mc{O}$-algebra $A \in \mc{C}$, as established in \cite[3.4.1.7]{HA} and \cite[3.4.1.9]{HA}.

\begin{construction} In the special case when $\mathcal{A}$ and $\mathcal{C}$ are the constant crystals, i.e., $\mathcal{A} = \IndCoh(W)\boxtimes \mathbf{1}$ and $\mathcal{C}=\IndCoh(V)\boxtimes \mathbf{1}$, the functor
$$ \on{res}: \mathbf{\Gamma}( \on{Id}_\dR, \mathcal{C} )  \to \mathbf{\Gamma}( \on{p}_{t,\dR}, \mathcal{C} ) $$
is just
$$p_t^!: \IndCoh(V)\otimes \DMod( \Ran_{x,\on{un}} ) \to \IndCoh(V)\otimes \DMod(\on{Arrow}( \Ran_{x,\on{un}}) ).$$
Hence we obtain a natural transformation
\begin{equation} \label{eqn-prop-A'-nt}
\xymatrix{
	B\on{-FactMod}_{\on{un}}( \IndCoh(V) \boxtimes \mathbf{1} ) \ar[d]^-{\on{prop}_{A}} \ar[r]^-{\oblv}
	& \IndCoh(V)\otimes \DMod( \Ran_{x,\on{un}} )  \\
	A\on{-FactMod}_{\on{un}}( \IndCoh(V) \boxtimes\mathbf{cores}_{p_{t,\dR}}(\mathbf{1}) )\ar@{=>}[ru] \ar[r]^-{\oblv}
	&
	\IndCoh(V)\otimes \DMod(\on{Arrow}( \Ran_{x,\on{un}}) ) \ar[u]^-{p_{t,!}},
}
\end{equation}
\end{construction}

The diagram (\ref{eqn-prop-A'-nt}) does not commute in general. But we have the following results:

\begin{lemma}[Propogation-restriction Lemma for Chiral Homology] \label{lem-AB-lemma} 
The diagram (\ref{eqn-prop-A'-nt}) commutes after taking chiral homology, i.e., after composing with
$$\on{C}_{x,\on{un}}^{\on{ch}}: \IndCoh(V)\otimes \DMod( \Ran_{x,\on{un}} ) \to \IndCoh(V).$$
\end{lemma}

Recall $\on{C}_{x,\on{un}}^{\on{ch}}$ is given by $!$-pushforward along $p:  \Ran_{x,\on{un}}\to \on{pt}$. In the case when $B = A$, the following stronger result is true:

\begin{lemma}[Propogation Lemma for Chiral Homology] \label{lem-AA-lemma} The composition
$$ p_! \circ p_{s,!}\circ p_t^!  \simeq p_! \circ p_{t,!}\circ p_t^! \to p_!  $$
in $\on{Funct}( \DMod( \Ran_{x,\on{un}} ) ) \to \on{Vect}$ is an equivalence. Moreover, the natural transformation
$$ p_{s,!}\circ p_t^! \to p^!\circ p_! $$
is an equivalence.
\end{lemma}

\begin{remark}
Before we prove, let us explain the names.
For a set of unital $A$-factorization modules $N_1, \ldots, N_k$ inserted at distinct points $x_1, \ldots, x_k \in X(\mb{C})$, let us write the chiral homology as
\[\langle N_1, \ldots, N_k \rangle^{A}_{x_1, \ldots, x_k};\]
Given a unital factorization $A$-module $M$ supported at $x$, the equality between fibers of $p_{s,!}\circ p_t^!(M)$ and $p^!\circ p_!(M)$ at some $(x, x_1, \ldots, x_n)$, given by the propagation lemma for chiral homology, translates to
\[\langle M, A, \ldots, A \rangle^{A}_{x, x_1, \ldots, x_n} \simeq \langle M \rangle^A_{x},\]
which says that when computing $A$-chiral homology, extra copies of $A$ would not affect the result. The propagation-restriction lemma, on the other hand, roughly states that if $M$ is a $B$-factorization module and $A \to B$ is a map between unital factorization algebras, then we have
\[\left\langle\langle M \rangle_x^A, \langle M, B \rangle_{x, x_1}^A, \langle M, B, B \rangle_{x, x_1, x_2}^A, \ldots \right\rangle \simeq \langle M \rangle^B_x\]
where the LHS now ranges over all of $\Ran_{x, \un}$.
\end{remark}

\begin{proof}[Proof for Lemma \ref{lem-AA-lemma}]
Consider the pullback diagram
	\[\xymatrix{
	\Ran_{x, \un} \times_{\{x\}} \Ran_{x, \un} \ar^{q_2}[r] \ar^{q_1}[d] & \Ran_{x, \un} \ar^{p}[d] \\
	\Ran_{x, \un} \ar^{p}[r] & \{x\}
	}\]
	We claim that $(q_1)_!$ is well-defined, and the natural transformation $(q_1)_! \circ (q_2)^! \Rightarrow p^! \circ p_!$ is an isomorphism. Indeed, this can be checked by by first establishing the corresponding facts for the following diagram
	\[\xymatrix{
	\Ran_{x, \un} \times_{R} \Ran_{x} \ar[r] \ar[d] & \Ran_{x} \ar[d] \\
	\Ran_{x, \un} \ar[r] & \{x\}
	}\]
	using the fact that $\Ran_x$ is pseudo-proper, then use the universal homological left cofinality of $\Ran_{x} \to \Ran_{x, \un}$, provable as in \cite[4.6.2]{gaitsgory2015atiyah}. There is a map of lax prestacks
	\[\iota: \tx{Arrow}(\Ran_{x, \un}) \to \Ran_{x, \un} \times_{\{x\}} \Ran_{x, \un} \hspace{1em} (I_0 \subseteq I_s \subseteq I_t) \mapsto ((I_0 \subseteq I_s), (I_0 \subseteq I_t))\]
	\[\tx{such that} \hspace{1em} q_2 \circ \iota \simeq p_t \hspace{1em} \tx{and} \hspace{1em} q_1 \circ \iota \simeq p_s.\]
	It suffices to establish that $\iota$ is universally homologically left cofinal. As in \cite[4.6.5]{gaitsgory2015atiyah}, it suffices to check that for every affine $S$, every $S$-point
	\[((I_0 \subseteq I_1), (I_0 \subseteq I_2)) \in \Ran_{x, \un} \times_{\{x\}} \Ran_{x, \un}(S),\]
	and every $t: T \to S$ where $T$ is affine, the category
	\[\{J \in \Ran(T), K \in \Ran(T), \alpha: I_0|_t \subseteq J \subseteq K~\tx{such that}~I_1|_t \subseteq J, I_2|_t \subseteq K\}\]
	is a contractible category; but this category has an initial element given by $(I_1|_t, I_1|_t \cup I_2|_t)$.
\end{proof}

\begin{proof}[Proof for Lemma \ref{lem-AB-lemma}]
We first prove the special case $A \simeq \tx{unit}_\mc{A}$, for $\eta_B: \tx{unit}_\mc{A} \to B$ the canonical map attached to any factorization algebra $B$. Let $\iota_s: \Ran_{x, \un} \to \tx{Arrow}(\Ran_{x, \un})$ be the inclusion map $(I_0 \subseteq I_1) \mapsto (I_0 \subseteq I_1 \subseteq I_1)$. It follows from Proposition \ref{prop:restriction-functor} and Proposition \ref{prop:unital-basic-properties} that
\[\tx{res}_{\eta_B}: B\tx{-FactMod}_\un(\mc{C}) \to \tx{unit}_\mc{A}\tx{-FactMod}_\un(\mc{C})\]
is isomorphic to the functor on $\mc{C}$ induced by
	\[p_s^! \circ \iota_{s}^!: \DMod(\tx{Arrow}(\Ran_{x, \un})) \to \DMod(\tx{Arrow}(\Ran_{x, \un})),\]
	So it suffices to establish $p_! \circ (p_s)_! \circ p_s^! \circ \iota_s^! \circ p_t^! \simeq p_! \circ (p_s)_! \circ p_s^! \simeq p_!$, which follows from the homological contractibility of $p_s$.
	
	Each $f: A \to B$ comes with a natural map $\eta_B \to f$, which by functoriality of $\tx{res}_f$ induces a natural transformation
	\[\alpha: (\tx{id} \times p_!) \circ \tx{prop}_{\tx{unit}_\mc{A}} \Rightarrow (\tx{id} \times p_!) \circ \tx{prop}_{A};\]
	unfolding the construction of $\tx{res}_f$ we see that $\alpha$ is induced by the following natural transformation of functors $\DMod(\tx{Arrow}(\Ran_{x, \un})) \to \Vect$:
	\[\alpha_{\tx{univ}}: p_! \circ (p_s)_! \circ (i \circ p_s)^! \Rightarrow p_! \circ (p_s)_!\]
	which is in turn induced the natural transformation $i \circ p_s \Rightarrow \tx{id}$ (here $i: \Ran_{x, \un} \to \tx{Arrow}(\Ran_{x, \un})$ is the map $I_0 \subseteq I_1 \mapsto I_0 \subseteq I_1 \subseteq I_1$). Thus it suffices to show that $\alpha_{\tx{univ}}$ is an equivalence. We get this by observing that $i \circ p_s$ is universally homologically left cofinal, which follows from the same argument as in the previous lemma.
\end{proof}

\begin{remark} \label{rem-RanRanRan}
One can also define $2\on{Chain}( \on{Ran}_{x,\on{un}} )$ classifying $I_0\subset I_1\subset I_2 \subset I_3$ and propogation functor
$$\on{prop}_A: A\on{-FactMod}_{\on{un}}( \mathbf{cores}_{p_{t,\dR}}(\mathcal{C}) ) \to A\on{-FactMod}_{\on{un}}( \mathbf{cores}_{p_{3,\dR}}(\mathcal{C}) ),$$
where $p_3: 2\on{Chain}( \on{Ran}_{x,\on{un}} )\to  \on{Ran}_{x,\on{un}}$ ``remembers'' $I_0\subset I_3$. One can similarly prove variants of the lemmas.
\end{remark}

\begin{construction} We apply Construction \ref{constr-propogation-fact-mod} to the factorization category $\mathcal{A}:=\mathcal{KL}(H)_{\mathscr{L}}$ and its factorization module category $\mathcal{C}:=  \widetilde{\mathcal{L}\mathfrak{h}}\on{-}mod^{\mathcal{L}^+H\times_{H_x} K_x }_{\mathscr{L}}$. Using the fact that $\on{Loc}$ is a morphism between crystals of categories, the natural transformation (\ref{eqn-nt-prop-res-oblv}) induces a natural transformation
\begin{equation} \label{eqn-nt-loc-prop}
\xymatrix{
	B\on{-FactMod}_{\on{un}}( \widetilde{\mathcal{L}\mathfrak{h}}\on{-}mod^{\mathcal{L}^+H\times_{H_x} K_x }_{\mathscr{L}} ) \ar[d]^-{\on{prop}_{A}} \ar[r]^-{\on{Loc}}
	& \DMod(\on{Bun}_H ^{K\on{-level}_x})_{\mathscr{T}} \otimes \DMod( \on{Ran}_{x,\on{un,dR}} )  \\
	A\on{-FactMod}_{\on{un}}( \mathbf{cores}_{p_{t,\dR}}( \widetilde{\mathcal{L}\mathfrak{h}}\on{-}mod^{\mathcal{L}^+H\times_{H_x} K_x }_{\mathscr{L}} ) )\ar@{=>}[ru]\ar[r]^-{\on{Loc}}
	&
	\DMod(\on{Bun}_H ^{K\on{-level}_x})_{\mathscr{T}} \otimes \DMod(\on{Arrow}(\on{Ran}_{x,\on{un,dR}}) ) \ar[u]^-{p_{t,!}}.
}
\end{equation}
\end{construction}

The following is the main result of this subsection.

\begin{lemma} \label{lem-loc-compatible-with-prop}
The diagram (\ref{eqn-nt-loc-prop}) commutes after taking chiral homology, i.e., after composing with
$$ \on{C}_{x,\on{un}}^{\on{ch}}: \DMod(\on{Bun}_H ^{K\on{-level}_x})_{\mathscr{T}} \otimes \DMod( \on{Ran}_{x,\on{un,dR}} ) \to  \DMod(\on{Bun}_H ^{K\on{-level}_x})_{\mathscr{T}}.$$
The similar claim for $\on{Loc}_{\on{co}}$ is also true.
\end{lemma}

\begin{warning} The lemma does \emph{not} follow immediately from Propogation-restriction Lemma for Chiral homology (see Lemma \ref{lem-AB-lemma}) because $\on{Loc}$ has no factorization structure. See Warning \ref{warn-no-fact-structure-on-Loc}.

\end{warning}

The rest of this subsection is devoted to the proof of the lemma. We only prove the claim for $\on{Loc}$ because the proof for $\on{Loc}_{\on{co}}$ is similar. Our strategy is to reduce to problem to a similar claim about another morphism
$$ \on{Loc}^+: \mathcal{L}^+\mathfrak{h}\on{-}mod^{\mathcal{L}^+H\times_{H_x} K_x } \to  \DMod(\on{Bun}_H ^{K\on{-level}_x})_{/\on{Bun}_H}\boxtimes \mathbf{1}$$
which is equipped with factorization structures such that we can deduce the claim from Propogation-restriction Lemma for Chiral homology (see Lemma \ref{lem-AB-lemma}). We first explain the construction of $\on{Loc}^+$.

\begin{construction} \label{constr-relative-crystal-localization}
By assumption, the twisting $\mathscr{T}$ on $\on{Bun}_H ^{K\on{-level}_x}$ is the restriction of a twisting on $\on{Bun}_H$. Hence we have an adjoint pair
$$\on{ind}: \DMod(\on{Bun}_H ^{K\on{-level}_x})_{/\on{Bun}_H} \adjoint \DMod(\on{Bun}_H ^{K\on{-level}_x})_{\mathscr{T}} : \on{oblv}, $$
where 
$$\DMod(\on{Bun}_H ^{K\on{-level}_x})_{/\on{Bun}_H} \simeq \IndCoh( \on{Bun}_H\times_{ \mathbb{B}H_{x,\dR} }  \mathbb{B}K_{x,\dR} ) $$
is the category of relative D-modules.

Let $U$ be a quasi-compact open substack of $\on{Bun}_H$. Consider the Cech nerve $U^{K\on{-level}_x,\bullet}$ of the projection 
$$
U^{K\on{-level}_x} = U\times_{ \mathbb{B}H_{x} }  \mathbb{B}K_{x} \to
U\times_{ \mathbb{B}H_{x,\dR} }  \mathbb{B}K_{x,\dR}. $$
It is easy to see the colimit of this simplicial laft prestack is just $U\times_{ \mathbb{B}H_{x,\dR} }  \mathbb{B}K_{x,\dR}$. Hence 
$$ \DMod(U ^{K\on{-level}_x})_{/U}  \simeq \lim_{!\on{-pull}} \IndCoh( U^{K\on{-level}_x,\bullet} ) \simeq \colim_{!\on{-push}} \IndCoh( U^{K\on{-level}_x,\bullet} ) $$
We have a simplicial map $ U^{K\on{-level}_x,\bullet} \times \Ran_{x,\on{un,dR}} \to \widehat{ \on{Hecke} }_{H,U}^{K\on{-level}_x,\bullet} $. Since the gerbe in Remark \ref{rem-compatible-pair-two-gerbe} is equipped with a trivialization on $U\times_{ \mathbb{B}H_{x,\dR} }  \mathbb{B}K_{x,\dR}$, the twisted $!$-pullback functors induce a cosimplicial lax-unital morphism
$$  \mathcal{I}nd\mathcal{C}oh(   \widehat{\on{Hecke}}_{H,U}^{K\on{-level}_x,\bullet}  )_{\mathscr{L},\Delta} \to  \IndCoh( U^{K\on{-level}_x,\bullet}) \boxtimes \mathbf{1},$$
and the twisted $!$-pushforward functors induce a simplicial unital morphism
$$\IndCoh( U^{K\on{-level}_x,\bullet}) \boxtimes \mathbf{1}\to  \mathcal{I}nd\mathcal{C}oh(   \widehat{\on{Hecke}}_{H,U}^{K\on{-level}_x,\bullet}  )_{\mathscr{L},\Delta^\op}.$$
Passing to limits and colimits, we obtain an adjoint pair
$$\on{ind}:  \DMod(U ^{K\on{-level}_x})_{/U}\boxtimes \mathbf{1}\adjoint \mathscr{L}_{\widehat{\on{Hecke}}_H^{K\on{-level}_x} } \on{-mod}( \mathcal{I}nd\mathcal{C}oh( U^{K\on{-level}_x})  ) :\on{oblv} . $$

By construction, the composition
$$  \DMod(U ^{K\on{-level}_x})_{/U}\boxtimes \mathbf{1}\adjoint \mathscr{L}_{\widehat{\on{Hecke}}_H^{K\on{-level}_x} } \on{-mod}( \mathcal{I}nd\mathcal{C}oh( U^{K\on{-level}_x})  ) \adjoint  \DMod( U^{K\on{-level}_x} )_{\mathscr{T}}\boxtimes \mathbf{1} $$
is just the the adjoint pair between categories of relative D-modules and (twisted) absolute D-modules.

\end{construction}

\begin{construction} Similar to Construction \ref{const-local-to-global-hecke}, we can construct a unital morphism
\begin{equation} \label{eqn-loc-relative-crystal}
 \on{Loc}_U^+: \mathcal{L}^+\mathfrak{h}\on{-}mod^{\mathcal{L}^+H\times_{H_x} K_x } \to  \DMod(U ^{K\on{-level}_x})_{/U}\boxtimes \mathbf{1} 
 \end{equation}
by writting both sides as simplicial colimits (or simplicial limits) and using $!$-pullback functors to connect them. Moreover, when $x$ corresponds to the empty set, the unital morphism
\begin{equation} \label{eqn-loc-relative-crystal-alg}
 \on{Loc}_U^+: \mathcal{R}ep(\mathcal L^+ H) \to  \IndCoh(U )\boxtimes \mathbf{1} 
 \end{equation}
can be upgraded to a unital factorization functor, where the LHS is a factorization category over $\on{Ran}_{\on{un,dR}}$, while the RHS is over $U\times \on{Ran}_{\on{un,dR}}$ (or equivalently, a $U$-family of factorization categories). Also, the morphism (\ref{eqn-loc-relative-crystal}) can be upgraded to a morphism between factorization module categories intertwining (\ref{eqn-loc-relative-crystal-alg}), where the RHS is viewed as a factorization module category over $ (U ^{K\on{-level}_x})_\dR \times_{U_\dR} U \times \Ran_{x,\on{un,dR}}$.

It follows we have a commutative diagram
\begin{equation} \label{eqn-loc-full-sph}
\xymatrix{
	\widetilde{\mathcal{L}\mathfrak{h}}\on{-}mod^{\mathcal{L}^+H\times_{H_x} K_x }_{\mathscr{L}} \ar[r] \ar[d]^-{\on{oblv}} &
	\mathscr{L}_{\widehat{\on{Hecke}}_H^{K\on{-level}_x} } \on{-mod}( \mathcal{I}nd\mathcal{C}oh( U^{K\on{-level}_x})  ) \ar[d]^-{\on{oblv}} \\
	\mathcal{L}^+\mathfrak{h}\on{-}mod^{\mathcal{L}^+H\times_{H_x} K_x } \ar[r]^-{\on{Loc}_U^+} &
	  \DMod(U ^{K\on{-level}_x})_{/U}\boxtimes \mathbf{1}.
}
\end{equation}
Moreover, using the base-change isomorphisms, one can show that this square is left adjointable along the vertical direction, i.e., it induces a commutative diagram
\begin{equation} \label{eqn-loc-full-sph-2}
\xymatrix{
	\widetilde{\mathcal{L}\mathfrak{h}}\on{-}mod^{\mathcal{L}^+H\times_{H_x} K_x }_{\mathscr{L}} \ar[r]  &
	\mathscr{L}_{\widehat{\on{Hecke}}_H^{K\on{-level}_x} } \on{-mod}( \mathcal{I}nd\mathcal{C}oh( U^{K\on{-level}_x})  ) \\
	\mathcal{L}^+\mathfrak{h}\on{-}mod^{\mathcal{L}^+H\times_{H_x} K_x } \ar[r]^-{\on{Loc}_U^+}\ar[u]^-{\on{ind}} &
	  \DMod(U ^{K\on{-level}_x})_{/U}\boxtimes \mathbf{1}\ar[u]^-{\on{ind}},
}
\end{equation}
which then induces a commutative diagram
\begin{equation} \label{eqn-loc-full-sph-3}
\xymatrix{
	\widetilde{\mathcal{L}\mathfrak{h}}\on{-}mod^{\mathcal{L}^+H\times_{H_x} K_x }_{\mathscr{L}} \ar[r]^-{\on{Loc}_U}  &
	\DMod(U ^{K\on{-level}_x})\boxtimes \mathbf{1} \\
	\mathcal{L}^+\mathfrak{h}\on{-}mod^{\mathcal{L}^+H\times_{H_x} K_x } \ar[r]^-{\on{Loc}_U^+}\ar[u]^-{\on{ind}} &
	  \DMod(U ^{K\on{-level}_x})_{/U}\boxtimes \mathbf{1}\ar[u]^-{\on{ind}}.
}
\end{equation}
Passing to right adjoints along vertical direction (and allowing $U$ to change), we obtain a 2-morphism
\begin{equation} \label{eqn-natural-transformation-loc-oblv}
\xymatrix{
	\widetilde{\mathcal{L}\mathfrak{h}}\on{-}mod^{\mathcal{L}^+H\times_{H_x} K_x }_{\mathscr{L}} \ar[r]^-{\on{Loc}} \ar[d]^-{\on{oblv}} &
	\DMod( \on{Bun}_H^{K\on{-level}_x} )_{\mathscr{T}}\boxtimes \mathbf{1}\ar[d]^-{\on{oblv}} \\
	\mathcal{L}^+\mathfrak{h}\on{-}mod^{\mathcal{L}^+H\times_{H_x} K_x } \ar[r]^-{\on{Loc}^+} \ar@{=>}[ru]&
	  \DMod(\on{Bun}_H ^{K\on{-level}_x})_{/\on{Bun}_H}\boxtimes \mathbf{1}.
}
\end{equation}
\end{construction}

This 2-morphism is \emph{not} invertible in general. However, we are going to prove the following result and use it to deduce Lemma \ref{lem-loc-compatible-!-pull}.

\begin{lemma} \label{lem-loc-compatible-oblv}
Up to chiral homology, the digram (\ref{eqn-natural-transformation-loc-oblv}) commutes for \emph{factorization objects}. More precisly, for any unital factorization algebra $A$ in $ \mathcal{KL}(H)_{\mathscr{L}}$, the natural transformation
\begin{equation} \label{eqn-lem-loc-compatible-oblv}
	\xymatrix{
	A\on{-FactMod}_{\on{un}}( \widetilde{\mathcal{L}\mathfrak{h}}\on{-}mod^{\mathcal{L}^+H\times_{H_x} K_x }_{\mathscr{L}} ) \ar[d]^-{\on{oblv}} 
	\ar[r]^-{ \on{Loc} } 
	& \DMod(\on{Bun}_H ^{K\on{-level}_x})_{\mathscr{T}} \otimes \DMod( \on{Ran}_{x,\on{un,dR}}) \ar[d]^-{\on{oblv}} \\
	\on{oblv}(A)\on{-FactMod}_{\on{un}}( \mathcal{L}^+\mathfrak{h}\on{-}mod^{\mathcal{L}^+H\times_{H_x} K_x }  ) 
	\ar[r]^-{\on{Loc}^+ } \ar@{=>}[ur]
	& \DMod(\on{Bun}_H ^{K\on{-level}_x})_{/\on{Bun}_H}  \otimes \DMod( \on{Ran}_{x,\on{un,dR}}) ,
	}
\end{equation}
is invertible after composing with 
$$ \on{C}^{\on{ch}}_{x,\on{un}}: \DMod(\on{Bun}_H ^{K\on{-level}_x})_{/\on{Bun}_H}  \otimes \DMod( \on{Ran}_{x,\on{un,dR}}) \to \DMod(\on{Bun}_H ^{K\on{-level}_x})_{/\on{Bun}_H}.$$
The similar claim for $\on{Loc}_{\on{co}}$ is also true.
\end{lemma}

We also need the following variant, which can be proved similarly\footnote{One needs to use the construction in Remark \ref{rem-RanRanRan}}.

\begin{variant} \label{lem-loc-compatible-oblv-variant}
Lemma \ref{lem-loc-compatible-oblv} remains valid if we 
\begin{itemize}
	\item replace $\on{Ran}_{x,\on{un,dR}}$ by $\on{Arrow}(\on{Ran}_{x,\on{un,dR}})$;
	\item replace $\widetilde{\mathcal{L}\mathfrak{h}}\on{-}mod^{\mathcal{L}^+H\times_{H_x} K_x }_{\mathscr{L}}$ and $\mathcal{L}^+\mathfrak{h}\on{-}mod^{\mathcal{L}^+H\times_{H_x} K_x }  $ by their (co)restriction along $p_{t,\dR}$;
	\item replcae the chiral homology functor $\on{C}^{\on{ch}}_{x,\on{un}}$ by $p_{s,!}$.
\end{itemize}
\end{variant}

We first deduce Lemma \ref{lem-loc-compatible-with-prop} from Lemma \ref{lem-loc-compatible-oblv}.

\noindent \emph{Proof of Lemma \ref{lem-loc-compatible-with-prop}}: Note that we have an obvious commutative diagram
$$
\xymatrix{
	B\on{-FactMod}_{\on{un}}( \widetilde{\mathcal{L}\mathfrak{h}}\on{-}mod^{\mathcal{L}^+H\times_{H_x} K_x }_{\mathscr{L}} ) \ar[d]^-{\on{prop}_{A}} \ar[r]^-{\on{oblv}}
	& \on{oblv}(B)\on{-FactMod}_{\on{un}}( \mathcal{L}^+\mathfrak{h}\on{-}mod^{\mathcal{L}^+H\times_{H_x} K_x } ) \ar[d]^-{\on{prop}_{\on{oblv}(A)}}  \\
	A\on{-FactMod}_{\on{un}}( \mathbf{cores}_{p_{t,\dR}}( \widetilde{\mathcal{L}\mathfrak{h}}\on{-}mod^{\mathcal{L}^+H\times_{H_x} K_x }_{\mathscr{L}} ) ) \ar[r]^-{\on{oblv}}
	&
	 \on{oblv}(A)\on{-FactMod}_{\on{un}}(\mathbf{cores}_{p_{t,\dR}}(  \mathcal{L}^+\mathfrak{h}\on{-}mod^{\mathcal{L}^+H\times_{H_x} K_x }) ) .
}
$$
Hence by Lemma \ref{eqn-lem-loc-compatible-oblv} and Variant \ref{lem-loc-compatible-oblv-variant}, we can reduce the problem to a similar claim about $\on{Loc}^+$. But then it follows immediately from Propogation-restriction Lemma for Chiral homology (see Lemma \ref{lem-AB-lemma}) because $\on{Loc}^+$ has a factorization structure.

\qed

The rest of this subsection is devoted to the proof of Lemma \ref{lem-loc-compatible-oblv}. The strategy is to use Propogation-restriction Lemma for Chiral homology (see Lemma \ref{lem-AB-lemma}) again to reduce to the case when $A$ is the unit factorization algebra.

\begin{notation} To avoid overboxes, we use the following shorthands:
\begin{equation*}
\begin{aligned}
	\mathcal{A}^+:= \mathcal{R}ep(\mathcal{L}^+H),\;\mathcal{A}:=\mathcal{KL}(H)_{\mathscr{L}},\; 
 \mathcal{C}^+:= \mathcal{L}^+\mathfrak{h}\on{-}mod^{\mathcal{L}^+H\times_{H_x} K_x },\; \mathcal{C}:=  \widetilde{\mathcal{L}\mathfrak{h}}\on{-}mod^{\mathcal{L}^+H\times_{H_x} K_x }_{\mathscr{L}} ,\\ D^+:= \DMod(U ^{K\on{-level}_x})_{/U} ,\; \mathcal{D}:=\mathscr{L}_{\widehat{\on{Hecke}}_H^{K\on{-level}_x} } \on{-mod}( \mathcal{I}nd\mathcal{C}oh( U^{K\on{-level}_x})  ),\; D':=  \DMod(U ^{K\on{-level}_x})_{\mathscr{T}}. 
 \end{aligned}
 \end{equation*}
We write $\mathbb{V}$ for the unit factorization algebra of $\mathcal{A}$. 
\end{notation}

\begin{construction} Consider the unital morphism $\on{Loc}_U:\mathcal{C} \to \mathcal{D} \xrightarrow{\on{ind}}  D'\boxtimes \mathbf{1}$. By Construction \ref{const-prop-unit=insert-unit}, it induces the following commutative diagram
$$ 
\xymatrix{
	A\on{-FactMod}_{\on{un}}(\mathcal{C}) \ar[d]^-{\on{prop}_{\mathbb{V}}}
	\ar[rr]^-{[\on{Loc}_U]_{\on{Id}_{\dR}}} & &
	D'\otimes \DMod( \on{Ran}_{x,\on{un}} ) \ar[d]^-{\on{ins}_{\on{unit}}\simeq p_s^!}
	\\	\mathbb{V}\on{-FactMod}_{\on{un}}(\mathbf{cores}_{p_{t,\dR}}(\mathcal{C}) ) \ar[rr]^-{[\on{Loc}_U]_{p_{t,\dR}}} & &
	D'\otimes \DMod( \on{Arrow}_{ \on{Ran}_{x,\on{un}}} ).
}
 $$
 It is easy to see $p_s$ is homologically contractible, i.e., $p_s^!$ is fully faithful. Hence we have a commutative diagram
\begin{equation} \label{eqn-prop-unit-loc-chiral-homology}
\xymatrix{
	A\on{-FactMod}_{\on{un}}(\mathcal{C}) \ar[d]^-{\on{prop}_{\mathbb{V}}}
	\ar[rr]^-{[\on{Loc}_U]_{\on{Id}_{\dR}}} & &
	D'\otimes \DMod( \on{Ran}_{x,\on{un}} ) 
	\\	\mathbb{V}\on{-FactMod}_{\on{un}}(\mathbf{cores}_{p_{t,\dR}}(\mathcal{C}) ) \ar[rr]^-{[\on{Loc}_U]_{p_{t,\dR}}} & &
	D'\otimes \DMod( \on{Arrow}_{ \on{Ran}_{x,\on{un}}} ) \ar[u]^-{p_{s,!}}.
}
 \end{equation}
 \end{construction}

\begin{construction}
On the other hand, the factorization structures on $\on{Loc}_{U}^+: \mathcal{C}^+ \to D^+\boxtimes \mathbf{1}$ induces a commutative diagram
$$ 
\xymatrix{
	\on{oblv}(A)\on{-FactMod}_{\on{un}}(\mathcal{C}^+) \ar[d]^-{\on{prop}_{\oblv(\mathbb{V})}}
	\ar[rr]^-{[\on{Loc}^+_U]_{\on{Id}_{\dR}}} & &
	\on{Loc}^+_U\circ \on{oblv}(A)\on{-FactMod}_{\on{un}}( D^+\boxtimes \mathbf{1} )  \ar[d]^-{\on{prop}_{\on{Loc}^+_U\circ \on{oblv}(\mathbb{V})}}
	\\	\oblv(\mathbb{V})\on{-FactMod}_{\on{un}}(\mathbf{cores}_{p_{t,\dR}}(\mathcal{C}^+) ) \ar[rr]^-{[\on{Loc}^+_U]_{p_{t,\dR}}} & &
	\on{Loc}^+_U\circ \on{oblv}(\mathbb{V})\on{-FactMod}_{\on{un}}( D^+\boxtimes \mathbf{cores}_{p_{t,\dR}}(\mathbf{1}) ).
}
 $$
By Propogation-restriction Lemma for Chiral homology (see Lemma \ref{lem-AB-lemma}), the natural transformation
$$
\xymatrix{
	\on{Loc}^+_U\circ \on{oblv}(A)\on{-FactMod}_{\on{un}}( D^+\boxtimes \mathbf{1} )  \ar[d]^-{\on{prop}_{\on{Loc}^+_U\circ \on{oblv}(\mathbb{V})}} \ar[r]^-{\on{oblv}}
	& D^+ \otimes \DMod( \Ran_{x,\on{un}})
	\\	
	\on{Loc}^+_U\circ \on{oblv}(\mathbb{V})\on{-FactMod}_{\on{un}}( D^+\boxtimes \mathbf{cores}_{p_{t,\dR}}(\mathbf{1}) ) \ar[r]^-{\on{oblv}} \ar@{=>}[ru]
	& D^+ \otimes \DMod( \on{Arrow}(\Ran_{x,\on{un}}))  \ar[u]^-{p_{t,!}}.
}
 $$
 becomes invertible after taking chiral homology. Hence
$$ 
\xymatrix{
	\on{oblv}(A)\on{-FactMod}_{\on{un}}(\mathcal{C}^+) \ar[d]^-{\on{prop}_{\oblv(\mathbb{V})}}
	\ar[rr]^-{[\on{Loc}^+_U]_{\on{Id}_{\dR}}} & &
	D^+ \otimes \DMod( \Ran_{x,\on{un}})
	\\	\oblv(\mathbb{V})\on{-FactMod}_{\on{un}}(\mathbf{cores}_{p_{t,\dR}}(\mathcal{C}^+) ) \ar[rr]^-{[\on{Loc}^+_U]_{p_{t,\dR}}}  \ar@{=>}[rru] & &
	D^+ \otimes \DMod( \on{Arrow}(\Ran_{x,\on{un}})) \ar[u]^-{p_{t,!}}
}
 $$
  becomes invertible after taking chiral homology.
\end{construction}

\begin{construction}
 Finally, similar to (\ref{eqn-lem-loc-compatible-oblv}), we have a natural transformation
 \begin{equation} \label{eqn-lem-loc-compatible-oblv-prop}
	\xymatrix{
	\mathbb{V}\on{-FactMod}_{\on{un}}(\mathbf{cores}_{p_{t,\dR}}(\mathcal{C}) )\ar[d]^-{\on{oblv}} 
	\ar[rr]^-{p_{s,!}\circ [\on{Loc}_U]_{p_{t,\dR}} } 
	& & D'\otimes \DMod( \on{Ran}_{x,\on{un}} )   \ar[d]^-{\on{oblv}} \\
	\oblv(\mathbb{V})\on{-FactMod}_{\on{un}}(\mathbf{cores}_{p_{t,\dR}}(\mathcal{C}^+) ) 
	\ar[rr]^-{p_{s,!}\circ [\on{Loc}^+_U]_{p_{t,\dR}} }  \ar@{=>}[urr]
	& & D^+\otimes \DMod( \on{Ran}_{x,\on{un}} ).
	}
\end{equation}
\end{construction}

\noindent \emph{Proof of Lemma \ref{lem-loc-compatible-oblv}}:
Combining the above three constructions, we have reduced Lemma \ref{lem-loc-compatible-oblv} to the following result.

\qed

\begin{lemma} \label{lem-loc-compatible-oblv-unit}
The natural transformation (\ref{eqn-lem-loc-compatible-oblv-prop}) is invertible.
\end{lemma}

\proof
Recall
$$ \mathbf{\Gamma}(\on{Id}_\dR,\mathcal{C}) \xrightarrow{\on{ins}_{\on{unit}}} \mathbb{V}\on{-FactMod}_{\on{un}}(\mathbf{cores}_{p_{t,\dR}}(\mathcal{C}) )$$
is an equivalence. Also recall the image of
$$ \mathbf{\Gamma}(\on{Id}_\dR,\mathcal{C}^+) \xrightarrow{\on{ind}} \mathbf{\Gamma}(\on{Id}_\dR,\mathcal{C}) $$
generates the target category because its right adjoint is conservative. Hence we only need to show (\ref{eqn-lem-loc-compatible-oblv-prop}) is invertible after precomposing with 
$$ \mathbf{\Gamma}(\on{Id}_\dR,\mathcal{C}^+) \xrightarrow{\on{ind}} \mathbf{\Gamma}(\on{Id}_\dR,\mathcal{C}) \xrightarrow{\on{ins}_{\on{unit}}} \mathbb{V}\on{-FactMod}_{\on{un}}(\mathbf{cores}_{p_{t,\dR}}(\mathcal{C}) ). $$
In other words, we only need to show
$$
\xymatrix{
	 \mathbf{\Gamma}(\on{p}_{t,\dR},\mathcal{C}) \ar[d]^-{\on{oblv}} 
	\ar[rr]^-{p_{s,!}\circ [\on{Loc}_U]_{p_{t,\dR}}} 
	& & D'\otimes \DMod( \on{Ran}_{x,\on{un}} )   \ar[d]^-{\on{oblv}} \\
	 \mathbf{\Gamma}(\on{p}_{t,\dR},\mathcal{C}^+)
	\ar[rr]_-{p_{s,!}\circ [\on{Loc}^+_U]_{p_{t,\dR}}} \ar@{=>}[urr]
	& & D^+\otimes \DMod( \on{Ran}_{x,\on{un}} ) ,
	}
$$
is invertible after precomposing with
$$ \mathbf{\Gamma}(\on{Id}_\dR,\mathcal{C}^+) \xrightarrow{\on{ind}} \mathbf{\Gamma}(\on{Id}_\dR,\mathcal{C}) \xrightarrow{\on{ins}_{\on{unit}}} \mathbf{\Gamma}(\on{p}_{t,\dR},\mathcal{C}).$$
Since $\on{ind}:\mathcal{C}^+ \to \mathcal{C} $ is unital, the above functor is equivalent to
$$\mathbf{\Gamma}(\on{Id}_\dR,\mathcal{C}^+) \xrightarrow{\on{ins}_{\on{unit}}} \mathbf{\Gamma}(\on{p}_{t,\dR},\mathcal{C}^+) \xrightarrow{\on{ind}} \mathbf{\Gamma}(\on{p}_{t,\dR},\mathcal{C}) 
  $$
Note that we have a commutative diagram:
$$
\xymatrix{
	\mathbf{\Gamma}(\on{Id}_\dR,\mathcal{C}^+) 
	\ar[r]^-{\on{ins}_{\on{unit}}}  \ar[d]^-{ \on{Loc}^+ }
	& \mathbf{\Gamma}(\on{p}_{t,\dR},\mathcal{C}^+)\ar[r]^-{\on{ind}}  \ar[d]^-{ \on{Loc}^+ }
	& \mathbf{\Gamma}(\on{p}_{t,\dR},\mathcal{C}) \ar[r]^-{\on{oblv}}  \ar[d]^-{ \on{Loc} }
	&  \mathbf{\Gamma}(\on{p}_{t,\dR},\mathcal{C}^+) \ar[d]^-{ \on{Loc}^+ } \\
	\mathbf{\Gamma}(\on{Id}_\dR,D^+\boxtimes \mathbf{1}) \ar[r]^-{\on{ins}_{\on{unit}}} 
	& \mathbf{\Gamma}(\on{p}_{t,\dR},D^+\boxtimes \mathbf{1})  \ar[r]^-{\on{ind}} 
	& \mathbf{\Gamma}(\on{p}_{t,\dR},\mathcal{D})  \ar[r]^-{\on{oblv}}
	&  \mathbf{\Gamma}(\on{p}_{t,\dR},D^+\boxtimes \mathbf{1})  ,
}
$$
where the left square is because $\on{Loc}^+:\mathcal{C}^+\to D^+\boxtimes \mathbf{1} $ is unital, and the right two squares follow from (\ref{eqn-loc-full-sph}) and (\ref{eqn-loc-full-sph-2}) (by taking $\mathbf{\Gamma}(\on{p}_{t,\dR},-)$). Hence we only need to show the natural transformation from
\begin{equation*}
 	\mathbf{\Gamma}(\on{Id}_\dR,D^+\boxtimes \mathbf{1}) \xrightarrow{\on{ins}_{\on{unit}}} 
	 \mathbf{\Gamma}(\on{p}_{t,\dR},D^+\boxtimes \mathbf{1})  \xrightarrow{\on{ind}}  \mathbf{\Gamma}(\on{p}_{t,\dR},D'\boxtimes \mathbf{1})  \xrightarrow{\on{oblv}}
	 \mathbf{\Gamma}(\on{p}_{t,\dR},D^+\boxtimes \mathbf{1}) \xrightarrow{p_{s,!}} \mathbf{\Gamma}(\on{Id}_\dR,D^+\boxtimes \mathbf{1}) 
\end{equation*}
to
$$ 	\mathbf{\Gamma}(\on{Id}_\dR,D^+\boxtimes \mathbf{1}) \xrightarrow{\on{ins}_{\on{unit}}} 
	 \mathbf{\Gamma}(\on{p}_{t,\dR},D^+\boxtimes \mathbf{1})  \xrightarrow{\on{ind}}  \mathbf{\Gamma}(\on{p}_{t,\dR},\mathcal{D})  \xrightarrow{\on{oblv}}
	 \mathbf{\Gamma}(\on{p}_{t,\dR},D^+\boxtimes \mathbf{1}) \xrightarrow{p_{s,!}} \mathbf{\Gamma}(\on{Id}_\dR,D^+\boxtimes \mathbf{1}) $$
is invertible. 

Note that we also have a natural transformation from
\begin{equation*}
 	\mathbf{\Gamma}(x,D^+\boxtimes \mathbf{1}) \xrightarrow{\on{ins}_{\on{unit}}} 
	 \mathbf{\Gamma}(\on{Id}_\dR,D^+\boxtimes \mathbf{1})  \xrightarrow{\on{ind}}  \mathbf{\Gamma}(\on{Id}_\dR,D'\boxtimes \mathbf{1})  \xrightarrow{\on{oblv}}
	 \mathbf{\Gamma}(\on{Id}_\dR,D^+\boxtimes \mathbf{1}) \xrightarrow{p_{s,!}} \mathbf{\Gamma}(x,D^+\boxtimes \mathbf{1}) 
\end{equation*}
to
$$ \mathbf{\Gamma}(x,D^+\boxtimes \mathbf{1}) \xrightarrow{\on{ins}_{\on{unit}}} 
	 \mathbf{\Gamma}(\on{Id}_\dR,D^+\boxtimes \mathbf{1})  \xrightarrow{\on{ind}}  \mathbf{\Gamma}(\on{Id}_\dR,\mathcal{D})  \xrightarrow{\on{oblv}}
	 \mathbf{\Gamma}(\on{Id}_\dR,D^+\boxtimes \mathbf{1}) \xrightarrow{p_{s,!}} \mathbf{\Gamma}(x,D^+\boxtimes \mathbf{1}) .$$
By the propogation lemma for chiral homology (see Lemma \ref{lem-AA-lemma}), the natural transformation in the last paragraph can be obtained from this one by tensoring up with $\DMod(\on{Ran}_{x,\on{un,dR}})$. Hence we only need to show this natural transformation is invertible.

In other words, we only need to show the chiral homology of the monad
$$ \DMod(U ^{K\on{-level}_x})_{/U}\boxtimes \mathbf{1}\to \mathscr{L}_{\widehat{\on{Hecke}}_H^{K\on{-level}_x} } \on{-mod}( \mathcal{I}nd\mathcal{C}oh( U^{K\on{-level}_x})  )  \to \DMod(U ^{K\on{-level}_x})_{/U}\boxtimes \mathbf{1}$$
is equivalent to the monad
$$ \DMod(U ^{K\on{-level}_x})_{/U} \to \DMod(U ^{K\on{-level}_x} )_{\mathscr{T}} \to \DMod(U ^{K\on{-level}_x})_{/U}. $$
When $\mathscr{T}$ is trivial and $x$ corresponds to the empty subset, this is the main theorem of \cite{rozenblyum2021connections}. The general case can be proved using the same method\footnote{Details will be provided in the second author's thesis.}.

\subsection{Compatibility with Restriction}
\label{ssec-loc-compatible-restriction}

\begin{notation}
In this and the following subsection, $U$ (resp. $V$) is a quasi-compact open substack of $\on{Bun}_H$ (resp. $\on{Bun}_{H'}$) such that the image of $V$ under the map $\on{Bun}_{H'} \to \on{Bun}_{H}$ is contained in $U$.
\end{notation}

We first construct a natural transformation unerlying the desired commutative square.

\begin{construction}
By Construction \ref{constr-Hecke-indcoh-vs-dmod}, we have a commutative diagram induced by twisted $!$-pullback functors
$$
\xymatrix{
		\DMod( U^{K\on{-level}_x} )_{\mathscr{T}}\boxtimes \mathbf{1} \ar[r]^-{\on{oblv}} \ar[d]^-{!\on{-pull}}
		& \mathscr{L}_{\widehat{\on{Hecke}}_H^{K\on{-level}_x} } \on{-mod}( \mathcal{I}nd\mathcal{C}oh( U^{K\on{-level}_x})  ) \ar[d]^-{!\on{-pull}}\\
		\DMod( V^{K'\on{-level}_x} )_{\mathscr{T}'}\boxtimes \mathbf{1} \ar[r]^-{\on{oblv}}
		& \mathscr{L}'_{\widehat{\on{Hecke}}_{H'}^{K'\on{-level}_x} } \on{-mod}( \mathcal{I}nd\mathcal{C}oh( V^{K'\on{-level}_x})  ).
}
$$
Passing to left adjoints, we obtain a $2$-morphism
$$
\xymatrix{
		\mathscr{L}_{\widehat{\on{Hecke}}_H^{K\on{-level}_x} } \on{-mod}( \mathcal{I}nd\mathcal{C}oh( U^{K\on{-level}_x})  ) \ar[r]^-{\on{ind}} \ar[d]^-{!\on{-pull}}
		& \DMod( U^{K\on{-level}_x} )_{\mathscr{T}}\boxtimes \mathbf{1} \ar[d]^-{!\on{-pull}}\\
		\mathscr{L}'_{\widehat{\on{Hecke}}_{H'}^{K'\on{-level}_x} } \on{-mod}( \mathcal{I}nd\mathcal{C}oh( V^{K'\on{-level}_x})  )
		 \ar[r]^-{\on{ind}} \ar@{=>}[ur]
		&\DMod( V^{K'\on{-level}_x} )_{\mathscr{T}'}\boxtimes \mathbf{1} .
}
$$
On the other hands, by Construction \ref{const-local-to-global-hecke} we have a commutative diagram induced by twisted $!$-pullback functors:
$$
\xymatrix{
	\widetilde{\mathcal{L}\mathfrak{h}}\on{-}mod^{\mathcal{L}^+H\times_{H_x} K_x }_{\mathscr{L}} \ar[r]^-{!\on{-pull}} \ar[d]^-{\on{res}}  &
	\mathscr{L}_{\widehat{\on{Hecke}}_H^{K\on{-level}_x} } \on{-mod}( \mathcal{I}nd\mathcal{C}oh( U^{K\on{-level}_x})  ) \ar[d]^-{!\on{-pull}}  \\
	\widetilde{\mathcal{L}\mathfrak{h}'}\on{-}mod^{\mathcal{L}^+H'\times_{H'_x} K'_x }_{\mathscr{L}'} \ar[r]^-{!\on{-pull}}  &
	\mathscr{L}'_{\widehat{\on{Hecke}}_{H'}^{K'\on{-level}_x} } \on{-mod}( \mathcal{I}nd\mathcal{C}oh( V^{K'\on{-level}_x})  ).
}
$$
Hence we obtain a $2$-morphism
\begin{equation} \label{eqn-natural-transformation-loc-restriction}
\xymatrix{
		\widetilde{\mathcal{L}\mathfrak{h}}\on{-}mod^{\mathcal{L}^+H\times_{H_x} K_x }_{\mathscr{L}}  \ar[r]^-{\on{Loc}_U}  \ar[d]^-{\on{res}} 
		& \DMod( U^{K\on{-level}_x} )_{\mathscr{T}}\boxtimes \mathbf{1} \ar[d]^-{!\on{-pull}}\\
		\widetilde{\mathcal{L}\mathfrak{h}'}\on{-}mod^{\mathcal{L}^+H'\times_{H'_x} K'_x }_{\mathscr{L}'} 
		 \ar[r]^-{\on{Loc}_V} \ar@{=>}[ur]
		&\DMod( V^{K'\on{-level}_x} )_{\mathscr{T}'}\boxtimes \mathbf{1}.
}
\end{equation}
\end{construction}

This $2$-morphism is \emph{not} invertible in general. However, we are going to prove the following result, which implies Theorem \ref{thm-localization-KL}(3):

\begin{lemma} \label{lem-loc-compatible-!-pull}
Up to chiral homology, the digram (\ref{eqn-natural-transformation-loc-restriction}) commutes for \emph{factorization objects}. More precisly, for any unital factorization algebra $A$ in $ \mathcal{KL}(H)_{\mathscr{L}}$, the natural transformation
\begin{equation} \label{eqn-lem-loc-compatible-!-pull}
	\xymatrix{
	A\on{-FactMod}_{\on{un}}( \widetilde{\mathcal{L}\mathfrak{h}}\on{-}mod^{\mathcal{L}^+H\times_{H_x} K_x }_{\mathscr{L}}  ) \ar[d]^-{\on{res}} 
	\ar[r]^-{[\on{Loc}]_x^{\on{ch}} } 
	& \DMod(\on{Bun}_H ^{K\on{-level}_x})_{\mathscr{T}}  \ar[d]^-{!\on{-pull}} \\
	\on{res}(A)\on{-FactMod}_{\on{un}}( \widetilde{\mathcal{L}\mathfrak{h}'}\on{-}mod^{\mathcal{L}^+H'\times_{H'_x} K'_x }_{\mathscr{L}'}   ) 
	\ar[r]^-{[\on{Loc}]_x^{\on{ch}} } \ar@{=>}[ur]
	& \DMod(\on{Bun}_{H'}^{K'\on{-level}_x})_{\mathscr{T}'},
	}
	\end{equation}
induced by it is invertible.
\end{lemma}

\proof Since
$$\on{oblv}:\DMod(\on{Bun}_{H'} ^{K'\on{-level}_x})_{\mathscr{T}'}\to \DMod(\on{Bun}_{H'} ^{K'\on{-level}_x})_{/\on{Bun}_{H'}}$$
is conservative, we only need to show (\ref{eqn-lem-loc-compatible-!-pull}) becomes invertible after composing with this functor. But this follows from  Lemma \ref{lem-loc-compatible-oblv} (applied to both $(K\to H)$ and $(K'\to H)$) and the following obvious commutative diagrams (where all the arrows are induced by twisted $!$-pullback functors):
\begin{equation*}
\xymatrix{
	\widetilde{\mathcal{L}\mathfrak{h}}\on{-}mod^{\mathcal{L}^+H\times_{H_x} K_x }_{\mathscr{L}} \ar[d]^-{\on{res}} \ar[r]^-{\on{oblv}} 
	& \mathcal{L}^+\mathfrak{h}\on{-}mod^{\mathcal{L}^+H\times_{H_x} K_x }  \ar[d]^-{\on{res}} 
	&
	\DMod(\on{Bun}_H ^{K\on{-level}_x})_{\mathscr{T}}  \ar[r]^-{\on{oblv}} \ar[d]^-{!\on{-pull}} &
	\DMod(\on{Bun}_H ^{K\on{-level}_x})_{/\on{Bun}_H} \ar[d]^-{!\on{-pull}} \\
	\widetilde{\mathcal{L}\mathfrak{h}'}\on{-}mod^{\mathcal{L}^+H'\times_{H'_x} K'_x }_{\mathscr{L}'} \ar[r]^-{\on{oblv}} 
	& \mathcal{L}^+\mathfrak{h}'\on{-}mod^{\mathcal{L}^+H'\times_{H'_x} K'_x },
	& 
	\DMod(\on{Bun}_{H'} ^{K'\on{-level}_x})_{\mathscr{T}'}  \ar[r]^-{\on{oblv}} &
	\DMod(\on{Bun}_{H'} ^{K'\on{-level}_x})_{/\on{Bun}_{H'}},
}
\end{equation*}

$$
\xymatrix{
	 \mathcal{L}^+\mathfrak{h}\on{-}mod^{\mathcal{L}^+H\times_{H_x} K_x }  \ar[d]^-{\on{res}}  \ar[r]^-{\on{Loc}_U^+}
	&
	\DMod(U ^{K\on{-level}_x})_{/U} \boxtimes\mathbf{1} \ar[d]^-{!\on{-pull}} \\
	\mathcal{L}^+\mathfrak{h}'\on{-}mod^{\mathcal{L}^+H'\times_{H'_x} K'_x }
	\ar[r]^-{\on{Loc}_{U'}^+}
	& 
	\DMod( (U') ^{K'\on{-level}_x})_{/U'}\boxtimes\mathbf{1} .
}
$$

\qed

\subsection{Compatibility with Semi-infinite Cohomology}

We first construct a natural transformation unerlying the desired commutative square.

\begin{construction}

It follows from Construction \ref{const-local-to-global-hecke} and Construction \ref{constr-Hecke-indcoh-vs-dmod} that we have a commutative diagram induced by twisted $*$-pushforward functors
$$
	\xymatrix{
	\mathcal{KL}(H')_{\mathscr{L}',\on{co}}   \ar[d]^-{\on{C}^{\frac{\infty}{2}}_{\on{co}}} 
	& \mathscr{L}'_{\widehat{\on{Hecke}}_{H'} } \on{-mod}( \mathcal{I}nd\mathcal{C}oh( V ) ) \ar[l] \ar[d]^-{*\on{-push}} \ar[r] &
	\DMod(V)_{\mathscr{T}'}\boxtimes 1 \ar[d]^-{*\on{-push}} \\
	 \mathcal{KL}(H)_{\mathscr{L},\on{co}}  
	& \mathscr{L}_{\widehat{\on{Hecke}}_{H} } \on{-mod}( \mathcal{I}nd\mathcal{C}oh( U ) )  \ar[l] \ar[r] &
	\DMod(U)_{\mathscr{T}}\boxtimes 1.
	}
	$$
Passing to left adjoints, we obtain
$$
	\xymatrix{
	\mathcal{KL}(H')_{\mathscr{L}',\on{co}}   \ar[d]^-{\on{C}^{\frac{\infty}{2}}_{\on{co}}} \ar[r]
	& \mathscr{L}'_{\widehat{\on{Hecke}}_{H'} } \on{-mod}( \mathcal{I}nd\mathcal{C}oh( V ) )  \ar[d]^-{*\on{-push}} \ar[r] &
	\DMod(V)_{\mathscr{T}'}\boxtimes 1 \ar[d]^-{*\on{-push}} \\
	 \mathcal{KL}(H)_{\mathscr{L},\on{co}}  \ar[r] \ar@{=>}[ru]
	& \mathscr{L}_{\widehat{\on{Hecke}}_{H} } \on{-mod}( \mathcal{I}nd\mathcal{C}oh( U ) )  \ar[r] &
	\DMod(U)_{\mathscr{T}}\boxtimes 1.
	}
	$$
The outer square is just
\begin{equation} \label{eqn-natural-transformation-loc-coh}
	\xymatrix{
	\mathcal{KL}(H')_{\mathscr{L}',\on{co}}   \ar[d]^-{\on{C}^{\frac{\infty}{2}}_{\on{co}}} \ar[r]^-{\on{Loc}_{V,\on{co}}}
	&
	\DMod(V)_{\mathscr{T}'}\boxtimes 1 \ar[d]^-{*\on{-push}} \\
	 \mathcal{KL}(H)_{\mathscr{L},\on{co}}   \ar[r]^-{\on{Loc}_{U,\on{co}}} \ar@{=>}[ru]
	&
	\DMod(U)_{\mathscr{T}}\boxtimes 1.
	}
	\end{equation}
\end{construction}

This $2$-morphism is \emph{not} invertible in general. However, we are going to prove the following result, which implies Theorem \ref{thm-localization-KL}(4):

\begin{lemma} \label{lem-loc-compatible-*-push}
Up to chiral homology, the digram (\ref{eqn-natural-transformation-loc-coh}) commutes for \emph{factorization objects}. More precisly, for any unital factorization algebra $A'$ in $ \mathcal{KL}(H')_{\mathscr{L}',\on{co}}$, the natural transformation
\begin{equation*} 
	\xymatrix{
	A'\on{-FactMod}_{\on{Ran}_{x,\on{un,dR}}}   \ar[d]^-{\on{C}^{\frac{\infty}{2}}_{\on{co}}} \ar[r]^-{[\on{Loc}]^{\on{ch}}_{x,\on{co}}}
	&
	\DMod(\on{Bun}_{H'})_{\mathscr{T}'} \ar[d]^-{*\on{-push}} \\
	 	\on{C}^{\frac{\infty}{2}}_{\on{co}}(A')\on{-FactMod}_{\on{Ran}_{x,\on{un,dR}}}  \ar[r]^-{[\on{Loc}]^{\on{ch}}_{x,\on{co}}}\ar@{=>}[ru]
	&
	\DMod(\on{Bun}_H)_{\mathscr{T}}.
	}
	\end{equation*}
induced by it is invertible.
\end{lemma}

\proof Let $\mathbb{V}'_{\on{co}}$ be the unit factorization algebra of $ \mathcal{KL}(H')_{\mathscr{L}',\on{co}}$. By Lemma \ref{lem-loc-compatible-with-prop}, we can reduce the lemma to the following result.

\qed

\begin{lemma} \label{lem-loc-compatible-*-push-prop}
The natural transformation
\begin{equation} \label{eqn-lem-loc-compatible-*-push-prop}
	\xymatrix{
	\mathbb{V}'_{\on{co}}\on{-FactMod}_{ \on{Arrow}(\on{Ran}_{x,\on{un,dR}})}   \ar[d]^-{\on{C}^{\frac{\infty}{2}}_{\on{co}}} \ar[r]^-{\on{Loc}_{\on{co}}}
	&
	\DMod(\on{Bun}_{H'})_{\mathscr{T}'} \otimes \DMod( \on{Arrow}(\on{Ran}_{x,\on{un,dR}}) ) \ar[d]^-{*\on{-push}} \\
	 	\on{C}^{\frac{\infty}{2}}_{\on{co}}(\mathbb{V}'_{\on{co}})\on{-FactMod}_{ \on{Arrow}(\on{Ran}_{x,\on{un,dR}})}   \ar[r]^-{\on{Loc}_{\on{co}}}\ar@{=>}[ru]
	&
	\DMod(\on{Bun}_H)_{\mathscr{T}} \otimes \DMod( \on{Arrow}(\on{Ran}_{x,\on{un,dR}}) ).
	}
	\end{equation}
is invertible after composing with
$$  p_{s,!}:  \DMod( \on{Arrow}(\on{Ran}_{x,\on{un,dR}}) )\to  \DMod( \on{Ran}_{x,\on{un,dR}}). $$
\end{lemma}

\begin{notation} In the proof of the lemma, we use the following shorthands:
\begin{itemize}
	\item $\pi: \on{Ran}_{x,\on{un}} \to \on{Ran}_{\on{un}}$ is the obvious preojection;
	\item $\mathcal{R}ep:= \mathcal{R}ep(\mathcal{L}^+ H),\; \mathcal{R}ep':= \mathcal{R}ep(\mathcal{L}^+ H'),\; \mathcal{KL}_{\on{co}}:= \mathcal{KL}(H)_{\mathscr{L},\on{co}},\; \mathcal{KL}'_{\on{co}}:= \mathcal{KL}(H')_{\mathscr{L}',\on{co}}$ are the unital factorization categories over $\on{Ran}_{\on{un}}$;
	\item $\mathcal{I}nd\mathcal{C}oh: \IndCoh(\on{Bun}_H) \boxtimes \mathbf{1},\;\mathcal{I}nd\mathcal{C}oh': \IndCoh(\on{Bun}_H') \boxtimes \mathbf{1} $
	are the crystal of categories over $\on{Ran}_{\on{un}}$, similarly for $\mathcal{DM}od$ and $\mathcal{DM}od'$.
\end{itemize}
\end{notation}

\proof Similar to the proof of Lemma \ref{lem-loc-compatible-oblv-unit}, we only need to show (\ref{eqn-lem-loc-compatible-*-push-prop}) is invertible after composing with
$$ \mathbf{\Gamma}(\pi_\dR, \mathcal{R}ep') \xrightarrow{\on{ind}} \mathbf{\Gamma}(\pi_\dR, \mathcal{KL}')   \xrightarrow{\on{ins}_{\on{unit}}} \mathbb{V}'_{\on{co}}\on{-FactMod}_{ \on{Arrow}(\on{Ran}_{x,\on{un,dR}})}. $$

We have commutative squares
$$
\xymatrix{
	\mathcal{R}ep' \ar[r]^-{\on{Loc}^+_{\on{co}}} \ar[d]^-{\on{ind}}
	&
	\mathcal{I}nd\mathcal{C}oh'  \ar[d]^-{\on{ind}} &
	\mathcal{R}ep \ar[r]^-{\on{Loc}^+_{\on{co}}} \ar[d]^-{\on{ind}}
	&
	\mathcal{I}nd\mathcal{C}oh   \ar[d]^-{\on{ind}} 
	& \mathcal{R}ep' \ar[r]^-{\on{ind}} \ar[d]^-{\on{inv}}
	&
	\mathcal{KL}'_{\on{co}}    \ar[d]^-{\on{C}^{\frac{\infty}{2}}_{\on{co}}} \\
	\mathcal{KL}'_{\on{co}}   \ar[r]^-{\on{Loc}_{\on{co}}}
	&
	\mathcal{D}mod' , &
		\mathcal{KL}_{\on{co}}   \ar[r]^-{\on{Loc}_{\on{co}}}
	&
	\mathcal{D}mod ,
		& \mathcal{R}ep \ar[r]^-{\on{ind}} 
	&
	\mathcal{KL}_{\on{co}}, 
}
$$
where the last one is induced by twisted $*$-pushforward morphisms, and the other two squares are just (the co-version of) (\ref{eqn-loc-full-sph-3}). Note that the horizontal morphisms are strictly unital. Hence we only need to show
\begin{equation} \label{eqn-proof-lem-loc-compatible-*-push-prop-1}
\xymatrix{
 \mathbf{\Gamma}((\pi\circ p_t)_\dR, \mathcal{R}ep')  \ar[d]^-{\on{inv}} \ar[r]^-{\on{Loc}^+_{\on{co}}} &
\mathbf{\Gamma}((\pi\circ p_t)_\dR, \mathcal{I}nd\mathcal{C}oh')  \ar[d]^-{*\on{-push}}
\\
 \mathbf{\Gamma}((\pi\circ p_t)_\dR, \mathcal{R}ep) \ar[r]^-{\on{Loc}^+_{\on{co}}} \ar@{=>}[ru] &
\mathbf{\Gamma}((\pi\circ p_t)_\dR, \mathcal{I}nd\mathcal{C}oh)
}
\end{equation}
is invertible after precomposing with
$$ \on{ins}_{\on{unit}}:\mathbf{\Gamma}(\pi_\dR, \mathcal{R}ep')\to  \mathbf{\Gamma}((\pi\circ p_t)_\dR, \mathcal{R}ep')$$
and composing with
$$p_{s,!} :\mathbf{\Gamma}((\pi\circ p_t)_\dR, \mathcal{I}nd\mathcal{C}oh)\to\mathbf{\Gamma}(\pi_\dR, \mathcal{I}nd\mathcal{C}oh).$$
Since $H'\to H$ is a surjection with unipotent kernel, the image of the restriction functor
$$ \on{res}:  \mathbf{\Gamma}(\pi_\dR, \mathcal{R}ep) \to \mathbf{\Gamma}(\pi_\dR, \mathcal{R}ep')$$
generates the target. Hence we can further precompose with this functor. Note that $$\mathbf{\Gamma}(\pi_\dR, \mathcal{R}ep)  \xrightarrow{\on{res}} \mathbf{\Gamma}(\pi_\dR, \mathcal{R}ep') \xrightarrow{ \on{ins}_{\on{unit}}}  \mathbf{\Gamma}((\pi\circ p_t)_\dR, \mathcal{R}ep')$$
is equivalent to
$$\mathbf{\Gamma}(\pi_\dR, \mathcal{R}ep) \xrightarrow{ \on{ins}_{\on{unit}}}  \mathbf{\Gamma}((\pi\circ p_t)_\dR, \mathcal{R}ep) \xrightarrow{\on{res}}  \mathbf{\Gamma}((\pi\circ p_t)_\dR, \mathcal{R}ep')$$
because $\on{res}: \mathcal{R}ep\to \mathcal{R}ep'$ is strictly unital.

Recall the pointwise-tensor products induce a symmetric monoidal structure on $\mathcal{R}ep$, and the monoidal unit of $\mathbf{\Gamma}( \on{Id}_\dR, \mathcal{R}ep)$ (which is the trivial representation) is also the factorization unit. Note that all the relevant prestacks in the above statement
$$ [\mathbb{B}\mathcal{L}^+H]_{\on{Ran}_{\on{un}}} ,\; [\mathbb{B}\mathcal{L}^+H']_{\on{Ran}_{\on{un}}}, \on{Bun}_H \times \on{Ran}_{\on{un}},\; \on{Bun}_{H'} \times \on{Ran}_{\on{un}}$$
are defined over $[\mathbb{B}\mathcal{L}^+H]_{\on{Ran}_{\on{un}}}$. Using the base-change isomorphisms and their corollary, the projection formulae, we can reduce to show (\ref{eqn-proof-lem-loc-compatible-*-push-prop-1}) is invertible when precomposing with $\on{res} \circ \on{ins}_{\on{unit}}$, composing with $p_{s,!}$, and evaluating at the monoidal unit of $\mathbf{\Gamma}(\pi_\dR, \mathcal{R}ep)$. In other words, we only need to show (\ref{eqn-proof-lem-loc-compatible-*-push-prop-1}) is invertible when evaluating at the monoidal unit of $ \mathbf{\Gamma}((\pi\circ p_t)_\dR, \mathcal{R}ep')$ and composing with $p_{s,!}$. 

Using Propogation Lemma for Chiral homology (see Lemma \ref{lem-AA-lemma}), we can further reduce to the following lemma.

\qed

\begin{lemma} \label{lem-conformal-block-chiral-homology}
The Beck-Chevalley natural transformation
$$
\xymatrix{
[\mathcal{R}ep(\mathcal{L}^+ H')]_{\on{Ran}_{\on{un,dR}}} \ar[r]^-{ \on{Loc}^+_{\on{co}}} \ar[d]^-{\on{inv}} &
\IndCoh(\on{Bun}_{H'})\otimes \DMod(\on{Ran}_{\on{un}} )\ar[d]^-{*\on{-push}} \\
[\mathcal{R}ep(\mathcal{L}^+ H)]_{\on{Ran}_{\on{un,dR}}} \ar[r]^-{ \on{Loc}^+_{\on{co}}} \ar@{=>}[ru] &
\IndCoh(\on{Bun}_{H})\otimes \DMod(\on{Ran}_{\on{un}} )
}
$$
is invertible when evaluating at the monoidal unit and taking chiral homology, i.e., composing with
$$ \on{C}_{\on{un}}^{\on{ch}}: \IndCoh(\on{Bun}_{H})\otimes \DMod(\on{Ran}_{\on{un}} )\to  \IndCoh(\on{Bun}_{H}). $$
The claim remains correct if we replace $\on{Loc}^+_{\on{co}}$ by $\on{Loc}$ and the $*$-pushforward functor by the $?$-pushforward functor.

\end{lemma}

\proof Recall $\on{Loc}^+_{\on{co}}$ is defined as the $*$-pullback functor for $\mathcal{I}nd\mathcal{C}oh_{*,\ren}$-theory along
$$ \on{Bun}_H \times \Ran_{\on{un,dR}} \to [\mathbb{B} \mathcal{L}^+ H)]_{\Ran_{\on{un,dR}}} $$
while $\on{Loc}^+$ is defined as the $!$-pullback functor for $\mathcal{I}nd\mathcal{C}oh^{!,\ren}$-theory along the same map. It follows that
$$\on{Loc}^+ \simeq  K_{\on{Bun}_H}\otimes \on{Loc}^+_{\on{co}}.$$
Then it is clear that the claims for $\on{Loc}$ and $\on{Loc}_{\on{co}}$ are equivalent. We will freely switch between these two versions.

Let $q: \on{Bun}_{H'} \to \on{Bun}_{H}$ be the projection. Let $I:=\on{inv}(\on{unit})\in \mathcal{R}ep(\mathcal{L}^+ H)$ be the unital factorization algebra. We want to show
$$\on{Loc}^+ (I) \to  q_?(\omega_{\on{Bun}_{H'} }) \boxtimes \omega_{ \on{Ran}_{\on{un}} } $$
is sent to an isomorphism by $\on{C}_{\on{un}}^{\on{ch}}$.

We need to use the theory of commutative factorization categories and commutative factorization algebras. We first review what these mean.

Recall the factorization structure on $\mathcal{R}ep(\mathcal{L}^+ H)$ is induced by the factorization structure on $[\mathcal{L}^+ H]_{\on{Ran}_{\on{un,dR}}}$. Moreover, the latter factorization structure on $[\mathcal{L}^+ H]_{\on{Ran}_{\on{un,dR}}}$ can be upgraded to a \emph{cocommutative} one. Roughly speaking, this means the isomorphism
$$[\mathcal{L}^+ H]_{ (\on{Ran}_{\on{un,dR}}\times\on{Ran}_{\on{un,dR}})_{\on{disj}} }  \simeq ([\mathcal{L}^+ H]_{ \on{Ran}_{\on{un,dR}}} \times [\mathcal{L}^+ H]_{ \on{Ran}_{\on{un,dR}}} )_{\on{disj}}$$
can be extended to a morphism
$$ [\mathcal{L}^+ H]_{ \on{Ran}_{\on{un,dR}}\times\on{Ran}_{\on{un,dR}} }  \to [\mathcal{L}^+ H]_{ \on{Ran}_{\on{un,dR}}} \times [\mathcal{L}^+ H]_{ \on{Ran}_{\on{un,dR}}}  $$
defined over $\on{Ran}_{\on{un,dR}}\times\on{Ran}_{\on{un,dR}}$ subject to higher compatibilities (which are actually structures). This implies the factorization structure on $\mathcal{R}ep(\mathcal{L}^+ H)$ can be upgraded to a \emph{commutative} one. This means the equivalence
$$   (\mathcal{R}ep(\mathcal{L}^+ H) \boxtimes \mathcal{R}ep(\mathcal{L}^+ H))|_{ (\on{Ran}_{\on{un,dR}}\times\on{Ran}_{\on{un,dR}})_{\on{disj}} } \simeq \mathcal{R}ep(\mathcal{L}^+ H) |_{ (\on{Ran}_{\on{un,dR}}\times\on{Ran}_{\on{un,dR}})_{\on{disj}} } $$
can be extended to a morphism
$$ \on{mult}:  (\mathcal{R}ep(\mathcal{L}^+ H) \boxtimes \mathcal{R}ep(\mathcal{L}^+ H))|_{ \on{Ran}_{\on{un,dR}}\times\on{Ran}_{\on{un,dR}} } \to \mathcal{R}ep(\mathcal{L}^+ H) |_{ \on{Ran}_{\on{un,dR}}\times\on{Ran}_{\on{un,dR}}} $$
between crystal of categories on $\on{Ran}_{\on{un,dR}}\times\on{Ran}_{\on{un,dR}}$ subject to higher compatibilities. 

The factorization functor $\on{inv}: \mathcal{R}ep(\mathcal{L}^+ H')\to \mathcal{R}ep(\mathcal{L}^+ H)$
can be upgrated to a \emph{right-lax commutative} factorization functor. Roughly speaking, this means we have a $2$-morphism
$$
\xymatrix{
	(\mathcal{R}ep(\mathcal{L}^+ H') \boxtimes \mathcal{R}ep(\mathcal{L}^+ H'))|_{ \on{Ran}_{\on{un,dR}}\times\on{Ran}_{\on{un,dR}} } 
	\ar[r]^-{\on{mult}} \ar[d]^-{ \on{inv}\boxtimes \on{inv} } &
	\mathcal{R}ep(\mathcal{L}^+ H') |_{ \on{Ran}_{\on{un,dR}}\times\on{Ran}_{\on{un,dR}}} \ar[d]^-{ \on{inv} }\\
	(\mathcal{R}ep(\mathcal{L}^+ H) \boxtimes \mathcal{R}ep(\mathcal{L}^+ H))|_{ \on{Ran}_{\on{un,dR}}\times\on{Ran}_{\on{un,dR}} } 
	\ar[r]^-{\on{mult}} \ar@{=>}[ru] &
	\mathcal{R}ep(\mathcal{L}^+ H) |_{ \on{Ran}_{\on{un,dR}}\times\on{Ran}_{\on{un,dR}}}
}
$$
subject to higher compatibilities. Namely, in the $\mathcal{I}nd\mathcal{C}oh^{!,\ren}$-theory, the horizontal morphisms are induced by $!$-pullback functors, the vertical morphisms are induced by $?$-pushforward fucntors (which are the right adjoints of $!$-pullback functors), and the $2$-morphism is induced by the Bech-Chevalley natural transformations. It follows that $I=\on{inv}(\on{unit})$ can be upgraded to a \emph{commutative} (unital) factorization algebra. This means we have a morphism
$$ \on{mult}( I\boxtimes I ) \to I|_{ \on{Ran}_{\on{un,dR}}\times\on{Ran}_{\on{un,dR}}} $$
in $[\mathcal{R}ep(\mathcal{L}^+ H)]_{ \on{Ran}_{\on{un,dR}}\times\on{Ran}_{\on{un,dR}}}$ subject to higher compatibilities.

Also, the commutative factorization algebra structure on $I$ can upgrade $\on{Loc}^+(I)$ to a commutative factorization algebra inside $\IndCoh(\on{Bun}_H)\boxtimes \mathbf{1}$, which is a $\on{Bun}_H$-family of commutative factorization categories. This means we have a canonical morphism
$$ \on{Loc}^+(I) \otimes^!_{\on{Bun}_H} \on{Loc}^+(I) \to \on{Loc}^+(I)|_{ \on{Ran}_{\on{un,dR}}\times\on{Ran}_{\on{un,dR}} } $$
subject to higher compatibilities.

It is clear that $q_?(\omega_{\on{Bun}_{H'} }) \boxtimes \omega_{ \on{Ran}_{\on{un}} }$ is also a commutative factorization algebra inside $\IndCoh(\on{Bun}_H)\boxtimes \mathbf{1}$ and the map
\begin{equation}\label{eqn-proof-lem-conformal-block-chiral-homology-1}
 \on{Loc}^+ (I) \to  q_?(\omega_{\on{Bun}_{H'} }) \boxtimes \omega_{ \on{Ran}_{\on{un}} } 
\end{equation}
is a morphism between commutative factorization algebras. 

By Proposition \ref{prop:comm-fact-main-diagonal-same}, the restriction functor
$$  [\mathcal{R}ep(\mathcal{L}^+ H)]_{ \on{Ran}_{\on{un,dR}} } \to [\mathcal{R}ep(\mathcal{L}^+ H) ]_{X_\dR} $$
induces an equivalence
$$ \on{CommFactAlg}_{\on{un}}(\mathcal{R}ep(\mathcal{L}^+ H)  ) \simeq \on{CommAlg}_{\on{un}}( [\mathcal{R}ep(\mathcal{L}^+ H) ]_{X_\dR} ).$$
Similarly we have
$$ \Delta^!: \on{CommFactAlg}_{\on{un}}(\IndCoh(\on{Bun}_H)\boxtimes \mathbf{1}  ) \simeq \on{CommAlg}_{\on{un}}( \IndCoh(\on{Bun}_H\times X_\dR) ).$$
Also, the symmtric monoidal functor 
$$\on{pr}_1^!:\IndCoh(\on{Bun}_H)\to  \IndCoh(\on{Bun}_H\times X_\dR) $$ induces a functor
$$\on{pr}_1^{!,\on{CommAlg}}: \on{CommAlg}_{\on{un}}( \IndCoh(\on{Bun}_H) ) \to \on{CommAlg}_{\on{un}}( \IndCoh(\on{Bun}_H\times X_\dR) ) $$
such that the corresponding functor
$$ \on{CommAlg}_{\on{un}}( \IndCoh(\on{Bun}_H) ) \to  \on{CommFactAlg}_{\on{un}}(\IndCoh(\on{Bun}_H)\boxtimes \mathbf{1}  ) $$
sends $\mathcal{M}$ to $\mathcal{M}\boxtimes \omega_{\on{Ran}_{\on{un}}}$. Passing to left adjoints, we obtain a commutative diagram
$$
\xymatrix{
\on{CommFactAlg}_{\on{un}}(\IndCoh(\on{Bun}_H)\boxtimes \mathbf{1}  ) \ar[r]_-{\Delta^!}^-\simeq\ar[d]^-{ \on{C}^{\on{ch}}_{\on{un}} }  &
\on{CommAlg}_{\on{un}}( \IndCoh(\on{Bun}_H\times X_\dR) )\ar[ld]^-{ ( \on{pr}_1^{!,\on{CommAlg}})^L } \\
\on{CommAlg}_{\on{un}}( \IndCoh(\on{Bun}_H) ).
}
$$
Hence we only need to show (\ref{eqn-proof-lem-conformal-block-chiral-homology-1}) is sent to an isomorphism by $( \on{pr}_1^{!,\on{CommAlg}})^L\circ \Delta^!$.

Unwinding the definitions, we have reduce the problem to the following claim. Consider the Beck-Chevalley natural transformation
$$
\xymatrix{
[\mathcal{R}ep(\mathcal{L}^+ H')]_{X_\dR} \ar[r]^-{ \on{Loc}^+_{\on{co}}} \ar[d]^-{\on{inv}} &
\IndCoh(\on{Bun}_{H'} \times X_\dR)\ar[d]^-{*\on{-push}} \\
[\mathcal{R}ep(\mathcal{L}^+ H)]_{X_\dR} \ar[r]^-{ \on{Loc}^+_{\on{co}}} \ar@{=>}[ru] &
\IndCoh(\on{Bun}_{H}\times X_\dR )
}
$$
and its value on the monoidal unit, which is a morphism in $\on{CommAlg}_{\on{un}}( \IndCoh(\on{Bun}_H\times X_\dR) )$. We only need to show this morphism is sent to an isomorphism by $( \on{pr}_1^{!,\on{CommAlg}})^L$. This claim can be rewritten using $\on{QCoh}$-theory as follows. The value of\footnote{For a symmetric monoidal DG category $\mathcal{C}$ equipped with a compatible $t$-structure, the category $\on{CommAlg}_{\on{un}}^{\on{cc}}(\mathcal{C})$ is the full subcategory of $\on{CommAlg}_{\on{un}}(\mathcal{C})$ containing those objects $A$ that are \emph{coconnective}, i.e., such that $A\in \mathcal{C}^{\ge 0}$ and $H^0(A)\in \mathcal{C}^{\heartsuit}$ is the monoidal unit.}
$$
\xymatrix{
\on{CommAlg}_{\on{un}}^{\on{cc}}(\on{QCoh}( [\mathbb{B}\mathcal{L}^+ H']_{X_\dR} )) \ar[r]^-{ *\on{-pull}} \ar[d]^-{*\on{-push}} &
\on{CommAlg}_{\on{un}}^{\on{cc}}(\QCoh(\on{Bun}_{H'} \times X_\dR))\ar[d]^-{*\on{-push}} \\
\on{CommAlg}_{\on{un}}^{\on{cc}}(\on{QCoh}( [\mathbb{B}\mathcal{L}^+ H]_{X_\dR} )) \ar[r]^-{ *\on{-pull}} \ar@{=>}[ru] &
\on{CommAlg}_{\on{un}}^{\on{cc}}(\QCoh(\on{Bun}_{H}\times X_\dR )) 
}
$$
at the monoidal unit is sent to an isomorphism by
$$ \on{CommAlg}_{\on{un}}^{\on{cc}}(\QCoh(\on{Bun}_{H}\times X_\dR )) \subset \on{CommAlg}_{\on{un}}(\QCoh(\on{Bun}_{H}\times X_\dR )) \xrightarrow{ (\on{pr}_1^{*,\on{CommAlg}})^L }  \on{CommAlg}_{\on{un}}(\QCoh(\on{Bun}_{H})). $$

Let $O_1 \to O_2$ be the value of the above natural transformation at the monoidal unit. Note that
$$(\on{pr}_1^{*,\on{CommAlg}})^L(O_2) \simeq q_*( O_{\on{Bun}_{H'}} ),$$
where $q: \on{Bun}_{H'} \to \on{Bun}_H$ is the projection map. Indeed, this can be seen by translating back to the chiral homology description for $(\on{pr}_1^{*,\on{CommAlg}})^L$. In particular, $(\on{pr}_1^{*,\on{CommAlg}})^L(O_2)$ is coconnective.

We claim $(\on{pr}_1^{*,\on{CommAlg}})^L(O_1)$ is also coconnective. We first finish the proof assuming this claim. We only need to show $O_1\to O_2$ is sent to an isomorphism by the partially defined left adjoint
$$L:=(\on{pr}_1^{*,\on{CommAlg}^{\on{cc}}})^L:\on{CommAlg}_{\on{un}}^{\on{cc}}(\QCoh(\on{Bun}_{H}\times X_\dR )) \to \on{CommAlg}_{\on{un}}^{\on{cc}}(\QCoh(\on{Bun}_{H} )).$$
Using the assumption that $H'\to H$ is surjective with unipotent kernel, one can check that the maps
$$  [\mathbb{B}\mathcal{L}^+ H']_{X_\dR}\to [\mathbb{B}\mathcal{L}^+ H]_{X_\dR},\;\on{Bun}_{H'}\to \on{Bun}_{H} $$
are \emph{co-affine}\footnote{A map $Y\to Z$ is co-affine if for any affine test scheme $S\to Z$, the fiber product $S\times_Z Y$ is co-affine in the sense of \cite{DAGVIII}.}. Hence we have
$$  [\mathbb{B}\mathcal{L}^+ H']_{X_\dR}\times_{[\mathbb{B}\mathcal{L}^+ H]_{X_\dR}  }  (\on{Bun}_{H}\times X_\dR) \simeq \on{cSpec}(O_1)  ,\; \on{Bun}_{H'}\times X_\dR \simeq \on{cSpec}(O_2),$$
where 
$$\on{cSpec}: \on{CommAlg}_{\on{un}}^{\on{cc}}(\QCoh(\on{Bun}_{H}\times X_\dR )) \to \on{PreStk}_{ / \on{Bun}_{H}\times X_\dR}
$$
is the (relative) co-affine spectrum functor (see \cite[Section 4.4]{DAGVIII}), it is charaterized by the functorial equivalences
$$ \on{Maps}_{\on{Bun}_{H}\times X_\dR}( \on{Spec}(A), \on{cSpec}(O) ) \simeq \on{Maps}_{ \on{CommAlg}_{\on{un}}(\QCoh(\on{Bun}_{H}\times X_\dR )) }(O,A) $$
for any \emph{connective} $A\in \on{CommAlg}_{\on{un}}^{\on{cn}}(\QCoh(\on{Bun}_{H}\times X_\dR ))$. Consider the Weyl retriction functor
$$W: \on{PreStk}_{ / \on{Bun}_{H}\times X_\dR} \to \on{PreStk}_{ / \on{Bun}_{H}}.$$
For any affine test scheme $Y$ over $\on{Bun}_H$, we have
\begin{eqnarray*}
& & \on{Maps}_{\on{Bun}_H}( Y,W([\mathbb{B}\mathcal{L}^+ H']_{X_\dR}\times_{[\mathbb{B}\mathcal{L}^+ H]_{X_\dR}  }  (\on{Bun}_{H}\times X_\dR) )  ) \\
&\simeq & \on{Maps}_{\on{Bun}_H\times X_\dR}( Y\times X_\dR,[\mathbb{B}\mathcal{L}^+ H']_{X_\dR}\times_{[\mathbb{B}\mathcal{L}^+ H]_{X_\dR}  }  (\on{Bun}_{H}\times X_\dR) )   \\
&\simeq & \on{Maps}_{[\mathbb{B}\mathcal{L}^+ H]_{X_\dR} }( Y\times X_\dR,[\mathbb{B}\mathcal{L}^+ H']_{X_\dR}  ) \\
& \simeq & \on{Maps}_{\mathbb{B}H\times X }( Y\times X,\mathbb{B}H'\times X )\\
& \simeq &  \on{Maps}_{\on{Bun}_H}( Y, \on{Bun}_{H'}),
\end{eqnarray*}
where we use the fact that the Weyl restriction functor $\on{PreStk}_{/X} \to \on{PreStk}_{/X_\dR}$ sends $\mathbb{B}H\times X$ to $[\mathbb{B}\mathcal{L}^+ H]_{X_\dR}$. Hence
$$W([\mathbb{B}\mathcal{L}^+ H']_{X_\dR}\times_{[\mathbb{B}\mathcal{L}^+ H]_{X_\dR}  }  (\on{Bun}_{H}\times X_\dR) ) \simeq \on{Bun}_{H'}.$$
Since $\on{cSpec}$ is fully faithful (see \cite[Theorem 4.4.1]{DAGVIII}), we obtain
$$L(O_1) \simeq q_*( O_{\on{Bun}_{H'}} ) \simeq L(O_2)$$
as desired.

It remains to show $(\on{pr}_1^{*,\on{CommAlg}})^L(O_1)$ is coconnective. Let $\mathfrak{n}$ be the Lie algebra of the kernel $H'\to H$. Recall we have a $\on{Lie}^*$-algebra $\mathfrak{n}_D:= \mathfrak{n}\otimes D_X \in \DMod(X)$. Similarly, we have a $\on{Bun}_H$-family of $\on{Lie}^*$-algebras $\widetilde{\mathfrak{n}_D}\in \on{QCoh}(\on{Bun}_H) \otimes \DMod(X)$ such that its restriction at a $H$-torsor $\mathcal{P}\in \on{Bun}_H$ is the $\on{Lie}^*$-algebra $(\mathfrak{n}_D)_{\mathcal{P}}$ twisted by $\mathcal{P}$ (defined via the adjoint action of $H$ on $\mathfrak{n}$). Let $C(\widetilde{\mathfrak{n}_D} )\in \on{QCoh}(\on{Bun}_H) \otimes \DMod(X)$ be the $\on{Bun}_H$-family of commutative algebras that is Koszul dual to $\widetilde{\mathfrak{n}_D} $ (see Proposition \ref{prop:chiral-koszul-duality}). By definition, $O_1\simeq C(\widetilde{\mathfrak{n}_D} )$. By \cite[Proposition 6.4.4]{francis2012chiral}, we have
$$(\on{pr}_1^{*,\on{CommAlg}})^L(C(\widetilde{\mathfrak{n}_D} )) \simeq C \circ \Gamma(X_\dR, \widetilde{\mathfrak{n}_D}),  $$
where $\Gamma(X_\dR, \widetilde{\mathfrak{n}_D}) \in \on{QCoh}(\on{Bun}_H)$ is a $\on{Bun}_H$-family of $\on{Lie}$-algebras, i.e., a Lie-algebra object in $\on{QCoh}(\on{Bun}_H)$ and $C \circ \Gamma(X_\dR, \widetilde{\mathfrak{n}_D})$ is its Koszul dual. Then we are done because $C \circ \Gamma(X_\dR, \widetilde{\mathfrak{n}_D})$ is clearly coconnective.

\qed

\section{Toric Fundamental Local Equivalence}
\label{sect:tori}
The starting point is the following classical result due to \cite[Theorem 1.6.6]{contou2013jacobienne}:
\begin{proposition}
	For any classical Noetherian scheme $S_0$, we have an isomorphism of $S$-group schemes
	\[\underline{\Hom}_{\tx{Grp}/S_0}(\mc{L} \mb{G}_m, \underline{\mb{G}_m}) \simeq \mc{L} \mb{G}_m\]
	where $\underline{\Hom}_{\tx{Grp}/S_0}$ denotes the relative group homomorphism prestack, $\mc{L} \mb{G}_m$ is the relative loop space over $S_0$, and $\underline{\mb{G}_m})$ is the constant group scheme.
\end{proposition}
In particular this gives a map of group prestacks known as the \emph{Contou-Carr\`ere symbol}
\[\langle \bullet, \bullet \rangle: \mc{L} \mb{G}_m \times_{S_0} \mc{L} \mb{G}_m \to \underline{\mb{G}_m}.\]
At one point it can be described as follows (\cite{beilinson2001epsilon}). The determinant line bundle on $\mc{L} \mb{G}_m$ determines a central extension $\widetilde{\mc{L} \mb{G}_m}$, and we have, for any test affine scheme $A$,
\[\langle f, g \rangle^{-1} = \tilde{f} \cdot \tilde{g} \cdot \tilde{f}^{-1} \cdot \tilde{g}^{-1} \in \mb{G}_m(A)\]
for any $f, g \in \mb{G}_m(K)(A)$ and $\tilde{f}, \tilde{g}$ any lifting to the central extension. Since the central extension canonically splits on $\mc{L}^+ \mb{G}_m$, the pairing is canonically trivial on $\mc{L}^+ \mb{G}_m$ and thus induces another group homomorphism
\[\mc{L}^+ \mb{G}_m \times_{S_0} (\mc{L} \mb{G}_m / \mc{L}^+ \mb{G}_m) \to \mb{G}_m;\]
In fact, this pairing is perfect in the following sense (\cite[Theorem 1.4.4]{contou2013jacobienne}):
\begin{proposition}
We have an isomorphism of $S_0$-group schemes
\[\underline{\Hom}_{\tx{Grp}/S_0}(\mc{L}^+ \mb{G}_m, \underline{\mb{G}_m}) \simeq \mc{L}\mb{G}_m / \mc{L}^+ \mb{G}_m.\]
\end{proposition}

\begin{remark}
An important property of this symbol is that it is \emph{bi-multiplicative}, i.e. for every $S$-point $x, y, z \in \mc{L}\mb{G}_m(S)$ we have:
\[\langle x y, z \rangle = \langle x, z \rangle \langle y, z \rangle, \langle x, y z \rangle = \langle x, y \rangle \langle x, z \rangle.\]
\end{remark}

As observed in \cite{beilinson2006langlands}, the Contou-Carr\`ere symbol is furthermore compatible with the flat connection on the loop group and is factorizable because the determinant line bundle is; therefore we have a (non-unital) \emph{factorizable} version of the symbol:
\[[\mc{L}^+ \mb{G}_m]_{\Ran_\dR} \times_{\Ran_\dR} [\tx{Gr}_{\mb{G}_m}]_{\Ran_\dR} \to \mb{G}_m \times \Ran_\dR\]
Delooping on the first variable, we have a pairing between $\tx{Gr}_{\mb{G}_m}$ and its 1-dual, i.e. a map
\[[\mb{B} \mc{L}^+ \mb{G}_m]_{\Ran_\dR} \times_{\Ran_\dR} [\tx{Gr}_{\mb{G}_m}]_{\Ran_\dR} \to \mb{B} \mb{G}_m \times \Ran_\dR\]
which is explicitly describable as follows. An $S$-point of the LHS is given by
\[(x, \mc{P}_1, \mc{P}_2, \alpha)\]
where $x = (x_i)_{i \in I}$ is an $S$-point of $\Ran_\dR$, $\mc{P}_1$, $\mc{P}_2$ two line bundles on the adic disk $\mc{D}'_{\Gamma_x}$, and $\alpha$ a trivialization of $\mc{P}_2$ on $\mc{D}^\circ_{\Gamma_x}$. Then the pairing sends it to
\[(x, \det(\Gamma(\mc{D}'_{\Gamma_x}, \mc{P}_1), \Gamma(\mc{D}'_{\Gamma_x}, \mc{P}_1 \tensor \mc{P}_2))^{-1} \tensor \det(\Gamma(\mc{D}'_{\Gamma_x}, \mc{O}), \Gamma(\mc{D}'_{\Gamma_x}, \mc{P}_2)))\]
where $\Gamma(\mc{D}'_{\Gamma_x}, \mc{P}_1)$ sends $\mc{P}_1$ to the corresponding lattice in the Tate vector space $\Gamma(\mc{D}^\circ_{\Gamma_x}, \mc{P}_1)$, which we identify with $\Gamma(\mc{D}^\circ_{\Gamma_x}, \mc{P}_1 \tensor \mc{P}_2)$ via $\alpha$; $\det$ stands for the determinant line bundle.

The following was shown in \cite[Theorem 5.2.1]{campbell2017geometric}:
\begin{proposition}
\label{prop:local-cartier-dual}
This pairing is again perfect. More precisely, we have an isomorphism of \emph{counital} factorization spaces
\[\underline{\Hom}_{\tx{Grp}/X^I}([\mb{B} \mc{L}^+ \mb{G}_m]_{X^I}, \underline{\mb{G}_m}) \simeq [\tx{Gr}_{\mb{G}_m}]_{X^I}\]
after base-changing along $X^I \to \Ran_\dR$.
\end{proposition}
The counital structure mentioned in the proposition above comes from a \emph{corr-unital} structure of the pairing, given by a map of corr-unital factorization spaces
\[[\mb{B}\mc{L}^+ \mb{G}_m]_{\Ran_{\un, \dR}} \times_{\Ran_{\un, \dR}} \tx{Gr}_{\mb{G}_m, \Ran_{\un, \dR}} \to \mb{B} \mb{G}_m \times \Ran_{\un, \dR}\]
here the fiber product is taken as corr-unital factorization spaces (note this is really necessary: $\mb{B}\mc{L}^+ \mb{G}_m$ is \emph{co-unital} but $\tx{Gr}_{\mb{G}_m}$ is \emph{unital}). Unfolding definition, we see the corresponding unital structure, for any fixed injection $\alpha: I \to J$ of finite sets, is given by
\[\xymatrix{
 & [\mb{B} \mc{L}^+ \mb{G}_m]_{X_J} \times_{X^J} X^J \times_{X^I} [\tx{Gr}_{\mb{G}_m}]_{X^I} \ar[ld] \ar[rd] \\
([\mb{B}\mc{L}^+ \mb{G}_m]_{X^I} \times_{X^I} [\tx{Gr}_{\mb{G}_m}]_{X^I}) \times_{X^I} X^J & & [\mb{B}\mc{L}^+ \mb{G}_m]_{X^J} \times_{X^J} [\tx{Gr}_{\mb{G}_m}]_{X^J}
}\]
the leftward map is given by restriction while the rightward one uses Beauville-Laszlo. Then the two maps to $\mb{B} \mb{G}_m$ (via either leg) canonically give the same result, so we have a well-defined corr-unital map.

Using the pairing $\coweightLat \tensor_{\mb{Z}} \weightLat \to \mb{Z}$, the above discussion immediately generalizes to give a pairing
\[[\mb{B}\mc{L}^+ T]_{\Ran_{\un, \dR}} \times_{\Ran_{\un, \dR}} \tx{Gr}_{\check{T}, \Ran_{\un, \dR}} \to \mb{B} \mb{G}_m \times \Ran_{\un, \dR}\]
uniquely determined by the requirement that, if $\mc{P}_1, \mc{P}_2$ are line bundles, $\lambda: \mb{G}_m \to T, \cmu: \mb{G}_m \to \check{T}$, then the pairing on the induced $T$ (resp.\ $\check{T}$) bundles should be the pairing of $\mc{P}_1$ and $\mc{P}_2$ raise to the $\lambda(\cmu)$-th power.

\subsection{Levelwise Cartier Duality}
Recall that we have a placid presentation $\mc{L}^+T_{X^I} \simeq \limit_i \mc{L}^{(n)}T_{X^I}$ by the $n$th jet schemes, where each $\mc{L}^{(n)}T_{X^I}$ is a smooth affine $X^I$-group scheme. Consider the Cartier dual
\[Z_n := \underline{\Hom}_{\tx{Grp}/X^I}(\mc{L}^{(n)}T_{X^I}, \underline{\mb{G}_m})\]
The goal of this subsection is to prove the following statement:
\begin{proposition}
\label{prop:GrT-as-nice-as-it-gets}
Each $Z_n$ is an ind-flat ind-finite group ind-scheme, and for each $n < m$, the map $Z_n \to Z_m$ is an ind-closed embedding, and we have $\tx{Gr}_{\check{T}, X^I} \simeq \colimit_n Z_n$. In other words, this provides an ind-flat placid presentation of $\tx{Gr}_{\check{T}, X^I}$.
\end{proposition}

\begin{remark}
	An ind-flat placid presentation of $\tx{Gr}_{G, X^I}$ is expected to exist for all reductive $G$. However, we are currently unaware of such a construction.
\end{remark}

\begin{proof}
	One readily reduces to the case of $T = \mb{G}_m$. Let $r: (\mc{L}^{(n)} \mb{G}_m)_{X^I} \to X^I$ be the projection, and set
	\[\mc{O}_n := r_*(\mc{O}) \in \QCoh^{\heartsuit}(X^I);\]
	it is a Hopf algebra object in said category. It is clearly flat as an $\mc{O}_{x^I}$-module, and is in fact locally free (\cite[Corollary I.3.3.12]{Gruson1971}), so we reduce to an immediate check on geometric fibers.

More explicitly, let $\Gamma_n$ be the scheme-theoretic sum of $n$th neighborhood of the graphs in $X^I \times X$, and let $p: \Gamma_n \to X^I$ be the projection map. If we consider $\mb{G}_m \simeq \mb{A}^2 \times_{\mb{A}^1} \{1\}$ as embedded into $\mb{A}^2$ via $x \mapsto (x, x^{-1})$, then we can write
\begin{equation}
\label{eq:sym-expression}
\mc{O}_n \simeq \tx{Sym}_{\mc{O}_{X^I}}(p_*(\mc{V}_2^\vee \tensor K)) \tensor_{\tx{Sym}_{\mc{O}_{X^I}}(p_*\tensor K)} \mc{O}_{X^I}
\end{equation}
where $K := p^!(\mc{O}_{X^I})$ is the relative dualizing sheaf, and $\mc{V}_r^\vee := \underline{\Hom}_{\mc{O}_{\Gamma_n}}(\mc{O}_{\Gamma_n}^{\oplus r}, \mc{O}_{\Gamma_n}) \in \QCoh^\heartsuit(\Gamma_n)$.

The algebra map $\tx{Sym}_{\mc{O}_{X^I}}(p_* K) \to \mc{O}_{X^I}$ is induced by the trace map $p_* p^!(\mc{O}) \to \mc{O}$, and the algebra map
\[\tx{Sym}_{\mc{O}_{X^I}}(p_* K) \to \tx{Sym}_{\mc{O}_{X^I}}(p_*(\mc{V}_2^\vee \tensor K))\]
is induced by a sheaf map
\[p_* K \to p_*(\mc{V}_2^\vee \tensor K) \tensor_{\mc{O}_{X^I}} p_*(\mc{V}_2^\vee \tensor K) \to \tx{Sym}^{2}_{\mc{O}_{X^I}}(p_*(\mc{V}_2^\vee \tensor K)) \inj \tx{Sym}_{\mc{O}_{X^I}}(p_*(\mc{V}_2^\vee \tensor K))\]
where all but the first map is canonical; the first map is in turn given by a coalgebra map
\[p_*(K) \to p_*(K) \tensor_{\mc{O}_{X^I}} p_*(K)\]
followed by the two maps induced by $K \simeq \mc{V}_1^\vee \tensor K \to \mc{V}_2^\vee \tensor K$ given by projections $\mc{V}_2 \to \mc{V}_1$. To get this coalgebra map, consider the diagram
\[\xymatrix@R=0.3em{
\Gamma_n \ar^{\Delta}[rd] \ar_{p}[rdddd] \\
& \Gamma_n \times_{X^I} \Gamma_n \ar^{u}[rr] \ar^{\pi}[ddd] & & \Gamma_n \times \Gamma_n \ar^{p \times p}[ddd] \\ \\ \\
& X^I \ar^{\Delta_X}[rr] & & X^I \times X^I
}\]
Then we have
\[p_*(K) \simeq p_* p^! \mc{O}_{X^I} \simeq \pi_* \Delta_* \Delta^! \pi^! \mc{O}_{X^I} \xrightarrow{\tx{counit}} \pi_* \pi^! \mc{O}_{X^I} \simeq \pi_* \pi^! \Delta_X^* \mc{O}_{X^I \times X^I} \simeq \pi_* u^* (p \times p)^!(\mc{O}_{X^I \times X^I})\]
\[\simeq \Delta_X^* (p \times p)_* (p \times p)^!(\mc{O}_{X^I \times X^I}) \simeq p_*(K) \tensor_{\mc{O}_{X^I}} p_*(K)\]
Note that each individual term on RHS of Equation~(\ref{eq:sym-expression}) has a coalgebra structure (given by the obvious diagonal maps), which induces a coalgebra structure on the RHS that agrees with that of $\mc{O}_n$.
In particular, the surjective map $\tx{Sym}_{\mc{O}_{X^I}}(p_*(\mc{V}_2^\vee \tensor K)) \surj \mc{O}_n$ is a coalgebra map. Let
\[\tx{Sym}^{\le m}_{\mc{O}_{X^I}}(p_*(\mc{V}_2^\vee \tensor K)) \in \tx{Coh}^{\heartsuit}(X^I)\]
denote the subcoalgebra corresponding to the $(\le m)$-th symmetric power for $m \in \mb{N}$, and let $\mc{O}_{n, m}$ denote the sheaf-theoretic image within $\mc{O}_{n}$.
	\begin{lemma}
	\label{lemma:Onm-locally-free}
	The coherent sheaves $\mc{O}_{n, m}$ are \emph{locally free subcoalgebras} of $\mc{O}_n$, and the maps $\mc{O}_{n, m} \subseteq \mc{O}_{n, m'}$ are bundle maps i.e. admit local splittings.
	\end{lemma}
	
\begin{proof}
	Introduce the notation
	\[\mathscr{F} := \tx{Sym}_{\mc{O}_{X^I}}(p_*(\mc{V}_2^\vee \tensor K)),\]
	and similarly $\mathscr{F}^{\le m}$ forr $\tx{Sym}^{\le n}$. Denote the map $\mathscr{F} \to \mc{O}_n$ by $\pi$, and let $\mathscr{K} := \ker(\pi)$, so that
	\[\mc{O}_{n, m} \simeq \mathscr{F}^{\le m} / (\mathscr{F}^{\le m} \cap \mathscr{K}).\]
	First note that $\mathscr{K} \inj \mathscr{F} \xrightarrow{\Delta} \mathscr{F} \tensor \mathscr{F}$ factors through $\mathscr{K} \tensor \mathscr{K}$, so $\mathscr{K}$ is a subcoalgebra of $\mathscr{F}$. It follows that $\mathscr{F}^{\le m} \cap \mathscr{K}$ is a subcoalgebra of $\mathscr{F}^{\le m}$ and thus $\mc{O}_{n, m}$ is a subcoalgebra of $\mc{O}_n$.
	
	Since our base scheme is reduced, to show $\mc{O}_{n, m}$ is locally free it suffces to check that fiber ranks are constant. Assume the collision pattern at a point $x \in X^I$ is
	\[(i_1, \ldots, i_N) \hspace{1em} i_1 + \ldots + i_N = |I|;\]
	in other words, the fiber of $\Gamma_n$ at $x$ is a product of $N$ fat points given by $\tx{Spec}(k[t] / t^{n i_j})$ for $1 \le j \le N$. Since $\Gamma_n \to X^I$ is Gorenstein, we know that
	\[\mc{O}_{n}|_{x} \simeq \tx{Sym}(A^\vee \times A^\vee) \tensor_{\tx{Sym}(A^\vee)} k, \hspace{1em} A := \prod_{j = 1}^{N} k[t]/t^{n i_j}.\]
	Unfolding definition, we see we can rewrite it as
	\[\mc{A} := k[x_{j, c_j}, y_{j', c_{j'}}] / \mc{I} \hspace{1em} 1 \le j, j' \le k, 0 \le c_j < n i_j, 0 \le c_{j'} < n i_{j'}\]
	where the quotient ideal $\mc{I}$ is generated by the relation
	\[\sum_{z = 0}^{r_j} x_{j, z} y_{j, r_j - z} = \begin{cases}
	1 & r_j = 0 \\
	0 & \tx{otherwise}
\end{cases} \hspace{1em} \forall 1 \le j \le k, 0 \le r_j < n i_j\]
The Hilbert series (as a filtered algebra) of $\mc{A}$ is the same as that of $k[x_{j, c_j}, y_{j', c_{j'}}] / \mc{I}'$, where $\mc{I}'$ is the same as $\mc{I}$ except we take $x_{j, 0} y_{j, 0}$ instead of $(x_{j_0} y_{j, 0} - 1)$. The generators provided above are obviously a regular sequence, so the Hilbert series is given by $\frac{(1 - t^2)^{\sum_j n i_j}}{(1 - t)^{2n|I|}} = \frac{(1 - t^2)^{n|I|}}{(1 - t)^{2n|I|}} = \left(\frac{1 + t}{1 - t}\right)^{n|I|}$, which is independent of the collision pattern.
This computation also shows that the quotient $\mc{O}_{n, m'} / \mc{O}_{n, m}$ is constant rank, therefore the inclusion is a bundle map.
	\end{proof}
	
	Introduce
	\[\mc{O}^\vee_{n, m} := \underline{\Hom}_{\mc{O}_{X^I}}(\mc{O}_{n, m}, \mc{O}_{X^I})\]
	which is a commutative algebra object in $\tx{Coh}^{\heartsuit}(X^I)$, and is locally free because $\mc{O}_{n, m}$ is. From the lemma above, the connecting morphisms of $\mc{O}_{n, m}^\vee$ are all surjective; consequently, they form a Mittag-Leffler system and thus $\colimit_m \underline{\tx{Spec}}_{/X^I}(\mc{O}^\vee_{n, m})$ is an ind-affine ind-scheme.
	
	By construction, to check the ind-closed embedding statement it suffices to show that $\mc{O}_{n, m} \inj \mc{O}_{n + 1, m}$ is a bundle map, which holds as the quotient has constant rank (which can be readily deduced from the Hilbert series calculation above).	We claim that $Z_n \simeq \colimit_m \underline{\tx{Spec}}_{/X^I}(\mc{O}^\vee_{n, m})$. By Zariski descent it suffices to demonstrate on a Zariski covering $p: \tx{Spec}(A) \to X^I$. Let us write
	\[B_{A} = \Gamma(\tx{Spec}(A), p^* \mc{O}_n)) \hspace{1em} B_{m, A} = \Gamma(\tx{Spec}(A), p^* \mc{O}_{n, m})),\]
	so we have $B_A = \colimit_m B_{m, A}$ as $A$-algebras. Setting $B_{m, A}^\vee := \Hom_{A}(B_m, A)$, we have
	\[\tx{Spec}(A) \times_{X^I} \underline{\Hom}_{\tx{Grp}/X^I}(\mc{L}^{(n)}T_{X^I}, \underline{\mb{G}_m})) \simeq \underline{\Hom}_{\tx{Grp}/\tx{Spec}(A)}(\mc{L}^{(n)}T_{\tx{Spec}(A)}, \underline{\mb{G}_m})\]
	\[\tx{Spec}(A) \times_{X^I} \colimit_m \underline{\tx{Spec}}_{/X^I}(\mc{O}^\vee_{n, m}) \simeq \colimit_m \underline{\tx{Spec}}_{/\tx{Spec}(A)}(B_{m, A}^\vee)\]
It is clear that both are affine schemes over $\tx{Spec}(A)$, so it suffices to match their $\tx{Spec}(A)$-sections, i.e. we want
\[\Hom_{\tx{HopfAlg}(A\tx{-mod})}(A[x, x^{-1}], B_A) \simeq \colimit_m \Hom_{\tx{CAlg}(A\tx{-mod})}(\Hom_{A\tx{-mod}}(B_{m, A}, A), A).\]
Since the algebras in consideration are in the heart, we can rewrite the RHS as
\begin{align*}
& \colimit_m \Hom_{\tx{CAlg}(A\tx{-mod})}(\Hom_{A\tx{-mod}}(B_{m, A}, A), A) \\
\simeq & \colimit_m \tx{Eq}(\Hom_{A\tx{-mod}}(B_{m, A}^\vee, A) \rightrightarrows \Hom_{A\tx{-mod}}(B_{m, A}^\vee \tensor_A B_{m, A}^\vee, A)) \\
\simeq & \tx{Eq}\left(\colimit_m \Hom_{A\tx{-mod}}(B_{m, A}^\vee, A) \rightrightarrows \colimit_m \Hom_{A\tx{-mod}}(B_{m, A}^\vee \tensor_A B_{m, A}^\vee, A)\right) \\
\simeq & \tx{Eq}\left(\colimit_m B_{m, A} \rightrightarrows \colimit_m B_{m, A} \tensor_A B_{m, A}\right) \\
\simeq & \colimit_m \tx{Eq}\left(B_{m, A} \rightrightarrows B_{m, A} \tensor_A B_{m, A}\right) \\
\simeq & \colimit_m \Hom_{\tx{HopfAlg}(A\tx{-mod})}(A[x, x^{-1}], B_{m, A}) \\
\simeq & \Hom_{\tx{HopfAlg}(A\tx{-mod})}(A[x, x^{-1}], B_A)
\end{align*}
where in 4th line we use local freeness, in 5th line the two maps are comultiplication and diagonal, respectively, and the last line is because $A[x, x^{-1}]$ is a compact object in $\tx{HopfAlg}(A\tx{-mod})$.

Recall that by Noetherian approximation, we have, for any classical $\tx{Spec}(A_0) \to X^I$ and $\tx{Spec}(A_r)$ is any Noetherian approximation of $\tx{Spec}(A_0)$,
	\[\Hom_{\tx{Sch}^\tx{cl}_{/X^I}}(\tx{Spec}(A_0) \times_{X^I} \mc{L}^+ T_{X^I}, \underline{\mb{G}_m}) \simeq \Hom_{\tx{Pro}(\tx{Sch}^\tx{cl,ft}_{/X^I})}(\tx{``}\limit_{n, r}\tx{''} \tx{Spec}(A_r) \times_{X^I} \mc{L}^{(n)} T_{X^I}, \underline{\mb{G}_m})\]
	\[\simeq \colimit_{n, r} \Hom_{\tx{Sch}^\tx{cl,ft}_{/X^I}}(\tx{Spec}(A_r) \times_{X^I} \mc{L}^{(n)} T_{X^I}, \underline{\mb{G}_m}) \simeq \colimit_n \Hom_{\tx{Sch}^\tx{cl}_{/X^I}}(\tx{Spec}(A_0) \times_{X^I} \mc{L}^{(n)} T_{X^I}, \underline{\mb{G}_m})\]
Since for any classical scheme $S$ and any classical $S$-group schemes $G, H$, we have
\[\Hom_{\tx{Grp}_{/S}}(G, H) = \tx{Eq}(\Hom_{/S}(G, H) \rightrightarrows \Hom_{/S}(G \times G, H))\]
we also have
	\[\Hom_{\tx{Grp}(\tx{Sch}^\tx{cl}_{/X^I})}(\tx{Spec}(A_0) \times_{X^I} \mc{L}^+ T_{X^I}, \underline{\mb{G}_m}) \simeq \colimit_n \Hom_{\tx{Grp}(\tx{Sch}^\tx{cl}_{/X^I})}(\tx{Spec}(A_0) \times_{X^I} \mc{L}^{(n)} T_{X^I}, \underline{\mb{G}_m})\]
	from which it follows that
	\[\underline{\Hom}_{\tx{Grp}/X^I}(\mc{L}^+ T_{X^I}, \underline{\mb{G}_m}) \simeq \colimit_n \underline{\Hom}_{\tx{Grp}/X^I}(\mc{L}^{(n)} T_{X^I}, \underline{\mb{G}_m}).\]
	By Proposition~\ref{prop:local-cartier-dual}, it follows that $\tx{Gr}_{\check{T}, X^I} \simeq \colimit_n Z_n$.
\end{proof}

\subsection{Fourier-Mukai Transform}
Introduce the notation
\[B_n := [\mb{B} \mc{L}^{(n)} T]_{X^I}.\]
We have a pairing map as before:
\[F_n: Z_n \times_{X^I} B_n \to \mb{B} \mb{G}_m.\]
Let $p_1, p_2$ denote the two projection maps from LHS to $Z_n$ and $B_n$ respectively. Let $\mathscr{L}_n$ denote the Poincare line bundle, then we can use $\mathscr{L}_n$ as the Fourier-Mukai kernel to define
\[\tx{FM}_n:= (p_1)_*(p_2^!(\bullet) \tensor \mathscr{L}_n): \IndCoh(B_n) \to \IndCoh(Z_n),\]
where the tensor product indicates the action of $\QCoh$ on $\IndCoh$.

 The following property of the Fourier-Mukai kernel is the main goal of this section:
	\begin{proposition}
	\label{prop:convolution-cartier}
	$\Upsilon((p_1)_*(\mathscr{L}_n)) \in \IndCoh(Z_n)$ and $(p_2)_*(\Upsilon(\mathscr{L}_n)) \in \IndCoh(B_n)$ are invertible objects in the convolution monoidal structure $\star$ given by
	\[\star: \IndCoh(Z_n) \tensor \IndCoh(Z_n) \to \IndCoh(Z_n) := m_*(p_1^!(\bullet) \boxtimes p_2^!(\bullet))\]
	where $p_1, p_2: Z_n \times_{X^I} Z_n \to Z_n$ are projections and $m: Z_n \times_{X^I} Z_n \to Z_n$ is the multiplication map (and similarly defined for $B_n$).
	\end{proposition}
	
\begin{corollary}
	Each $\tx{FM}_n$ is an equivalence of categories.
\end{corollary}
\begin{proof}
	The proof is standard and can be found in e.g. \cite[Section 5.1]{Elliptic-I}, using Proposition~\ref{prop:convolution-cartier} as the main input. The necessary modifications are as follows:
	\begin{itemize}
	\item Replace the use of biextensions in 4.2.2 of \emph{loc.cit.} with the bi-multiplicativity of $F_n$, which follows from definition;
	\item For the inverse Fourier transformation, use $(p_2)_!(p_1^{*, \tx{IndCoh}}(\bullet) \tensor \mathscr{L}_n)$ instead. (This accounts for the asymmetry of the statement above.) Tracing through the proof one sees that for the reverse direction, one needs the projection formula
	\[f_*(\mc{E}) \tensor \mc{F} \simeq f_*(\mc{E} \tensor f^{*, \IndCoh}(\mc{F}))\]
	for $f: Z_n \times_{X^I} B_n \times_{X^I} Z_n \to Z_n \times_{X^I} Z_n$, $\mc{E} \in \QCoh(Z_n \times_{X^I} B_n \times_{X^I} Z_n)$, as well as the QCoh base change formula for the next diagram
	\[\xymatrix{
	Z_n \times_{X^I} B_n \times_{X^I} Z_n \ar[d] \ar[r] & Z_n \times_{X^I} B_n \ar[d] \\
	Z_n \times_{X^I} Z_n \ar[r] & Z_n
	}\]
	both of which follow from the fact that $B_n$ is a QCA stack of bounded Tor dimension.
\end{itemize}
The rest of the proof proceeds as in \emph{loc.cit.}.
\end{proof}	
	
	\begin{proof}[Proof of Proposition \ref{prop:convolution-cartier}]
	First we note that $\mathscr{L}_n$ is canonically bi-rigidified, i.e. we have isomorphisms
	\[\mathscr{L}_n|_{\{e_Z\} \times B_n} \simeq \mc{O}_{B_n} \hspace{1em} \mathscr{L}_n|_{Z_n \times \{e_B\}} \simeq \mc{O}_{Z_n},\]
	where $e_Z: X^I \to Z_n, e_B: X^I \to B_n$ are the unit sections. Indeed, $e_B$ is the trivial $\mc{L}^{(n)} T$-bundle on $X^I$, so the map $Z_n \times \{e_B\} \to Z_n \times_{X^I} B_n \to \mb{B} \mb{G}_m$ factors through $\tx{pt} \inj \mb{B} \mb{G}_m$. Likewise, identity for $Z_n$ is the constant group homomorphism $\mc{L}^{(n)} T \to \{1\} \to \underline{\mb{G}_m}$, so the induced bundle is again always trivial.
	
	Now consider the following two Cartesian diagrams:
	\[\xymatrix{
	Z_n \ar^{\pi}[d] \ar^{\widetilde{e_B}}[r] & Z_n \times_{X^I} B_n \ar^{p_2}[d] \\
	X^I \ar^{e_B}[r] & B_n
	} \hspace{1em} \xymatrix{
	B_n \ar^{\pi}[d] \ar^{\widetilde{e_Z}}[r] & Z_n \times_{X^I} B_n \ar^{p_1}[d] \\
	X^I \ar^{e_Z}[r] & Z_n
	}\]
	We obtain the following from the left diagram:
	\begin{align*}
	& \Hom_{\IndCoh(Z_n)}(\omega_{Z_n}, \omega_{Z_n})
	\simeq \Hom_{\IndCoh(Z_n)}(\omega_{Z_n}, \pi^!(\omega_{X^I})) \\
	\simeq & \Hom_{\IndCoh(X^I)}(\pi_* \omega_{Z_n}, \omega_{X^I}) \simeq \Hom_{\IndCoh(X^I)}(\pi_* \widetilde{e_B}^!(\omega_{Z_n \times_{X^I} B_n}), \omega_{X^I}) \\
	\simeq & \Hom_{\IndCoh(X^I)}(\pi_* \widetilde{e_B}^!(\Upsilon(\mc{O}_{Z_n \times_{X^I} B_n})), \omega_{X^I}) \simeq \Hom_{\IndCoh(X^I)}(\pi_* \Upsilon(\widetilde{e_B}^*(\mc{O}_{Z_n \times_{X^I} B_n})), \omega_{X^I}) \\
	\simeq & \Hom_{\IndCoh(X^I)}(\pi_* \Upsilon(\widetilde{e_B}^*(\mathscr{L}_n), \omega_{X^I}) \simeq \Hom_{\IndCoh(X^I)}(\pi_* \widetilde{e_B}^!(\Upsilon(\mathscr{L}_n)), \omega_{X^I}) \\
	\simeq & \Hom_{\IndCoh(X^I)}(e_B^! (p_2)_*(\Upsilon(\mathscr{L}_n)), \omega_{X^I}) \simeq \Hom_{\IndCoh(X^I)}(K_{X^I / B_n} \tensor e_B^{*,\IndCoh} (p_2)_*(\Upsilon(\mathscr{L}_n)), \omega_{X^I}) \\
	\simeq & \Hom_{\IndCoh(X^I)}(e_B^{*,\IndCoh} (p_2)_*(\Upsilon(\mathscr{L}_n)), K_{X^I / B_n}^{-1} \tensor \omega_{X^I}) \\
	\simeq & \Hom_{\IndCoh(B_n)}((p_2)_*(\Upsilon(\mathscr{L}_n)), (e_B)_*(\Upsilon(K_{X^I / B_n}^{-1})))
	\end{align*}
	where we used the fact that the relative dualizing complex $K_{X^I / B_n} := e_B^{!}(\mc{O}_{B_n})$ is a line bundle, which is because the morphism $X^I \to B_n$ is smooth. Thus from the identity arrow $\omega_{Z_n} \to \omega_{Z_n}$ we obtain a \emph{root morphism}
	\[\rho_Z: (p_2)_*(\Upsilon(\mathscr{L}_n)) \to (e_B)_*(\Upsilon(K_{X^I / B_n}^{-1}))\]
	For the other root morphism, note that similarly we have
	\begin{align*}
	& \Hom_{\QCoh(B_n)}(\mc{O}_{B_n}, \mc{O}_{B_n}) \\
	\simeq & \Hom_{\QCoh(X^I)}(\pi^* \Upsilon^{-1}(\omega_{X^I}), \widetilde{e_Z}^*(\mc{O}_{Z_n \times_{X^I} B_n})) \\
	\simeq & \Hom_{\QCoh(X^I)}(\Upsilon^{-1}(\omega_{X^I}), \pi_* \widetilde{e_Z}^*(\mc{O}_{Z_n \times_{X^I} B_n})) \\
	\simeq & \Hom_{\IndCoh(X^I)}(\omega_{X^I}, \Upsilon(\pi_* \widetilde{e_Z}^*(\mc{O}_{Z_n \times_{X^I} B_n}))) \\
	\simeq & \Hom_{\IndCoh(X^I)}(\omega_{X^I}, \Upsilon(\pi_* \widetilde{e_Z}^*(\mathscr{L}_n))) \\
	\simeq & \Hom_{\IndCoh(X^I)}(\omega_{X^I}, \Upsilon(e_Z^* (p_1)_*(\mathscr{L}_n))) \\
	\simeq & \Hom_{\IndCoh(X^I)}(\omega_{X^I}, e_Z^! \Upsilon((p_1)_*(\mathscr{L}_n))) \\
	\simeq & \Hom_{\IndCoh(Z_n)}((e_Z)_!(\omega_{X^I}), \Upsilon((p_1)_*(\mathscr{L}_n)))
	\end{align*}
	where in the middle we used QCoh base change for the second diagram above. From the identity of $\mc{O}_{B_n}$ we get our second root morphism
	\[\rho_B: (e_Z)_!(\omega_{X^I}) \to \Upsilon((p_1)_*(\mathscr{L}_n)).\]

	Observe that $\Upsilon$ on the right is symmetric monoidal, and $e_B$ is symmetric monoidal with respect to the standard monoidal structure on $\IndCoh(X^I)$ and the convolution monoidal structure on $\IndCoh(B_n)$. It follows that if $\rho_Z$ is an isomorphism, then $(p_2)_*(\Upsilon(\mathscr{L}_n))$ is invertible in the convolution monoidal structure; same goes for $\rho_B$.	
	
	It is clear that the root morphisms are compatible with base change on $X^I$. Using the usual open-closed triangle, it suffices to prove them being isomorphism on each strata of $X^I$ (where we take the diagonal stratification). Since $\mc{L}^{(n)} T$ is \'etale locally trivial on each strata, all sheaves involved are \'etale locally constant; hence it suffices to prove the isomorphism after base-changing to each geometric point.
	
	We again reduce to $T = \mb{G}_m$. Suppose we are at a point $x \in X^I$ with collision pattern $(i_1, \ldots, i_N)$, then
	\[\mc{L}^{(n)}|_{x} \simeq \mb{G}_m^{N} \times \mb{G}_a^{n|I| - N} \hspace{1em} Z_n|_x \simeq \mb{Z}^{N} \times (\mb{A}^{n|I| - N})^\wedge_{0};\]
	so we see it suffices to check the root morphisms are isomorphisms for each collision type. The rest of the proof is now explicit computation, which immediately reduces to the case $(\mb{G}_m, \mb{Z})$ and $(\mb{G}_a, (\mb{A}^1)_0^\wedge)$.
	
	The $\mb{G}_m$ case is simple: the Poincare bundle places a weight $j$ representation on $\{j\} \times \mb{B} \mb{G}_m \subseteq \mb{Z} \times \mb{B} \mb{G}_m$. From this one readily computes that both sides of $\rho_Z$ are regular representation
	\[\mc{O}_{\mb{G}_m} \in \QCoh(\mb{B} \mb{G}_m)\]
	and both sides of $\rho_B$ are
	\[i_*(k) \in \IndCoh(\mb{Z}) \hspace{1em} i: \{0\} \inj \mb{Z}.\]
	For $\mb{G}_a$, let $Z := (\mb{A}^1)_{0}^\wedge$ and $f: \mb{G}_a \times Z \to \mb{G}_m$ be the pairing map. On the level of algebras, $f$ is given by (truncations) of the map
	\[k[x, x^{-1}] \to k[y][[t]] \hspace{1em}x \mapsto e^{y t};\]
	Our first job is to identify the Poincare bundle. By descent we have
	\[\QCoh(B \mb{G}_a \times Z) \simeq \tx{Tot}(\QCoh(\mb{G}_a^\bullet \times Z)),\]
	unwinding definition, we see that $\mathscr{L}$ can understood as $\mc{O}_Z \in \QCoh(Z)$ equipped with the equivariance given by multiplying by $e^{t y}$.
	Thus, as plain vector space, $(p_2)_*(\Upsilon(\mathscr{L}))$ is isomorphic to
	\[\pi_*(\omega_Z) \simeq \colimit (k \xrightarrow{t} k[t]/t^2 \xrightarrow{t} k[t]/t^3 \ldots) \simeq k[t^{-1}]\]
	where the map $k[t]/t^r \to k[t^{-1}]$ is by multiplying by $t^{-r + 1}$; this is additionally a colimit of cocommutative coalgebras. Each term admits an $k[y]$-coaction by truncations of $e^{t y}$, which are compatible with these transition maps, so $k[t^{-1}]$ itself admits a $k[y]$-coaction given by $t^{-r} \mapsto t^{-r}\left(1 + \ldots + (yt)^r / (r)!\right)$, which is isomorphic to the regular representation of $k[y]$ via $t^{-r} \mapsto \sum_{i = 1}^{r} \frac{y^i}{i!}$.

	For the other root morphism $\rho_B$, note that $!$-pullback along $e_Z$ is conservative, so it suffices to compute
\[e_Z^! \circ (e_Z)_!(k) \simeq k \oplus k[1] \simeq C^*(\mb{B}\mb{G}_a) \simeq \Upsilon(\pi_* \circ \widetilde{e_Z}^*(\mathscr{L})) \simeq e_Z^! \circ \Upsilon \circ (p_1)_*(\mathscr{L}).\]
This finishes the computation.
	\end{proof}

Now we move towards factorization.
\newcommand{\ICR}{\mc{I}nd\mc{C}oh_{*, \tx{ren}}}

Recall from Definition \ref{def:kl-factorizable} that we have a unital factorization category $\mathcal{R}ep( \mathcal{L}^+ T )$ over $\Ran_{\un, \dR}$, defined using renormalized ind-coherent sheaves. On the other hand, by Construction \ref{constr-unit-fact-prestack-indcoh-fact-cat}, the unital structure of $[\tx{Gr}_{\check{T}}]_{\un, \dR}$ equips the category
\[\ICR([\tx{Gr}_{\check{T}}]_{\Ran_{\un, \dR}})\]
with the structure of a unital factorization category over $\Ran_{\un, \dR}$ as well. (The condition from Remark \ref{rem-constr-unit-fact-prestack-indcoh-fact-cat} is easily verifiable.)

\begin{proposition}
\label{prop:fourier-mukai-indcoh}
	There exists an \emph{unital} equivalence of unital factorization categories over $\Ran_\dR$
	\[\tx{FM}: \mathcal{R}ep( \mathcal{L}^+ T ) \simeq \mc{I}nd\mc{C}oh_{*, \tx{ren}}([\tx{Gr}_{\check{T}}]_{\Ran_{\un, \dR}})\]
	whichi induces a pointwise equivalence
	\[\tx{FM}: \Rep(\mc{L}^+ T) \simeq \IndCoh_{*, \tx{ren}}([\tx{Gr}_{\check{T}}])\]
	that satisfies the following additional constraint: for each dominant weight $\clambda$, let $k^\clambda \in \Rep(\mc{L}^+ T)$ be the restriction of $k^\clambda \in \Rep(T)$ along $\Arc T \xrightarrow{t = 0} T$, then $\tx{FM}(k^\lambda)$ is the $\delta$-sheaf supported at $\clambda$.
\end{proposition}
\begin{proof}
For simplicity, introduce
\[[\tx{Pair}]_{\Ran_{\un, \dR}} := [\mb{B}\mc{L}^+ T]_{\Ran_{\un, \dR}} \times_{\Ran_{\un, \dR}} \tx{Gr}_{\check{T}, \Ran_{\un, \dR}}\]
and let $[\tx{Pair}]_{X^I}$ denote its base change along any $X^I \to \Ran_{\un, \dR}$. Consider the diagram
\[\xymatrix{
	& [\tx{Pair}]_{\Ran_{\un, \dR}} \ar_{p}[ld] \ar^{q}[rd] \\
[\mb{B}\mc{L}^+ T]_{\Ran_{\un, \dR}} & & \tx{Gr}_{\check{T}, \Ran_{\un, \dR}}
}\]
The top space again satisfies the condition from Remark \ref{rem-constr-unit-fact-prestack-indcoh-fact-cat}, so that  $\ICR([\tx{Pair}]_{\Ran_{\un, \dR}})$ is a unital factorization category.
Note that the map $p$ is $(?, \tx{ren})$-pullable and strictly unital, $q$ is $*$-pullable and strictly co-unital.

From the pairing map of corr-unital factorization spaces
\[[\tx{Pair}]_{\Ran_{\un, \dR}} \to \mb{B} \mb{G}_m \times \Ran_{\un, \dR}\]
we obtain an universal Poincare bundle
\[\Upsilon(\mathscr{L}_\tx{univ}) \in \ICR([\tx{Pair}]_{\Ran_{\un, \dR}})\]
using which we define the Fourier-Mukai morphism and its left adjoint
\[\tx{FM} := q_{*, \tx{ren}} \circ (- \tensor \mathscr{L}_\tx{univ}) \circ p^{?}_{\tx{ren}}: \mathcal{R}ep( \mathcal{L}^+ T ) \to \mc{I}nd\mc{C}oh_{*, \tx{ren}}([\tx{Gr}_{\check{T}}]_{\Ran_{\un, \dR}})\]
here by $\tensor^!$ we mean the $!$-pullback in the $!$-theory; we note that in the present setting, the $!$-theory and $*$-theory are canonically equivalent.
The general setup in Construction \ref{constr-functorial-unital-fact-indcoh} tells us that $\tx{FM}$ is a \emph{lax} unital map; however, explicit computation performed in the previous Proposition demonstrates that the unital map is actually always an isomorphism.
Thus we obtain a \emph{strictly unital} functor.

To see that it is an equivalence, we note that the restriction functors
\[\tx{FactCat}_\un \to \tx{ShvCat}(\Ran_{\un, \dR}) \to \tx{ShvCat}(\Ran_\dR) \to \tx{ShvCat}(\Ran)\]
are all conservative; indeed, the first two are general categorical facts, and the third one is because we can always use smoothness of $X^I$ to \emph{non-canonically} lift any $S$-point of $X^I_\dR$ to one of $X^I$.
It thus suffices to establish the map is an equivalence on each $X^I$. Using result from Section \ref{sect:indcoh-crys} we can write
\[\ICR([\tx{Gr}_{\check{T}}]_{X^I}) \simeq \colimit_{i, !\tx{-push}} \IndCoh(Z_n) \hspace{1em} \ICR([\mb{B} T]_{X^I}) \simeq \colimit_{i, *\tx{-pull}} \IndCoh(B_n)\]
with the stratification given above, and $\tx{FM}_{X^I}$ restricts to the system $\tx{FM}_n$ above; the compatibility of the system is given by base change. In particular, the additional constraint is met since the representations $k^\clambda$ are pulled back from $B_1$.
\end{proof}

\subsection{Monad Comparison}
Let $[\kappa, E]$ be a non-degenerate quantum parameter for the maximal torus $T$, and $[\check{\kappa}, \check{E}]$ be the Langlands dual parameter. By \cite{zhao2017quantum}, to such data we can attach a unital factorization line bundle $\mathscr{L}_{\widehat{\tx{Hecke}}_T}$ over $[\widehat{\tx{Hecke}}_T]_{\Ran_{\un, \dR}}$ as in Section \ref{sect:km-fact-twisted}. On the other hand, let
\[[\tx{Gr}_{\check{T}, \tx{rel}}^\bullet]_{\Ran_{\un, \dR}}\]
denote the Cech nerve of $[\tx{Gr}_{\check{T}}]_{\Ran_{\un, \dR}} \to [(\tx{Gr}_{\check{T}})_\dR]_{\Ran_{\un, \dR}}$; this is a cosimplicial object in the category of unital factorization spaces. The data $[\check{\kappa}, \check{E}]$ also gives rise to a multiplicative factorization line bundle $\mathscr{L}_{\tx{Gr}_{\check{T}, \tx{rel}}}$ over $[\tx{Gr}_{\check{T}, \tx{rel}}]_{\Ran_{\un, \dR}}$.
\begin{remark}
The compatibility between $\mathscr{L}_{\tx{Gr}_{\check{T}, \tx{rel}}}$ and $\mathscr{L}_{\widehat{\tx{Hecke}}_{T}}$ is as follows. Because $\kappa$ is non-degenerate, it induces an isomorphism of the Lie algebras $\mf{t} \simeq \check{\mf{t}}$ and thus an isomorphism of corr-unital factorization prestacks
\[\Psi_{[\kappa, E]}: [\widehat{\tx{Gr}}_T]_{\Ran_{\un, \dR}} \simeq [\widehat{\tx{Gr}}_{\check{T}}]_{\Ran_{\un, \dR}};\]
Let $p_0, p_1: [\widehat{\tx{Gr}}_T]_{\Ran_{\un, \dR}} \to [\widehat{\tx{Hecke}}_T]_{\Ran_{\un, \dR}}$ denote the two projection coming from base changes of the two maps $[\widehat{\tx{Hecke}}_T]_{\Ran_{\un, \dR}} \to [\mb{B}\mc{L}^+T]_{\Ran_{\un, \dR}}$ along $\tx{pt} \to [\mb{B}\mc{L}^+T]_{\Ran_{\un, \dR}}$ respectively. On the other hand, by forgetting one of the two $\check{T}$-bundles we have two evaluation maps $e_0, e_1: [\tx{Gr}_{\check{T}, \tx{rel}}]_{\Ran_{\un, \dR}} \to [\tx{Gr}_{\check{T}}]_{\Ran_{\un, \dR}}$, both of which are (split) surjections of relative group indschemes, whose kernel is isomorphic to $[\widehat{\tx{Gr}}_{\check{T}}]_{\Ran_{\un, \dR}}$. In other words we have two inclusions $q_0, q_1: [\widehat{\tx{Gr}}_{\check{T}}]_{\Ran_{\un, \dR}} \to [\tx{Gr}_{\check{T}, \tx{rel}}]_{\Ran_{\un, \dR}}$. The compatibility is then given by
\[\Psi_{[\kappa, E]}^! \circ q_i^!(\mathscr{L}_{\tx{Gr}_{\check{T}, \tx{rel}}}) \simeq p_i^!(\mathscr{L}_{\widehat{\tx{Hecke}}_T}) \hspace{1em} i = 0, 1;\]
We note that this pullback actually determines the entire line bundle, in the following sense. In the case of Hecke stack, the map $p_0, p_1$ admit splittings $r_0, r_1$, given by taking\footnote{Here the commutativity of $T$ is used.} the difference of the two bundles (resp. the inverse of the difference); the kernel of each $r_i$ is isomorphic to $[\mb{B} \mc{L}^+ T]_{\Ran_{\un, \dR}}$, giving rise to two isomorphisms
\[s_0, s_1: [\widehat{\tx{Hecke}}_T]_{\Ran_{\un, \dR}} \simeq [\mb{B}\mc{L}^+T]_{\Ran_{\un, \dR}} \times_{\Ran_{\un, \dR}} [\widehat{\tx{Gr}}_T]_{\Ran_{\un, \dR}};\]
if we let $\tx{pr}_{\widehat{\tx{Gr}}_T}$ be the projection to the second factor, then unfolding definition one sees that
\[s_i^! \circ \tx{pr}_{\widehat{\tx{Gr}}_T}^! \circ p_i^!(\mathscr{L}_{\widehat{\tx{Hecke}}_T}) \simeq \mathscr{L}_{\widehat{\tx{Hecke}}_T} \hspace{1em} i = 0, 1;\]
the situation for $\tx{Gr}_{\check{T}}$ is completely analogous.
\end{remark}

Recall from Construction \ref{constr:hecke-line-bundle-is-algebra} that $\mathscr{L}_{\widehat{\tx{Hecke}}_T}$ has additional structure of an associative algebra object (internal to the category of unital factorization algebras) under the convolution monoidal structure\footnote{We note that for Torus the Tate bundle is canonically trivial, so we are allowed to freely switch between $!$ and $*$ theories; for notational unity we stick with the $*$-theory everywhere.} of $\ICR([\widehat{\tx{Hecke}}_T]_{\Ran_{\un, \dR}}) \in \tx{Alg}(\tx{FactCat}_{\un})$. In an completely analogous manner, convolution upgrades $\ICR([\tx{Gr}_{\check{T}, \tx{rel}}]_{\Ran_{\un, \dR}})$ to an element of $\tx{Alg}(\tx{FactCat}_{\un})$, and $\mathscr{L}_{\tx{Gr}_{\check{T}, \tx{rel}}}$ obtains the structure of an associative algebra object within the category of unital factorization algebras in $\ICR([\tx{Gr}_{\check{T}, \tx{rel}}]_{\Ran_{\un, \dR}})$.

Through convolution $\ICR([\tx{Gr}_{\check{T}, \tx{rel}}]_{\Ran_{\un, \dR}})$ acts on $\ICR([\tx{Gr}_{\check{T}}]_{\Ran_{\un, \dR}})$ (within the symmetric monoidal category $\tx{FactCat}_{\un}$), and we shall define twisted D-modules on affine Grassmannian via this action:
\begin{definition}
We define
\[\DMod_{[\check{\kappa}, \check{E}]}([\tx{Gr}_{\check{T}}]_{\Ran_{\un, \dR}}) := \mathscr{L}_{\tx{Gr}_{\check{T}, \tx{rel}}}\mod(\ICR([\tx{Gr}_{\check{T}}]_{\Ran_{\un, \dR}})) \in \tx{FactCat}_\un.\]
\end{definition}

\begin{remark}
Thanks to results from Section \ref{sect:indcoh-crys}, one can readily verify that over each $X^I_\dR$, this recovers the usual category of $\mathscr{L}_{\tx{Gr}_{\check{T}, \tx{rel}}}$-twisted D-modules on the Beilinson-Drinfeld grassmannian. Over each $X^I_\dR$, The effect of the monad induced by $\mathscr{L}_{\tx{Gr}_{\check{T}, \tx{rel}}}$ is the composition of the induction functor
\[\tx{ind}_{\DMod_{[\check{\kappa}, \check{E}]}}: \IndCoh_{*, \tx{ren}}([\tx{Gr}_{\check{T}}]_{X^I_\dR}) \to \DMod_{[\check{\kappa}, \check{E}]}([\tx{Gr}_{\check{T}}]_{X^I_\dR})\]
and the forgetful functor (which is the right adjoint of). In particular, the unit factorization algebra is given by $\delta$-D-module supported at the unit section.
\end{remark}

On the other hand, from Definition \ref{def:kl-factorizable} we have the factorizable Kazhdan--Lusztig category, defined as
\[\mc{KL}(T)_{[\kappa, E]} := \mc{KL}(T)_{\mathscr{L}_{\widehat{\tx{Hecke}}_T}} := \mathscr{L}_{\widehat{\tx{Hecke}}_T}\mod(\mc{R}ep(\mc{L}^+ T)).\]

\begin{proposition}[Toric Fundamental Local Equivalence]
\label{prop:toric-FLE}
	We have an equivalence of unital factorization categories
	\[\tx{FLE}_T: \mc{KL}(T)_{[\kappa, E]} \simeq \DMod_{[\check{\kappa}, \check{E}]}([\tx{Gr}_{\check{T}}]_{\Ran_{\un, \dR}})\]
	compatible with the (strictly unital factorizable) induction functors and the (lax-unital factorizable) forgetful functors to the two sides of Proposition \ref{prop:fourier-mukai-indcoh}.
\end{proposition}

\begin{proof}
We note that there is a third cosimplicial object in the category of corr-unital factorization spaces
\[[\tx{Prod}^\bullet]_{\Ran_{\un, \dR}}\]
mapping to both $[\tx{Gr}_{\check{T}, \tx{rel}}^\bullet]_{\Ran_{\un, \dR}}$ and $[\widehat{\tx{Hecke}}_T^\bullet]_{\Ran_{\un, \dR}}$. On the zeroth level it is
\[[\tx{Prod}^0]_{\Ran_{\un, \dR}} := [\mb{B} \mc{L}^+ T]_{\Ran_{\un, \dR}} \times_{\Ran_{\un, \dR}} [\tx{Gr}_{\check{T}}]_{\Ran_{\un, \dR}};\]
here the fiber product is taken within the category of corr-unital factorization spaces. It comes with two maps
\[\pi^0_L: [\tx{Prod}^0]_{\Ran_{\un, \dR}} \to [\mb{B} \mc{L}^+ T]_{\Ran_{\un, \dR}} \hspace{1em} \pi^0_R: [\tx{Prod}^0]_{\Ran_{\un, \dR}} \to [\tx{Gr}_{\check{T}}]_{\Ran_{\un, \dR}};\]
On the first level it is given by
\[[\tx{Prod}^1]_{\Ran_{\un, \dR}} := [\tx{Prod}^0]_{\Ran_{\un, \dR}} \times^{\pi_L^0, p^{\widehat{\tx{Hecke}}_T}_0}_{[\mb{B} \mc{L}^+ T]_{\Ran_{\un, \dR}}} [\widehat{\tx{Hecke}}_T]_{\Ran_{\un, \dR}}\]
where $p^{\widehat{\tx{Hecke}}_T}_0$ is the zeroth structural map $[\widehat{\tx{Hecke}}_T]_{\Ran_{\un, \dR}} \to [\mb{B} \mc{L}^+ T]_{\Ran_{\un, \dR}}$. We let $p^{\tx{Prod}^1}_0$ and $\pi^1_L$ denote the projection to $[\tx{Prod}^0]_{\Ran_{\un, \dR}}$ and $[\widehat{\tx{Hecke}}_T]_{\Ran_{\un, \dR}}$ respectively.
As in the previous remark, the short exact sequence
\[\Ran_{\un, \dR} \to [\tx{Gr}_{\check{T}}]_{\Ran_{\un, \dR}} \to [\tx{Gr}_{\check{T}, \tx{rel}}]_{\Ran_{\un, \dR}} \to [\widehat{\tx{Gr}}_{\check{T}}]_{\Ran_{\un, \dR}} \to \Ran_{\un, \dR}\]
admits two splittings, giving rise to two isomorphisms
\[t_0, t_1: [\tx{Gr}_{\check{T}, \tx{rel}}]_{\Ran_{\un, \dR}} \simeq [\tx{Gr}_{\check{T}}]_{\Ran_{\un, \dR}} \times_{\Ran_{\un, \dR}} [\widehat{\tx{Gr}}_{\check{T}}]_{\Ran_{\un, \dR}}\]
We thus have a map
\[[\tx{Prod}^1]_{\Ran_{\un, \dR}} \simeq [\tx{Prod}^0]_{\Ran_{\un, \dR}} \times^{\pi_L^0, p^{\widehat{\tx{Hecke}}_T}_0}_{[\mb{B} \mc{L}^+ T]_{\Ran_{\un, \dR}}} [\widehat{\tx{Hecke}}_T]_{\Ran_{\un, \dR}}\]
\[\xrightarrow[s_0]{\simeq} [\mb{B} \mc{L}^+ T]_{\Ran_{\un, \dR}} \times_{\Ran_{\un, \dR}} [\tx{Gr}_{\check{T}}]_{\Ran_{\un, \dR}} \times_{\Ran_{\un, \dR}} [\widehat{\tx{Gr}}_{T}]_{\Ran_{\un, \dR}}\]
\[\xrightarrow[\Psi_{[\kappa, E]}]{\simeq} [\mb{B} \mc{L}^+ T]_{\Ran_{\un, \dR}} \times_{\Ran_{\un, \dR}} [\tx{Gr}_{\check{T}}]_{\Ran_{\un, \dR}} \times_{\Ran_{\un, \dR}} [\widehat{\tx{Gr}}_{\check{T}}]_{\Ran_{\un, \dR}}\]
\[\xrightarrow[t_0^{-1}]{\simeq} [\tx{Prod}^0]_{\Ran_{\un, \dR}} \times^{\pi_R^0, p^{\tx{Gr}_{\check{T}, \tx{rel}}}_0}_{[\tx{Gr}_{\check{T}}]_{\Ran_{\un, \dR}}} [\tx{Gr}_{\check{T}, \tx{rel}}]_{\Ran_{\un, \dR}}\]
and we let $\pi^1_R$ denote the projection to the second factor. The following diagram of corr-unital factorization spaces, in which all squares are Cartesian, summarizes the construction above:
\[		
	\xymatrix{
	[\mb{B} \mc{L}^+ T]_{\Ran_{\un, \dR}} & 
[\widehat{\tx{Hecke}}_T]_{\Ran_{\un, \dR}}	
	\ar_{p^{\widehat{\tx{Hecke}}_T}_0}[l] \ar^{p^{\widehat{\tx{Hecke}}_T}_1}[r] & [\mb{B} \mc{L}^+ T]_{\Ran_{\un, \dR}} \\
	[\tx{Prod}^0]_{\Ran_{\un, \dR}} \ar^{\pi_L^0}[u] \ar_{\pi_R^0}[d] & [\tx{Prod}^1]_{\Ran_{\un, \dR}} \ar^{p^{\tx{Prod}^1}_0}[l] \ar_{p^{\tx{Prod}^1}_1}[r] \ar^{\pi_L^1}[u] \ar_{\pi_R^1}[d] & [\tx{Prod}^0]_{\Ran_{\un, \dR}} \ar_{\pi_L^0}[u] \ar^{\pi_R^0}[d] \\
	[\tx{Gr}_{\check{T}}]_{\Ran_{\un, \dR}} & [\tx{Gr}_{\check{T}, \tx{rel}}]_{\Ran_{\un, \dR}} \ar^{p^{\tx{Gr}_{\check{T}, \tx{rel}}}_0}[l] \ar_{p^{\tx{Gr}_{\check{T}, \tx{rel}}}_1}[r] & [\tx{Gr}_{\check{T}}]_{\Ran_{\un, \dR}}
	}\]
By similarly carrying it out over each level, we obtain a cosimplicial object $[\tx{Prod}^\bullet]_{\Ran_{\un, \dR}}$ equipped with two maps $\pi_L^\bullet, \pi_R^\bullet$ of cosimplicial objects in corr-unital factorization spaces. Furthermore, chasing through definitions one sees that the two maps
\[[\tx{Prod}^1]_{\Ran_{\un, \dR}} \xrightarrow{p^{\tx{Prod}^1}_0} [\tx{Prod}^0]_{\Ran_{\un, \dR}} \to \mb{B} \mb{G}_m \hspace{1em} [\tx{Prod}^1]_{\Ran_{\un, \dR}} \xrightarrow{p^{\tx{Prod}^1}_1} [\tx{Prod}^0]_{\Ran_{\un, \dR}} \to \mb{B} \mb{G}_m\]
are canonically isomorphic; same procedure for higher levels equip $\mathscr{L}_{\tx{univ}}$ with an equivariant structure with respect to the groupoid $[\tx{Prod}^\bullet]_{\Ran_{\un, \dR}}$.
Noting that the bottom-left and the top-right square satisfy base change for $\ICR$-theory, the necessary homotopy data for upgrading to an equivalence of monadic modules now follow from base change.
\end{proof}

\begin{remark}
The additional constringt from Proposition \ref{prop:fourier-mukai-indcoh}, together with the compatibility with induction, means we have the following pointwise compatibility (see Section \ref{sect:prep_km} for notations): the \emph{induced representation} $\mb{V}_\kappa^{\clambda} \in \tx{KL}_\kappa(T)$ is sent to the $\delta$ D-module supported at $\clambda$ on $\tx{Gr}_{\check{T}}$.
\end{remark}

\subsection{Global Compatibility}
\label{sect:tori-glob-compatibility}
Let $X$ be any smooth complete curve. Recall that we have the global Langlands equivalence
\[\tx{FM}_{\tx{glob}, \kappa}: \DMod_{[\kappa, 0]}(\tx{Bun}_T) \simeq \DMod_{[\check{\kappa}, 0]}(\tx{Bun}_{\check{T}}),\]
which can be constructed as follows. We have a universal Weil pairing
\[\tx{Bun}_T \times \tx{Bun}_{\check{T}} \to \mb{B}\mb{G}_m\]
induced by the pairings $\mb{B} T \times \weightLat \to \mb{B} \mb{G}_m, \mb{B} \check{T} \times \coweightLat \to \mb{B} \mb{G}_m$ and the self-duality pairing $\tx{Pic} \times \tx{Pic} \to \mb{B} \mb{G}_m$. Let $\mc{P}$ denote the resulting universal Poincare bundle on $\tx{Bun}_T \times \tx{Bun}_{\check{T}}$. In the case of $T \simeq \check{T} \simeq \mb{G}_m$, its fiber at a pair of line bundles $\mc{L}, \mc{L}'$ is explicitly given by
\[\det \Gamma(X, \mc{O}) \tensor \det \Gamma(X, \mc{L} \tensor \mc{L}') \tensor (\det \Gamma(X, \mc{L}))^{-1} \tensor (\det \Gamma(X, \mc{L}'))^{-1}.\]
Set $\pi_1: \tx{Bun}_T \times \tx{Bun}_{\check{T}} \to \tx{Bun}_T$, $\pi_2: \tx{Bun}_T \times \tx{Bun}_{\check{T}} \to \tx{Bun}_{\check{T}}$, so that
\[\tx{FM}_\IndCoh: (\pi_2)_*(\pi_1^!(-) \tensor \mc{P}): \IndCoh(\tx{Bun}_T) \to \IndCoh(\tx{Bun}_{\check{T}})\]
is an equivalence. 
\[\tx{FM}_\IndCoh^{-1} \simeq (\pi_1)_*(\pi_2^!(-) \tensor \mc{P}^\tx{inv}).\]

Let $\widehat{\tx{Bun}}_{T}$ denote the formal completion of the diagonal within $\tx{Bun}_{T} \times \tx{Bun}_{T}$, i.e. we have the Cartesian diagram
\[\xymatrix{
\widehat{\tx{Bun}}_{T} \ar^{\widehat{\Delta}_{T}}[r] \ar[d] & \tx{Bun}_{T} \times \tx{Bun}_{T} \ar[d] \\
(\tx{Bun}_{T})_\dR \ar^(0.4){\Delta_{T}}[r] & (\tx{Bun}_{T})_\dR \times (\tx{Bun}_{T})_\dR
}\]
Let $\mc{L}_{\kappa}$ be the line bundle on $\widehat{\tx{Bun}}_{T}$ satisfying
\[\mc{L}_{\kappa}\mod(\IndCoh(\tx{Bun}_{T})) \simeq \DMod_{[\check{\kappa}, 0]}(\tx{Bun}_{T});\]
As in the local setting, taking differences of the two $T$-bundles gives two splittings
\[s_0, s_1: \widehat{\tx{Bun}}_T \simeq \tx{Bun}_T \times (\tx{Bun}_T)^\wedge_{\tx{triv}}\]
where the last factor is the formal completion around the trivial $T$-bundle, and $\mc{L}_\kappa$ is pulled back (via either map) from a line bundle $\mc{L}_\kappa^\wedge$ on $(\tx{Bun}_T)^\wedge_{\tx{triv}}$. The situation is exactly the same for $\mc{L}_{\check{\kappa}}$ over $\widehat{\tx{Bun}_{\check{T}}}$, and under the isomorphism $(\tx{Bun}_T)^\wedge_{\tx{triv}} \simeq (\tx{Bun}_{\check{T}})^\wedge_{\tx{triv}}$ induced by $\kappa$, $\mc{L}_\kappa^\wedge$ and $\mc{L}_{\check{\kappa}}^\wedge$ correspond to each other. Using the same simplicial construction as in the proof of Proposition \ref{prop:toric-FLE}, we see that $\tx{FM}_{\IndCoh}$ upgrades to an equivalence $\tx{FM}_{\tx{glob}, \kappa}$ between twisted D-modules.
\begin{remark}
Let us explain how to relate the above construction to the classical one. For any three laft prestacks $X, Y, Z$, we have the natural convolution functor
\[\star: \IndCoh(X \times Y) \tensor \IndCoh(Y \times Z) \to \IndCoh(X \times Z);\]
In the case of $X = \tx{Bun}_{\check{T}}, Y = \tx{Bun}_T, Z = \tx{pt}$, the operation $\mc{P} \star (-)$ becomes the Fourier-Mukai transform for IndCoh. In the case of $X = Y = \tx{Bun}_{T}, Z = \tx{pt}$, we obtain the usual convolution action of $\IndCoh(\tx{Bun}_T \times \tx{Bun}_T)$ on $\IndCoh(\tx{Bun}_T)$, and a $[\kappa, 0]$-twisted D-module is exactly a module over $(\widehat{\Delta}_{T})_!(\mc{L}_\kappa)$ under this action. Now, the isomorphism between $\mc{L}_\kappa$ and $\mc{L}_{\check{\kappa}}$ over $(\tx{Bun}_T)^\wedge_{\tx{triv}} \simeq (\tx{Bun}_{\check{T}})^\wedge_{\tx{triv}}$ gives an identification
\[\mc{P} \star (\widehat{\Delta}_{T})_!(\mc{L}_\kappa) \simeq (\widehat{\Delta}_{\check{T}})_!(\mc{L}_{\check{\kappa}}) \star \mc{P},\]
from which it follows that $\mc{P} \star (-)$ upgrades to a functor
\[ (\widehat{\Delta}_{T})_!(\mc{L}_\kappa)\mod(\IndCoh(\tx{Bun}_T), \star) \xrightarrow[\simeq]{\mc{P} \star (-)} (\widehat{\Delta}_{\check{T}})_!(\mc{L}_{\check{\kappa}})\mod(\IndCoh(\tx{Bun}_{\check{T}}), \star)\]
This is the functor constructd in \cite[Theorem 6.5]{polishchuk2001fourier}; it tautologically agrees with our construction above.
\end{remark}

The compatibility with our construction above is summarized in the following diagram (c.f. \cite[Section 6.4.2]{gaitsgory2016eisenstein}):
\begin{proposition}
For $[\kappa, 0]$ positive, the following diagram
\[\xymatrix{
\mc{KL}(T)_{[\kappa, 0]} \ar^(0.4){\tx{FLE}_T}[r] \ar^{\tx{Loc}}[d] & \DMod_{[\check{\kappa}, 0]}([\tx{Gr}_{\check{T}}]_{\Ran_{\un, \dR}}) \ar^{\tx{Poinc}_*}[d] \\
\DMod_{[\kappa, 0]}(\tx{Bun}_T) \ar^{\tx{FM}_{\tx{glob}, \kappa}}[r] & \DMod_{[\check{\kappa}, 0]}(\tx{Bun}_{\check{T}})
}\]
commutes. Here $\tx{Loc}$ is the chiral localization functor as defined in Theorem \ref{thm-localization-KL}, and $\tx{Poinc}_*$ is the $*$-pushforward along $\pi_{\tx{Gr}}: [\tx{Gr}_{\check{T}}]_{\Ran_{\un, \dR}} \to \tx{Bun}_{\check{T}}$.
\end{proposition}

\begin{proof}
It suffices to construct the commutativity of the diagram obtained by passing the two vertical arrows to right adjoints, i.e.
\[\xymatrix{
\mc{KL}(T)_{[\kappa, 0]} \ar^(0.4){\tx{FLE}_T}[r] & \DMod_{[\check{\kappa}, 0]}([\tx{Gr}_{\check{T}}]_{\Ran_{\un, \dR}}) \\
\DMod_{[\kappa, 0]}(\tx{Bun}_T) \ar^{\tx{Loc}^R}[u] \ar^{\tx{FM}_{\tx{glob}, \kappa}}[r] & \DMod_{[\check{\kappa}, 0]}(\tx{Bun}_{\check{T}}) \ar^{\pi_{\tx{Gr}}^!}[u]
}\]
It further suffices to restrict to each (quasi-compact) connected component $U$ of $\tx{Bun}_T$, i.e. \[\xymatrix{
\mc{KL}(T)_{[\kappa, 0]} \ar^(0.4){\tx{FLE}_T}[r] & \DMod_{[\check{\kappa}, 0]}([\tx{Gr}_{\check{T}}]_{\Ran_{\un, \dR}}) \\
\DMod_{[\kappa, 0]}(U) \ar^{\tx{Loc}^R}[u] \ar^{\tx{FM}_{\tx{glob}, \kappa}}[r] & \DMod_{[\check{\kappa}, 0]}(\tx{Bun}_{\check{T}}) \ar^{\pi_{\tx{Gr}}^!}[u]
}\]
where the left vertical arrow is the functor constructed in Diagram \eqref{eq:coloc-construction} during Construction \ref{const-local-to-global-hecke}, which is induced by the $*$-pushforward along $\widehat{\tx{Hecke}}_{T, U}^\bullet \to \widehat{\tx{Hecke}}_T^\bullet$. At level $[0]$, the functor
\[\tx{Loc}^R_{[0]}: \IndCoh(U) \to \IndCoh([\tx{Gr}_{\check{T}}]_{\Ran_\dR})\] is given by $*$-pushing along $U \to [\mb{B} \Arc T]_{\Ran_\dR}$; routine base change shows that the top circuit of the diagram above becomes $q_* \circ \varphi_* (\mc{P}_\tx{glob} \tensor p^!(-))$, where
\[\xymatrix{
U \times [\tx{Gr}_{\check{T}}]_{\Ran_\dR} \ar^{p}[d] \ar^(0.35){\varphi}[r] & [\mb{B} \Arc T]_{\Ran_\dR} \times_{\Ran_\dR} [\tx{Gr}_{\check{T}}]_{\Ran_\dR} \ar^(0.68){q}[r] & [\tx{Gr}_{\check{T}}]_{\Ran_\dR} \\
U
}\]
and $\mc{P}_\tx{glob}$ is the pullback of the global Poincare bundle. On the other hand, the bottom circuit becomes $q_* (\mc{P}_\tx{loc} \tensor \varphi_* \circ p^!(-))$, where $\mc{P}_\tx{loc}$ is the pullback of the local Poincare bundle. The equivalence at levele $[0]$ thus follows from the fact that $\varphi^*(\mc{P}_\tx{loc}) \simeq \mc{P}_{\tx{glob}}$. Because the local twisting is pulled back from the global twisting, similar analysis yields a commutative diagram for each level, compatible in an obvious manner.
\end{proof}

The Weil pairing
\[\langle -, - \rangle_{\tx{Weil}}: \tx{Bun}_T \times \tx{Bun}_{\check{T}} \to \mb{G}_m\]
defines, for every $\check{T}$ bundle $\check{\mc{Q}}$, a line bundle $\mc{L}_{\check{\mc{Q}}}$ on $\tx{Bun}_T$, whose fiber at $\mc{P}$ is given by $\langle \mc{P}, \check{\mc{Q}} \rangle_{\tx{Weil}}$; more precisely, this is the (quasi-coherent) Fourier-Mukai dual of the skyscraper sheaf at $\mc{Q}$.

Considering $\crho$ as a coweight for $\check{T}$, it induces from $\omega_X$ a $\check{T}$-bundle $\omega_X^{\crho}$ on $X$, so we have a line bundle $\mc{L}_{\omega_X^{\crho}}$. Tracing through definition, one sees that the twisting corresponding to $[0, \AnomalyTerm]$ is precisely an integral twisting induced by $\mc{L}_{\omega_X^{\crho}}$.
It follows that the same pull-push procedure gives rise to a twisted version
\[\tx{FM}_{\tx{glob}, \kappa}^{\omega^{\crho}}: \DMod_{[\kappa, \AnomalyTerm]}(\tx{Bun}_T) \simeq \DMod_{\DualTwisting}(\tx{Bun}_{\check{T}})\]
here $\DualTwisting$ is the twisting pulled back from the twisting $[\check{\kappa}, 0]$ along the map $\tx{Bun}_{\check{T}} \to \tx{Bun}_{\check{T}}$ given by translation by $(\omega_X^{\crho})^{-1}$.

\begin{definition}
\label{def:shifted-FLE}
We define the \emph{$\omega^{\crho}$-shifted torus FLE} to be variant of the above construction, where the kernel transformation in Proposition \ref{prop:fourier-mukai-indcoh} is replaced with the kernel transformation (using the same Poincare bundle as before) along the following diagram
\[\xymatrix{
	& [\tx{Pair}]_{\Ran_{\un, \dR}} \ar_{p}[ld] \ar^{q}[rd] \\
[\mb{B}\mc{L}^+ T]_{\Ran_{\un, \dR}} & & [\tx{Gr}_{\check{T}}]_{\Ran_{\un, \dR}} \ar^{\tensor \omega^{\crho}}_{\simeq}[rd] \\
& & & [\tx{Gr}_{\check{T}}^{\omega^{\crho}}]_{\Ran_{\un, \dR}}
}\]
Unfolding definition, one sees that the appropriate twisting signature is
\[\tx{FLE}^{\omega^{\crho}}_T: \mc{KL}(T)_{[\kappa, \AnomalyTerm]} \simeq \DMod_{\DualTwisting}([\tx{Gr}_{\check{T}}^{\omega^{\crho}}]_{\Ran_{\un, \dR}}).\]
\end{definition}

It follows by unwinding definition that the following diagram commutes, again for $[\kappa, 0]$ positive:
\[\xymatrix{
\mc{KL}(T)_{[\kappa, \AnomalyTerm]} \ar^(0.4){\tx{FLE}^{\omega^{\crho}}_T}[rr] \ar^{\tx{Loc}}[d] & & \DMod_{\DualTwisting}([\tx{Gr}_{\check{T}}^{\omega^{\crho}}]_{\Ran_{\un, \dR}}) \ar^{\tx{Poinc}_*}[d] \\
\DMod_{[\kappa, \AnomalyTerm]}(\tx{Bun}_T) \ar^{\tx{FM}_{\tx{glob}, \kappa}^{\omega^{\crho}}}[rr] & & \DMod_{\DualTwisting}(\tx{Bun}_{\check{T}})
}\]

\begin{remark}
This $\crho$-shift is also reflected in the cotangent fiber shifts on the Kac-Moody side, which we opted to ignore in this paper (c.f. Warning \ref{warn:cotangent-fiber}).
\end{remark}

\section{Categorical Verdier Duality}
\label{sect:factvse2}

\begin{notation}
	For any $\mc{C}$ symmetric monoidal and $M$ a smooth manifold, we will let $\Psi_\mc{C}$ denote the fully faithful embedding $\mb{E}_M\tx{-Alg}_\tx{nu}^{\tx{loc}}(\mc{C}) \to \tx{Disk}(M)\tx{-Alg}(\mc{C}) \to \tx{Fact}(M)\tx{-Alg}(\mc{C})$ which is the combination of \cite[5.4.5.15]{HA} and \cite[5.5.4.10]{HA}. We let $\tx{Fact}(M)\tx{-Alg}^\tx{fact}(\mc{C})$ denote the essential image (i.e. the $\mc{C}$-valued factorizable cosheaves) and let $\Psi_\mc{C}^{-1}$ denote the inverse functor defined on it. When $\mc{C}$ is clear from the context, we will omit the subscript $\mc{C}$.
\end{notation}

\begin{convention}
\label{conv:C-is-R2}
We identify the topological manifold underlying $\mb{A}^1(\mb{C})$ as an oriented real manifold with $\mb{R}^2$ in the standard way. This in particular gives two distinguished embeddings $\mb{R}^1 \subseteq \mb{R}^2$ (thus two ways of forgetting $\mb{E}_2$-structures to $\mb{E}_1$-structures; for our purposes it would not matter which one to choose) and a choice of generator for $\pi_1(\mb{C}^\times)$.
\end{convention}

\subsection{Constructible Factorization Categories}
We first make precise of what ``topological factorization category'' means. For technical reasons, it is easier to work with a slight variant which we call \emph{constructible} factorization categories.

Let $\tx{StratTop}$ denote the category of conically stratified $D_\omega$ topological spaces $S \to A$ (where $A$ is the stratification poset) as described in \cite{lejay2021constructible}. A morphism is a square
\[\xymatrix{S \ar[r] \ar[d] & S' \ar[d] \\ A \ar[r] & A'}\]
i.e. the induced stratification on $S$ needs to be coarser than $A$. $(S \to A) \times (S' \times A') := (S \times S \to A \times A')$ gives it a symmetric monoidal structure.

\begin{definition}
	For any $X \to A \in \tx{StratTop}$ and any symmetric monoidal category $\mc{C}$, we let $\tx{Shv}_A(X; \mc{C})$ (resp. $\tx{cShv}_A(X; \mc{C})$) denote the category of $\mc{C}$-valued \emph{hyper(co)complete} $A$-constructible (co)sheaves (a notion defined in \cite{lejay2021constructible}). The functoriality is given by (usual) pushforward and \emph{hyper}pullbacks\footnote{that it preserves constructibility is Corollary 2.18 of \emph{loc.cit.}.}. We will abuse language by referring to such objects simply as \emph{$A$-constructible (co)sheaves}.
\end{definition}

\begin{notation}
We will let $(f_\dagger, f^\dagger)$ denote the adjoint pair of cosheaf pushforward / (hyper)pullback.
\end{notation}

\begin{remark}
Since the problem of hypercompleteness only occurs when we consider Ran space itself, we have opted to abuse notation in favor of brevity. We hope this does not cause too much confusion.
\end{remark}

Recall that the category of all $\mc{C}$-valued hypercocomplete cosheaves $\tx{cShv}(X; \mc{C})$ has a symmetric monoidal structure for any symmetric monoidal $\mc{C}$ whose tensor commutes with colimit in each variable: by \cite[2.2.1.9]{HA}, there is an unique such structure such that the cosheafify-then-hypercocomplete functor
\[\tx{coShfycoHyp}: \tx{Fun}(\mc{U}(X), \mc{C}) \to \tx{cShv}(X; \mc{C})\]
upgrades to a symmetric monoidal functor, where we equip the former with pointwise tensor product.
\begin{remark}
\label{remark:cShvA-sym-monoidal}
If $\varphi: \mc{C} \to \mc{D}$ is a symmetric monoidal functor, then so is the induced map $\tilde{\varphi}: \tx{cShv}(X; \mc{C}) \to \tx{cShv}(X; \mc{D})$. The category of $A$-constructible cosheaves is stable under this tensor product and $\tilde{\varphi}$ restricts to a symmetric monoidal functor $\tx{cShv}_A(X; \mc{C}) \to \tx{cShv}_A(X; \mc{D})$.
\end{remark}

\paragraph{Enter-Path Model}
Let us introduce a different model for $A$-constructible $\mc{C}$-valued cosheaves. By \cite{lejay2021constructible} (extending \cite[A.9.3]{HA}), for $(X \to A) \in \tx{StratTop}$, we have
\[\tx{cShv}_{A}(X; \mc{C}) \simeq \tx{Shv}_A(X; \mc{C}^\op)^\op \simeq \tx{Fun}(\tx{Exit}_A(X), \mc{C}^\op)^\op \simeq \tx{Fun}(\tx{Entr}_A(X), \mc{C})\]
where $\tx{Entr}_A(X) := \tx{Exit}_A(X)^\op$. For every map $X \to Y$ of $A$-stratified space we have an obvious functor $\tx{Entr}_A(X) \to \tx{Entr}_A(Y)$, and precomposition gives the pullback functor; it then follows that $\dagger$-pushforward is given by left Kan extension along $\tx{Entr}_A(X) \to \tx{Entr}_A(Y)$.

A special case is when $A$ is trivial. In such case, $\tx{Entr}_A(X) \simeq X^\tx{Spc}$ is a Kan complex, so we have
\[\tx{cShv}_A(X; \mc{C}) \simeq \tx{Fun}(X^{\tx{Spc}, \op}, \mc{C}) \simeq \limit_{X^\tx{Spc}} \mc{C} =: \ms{LS}(X; \mc{C}).\]
Note that we also have $\ms{LS}(X; \mc{C}) \simeq \colimit_{X^{\tx{Spc}, \op}} \mc{C}$.

The assignments
\[\tx{ShvMod}: \tx{StratTop}^\op \to 1\tx{-Cat} \hspace{1em} (S \to A) \mapsto \tx{cShv}_A(S)\mod(\DGCat)\]
is therefore well-defined, sending $f: (S \to A) \to (S' \to A')$ to $M \mapsto M \tensor_{\tx{cShv}_{A'}(S')} \tx{cShv}_A(S)$. Similarly we have another map
\[\tx{StratTop}^\op \to 1\tx{-Cat} \hspace{1em} (S \to A) \mapsto \{M \in \tx{cShv}_A(S)\mod(\DGCat), m \in M\}\]
these two functors give rise to a map of two Cartesian fibrations
\[\tx{ShvModSect}: \tx{StratTop}^{\tx{ShvModSect}} \to \tx{StratTop}^{\tx{ShvMod}}\]
over $\tx{StratTop}$. Let $\tx{StratTop}^{\op, \tx{ShvMod}, \vee} \to \tx{StratTop}^\op$ and $\tx{StratTop}^{\op, \tx{ShvModSect}, \vee} \to \tx{StratTop}^\op$ denote the dual coCartesian fibrations.
Taking box product then inducing along
\[\tx{cShv}_A(S) \tensor \tx{cShv}_{A'}(S') \xrightarrow{\boxtimes} \tx{cShv}_{A \times A'}(S \times S')\]
gives both of them symmetric monoidal structures, such that the map (and their maps to $\tx{StratTop}$) are symmetric monoidal.

Note that the $\mb{C}$-points of $\Ran(X)$ gives an oplax symmetric monoidal functor $\RanBetti: \fSetnu \to \tx{StratTop}$, where the stratification is by number of distinct points; we let $A_\Ran$ denote this stratification.
Similarly, the $\mb{C}$-points of $U_\Ran$ gives an oplax symmetric monoidal functor $\BettiU: \ECOp^\op \to \tx{StratTop}$. We let $A(p)$ denote its stratification for $p \in \ECOp^\op$.

\begin{definition}
We let $\tx{FactCat}^{\CoMark}_{\nonun, \ConstructibleMark, \LWeakMark}$ denote the full subcategory of
\[\tx{Fun}^{\tx{lax}, \tensor}(\fSetnu^\op, \tx{StratTop}^{\op, \tx{ShvMod}, \vee}) \times_{\tx{Fun}^{\tx{lax}, \tensor}(\fSetnu^\op, \tx{StratTop}^\op)} \{\RanBetti^\op\}\]
consisting of functors $F$ such that:
\begin{itemize}
\item[(Condition R)]For every $I \surj J$, the arrow $F(I) \to F(J)$ is coCartesian with respect to $\tx{StratTop}^{\op, \tx{ShvMod}, \vee} \to \tx{StratTop}^\op$;
\end{itemize}
\end{definition}

\begin{definition}
We let $\tx{FactCat}^{\CoMark}_{\nonun, \ConstructibleMark, \PartitionMark, \LWeakMark}$ denote the full subcategory of
\[\tx{Fun}^{\tx{lax}, \tensor}(\ECOp, \tx{StratTop}^{\op, \tx{ShvMod}, \vee}) \times_{\tx{Fun}^{\tx{lax}, \tensor}(\ECOp, \tx{StratTop}^\op)} \{\BettiU^\op\}\]
consisting of functors $F$ such that:
\begin{itemize}
\item[(Condition R')]For every $a \to b \in \ECOp$, the arrow $F(a) \to F(b)$ is coCartesian with respect to $\tx{StratTop}^{\op, \tx{ShvMod}, vee} \to \tx{StratTop}^\op$.
\end{itemize}
\end{definition}

The argument of Proposition \ref{prop:fact-cat-alg-partition-description} shows that we have a pair of functors
\[\tx{Conv}^\CoMark_{\PartitionMark \to \Ran}: \tx{FactCat}^\CoMark_{\nonun, \ConstructibleMark, \PartitionMark, \LWeakMark} \adjoint \tx{FactCat}^\CoMark_{\nonun, \ConstructibleMark, \LWeakMark} :\tx{Conv}^\CoMark_{\Ran \to \PartitionMark}\]
where, if $F$ is an element of the latter, then it gets sent to
\[p = (I \surj J) \mapsto F(J) \tensor_{\tx{cShv}_{A_\Ran}(\Ran_\nonun(J))} \tx{cShv}_{A(p)}(\BettiU(p));\]
and if $T$ is an element of the former, it gets sent to
\[J \mapsto \limit_{p = (J' \surj J)} T(p) \in \tx{cShv}_{A_\Ran}(\Ran_\nonun(J))\mod;\]
where each $T(p)$ is considered as a $\tx{cShv}_{A_\Ran}(\Ran_\nonun(J))$-module via $\dagger$-restriction.
\begin{proposition}
The two functors above are equivalences of categories.
\end{proposition}
\begin{proof}
Following the argument of Proposition \ref{prop:fact-cat-alg-partition-description} it suffices to check that the (co)unit yield identities. By conservativity it suffices to check it on the underlying (plain) categories.
For simplicity, set $C := \tx{cShv}_{A_\Ran}(\Ran_\nonun([1]))$ and $C_J := \tx{cShv}_{A(J \surj 1)}(\BettiU(J \surj 1))$.

On one hand, we need to check for every $F \in \tx{FactCat}_{\nonun, \ConstructibleMark}$, we have
\[F(1) \simeq \limit_{J} F(1) \tensor_{C} C_J;\]
passing to left adjoints and use the projection formula for closed embeddings, we reduce to showing that $C \to \limit_{J, \dagger\tx{-pull}} C_J$ is an equivalence. Let $C_{\le n} := \tx{cShv}_{A_{\Ran^{\le n}}}(\Ran_{\nonun}^{\le n})$, where $\Ran_{\nonun}^{\le n}$ is the topological subspace of $\Ran_{\nonun}$ consisting of points of cardinality at most $n$, and $A_{\Ran^{\le n}}$ is the obvious stratification. It was already proven by devissage in \cite[Theorem 2.19]{lejay2021constructible} that $C \simeq \limit_{n, \dagger\tx{-pull}} C_{\le n}$; by passing to left adjoints and combining colimit diagrams we see that it suffices to prove
\begin{equation}
\label{eq:cshv-induction-diagram}
C_{\le n} \to \limit_{J, |J| \le n, \dagger\tx{-pull}} C_J
\end{equation}
is an equivalence for each $n$. We run induction on $n$, the $n = 1$ case being evident. Let $C_{=n} := \ms{LS}(\Ran_\nonun^{=n})$, where $\Ran_{\nonun}^{=n} \subseteq \Ran_{\nonun}^{\le n}$ is the open subspace consisting of points of cardinality exactly $n$. For each $J$, let $C_{J}^{=} := \ms{LS}(U_{\BettiU(J \surj 1)}^{=})$ where $\BettiU(J \surj 1)^{=}$ is the open stratum of $U_{\BettiU(J \surj 1)}$, and let $\BettiU(J \surj 1)^{<}$ be the union of all closed proper strata of $\BettiU(J \surj 1)$, and set $C_J^{<} := \tx{cShv}_{A(J \surj 1)}(\BettiU(J \surj 1)^{<})$. Then we have
\[C_{\le n} \simeq \tx{lax} \limit [C_{=n} \to C_{\le n - 1}] \hspace{1em} C_J \simeq \tx{lax} \limit [C_J^{=} \to C_{J}^{<}]\]
where both functors are nearby cycles, and lax limits are taken with respect to the diagram $\bullet \to \bullet$. We also have $C_{=n} \simeq \limit_{J, |J| = n} C_J^{=}$ since $\Ran_\nonun^{=n}$ is the colimit of $\BettiU(J \surj 1)^{=}$. Noting that the diagram on the RHS of Equation \eqref{eq:cshv-induction-diagram} only goes from $|J| = n$ to $|J| < n$, we can also rewrite it as
\[\limit [(\limit_{J, |J| = n} C_J) \to (\limit_{K, |K| \le n - 1} C_K)]\]
Thus we rewrite Equation \eqref{eq:cshv-induction-diagram}, using induction hypothesis, as
\[\tx{lax} \limit [(\limit_{J, |J| = n} C_J^{=}) \to \limit_{K, |K| \le n - 1} C_K] \to \limit \left[\left(\limit_{J, |J| = n} (\tx{lax} \limit [C_J^{=} \to C_{J}^{<}]) \right) \to (\limit_{K, |K| \le n - 1} C_K)\right]\]
which can be seen to be an equivalence by unfolding definition.

On the other hand, we need to check for every $G \in \tx{FactCat}_{\nonun, \ConstructibleMark, \PartitionMark}$ and $I$, we have
\[G(I \surj 1) \simeq (\limit_{J} G(J \surj 1)) \tensor_{C} C_J\]
Again by projection formula, $C_J$ is dualizable as a $C$-module category, so we can swap the tensor inside. Similar to Proposition \ref{prop:fact-cat-alg-partition-description}, we now claim that the limit stabilizes to a constant as long as $|J| > |I|$. Indeed, since each $C_J$ fully faithfully embeds into $C$, the tensor product $G(J \surj 1) \tensor_{C} C_I$ in such range becomes $G(J \surj 1) \tensor_{C_J} C_I \simeq G(I \surj 1)$.
\end{proof}

\begin{definition}
\label{def:constructible-fact-objects}
We will use the name (non-unital) \emph{weak constructible} factorization category over $X$ to denote an element of
\[\tx{FactCat}^\CoMark_{\nonun, \ConstructibleMark, \PartitionMark, \LWeakMark}.\] (The category $\tx{FactCat}^\CoMark_{\nonun, \ConstructibleMark, \LWeakMark}$ will not play a role in our proof.) Similarly,
\begin{itemize}
\item for a fixed constructible factorization category $\mc{C}$, we can define the category
\[\tx{FactAlg}^\CoMark_{\nonun, \ConstructibleMark, \PartitionMark, \LWeakMark}(\mc{C})\]
of weak constructible factorization algebras internal to $\mc{C}$;
\item within said category, we can define the full subcategory of (strong) \emph{factorization} algebra internal to $\mc{C}$, where we require 
an additional condition:
\begin{center}
(condition $L_\tx{dual}$) For every $I \surj J \surj K \in \fSetnu^\op$, the arrow $\bigotimes_{k} F(I_k \surj J_k) \to F(I \surj J)$ is coCartesian with respect to $\tx{StratTop}^{\op, \tx{ShvModSect}, \vee} \to \tx{StratTop}^{\op, \tx{ShvMod}, \vee}$;
\end{center}
We denote this category by
\[\tx{FactAlg}^\CoMark_{\nonun, \ConstructibleMark, \PartitionMark}(\mc{C});\]
\item for a fixed point $x \in X(\mb{C})$, we also have the category of $\mc{C}$-constructible weak factorization module categories over $\Ran_x(\mb{C})$, to be denoted
\[\tx{FactModCat}^\CoMark_{\nonun, \ConstructibleMark, \PartitionMark, \LWeakMark, x}(\mc{C});\]
\item finally, fix $A \in \tx{FactAlg}^\CoMark_{\nonun, \ConstructibleMark, \PartitionMark}(\mc{C})$ and $\mc{M} \in \tx{FactModCat}^\CoMark_{\nonun, \ConstructibleMark, \PartitionMark, \LWeakMark, x}(\mc{C})$ we have the category of constructible weak factorization $A$-modules within $\mc{M}$, denoted
\[A\tx{-FactMod}^\CoMark_{\nonun, \ConstructibleMark, \PartitionMark, \LWeakMark, x}(\mc{M});\]
imposing an analogous condition ($L_\tx{dual}$), we obtain the category of (strong) factorization $A$-modules, denoted
\[A\tx{-FactMod}^\CoMark_{\nonun, \ConstructibleMark, \PartitionMark, x}(\mc{M}).\]
\end{itemize}
\end{definition}

\medskip

\subsection{How to Glue Cosheaves}

Recall the notion of \emph{lax} limits of categories and \emph{lax} functors between them (Section \ref{sect:lax-limits}). Constructible cosheaves constitute an important example of such:

\begin{lemma}
\label{lemma:cShv-glue}
	Let $\mc{C}$ a symmetric monoidal $\infty$-category, $(X \to A) \in \tx{StratTop}$, then $\tx{cShv}_A(X; \mc{C})$ is equivalent to the lax limit corresponding to $\{\ms{LS}(X_a; \mc{C})\}_{a \in A}$ where the connecting functors are $(j_b)^\dagger \circ (j_a)_\dagger$ for $a \to b$.
\end{lemma}
\begin{proof}
We check the conditions above are satisfied. Fully faithfulness is \cite[B.2.4]{knudsen2018higherv1}. Joint conservativeness is because isomorphism can be checked on costalks (thanks to hypercocompleteness). Assumption on $F_j^R \circ F_i$ is \cite[B.2.6]{knudsen2018higherv1}.
It remains to check that $j_a^\dagger: \tx{cShv}_A(X; \mc{C}) \to \ms{LS}(X_a; \mc{C})$ is continuous. By cosheaf condition, it suffices to compare
	\[\colimit_{\alpha} j_a^\dagger \mc{F}_\alpha(W) \to j_a^\dagger (\colimit_\alpha \mc{F}_\alpha)(W)\]
	for each open $F_\sigma$ set $W \subseteq X_\alpha$. By \cite[A.5.12]{HA} among all open sets of $X$ containing $W$ there exists a filtered cofinal subcategory consisting of spaces $\{Q_\beta\}$ which are homeomorphic to a fixed conical neighborhood $W \times C(L)$ of $W$. By homotopy invariance of constructible (co)sheaves \cite{haine2020homotopyinvariance}, the map above can then be identified with $\colimit_\alpha \mc{F}_\alpha(Q_\beta) \xrightarrow{\simeq} (\colimit_\alpha \mc{F}_\alpha)(Q_\beta)$ for any such $Q_\beta$.
\end{proof}

Recall from \cite[5.5.3.3]{HTT} the class of \emph{presentable} coCartesian fibrations can be characterized as those whose corresponding functors factor through $\tx{Pr}^L \subseteq \tx{Cat}$.
A morphism between two presentable coCartesian fibrations $p_1: E_1 \to S$ and $p_2: E_2 \to S$ is a morphism of coCartesian fibrations, such that its induced maps $(E_1)|_s \to (E_2)|_s$ are all continuous.

The following curious phenomenon is the basis of our construction:

\begin{proposition}
\label{prop:cShvA-presentable-coCart}
	Let $p: E \to D$ be a presentable coCartesian fibration. Then the canonical functor
	\[\tx{cShv}_A(X; E) \to \tx{cShv}_A(X; D)\]
	is a presentable coCartesian fibration.
\end{proposition}
\begin{proof}
We note that a morphism $\widetilde{F} \to \widetilde{G}$ upstairs is $\pi$-coCartesian if and only if for every $U \in \mc{U}(X)$, the morphism $\widetilde{F}(U) \to \widetilde{G}(U) \in \tx{Fun}([1], E)$ is a $p$-coCartesian arrow: by \cite[2.4.1.8]{HTT} this is evidently true for $\pi'$ in the above diagram, and thus is true for $\pi$. Now fix $e: (F \to G) \in \tx{Fun}([1], \tx{cShv}_A(X; D))$ and $\widetilde{F}$ over $F$. Our goal is to construct a coCartesian lift $\widetilde{e}$ above $e$.
	
	\medskip
	
	First we consider the case when $A$ is trivial. Then $\pi$ becomes the map $\limit_{X^\tx{Spc}} E \to \limit_{X^\tx{Spc}} D$. Then $e$ becomes a compatible family $\{e_x: F_x \to G_x\}_{x \in X^\tx{Spc}}$ of objects in $D$, and $\widetilde{F}$ becomes $\{\widetilde{F}_x\}$ a family of objects in $E$. There is an obvious map $\widetilde{e}: \{\widetilde{F}_x\} \to \{(e_x)_!(\widetilde{F}_x)\}$; we now prove this is $\pi$-coCartesian.
	
	Elements in $E$ can be expressed as $(x, d)$ where $d \in D$ and $x \in E|_{d}$. Recall that within $E$ (which is a lax limit), colimits are computed as 
	\[\colimit_i (x_i, d_i) \simeq (\colimit_i \tx{ins}_i(x_i), \colimit_i d_i)\]
	where $\tx{ins}_i$ is the coCartesian pushforward along $d_i \to \colimit_i d_i$. It follows if $\{(x_i, d_i)\}$ is a diagram in $E$ and $\{\varphi_i: d_i \to d'_i\}$ is a diagram of arrows in $D$, then the coCartesian lift of $\colimit_i \varphi_i$ above $\colimit_i (x_i, d_i)$ is
	\[(\colimit_i \tx{ins}^d_i(x_i), \colimit_i d_i) \to (\colimit_i \tx{ins}^{d'}_i((\varphi_i)_!(x_i)), \colimit_i d'_i);\]
	It thus follows from cosheaf condition that $\pi$-coCartesian-ness in $\tx{cShv}(X; E)$ can be checked just on disks in $X$. Since $\dagger$-fibers can be computed as a constant filtered diagram of sections on such disks, the same condition can be checked on $\dagger$-fibers, whence the claim.
	
	We need to check the fibration is presentable. To check fibers contain filtered colimits, it suffices to check that, if $c_i$ is a diagram upstairs that map to constant diagram $\underline{c}$ downstairs, then $\colimit_i c_i$ maps to the same $\underline{c}$, which follows from the fact that $\pi$ is continuous. That the connecting functors are continuous (i.e. coCartesian liftings preserve colimits) follows from the construction of $\widetilde{e}$.
	
	\medskip
	
	Now we consider general $A$-constructible cosheaves. Let $\Phi$ (resp. $\widetilde{\Phi}$) denote the lax-functors encoding $A$-constructible cosheaves valued in $D$ (resp. $E$) on $X$. There is an evident lax natural transformation $\ms{P}: \widetilde{\Phi} \to \Phi$ such that $\pi \simeq \ms{Glue}(\ms{P})$. Unfolding definition, we see that there is a family of functors
	\[\{\ms{P}_s: \mf{G}(\widetilde{\Phi})_s \to \mf{G}(\Phi)_s\}_{s \in \ms{sd}(A)}\]
	such that $\ms{Glue}(P) \simeq \limit_{s \in \ms{sd}(A)} \ms{P}_s$. Now, each $\ms{P}_s$ is of the form
	\[\tx{Fun}(K_s, \ms{LS}(X_{\ms{max}(s)}; E)) \to \tx{Fun}(K_s, \ms{LS}(X_{\ms{max}(s)}; D))\]
	and is component-wise induced by $\pi$ from the locally constant case. Given $e: (\{f_s\}_{s \in \ms{sd}(A)} \to \{g_s\}_{s \in \ms{sd}(A)}) \in \tx{Fun}([1], \ms{Glue}(\Phi))$ and a lift of $\{\tilde{f}_s\}$ $\{f_s\}$, then, we have a canonical lift $\widetilde{e}: \{\tilde{f}_s\} \to \{\tilde{g}_s\}$ given by component-wise coCartesian lift, using the locally constant case handled above. We declare this is a $\pi$-coCartesian lift. Replacing disks by the sets $Q_\beta$ considered in Lemma \ref{lemma:cShv-glue}, we see that $\pi$-coCartesian-ness can still be checked at $\dagger$-fibers, thus the claim follows from the construction of $\widetilde{e}$ and the previous case. That it is a presentable fibration again follows from unwinding the construction of $\widetilde{e}$.
\end{proof}

Let $\DGCatMarked$ denote the category $\DGCat_{\Vect/}$, i.e. the pair $(c, C)$ where $c \in C \in \DGCat$. The forgetful functor $p: \DGCatMarked \to \DGCat$ is continuous and symmetric monoidal, and is in fact the presentable coCartesian fibration associated with the tautological map $\DGCat \to \tx{Pr}^L$. We let
\[\Theta_{X, A}: \tx{cShv}_A(X; \DGCat) \to \tx{Pr}^L\]
denote functor corresponding to $\tx{cShv}_A(X; \DGCatMarked) \to \tx{cShv}_A(X; \DGCat)$.

\begin{example}
\label{example:LS-theta-value}
If $A$ is trivial and $\mc{F} \in \ms{LS}(X; \DGCat)$ is a locally constant cosheaf, then we have $\Theta_{X, \tx{triv}}(\mc{F}) \simeq \Gamma(X, \mc{F}) \in \DGCat$. Indeed, switching to the LS model, we can identify the fiber as
	\[\limit_{X^\tx{Spc}} \DGCatMarked \times_{\limit_{X^\tx{Spc}} \DGCat} \{x^\dagger(\mc{F})\} \simeq \limit_{x \in X^\tx{Spc}} x^\dagger(\mc{F}) \simeq \colimit_{x \in X^{\tx{Spc}, \op}} x^\dagger(\mc{F}),\]
	which is by definition the $\dagger$-pushforward to point.
\end{example}

\begin{example}	
In particular, let $\ms{LS}$ denote the cosheafification of the constant copresheaf $\bullet \mapsto \Vect$ (it is easy to see that $\ms{LS}(U) \simeq \tx{LS}(U)$, whence the name). Then $\Theta_{X, \tx{triv}}(\ms{LS}) \simeq \tx{LS}(X)$.
\end{example}

To describe $\Theta_{X, A}$ in more general cases we need some preparation.

\begin{lemma}
\label{lemma:fibration-pull-push-commute}
	In the above setting (in fact, for any map $E \to D$ of categories), for every $A', A'' \subseteq A$, we have a commutative squares
	\[\xymatrix{
\tx{cShv}_{A'}(X_{A'}; E) \ar^{(\widetilde{j_{A'}})_\dagger}[rr] \ar^{\pi_a}[d] & & \tx{cShv}_{A}(X; E) \ar^{\pi}[d] \ar^{\pi}[d] \ar^{(\widetilde{j_{A''}})^\dagger}[rr]  & & \tx{cShv}_{A''}(X_{A''}; E) \ar^{\pi_b}[d] \\
\tx{cShv}_{A'}(X_{A'}; D) \ar^{(j_{A'})_\dagger}[rr] & & \tx{cShv}_\tx{A}(X; D) \ar^{(j_{A''})^\dagger}[rr] & & \tx{cShv}_{A''}(X_{A''}; D)
	}\]
	such that the top arrows preserve coCartesian arrows.
\end{lemma}

\begin{proof}
	The first square exists and commutes by definition. Passing to the adjoint horizontally, we have a natural transformation and we reduce to a check on $\dagger$-fibers; the homotopy invariance argument from Lemma \ref{lemma:cShv-glue} then shows that the diagram actually commutes. The second claim follow similarly.%
\end{proof}

It follows that we have a natural transformation $\eta_{ab}: \Theta_{X_a, \tx{triv}} \Rightarrow \Theta_{X_b, \tx{triv}} \circ (j_b)^\dagger (j_a)_\dagger$ (note that $\dagger$-pullback is continuous for constructible sheaves).
Now, let $F = \{F_a, \varphi_{ab}\}$ denote an element of $\ms{Glue}(\Phi) \simeq \tx{cShv}_A(X; D)$ where $F_i \in \tx{cShv}_\tx{lc}(X_a; D)$ and $\varphi_{ab}: j_b^\dagger \circ (j_a)_\dagger(F_a) \to F_b$. We have a continuous functor
\[\theta_{ab} := \Theta_{X_a, \tx{triv}}(F_a) \xrightarrow{\eta_{ab}} \Theta_{X_b, \tx{triv}}((j_b)^\dagger (j_a)_\dagger F_a) \xrightarrow{\varphi_{ab}} \Theta_{X_b, \tx{triv}}(F_b).\]

\begin{proposition}
\label{prop:theta-from-gluing}
$\{\Theta_{X_a, \tx{triv}}(F_a), \theta_{ab}\}$ combine into a lax functor $\Xi_F: A \dashrightarrow \tx{Cat}$ such that $\Theta_{X, A}(F) \simeq \ms{Glue}(\Xi_F)$.
\end{proposition}
\begin{proof}
	By Proposition \ref{prop:glue-cat-reconstruction}, we need to check the following conditions:
	\begin{enumerate}
	\item There exists fully faithful continuous functors $\phi_a: \Theta_{X_a, \tx{triv}}(F_a) \to \Theta_{X, A}(F)$, whose right adjoints are continuous and jointly conservative, and an equivalence $\theta_{ab} \simeq \phi_b^R \circ \phi_a$;
	\item We have $\theta_{ab}(\bullet) \simeq \tx{init}_{\Theta_{X_b, \tx{triv}}(F_b)}$ for every $i \not\to j$ in $I$.
	\end{enumerate}
	The second condition is evident because $(\widetilde{j_b})^\dagger (\widetilde{j_a})_\dagger$ preduces the initial precosheaf. If we consider $F$ as an element of $\tx{cShv}_A(X; D)$, then each $F_a \simeq (j_a)^\dagger F$. The first diagram above then gives a functor
	\[\phi_a: \Theta_{X_a, \tx{triv}}((j_a)^\dagger F) \xrightarrow{\eta_\dagger} \Theta_{X, A}((j_a)_\dagger (j_a)^\dagger F) \xrightarrow{\tx{counit}} \Theta_{X, A}(F)\]
	From the second square above, we also obtain a natural transformation $\eta^\dagger: \Theta_{X, A} \Rightarrow \Theta_{X_a, \tx{triv}} \circ (j_a)^\dagger$, hence a (continuous) functor
	\[\psi_a := \Theta_{X, A}(F) \xrightarrow{\eta^\dagger} \Theta_{X_a, \tx{triv}}((j_a)^\dagger F);\]
	We claim that $\psi_a$ is the right adjoint of $\phi_a$, and $\psi_a \circ \phi_a \simeq \tx{id}$.
	This would prove the original claim, since joint conservativity follows from Proposition \ref{prop:lax-limit-properties}, and the equivalence $\theta_{ab} \simeq \varphi_b^R \circ \varphi_a$ follows from the functoriality of $\eta^\dagger$.
	
	First, note that the composition $\psi_a \circ \phi_a$ is computed by first send an object in the fiber over $F_a$ along $(\widetilde{j_a})_\dagger (\widetilde{j_a})^\dagger$ then coCartesian push along the arrow $(j_a)^\dagger \to (j_a)^\dagger (j_a)_\dagger (j_a)^\dagger F$ in the base; however, the latter arrow is identity, and the first operation is fully faithful.
	To prove adjointness, it remains to exhibit a counit $\phi_a \circ \psi_a \Rightarrow \tx{id}$ and check the adjoint axioms. If we denote an object in $\Theta_{X, A}(F)$ by $(F, x)$, then this is a functorial morphism
	\[(F, \tx{counit}_!((\widetilde{j_a})_\dagger (\widetilde{j_a})^\dagger x)) \to (F, x)\]
	where $\tx{counit}_!$ is the coCartesian push along $(j_a)_\dagger (j_a)^\dagger F \to F$. By coCartesian-ness this is simply a morphism $(\widetilde{j_a})_\dagger (\widetilde{j_a})^\dagger x) \to x$ within the total category, which is given by the counit upstairs. The adjoint axioms are checked by unfolding definitions.
\end{proof}

\begin{lemma}
	\label{lemma:transition-theta-nearby-cycle}
	For $(X \to A) \in \tx{StratTop}$, the map $\theta_{ab}: \Theta_{X_a, \tx{triv}}(\ms{LS}) \to \Theta_{X_b, \tx{triv}}(\ms{LS})$ is given by
	\[\ms{LS}(X_a) \xrightarrow{(j_a)_\dagger} \tx{cShv}_A(X; \Vect) \xrightarrow{(j_b)^\dagger} \ms{LS}(X_b).\]
	From Proposition \ref{prop:glue-cat-reconstruction} we thus have $\Theta_{X, A}(\ms{LS}) \simeq \tx{cShv}_A(X; \Vect)$.
\end{lemma}

\begin{proof}
Unfolding definition using the enter-path model, we see that the functor
\[\ms{LS}(X_a) \xrightarrow{\eta_{ab}} \Theta_{X_b, \tx{triv}}((j_b)^\dagger (j_a)_\dagger(\ms{LS}|_{X_a})) \xrightarrow{\varphi_{ab}} \ms{LS}(X_b)\]
identifies with
\[\xymatrix{
\tx{Fun}_{/X_a^{\tx{Spc}, \op}}(X_a^{\tx{Spc}, \op}, X_a^{\tx{Spc}, \op} \times_{\tx{Entr}_A(X)} \tx{Un}(\tx{LKE}_{\iota_a} \tx{Res}_{\iota_a} C_\Vect)) \ar^(0.55){\tx{LKE}}[r] & \tx{Fun}_{/\tx{Entr}_{A}(X)}(\tx{Entr}_{A}(X), \tx{Un}(\tx{LKE}_{\iota_a} \tx{Res}_{\iota_a} C_\Vect)) \ar^{\eta_a}[d] \\
 \tx{Fun}_{/X_a^{\tx{Spc}, \op}}(X_a^{\tx{Spc}, \op}, \tx{Un}(\tx{Res}_{\iota_a} C_\Vect)) \ar@{-->}[r] \ar^{\simeq}[u]
& \tx{Fun}_{/\tx{Entr}_{A}(X)}(\tx{Entr}_{A}(X), \tx{Un}(C_\Vect)) \ar^{\tx{res}}[d] \\
& \tx{Fun}_{/X_b^\tx{Spc, \op}}(X_b^{\tx{Spc}, \op}, \tx{Un}(\tx{Res}_{\iota_b} C_\Vect))
}\]
where $C_\Vect: \tx{Entr}_A(X) \to \DGCat$ is the constant functor mapping everything to $\Vect$, $\iota_a: X_a^{\tx{Spc}, \op} \to \tx{Entr}_A(X)$ the inclusion, $\tx{Un}$ denotes the corresponding coCartesian fibration, and $\eta_a$ is given by the universal property of LKE. It suffices to show that the dashed arrow is given by $(j_a)_\dagger$. Flipping the left vertical arrow, it suffices to prove the commutativty of
\[\xymatrix{
\tx{Fun}_{/\tx{Entr}_A(X)}(X_a^{\tx{Spc}, \op}, \tx{Un}(\tx{LKE}_{\iota_a} \tx{Res}_{\iota_a} C_\Vect)) \ar^{\tx{LKE}}[r] \ar^{\eta_a}[d] & \tx{Fun}_{/\tx{Entr}_{A}(X)}(\tx{Entr}_{A}(X), \tx{Un}(\tx{LKE}_{\iota_a} \tx{Res}_{\iota_a} C_\Vect)) \ar^{\eta_a}[d] \\
 \tx{Fun}_{/\tx{Entr}_A(X)}(X_a^{\tx{Spc}, \op}, \tx{Un}(C_\Vect)) \ar^{\tx{LKE}}[r] 
& \tx{Fun}_{/\tx{Entr}_{A}(X)}(\tx{Entr}_{A}(X), \tx{Un}(C_\Vect))
}\]
which is because $\eta_a$ is pointwise colimit-preserving.
\end{proof}

\subsection{Constructible Action}

Recall that a symmetric monoidal coCartesian fibration is a coCartesian fibration that is also a symmetric monoidal functor, and such that the coCartesian arrows are stable under tensor product (\cite[3.3.4.6]{gaitsgory2019weil}). By \cite[3.3.5.17]{gaitsgory2019weil}, such a data $p: \mc{E} \to \mc{B}$ corresponds to a \emph{lax} symmetric monoidal functor $\mc{B} \to \tx{Cat}_\infty$, where the latter is equipped with the Cartesian monoidal structure. It is easy to see that if $p$ is a presentable fibration, then this functor descends to a lax symmetric monoidal functor $\mc{B} \to \tx{Pr}^L$. %

The map $p: \DGCatMarked \to \DGCat$ is a symmetric monoidal presentable coCartesian fibration, hence by Proposition \ref{prop:cShvA-presentable-coCart} and Remark \ref{remark:cShvA-sym-monoidal} the situation above applies. Recall that the unit for $\tx{cShv}(X; \DGCat)$ is $\ms{LS}$, being the cosheafification of the precosheaf unit. Thus we see that for any $F \in \tx{cShv}_A(X; \DGCat)$, $\Theta_{X, A}(F)$ has a natural $\tx{cShv}_A(X; \Vect)$-module structure, i.e. we have an enhanced functor
\[\Theta_{X, A}: \tx{cShv}_A(X; \DGCat) \to \tx{cShv}_A(X)\mod(\DGCat)\]
which we denote by the same symbol.

In general, the $\dagger$-pullback for $\DGCat$-valued constructible cosheaves is \emph{incompatible} with base-change as $\tx{Shv}_A$-modules. The following is thus a nontrivial property of this action:

\begin{proposition}
\label{prop:ShvA-action-base-change}
	Let $(X \to A) \in \tx{StratTop}$, and let $A' \subseteq A$ be either upward-closed or an initial element $\{0\}$, then for any $F \in \tx{cShv}_A(X; \DGCat)$ we have
	\[\Theta_{X, A}(F) \tensor_{\tx{cShv}_A(X)} \tx{cShv}_{A'}(X_{A'}) \simeq \Theta_{X, A'}(j_{A'}^\dagger F).\]
\end{proposition}

\begin{proof}

For simplicity let $X' := X_{A'}$. First we need a fact:
\begin{lemma}
\label{lemma:cShv-sub-self-dual}
$\tx{cShv}_{A'}(X')$ is self-dual as a $\tx{cShv}_A(X)$-module category.
\end{lemma}
\begin{proof}
	Note that $\tx{cShv}_A(X)$ is self-dual as a module over itself, so it suffices to demonstrate $\tx{cShv}_{A'}(X')$ as a retract of $\tx{cShv}_A(X)$ within $\tx{cShv}_A(X)\mod$. If $A'$ is an upward-closed subposet i.e. $i: X' \subseteq X$ is a closed embedding, we check that the adjoint pair $i_\dagger: \tx{cShv}_{A'}(X') \adjoint \tx{cShv}_A(X): i^\dagger$ (where the left adjoint is fully faithful by Proposition \ref{prop:lax-limit-properties}) is $\tx{cShv}_A(X)$-linear, which follows from the explicit description given there. When $A' = \{0\}$ i.e. $j: X' \subseteq X$ is an open embedding), since we are in the presentable stable setting, we can instead use the existence of a \emph{right} adjoint $j_*$ to $j^\dagger$ (which can be found in \cite[Proposition 6.1.6]{ayala2021stratified}), which sets every other component of the lax limit to the zero object; the rest is the same.
\end{proof}
By Proposition \ref{prop:theta-from-gluing} we can rewrite the equation as
\[\ms{Glue}_{a \in A}(\Theta_{X_a, \tx{triv}}(j_a^\dagger(F))) \tensor_{\tx{cShv}_A(X)} \tx{cShv}_{A'}(X') \xrightarrow{?} \ms{Glue}_{a' \in A'}(\Theta_{X_{a'}, \tx{triv}}(j_{a'}^\dagger(F)))\]
Using duality we further expand the left as a cosimplicial limit
\[\limit_\bullet \tx{LFun}(\tx{cShv}_{A'}(X') \tensor \tx{cShv}_A(X)^{\tensor \bullet}, \ms{Glue}_{a \in A}(\Theta_{X_a, \tx{triv}}(j_a^\dagger(F))))\]
Combining Lemma \ref{lemma:lax-limit-functorialities} and Lemma \ref{lemma:glue-as-strict-limit} we rewrite it as
\[\ms{Glue}_{a \in A}(\limit_\bullet \tx{LFun}(\tx{cShv}_{A'}(X') \tensor \tx{cShv}_A(X)^{\tensor \bullet}, \Theta_{X_a, \tx{triv}}(j_a^\dagger(F))))\]
which is
\[\ms{Glue}_{a \in A}(\tx{cShv}_{A'}(X')_{\tensor \tx{cShv}_A(X)} \Theta_{X_a, \tx{triv}}(j^\dagger_a(F))).\]
Thus we reduce to proving that $\tx{cShv}_{A'}(X') \tensor_{\tx{cShv}_A(X)} \Theta_{X_a, \tx{triv}}(j^\dagger_a(F))$ is $\Theta_{X_a, \tx{triv}}(j_a^\dagger(F))$ if $a \in A'$, and the empty category otherwise. Note that we can rewrite the expression as
\[\tx{cShv}_{A'}(X') \tensor_{\tx{cShv}_A(X)} \ms{LS}(X_a) \tensor_{\ms{LS}(X_a)} \Theta_{X_a, \tx{triv}}(j^\dagger_a(F))\]
thus we reduce to the case of $\ms{LS}(X_a)$, which is elementary.
\end{proof}

For $p = (I \surj J) \in \ECOp^\op$, let $j_{p}: \BettiU(p) \to \Ran(\mb{C})$ denote the inclusion. We are now ready to prove

\begin{proposition}
\label{prop:E2-to-fact-cat}
	Given $F \in \tx{Fact}(X(\mb{C}))\tx{-Alg}(\DGCat)$, the assignment
	\[p \in \ECOp^\op \mapsto \Theta_{\BettiU(p), A(p)}(j_{p}^\dagger F)\]
	is a constructible factorization category over $X(\mb{C})$. Combined with \cite[5.5.4.10]{HA}, we obtain a functor
	\[\tx{Fact}: \mb{E}_{X(\mb{C}), \nonun}\tx{-Alg}(\DGCat) \to \tx{FactCat}^\CoMark_{\nonun, \ConstructibleMark, \PartitionMark, \WeakMark}.\]
\end{proposition}

\begin{proof}

For any symmetric monoidal category $\mc{C}$ and any smooth manifold $M$, we let
\[\tx{Fact}(M)\tx{-Alg}^\tx{fact}(\mc{C}) \subseteq \tx{Fact}(M)\tx{-Alg}(\mc{C})\]
denote the essential image of $\Psi$ introduced in \cite[5.5.4.10]{HA}. Elements of this subcategory will be called \emph{factorizable} $\mc{C}$-valued $\tx{Fact}(M)$-algebra.

We will construct (A) and (B) below
\[\mb{E}_{X(\mb{C}), \nonun}\tx{-Alg}(\DGCat) \xrightarrow[\simeq]{\Psi} \tx{Fact}(X(\mb{C}))\tx{-Alg}^\tx{fact}(\DGCat) \xrightarrow[\simeq]{(A)} \tx{FactCat}^{\tx{cShvCat}, \CoMark}_{\nonun, \ConstructibleMark, \WeakMark}\]
\[\xrightarrow[\simeq]{\tx{Conv}^\CoMark_{\Ran \to \PartitionMark}} \tx{FactCat}^{\tx{cShvCat}, \CoMark}_{\nonun, \ConstructibleMark, \PartitionMark, \WeakMark} \xrightarrow{(B)} \tx{FactCat}^\CoMark_{\nonun, \ConstructibleMark, \PartitionMark, \WeakMark}.\]

For (A), let $\tx{Fact}(X(\mb{C}))_{n \to 1}$ denote the full subcategory of coCartesian arrows in $\tx{Fact}(X(\mb{C}))^{\tensor}$ that lie above the active arrow $[n] \to [1] \in \fSetnu$. Consider the functor
\[\mathscr{P}: \fSetnu^\op \to 1\tx{-Cat} \hspace{1em} [n] \mapsto \tx{Fun}(\tx{Fact}(X(\mb{C}))_{n \to 1}, \DGCat)\]
where for $[n] \to [m] \in \fSetnu^\op$, the map is given by composing with $\tx{Fact}(X(\mb{C}))_{m \to 1} \to \tx{Fact}(X(\mb{C}))_{n \to 1}$, given by the factorization of coCartesian arrows. Let $\tx{Un}(\mathscr{P})$ denote the corresponding coCartesian fibration. Moreover, every $[n_1] \sqcup [n_2] \simeq [n]$ gives rise to a map
\[\tx{Fact}(X(\mb{C}))_{n \to 1} \to \tx{Fact}(X(\mb{C}))_{n_1 \to 1} \times \tx{Fact}(X(\mb{C}))_{n_2 \to 1}\]
so box product equips $\tx{Un}(\mathscr{P})$ a symmetric monoidal structure.
We have a naturally defined functor
\[\eta: \tx{Fact}(X(\mb{C}))\tx{-Alg}(\DGCat) \to \tx{Fun}^{\tx{lax}, \tensor}(\fSetnu^\op, \tx{Un}(\mathscr{P}))\]
given by sending a $\tx{Fact}(X(\mb{C}))$-algebra $F$ and $[n] \in \fSetnu^\op$ to $(n, (p_0 \to p_1) \mapsto F(p_1))$; the lax symmetric monoidal structure is given by the map $F(U) \tensor F(V) \to F(U \star V)$ given by the algebra structure of $F$.
On the other hand, note that $\tx{FactCat}^{ \tx{cShvCat}, \CoMark}_{\nonun, \ConstructibleMark, \WeakMark}$ embeds into the RHS, because each $\tx{Fact}(X(\mb{C}))_{n \to 1}$ provides a basis for the topology of $\Ran([n])$. It is also clear that its essential image coincides with that of $\eta$: both consists of sections that satisfy cosheaf and factorization conditions. This provides the desired equivalence (A).

For (B), we first need the fact that $\Theta_{X, A}$ is \emph{lax} functorial in $(X \to A)$; i.e. for every $f: (X \to A) \to (Y \to B)$ we have a map of $\tx{cShv}_A(X)$-modules
\begin{equation}
\label{eq:theta-base-change}
\Theta_{Y, B}(-) \tensor_{\tx{cShv}_B(Y)} \tx{cShv}_A(X) \to \Theta_{X, A}(f^\dagger(-)).
\end{equation}
First, we note that the diagram
\[\xymatrix@C=1em{
\tx{cShv}_B(Y; \DGCatMarked) \ar[rr] \ar[rd] & & \tx{cShv}_A(X; \DGCatMarked) \times_{\tx{cShv}_A(X; \DGCat)} \tx{cShv}_B(Y; \DGCat) \ar[ld] \\
 & \tx{cShv}_B(Y; \DGCat)
}\]
coming from the left square in Lemma \ref{lemma:fibration-pull-push-commute} upgrades to a diagram of symmetric monoidal functors, where the left arrow is given by Remark \ref{remark:cShvA-sym-monoidal}. (Indeed, pullback of cosheaves is symmetric monoidal, so the right top corner is a limit of symmetric monoidal categories.) It follows that we have a map of $\tx{cShv}_B(Y)$-module categories $\Theta_{Y, B}(-) \to \Theta_{X, A}(f^\dagger(-))$, where $\tx{cShv}_B(Y)$ acts on the right via the map $f^\dagger: \tx{cShv}_B(Y) \to \tx{cShv}_A(X)$; the claim above then follows by passing to left adjoint.

Thus we have a map $\tilde{\Theta}: \tx{StratSpc}^{\tx{cShvCat}, \vee} \to \tx{StratSpc}^{\tx{ShvMod}, \vee}$ over $\tx{StratSpc}^\op$. Moreover, for $X \to A, Y \to B$, the diagram
\[\xymatrix{
\tx{cShv}_A(X; \DGCatMarked) \times \tx{cShv}_B(Y; \DGCatMarked) \ar[r] \ar[d] & \tx{cShv}_{A \times B}(X \times Y; \DGCatMarked) \ar[d] \\
\tx{cShv}_A(X; \DGCat) \times \tx{cShv}_B(Y; \DGCat) \ar[r] & \tx{cShv}_{A \times B}(X \times Y; \DGCat)
}\]
is commutative; it follows that $\tilde{\Theta}$ is \emph{lax} symmetric monoidal. The existence of (B) now follows Proposition \ref{prop:ShvA-action-base-change} above.
\end{proof}

\subsubsection{Variant: Module Categories}
The concept of operadic module is defined in \cite[Chapter 3.3]{HA}. Fix $\mc{C}$ an $\mb{E}_M$-category and $A \in \mb{E}_M\tx{-Alg}_\nonun(\mc{C})$. Fix a point $x \in M$ and a particular disk embedding $(D_x \subseteq M) \in (\mb{E}_M^{\tensor})_{\langle 1 \rangle}$, and let $\mb{E}_{D_x}$ denote the full subcategory of $\mb{E}_M^{\tensor}$ consisting of copies of $(D_x \subseteq M)$. Let $\mc{K}'$ be the category defined as
\[\xymatrix{
\mc{K}' \ar[r] \ar[d] & \mc{K}_{\mb{E}_M^{\tensor}} \ar^{e_0, e_1}[d] \\
\mb{E}_{D_x} \times \mb{E}_M^{\nonun, \tensor} \ar[r] & \mb{E}_M^{\tensor} \times \mb{E}_M^{\tensor}
}\]
and let
\[\tx{Mod}_{A, x}^{\mb{E}_M, \nonun}(\mc{C}) \subseteq \tx{Fun}_{\mb{E}_M^{\nonun, \tensor}}(\mc{K}', \mc{C}^{\tensor}) \times_{\tx{Fun}_{\mb{E}_M^{\nonun, \tensor}}(\mb{E}_M^{\nonun, \tensor}, \mc{C}^{\tensor})} \{A\}\]
be the full subcategory of maps that preserve inert arrows (the map $\mb{E}_M^{\nonun, \tensor} \inj \mc{K}'$ is the inclusion of null arrows c.f. \cite[3.3.3.2]{HA}). By definition, the category is an operad over $N(\tx{Surj}_*)$, whose fiber over $[k]$ consists of all $k$-tuples of non-unital $A$-$\mb{E}_M$-modules inside $\mc{C}$ and $k$-tuples of non-unital $\mb{E}_M$-module morphisms between them.

Define a category $\tx{Fact}^\ms{Mod}(M)_x^{\tensor}$ as follows:
\begin{itemize}
\item Its objects are pairs
\[(\varphi: I \subseteq [n] \in \Setli, (U_1, \ldots, U_n))\]
where $I$ is a (potentially empty) subset, each $U_j$ such that $j \in \im(\varphi)$ is an open subseteq of $\Ran(M)_x$, and is otherwise an open subset of $\Ran(M)$;
\item A morphism
\[(\varphi: I \subseteq [n], (U_1, \ldots, U_n)) \to (\psi: J \subseteq [m], (V_1, \ldots, V_m))\]
only exists when $|I| = |J|$, in which case it is a surjection $\alpha: [n] \surj [m]$ that restricts to a bijection $I \simeq J$, and such that for each $j \in J$, $\{U_i\}_{\alpha(i) = j}$ are pairwise disjoint, and that $\star_{\alpha(i) = j} U_i \subseteq V_j$.
\end{itemize}
The category $N(\tx{Fact}^\ms{Mod}(M)_x^{\tensor})$ is then an $\infty$-operad over $\ms{CM}_{\nonun}^{\tensor}$. A straightforward extension of \cite[5.5.4.10]{HA} (where one uses \cite[Lemma 6.6 and Proposition 2.38]{ayala2019factorization} to prove the fully faithfulness) allows us to conclude:

\begin{proposition} %
	The category $\tx{Mod}_{A, x}^{\mb{E}_M, \nonun}(\mc{C})$ identifies with the  full subcategory of
	\[\tx{Alg}_{\tx{Fact}^\ms{Mod}(M)_x}(\mc{C}) \times_{\tx{Alg}_{\tx{Fact}(M)}(\mc{C})} \{\Psi(A)\}\]
	spanned by objects $F$ satisfying the following conditions:
\begin{enumerate}
\item the restriction to the objects in $\tx{Fact}^\ms{Mod}(M)_{x, \langle 1 \rangle}$ containing a marked disk gives a constructible cosheaf on $\Ran(M)_x$; and
\item the map $F(U) \tensor F(V) \to F(U \star V)$ is an isomorphism for all $U, V$ independent.
\end{enumerate}
\end{proposition}
The same argument as before then gives that
\begin{corollary}
\label{cor:E2Mod-to-FactModCat}
	For any $\mb{E}_{X(\mb{C})}$ DG category $\mc{C}$, There exists a canonically defined functor
	\[\tx{Fact}_\tx{ModCat}: \mc{C}\mod^{\mb{E}_{X(\mb{C})}, \nonun}(\DGCat) \to \tx{Fact}(\mc{C})\tx{-FactModCat}^\CoMark_{\nonun, \ConstructibleMark, \PartitionMark, \WeakMark, x}.\]
\end{corollary}

\subsubsection{Algebras and Modules}
Now we compare internal $\mb{E}_{X(\mb{C})}$-algebras and modules with the corresponding constructible factorization concepts. Fix an $\mb{E}_{X(\mb{C})}$-DG category $\mc{C}$ and let its corresponding factorizable cosheaf of DG categories be $\Psi(\mc{C})$. Note that by definition, the category of $\mb{E}_{X(\mb{C})}$-algebras internal to $\mc{C}$ is equivalent to the category 
\[\tx{Fact}(X(\mb{C}))\tx{-Alg}^\tx{fact}(\DGCatMarked) \times_{\tx{Fact}(X(\mb{C}))\tx{-Alg}^\tx{fact}(\DGCat)} \{\Psi(\mc{C})\}.\]
The argument in the proof of Proposition \ref{prop:E2-to-fact-cat} shows that this is equivalent to a system of objects $\{c_p\}$ in $\tx{cShv}_{A(p)}(\BettiU(p); \DGCatMarked)$ above the elements in $\tx{cShv}_{A(p)}(\BettiU(p); \DGCat)$ corresponding to $j_p^\dagger(\Psi(\mc{C}))$, which are just objects in each $\Theta_{\BettiU(p), A(p)}(j^\dagger_p(\Psi(\mc{C})))$. Moreover, that these $c_p$ correspond to each other under $\dagger$-pullbacks is tautologically equivalent to the statement that they map to each other under the map in Equation \ref{eq:theta-base-change}. The module case is analogous. We are thus led to the following fact:
\begin{proposition}
\label{prop:E2Alg-Mod-to-FactAlg-Mod}
	For any $\mb{E}_{X(\mb{C})}$ DG category $\mc{C}$, we have an equivalence of categories
	\[\tx{Fact}_\tx{Alg}: \mb{E}_{X(\mb{C})}\tx{-Alg}_\nonun(\mc{C}) \simeq \tx{FactAlg}^\CoMark_{\nonun, \ConstructibleMark, \PartitionMark}(\tx{Fact}(\mc{C}));\]
	for any $\mb{E}_{X(\mb{C})}$-algebra $A$ inside $\mc{C}$, a fixed point $x \in X(\mb{C})$ and any $\mb{E}_{X(\mb{C})}$ module category $\mc{M}$ of $\mc{C}$, we have an equivalence of categories
	\[\tx{Fact}_\tx{Mod}: A\mod^{\mb{E}_{X(\mb{C}), \nonun}}(\mc{M}) \simeq \tx{Fact}_\tx{Alg}(A)\tx{-FactMod}^\CoMark_{\nonun, \ConstructibleMark, \PartitionMark, x}(\tx{Fact}_{\tx{ModCat}}(\mc{M})).\]
\end{proposition}

\begin{remark}[Case of $\mb{A}^1$]
\label{remark:A1-sees-no-ribbon}
When our curve is $\mb{A}^1$ (which is what we use in the main paper), the operad $\mb{E}_{\mb{A}^1(\mb{C})}$ is isomorphic to $\mb{E}_2$ (\cite[5.4.5.3]{HA}), so LHS of Proposition \ref{prop:E2Alg-Mod-to-FactAlg-Mod} becomes (non-unital) $\mb{E}_2$-algebras and $\mb{E}_2$-modules in the usual sense. However, we avoid utilizing this trivialization until the comparison with Drinfeld center (Proposition \ref{prop:HH-is-E_n-mod}), as it obfuscates the $SO(2)$-equivariance necessary for \emph{combinatorially} constructing the relevant factorization algebras.
\end{remark}

\section{Quantum Torus}

\label{sect:q-torus}

\begin{notation}
In this subsection, we use the notation $\tb{K}(G, n)$ to denote the homotopy type of the Eilenberg-Maclane space $K(G, n)$.
\end{notation}

Let $T$ be a torus over $\mb{C}$ with weight lattice $\weightLat$ and coweight lattice $\coweightLat$; let $\cT$ be the Langlands dual torus. Consider the mapping space
\[\mc{P}_\tx{torus} := \tx{Maps}_{\mb{E}_0\tx{-Alg}(\tx{Spc})}(B(\cT(\mb{C}))^\tx{Spc}, \tb{K}(\mb{C}^\times, 4))\]
This is a strict 2-groupoid classifying braided monoidal 1-categories known as \emph{quantum tori} associated with $T$. It has two homotopy invariants (\cite[Appendix C]{gaitsgory2018parameters}):
\[\pi_0(\mc{P}_\tx{torus}) = \tx{Quad}(\weightLat), \pi_2(\mc{P}_\tx{torus}) = \Hom(\weightLat, \mb{C}^\times)\]
where the former is the set of quadratic forms $\weightLat \times \weightLat \to \mb{C}^\times$.

Let us fix a quadratic form $q \in \pi_0(\mc{P}_\tx{torus})$ and let $\Rep_q(T)$ denote the $\mb{E}_2$-DG Category defined in Definition \ref{def:quantum-torus}.
By definition, as an $\mb{E}_1$ DG category, $\Rep_q(T)$ is isomorphic to the $\mb{E}_1$ DG category underlying the symmetric monoidal DG category $\Rep(T)$. It follows from \cite[5.1.3.2]{HA} that the $\mb{E}_2$ operations are $t$-exact with respect to the standard $t$-structure, so that the $\mb{E}_2$ category $\Rep_q(T)^\heartsuit$ form a full $\mb{E}_2$ DG subcategory. By \cite[5.1.2.4]{HA}, a $\mb{E}_2$-structure on a Grothendieck abelian category is a (classical) braided monoidal structure on the underlying 1-category (here Convention \ref{conv:C-is-R2} is invoked), and the corresponding braided monoidal structure on $\Rep_q(T)^\heartsuit$ is the same as what we described above.

Crucially, the above conditions actually uniquely identifies $\Rep_q(T)$:

\begin{proposition}
\label{prop:RepqT-is-unique}
Let $\mc{E}$ be an $\mb{E}_2$ DG category. Suppose we are given an equivalence of $\mb{E}_1$ DG categories $\alpha: \mc{E} \simeq \Rep(T)$. Then $\mc{E}$ inherits a $t$-structure from $\Rep(T)$ and, by above, $\mc{E}^{\heartsuit}$ inherits a braided monoidal structure. Suppose we are also given a braided monoidal equivalence $\beta: \mc{E}^{\heartsuit} \simeq \Rep_q(T)^{\heartsuit}$. Then to $\alpha$ and $\beta$ we can canonically attach an equivalence of $\mb{E}_2$ DG categories $\mc{E} \simeq \Rep_q(T)$.
\end{proposition}

\begin{proof}
Let us introduce some notations:
\begin{itemize}
\item Let $\tx{Groth}_\infty$ denote the category of all Grothendieck prestable categories (defined in \cite[C.3.0.5]{lurie2018spectral}); it comes with a symmetric monoidal structure by \cite[C.4.2.1]{lurie2018spectral}. Let $\tx{Groth}_\infty^\tx{lex}$ denote the 1-full subcategory where we require morphisms to be left exact, and let $\tx{Groth}_\infty^{\tx{sep}, \tx{lex}}$ denote the full subcategory consisting of separated prestable categories (\cite[C.1.2.12]{lurie2018spectral});
\item Let $\tx{Groth}_\tx{ab}^\tx{ex}$ denote the category of Grothendieck abelian categories and exact functors between them;
\item Let $\tx{Groth}_\infty^{+, \tx{lex}}$ be the category of pairs $(\mc{C}, \mc{C}_{\ge 0})$, where $\mc{C} \in \tx{Pr}^{\tx{st}}$ and $\mc{C}_{\ge 0}$ is a full subcategory closed under colimits and extensions. A morphism $(\mc{C}, \mc{C}_{\ge 0}) \to (\mc{C}', \mc{C'}_{\ge 0})$ is a continuous functor $\mc{C} \to \mc{C'}$ that is left $t$-exact, and restricts to a left exact functor $\mc{C}_{\ge 0} \to \mc{C}'_{\ge 0}$.
\end{itemize}
	By \cite[C.5.4.5]{lurie2018spectral}, we have a \emph{fully faithful} map
	\[\mc{D}_{\ge 0}: \tx{Groth}_\tx{ab}^\tx{ex} \to \tx{Groth}_\infty^{\tx{sep}, \tx{lex}},\]
	given by sending an Grothendieck abelian category $\mc{H}$ to $\mc{D}(\mc{H})_{\ge 0}$ defined in \cite[C.1.4.5]{lurie2018spectral}. Moreover, the essential image consists of those prestable categories which are separated and 0-complicial (\cite[C.5.3.1]{lurie2018spectral}). By \cite[C.3.1.4 and C.3.2.1]{lurie2018spectral}, the stabilization functor
	\[\tx{Sp}: \tx{Groth}_\infty^\tx{lex} \to \tx{Groth}_\infty^{+, \tx{lex}} \hspace{1em} \mc{C} \mapsto (\tx{Sp}(\mc{C}), \tx{Sp}(\mc{C})_{\ge 0})\]
	is a fully faithful embedding; here $\tx{Sp}(\mc{C})_{\ge 0}$ is the essential image of $\mc{C} \subseteq \tx{Sp}(\mc{C})$. By \cite[1.3.5.21]{HA}, the combination $\tx{Sp} \circ \mc{D}_{\ge 0}$ is the functor of derived category $\mc{D}(-)$, so we have a fully faithful embedding
	\[\mc{D}: \tx{Groth}_\tx{ab}^\tx{ex} \to \tx{Groth}_\infty^{+, \tx{lex}}.\]
	Let $\tx{Pr}^{\tx{t-ex}}$ denote the essential image. By \cite[C.3.1.5]{lurie2018spectral}, this category is the category of stable presentable categories $\mc{C}$ equipped with an accessible $t$-structure $(\mc{C}_{\ge 0}, \mc{C}_{\le 0})$ that is right complete and compatible with filtered colimits, and such that $\mc{C}_{\ge 0}$ is separated and 0-complicial. Morphisms in this category are the $t$-exact continuous functors. In other words, we obtain an equivalence
	\[\mc{D}: \tx{Groth}_\tx{ab}^\tx{ex} \simeq \tx{Pr}^{\tx{t-ex}}: \tx{Heart}\]

Both sides of $\mc{D}$ admit finite products, given by product of $\infty$-categories and $(\mc{C} \times \mc{C}')_{\ge 0} := \mc{C}_{\ge 0} \times \mc{C}'_{\ge 0}$. We claim $\mc{D}$ commutes with finite products. It suffices to check this on the underlying pair of $\infty$-categories, i.e. that for any $A, B \in \tx{Groth}_\tx{ab}^\tx{ex}$ we have
\[\mc{D}(A \times B) \simeq \mc{D}(A) \times \mc{D}(B) \hspace{1em} \mc{D}(A \times B)_{\ge 0} \simeq \mc{D}(A)_{\ge 0} \times \mc{D}(B)_{\ge 0};\]
which is immediate from the description of $\mc{D}(-)$ as a localization via the injective model structure.
Let $(\tx{Groth}_\tx{ab}^\tx{ex})^{\times}$ and $(\tx{Pr}^{\tx{t-ex}})^\times$ denote the Cartesian symmetric monoidal structures on them (defined in \cite[Chapter 2.4.1]{HA}), then $\mc{D}(-)$ is symmetric monoidal functor between them.

Let $\ms{CM}$ denote the $\infty$-operad encoding a symmetric monoidal algebra along with its left module. The data of an $\ms{CM} \tensor \mb{E}_2$-algebra in $\tx{Pr}^L$, whose underlying $\mb{E}_\infty \tensor \mb{E}_2 \simeq \mb{E}_\infty$-algebra is $\Vect$, is precisely an $\mb{E}_2$ DG category; likewise, a morphism of $\mb{E}_2$ DG categories is a map of such $\ms{CM} \tensor \mb{E}_2$-algebras that restricts to identity on the algebra part.
We have a pair of equivalence functors
\[\mc{D}^\tx{enh}: (\ms{CM} \tensor \mb{E}_2)\tx{-Alg}((\tx{Groth}_\tx{ab}^\tx{ex})^{\times}) \simeq (\ms{CM} \tensor \mb{E}_2)\tx{-Alg}((\tx{Pr}^{\tx{t-ex}})^\times): \tx{Heart}^\tx{enh}\]
both computed by the underlying functors. The symmetric monoidal forgetful functor $(\tx{Pr}^{\tx{t-ex}})^\times \to 1\tx{-Cat}^\times$ then gives a functor $\oblv_t: (\ms{CM} \tensor \mb{E}_2)\tx{-Alg}((\tx{Pr}^{\tx{t-ex}})^\times) \to (\ms{CM} \tensor \mb{E}_2)\tx{-Alg}(1\tx{-Cat}^\times)$. On the other hand, by \cite[4.8.1.9 and 4.8.1.15]{HA}, there is a \emph{fully faithful} embedding
\[\oblv_\tx{Pr}: (\ms{CM} \tensor \mb{E}_2)\tx{-Alg}(\tx{Pr}^L) \subseteq (\ms{CM} \tensor \mb{E}_2)\tx{-Alg}(1\tx{-Cat}^\times),\]
whose image consists of those algebras whose underlying categories are presentable, and whose action maps are colimit-preserving.

Let $(\Vect, \Rep_q(T))$ be the $(\ms{CM} \tensor \mb{E}_2)$-algebra in $\tx{Pr}^L$ observing the $\mb{E}_2$ DG category structure of $\Rep_q(T)$. On the other hand, let $(\Vect^\heartsuit, \Rep_q(T)^\heartsuit) \in (\ms{CM} \tensor \mb{E}_2)\tx{-Alg}((\tx{Groth}_\tx{ab}^\tx{ex})^{\times})$ denote the algebra encoding the $k$-linear braided monoidal category structure of $\Rep_q(T)^{\heartsuit}$. It follows from unfolding definition that we have
\[\oblv_t \circ \mc{D}^\tx{enh}((\Vect^\heartsuit, \Rep_q(T)^\heartsuit)) \simeq \oblv_{\tx{Pr}}((\Vect, \Rep_q(T))).\]

Recall the category $\mc{E}$ assumed in the statement of the lemma. It now suffices to lift it to an element $(\mc{E}, \mc{E}_{\ge 0})$ of $(\tx{Pr}^{\tx{t-ex}})^\times$ and provide\ an isomorphism of $\ms{CM} \tensor \mb{E}_2$-algebras $\tx{Heart}^\tx{enh}((\mc{E}, \mc{E}_{\ge 0})) \simeq (\Vect^\heartsuit, \Rep_q(T)^\heartsuit)$ that is identity on the algebra part. This is precisely the data of $\alpha$ and $\beta$.
\end{proof}

\subsection{Geometric Construction}
\begin{notation}
Because $\cT(\mb{C})$ is commutative, $B \cT(\mb{C})^{\tx{Spc}}$ is an infinite loop space. For notational simplicity, let us introduce
\[\doubleLoop := \Omega^2(B \cT(\mb{C})^\tx{Spc}) \in \mb{E}_\infty\tx{-Alg}(\tx{Spc});\]
this notation is justified because there is indeed an identification of both sides as $\mb{E}_\infty$-spaces, when the lattice is understood as a discrete topological group. We make this distinction so that $\weightLat$ by itself always denotes the lattice.
\end{notation}

Start with a fixed\footnote{The ambiguity in $\pi_2$ does not pose a problem, since the homotopy coherence can be checked by hand.} factorizable Betti gerbe $\TheGerbe$ on $\tx{Gr}_{\cT}(\mb{A}^1)_\Ran$; more precisely, this is an element of $\Theta(\weightLat, \tb{Ge}_\tx{an})(\mb{A}^1)$, the strict 2-groupoid classfying factorizable Betti gerbes on $\tx{Gr}_{\check{T}}$ over $\mb{A}^1$ (\cite[5.3.5]{zhao2020tame}).

The $\mb{C}^*$-version of \cite[3.1.9]{gaitsgory2018parameters} is that\footnote{This is not a formal consequence of the torsion coefficient version, but follows the same proof as in \emph{loc.cit.}.}, for any smooth curve $X$, the data of factorizable Betti gerbes on $\tx{Gr}_{\check{T}}(X)_\Ran$ is given by
\[\tx{Maps}_{\tx{Spc}}(B \check{T}(\mb{C}) \times X(\mb{C}), \tb{K}(\mb{C}^\times, 4));\]
In particular, for $\mb{A}^1$ we thus obtain a map
\[\Theta(\weightLat, \tb{Ge}_\tx{an})(\mb{A}^1) \to \mc{P}_{\tx{torus}}.\]

Delooping the image of $\TheGerbe$ twice, we obtain a map of $\mb{E}_2$-spaces
\[g: \doubleLoop \to \tb{K}(\mb{C}^\times, 2)\]

\begin{lemma}
\label{lemma:LS-tw-univ}
On the classifying space $\tb{K}(\mb{C}^\times, 2) \in \mb{E}_\infty\tx{-Alg}(\tx{Spc})$ there exists a universal locally constant cosheaf of categories $\ms{LS}^\tx{tw}_{\tx{univ}}$ corresponding to a map of $\mb{E}_\infty$-spaces $u: \tb{K}(\mb{C}^\times, 2) \to \DGCat^\tx{Spc}$.
\end{lemma}
\begin{proof}
Let $\DGCat_1$ be the full $\infty$-subcategory of $\DGCat$ consisting of the unit $\Vect$; it is a $\mb{E}_\infty$-subalgebra of $\DGCat$ within $\VVLCat$. The colocalization $\VVLCat \to \VVLSpc$ then sends it to an element $\DGCat_1^\tx{Spc} \in \mb{E}_\infty\tx{-Alg}(\VVLSpc)$. It is evidently grouplike \cite[5.2.6.16]{HA}, so to achieve our purpose it suffices to produce an $\mb{E}_\infty$ map $\mb{C}^\times \to \Omega^2(\DGCat_1^\tx{Spc})$. The grouplike $\mb{E}_\infty$ space $\Omega(\DGCat_1^\tx{Spc})$ is the full subgroupoid of $\tx{Maps}_{\DGCat}(\Vect, \Vect) \simeq \Vect^{\tx{Spc}} \in \mb{E}_\infty\tx{-Alg}(\tx{Spc})$ consisting of invertible elements. Its loop space at base point $k \in \Vect$, then, is the full subgroupoid of $\tx{Maps}_\Vect(k, k) \simeq \tx{Dold-Kan}(k)$ consisting of invertible elements, which is the (discrete!) group of invertible elements of $k = \mb{C}$.
\end{proof}

Since $\tx{Spc} \subseteq 1\tx{-Cat}$ and $\DGCat^\tx{Spc} \subseteq \DGCat$ are both symmetric monoidal, we can consider $u$ as a map between $\mb{E}_\infty$-categories.
Thus $u \circ g$ is a (strict) $\mb{E}_2$ functor $\doubleLoop \to \DGCat$. Let $\pi: \doubleLoop \to \tx{pt}$ be the projection. Because left Kan extension along $\pi$ is \emph{symmetric monoidal} with respect to the Day convolution structure on functor spaces, the result $\pi_{\dagger} g^\dagger(\ms{LS}_\tx{univ}^\tx{tw}) \simeq \colimit_{\doubleLoop} u \circ g$ is a lax $\mb{E}_2$-functor $\tx{pt} \to \DGCat$, i.e. a $\mb{E}_2$ DG category.

As a map between $\mb{E}_1$-spaces, there exists a homotopy $\alpha$ between $g$ and the trivial map, because
\[\pi_i\tx{Maps}_{\mb{E}_1\tx{-Alg}(\tx{Spc})}(\doubleLoop, \tb{K}(\mb{C}^\times, 2)) = \check{H}^{3 - i}(\check{T}(\mb{C}), \mb{C}^\times) = \begin{cases}
\Hom(\weightLat, \mb{C}^\times) & i = 2 \\
0 & \tx{else}
\end{cases}\]
so (up to ambiguity in $\pi_2$) we obtain a canonical identification $\alpha: \pi_{\dagger} g^\dagger(\ms{LS}_\tx{univ}^\tx{tw}) \simeq \Rep(T)$ as $\mb{E}_1$ DG categories.
By discussion above, $\pi_{\dagger} g^\dagger(\ms{LS}_\tx{univ}^\tx{tw})$ now inherits a $t$-structure such that its heart is a braided monoidal category, and we only need to identify the quadratic form associated with it.

Let $\Vect_{\check{\lambda}}$ denote the image of $\check{\lambda}$ under $g \circ u$. If we fix any $e \in E := (\mb{E}_2^{\tensor})_{\langle 2 \rangle}$ and let $e_!$ denote the coCartesian pushforward, then the value of the quadratic form in question at $\check{\mu}$ is the effect of the 2-morphism induced by the generator of $\pi_1(E^\tx{Spc}) \simeq \mb{Z}$ (in accordance with Convention \ref{conv:C-is-R2})
\[\xymatrix{
e_!(\colimit_{\check{\lambda}_1} \Vect_{\check{\lambda}_1}, \colimit_{\check{\lambda}_2}\Vect_{\check{\lambda}_2}) \ar[r] \ar[d] \drrtwocell<\omit> & \colimit_{\check{\lambda}_1, \check{\lambda}_2} e_!(\Vect_{\check{\lambda}_1}, \Vect_{\check{\lambda}_2}) \ar[r] & \colimit_{\check{\lambda}} \Vect_{\check{\lambda}} \ar[d] \\
e_!(\colimit_{\check{\lambda}_1} \Vect_{\check{\lambda}_1}, \colimit_{\check{\lambda}_2}\Vect_{\check{\lambda}_2}) \ar[r] & \colimit_{\check{\lambda}_1, \check{\lambda}_2} e_!(\Vect_{\check{\lambda}_1}, \Vect_{\check{\lambda}_2}) \ar[r] & \colimit_{\check{\lambda}} \Vect_{\check{\lambda}}
}\]
at $(k_{\check{\mu}}, k_{\check{\mu}})$; unfolding definition, this is the effect of the 2-morphism
\[\xymatrix{
\Vect_{\check{\mu}} \tensor \Vect_{\check{\mu}} \ar[r] \ar[d] \drtwocell<\omit> & \Vect_{\check{\mu} + \check{\mu}} \ar[d] \\
\Vect_{\check{\mu}} \tensor \Vect_{\check{\mu}} \ar[r] & \Vect_{\check{\mu} + \check{\mu}}
}\]
which is by construction the value $q_g(\check{\mu})$, where $q_g$ denote the image of $g$ under
\[\pi_0(\tx{Maps}_{\mb{E}_2\tx{-Alg}(\tx{Spc})}(\doubleLoop, \tb{K}(\mb{C}^\times, 2)) \xrightarrow[\tx{delooping}]{\simeq} \pi_0(\mc{P}_{\tx{torus}}) \simeq \tx{Quad}(\weightLat),\]
i.e. the quadratic form inside $\TheGerbe$. By Proposition \ref{prop:RepqT-is-unique}, this gives us an identification of $\mb{E}_2$ DG categories
\[\pi_{\dagger} g^\dagger(\ms{LS}_\tx{univ}^\tx{tw}) \simeq \Rep_{q_g}(T).\]

\medskip

We now obtain a more explicit description of $\Psi(\Rep_{q_g}(T))$ as a factorizable cosheaf of spaces.

In general, let $\mc{C}$ be a non-unital $\mb{E}_2$-category, and $\mc{D}$ a \emph{symmetric monoidal} category. Suppose we are given a $\mb{E}_2^\tx{nu}$-functor $F: \mc{C} \to \mc{D}$, then by \cite[5.5.4.10]{HA} we obtain a morphism of $\tx{Fact}(\mb{R}^2)$-algebras in $1\tx{-Cat}$ $\Psi(F): \Psi(\mc{C}) \to \Psi(\mc{D})$, which is (\cite[2.4.2.4]{HA}) equivalent to a map $\Psi(F)^{\tensor}: \Psi(\mc{C})^{\tensor} \to \Psi(\mc{D})^{\tensor}$ of coCartesian fibrations over $\tx{Fact}(\mb{R}^2)^{\tensor}$. Now, $\mc{D}$ being symmetric monoidal means we have a well-defined map
\[\Psi(\mc{D})^{\tensor} \to \tx{Fact}(\mb{R}^2)^{\tensor} \times_{\CommOperad} \mc{D}^{\tensor}\]
of coCartesian fibrations over $\tx{Fact}(\mb{R}^2)^{\tensor}$, given by the \emph{unique} push along each $\Psi(\mc{D})(U) \to \Psi(\mc{D})(\Ran(\mb{R}^2)) \simeq \mc{D}$. Precomposing with $\Psi(F)$ we obtain a map $\Psi(F)_{\Gamma}: \Psi(\mc{C})^{\tensor} \to \mc{D}^{\tensor}$ of coCartesian fibrations over $\CommOperad$ that preserve inert maps. Taking relative left Kan extension along $\Psi(\mc{C})^{\tensor} \to \tx{Fact}(\mb{R}^2)^{\tensor}$ with respect to $\CommOperad$, we arrive at a map
\[\tx{LKE}(\Psi(F)_{\Gamma}): \tx{Fact}(\mb{R}^2)^{\tensor} \to \mc{D}^{\tensor}\]
over $\CommOperad$. On the other hand $\tx{LKE}(F) := \colimit F \in \tx{Fun}(\tx{pt}, \mc{D}) \simeq \mc{D}$ is a $\mb{E}_2^\tx{nu}$-algebra within $\mc{D}$, then we can attach to it an object $\Psi_\mc{D}(\tx{LKE}(F))$.
Then we have $\tx{LKE}(\Psi(F)_{\Gamma}) \simeq \Psi_{\mc{D}}(\tx{LKE}(F))$: indeed it is easy to see that LHS is a factorizable cosheaf, thus it suffices to produce an equivalence $\tx{LKE}(\Psi(F)_{\Gamma})|_{\tx{Disk}(M)} \simeq \tx{LKE}(F)$, which now follows from definition.

Apply the above paradigm to $u \circ g: \doubleLoop \to \DGCat$, we obtain an identification
\[\Psi(\Rep_{q_g}(T)) \simeq \tx{LKE}(\Psi(u \circ g)_{\Gamma}).\]
On the other hand, from $g: \doubleLoop \to \tb{K}(\mb{C}^\times, 2)$ we can obtain a map over $\CommOperad$ given by
\[\Psi(\doubleLoop)^{\tensor} \xrightarrow{\Psi(g)_{\Gamma}} \tb{K}(\mb{C}^\times, 2)^{\tensor} \xrightarrow{u^{\tensor}} \DGCat^{\tensor}\]
Again taking relative LKE we obtain
\[\tx{LKE}(u^{\tensor} \circ \Psi(g)_\Gamma): \tx{Fact}(\mb{R}^2)^{\tensor} \to \DGCat^{\tensor};\]
Note that $\Psi(u \circ g)_\Gamma \simeq u^{\tensor} \circ \Psi(g)_\Gamma$, so we have
\[\Psi(\Rep_{q_g}(T)) \simeq \tx{LKE}(\Psi(u \circ g)_\Gamma) \simeq \tx{LKE}(u^{\tensor} \circ \Psi(g)_\Gamma).\]
The RHS can be understood as assignment to each $U \subseteq \Ran(\mb{R}^2)$ the DG category of $q$-twisted local systems on $\Psi(\doubleLoop)(U)$; the extra structure here is the factorization property of this assignment, which we shall discuss after a small digression.

\subsubsection{Comparison with Affine Grassmannian}
\newcommand{\GrTcShv}{\ms{Gr}\check{\ms{T}}}
The factorization space $\Psi(\doubleLoop)$ admits another description as follows.
Recall the non-unital factorization space $\tx{Gr}_{\cT}(\mb{A}^1)_\Ran$, whose $\mb{C}$-points give a map of topological spaces
\[\tx{Gr}_{\cT}(\mb{A}^1)_\Ran(\mb{C}) \to \Ran(\mb{C}) \simeq\Ran(\mb{R}^2).\]
Consider the copresheaf of spaces over $\Ran(\mb{R}^2)$ given by
\[\GrTcShv(U) := [U \times_{\Ran(\mb{R}^2)} \tx{Gr}_{\cT}(\mb{A}^1)_\Ran(\mb{C})]^\tx{Spc} \in \tx{Spc};\]
The natural inclusion of subspaces then upgrades it to a $\tx{Fact}(\mb{R}^2)$-algebra. We claim that it is a $\tx{Spc}$-valued factorizable cosheaf: indeed, before passing through $\tx{Top} \to \tx{Spc}$ this claim was verified in \cite[Section 4]{nocera2021model}; since the colimits in the cosheaf condition in this case are homotopy colimits (the relevant maps being monomorphisms of topological spaces), the cosheaf condition still holds as a $\tx{Spc}$-valued cosheaf.

It then follows that $\Psi^{-1}(\GrTcShv)$ is a well-defined $\mb{E}_2$-algebra in $\tx{Spc}$. Since we have an identification of weak homotopy types $\tx{Gr}_{\cT}(\mb{A}^1)(\mb{C})^\tx{Spc} \simeq \Omega^2(B \cT(\mb{C})^\tx{Spc})$, one is naturally led to believe that $\Psi^{-1}(\GrTcShv)$ carries the $\mb{E}_2$-structure from being a double loop space. However, this is not immediate from the factorization perspective.

\begin{remark}
For a general reductive $G$, the identification of homotopy types $\tx{Gr}_G(\mb{C})^\tx{Spc} \simeq \Omega^2(BG(\mb{C})^\tx{Spc})$ associates to the former an $\mb{E}_2$-structure. By looking at sections, one expects that the corresponding cosheaf of spaces $\Psi(\tx{Gr}_G(\mb{C})^\tx{Spc})$ is the $\mb{C}$-points of the Beilinson-Drinfeld Grassmannian $\tx{Gr}_G(\mb{A}^1)_\Ran$. We do not need, nor do we prove, this identification for general $G$.
\end{remark}

Luckily, we know a bit more about this space, namely that $\tx{Gr}_{\cT}(\mb{A}^1)_\Ran$ is a \emph{commutative} factorization space under the fiber product symmetric monoidal structure on non-unital factorization spaces. This structure induces an upgrade of $\GrTcShv$ as an element of $\mb{E}_\infty\tx{-Alg}(\tx{Fact}(\mb{R}^2)\tx{-Alg}(\tx{Spc}))$ with respect to the section-wise fiber product (recall that in $\tx{Spc}$ colimits distributes over fiber products). Because $\Psi$ is symmetric monoidal (this follows from unwinding construction), we end up producing an element
\[\Psi^{-1}(\GrTcShv) \in \mb{E}_\infty\tx{-Alg}(\mb{E}_2\tx{-Alg}(\tx{Spc})) \simeq \mb{E}_\infty\tx{-Alg}(\tx{Spc}).\]
We claim that we have an $\mb{E}_\infty$ equivalence $\Psi^{-1}(\GrTcShv) \simeq \doubleLoop$, where the latter is given the evident symmetric monoidal structure: indeed, by recognition principle \cite[5.2.6.26]{HA} it suffices to compute the connective spectrum corresponding to $\Psi^{-1}(\GrTcShv)$. The underlying space has discrete homotopy type, and the abelian group structure is forced by the additional $\mb{E}_\infty$ structure mentioned above, so the claim follows from \cite[1.4.3.6]{HA}. It follows that
\begin{corollary}
There is an unique identification of $\tx{Spc}$-valued factorizable cosheaves $\Psi(\doubleLoop) \simeq \GrTcShv$.
\end{corollary}

\begin{remark}
\label{rmk:rho-shift-does-nothing-on-shvc}
In the main proof, it is not $\GrTcShv$, but $\ms{Gr}\check{\ms{T}}^{\omega^{\crho}}$, the $\mb{C}$-points of the $\crho$-shifted Grassmannian (c.f. Definition \ref{def:rho-shifted-gr}) that we want to compare with $\mb{E}_2$. However, we note that over $\mb{A}^1$ their analytifications are canonically isomorphic.
\end{remark}

\subsubsection{Ind-Lisse Sheaves}

We begin with a basic fact relating local systems as defined above and \emph{ind-lisse} constructible sheaves. Let $Z$ be an aft scheme over $\mb{C}$ and let $\TheGerbe_0 \in \tx{Maps}(Z, \tx{\tb{Ge}}_\tx{an})$ be an analytic gerbe. Let $g_0: Z(\mb{C})^\tx{Spc} \to \tb{K}(\mb{C}^\times, 2)$ be the map induced by $\TheGerbe_0$. Let us define
\begin{equation}
\label{eq:categorical-twisted-ls-def}
\ms{LS}_{g_0}(Z) := \ms{LS}_{g_0}(Z(\mb{C})^\tx{Spc}) := \limit_{Z(\mb{C})^{\tx{Spc}}}\left(Z(\mb{C})^{\tx{Spc}} \xrightarrow{g_0} \tb{K}(\mb{C}^\times, 2) \xrightarrow{u} \DGCat\right);
\end{equation}
The same argument in \cite[Section 1.4]{gaitsgory2020toy} still gives us the $\dagger$-pull and $\dagger$-push functors. On the other hand, let $\tx{Shv}^!_{\tx{ind-lisse}, \TheGerbe_0}(Z)$ denote the full subcategory of $\tx{Shv}^!_{c, \TheGerbe_0}(Z)$ consisting of objects whose value lands in the ind-lisse subcategory (i.e. filtered colimits of lisse objects) when pulled back to each $V \in \tx{Split}(\TheGerbe_0, Z)$. We will refer to these two categories as \emph{$\TheGerbe_0$-twisted local systems} and \emph{$\TheGerbe_0$-twisted ind-lisse sheaves} respectively.

The following is an immediate generalization of \cite[Lemma A.4.2]{gaitsgory2020toy}:
\begin{lemma}
	\label{lemma:twisted-LS-is-twisted-lisse}
	We have $\ms{LS}_{g_0}(Z) \simeq \tx{Shv}^!_{\tx{ind-lisse}, \TheGerbe_0}(Z)$.
\end{lemma}
\begin{proof}
	First note that \cite[Theorem 3.21]{blanc2015topological}, which says that the homotopy type of $X(\mb{C})$ can be computed via the Cech nerve of an h-cover, implies that $\ms{LS}_{g_0}$ satisfies h-descent. On the other hand, $\tx{Shv}^!_{\tx{ind-lisse}, \TheGerbe_0}(Z)$ satisfies \'{e}h-descent \emph{by definition}. Moreover, if (the pullback of) $\TheGerbe_0$ is trivializable on some cover $V$, then we have
	\[\ms{LS}_{g_0}(V) \simeq \left(\vcenter{\xymatrix{\bullet \ar[r] \ar[d] & \tx{pt} \ar^{\tx{triv}}[d] \\
\tx{pt} \ar^(0.3){g_0}[r] & \tx{Maps}(V, \tb{K}(\mb{C}^\times, 2))}}\right) \times^{\tx{Maps}(V, \tb{K}(\mb{C}^\times, 1))} \ms{LS}(V) \simeq \tx{Shv}^!_{\tx{ind-lisse}, \mathscr{G}_0}(V);\]
The claim now follows from descent.
\end{proof}

For future use, we remark that, if $\TheGerbe_0$ admits an \emph{\'{e}tale} trivialization, we can equip both sides of the above lemma with $t$-structures such that the equivalence is $t$-exact:

\begin{itemize}
\item First, being an ind-category, $\tx{Shv}^!_c(Z) := \tx{Ind}(\tx{Shv}^!_{\tx{constr}}(Z))$ admits a unique $t$-structure (\cite[C.2.4.3]{lurie2018spectral}) specified by $\tx{Shv}^!_c(Z)^{\le 0} := \tx{Ind}(\tx{Shv}^!_{\tx{constr}}(Z)^{\tx{perv}\le 0})$. It guarantees, by design, that $\tx{Shv}^!_{\tx{constr}}(Z) \subseteq \tx{Shv}^!_c(Z)$ is $t$-exact;
\item Because the lisse condition is defined on each cohomology sheaf, the lisse subcategory is closed under extensions. It follows (\cite[1.4.4.1]{HA}) that $\tx{Shv}^{!}_{\tx{ind-lisse}}(Z)^{\le 0} := \tx{Ind}(\tx{Shv}^!_{\tx{constr}}(Z)^{\tx{perv}\le 0} \cap \tx{Shv}^!_\tx{lisse}(Z))$ generates a $t$-structure for $\tx{Shv}^!_\tx{ind-lisse}(Z)$. Moreover, since $\tx{Shv}^!_\tx{lisse}(Z)$ is closed under truncation $\tau^{\ge 0}$, it follows that $\tx{Shv}^!_\tx{ind-lisse}(Z) \subseteq \tx{Shv}^!_c(Z)$ is $t$-exact;
\item By assumption, there exists an \'{e}tale covering $V / Z$ such that \[\tx{Shv}^!_{\tx{ind-lisse}, \TheGerbe_0}(Z) \simeq \tx{Tot}_{!\tx{-pull}} \tx{Shv}^!_{\tx{ind-lisse}, \TheGerbe_0}(V^\bullet/Z);\]
each entry on the RHS has a $t$-structure such that the connecting morphisms are $t$-exact, so it equips LHS with an unique $t$-structure such that $!$-pulling back to each $V^\bullet$ is $t$-exact;
\item On the other hand, $\ms{LS}(Z)$, being a limit of categories, has an unique $t$-structure such that each $\dagger$-fiber map is $t$-exact. We also have a descent isomorphism $\ms{LS}_{g_0}(Z) \simeq \tx{Tot}_{\dagger\tx{-pull}} \ms{LS}_{g_0}(V^\bullet/Z)$, in which each evaluation is $t$-exact, so we can equip it with a $t$-structure such that the previous lemma is $t$-exact.
\end{itemize}

\subsection{Comparison with Constructible Sheaves}
We now return to the question of understanding $\Psi(\Rep_{q_g}(T))$ explicitly.

\begin{notation}
	During this subsection we will set $\ms{GR}_{p} := \tx{Gr}_{\cT}(\mb{A}^1)_{U(p)}$ for simplicity of notation.
\end{notation}

From recollement it is easy to see that, for any analytic factorizable gerbe $\TheGerbe_0$, the system of categories
\[p \in \ECOp^\op \mapsto \tx{Shv}^!_{\tx{Ran-c}, \TheGerbe_0}(\ms{GR}_p)\]
defined in Section \ref{sect:rh} form a weak constructible factorization category, where cosheaves on the base act via the previous lemma. Recall that $\TheGerbe$ is our fixed factorizable Betti gerbe, and $g$ is the associated map of $\mb{E}_2$-spaces. Our last step of comparison is:

\begin{proposition}
\label{prop:twisted-shv-is-RepqT}
We have an identification of weak constructible factorization categories
\[\tx{Fact}(\tx{LKE}(u^{\tensor} \circ \Psi(g)_\Gamma)) \simeq \tx{Shv}^!_{\tx{Ran-c}, \TheGerbe}(\ms{GR}_\bullet);\]
here $\tx{Fact}$ is the functor constructed in Proposition \ref{prop:E2-to-fact-cat}.
\end{proposition}

\begin{proof}
For each $p \in \ECOp^\op$, let $A(p)$ denote the stratification poset for $\BettiU(p)$.
We slightly abuse notation and let $g$ denote the map $\ms{GR}_p^\tx{Spc} \to \tx{Gr}_{\check{T}}(\mb{A}^1)^\tx{an}_{\Ran} \xrightarrow{\Psi(g)_\gamma} \tb{K}(\mb{C}^\times, 2)$, and let $\tx{Un}(\varphi_p)$ denote the coCartesian fibration over $\tx{Entr}_{A(p)}(\ms{GR}_p)$ corresponding to the map
\[\varphi_p := \tx{Entr}_{A(p)}(\ms{GR}_p) \to \ms{GR}_p^{\tx{Spc}, \op} \simeq \ms{GR}_p^{\tx{Spc}} \xrightarrow{u \circ g} \DGCat.\]
Define 
\[\tx{cShv}_{A(p), g}(\ms{GR}_p) := \tx{Fun}_{/\tx{Entr}_{A(p)}(\ms{GR}_p)}(\tx{Entr}_{A(p)}(\ms{GR}_p), \tx{Un}(\varphi_p));\]
then the $p$-component of LHS is the lax limit of the system
$\{\ms{L}_a := \ms{LS}_{g}(\ms{GR}_a)\}_{a \in A(p)}$ where the transition functor for $a > b$ is, by the same argument as Lemma \ref{lemma:transition-theta-nearby-cycle}, given by the nearby cycle functors
\[\ms{LS}_g(\ms{GR}_a) \xrightarrow{(j_b)_\dagger} \tx{cShv}_{A(p), g}(\ms{GR}_p) \xrightarrow{(j_b)^\dagger} \ms{LS}_g(\ms{GR}_b)\]
for appropriately defined $\dagger$-pull and $\dagger$-push of twisted constructible sheaves. By recollement, we see that the $p$-component of LHS is isomorphic to $\tx{cShv}_{A(p), g}(\ms{GR}_p)$.
On the other hand, the $p$-component of the RHS is the category $\tx{Shv}^!_{c, A(p), \TheGerbe}(\ms{GR}_p)$ of ind-$A(p)$-constructible sheaves.

In general, assume $Z$ is a prestack such that $Z(\mb{C})$ is stratified by some poset $A$. Covariant Verdier duality gives a functor $D^c: \tx{Shv}^!_{\tx{constr}, A}(Z) \to \tx{cShv}_{A}(Z(\mb{C})^\tx{Spc})$ which lands in the compact subcategory of the latter (consisting of those sheaves that are fiber-wise finite-dimensional) and is fully-faithful. Ind-extending, we obtain a fully faithful functor $D: \tx{Shv}^!_{c, A}(Z) \to \tx{cShv}_A(Z(\mb{C})^{\tx{Spc}})$. Furthermore, this equivalence is compatible with pullbacks and pushforwards from individual strata; it follows that $D$ is also essentially surjective, since pushforwards from individual strata generate the entire category.

Let us fix a \emph{Zariski}\footnote{After proving a smooth base change theorem for constructible cosheaves, this can be relaxed to ``\'{e}tale'' by realizing twisted constructible cosheaves as lax limits of twisted local systems, and use \'{e}tale descent of the latter. However, proper-surjective descent may not be satisfied by twisted constructible cosheaves.} cover $U^\bullet / \ms{GR}_p$ over which $\TheGerbe$ is trivial; this is always achievable for $\mb{A}^1$. By definition, $\tx{Shv}^!_{c, A(p), \mathscr{G}}(\ms{GR}_p)$ satisfies descent, i.e. we have
\[\tx{Shv}^!_{c, A(p), \mathscr{G}}(\ms{GR}_p) \simeq \tx{Tot}\left(\tx{Shv}^!_{c, A(p), \mathscr{G}}(U^\bullet / \ms{GR}_p)\right) \simeq \tx{Tot}_{\mathscr{G}}\left(\tx{Shv}^!_{c, A(p)}(U^\bullet / \ms{GR}_p)\right)\]
where $\mathscr{G}$ in the RHS denotes the twist on the cosimplicial system. On the other hand, we claim we also have
\[\tx{cShv}_{A(p), g}(\ms{GR}_p) \simeq \tx{Tot}\left(\tx{cShv}_{A(p), g}(U^\bullet / \ms{GR}_p)\right) \simeq \tx{Tot}_g\left(\tx{cShv}_{A(p)}(U^\bullet / \ms{GR}_p)\right)\]
where the $g$ subscript in RHS denotes a twist defined using the action\footnote{Via the universal map $\tb{K}(\mb{C}^\times, 1) \to \Vect$, any element of $\tx{Maps}(\tx{Entr}(V), \tb{K}(\mb{C}^\times, 1))$ yields a constructible cosheaf on $V$.} of $\tx{Maps}(\tx{Entr}(V), \tb{K}(\mb{C}^\times, 1))$ on $\tx{cShv}_{A(p)}(V)$, similar to the case in Lemma \ref{lemma:twisted-LS-is-twisted-lisse}. Indeed, noting that $\tx{cShv}_{A(p), g}(\ms{GR}_p)$ consists of lifts of $\varphi_p$ to a map $\tx{Entr}_{A(p)}(\ms{GR}_p) \to \DGCatMarked$, the claim follows from \cite[A.7.1]{HA} which gives $\tx{Entr}_{A(p)}(\ms{GR}_p) \simeq |\tx{Entr}_{A(p)}(U^\bullet)|$.

We have thus constructed an equivalence $\tx{Shv}^!_{c, A(p), \mathscr{G}}(\ms{GR}_p) \simeq \tx{cShv}_{A(p), g}(\ms{GR}_p)$, which is also compatible with $!/\dagger$-pullbacks to individual strata. The compatibility with actions of constructible (co)sheaves on the base and the factorization structure can be established routinely.
\end{proof}

Combined with discussion from the previous section, we finally arrive at:
\begin{corollary}
\label{cor:RepqT-is-GrcT}
	We have an identification of weak constructible factorization categories
	\[\tx{Fact}(\Psi(\Rep_{q_g}(T))) \simeq \tx{Shv}^!_{\tx{Ran-c}, \TheGerbe}(\tx{Gr}_{\cT}(\mb{A}^1)_{\Ran}),\]
	where the RHS is defined in \ref{def:shv-ranc-weak-constr-factcat}.
\end{corollary}

\begin{remark}
We emphasize that the parameter $q_g$ depends only on the bilinear form because we work over $\mb{A}^1$; for general curve, the line bundle term of the quantum parameter would result in a $\mb{E}_{X(\mb{C})}$-category that is not induced from an $\mb{E}_2$-category.
\end{remark}

\subsection{Koszul/Verdier Duality}
Recall the notion of $\mb{E}_2$-Koszul duality described in Section \ref{sect:E2-koszul-duality}. The geometric incarnation of Koszul duality is that of \emph{Verdier duality on the Ran space}. For technical simplicity, we will restrict ourselves to consider Verdier duality on the configuration space instead; algebraically, this corresponds to considering negatively-weighted non-unital $\mb{E}_2$-algebras.

\begin{definition}
	Let $\mathscr{G}$ be an analytic gerbe on $[\tx{Gr}_{\check{T}}]_\Ran$ and let $\mathscr{G}^{-1}$ be the inverse gerbe (i.e. the gerbe corresponding to the inverse bilinear form). For $\mc{F} \in \tx{Shv}^!_{c, \mathscr{G}}(\tx{Conf})$, the \emph{Verdier dual} of $\mc{F}$, to be denoted $\mb{D}^\tx{Verdier}(\mc{F})$, is the object in $\tx{Shv}^!_{c, \mathscr{G}^{-1}}(\tx{Conf})$ representing the functor
	\[\mc{G} \in \tx{Shv}^!_{c, \mathscr{G}^{-1}}(\tx{Conf}) \mapsto \tx{Maps}_{\tx{Shv}^!_{c, \mathscr{G} \times \mathscr{G}^{-1}}(\tx{Conf} \times_{\Ran} \tx{Conf})}(\mc{F} \boxtimes \mc{G}, \Delta_*(\omega_{\tx{Conf}})) \in \tx{Spc}\]
Here $\omega_\tx{Conf}$ is the dualizing sheaf (i.e. the $!$-pullback of $k \in \Vect$) on the configuration space.
\end{definition}

Suppose now that $\mathscr{G}$ is a \emph{factorizable gerbe} and $\mc{F}$ is a factorizable algebra on the configuration space; then it is immediate that the Verdier dual of $\mc{F}$ underlies another factorization algebra (with respect to the dual gerbe), and we will again denote it by $\mb{D}^\tx{Verdier}(\mc{F})$. On compact objects, we therefore have a well-defined functor
\[\mb{D}^\tx{Verdier}: \tx{FactAlg}(\tx{Shv}^!_{\tx{constr}, \mathscr{G}}(\tx{Conf})) \to \tx{FactAlg}(\tx{Shv}^!_{\tx{constr}, \mathscr{G}^{-1}}(\tx{Conf}))^\op.\]

The enhanced Verdier duality functor has the following universal property, which is easy to verify by unfolding definition: suppose $\mc{F}$ is a $\mathscr{G}$-twisted factorizable algebra on $\tx{Conf}(\mb{A}^1)$ whose underlying sheaf object is compact, and $\mc{G}$ is any $\mathscr{G}^{-1}$-twisted factorization algebra on $\tx{Conf}(\mb{A}^1)$. Then we have
	\[\tx{Maps}_{\tx{FactAlg}(\tx{Shv}^!_{c, \mathscr{G}^{-1}}(\tx{Conf}))}(\mc{G}, \mb{D}^\tx{Verdier}(\mc{F})) \simeq \tx{Maps}_{\tx{FactAlg}(\tx{Shv}^!_{c, \mathscr{G} \times \mathscr{G}^{-1}}(\tx{Conf} \times_{\Ran} \tx{Conf}))}(\mc{F} \boxtimes \mc{G}, \Delta_*(\omega_{\tx{Conf}}))\]
	here we note that both $\Delta_*(\omega_\tx{Conf})$ and $\mc{F} \boxtimes \mc{G}$ have evident upgrades to factorization algebras.

We can now state the following compatibility result, first claimed in \cite[3.4.5]{gaitsgory2021factorization}:

\begin{proposition}
\label{prop:fact-koszul-verdier-compatibility}
	Let $q$ be the quadratic form corresponding to a factorization gerbe $\mathscr{G}$ on $[\tx{Gr}_{\check{T}}]_\Ran$. Suppose $A$ is a non-unital $\mb{E}_2$-algebra in $\Rep_q(T)$ concentrated in negative weights which is weight-wise finite dimensional (i.e. $\tx{FactAlg}(A)$ lands in the compact subcategory $\tx{Shv}^!_{\tx{constr}, \mathscr{G}}(\tx{Conf})$). Then we have
	\[\tx{FactAlg}^\op \circ \tx{inv}^\op \circ \tx{KD}^\nonun_{\mb{E}_2}(A) \simeq \mb{D}^\tx{Verdier} \circ \tx{FactAlg}(A),\]
	where $\tx{inv}: \Rep_{q^{-1}}(T) \to \Rep_{q^{-1}}(T)$ is the automorphism induced by the Cartan involution of $T$ (i.e. one which sends each $\clambda$-weight space to $-\clambda$).
\end{proposition}

\begin{proof}
We start by giving alternative description of Koszul pairing for $\Rep_q(T)$.
Note that $\tx{Alg}(\Rep_q(T))$ admits a monoidal structure, and the functor
\[\tx{biact}: [\tx{Alg}(\Rep_q(T)) \times \tx{Alg}(\Rep_{q^{-1}}(T))] \times \tx{Alg}(\Rep_q(T)) \to \tx{Alg}(\Rep_q(T))\]
given by two-sided multiplication, observes the structure of $\mc{M} := \tx{Alg}(\Rep_q(T))$ as being left tensored over the monoidal category $\mc{A} := \tx{Alg}(\Rep_q(T)) \times \tx{Alg}(\Rep_{q^{-1}}(T))$. Fix $A \in \mb{E}_2\tx{-Alg}(\Rep_q(T)), B \in \mb{E}_2\tx{-Alg}(\Rep_{q^{-1}}(T))$.
Consider the category
\[\mc{C} := \{a \in \mc{A}, \varphi: \tx{biact}(a, k) \to k \in \mc{M}\};\]
as in \cite[4.7.1]{HA}, $\mc{C}$ has the structure of a monoidal category, and an $(A, B)$-bimodule structure on $k$ is precisely an element of $\tx{Alg}(\mc{C})$ lifting $(A, B) \in \tx{Alg}(\mc{A})$.

Let $\tx{prod}: \Rep_{q \times q^{-1}}(T \times T) \to \Rep_q(T)$ be the monoidal functor of tensor product of graded vector spaces; it induces a plain functor \[\tx{prod}: \tx{Alg}(\Rep_{q \times q^{-1}}(T \times T)) \to \tx{Alg}(\Rep_q(T));\]
 Define
\[\mc{A}' := \tx{Alg}(\Rep_{q \times q^{-1}}(T \times T)) \hspace{1em} \mc{C}' := \{a' \in \mc{A}', \varphi: \tx{prod}(a') \to k \in \mc{M}\}.\]
Note that we have a box product $\boxtimes: \mc{A} \to \mc{A}'$ which is monoidal; also note that, as a plain functor, $\tx{biact}(-, k): \mc{A} \to \mc{M}$ factors as
\[\mc{A} \xrightarrow{\boxtimes} \mc{A}' \xrightarrow{\tx{prod}} \mc{M};\]
It follows that an $(A, B)$-bimodule structure on $k$ is the same as an element of $\tx{Alg}(\mc{C}')$ lifting $A \boxtimes B \in \tx{Alg}(\mc{A}')$. Next note that the plain functor $\tx{prod}$ admits a right adjoint $\tx{prod}^R$, and thus we have an induced adjoint pair
\[\tx{prod}: \mc{A}' \adjoint \mc{M}: \tx{prod}^R;\]
As in \cite[4.7.1.41]{HA}, $\tx{prod}^R(k) \in \mc{A}'$ equipped with the obvious map $\tx{prod} \circ \tx{prod}^R(k) \to k$ is the final object in $\mc{C}'$, and thus admits an unique lifting to an element in $\tx{Alg}(\mc{C}')$; furthermore, an $(A, B)$-bimodule structure on $k$ is the same as a map
\[A \boxtimes B \to \tx{prod}^R(k) \in \tx{Alg}(\mc{A}') \simeq \mb{E}_2\tx{-Alg}(\Rep_{q \times q^{-1}}(T \times T)).\]
We shall explicitly identify $\tx{prod}^R(k)$. Consider the following \emph{group homomorphism}
\[\tx{inv-mult}: T \times T \xrightarrow{\tx{id} \times (-)^{-1}} T \times T \xrightarrow{\tx{mult}} T\]
We observe that restriction along $\tx{inv-mult}$ gives a $\mb{E}_2$-map
\[\tx{res}_{\tx{inv-mult}}: \Rep(T) \to \Rep_{q \times q^{-1}}(T \times T);\]
indeed, this claim can be checked on the level of braided monoidal abelian categories. Consider the $\QCoh$-adjunction
\[(p_T)^*: \Vect \adjoint \Rep(T): (p_T)_*\]
the left adjoint is a fortiori $\mb{E}_2$, thus the right adjoint $(p_T)_*$ lifts to a lax $\mb{E}_2$-map, so that $\mc{O}_T := (p_T)_*(k)$ obtains an $\mb{E}_2$-algebra structure, and we set
\[\widetilde{\mc{O}_T} := \tx{res}_{\tx{inv-mult}}(\mc{O}_T) \in \mb{E}_2\tx{-Alg}(\Rep_{q \times q^{-1}}(T \times T)).\]
Because the object lies in the heart, it is easy to see that $\tx{prod}^R(k) \simeq \widetilde{\mc{O}_T}$ as $\mb{E}_2$-algebras.

If $A$ is an augmented algebra in $\Rep_q(T)$, then we have a naturally defined map of $\mb{E}_2$-algebras
\[\epsilon_A := A \boxtimes k \to k \boxtimes k \simeq k \to (p_T)_* (p_T)^* (k) \simeq \widetilde{\mc{O}_T}.\]
The discussion above can be summarized as
\begin{corollary}
	For every augmented $\mb{E}_2$-algebra $A$, its $\mb{E}_2$ Koszul dual is the universal augmented $\mb{E}_2$-algebra in $\Rep_{q^{-1}}(T)$ representing the functor
	\[B \in \mb{E}_2\tx{-Alg}_{\tx{aug}}(\Rep_{q^{-1}}(T)) \mapsto \tx{Maps}_{\mb{E}_2\tx{-Alg}(\Rep_{q \times q^{-1}}(T \times T))}(A \boxtimes B, \widetilde{\mc{O}_T})\]
	\[\times_{\tx{Maps}_{\mb{E}_2\tx{-Alg}(\Rep_{q \times q^{-1}}(T \times T))}(A \boxtimes k, \widetilde{\mc{O}_T}) \times \tx{Maps}_{\mb{E}_2\tx{-Alg}(\Rep_{q \times q^{-1}}(T \times T))}(k \boxtimes B, \widetilde{\mc{O}_T})} \{\epsilon_A, \epsilon_{B}\}.\]
\end{corollary}

We observe that when we restricts to negatively-weighted non-unital $\mb{E}_2$-algebras, the pairing admits a simpler description:
\begin{lemma}
	For every non-unital negatively-weighted $\mb{E}_2$-algebra $A \in \Rep_q(T)$, its $\mb{E}_2$-Koszul dual is the universal positively-weighted non-unital $\mb{E}_2$-algebra in $\Rep_{q^{-1}}(T)$ representing the functor
	\[B \in \mb{E}_2\tx{-Alg}_\nonun(\Rep_{q^{-1}}(T))^{>0} \mapsto \tx{Maps}_{\mb{E}_2\tx{-Alg}_\nonun(\Rep_{q \times q^{-1}}(T \times T))}(A \boxtimes \tx{inv}(B), \widetilde{\mc{O}_T}_\Delta^{<0})\]
	where
	\begin{itemize}
	\item $\mb{E}_2\tx{-Alg}_\nonun(\Rep_{q^{-1}}(T))^{>0}$ denotes the full subcategory of $\mb{E}_2$-algebras whose underlying object lies in the subcategory of positively-weighted objects (same for $(-)^{>0}$);
	\item $\tx{inv}(B)$ is the $\mb{E}_2$-functor of inverting the weight lattice;
	\item Let $\widetilde{\mc{O}_T}_\Delta$ be the $\mb{E}_2$-algebra whose underlying object is the restriction of $\mc{O}_T$ along $\tx{mult}$, then $\widetilde{\mc{O}_T}_{\Delta}^{<0}$ is the image of $\widetilde{\mc{O}_T}_\Delta$ under the lax (non-unital) $\mb{E}_2$ functor
	\[\Rep_{q \times q^{-1}}(T \times T) \to \Rep_{q \times q^{-1}}(T \times T)^{<0, <0} \to \Rep_{q \times q^{-1}}(T \times T)\]
	where the first (lax $\mb{E}_2$) map is projection and second ($\mb{E}_2$) map is inclusion.
	\end{itemize}
\end{lemma}

\begin{proof}
	The explicit description of the $\mb{E}_2$-Koszul dual via double Bar construction guarantees that it is positively-weighted, so by Yoneda, it suffices to specify how other positively-weighted non-unital $\mb{E}_2$-algebras map to it. The map $\tx{id} \times \tx{inv}$ is a $\mb{E}_2$-functor on $\Rep_{q \times q^{-1}}(T \times T)$, so the space in the statement is the same as
	\[\mc{S}_\tx{red} := \tx{Maps}_{\mb{E}_2\tx{-Alg}_\nonun(\Rep_{q \times q^{-1}}(T \times T))}(A \boxtimes B, \widetilde{\mc{O}_T}^{<0, >0})\]
	where $\widetilde{\mc{O}_T}^{<0, >0}$ is the non-unital $\mb{E}_2$-algebra obtained from $\widetilde{\mc{O}_T}$ as $\tx{inc}_{<0, >0} \circ \tx{proj}_{<0, >0}$, where
	\[\tx{inc}_{<0, >0}: \mb{E}_2\tx{-Alg}_{\nonun}(\Rep_{q \times q^{-1}}(T \times T)^{<0, >0}) \adjoint \mb{E}_2\tx{-Alg}_{\nonun}(\Rep_{q \times q^{-1}}(T \times T)) :\tx{proj}_{<0, >0}\]
	is the adjoint pair on non-unital $\mb{E}_2$-algebras induced by the inclusion-projection adjoint pair. Set
	\[\mc{S} := \tx{Maps}_{\mb{E}_2\tx{-Alg}(\Rep_{q \times q^{-1}}(T \times T))}((A \oplus k) \boxtimes (B \oplus k), \widetilde{\mc{O}_T})\]
	\[\times_{\tx{Maps}_{\mb{E}_2\tx{-Alg}(\Rep_{q \times q^{-1}}(T \times T))}((A \oplus k) \boxtimes k, \widetilde{\mc{O}_T}) \times \tx{Maps}_{\mb{E}_2\tx{-Alg}(\Rep_{q \times q^{-1}}(T \times T))}(k \boxtimes (B \oplus k), \widetilde{\mc{O}_T})} \{\epsilon_{A \oplus k}, \epsilon_{B \oplus k}\};\]
	Since $\Rep_q(T)$ is stable, $(A \oplus k) \boxtimes (B \oplus k) \simeq A \boxtimes B \oplus A \boxtimes k \oplus k \boxtimes B \oplus k$ is the product of its summands (as $\mb{E}_2$-algebras), so we see that $\mc{S}$ is isomorphic to $\tx{Maps}_{\mb{E}_2\tx{-Alg}(\Rep_{q \times q^{-1}}(T \times T))}(A \boxtimes B \oplus k, \widetilde{\mc{O}_T})$, which is $\tx{Maps}_{\mb{E}_2\tx{-Alg}_\nonun(\Rep_{q \times q^{-1}}(T \times T))}(A \boxtimes B, \widetilde{\mc{O}_T})$ since $(-) \oplus k$ is the left adjoint of the functor of forgetting the unital structure. Thus the adjunction $(\tx{inc}_{<0, >0}, \tx{proj}_{<0, >0})$ gives an isomorphism $\mc{S}_{\tx{red}} \simeq \mc{S}$.
\end{proof}
Under the correspondence between factorizable sheaves and non-unital $\mb{E}_2$-algebras, box product matches with box product, and $\Delta_*(\omega_{\tx{Conf}})$ matches with $\widetilde{\mc{O}_T}_\Delta^{<0}$ (it is easy to see the underlying object matches; to upgrade to $\mb{E}_2$-structure, use the fact that the latter lies in the heart). The original claim follows.
\end{proof}

\section{Main Computations}

\label{sect:main-computation}

\subsection{Proof of Theorem \ref{thm:jacquet-is-equivalence}}
We proceed in steps.

\paragraph{Preparation}
For each $\clambda \in \weightLat$ let $\iota_{\clambda \cdot x_0}$ denote the inclusion of the point $(\clambda \cdot x_0)$ into the marked configuration space. The costandard object $\CoStObj_{\tx{Conf}, \kappa}^{\cmu}$ is uniquely characterized by its right orthogonal property with the standard objects $\StObj_{\tx{Conf}, \kappa}^{\clambda}$; since the latter corepresent $!$-fiber at $\clambda \cdot x_0$, we see that $\CoStObj_{\tx{Conf}, \kappa}^{\cmu}$ is the unique factorization module $M$ that satisfies
\[\iota_{\clambda \cdot x_0}^{!}(M) \simeq \begin{cases} \mb{C} & \clambda = \cmu \\ 0 & \tx{else}\end{cases};\]
However, by construction, $\tb{J}_{\kappa}(\CoStObj_{\tx{KM}, \kappa}^{\clambda})$ exactly satisfies this property.

Now recall that the standard objects form a set of compact generators of the source category (Lemma \ref{lemma:wakimoto-generates}). On the other hand, $\StObj_{\tx{Conf}, \kappa}^{\clambda}$ also compactly generates. It therefore suffices to show that $\tb{J}^\tx{KM}_\kappa$ sends standard to standard, and is fully faithful on these standard objects.
The latter fact is, in fact, obvious, as we have
\[\Hom_{\kmgmod^I}(\StObj_{\tx{KM}, \kappa}^{\clambda}, \StObj_{\tx{KM}, \kappa}^{\cmu}) \simeq C^{\semiinf}(\Loop \mf{n}, \Arc \mf{n}, \StObj_{\tx{KM}, \kappa}^{\cmu})^{\clambda}\]
\[\simeq \Hom_{\TheFactAlg\tx{-FactMod}_\nonun(\DMod_{\DualTwisting}(\tx{Conf}_{\infty \cdot x_0}))}(\StObj_{\tx{Conf}, \kappa}^{\clambda}, \StObj_{\tx{Conf}, \kappa}^{\cmu}).\]

Because of the right orthogonality characterization of the standard objects, it suffices to prove that
\[\iota_{\mu \cdot x_0}^{*}(\tb{J}_\kappa(\StObj_{\tx{KM}, \kappa}^{\clambda})) \simeq \begin{cases} \mb{C} & \clambda = \cmu \\ 0 & \tx{else}\end{cases}\]

\paragraph{Adding One Point}
To do this we invoke the local-global compatibility by switching to the complete curve $\overline{X} \simeq \mb{P}^1$. We let $x_\infty \in \mb{P}^1(\mb{C})$ denote the added point, and set $\modulePoints := (x_0, x_\infty)$. The definition of Jacquet functor proceeds in the exact same way to give
\[\tb{J}_{\kappa, \modulePoints}: (\kmgmod^I)_{x_0} \tensor (\kmgmod^I)_{x_\infty} \to \TheFactAlg\tx{-FactMod}_\nonun(\DMod_{\DualTwisting}(\tx{Conf}_{\infty \cdot \modulePoints}(\mb{P}^1))).\]
Let
\[\tx{AJ}: \tx{Conf}_{\infty \cdot \modulePoints}(\mb{P}^1) \subseteq [\tx{Gr}^{\omega^{\crho}}_{\check{T}}(\mb{P}^1)]_{\Ran_{\modulePoints, \dR}} \to \tx{Bun}_{\check{T}}(\mb{P}^1)\]
denote the Abel-Jacobi map, given by\footnote{Note that $\tx{Pic}(\mb{P}^1) \tensor_{\mb{Z}} \weightLat \simeq \tx{Bun}_{\check{T}}(\mb{P}^1)$ requires choosing a base point; here we use $x_\infty$.}
\[\sum_{i} \clambda_i \cdot x_i \mapsto \sum_i \mc{O}(-x_i) \tensor \clambda_i \in \tx{Pic}(\mb{P}^1) \tensor_{\mb{Z}} \weightLat \simeq \tx{Bun}_{\check{T}}(\mb{P}^1).\]
The gerbe $\DualTwisting$ is, by design, canonically trivialized at $\tx{Bun}_{\check{T}}(\mb{P}^1)_{\omega^{\crho}} \simeq \mb{B} \check{T}$, the component of $\tx{Bun}_{\check{T}}(\mb{P}^1)$ corresponding to the $\check{T}$-bundle $\omega^{\crho}$; noting that the component $\tx{Conf}_{\infty \cdot \modulePoints}(\mb{P}^1)^{=2 \crho}$ gets sent to this component via the Abel-Jacobi map, we see that the gerbe $\DualTwisting$ admits a canonical trivialization on that component as well.
We let
\[\tx{Triv}_{\DualTwisting}^{\tx{Bun}}: \DMod_{\DualTwisting}(\tx{Bun}_{\check{T}}(\mb{P}^1)_{\omega^{\crho}}) \to \DMod(\tx{Bun}_{\check{T}}(\mb{P}^1)_{\omega^{\crho}})\]
\[\tx{Triv}_{\DualTwisting}: \DMod_{\DualTwisting}(\tx{Conf}^{=2 \crho}(\mb{P}^1)) \to \DMod(\tx{Conf}^{=2 \crho}(\mb{P}^1))\]
denote the resulting isomorphisms.

For any $\tx{Conf}^{=\mu}(\mb{A}^1)$, consider the open embedding $j_{\mu \to 2 \crho}: \tx{Conf}_{\infty \cdot x_0}^{=\mu}(\mb{A}^1) \inj \tx{Conf}_{\infty \cdot \modulePoints}^{=2 \crho}(\mb{P}^1)$. The gerbe $\DualTwisting$ on the former is isomorphic to the pullback from the latter, so by above, we see that there is a canonical trivialization of this gerbe on any $\tx{Conf}^{=\mu}(\mb{A}^1)$; let
\[\tx{Triv}_{\DualTwisting}: \DMod_{\DualTwisting}(\tx{Conf}^{=\mu}(\mb{A}^1)) \to \DMod(\tx{Conf}^{=\mu}(\mb{A}^1))\]
denote the resulting isomorphism. It suffices, then, to compute the $*$-fiber on the RHS using this trivialization.

Note that the $\mb{G}_m$ action on the curve $\mb{A}^1$ induces an $\mb{G}_m$ action on $\tx{Conf}_{\infty \cdot x_0}^{=\mu}(\mb{A}^1)$ with the point $\mu \cdot 0$ being the attractor locus (we refer readers to \cite{drinfeld2016geometric} for the meaning of these terms), and the sheaf $\tx{Triv}_{\DualTwisting}(\tb{J}_\kappa(\StObj_{\tx{KM}, \kappa}^{\clambda}))$ is $\mb{G}_m$-monodromic with respect to this action. By contraction principle (\cite[4.1.5]{drinfeld2016geometric}), we have
\[i_{\mu \cdot x_0}^*(\tx{Triv}_{\DualTwisting}(\tb{J}_\kappa(\StObj_{\tx{KM}, \kappa}^{\clambda}))) \simeq \Gamma_\dR(\tx{Conf}^{=\mu}(\mb{A}^1), \tx{Triv}_{\DualTwisting}(\tb{J}_\kappa(\StObj_{\tx{KM}, \kappa}^{\clambda})))\]
where the latter means $*$-pushforward to point; the latter can also be rewritten as
\begin{equation}
\label{eq:main-eqline-1}
\Gamma_\dR(\tx{Conf}^{=\mu}(\mb{P}^1), \tx{Triv}_{\DualTwisting}((j_{\mu \to 2 \crho})_*(\tb{J}_\kappa(\StObj_{\tx{KM}, \kappa}^{\clambda}))))
\end{equation}
Since we already showed $\tb{J}_{\kappa}(\CoStObj_{\tx{KM}, \kappa}^{\clambda}) \simeq \CoStObj_{\tx{Conf}, \kappa}^{\clambda}$, by factorization it is easy to see that
\[(j_{\mu \to 2 \crho})_*(\tb{J}_\kappa(\StObj_{\tx{KM}, \kappa}^{\clambda}))) \simeq \oblv_{\tx{FactMod}}\left(\tb{J}_{\kappa, \modulePoints}(\StObj_{\tx{KM}, \kappa}^{\clambda} \boxtimes \CoStObj_{\tx{KM}, \kappa}^{2 \crho - \cmu})\right).\]
The trivialization above also gives a trivialization for the inverse gerbe $-\DualTwisting$ on $\tx{Bun}_{\check{T}}(\mb{P}^1)_{\omega^{\crho}}$, and we set
\[\tx{Delta}_{\omega^{\crho}, -\DualTwisting} := (\j_\omega^{\crho})_! (\tx{Triv}_{-\DualTwisting}^{\tx{Bun}})^{-1}(\omega_{\mb{B} \check{T}}) \in \DMod_{-\DualTwisting}(\tx{Bun}_{\check{T}}(\mb{P}^1))\]
where $\j_\omega^{\crho}: \tx{Bun}_{\check{T}}(\mb{P}^1)_{\omega^{\crho}} \subseteq \tx{Bun}_{\check{T}}(\mb{P}^1)$ is the inclusion of the component.
Following definition we can rewrite Equation \eqref{eq:main-eqline-1} as
\[\Gamma_\dR(\tx{Conf}^{=\mu}(\mb{P}^1), \tx{Triv}_{\DualTwisting}(\oblv_{\tx{FactMod}}\left(\tb{J}_{\kappa, \modulePoints}(\StObj_{\tx{KM}, \kappa}^{\clambda} \boxtimes \CoStObj_{\tx{KM}, \kappa}^{2 \crho - \cmu})\right)))\]
\[\simeq \left\langle \tx{AJ}_* \circ \oblv_{\tx{FactMod}}\left(\tb{J}_{\kappa, \modulePoints}(\StObj_{\tx{KM}, \kappa}^{\clambda} \boxtimes \CoStObj_{\tx{KM}, \kappa}^{2 \crho - \cmu})\right), \tx{Delta}_{\omega^{\crho}, -\DualTwisting} \right\rangle_{\DMod(\tx{Bun}_{\check{T}}(\mb{P}^1))}\]
here we are using the usual duality pairing for $\tx{Bun}_T$ given by de Rham cohomology of the $!$-tensor.

Let $\tx{FM}_{\tx{glob}}^{\crho}$ be the $\crho$-shifted Fourier-Mukai equivalence described in Section \ref{sect:tori-glob-compatibility}. It follows from definition that
\[(\tx{FM}_{\tx{glob}}^{\crho})^{-1}(\tx{Delta}_{\omega^{\crho}, -\DualTwisting}) \simeq \tx{Delta}_{\tx{triv}, -[\kappa, \AnomalyTerm]} \in \DMod_{-[\kappa, \AnomalyTerm]}(\tx{Bun}_{T}(\mb{P}^1));\]
where the RHS is defined as before, using the canonical trivialization at the component of the trivial bundle.
Using compatibility between $\crho$-shifted torus FLE and $\crho$-shifted Fourier-Mukai shown in Section \ref{sect:tori-glob-compatibility}, we can further rewrite the previous expression as
\begin{equation}
\label{eq:main-eqline-2}
\left\langle [\tx{Loc}]^{\tx{ch}}_{\modulePoints}\left(\tb{J}_{\kappa}^{\tx{KM}}(\StObj_{\tx{KM}, \kappa}^{\clambda} \boxtimes \CoStObj_{\tx{KM}, \kappa}^{2 \crho - \cmu})\right), \tx{Delta}_{\tx{triv}, -[\kappa, \AnomalyTerm]} \right\rangle_{\DMod(\tx{Bun}_{T}(\mb{P}^1))}
\end{equation}
Using the diagram Equation \ref{eq:local-global-diagram}, we can rewrite Equation \eqref{eq:main-eqline-2} as
\begin{equation}
\label{eq:main-eqline-3}
\left \langle \tx{CT}_{*}^{\tx{shifted}} \circ [\tx{Loc}]_{\modulePoints}(\StObj_{\tx{KM}, \kappa}^{\clambda} \boxtimes \CoStObj_{\tx{KM}, \kappa}^{2 \crho - \cmu}), \tx{Delta}_{\tx{triv}, -[\kappa, \AnomalyTerm]} \right\rangle_{\DMod(\tx{Bun}_T(\mb{P}^1))}
\end{equation}
We have dualities
\[\DMod_{[\kappa, \AnomalyTerm]}(\tx{Bun}_T(\mb{P}^1))^\vee \simeq \DMod_{-[\kappa , \AnomalyTerm]}(\tx{Bun}_T(\mb{P}^1))\]
\[\DMod_{(\kappa, 0)}(\tx{Bun}_G(\mb{P}^1)^{\tx{level}_B, \modulePoints})^\vee \simeq \DMod_{(-\kappa, 0)}(\tx{Bun}_G(\mb{P}^1)^{\tx{level}_B, \modulePoints})_{\tx{co}}\]
We remind readers that the latter, known as the \emph{pseudo-identity functor} of $\tx{Bun}_G$, is a highly non-trivial construction and relies heavily on the geometry of $\tx{Bun}_G$ (see \cite{drinfeld2015compact} for details). We define
\[\tx{Eis}_{\tx{co}, *}^\tx{shifted} := (\tx{CT}_*^\tx{shifted})^\vee\]
as the dual functor of $\tx{CT}_*^\tx{shifted}$, and rewrite Equation \eqref{eq:main-eqline-3} formally as
\[\left \langle [\tx{Loc}]_{\modulePoints}(\StObj_{\tx{KM}, \kappa}^{\clambda} \boxtimes \CoStObj_{\tx{KM}, \kappa}^{2 \crho - \cmu}), \tx{Eis}_{\tx{co}, *}^\tx{shifted}(\tx{Delta}_{\tx{triv}, -[\kappa, \AnomalyTerm]}) \right\rangle_{\DMod(\tx{Bun}_G(\mb{P}^1)^{\tx{level}_B, \modulePoints})}\]
Now, the same argument as in \cite[Proposition 1.5.5]{gaitsgory2017strange} tells us that $\tx{Eis}_{\tx{co}, *}^\tx{shifted}(\tx{Delta}_{\tx{triv}, -[\kappa, \AnomalyTerm]})$ can be explicit computed as (we use the notation from \emph{loc.cit.})
\[(j_U)_{\tx{co}, *} \circ r_* (\omega_{\tx{pt} / N})[-\dim(\mf{n})]\]
here
\[r: \tx{pt} / N \simeq \tx{Bun}_{N}(\mb{P}^1) \to (\tx{pt} / G) \times_{\tx{Bun}_G(\mb{P}^1)} \tx{Bun}_G(\mb{P}^1)^{\tx{level}_B, \modulePoints} \simeq (\tx{pt} / B)_{0} \times_{\tx{pt} / G} (\tx{pt} / B)_{\infty}\]
is the evident projection from $\tx{Bun}_B$, and
\[j_U: (\tx{pt} / G) \times_{\tx{Bun}_G(\mb{P}^1)} \tx{Bun}_G(\mb{P}^1)^{\tx{level}_B, \modulePoints} \subseteq \tx{Bun}_G(\mb{P}^1)^{\tx{level}_B, \modulePoints}\]
is the inclusion of the open substack corresponding to the trivial bundle. By definition, for any $\mc{F}$ we have
\[\langle \mc{F}, (j_U)_{\tx{co}, *} \circ r_* (\omega_{\tx{pt} / N}) \rangle_{\DMod(\tx{Bun}_G(\mb{P}^1)^{\tx{level}_B, \modulePoints})}\]
\[\simeq \langle j_U^*(\mc{F}), r_* (\omega_{\tx{pt} / N}) \rangle_{\DMod((\tx{pt} / B)_{0} \times_{\tx{pt} / G} (\tx{pt} / B)_{\infty})}\]
\[\simeq \langle r^! \circ j_U^*(\mc{F}), \omega_{\tx{pt} / N} \rangle_{\DMod(\tx{pt} / N)} \simeq \iota_{\tx{triv}}^!(\mc{F})\]
for $\iota_{\tx{triv}}: \tx{pt} \to \tx{Bun}_G(\mb{P}^1)^{\tx{level}_B, \modulePoints}$ the inclusion of the trivial bundle. Thus we arrive at the expression
\begin{equation}
\label{eq:main-eqline-4}
\iota_{\tx{triv}}^!([\tx{Loc}]_{\modulePoints}(\StObj_{\tx{KM}, \kappa}^{\clambda} \boxtimes \CoStObj_{\tx{KM}, \kappa}^{2 \crho - \cmu}))[-\dim(\mf{n})].
\end{equation}
This fiber is provided by the (cohomologically shifted) chiral homology (over $\mb{P}^1$) of $(\StObj_{\tx{KM}, \kappa}^{\clambda},  \CoStObj_{\tx{KM}, \kappa}^{2 \crho - \cmu})$, considered as two modules over the Kac-Moody chiral algebra $\mb{V}_\kappa^0$ supported at $0$ and $\infty$. Note that this is not a tautology: instead it follows from applying Lemma \ref{lem-loc-compatible-oblv}.
As is classically known, this is the coinvariant
\[C_*(\mf{g}[t, t^{-1}], (\StObj_{\tx{KM}, \kappa}^{\clambda})_{0} \tensor (\CoStObj_{\tx{KM}, \kappa}^{2 \crho - \cmu})_{\infty}),\]
where $t$ is the coordinate around $0 \in \mb{P}^1$ (we use subscript to remind ourselves of the local coordinates). Thus we have reduced the computation of $*$-fiber to a representation-theoretic question.

\paragraph{Contragredient Duality}
Using weight space of the Sugawara $L_0$ operator one can define the \emph{contragredient dual} of $M \in (\kmgmod^I)^\heartsuit$ (c.f. \cite{arkhipov_gaitsgory}), which we denote by $M^*$. The key property of this contragredient dual is that
\[C_*(\mf{g}[t, t^{-1}], M_{0} \tensor (N^*)_{\infty}) \simeq \Hom_{\kmgmod^T}(M, N)^*.\]
The contragredient duality exchanges Verma module $\mb{M}^{\clambda}$ with the contragredient Verma module $\mb{M}^{-\clambda, I^-}$ \emph{for the opposite Iwahori $I^-$}. Likewise, it exchanges $\mb{W}_{\kappa}^{w_0, \clambda}$ with $\mb{W}_{\kappa}^{1, -\clambda, I^-}$. We thus have
\begin{align*}
& (C_*(\mf{g}[t, t^{-1}], (\StObj_{\tx{KM}, \kappa}^{\clambda})_{0} \tensor (\CoStObj_{\tx{KM}, \kappa}^{2 \crho - \cmu})_{\infty})[-\dim(\mf{n})])^* \\
\simeq & \Hom_{\kmgmod^T}((\StObj_{\tx{KM}, \kappa}^{\clambda})_{0}, (\CoStObj_{\tx{KM}, \kappa}^{2 \crho - \cmu})^*)_0)[\dim(\mf{n})] \\
\simeq & \Hom_{\kmgmod^T}(\StObj_{\tx{KM}, \kappa}^{\clambda}, (\mb{W}_{\kappa}^{w_0, 2 \crho - \cmu})^*)[\dim(\mf{n})] \\
\simeq & \Hom_{\kmgmod^T}(\StObj_{\tx{KM}, \kappa}^{\clambda}, \mb{W}_{\kappa}^{1, \cmu - 2\crho, I^-})[\dim(\mf{n})] \\
\simeq & \Hom_{\kmgmod^I}(\StObj_{\tx{KM}, \kappa}^{\clambda}, \tx{Av}_{*, I / T}(\mb{W}_{\kappa}^{1, \cmu - 2\crho, I^-}))[\dim(\mf{n})] \\
\simeq & \Hom_{\kmgmod^I}(\StObj_{\tx{KM}, \kappa}^{\clambda}, \mb{W}_{\kappa}^{w_0, \cmu}) \\
\simeq & C^{\semiinf}(\Loop \mf{n}, \Arc \mf{n}, \mb{W}_\kappa^{w_0, \cmu})^{\clambda} \simeq  \begin{cases}
\mb{C} & \tx{if}~\clambda = \cmu \\
0 & \tx{else}
\end{cases}
\end{align*}
which finishes the proof.

\subsection{Proof of Proposition \ref{prop:KM-E2-Lus-match}}

Our notations are as follows:
\begin{itemize}
\item Let $\mathscr{G}_{\DualTwisting}$ be the gerbe on $[\tx{Gr}_{\check{T}}(\mb{A}^1)]_{\Ran_\dR}$ coming from the ``$\crho$-shifted'' FLE, as well as its restriction to $\tx{Conf}$ (for the constructible case, use its Riemann-Hilbert counterpart);
\item Let $\tx{Conf}^\clambda$ denote the component of $\tx{Conf}$ where the total degree is $\clambda$, and let
\[\Delta_\clambda: X \xrightarrow{x \mapsto \clambda x} \tx{Conf}^\clambda\]
the main diagonal embedding, $j_\clambda$ the complement open embedding;
\item Let $\tx{Conf}^\circ$ the open part of $\tx{Conf}$, consisting of tuples whose coefficients are all negative simple roots.
\end{itemize}

\begin{remark}
The $\crho$-shift has the following effect. Suppose $\clambda = \sum_{i} (- n_i) \cdot \check{\alpha}_i$ is a negative combination of simple roots, then the restriction of $\mathscr{G}_{\DualTwisting}$ to $\tx{Conf}^\clambda$ is given by (see \cite[1.3.5]{gaitsgory2021factorization})
\[\left(\bigotimes_i \mc{O}(-\Delta_i')^{\check{\kappa}(\alpha_i, \alpha_i) / 2}\right) \bigotimes \left(\bigotimes_{i \neq j} \mc{O}(-\Delta_{i, j})^{\check{\kappa}(\check{\alpha}_i, \check{\alpha}_j)}\right);\]
here $\Delta_i'$ is the diagonal divisor in the symmetric power of the curve $X^{(n_i)}$ corresponding to $\clambda$, the second product runs over all unordered pairs $i \neq j$, and $\Delta_{i, j}$ is the incidence divisor in $X^{(n_i)} \times X^{(n_j)}$. Directly inspecting the formula, one sees that the gerbe is canonically trivial on the open part $\tx{Conf}^{\circ}$.
\end{remark}

Under the assumption from Condition \ref{cond:avoid-small-roots-of-unity}, \cite{gaitsgory2021factorization} defines a factorization \emph{perverse sheaf} $\Omega_q^\tx{Lus}$ via a combinatorial procedure, which we recall now:

\begin{definition}
The factorizable algebra $\Omega_{q, \tx{alg}}^\tx{Lus}$ (resp.\ $\Omega_{q, \tx{top}}^\tx{Lus}$) is the regular holonomic D-module (resp.\ perverse sheaf) defined via the following inductive procedure:
\begin{itemize}
\item using the canonical trivialization of the gerbe on $\tx{Conf}^\circ$, we have the \emph{sign local system} on $\tx{Conf}^\circ$;
\item by \emph{stipulating} our sheaf is a factorization algebra, we may assume that the sheaf has been defined away from the image of $\Delta_\clambda$;
\item on $\tx{Conf}^\clambda$, if $\clambda$ can be written as $w(\crho) - \crho$ for $\ell(w) = 2$, we take $!*$-extension along $j_\clambda$ (these encode the \emph{quantum Serre relations});
\item Otherwise, we take $H^0$ of the $!$-extension along $j_\clambda$.
\end{itemize}
\end{definition}

\begin{convention}
For arguments that apply to both algebraic and topological situation equally, we will simply write this factorization algebra as $\Omega_q^\tx{Lus}$.
\end{convention}

As our next criterion demonstrates, this object enjoys strong rigidity property:

\begin{proposition}
\label{prop:omega-lus-criterion}
$\Omega_{q}^\tx{Lus}$ is the \emph{unique} non-unital Ran-constructible factorization algebra within
\[\DMod_{\mathscr{G}_{\DualTwisting}}([\tx{Gr}_{\check{T}}]_{\Ran}) \hspace{0.5em} \tx{or} \hspace{0.5em} \tx{Shv}^!_{\tx{RH}(\mathscr{G}_{\DualTwisting})}([\tx{Gr}_{\check{T}}]_{\Ran})\] satisfying the following cohomological estimates: for every $x \in X(\mb{C})$ and every $\lambda \in \weightLat$,
\begin{enumerate}
\item if $\clambda \notin \weightLat^{<0}$, then the $!$-fiber at $\clambda \cdot x$ is zero;
\item the $!$-fiber at every $\clambda \cdot x$ has no negative cohomology;
\item if $\clambda$ is a simple negative root, then either the $*$-fiber at $\clambda \cdot x$ is $\mb{C}[1]$, or the $!$-fiber at $\clambda \cdot x$ is $\mb{C}[-1]$;
\item if $\clambda$ equals $w(\crho) - \crho$ for some $\ell(w) = 2, w \in W_\tx{fin}$, then the $!$-fiber at $\clambda \cdot x$ vanishes at $H^0$ and $H^1$, and $*$-fiber at $\clambda \cdot x$ vanishes at $H^{0}$ and $H^{-1}$;
\item otherwise, the $!$-fiber at $\clambda \cdot x$ vanishes at $H^0$, and $*$-fiber at $\clambda \cdot x$ vanishes at $H^{0}$, $H^{-1}$ and $H^{-2}$.
\end{enumerate}
\end{proposition}

\begin{proof}
	We first verify that $\Omega_q^\tx{Lus}$ satisfies the conditions listed above. Condition (1) and (2) are evident; condition (3) is \cite[3.5.9]{gaitsgory2021factorization}. For condition (4) and (5), (the Koszul dual of) \cite[3.7.3, (a) and (b)]{gaitsgory2021factorization} gives the $*-$ and $!$-estimate respectively, with the $!$-condition for (5) being evident.

	Now suppose $\mc{F}$ satisfies condition (1) to (5) above. Condition (1) guarantees that the sheaf is supported on the configuration space. Since the sheaf is locally constant along the diagonal stratification and has factorization property, condition (3) guarantees that it is the sign local system on the open locus. The vanishing condition on $H^0$ from (2), (4) and (5) guarantee that $\mc{F}$ is actually a perverse sheaf over the entire configuration space.
	
	We will now show that $\mc{F}$ is uniquely determined on each $\tx{Conf}^\clambda$ by running induction on the total degree of $\clambda$. The base case is $|\clambda| = 1$, which we already argued above. So we assume $|\clambda| \ge 2$. By factorization and induction, we can assume $\mc{F}$ is already uniquely determined away from the main diagonal. Thus the proposition would follow from the following claims:
	\begin{itemize}
	\item if $\clambda = w(\crho) - \crho$ for $\ell(w) = 2$, then $\mc{F} \simeq H^0((j_\clambda)_! (j_\clambda)^! \mc{F})$;
	\item otherwise, $\mc{F} \simeq (j_\clambda)_{!*} (j_\clambda)^! \mc{F}$.
	\end{itemize}
	For (b), note that we have an exact sequence
	\[H^0((j_\clambda)_! (j_\clambda)^! (\mc{F})) \to \mc{F} \to H^0((\Delta_{\clambda})_* (\Delta_{\clambda})^* (\mc{F}))\]
	and the last term vanishes by $*$-fiber condition in (4), which means the first arrow is surjective; at the same point, we also have the exact sequence
	\[H^0((\Delta_{\clambda})_! (\Delta_{\clambda})^! (\mc{F})) \to \mc{F} \to H^0((j_\clambda)_* (j_\clambda)^* (\mc{F}))\]
	and the first term vanishes by $!$-fiber condition in (4), so the second map is injective. So we have
	\[H^0((j_\clambda)_! (j_\clambda)^! (\mc{F})) \surj \mc{F} \inj H^0((j_\clambda)_* (j_\clambda)^* (\mc{F}))\]
	implying point (b). For (a), note that we have
	\[H^{-1}((\Delta_{\clambda})_* (\Delta_{\clambda})^* (\mc{F})) \to H^0((j_\clambda)_! (j_\clambda)^! (\mc{F})) \to \mc{F} \to H^0((\Delta_{\clambda})_* (\Delta_{\clambda})^* (\mc{F}))\]
	and the first and the last term vanish by the $*$-fiber estimate in (4).
\end{proof}

Our remaining job is to use this criterion to match the factorization algebras produced from both the Kac-Moody side and the quantum side.

\subsubsection{Kac-Moody Side}

\begin{definition}
\label{def:kappa-almost-admissible}
Fixed a bilinear form $\kappa$. For each simple factor $\mf{g}_i$ of $\mf{g}$, let $\kappa|_{\mf{g}_i} = \frac{c - h_i^\vee}{2 h_i^\vee} \tx{Kil}_{\mf{g}_i}$; we say $\kappa$ is \emph{almost admissible} on $\mf{g}_i$ if either $c$ is irrational, or $c = \frac{p}{q}$ for coprime $p, q \in \mb{Z}^{> 0}$, such that
\[p \ge \begin{cases}
h_i^\vee - 1 & \tx{if}~(q, d_i) = 1 \\
h_i - 1 & \tx{else}
\end{cases}\]
here $d_i$ is the lacing number of $\mf{g}_i$, and $h_i$ (resp.\ $h^\vee_i$) is the Coxeter (resp.\ dual Coxeter) number of $\mf{g}_i$. We say $\kappa$ is \emph{almost admissible} if it is so on each simple factor of $\mf{g}$.
\end{definition}

\begin{remark}
We call this \emph{almost admissible} because the lower bound on $p$ is precisely 1 below that of the classical notion of admissibility as defined in \cite{kac2008rationality}.
\end{remark}

First observe that, even though $w \cdot (-2 \crho)$ is not necessarily dominant, we still have
\begin{lemma}
\label{lemma:wa-is-verma-on-2crho}
	When $\kappa$ is almost admissible, the canonical map $\mb{M}_{-\kappa}^{w \cdot (-2\crho)} \to \mb{W}_{-\kappa}^{1, w \cdot (-2 \crho)}$ is an isomorphism for every $w \in W^\tx{fin}$.
\end{lemma}
	\begin{proof}
	Choose sufficiently large $\mu \gg 0$ such that $w \cdot (-2 \crho) + \kappa(\mu) \in \weightLat^+$. By \cite[Section 3.4.3]{gaitsgory2021conjectural}, we have
	\[\mb{M}_{-\kappa}^{w \cdot (-2 \crho) + \kappa(\mu)} \simeq \mb{W}_{-\kappa}^{1, w \cdot (-2 \crho) + \kappa(\mu)};\]
	by \cite[Section 2.5.9]{gaitsgory2021conjectural}, (and because $J_{\mu} \star \bullet$ for $\mu \gg 0$ is invertible,) the claim follows if we prove
	\[j_{\mu, !} \star \mb{M}_{-\kappa}^{w \cdot (-2 \crho)} \simeq \mb{M}_{-\kappa}^{w \cdot (-2 \crho) + \kappa(\mu)}.\]
	Note that $-2 \crho$ is $\kappa$-admissible weight by our assumption: the first condition of Definition \ref{def:admissible-weight} is trivial, and the second condition is guaranteed precisely when $\kappa$ is almost admissible. Thus Theorem \ref{thm:kashiwara-tanisaki} can be invoked, and the claim above reduces to the following claim about $-\kappa$-twisted D-modules on the affine flag variety:
	\[j_{\mu, !} \star j_{w, !} \simeq j_{\mu w, !};\]
	this is guaranteed if we have $\ell(\mu) + \ell(w) = \ell(\mu w)$, which in turn follows from the following general length formula for elements in the extended affine Weyl group: for $\mu_0 \in \coweightLat$ and $w_0 \in \Wfin$, we have
	\[\ell(\mu_0 w_0) = \sum_{\alpha \in \Delta^+, w_0^{-1}(\alpha) \in \Delta^+}|\langle \mu_0, \alpha \rangle| + \sum_{\alpha \in \Delta^+, w_0^{-1}(\alpha) \in \Delta^-} |1 + \langle \mu_0, \alpha \rangle|.\]
\end{proof}

\begin{corollary}
When $\kappa$ is positive and is almost admissible, there exists a BGG resolution
\[\left[0 \to \StObj_{\tx{KM}, \kappa}^{w_0 \cdot 0} \to \ldots \to \bigoplus_{w \in \Wfin, \ell(w) = 1} \StObj_{\tx{KM}, \kappa}^{w \cdot 0} \to \underline{\StObj_{\tx{KM}, \kappa}^{0}} \to 0\right] \simeq \mb{V}_\kappa^0\]
of the vacuum by standard objects at positive level. Here on LHS the underline marks cohomological degree $0$ (thus e.g. the object $\StObj_{\tx{KM}, \kappa}^{0}$ is placed at cohomological degree 0).
\end{corollary}

\begin{proof}
	At level $-\kappa$, we have a standard BGG resolution of the vacuum by Verma modules
	\[\left[0 \to \mb{M}_{-\kappa}^{-2 \crho} \to \bigoplus_{w \in \Wfin, \ell(w) = 1} \mb{M}_{-\kappa}^{w \cdot (-2\crho)} \to \ldots \to \underline{\mb{M}_{-\kappa}^{w_0 \cdot (-2\crho)}} \to 0\right] \simeq \mb{V}_{-\kappa}^0;\]
	by Lemma \ref{lemma:wa-is-verma-on-2crho}, this can be rewritten as
	\[\left[0 \to \mb{W}_{-\kappa}^{1, -2 \crho} \to \bigoplus_{w \in \Wfin, \ell(w) = 1} \mb{W}_{-\kappa}^{1, w \cdot (-2\crho)} \to \ldots \to \underline{\mb{W}_{-\kappa}^{1, w_0 \cdot (-2\crho)}} \to 0\right] \simeq \mb{V}_{-\kappa}^0;\]
	Given that $\mb{D}_{\tx{KM}}(\mb{V}_{-\kappa}^0) \simeq \mb{V}_{\kappa}^0$, the original claim now follows from applying $\mb{D}_{\tx{KM}}$ to the equation above (note that it flips the direction of arrows).
\end{proof}

The following estimates will be needed in light of Proposition \ref{prop:omega-lus-criterion}:

\begin{proposition}
\label{prop:KM-vacuum-verification}
	If $\kappa$ is positive and almost admissible, $\mb{V}_\kappa^0$ is the vacuum module, then:
	\begin{enumerate}
	\item We have
	\[\Hom_{\kmgmod^I}(\mb{V}_\kappa^0, \CoStObj_{\tx{KM}, \kappa}^{\clambda}) = \begin{cases}
\mb{C}[-\ell(w_f)] & \clambda = w_f(\crho) - \crho, w_f \in \Wfin \\
0 & \tx{else}
\end{cases}\]
	\item For every $\clambda \in \weightLat$ $C_*^{\semiinf}(\Loop \mf{n}, \Arc \mf{n}, \mb{V}_\kappa^0)^\clambda$ vanishes unless $\clambda \in \weightLat^{\le 0}$, is concentrated in cohomological degree $\ge 0$, and is nonzero in degree $0$ only if $\clambda = 0$;
	\item For $\clambda = w(\crho) - \crho$, $w \in \Wfin, \ell(w) = 2$, $H^1(C_*^{\semiinf}(\Loop \mf{n}, \Arc \mf{n}, \mb{V}_\kappa^0)^\clambda) = 0$.
	\end{enumerate}
\end{proposition}

\begin{proof}
	(1) follows directly from the previous corollary.
	
	For (2), let $C := C_*^{\SI}(\Loop \mf{n}, \Arc N, \mb{V}_{\kappa}^0)$ denote the semi-infinite complex for the vacuuum, and let $(\mf{n}_i)_{i \ge 0}$ denote a sequence of compact-open subalgebras
	\[\Arc \mf{n} \simeq \mf{n}_0 \subset \mf{n}_1 \subset \ldots \subset \bigcup_{i \ge 0} \mf{n}_i \simeq \Loop \mf{n};\]
	As in \cite[Appendix A.8]{raskin2021mathcal}, one may impose a $\mb{Z}^-$-filtered PBW-type filtration
	\[\ldots \subset F_{-1} C \subset F_0 C \simeq C\]
	on the semi-infinite chain complex, such that $\limit_{-i \to -\infty} F_{-i} C \simeq 0$ and
	\[\tx{gr}_{-i} C \simeq C^*(\Arc \mf{n}, \tx{Sym}^i(\Loop (\mf{g} / \mf{n}) / \Arc (\mf{g} / \mf{n}))).\]
	Note that for each $i$, $H^\bullet(\tx{gr}_{-i} C)^{\clambda}$ is zero unless $\clambda \in \weightLat^{< 0}$; we also have
	\[H^0(\tx{gr}_{-i} C) \subseteq H^0(\mf{n}, \tx{Sym}^i(\mf{g} / \mf{n})) \subseteq \tx{Sym}(\mf{g} / \mf{n})^{\mf{n}} \simeq \tx{Sym}(\mf{h})\]
	so we see that $H^0(\tx{gr}_{-i} C)^{\clambda} = 0$ unless $\clambda = 0$. Now we have, within $\KLT^{<0}$,
	\[0 \simeq \tau_{< 0}~0 \simeq \tau_{< 0} \limit_{-i} F_{-i} C \simeq \limit_{-i} \tau_{< 0} F_{-i} C;\]
	as we have $\tau_{<0}~\tx{gr}_{-i} C \simeq 0$, from the usual exact triangle we see that $\tau_{<0} F_{-i-1} C \to \tau_{<0} F_{-i} C$ is an isomorphism for each $i$, so that limit diagram is constant and each term is zero; in particular we have $\tau_{<0} C \simeq \tau_{<0} F_0 C \simeq 0$. The other two claims of (2) now follows the same argument.
	
	It remains to compute (3). We have
	\[\Hom_{\kmgmod^I}(\StObj_{\tx{KM}, \kappa}^{w_f \cdot 0}, \mb{V}_\kappa^0) \simeq \Hom_{\nkmgmod^I}(\mb{D}_{\tx{KM}}(\mb{V}_\kappa^0), \mb{D}_{\tx{KM}}(\StObj_{\tx{KM}, \kappa}^{w_f \cdot 0}))\]
	\[\simeq \Hom_{\nkmgmod^I}(\mb{V}_{-\kappa}^0, \mb{W}_{-\kappa}^{1, w_f \cdot (-2 \crho)}[\dim(\mf{n})]) \simeq \Hom_{\nkmgmod^I}(\mb{V}_{-\kappa}^0, \mb{M}_{-\kappa}^{w_f \cdot (-2 \crho)})[\dim(\mf{n})],\]
	where in the last step we used Lemma \ref{lemma:wa-is-verma-on-2crho} above.
	Because $\kappa$ is almost admissible, we can invoke Kashiwara-Tanisaki for the weight $(-2 \crho)$. Consider the closed embedding of the basic affine space
	\[\iota_{\tx{Fin}}: G / B \to \Arc G / I \simeq \tx{Fl}_\tx{G};\]
	the weight $(-2 \crho)$ determines a twisting on $G / B$ that is compatible with $(-\kappa, -2 \crho)$ via pullback along $\iota_{\tx{Fin}}$, so we have a well-defined functor $\iota_{\tx{Fin}, !}: \DMod_{(-2 \crho)}(G / B) \to \DMod_{(-\kappa, -2 \crho)}(\tx{Fl}_G)$, and an induced functor	
	\[\iota_{\tx{Fin}, !}^N: \DMod_{(-2 \crho)}(G / B)^N \to \DMod_{(-\kappa, -2 \crho)}(\tx{Fl}_G)^{I^\circ}.\]
	For each $w_f \in \Wfin$, $j_{\tx{Fin}, w_f, !}$ (resp.\ $j_{\tx{Fin}, w_f, *}$) denote the $!$ (resp.\ $*$) extension on $G / B$ of the IC sheaf of the corresponding $N$-orbit; these correspond to finite Verma (resp.\ dual Verma) modules under the usual Beilinson-Bernstein localization. We note that the $w_f$-orbits in the affine flag variety are isomorphic to these finite orbits under $\iota_\tx{Fin}$. Further observe that by base change we have
	\[\Gamma_{-\kappa}(\tx{Fl}_G, \tx{IC}_{G / B})^{I^\circ} \simeq \mb{V}_{-\kappa}^0;\]
	here in writing $\tx{IC}_{G / B}$ we note that $(-2 \crho)$ is integral. We claim that this means
	\[\Hom_{\DMod_{(-\kappa, -2 \crho)}(\tx{Fl}_G)^{I^\circ}}(\tx{IC}_{G/B}, j_{w_f, !}) \simeq \Hom_{\nkmgmod^{I^\circ}}(\mb{V}_{-\kappa}^0, \mb{M}_{-\kappa}^{w_f \cdot (-2 \crho)});\]
	indeed, $\tx{IC}_{G/B}$ admits a \emph{finite} resolution by standard objects, and $j_{w_f, !}$ a \emph{finite} resolution by costandard objects (\cite[Lemma 15]{arkhipov2009perverse}). Because everything commutes with finite colimits, the claim follows from the fact that standard and costandard objects are orthogonals in both $\DMod_{(-\kappa, -2 \crho)}(\tx{Fl}_G)^{I_0}$ and $\nkmgmod^{I_0}$. Next, note that
	\[\Hom_{\DMod_{-2 \crho}(G / B)^N}(\tx{IC}_{G/B}, j_{\tx{Fin}, w_f, !})[\dim(\mf{n})] \simeq \Hom_{\mf{g}\mod^N}(k, \mb{M}_\tx{fin}^{w_f \cdot (-2 \crho)})[\dim(\mf{n})]\]
	\[\simeq \Hom_{\mf{g}\mod^N}(\mb{M}_{\tx{fin}}^{w_f \cdot 0}, k) \simeq \mb{C}[-\ell(w_f)]\]
	where we use the finite analogue of $\mb{D}_{\tx{KM}}$ in the middle. Now the claim original claim follows from the fact that
	\[\Hom_{\nkmgmod^{I}}(\mb{V}_{-\kappa}^0, \mb{M}_{-\kappa}^{w_f \cdot (-2 \crho)}[\dim(\mf{n})]) \simeq H_T^*(\Hom_{\nkmgmod^{I^\circ}}(\mb{V}_{-\kappa}^0, \mb{M}_{-\kappa}^{w_f \cdot (-2 \crho)}[\dim(\mf{n})])).\]
\end{proof}

\begin{lemma}
	If $\kappa$ is positive and almost admissible, we have $\TheFactAlg \simeq \Omega_{q, \tx{alg}}^\tx{Lus}$.
\end{lemma}

\begin{proof}
Proposition \ref{cor:omega-km-constructible} allows us to invoke Proposition \ref{prop:omega-lus-criterion}. We already know that the $!$-fibers of $\TheFactAlg$ are given by semi-infinite cohomology of the vacuum; to interpret the $*$-fibers, we invoke Theorem \ref{thm:jacquet-is-equivalence}. Recall from Section \ref{sect:conf-hw} that the costandard object $\CoStObj_{\tx{Conf}, \kappa}^{\clambda}$ corepresent the functor of taking $*$-fibers at $\clambda \cdot x$ within the factorization module category. Let $\iota_{\clambda \cdot x}$ denote the inclusion, then we have
\[\iota_{\clambda \cdot x}^*(\TheFactAlg) \simeq \Hom_{\TheFactAlg\tx{-FactMod}_\nonun(\DMod_{\DualTwisting}(\tx{Conf}_{\infty \cdot x_0}))}(\TheFactAlg, \CoStObj_{\tx{Conf}, \kappa}^{\clambda})\]
\[\simeq \Hom_{\kmgmod^I}(\mb{V}_\kappa^0, \CoStObj_{\tx{KM}, \kappa}^{\clambda}).\]
Now we can invoke (1) of Proposition \ref{prop:KM-vacuum-verification} to finish the argument.
\end{proof}

\subsubsection{Quantum Side}

Recall that in \cite{gaitsgory2021factorization} the factorization algebra $\Omega_q^\tx{Lus}$ was related to $U_q^\tx{KD}(N^-) \in \tx{BiAlg}(\Rep_q(T))$ (defined in Definition \ref{def:lus-and-kd}) by appealing to a general construction in \cite{bezrukavnikov2006factorizable} that uses hyperbolic restriction to produce a factorizable perverse sheaf from a (graded) Hopf algebra in the heart. We take the following definition instead:
\begin{definition}
We define a topological factorization algebra $\Omega_{\mb{E}_2, q}^\tx{Lus}$ in the topological factorization category $\tx{Shv}^!_{\tx{RH}(\mathscr{G}_{\DualTwisting})}([\tx{Gr}_{\check{T}}]_{\Ran})$ to be the image of
\[\tx{AugIdeal}(\tx{coBar}(U_{q}^{\tx{KD}}(N^-))) \in \mb{E}_2\tx{-Alg}_{\tx{nu}}(\Rep_q(T))\]
under the equivalence in Proposition \ref{prop:E2Alg-Mod-to-FactAlg-Mod}.
\end{definition}

\begin{lemma}
	If $q$ is non-degenerate, avoids small torsion and almost admissible, we have $\Omega_{\mb{E}_2, q}^\tx{Lus} \simeq \Omega_{q, \tx{top}}^\tx{Lus}$.
\end{lemma}

\begin{proof}
	From Proposition \ref{prop:fact-koszul-verdier-compatibility}, using the fact that
	\[\Hom_{U_{q}^\tx{KD}(N^-)\tx{-comod}(\Rep_{q}(T))}(k, k) \simeq \Hom_{U_{q}^\tx{Lus}(N)\tx{-mod}_{\tx{loc.nilp}}(\Rep_{q}(T))}(k, k) \simeq \Hom_{U_{q}^\tx{Lus}(N)\tx{-mod}(\Rep_{q}(T))}(k, k)\]
	as in \cite[Lemma 4.3.5]{gaitsgory2021conjectural} (note that $k$ is \emph{by definition} compact in $U_{q}^\tx{Lus}(N)\tx{-mod}_{\tx{loc.nilp}}(\Rep_{q}(T))$), it suffices to check that if $q$ satisfies the conditions listed at the beginning of Section \ref{sect:prep-q}, then we have the following cohomological estimates:
\begin{enumerate}
\item $C^*(U_q^\tx{Lus}(N), k)^{\clambda}$ is zero unless $\clambda \in \weightLat^{<0}$, is always concentrated in degree $\ge 0$, and only has nonzero $H^0$ if $\clambda = 0$;
\item If $\clambda$ is a simple negative root, then $C_*(U_{q^{-1}}^\tx{KD}(N^-), k)^{\clambda} \simeq \mb{C}[1]$;
\item If $\clambda = w(\crho) - \crho$ for $\ell(w) = 2$, $w \in W_\tx{fin}$, then $H^1(U_q^\tx{Lus}(N), k)^{\clambda} \simeq 0$ and $H_1(U_{q^{-1}}^\tx{KD}(N^-), k)^{\clambda} \simeq 0$;
\item For every other $\clambda$, we have $H_2(U_{q^{-1}}^\tx{KD}(N^-), k)^{\clambda} \simeq 0$.
\end{enumerate}
Consider the resolution of $k$ given by
\[\ldots \xrightarrow{\delta_2} U \tensor \ker(\delta_0) \xrightarrow{\delta_1} U \xrightarrow{\delta_0} k \to 0\]
for either $U = U_q^\tx{Lus}(N)$ or $U = U_{q^{-1}}^{\tx{KD}}(N^-)$. Here each $\delta_i$ is the multiplication map (so $\delta_0$ is the augmentation ideal). (1) and (2) then obvious. Both $H^1$ in (3) are spanned by generators of the algebra, so (3) follows from the fact that none of the generators (for either algebra) has the prescribed form. For (4), note that $H_2$ is a quotient of $\ker(I \tensor I \xrightarrow{m} I)$ where $I$ is the augmentation ideal of $U_{q^{-1}}^\tx{KD}(N^-)$, i.e. the span of relations. As noted in \cite[Remark 3.3.7]{gaitsgory2021factorization}, when the parameter $q$ satisfies said conditions, this space is spanned by quantum Serre relations and thus has no other weight.
\end{proof}

\section{Apppendix A: Lax Limits} %
\label{sect:lax-limits}
In this section we develop some basic properties of lax limits parallel to those in \cite{ayala2021stratified}, which does not immediately follow from \emph{loc.cit.} because we need to work with the unstable non-presentable setting.

Recall from \cite{ayala2021stratified} the definition of a \emph{lax limit} $\ms{Glue}(\Phi)$ corresponding to a \emph{lax functor} $\Phi: I \dashrightarrow \tx{Cat}$. Intuitively, given a poset $(I, >)$, a lax functor $\Phi$ encodes the data
\[i \in \mc{I} \mapsto \Phi(i) \in \tx{Cat} \hspace{1em} (i > j) \leadsto \Phi_{ij}: \Phi(i) \to \Phi(j) \hspace{1em} (i > j > k) \leadsto \eta_{ijk}: \Phi_{jk} \circ \Phi_{ij} \Rightarrow \Phi_{ik}\]
and a lax limit encodes the data
\[x_i \in \Phi(i), i > j \leadsto \Phi_{ij}(x_i) \xrightarrow{t_{ij}} x_j\]
compatible withe the data of $\Phi_{ij}$ and $\eta_{ijk}$ above. Here's the formal definition:

\begin{definition}
A (right-)lax functor $\Phi: I \dashrightarrow \tx{Cat}$ is a locally Cartesian fibration $\mc{E}_\Phi$ over $I^\op$. Its (left-)lax limit $\ms{Glue}(\Phi)$ is the category of functors $\ms{sd}(I^\op)^\op \to \mc{E}_\Phi$ that preserve locally Cartesian\footnote{In Notation A.4.1 of \emph{loc.cit.} \emph{Cartesian} arrows is written instead, but this is likely a typo.} arrows over $I^\op$, where $\ms{sd}$ is the category of subdivisions defined in \emph{loc.cit.} %
\end{definition}
We also consider the oplax limit $\ms{opGlue}(\Psi)$ of an oplax functor $\Psi$, defined analogously.

\begin{convention}
Our usage of ``lax functor'' refers to \emph{right}-lax functors in \cite{ayala2021stratified}, and ``oplax functor'' refers to \emph{left}-lax functors; on the other hand, our usage of ``lax limits'' refer to their \emph{left}-lax limits, and ``oplax limits'' refer to their \emph{right}-lax limits. This sightly unfortunate hodgepodge is in alignment with the usage in \cite{GR-DAG1}.
To maintain some sanity we will only speak of lax-limits of lax-functors and oplax-limits of oplax-functors; the other two combinations will not occur in this text.
\end{convention}

\begin{remark}
	It follows from \cite[2.4.6.5]{HTT} that, because $I$ is a poset, any lax functor is automatically an isofibration (a.k.a. categorical fibration).
\end{remark}

Here is an example demonstrating this encoding. Let us suppose $I$ is the poset $[2]^\op := (2 > 1 > 0)$, and define $\Phi(i) := \mc{E}_{\phi}|_{i}, i \in I$. The category $\ms{sd}(I^\op)^\op$ looks like follows:
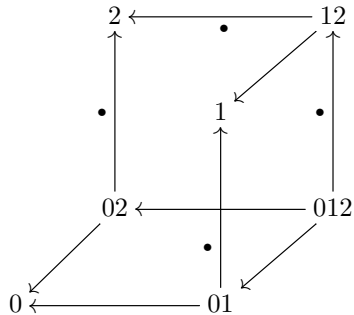
\[
\xymatrix{
 & 2 &  & 12 \ar^{\bullet}[ll] \ar[ld] \\
 &  & 1 &  \\
 & 02 \ar^{\bullet}[uu] \ar[ld] &  & 012 \ar[ll] \ar^{\bullet}[uu] \ar[ld] \\
0 &  & 01 \ar[ll] \ar^(0.3){\bullet}[uu] & 
}\]
and it is a locally Cartesian fibration over $[2]$ via the map $\ms{min}$ (the marked arrows are locally Cartesian). Thus, a typical element in the lax limit looks like
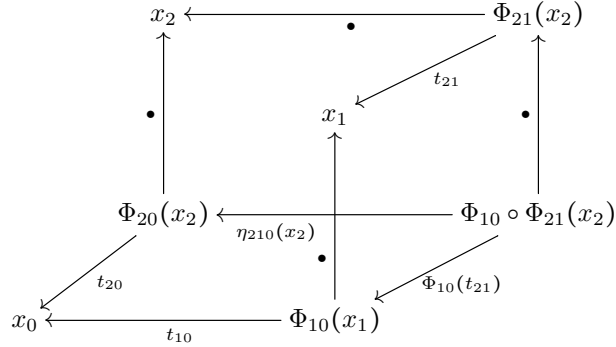
\[
\xymatrix{
 & x_2 &  & \Phi_{21}(x_2) \ar^{\bullet}[ll] \ar^{t_{21}}[ld] \\
 &  & x_1 &  \\
 & \Phi_{20}(x_2) \ar^{\bullet}[uu] \ar^{t_{20}}[ld] &  & \Phi_{10} \circ \Phi_{21}(x_2) \ar^(0.7){\eta_{210}(x_2)}[ll] \ar^{\bullet}[uu] \ar^{\Phi_{10}(t_{21})}[ld] \\
x_0 &  & \Phi_{10}(x_1) \ar^{t_{10}}[ll] \ar^(0.3){\bullet}[uu] & 
}
\]
Here is an alternative description, useful in its own right (\cite[A.6.6]{ayala2021stratified}):
\begin{lemma}
\label{lemma:glue-as-strict-limit}
	There exists a diagram $\mf{G}(\Phi): \ms{sd}(I) \to \tx{Cat}$ such that $\ms{Glue}(\Phi) \simeq \limit_{\ms{sd}(I)} \mf{G}(\Phi)$.
\end{lemma}
A typical element in this description looks like (note the diagram has been flipped):
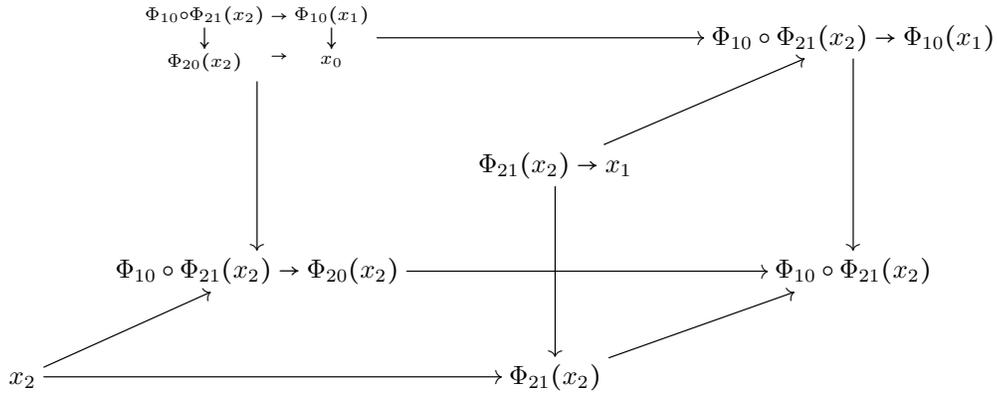
\[
\xymatrix{
 & \underset{\Phi_{20}(x_2)}{\overset{\Phi_{10} \circ \Phi_{21}(x_2)}{\downarrow}} \underset{\rightarrow}{\overset{\rightarrow}{\phantom{?}}} \underset{x_0}{\overset{\Phi_{10}(x_1)}{\downarrow}} \ar[dd] \ar[rr] &  & \Phi_{10} \circ \Phi_{21}(x_2) \to \Phi_{10}(x_1) \ar[dd] \\
 &  & \Phi_{21}(x_2) \to x_1 \ar[ru] \ar[dd] &  \\
 & \Phi_{10} \circ \Phi_{21}(x_2) \to \Phi_{20}(x_2) \ar[rr] &  & \Phi_{10} \circ \Phi_{21}(x_2) \\
x_2 \ar[ru] \ar[rr] &  & \Phi_{21}(x_2) \ar[ru] & 
}
\]
as an element of
\[
\xymatrix{
 & \tx{Fun}([1] \times [1], \Phi(0)) \ar^{(\tx{src}, \tx{id})}[dd] \ar^{(\tx{id}, \tx{src})}[rr] &  & \tx{Fun}([1], \Phi(0)) \ar^{\tx{src}}[dd] \\
 &  & \tx{Fun}([1], \Phi(1)) \ar^{\Phi_{10}}[ru] \ar^(0.3){\tx{src}}[dd] &  \\
 & \tx{Fun}([1], \Phi(0)) \ar^(0.3){\tx{src}}[rr] &  & \Phi(0) \\
\Phi(2) \ar^{\eta_{210}}[ru] \ar^{\Phi_{21}}[rr] &  & \Phi(1) \ar^{\Phi_{10}}[ru] & 
}
\]

\begin{remark}
\label{remark:lax-functor-fully-faithful}
Let $\Phi, \Phi': I \dashrightarrow \tx{Cat}$ be lax functors, and suppose we are given a fully faithful morphism $i: \Phi' \to \Phi$, by which we mean a morphism $\mc{E}_{\Phi'} \to \mc{E}_\Phi$ preserving locally Cartesian arrows, and which is a fiberwise fully faithful embedding of categories. Then we have an induced fully faithful embedding $\ms{Glue}(\Phi') \subseteq \ms{Glue}(\Phi)$, whose essential image consists of those objects in $\mf{G}(\Phi)$ which lands in $\Phi'$ in each component.
\end{remark}

We now summarize the main properties of lax limits:

\begin{proposition}
Let $\Phi: I \dashrightarrow \tx{Cat}$ be a lax functor such that $\Phi(i)$ all contain small colimits and $\Phi_{ij}$ all commutes with small colimits. Let $I_0 \subseteq I$ a full subposet. Then the following hold:
\label{prop:lax-limit-properties}
\begin{enumerate}
\item $\ms{Glue}(\Phi)$ contains small colimits, and the projection functor $\tx{ev}_{I_0}: \ms{Glue}(\Phi) \to \ms{Glue}(\Phi|_{I_0})$ preserve small colimits;
\item The canonical functor $\tx{ev}: \ms{Glue}(\Phi) \to \prod_{i \in I} \Phi(i)$ is conservative;
\item If $I_0$ is either upward closed\footnote{i.e. if $(a \to b) \in D$ and $a \in D$ then $b \in D$.} or an initial object $\{0\}$\footnote{In fact, the statement is true for \emph{any} full subposet $I_0$, but the proof is more involved so we opt out of it for now.}, then $\tx{ev}_{I_0}$ admits a fully faithful left adjoint which we denote by $\tx{ins}_{I_0}$. In the upward-closed case, the essential image consists of those objects $x$ such that $\tx{ev}_i(x)$ is the initial object in $\Phi(i)$ for $i \not\in I_0$; in the initial singleton case, the essential image consists of those objects where each $x_i \simeq \Phi_{0i}(x_0)$ and each $t_{ij}$ is the map induced by $\eta_{0ij}$.
\end{enumerate}
Further, if $\Phi(i)$ each admit small limits, then so does $\ms{Glue}(\Phi)$.
\end{proposition}

\begin{proof}
	The last statement follows directly from Lemma \ref{lemma:glue-as-strict-limit}. For (1), note that when $\Phi(i)$ all contain small colimits, so do all terms in $\mf{G}(\Phi)$, and the connecting functors in $\mf{G}(\Phi)$ are all continuous when $\Phi_{ij}$ are; then (1) follows from the fact that in such case, colimits in strict limits are computed component-wise. (2) is because every vertex in $\ms{sd}(I)^\op$ is connected by an (unique) path of locally Cartesian arrows to one in $I^\op$. It remains to prove (3).
	Let $\mc{E} \to I^\op$ denote the locally Cartesian fibration corresponding to $\Phi$, and $\mc{E}_0 := \mc{E} \times_{I^\op} I_0^\op \to I_0^\op$ denote the restriction to $I_0$.
	First we consider when $I_0$ is upward-closed. Unfolding the construction of $\ms{sd}$, one sees that $I_0$ being upward closed implies that the relative Kan extension $\tx{LKE}_i$ depicted by the diagram below
	\[\xymatrix{
	\ms{sd}(I_0^\op)^\op \ar[r] \ar[d] & \mc{E} \ar[d] \\
	\ms{sd}(I^\op)^\op \ar@{-->}[ru] \ar[r] & I^\op
	}\]	
	factors through the subcategory $\ms{Glue}(\Phi)$; in other words, we have uniquely defined $\tx{LKE}_i'$ and $\tx{res}_i'$ as below that fits into the following commutative squares:
\[
\xymatrix{
\tx{Fun}_{\tx{loc.Cart}_{/I^\op}}(\ms{sd}(I_0^\op)^\op, \mc{E}) \ar^{\tx{LKE}_i'}[r] \ar^{\oblv}[d] & \tx{Fun}_{\tx{loc.Cart}_{/I^\op}}(\ms{sd}(I^\op)^\op, \mc{E}) \ar^{\oblv}[d] \\
\tx{Fun}_{/I^\op}(\ms{sd}(I_0^\op)^\op, \mc{E}) \ar^{\tx{LKE}_i}[r] & \tx{Fun}_{/I^\op}(\ms{sd}(I^\op)^\op, \mc{E})
}
\]
\[
\xymatrix{
\tx{Fun}_{\tx{loc.Cart}_{/I^\op}}(\ms{sd}(I^\op)^\op, \mc{E}) \ar^{\tx{res}_i'}[r] \ar^{\oblv}[d] & \tx{Fun}_{\tx{loc.Cart}_{/I^\op}}(\ms{sd}(I_0^\op)^\op, \mc{E}) \ar^{\oblv}[d] \\
\tx{Fun}_{/I^\op}(\ms{sd}(I^\op)^\op, \mc{E}) \ar^{\tx{res}_i}[r] & \tx{Fun}_{/I^\op}(\ms{sd}(I_0^\op)^\op, \mc{E})
}
\]
where $(\tx{LKE}_i, \tx{res}_i)$ is the adjoint pair of relative Kan extensions (\cite[4.3.2.17]{HTT}), so that $(\tx{LKE}_i', \tx{res}_i')$ is again a pair of adjunction. Now we claim that we have an equivalence
\begin{equation}
\label{eq:locCart-bc}
\tx{Fun}_{\tx{loc.Cart}_{/I_0^\op}}(\ms{sd}(I_0^\op)^\op, \mc{E}_0) \simeq \tx{Fun}_{\tx{loc.Cart}_{/I^\op}}(\ms{sd}(I_0^\op)^\op, \mc{E})
\end{equation}
given by composing with $\mc{E}_0 \to \mc{E}$ and base change. Under this equivalence, $\tx{ev}_{I_0}$ becomes $\tx{res}'_i$, so we set $\tx{ins}_{I_0} := \tx{LKE}'_i$. The original proposition follows immediately.
To see the claim, note that we have an equivalence of categories
\[\tx{Fun}_{/I_0^\op}(\ms{sd}(I_0^\op)^\op, \mc{E}_0) \simeq \tx{Fun}_{/I^\op}(\ms{sd}(I_0^\op)^\op, \mc{E})\]
Indeed, there is a natural functor from LHS to RHS given by composition; to see that it is an equivalence, note that it suffices to check
\[\tx{Maps}(\Delta^n, \tx{Fun}_{/I_0^\op}(\ms{sd}(I_0^\op)^\op, \mc{E}_0)) \xrightarrow{\simeq} \tx{Maps}(\Delta^n, \tx{Fun}_{/I^\op}(\ms{sd}(I_0^\op)^\op, \mc{E}))\]
for all $\sigma: \Delta^n \to \tx{Fun}_{/I_0^\op}(\ms{sd}(I_0^\op)^\op, \mc{E}_0)$; This then becomes
\[\tx{Maps}_{/I_0^\op}(\Delta^n \times \ms{sd}(I_0^\op)^\op, \mc{E}_0) \xrightarrow{\simeq} \tx{Maps}_{/I^\op}(\Delta^n \times \ms{sd}(I_0^\op)^\op, \mc{E})\]
which is the universal property of $\mc{E}_0$ as a pullback. Equation~\eqref{eq:locCart-bc} is then immediate.

Now we treat the case of $I_0 \simeq \{0\}$ is the initial element. Let $\ms{sd}_0$ be the full subcategory of $I^op$ consisting of chains that end with $0$. It suffices to prove
\begin{itemize}
\item Projection to $0$ induces an equivalence $\tx{Fun}_{\tx{loc.Cart}_{/I^\op}}(\ms{sd}_0, \mc{E}) \simeq \Phi(0)$; and
\item The relative left Kan extension $\tx{LKE}_i'$ as above is again well-defined, with $\ms{sd}_0$ replacing $\ms{sd}(I_0^\op)^\op$.
\end{itemize}
For the second claim, suppose $f: \ms{sd}_0 \to \mc{E}$ is an element of $\tx{Fun}_{\tx{loc.Cart}_{/I^\op}}(\ms{sd}_0, \mc{E})$, and let $\tx{LKE}_i(f)$ be the result of the Kan extension. If $P \to Q$ is an arrow in $\ms{sd}(I^\op)^\op$ locally Cartesian above $I^\op$ (meaning $P$ is a suffix of $Q$), then the map $\tx{LKE}_i(f)(P) \to \tx{LKE}_i(f)(Q)$ is either $f(P) \to f(Q)$ (if $0 \in Q$) or $f(P0) \to f(Q0)$ (where $P0$ is the result of appending $0$ to $P$); in either case, we see that $\tx{LKE}_i(f)$ preserve locally Cartesian arrows. It remains to prove the first claim.

To simplify notation we set $A := \ms{sd}_0$ and let $C$ be any locally Cartesian fibration over $I^\op$. Let $0$ be the initial object of $I$, we shall show that projection induces an equivalence $\tx{Fun}_{\tx{loc.Cart}_{/I^\op}}(A, C) \simeq C|_{0}$.	We let $A_m$ denote the full subcategory spanned by totally ordered subsets of $I^\op$ that are of size $\le m$ and end with $0$. We note that $A_1 \simeq \{0\}, A \simeq \bigcup_{m \ge 1} A_m$, and that $A_m$ is obtained from $A_{m - 1}$ by attaching a list of cones, one for each object in $A_{m} \setminus A_{m - 1}$. Relative right Kan extension then gives a sequence of adjunctions
	\[\ldots \adjoint \tx{Fun}_{/I^\op}(A_2, C) \adjoint \tx{Fun}_{/I^\op}(A_1, C) \simeq C|_0\]
	Noting that $\tx{Fun}_{/I^\op}(A, C) \simeq \limit_m \tx{Fun}_{/I^\op}(A_m, C)$, this induces a pair of adjunctions $\tx{Fun}_{/I^\op}(A, C) \adjoint C|_0$. It suffices to show that each pair of adjunction
	\[\tx{res}_{m + 1}: \tx{Fun}_{/I^\op}(A_{m + 1}, C) \simeq \tx{Fun}_{/I^\op}(A_{m}, C): \tx{RKE}_{m}\] restricts to an equivalence
	\[\tx{res}'_{m + 1}: \tx{Fun}_{\tx{loc.Cart}/_{I^\op}}(A_{m + 1}, C) \simeq \tx{Fun}_{\tx{loc.Cart}/_{I^\op}}(A_{m}, C): \tx{RKE}'_{m}\]
	from which the original claim is immediate. Since $\tx{res}_{m + 1}$ evident preserve the subcategories, we first want to show that if $f \in \tx{Fun}_{\tx{loc.Cart}/_{I^\op}}(A_{m}, C)$ preserve locally Cartesian arrows, then so does $\tx{RKE}_m(f)$.
	
	Suppose $\underline{a0} := (a_1 \ldots a_{m} 0)$ is a vertex in $A_{m + 1} \setminus A_m$, then out of it there is an unique locally Cartesian arrow $e$ mapping to $(a_1 \ldots a_m) \in A_m \setminus A_{m - 1}$. The value of $\tx{RKE}_m(f)$ at $\underline{a0}$ is the limit of all Cartesian lifts of $f(\underline{a'})$ to $\underline{a0}$, where $\underline{a'}$ is a substring of $a_1 \ldots a_m$. Among these, the Cartesian lift of $f(a_1 \ldots a_m)$ is the initial object, and it follows that the arrow lifting $e$ is locally Cartesian.
	The adjunction pair $(\tx{res}'_{m + 1}, \tx{RKE}_m')$ is then immediately an equivalence of categories, since the value of any $f \in \tx{Fun}_{\tx{loc.Cart}/_{I^\op}}(A_{m + 1}, C)$ at any $(a_1 \ldots a_m 0)$ is fixed as the Cartesian lift along the chain $(a_1 \ldots a_m 0) \to (a_1 \ldots a_{m - 1} 0) \to \ldots \to (0)$.
\end{proof}

For an arbitrary $i \in I$, by considering $\{i\} \subseteq I_{i/} \subseteq I$ we see that

\begin{corollary}
	For the same setup as above, for \emph{any} $i \in I$, $\tx{ev}_i$ admits a fully faithful adjoint $\tx{ins}_i$.
\end{corollary}

Now we consider the functorial properties of lax limits:

\begin{lemma}
\label{lemma:lax-limit-functorialities}
	For any category $A$ and any lax functor $\Phi: I \dashrightarrow \tx{Cat}$ as in Proposition \ref{prop:lax-limit-properties}, we have:
	\begin{enumerate}
	\item $\tx{Fun}(A, \ms{Glue}(\Phi)) \simeq \ms{Glue}(\tx{Fun}(A, \Phi(-)))$;
	\item $\tx{LFun}(A, \ms{Glue}(\Phi)) \simeq \ms{Glue}(\tx{LFun}(A, \Phi(-)))$;
	\end{enumerate}	
If all connecting functors $\Phi_{ij}$ admit right adjoints, then the above further identify with
	\begin{enumerate}
	\setcounter{enumi}{2}
	\item $\tx{RFun}(A, \ms{Glue}(\Phi)) \simeq \ms{opGlue}(\tx{RFun}(A, \Phi(-)))$, where RHS is the \emph{oplax} functor $I^\op \dashrightarrow \tx{Cat}$ defined by
\[i \mapsto \tx{RFun}(A, \Phi(i)) \hspace{1em} (i \to j) \in I \mapsto \tx{Run}(A, \Phi(j)) \xrightarrow{\Phi_{ij}^R \circ (-)} \tx{RFun}(A, \Phi(i));\]
	\end{enumerate}
	And if the lax functor factors through $I \dashrightarrow \tx{Pr}^L$ and $A$ is presentable, then the above further identify with
	\begin{enumerate}
	\setcounter{enumi}{3}
	\item $A \tensor \ms{Glue}(\Phi) \simeq \ms{Glue}(A \tensor \Phi(-))$.
	\end{enumerate}
\end{lemma}

\begin{proof}
	(1) is evident from Lemma~\ref{lemma:glue-as-strict-limit}, since each component of $\mf{G}(\Phi)$ is a functor category and we can swap $\tx{Fun}(A, -)$ to the inside. For (2), note that we have a natural functor from left to right using the natural map
	\[\tx{LFun}(A, \tx{Fun}(S, T)) \to \tx{Fun}(S, \tx{LFun}(A, T));\]
	then (2) follows from the fact that both embed fully faithfully into (1) with the same essential image as in Remark \ref{remark:lax-functor-fully-faithful}. With the same argument we identify LHS of (3) with the full subcategory of $\ms{Glue}(\tx{Fun}(A, \Phi(-)))$ where each component lands in $\tx{RFun}(A, \Phi(-))$. Now, if all arrows in a lax functor $\Psi: A \dashrightarrow \tx{Cat}$ admits right adjoints, then it is also a locally \emph{Cartesian} fibration, and the corresponding glue category is also the op-glue category for this oplax functor. Finally, (4) follows from (3) and \cite[4.8.1.17]{HA}, once we identify the connecting functors on both sides; this follows from unfolding the proof of \emph{loc.cit}.
\end{proof}

The following general observation tells us that lax limits are ubiquitous in application:

\begin{proposition}
\label{prop:glue-cat-reconstruction}
Let $C$ be a category containing all small colimits, $\{C_i\}_{i \in I_\tx{set}}$ a small set of cocomplete full subcategories. Let $F_i: C_i \to C$ denote the inclusion functors, and suppose that their right adjoints $\{F_i^R\}$ are continuous and jointly conservative. Let $I$ be a poset whose underlying set of vertices is $I_\tx{set}$, such that for every $i \not\to j$ in $I$, we have $F_j^R \circ F_i(\bullet) \simeq \tx{init}_{C_j}$. Then:
\begin{enumerate}
	\item The data $\Phi(i) := C_i, \Phi_{ij} := F_j^R \circ F_i$ extends to a lax functor $\Phi: I \dashrightarrow \tx{Cat}$; and
	\item There exists an equivalence $\varphi: C \simeq \ms{Glue}(\Phi)$.
\end{enumerate}
\end{proposition}

\begin{proof}
We first build up $\mc{E}_\Phi$ by declaring it to be the full subcategory $\{(c, i) \in C \times I^\op \mid c \in C_i\} \subseteq C \times I^\op$; it follows from unfolding definition that this is the desired lax functor.
Let $\tx{Const}_C$ denote the constant (strict) functor corresponding to the projection $C \times I^\op \to I^\op$. The fiberwise maps $F_i$ induce a oplax\footnote{Note that such data is a map $\mc{E}_\Phi \to C \times I^\op$ over $I^\op$ that does \emph{not} preserve locally Cartesian arrows.} natural transformation $\eta_F: \Phi \to \tx{Const}_C$, so $F_i^R$ induces an lax natural transformation $\eta_{F^R}: \tx{Const}_C \to \Phi$ (\cite[A.8.1]{ayala2021stratified}). Since constant diagram functor is the left adjoint to taking lax limit (\cite[A.5.8]{ayala2021stratified}), we obtain a canonical functor $\varphi : C \to \ms{Glue}(\Phi)$ from adjunction.
By definition, we have $F_i^R \simeq \tx{ev}_i \circ \varphi$.

We note that $\varphi$ admits a left adjoint $\varphi^L$, which will by construction satisfy $F_i \simeq \varphi^L \circ \tx{ins}_i$; we stipulate it to be
\[\tx{Fun}_{\tx{loc.Cart}_{/I^\op}}(\ms{sd}(I^\op)^\op, \mc{E}_\Phi) \xrightarrow{\oblv} \tx{Fun}_{/I^\op}(\ms{sd}(I^\op)^\op, \mc{E}_\Phi) \xrightarrow{\theta_F \circ (-)} \tx{Fun}_{/I^\op}(\ms{sd}(I^\op)^\op, C \times I^\op)\]
\[\simeq \tx{Fun}(\ms{sd}(I^\op)^\op, C) \xrightarrow{\colimit} C\]
where $\theta_F$ is the map $\mc{E}_\Phi \to C \times I^\op$ corresponding to $\eta_F$ described above. That this is the left adjoint of $\varphi$ can be proven similar to \cite[2.5.12]{ayala2021stratified}.
	
	We check $\tx{Id} \Rightarrow \varphi \circ \varphi^L$ is an equivalence. We note that $\ms{Glue}(\Phi)$ is generated by objects of the form $\tx{ins}_i(x_i)$: indeed, let $T = \prod_i \tx{ins}_i \circ \prod_i \tx{ev}_i$, then there is a natural map $\colimit(\ldots TTx \rightrightarrows T x) \to x$; since $\prod_i \tx{ev}_i$ is both conservative and continuous, this is an isomorphism. Thus it suffices to check $\tx{ev}_j \circ \tx{ins}_i \Rightarrow \tx{ev}_j \circ \varphi \circ \varphi^L \circ \tx{ins}_i \simeq F_j^R \circ F_i$ is an equivalence for each $i, j$, which follows from Proposition \ref{prop:lax-limit-properties}. We now check $\varphi^L \circ \varphi \Rightarrow \tx{Id}$ is an equivalence. Similar as above, it suffices to check $F_j^R \circ \varphi^L \circ \varphi \circ F_i \Rightarrow F_j^R \circ F_i$ is an equivalence for each $i, j$. Unfolding definition, we see that this functor can be computed as
	\[x_i \in C_i \mapsto F_j^R (\colimit_{k \in A} F_k \circ F_k^R \circ F_i(x_i)) \simeq \colimit_{k \in A} F_j^R F_k F_k^R F_i(x_i)\]
	By assumption, this simplifies to
	\[\colimit_{(i \to k \to j) \in A} F_j^R F_k F_k^R F_i(x_i) \simeq F_j^R F_j F_j^R F_i(x_i) \simeq F_j^R F_i(x_i)\]
	because $j$ is the final object among all such $k$.
\end{proof}

\section{Appendix B: Factorization Structure}

\label{sect:factspc}

In this section we give a self-contained exposition of what \emph{factorization structure} is. It essentially repeats \cite{raskin2015chiral} and should be skipped by readers who are already familiar with the objects involved. The only thing new here is Section \ref{sect:fact-restriction}, where we describe how to perform the (intuitively obvious, but technically annoying) operation of \emph{restriction of factorization modules}.

\begin{remark}
In \cite{raskin2015chiral}, where most of the foundation was laid down, Raskin used correspondences to encode factorization structures. Here we prefer the (equivalent, but uglier) fibrational approach as it is slightly more convenient for technical purposes.
\end{remark}

\subsection{Some Abstract Nonsense}

A few things before we start; readers should skip ahead.

\paragraph{Notation}
Given a category $\mc{C}$ that is coCartesian over an $\infty$-operad $\baseOp$, we'll call $\mc{C}_{\langle 1 \rangle}$, its fiber over $\langle 1 \rangle$, an $\baseOp$-monoidal category (and avoid calling $\mc{C}$ itself a $\baseOp$-monoidal category). When we have two $\baseOp$-monoidal categories $\mc{C}$ and $\mc{D}$, we'll refer to an element of $\tx{Alg}_{\baseOp}(\mc{C}^{\tensor}, \mc{D}^{\tensor})$ as a lax $\baseOp$-monoidal functor and one that preserves coCartesian arrows as a (strict) $\baseOp$-monoidal functor. The category of lax $\baseOp$-monoidal functor will be denoted $\tx{Fun}^{\tx{lax}, \baseOp}(\mc{C}, \mc{D})$, and the full subcategory of strictly $\baseOp$-monoidal functors will be denoted $\tx{Fun}^{\baseOp}(\mc{C}, \mc{D})$. There is also the dual notion of \emph{oplax} $\baseOp$-monoidal functors, which we shall denote by $\tx{Fun}^{\tx{oplax}, \baseOp}(\mc{C}, \mc{D})$.

\paragraph{Twisted Arrows}
For any \emph{ordinary} category $\mc{C}$, we can assign a corresponding category $\tx{Tw}(\mc{C})$ whose objects are arrows in $\tx{C}$ and whose 1-morphisms between $A \xrightarrow{p} B$ and $A' \xrightarrow{q} B'$ is an invertible 2-morphism in $\tx{C}$ of the following shape:
\[\xymatrix{
A \ar[d] \ar^p[r] & B \\
A' \ar^q[r] & B' \ar[u]
}\]
This notion admits an $\infty$-categorical extension; we refer to \cite[5.2.1]{HA} for the precise definition. Let $C, D$ be $(\infty, 1)$-categories where $D$ has fiber products. Let us call a functor from $\tx{Tw}(C)$ to $D$ \emph{composable} if, for every composable $X \xrightarrow{f} Y \xrightarrow{g} Z$ in $C$, the square
\[
\xymatrix{
g \circ f \ar[r] \ar[d] & g \ar[d] \\
f \ar[r] & \tx{id}_Y
}
\]
maps to a Cartesian square. The following is \cite[B.9]{raskin2015chiral}:
\begin{lemma}
\label{lemma:twisted-arrow-maps}
We have an equivalence of categories:
\[\{\tx{Composable functors from Tw(C) to D}\} \simeq \{\tx{Functors from C to }D^\tx{corr}\};\]
Moreover, this is compatible with the tensor structure on $1\tx{-Cat}$, so in particular we get another equivalence of categories
\[\{\tx{Composable symmetric monoidal functors}~\tx{Tw}(C) \to D\} \simeq \{\tx{Symmetric monoidal functors}~C \to D^\tx{corr}\}.\]
\end{lemma}

\paragraph{Monoidal Envelope}

For every $\infty$-operad $\baseOp$ and every $\baseOp$-monoidal\footnote{There is in fact no need for $\mc{C}$ to be $\baseOp$-monoidal: we just need $\mc{C}^{\tensor} \to \baseOp$ to be a fibration of operads. However we do not need this generality here.} category $\mc{C}$, there exists another $\baseOp$-monoidal category $\tx{Env}_{\baseOp}(\mc{C})$ with an universal (lax $\baseOp$-monoidal) map $\mc{C} \to \tx{Env}_\baseOp(\mc{C})$ such that it induces an equivalence
\[\tx{Fun}^{\tx{lax}, \baseOp}(\mc{C}, \mc{D}) \simeq \tx{Fun}^{\baseOp}(\tx{Env}_{\baseOp}(\mc{C}), \mc{D})\]
for every $\baseOp$-monoidal category $\mc{D}$. We also have the dual notion of a \emph{comonoidal} envelope.

\begin{remark}
This notion, as well as many others below, makes sense for both unital and non-unital operads. We do not dwell on this distinction too much, as it should always be clear from context which one is needed.
\end{remark}

\begin{definition}
$\fSetnu$ is the nerve of the $(1, 1)$-category consisting of nonempty finite sets and surjections among them.
\end{definition}

\begin{example}
	The \emph{symmetric monoidal} envelope of the symmetric monoidal singleton category is the symmetric monoidal category $\fSetnu$. The symmetric monoidal envelope of $\fSetnu$ is the twisted arrow category $\tx{Tw}(\fSetnu)^\op$.
\end{example}

\begin{remark}
\label{rmk:tw-composability-weak}
	It follows that a \emph{commutative algebra object} in $\mc{C}^\tx{corr}$ (which is a lax symmetric monoidal functor $* \to \mc{C}^\tx{corr}$) is also a symmetric monoidal functor $\fSetnu \to \mc{C}^\tx{corr}$, thus a composable symmetric monoidal functor $\tx{Tw}(\fSetnu) \to \mc{C}$, which is in turn an \emph{oplax} functor $\fSetnu \to \mc{C}$ satisfying some composability condition.
\end{remark}

\paragraph{Monoidal Kan Extensions}
We will need the following version of monoidal Kan extension, whose proof follows exactly that of \cite[Lemma 2.16]{ayala2016factorization}:
\begin{lemma}
\label{lemma:sm-right-kan-extension}
Fix $\baseOp$ an $\infty$-operad. Let $G: \mc{C} \to \mc{D}$ be an $\baseOp$-monoidal functor between $\baseOp$-monoidal categories such that $\mc{C}$ is small and $\mc{D}$ locally small. Suppose we have another $\baseOp$-monoidal category $\mc{V}$, and suppose there exists a subclass $\mc{K}$ of $\tx{sSet}$ such that:
	\begin{enumerate}
	\item $\mc{V}$ admits $\mc{K}$-limits;
	\item each slice category $\mc{C}_{d/}$ belongs to $\mc{K}$;
	\end{enumerate}
	then the right Kan extension $\tx{RKE}_{F}: \tx{Fun}(\mc{C}, \mc{V}) \to \tx{Fun}(\mc{D}, \mc{V})$ exists. Furthermore, if
	\begin{enumerate}[resume]
	\item for each active arrow $p: x \to y$ in $\baseOp$ and every $d \in \mc{D}^{(x)}$, the functor
	\[\mc{C}^{(x)}_{d/} \to \mc{C}^{(y)}_{p(d)/}\]
	is initial, where $d \to p(d)$ is the coCartesian lift of $p$,
	\end{enumerate}
	then the forgetful functors between the category of lax $\baseOp$-monoidal functors
	\[\tx{Fun}^{\tx{lax}, \baseOp}(\mc{D}, \mc{V}) \to \tx{Fun}^{\tx{lax}, \baseOp}(\mc{C}, \mc{V})\]
	admits a right adjoint and is computed as right Kan extension of the underlying functor along $G$; finally, if further
	\begin{enumerate}[resume]
	\item the $\baseOp$-monoidal structure of $\mc{V}$ distributes over $\mc{K}$-limits;
	\end{enumerate}
	then the right Kan extension also provides a right adjoint to the forgetful functors between the category of $\baseOp$-monoidal functors
	\[\tx{Fun}^{\baseOp}(\mc{D}, \mc{V}) \to \tx{Fun}^{\baseOp}(\mc{C}, \mc{V}).\]	
An analogous statement holds true for left Kan extension as well.
\end{lemma}

\begin{example}
In the case of $\baseOp = \CommOperad$, the third condition is equivalent to the following:
\begin{itemize}
\item The tensor product functor $\mc{C}_{d_1 / } \times \mc{C}_{d_2 / } \to \mc{C}_{d_1 \tensor d_2 /}$ is initial; and
\item The functor $\mc{C}_{\tb{1}_\mc{C}/} \to \mc{C}_{\tb{1}_\mc{D}/}$ is initial.
\end{itemize}
\end{example}

\subsection{Non-unital Factorization Spaces}

For any given prestack $\mc{X}$, let $\Ran_{\tx{nu}}(\mc{X})$ denote the prestack associated with $\mc{X}$, defined as
\[\tx{Ran}_{\tx{nu}, \mc{X}} := \colimit_{I \in (\fSetnu)^\tx{op}} \mc{X}^I\]

We let $[\mc{X} \times \mc{X}]_{\tx{disj}}$ to be the subfunctor whose $S$-point is given by
\[\{x: S \to X^2, \alpha: S \times_{X^2} X \simeq \emptyset\};\]
then we let $[\mc{X}^I \times \mc{X}^J]_{\tx{disj}}$ denote the subfunctor consisting of those points of $(\mc{X}^I \times \mc{X}^J)(S)$ such that for every $i \in I, j \in J$, the corresponding $S$-point of $\mc{X}^2$ lies in $[\mc{X} \times \mc{X}]_{\tx{disj}}$. Finally, we set
\[[\Ran_{\nonun, \mc{X}} \times \Ran_{\nonun, \mc{X}}]_{\tx{disj}} := \colimit_{I, J \in (\fSetnu)^\op} [\mc{X}^I \times \mc{X}^J]_{\tx{disj}}.\]

\begin{convention}
We use $\Ran_\tx{nu}$ to denote $\Ran_{\tx{nu}, X}$ (where $X$ is our global curve) and $\Ran_{\tx{nu},\dR}$ for $\Ran_{\tx{nu}, X_\dR}$. When we use the latter, we refer to relevant objects as \emph{crystals} of factorizable objects.
\end{convention}

\begin{definition}
Let
\[\tx{WeakFactMonoid}_\nonun := \tx{Fun}^{\tx{oplax}, \tensor}_{\nonun}(\fSetnu, \tx{PreStk}).\]
An element of this category is called a non-unital \emph{weak factorization monoids}.
\end{definition}

\begin{remark}
For $\mc{C}$ an arbitrary $(\infty, 1)$-category, there is the $(\infty, 1)$-category of correspondences in $\mc{C}$, defined in \cite{raskin2015chiral}; we denote it by $\mc{C}_\tx{corr}$. By Remark \ref{rmk:tw-composability-weak}, a non-unital commutative algebra within $\tx{PreStk}_\tx{corr}$ is a weak factorization space satisfying an extra \emph{composability} condition. We will refer to such objects as \emph{factorization monoids}.
\end{remark}

\begin{convention}
	When there is no ambiguity on the factorization structure, we will sometimes abuse notation and use the prestack $F([1])$ to refer to the entire weak factorization monoid.
\end{convention}

There are two non-unital factorization monoids whose underlying prestack is $\tx{Ran}_\nonun$. First the standard one, given by
\[\Ran_\nonun(I) := [\underbrace{\Ran_\tx{nu} \times \ldots \times \Ran_\tx{nu}}_{|I|~\tx{times}}]_\tx{disj};\]
	the oplax structure is the obvious inclusion $\iota$, and the map corresponding to $I \surj J$ is the ``add'' map. In terms of correspondences, it can be depicted by the following diagram:
\[
\xymatrix{
 & [\Ran_{\tx{nu}} \times \Ran_{\tx{nu}}]_{\tx{disj}} \ar_{\iota}[ld] \ar^{\ms{add}}[rd] \\
 \Ran_{\tx{nu}} \times \Ran_{\tx{nu}} & & \Ran_{\tx{nu}}
}
\]
There is also $\tx{Ran}_{\tx{nu}}^{\tx{comm}}$, the \emph{commutative} Ran space, given by
\[\Ran_\nonun^\tx{Comm}(I) := \Ran_\tx{nu}^{\times |I|}.\]
Graphically it is depicted as
\[
\xymatrix{
 & \Ran_{\tx{nu}} \times \Ran_{\tx{nu}} \ar_{\tx{id}}[ld] \ar^{\ms{add}}[rd] \\
 \Ran_{\tx{nu}} \times \Ran_{\tx{nu}} & & \Ran_{\tx{nu}}
}
\]

\begin{definition}
\label{def:factspace-nu}
	We let $\tx{FactSpace}_\nonun$ denote the full subcategory of
	\[\tx{Fun}^{\tx{oplax}, \tensor}_\nonun(\fSetnu, \tx{PreStk})_{/\Ran_\nonun}\]
	consisting of oplax symmetric monoidal functors $F$ over $\Ran_\nonun$ satisfying the following conditions:
	\begin{itemize}
	\item[(Condition R)] Natural transformations as plain functors are Cartesian, i.e. we have the following Cartesian squares for all $I \surj J$:
	\[\xymatrix{
	F(I) \ar[r] \ar[d] & F(J) \ar[d] \\
	\Ran_\nonun(I) \ar[r] & \Ran_\nonun(J)
	}\]
	\item[(Condition L)] Natural transformations from the oplax structure are Cartesian, i.e. we have the following Cartesian squares for all $I \surj J$:
	\[\xymatrix{
	F(I) \ar[r] \ar[d] & \prod_j F(I_j) \ar[d] \\
	\Ran_\nonun(I) \ar[r] & \prod_j \Ran_\nonun(I_j)
	}\]
	\end{itemize}
	An element of this category is referred to as a non-unital factorization space.
\end{definition}
At the first level, condition R (resp. L) means that the right (resp. left) square of the following diagram is Cartesian:
\[
\xymatrix{
		& F(2) \ar[d] \ar[ld] \ar[rd] \\
	F(1) \times F(1) \ar_{}[d] & [\Ran_\nonun \times \Ran_\nonun]_\tx{disj} \ar[ld] \ar[rd] & F(1) \ar^{}[d] \\
	\Ran_\nonun \times \Ran_\nonun & & \Ran_\nonun
}
\]

\subsubsection{Partition Description}
Since the Ran space is a colimit over different $X^I$, we sometimes find it useful to describe factorization data as a collection of data living over $X^I$ satisfying compatibilities. Here we follow the construction of \cite[Section 8]{raskin2015chiral}.

\begin{notation}
We set $\ECOp := \tx{Tw}(\fSetnu)^\op$ and $\ETOp := \tx{Env}(\ECOp)$.
\end{notation}

\begin{definition}
\label{def:u-ran}
Introduce $U_\Ran \in \tx{Fun}^{\tx{oplax}, \tensor}_\nonun(\ECOp^\op, \tx{PreStk})$ such that
\[U_\Ran(p: I \surj J)(S) := \{(x_i)_{i \in I} \in X^I(S) \mid p(i_1) \neq p(i_2) \implies x_{i_1} \neq x_{i_2}~\forall i_1, i_2 \in I\}\]
This is the ``partition'' incarnation of the Ran space; it is again an oplax symmetric monoidal functor.
\end{definition}
Let $\pi: \ECOp^\op \to \fSetnu, (I \surj J) \mapsto J$ denote the projection. Operadic left Kan extension (Lemma \ref{lemma:sm-right-kan-extension}) preserves oplax symmetric monoidal functors, and we observe that $\Ran_\nonun \simeq \tx{LKE}_\pi(U_\Ran)$ because $\pi$ is a coCartesian fibration.

\begin{definition}
\label{def:factspace-part-nu}
We let $\tx{FactSpace}_{\nonun, \PartitionMark}$ denote the full subcategory of $\tx{Fun}^{\tx{oplax}, \tensor}_\nonun(\ECOp^{\tx{op}}, \tx{PreStk})_{/U_\Ran}$ consisting of functors $F$ over $U_\Ran$ such that:
\begin{itemize}
\item[(Condition R')] Natural transformations as plain functors are Cartesian, i.e. we have the following Cartesian squares for all $(I' \to J') \to (I \to J) \in \ECOp$:
	\[\xymatrix{
	F(I \surj J) \ar[r] \ar[d] & F(I' \surj J') \ar[d] \\
	U_\Ran(I \surj J) \ar[r] & U_\Ran(I' \surj J')	
	}\]
\item[(Condition L')] Natural transformations from the oplax structure are Cartesian, i.e. we have the following Cartesian squares for all $I \surj J \surj K$:
	\[\xymatrix{
	F(I \surj J) \ar[r] \ar[d] & \prod_{k} F(I_k \surj J_k) \ar[d] \\
	U_\Ran(I \surj J) \ar[r] & \prod_{k} U_\Ran(I_k \surj J_k)
	}\]
\end{itemize}
\end{definition}

The following fact is parallel to \cite[Section 8]{raskin2015chiral}

\begin{proposition}
\label{prop:fact-spc-formulations-eqv}
$\tx{FactSpace}_\nonun \simeq \tx{FactSpace}_{\nonun, \PartitionMark}$.
\end{proposition}

Modulo a quick combinatorial check, this reduces to a general claim about left Kan extensions described below, to which we devote the rest of this subsection. Since the proof is pure combinatorics, we recommend skipping it on the first pass.

\begin{definition}
	For a cospan $I \surj^{p} J' \twoheadleftarrow^{q} J$ in $\fSetnu$, define the (ordinary) category $\tx{Factor}(p, q)$:
	\begin{itemize}
	\item Its objects are surjections $I \surj^r J$ such that $q \circ r = p$; and
	\item A map $r \to r'$ is a bijection $a: I \simeq I$ such that $r \circ a = r'$, $p \circ a = p$.
	\end{itemize}
	when the relevant maps are clear sometimes we just write $\tx{Factor}(I \surj J' \twoheadleftarrow J)$.
\end{definition}

\begin{proposition}
\label{prop:fact-formulations-equiv}
	Suppose $\mc{C}$ is a symmetric monoidal $\infty$-topos, such that the monoidal product commutes with filtered colimits in each variable and distributes over pullbacks. Let $U: \ECOp^\op \to \mc{C}$ be a oplax symmetric monoidal functor, and $\tx{LKE}_\pi(U): \fSetnu \to \mc{C}$ be its (operadic) left Kan extension. Suppose that, for every diagram
	\[\xymatrix{ & & E \ar[d] \\
	A \ar[r] & B \ar[r] & C \ar[r] & D}\]
	in $\fSetnu$, the outer square (equivalently, both squares) of the diagram
\[\xymatrix{
\colimit_{\tx{Factor}(B \surj C \twoheadleftarrow E)} U(B \surj E) \ar[r] \ar[d] & \colimit_{\tx{Factor}(A \surj C \twoheadleftarrow E)} U(A \surj E) \ar[d] \\
U(B \surj C) \ar[r] \ar[d] & U(A \surj C) \ar[d] \\
\bigotimes_{d \in D} U(B_d \surj C_d) \ar[r] & \bigotimes_{d \in D} U(A_d \surj C_d)
}\]
is Cartesian. Then operadic Kan extension
	\[\tx{LKE}_\pi^{/U}: \tx{Fun}^{\tx{oplax}, \tensor}_\nonun(\ECOp^\op, \mc{C})_{/U} \to \tx{Fun}^{\tx{oplax}, \tensor}_\nonun(\fSetnu, \mc{C})_{/\tx{LKE}_\pi(U)}\]
	induces an equivalence between the full subcategory of LHS satisfying condition (R') and (L') as in Definition \ref{def:factspace-part-nu} and the full subcategory of RHS satisfying condition (R) and (L) as in Definition \ref{def:factspace-nu}.
\end{proposition}

\begin{proof}
Note that $\ECOp$ is the comonoidal envelope of $\fSetnu$ and $\ETOp^\op$ is the comonoidal envelope of $\ECOp^\op$. Then there is a commutative diagram
\[\xymatrix{
\tx{Fun}^{\tx{oplax}, \tensor}_\nonun(\ECOp^\op, \mc{C}) \ar^{\tx{LKE}_\pi}[r] \ar^{F \mapsto F_\tx{univ}}_{\simeq}[d] & \tx{Fun}^{\tx{oplax}, \tensor}_\nonun(\fSetnu, \mc{C}) \ar^{F \mapsto F_\tx{univ}}_{\simeq}[d] \ar@<5pt>^{\tx{Res}}[l] \\
\tx{Fun}^{\tensor}(\ETOp^\op, \mc{C}) \ar^{\tx{LKE}_\tau}[r] & \tx{Fun}^{\tensor}(\ECOp, \mc{C}) \ar@<5pt>^{\tx{Res}}[l]
}\]
where $\tau(I \surj J \surj K) = (J \surj K)$; note that $\tx{LKE}_\tau$ is again just computed by LKE on the underlying functor. Let
\[\mc{L} := \tx{LKE}_\tau^{/U_\tx{univ}}: \tx{Fun}^{\tensor}(\ETOp^\op, \mc{C})_{/U_\tx{univ}} \to \tx{Fun}^{\tensor}(\ECOp, \mc{C})_{/\tx{LKE}_\tau(U_{\tx{univ}})}\]
be the induced functor; its right adjoint $\mc{R}$ is then given by sending $F_\tx{univ}$ to the fiber product
\[U_\tx{univ} \times_{\tx{Res} \circ \tx{LKE}_\tau(U_\tx{univ})}\tx{Res}(F_\tx{univ})\]
within the category $\tx{Fun}^{\tensor}(\ETOp^\op, \mc{C})$. Since the category of symmetric monoidal functors embed fully faithfully into the category of \emph{lax} symmetric monoidal functor, i.e. commutative algebra objects in the Day convolution monoidal structure, and because tensor product distributes over pullbacks by assumption, this pullback can be computed in the underlying category
$\tx{Fun}(\ETOp^\op, \mc{C})$.

For brevity, set $R := \tx{LKE}_\pi(U)$, so that $R_\tx{univ} \simeq \tx{LKE}_\tau(U_\tx{univ})$. We let

\begin{itemize}
\item $\mathscr{L}$ denote the full subcategory of $\tx{Fun}^{\tensor}(\ETOp^\op, \mc{C})_{/U_\tx{univ}}$ consisting of those object $G_\tx{univ}$ such that, for every $(I \surj J \surj K) \to (I' \surj J' \surj K')$ in $\ETOp^\op$, the following square is Cartesian:
	\[
\xymatrix{
G_\tx{univ}(I \surj J \surj K) \ar[r] \ar[d] & G_\tx{univ}(I' \surj J' \surj K') \ar[d] \\
\bigotimes_{k} U(I_k \surj J_k) \simeq U_\tx{univ}(I \surj J \surj K) \ar[r] & U_\tx{univ}(I' \surj J' \surj K') \simeq \bigotimes_{k'} U(I'_{k'} \surj J'_{k'})
}	
	\]
\item $\mathscr{R}$ denote the full subcategory of $\tx{Fun}^{\tensor}(\ECOp, \mc{C})_{/R_\tx{univ}}$ consisting of those object $F_\tx{univ}$ such that, for every $(J \surj K) \to (J' \surj K')$ in $\ECOp$, the following square is Cartesian:
	\[
\xymatrix{
F_\tx{univ}(J \surj K) \ar[r] \ar[d] & F_\tx{univ}(J' \surj K') \ar[d] \\
\bigotimes_{k} R(J_k) \simeq R_\tx{univ}(J \surj K) \ar[r] & R_\tx{univ}(J' \surj K') \simeq \bigotimes_{k'} R(J'_{k'})
}	
	\]
\end{itemize}
It is immediate that the two subcategories originally considered correspond to $\mathscr{L}$ and $\mathscr{R}$ respectively. Our job then is to check that $\tx{Res}$ and $\tx{LKE}_\tau$ both preserve these subcategories, and on these subcategories the unit/counit induce isomorphisms on the level of underlying functors.

If $F_\tx{univ} \in \mathscr{R}$ then $\mc{R}(F_\tx{univ}) \in \mathscr{L}$ by unfolding definition. Conversely, for $G_\tx{univ} \in \mathscr{L}$, we have
\[\tx{LKE}_\tau^{/U_\tx{univ}}G_\tx{univ}(J \surj K) \simeq \colimit_{\substack{\xymatrix@=0.5em{I' \ar[r] & J' \ar[d] \ar[r] & K' \\
 & J \ar[r] & K \ar[u]}}} G_\tx{univ}(I' \surj J' \surj K') \simeq \colimit_{I \surj J} G_\tx{univ}(I \surj J \surj K)\]
by cofinality, so our job is to prove, for any fixed map $(J \surj K) \to (J' \surj K') \in \ECOp$, the diagram
 \[\xymatrix{
 \colimit_{(I \surj J) \in \mc{B}} G_\tx{univ}(I \surj J \surj K) \ar[r] \ar[d] & \colimit_{(I' \surj J') \in \mc{A}} G_\tx{univ}(I' \surj J' \surj K') \ar[d] \\
 \colimit_{(I \surj J) \in \mc{B}} U_\tx{univ}(I \surj J \surj K) \ar[r] & \colimit_{(I' \surj J') \in \mc{A}} U_\tx{univ}(I' \surj J' \surj K')
 }
 \]
is Cartesian, where $\mc{B} := (\fSetnu_{/J})^\op$ and $\mc{A} := (\fSetnu_{/J'})^\op$. Note that $\mc{B} \to \mc{A}$ is a coCartesian fibration and the corresponding functor sends $L \surj J'$ to $\tx{Factor}(L \surj J' \twoheadleftarrow J)$, so we have the diagram
\[\xymatrix@=1em{\bullet \ar[r] \ar[d] & \colimit_{a \in \mc{A}} G_\tx{univ}(a  \surj J' \surj K') \ar[d] \\
\colimit_{a \in \mc{A}} \colimit_{b \in \mc{B}_a} U_\tx{univ}(b \surj J \surj K) \ar[r] & \colimit_{a \in \mc{A}} U_\tx{univ}(a \surj J' \surj K')}\]
We claim that this is the same as the colimit of pullbacks
\[\colimit_{a \in \mc{A}} \xymatrix@=0.7em{\bullet \ar[r] \ar[d] & G_\tx{univ}(a  \surj J' \surj K') \ar[d] \\
\colimit_{b \in \mc{B}_a} U_\tx{univ}(b \surj J \surj K) \ar[r] & U_\tx{univ}(a  \surj J' \surj K')}\]
Assuming this, since colimits are universal, we can write the pullback as
\[\colimit_{a \in \mc{A}} \colimit_{b \in \mc{B}_a} \xymatrix@=0.7em{\bullet \ar[r] \ar[d] & G_\tx{univ}(a \surj J' \surj K') \ar[d] \\ U_\tx{univ}(b \surj J \surj K) \ar[r] & U_\tx{univ}(a \surj J' \surj K')}\]
\[\simeq \colimit_{a \in \mc{A}} \colimit_{b \in \mc{B}_a} G_\tx{univ}(b \surj J \surj K) \simeq \colimit_{\mc{B}} G_\tx{univ}\]
as desired. To check the claim, rewrite
\[x_a = G_\tx{univ}(a \surj J' \surj K'), y_a = U_\tx{univ}(a \surj J' \surj K'), z_a = \colimit_{b \in \mc{B}_a} U_\tx{univ}(b \surj J \surj K);\]
\[x = \colimit_a x_a, y = \colimit_a y_a, z = \colimit_a z_a.\]
The claim is that $\colimit_{a} x_a \times_{y_a} z_a \simeq x \times_y z$. Since colimits are universal, we always have $x \times_y z \simeq \colimit_a x_a \times_y z \simeq \colimit_a x_a \times_{y_a} y_a \times_y z$, so we only need $y_a \times_y z \simeq z_a$, i.e. we need for every fixed $I' \surj J'$, the diagram
\[\xymatrix{
\colimit_{\tx{Factor}(I' \surj J' \twoheadleftarrow J)} U_\tx{univ}(I' \surj J \surj K) \ar[r] \ar[d] & \tx{LKE}_\tau U_\tx{univ}(J \surj K) \ar[d] \\
\bigotimes_{k' \in K'} U(I'_{k'} \surj J'_{k'}) \simeq U_\tx{univ}(I' \surj J' \surj K') \ar[r] & \tx{LKE}_\tau U_\tx{univ}(J' \surj K') \simeq \bigotimes_{k' \in K'} \tx{LKE}_\pi U(J'_{k'})
}\]
One further reduces to $|K| = 1$; this then follows from descent and our assumption on $U$.

Universality of colimits guarantees $\mc{L} \circ \mc{R} \simeq \tx{id}$; to check $\tx{id} \simeq \mc{R} \circ \mc{L}$ on $\mathscr{L}$, note that we have a pullback square
 \[\xymatrix{\mc{R} \circ \mc{L}(G_\tx{univ})(I \surj J \surj K) \ar[r] \ar[d] & \colimit_{I' \surj J} G_\tx{univ}(I' \surj J \surj K) \simeq \tx{Res} \circ \tx{LKE}_\tau(G_\tx{univ})(I \surj J \surj K) \ar[d] \\
 U_\tx{univ}(I \surj J \surj K) \ar[r] & \colimit_{I' \surj J} U_\tx{univ}(I' \surj J \surj K) \simeq \tx{Res} \circ \tx{LKE}_\tau(U_\tx{univ})(I \surj J \surj K)}\]
 Descent and the condition on $\mathscr{L}$ then guarantees that the pullback is isomorphic to $G_\tx{univ}(I \surj J \surj K)$.
\end{proof}

\subsection{Unital Factorization Spaces}
\label{sect:unital-fact-spaces}
\newcommand{\unaug}{\tx{un,aug}}
Now we introduce the unital Ran space $\Ran_\un$.
\begin{definition}
	The \emph{unital} Ran space $\Ran_\un$ is the lax prestack whose $S$-point is the poset category of (potentially empty) finite subsets of $X(S)$ with the inclusion relation.
\end{definition}
As we can write $\Ran_\nonun = \colimit_{I \in \fSetnu^\tx{op}} X^I$, we can also write the unital Ran as a \emph{partially oplax colimit}.

\begin{definition}
Let $F: C \to 1\tx{-Cat}$ be a functor, and $I$ be a class of morphisms in $C$. Let $\tx{Cart}(F) \to C^\op$ be the corresponding Cartesian fibration. The partially oplax colimit, which we denote by $\tx{oplax} \colimit_{C, I} F$, is the localization of $\tx{Cart}(F)$ with respect to all Cartesian arrows that lies above arrows in $I$.
\end{definition}

\begin{definition}
Let $\Setli$ denote the category of all (potentially empty) finite sets and arbitrary morphisms among them; it contains $\fSetnu$ as a 1-full subcategory.
\end{definition}
Let $X^{(-)}$ denote the functor $\Setli^\op \to \tx{PreStk}$; it follows from unpacking definition (c.f. \cite{gaitsgory2019atiyahbott}) that
\[\Ran_\un \simeq \tx{oplax} \colimit_{\Setli^\op, \fSetnu^\tx{op}} X^{(-)}.\]
\begin{definition}
An \emph{unital} weak factorization monoid is an \emph{oplax} \emph{oplax-unital}\footnote{This means that for such a monoid $F$, $F(0)$ is \emph{not} required to be $\tx{pt}$.} symmetric monoidal functor from $\Setli$ to $\tx{LaxPreStk}$; the category is denoted as
\[\tx{WeakFactMonoid}_\un := \tx{Fun}^{\tx{oplax}, \tensor}(\Setli, \tx{LaxPreStk}).\]
\end{definition}
We have the following correspondence diagram:
\[\xymatrix{
 & [\Ran_\un \times \Ran_\un]_\tx{disj} \ar[ld] \ar[rd] \\
 \Ran_\un \times \Ran_\un & & \Ran_\un
} \hspace{1em} \xymatrix{
 & \tx{pt} \ar[ld] \ar[rd] \\
 \tx{pt} & & \Ran_\un
}\]
where $[\Ran_\un \times \Ran_\un]_\tx{disj}$ is the lax sub-prestack of $\Ran_\un \times \Ran_\un$ consisting of objects that lies in $[\Ran_\nonun \times \Ran_\nonun]_\tx{disj}$. It gives rise to a unital factorization monoid structure on $\Ran_\un$.

\begin{definition}
	We let $\tx{FactSpace}_\un$ denote the 1-full subcategory of
	\[\tx{Fun}^{\tx{oplax}, \tensor}(\Setli, \tx{LaxPreStk})_{/\Ran_\un}\]
	consisting of unital factorization monoids $F$ over $\Ran_\un$ satisfying the following conditions:
	\begin{itemize}
	\item[(Condition R)] Natural transformations as plain functors are Cartesian for $I \surj J$ surjective:
	\[\xymatrix{
	F(I) \ar[r] \ar[d] & F(J) \ar[d] \\
	\Ran_\un(I) \ar^{\tx{add}}[r] & \Ran_\un(J)
	}\]
	\item[(Condition L)] Natural transformations from the oplax structure are Cartesian for $I \surj J$ surjective:
	\[\xymatrix{
	F(I) \ar[r] \ar[d] & \prod_j F(I_j) \ar[d] \\
	\Ran_\un(I) \ar[r] & \prod_j \Ran_\un(I_j)
	}\]
	\item[(Condition F)] For each $I \in \fSetnu$, the map $F(I) \to \Ran_\un(I)$ is a left fibration\footnote{In the model-independent world, ``left fibration'' means ``coCartesian fibration fibered in spaces''.} (i.e. it is pointwise a left fibration in a functorial way);
	\item[(Condition P)] A map $F \Rightarrow G$ between two elements is required to pointwise preserve coCartesian arrows;
	\item[(Condition C)] For every composable pair of morphisms $I \to J \to K$ in $\Setli$, the commutative square
\[\xymatrix{
	\prod_{k = 1}^{K} F(I_k) \ar[r] \ar[d] & \prod_{k = 1}^{K} F(J_k) \ar[d] \\
	\prod_{j = 1}^{J} F(I_j) \ar[r] & \prod_{j = 1}^{J} F([1])
}\]
	is Cartesian.
	\end{itemize}
	A unital factorization space is then an element of $\tx{FactSpace}_\un$.
\end{definition}

\begin{remark}
	Condition C corresponds to composability in Lemma \ref{lemma:twisted-arrow-maps}; it was redundant in the non-unital case because it is implied by condition R and L, and the fact that $\Ran_\nonun$ itself satisfies condition C.
\end{remark}

\begin{variant}
We let $\tx{FactSpace}_{\un, \WeakMark}$ denote the category sharing the same definition as above but with condition (R) dropped. An element of this category is referred to as a \emph{weak} unital factorization space.
\end{variant}

\begin{variant}
If we drop condition P above, we obtain the notion of a \emph{lax-unital} morphism between factorization spaces. The category of factorization spaces and lax-unital morphisms between them is denoted $[\tx{FactSpace}_{\un}]_{\tx{lax-unital}}$.
\end{variant}

\begin{variant}
If we replace ``left'' in condition F by ``right'' and ``coCartesian'' in condition P by ``Cartesian'', we obtain the notion of a \emph{co-unital} factorization space (the category of which is denoted $\tx{FactSpace}_{\tx{co-un}}$). The lax-counital version is denoted $[\tx{FactSpace}_{\tx{co-un}}]_{ \tx{lax-counital}}$.
\end{variant}

Let us explain how this definition compare with the existing definitions in literature.
For every map $\varphi: S \to \Ran_\un$ from an arbitrary prestack $S$, let $C_\varphi$ denote the fiber product lax prestack; since condition F requires the fibration to be fibered in spaces, $C_\varphi$ will be an actual prestack.
The coCartesian part of condition F implies that, for every natural transformation $\varphi \Rightarrow \psi: S \to \Ran_\un$, we have a map $C_\varphi \to C_\psi$ as prestacks over $S$. Condition P further guarantees that this construction is functorial in $S$, and thus combines into a map
\[\tx{Fun}([1], \Ran_\un) \times_{\Ran_\un, \tx{source}} F(\tb{1}) \to F(\tb{1})\]
which is the structured described in e.g. \cite[Section 4]{gaitsgory2017semi}.
In particular, for an injective map $I \inj J$, the natural transformation
\[(X^J \to X^I \to \Ran \to \Ran_\un) \Rightarrow (X^J \to \Ran \to \Ran_\un)\]
gives a morphism
\[X^J \times_{X^I} (F(\tb{1}) \times_{\Ran_\un} X^I) \to (F(\tb{1}) \times_{\Ran_\un} X^J)\]
which is the usual unital data.

The following definition is reminiscent of the ``relative'' $\infty$-category theory developed in \cite{shah2017parametrized}:

\begin{definition}
Let $C_1 \xrightarrow{\varphi} C_2 \xrightarrow{\pi} D$ be maps between simplicial sets, such that $\pi \circ \varphi$ and $\pi$ are both flat inner fibrations. We say $\varphi$ has $D$-left anodyne lifting property if, for every map $A \xrightarrow{i} B \xrightarrow{p} D$ where $i$ is left anodyne, the induced map
\[(C_1^{B} \times_{D^{B}} \{p\}) \to (C_1^A \times_{D^A} \{p \circ i\}) \times_{(C_2^A \times_{D^A} \{p \circ i\})} (C_2^B \times_{D^B} \{p\})\]
is a trivial Kan fibration. The notion of $D$-right anodyne lifting property is defined analogously.
\end{definition}

\begin{example}
By \cite[2.1.2.9]{HTT} if $\varphi$ is itself a left/right fibration, then it has $D$-left/right anodyne lifting property for any $D$.
\end{example}

\begin{variant}
We can further weaken condition F to being a \emph{flat inner fibration}\footnote{The model-independent notion is ``exponentiable fibration'' c.f. \cite{ayala2020fibrations}.} fibered in spaces; we choose to call the resulting structure a \emph{corr-unital} factorization space.\footnote{This structure has also been referred to simply as \emph{unital} in literature, which to us seem like an unfortunate naming clash.} In the corr-unital case, we have four types of condition P:
\begin{enumerate}
\item The default notion is when we drop condition P altogether. The resulting category will be denoted $\tx{FactSpace}_{\tx{corr-un}}$;
\item We can restrict morphism spaces to those maps $\varphi: F \Rightarrow G$ such that, for every $I$ and $S$, $\varphi(I)(S)$ has the $\Ran_\un(I)(S)$-left anodyne lifting property. The resulting category will be called the category of corr-unital factorization spaces \emph{with strictly unital maps} and will be denoted $[\tx{FactSpace}_{\tx{corr-un}}]_{\tx{str-unital}}$;
\item We can replace ``left'' with ``right'' above and obtain a dual version of corr-unital factorization spaces \emph{with strictly counital maps}, which we denote by $[\tx{FactSpace}_{\tx{corr-un}}]_{\tx{str-counital}}$;
\item Finally we have the category of corr-unital factorization spaces \emph{with strictly bi-unital maps} where we require both lifting properties at the same time, defined as $[\tx{FactSpace}_{\tx{corr-un}}]_{\tx{str-biunital}}$.
\end{enumerate}
\end{variant}

In the corr-unital case, recall that to each map $\mc{Y} \to \mc{Z}$ of lax prestacks we can form a corresponding map of arrow prestack $\tx{Fun}([1], \mc{Y}) \to \tx{Fun}([1], \mc{Z})$. A natural transformation $\varphi \Rightarrow \psi$ gives a map $S \to \tx{Fun}([1], \Ran_\un)$, so we can form
\[\xymatrix{
 & S \times_{\tx{Fun}([1], \Ran_\un)} \tx{Fun}([1], F(\tb{1})) \ar^{\tx{proj}_0}[ld] \ar_{\tx{proj}_1}[rd] \\
 C_\varphi & & C_\psi}\]
which is the usual correspondence diagram for the (corr-)unital structure; being a flat inner fibration then guarantees that this is compatible with compositions of natural transformations.
For a morphism $F \Rightarrow G$ in $\tx{FactSpace}_{\tx{corr-un}}$, we obtain a commuting diagram:
\[\xymatrix{
 & S \times_{\tx{Fun}([1], \Ran_\un)} \tx{Fun}([1], F(\tb{1})) \ar[d] \ar^{\tx{proj}_0}[ld] \ar_{\tx{proj}_1}[rd] \\
 C_\varphi \ar[d] & S \times_{\tx{Fun}([1], \Ran_\un)} \tx{Fun}([1], F(\tb{1})) \ar^{\tx{proj}_0}[ld] \ar_{\tx{proj}_1}[rd] & C_\psi \ar[d] \\
D_\varphi & & D_\psi 
 }\]
where $D_\varphi$ and $D_\psi$ are the corresponding prestacks for $G$. The condition of being strictly unital (resp. strictly co-unital, strictly bi-unital) then corresponds to the condition of the left square (resp. right square, both squares) being Cartesian.

\begin{example}
	Any map between unital (resp. co-unital) factorization spaces, considered as one between corr-unital factorization spaces, is automatically strictly unital (resp. counital). A strictly \emph{co-unital} map between \emph{unital} factorization spaces (or vice versa) forces the non-dgenerate square in the above diagram to be Cartesian. Lax-(co)unital maps between (co)unital factorization spaces are not allowed when considering them as corr-unital factorization spaces.
\end{example}

\begin{remark}
	In the correspondence framework of \cite{raskin2015chiral} one can also write down a definition of corr-unital factorization spaces, namely as a $\tx{PreStk}_{\tx{corr}}$-valued lax prestack. To compare this with our definition one needs to compare the category of categories over $[1]$ fibered in spaces with the category of spans in spaces, which is part of Problem 0.17 of \cite{ayala2020fibrations}. We will not need such comparison here.
\end{remark}

\begin{remark}
	Using the lax colimit description one can formulate a description of the category of unital factorization spaces with partition data similar to Proposition \ref{prop:fact-spc-formulations-eqv}, but the description is more complicated. We leave it to the interested reader.
\end{remark}

\begin{convention}
	From now on, we will only introduce the construction of factorization structures for the non-unital version, with the implicit understanding that the corresponding unital version can be defined similarly. There are a few places where complications do occur, at which point we address them separately.
\end{convention}

\subsection{Factorization Module Spaces}

\paragraph{The CM Operad}
We have the colored operad $\tx{CMod}_\nonun$, which consists of two colors $\mf{a}$ and $\mf{m}$, and such that
\[\tx{Mul}(\{X_i\}, Y) = \begin{cases}
* & \tx{if}~Y = X_i = \mf{a}\tx{ for all }i \\
* & \tx{if}~Y = \mf{m}\tx{ and }X_i = \mf{m}\tx{ for exactly one }i \\
\emptyset & \tx{otherwise}
\end{cases}\]
Let $\CMnu^{\tensor}$ be the associated $\infty$-operad. Unwinding definitions we find that it is given by the following category:
\begin{itemize}
\item Its objects are pairs $(I, M)$ where $I$ is a nonempty finite set and $M \subseteq I$ is a (potentially empty) subset;
\item A morphism from $(I, M)$ to $(J, N)$ only exists when $|M| = |N|$, in which case it is a surjection $p: I \surj J$ that restricts to a bijection $M \simeq N$. We will often only write $M \subseteq I \surj J$ for the morphism when no confusion is likely to occur.
\end{itemize}

\begin{variant}
We also have the \emph{unital} CM operad $\CM^{\tensor}$, where we replace ``surjective map'' with ``arbitrary map'' in the definition above.
\end{variant}

\begin{remark}
Let $\tx{Mod}^\tx{Comm}(\mc{C})$ be the $\infty$-operad of $E_\infty$-algebras and its $E_\infty$-modules (defined in \cite[3.3.3.8]{HA}), and let $\tx{Mod}^\tx{Comm}_\nonun(\mc{C}) := \tx{Mod}^{\tx{Comm}}(\mc{C}) \times_{\CommOperad} N_*(\tx{Surj})$ be the base change such that we only allow surjective morphisms in $\CommOperad$. By \cite[4.5.1.4]{HA} (and its non-unital variant) we have an equivalence of categories $\tx{Alg}_{\CMnu}(\mc{C}) \simeq \tx{Mod}^\tx{Comm}_\nonun(\mc{C})$ and $\tx{Alg}_{\CM}(\mc{C}) \simeq \tx{Mod}^\tx{Comm}(\mc{C})$.
\end{remark}

\newcommand{\ECMOp}{\ms{ECM}^{\op}}

\begin{definition}
Taking opposite category over each arity, $(\CM^{\op})^{\tensor}$ is again a $\CM^{\tensor}$-monoidal operad. Let $\tx{Env}((\CM^{\op})^{\tensor})$ denote its $\CM$-monoidal envelope of $\CM^{\op}$, and set
\[\ECMOp := \tx{Env}((\CM^{\op})^{\tensor})_{\langle 1 \rangle}.\]
\end{definition}
This category can be explicitly described as follows:
\begin{itemize}
\item Its objects are pairs $(Z \subseteq I)$, where $I$ is a finite (potentially empty) set and $Z$ is a subset of size at most 1;
\item Those objects with $Z \simeq \emptyset$ are of color $\mf{a}$, and others are of color $\mf{m}$;
\item A morphism $(Z_1 \subseteq I_1) \to (Z_2 \subseteq I_2)$ only exists when $|Z_1| \simeq |Z_2|$, in which case it is a map $I_1 \gets I_2$ in $\Setli$ that maps $Z_2$ to $Z_1$ isomorphically;
\item For $I \in \Setli^\op$, we will use $I_\mf{a}$ to denote $(\emptyset \subseteq I) \in \ECMOp$, and $I_\mf{m}$ to denote $(\{*\} \subseteq I) \in \ECMOp$.
\end{itemize}

\begin{variant}
The non-unital version $\ECMOp_\nonun$ is defined analogously, by replacing $\Setli$ with $\fSetnu$ everywhere, and disallowing empty sets for $I$.
\end{variant}

\begin{remark}
By definition, a $\CM$-monoidal functor from $\ECMOp$ to another $\CM$-monoidal category $\mc{C}$ is precisely a \emph{lax} $\CM$-monoidal functor $\CM^\op \to \mc{C}$.
\end{remark}

Now we are ready to define factorization module spaces:

\begin{definition}
	Let $\tx{RanModMonoid}_{\nonun}$ denote the category
	\[\tx{Fun}^{\CM\tx{-lax}, \tensor}_{\nonun}(\ECMOp_{\nonun}, \tx{PreStk}^{\op}) \times_{\tx{Fun}^{\tx{lax}, \tensor}_{\nonun}(\fSetnu^\op, \tx{PreStk}^{\op})} \{\Ran_\nonun^\op\}\]
	We refer to this category as category of Ran-module monoids. If $Z$ is the image of $[1]_\mf{m}$, we will denote the corresponding element by $(\Ran_\nonun, Z)$.
\end{definition}

\begin{example}
	$\Ran_\nonun$ is itself a Ran-module monoid.
\end{example}

\begin{definition}
	For a fixed Ran-module monoid $Z$, let $\tx{FactModSpace}_{\nonun}(Z)$ denote the full subcategory of		
	\[\tx{Fun}^{\CM\tx{-lax}, \tensor}_{\nonun}(\ECMOp_{\nonun}, \tx{PreStk}^{\op})_{(\Ran_\nonun, Z)/}\]
	satisfying condition L (i.e. the natural transformations from the lax structure are Cartesian). We refer to this category as the category of non-unital factorization module spaces of type $Z$.
\end{definition}

\begin{variant}
The unital versions $\tx{RanModMonoid}_\un$ and $\tx{FactModSpace}_\un(Z)$ are defined analogously using $\ECMOp$ instead.
\end{variant}

\begin{example}
\label{example:Ran-module-spaces}
Let $\tx{Arr}(\Ran_\un)$ denote the arrow lax prestack, given by $S \mapsto \tx{Fun}([1], \Ran_\un(S))$. It admits two maps $\pi_s$ and $\pi_t$ to $\Ran_\un$, by mapping to the source and target respectively.

For any $\epsilon_R: R \to \Ran_\un$, we define
\[\Ran_{\un, \epsilon_R} := R_{\epsilon_R/} := R \times_{\epsilon_R, \Ran_\un, \pi_s} \tx{Arr}(\Ran_\un).\]
For any map $\eta: R \to \Ran$ to the non-unital Ran space, we define $\Ran_\eta$ as
\[\Ran_{\nonun, \eta} := \left(\Ran_{\un, (\Ran \inj \Ran_\un) \circ \eta}\right) \hspace{1em} \Ran_{\nonun, \eta} := \Ran_{\un, \eta}^{\tx{PreStk}}.\]
We sometimes also the notation $\tx{Ran}_R$ to denote this space, when the map is clear from the context. Explicitly, on the object-level it is given by
\[\Ran_{\nonun, R}(S) = \{\alpha \in R(S), \beta \in \Ran(S) \mid \epsilon_R \circ \alpha \subseteq \beta\}\]
We have evident maps
\[\tx{pr}_R: (\alpha, \beta) \mapsto \alpha \hspace{1em} \tx{pr}_\Ran: (\alpha, \beta) \mapsto \beta;\]
and the action of $\Ran_\un$ on $\Ran_{\un, R}$ is by union with $\beta$.

Cases that we will need include (all for the unital version):
\begin{itemize}
\item $R = \emptyset$, so we get $\Ran_\un$, the unmarked Ran space itself;
\item For $S$ affine and a fixed subgroupoid $\overrightarrow{x} = \{x_1, \ldots, x_n\} \subset X(S)$, we let $\Ran_{\un, \overrightarrow{x}}$ denote the lax prestack $\Ran_S$ for the map $S \xrightarrow{\overrightarrow{x}} X^n \to \Ran \to \Ran_{\un}$;
\item $R = X^I \to \Ran \to \Ran_{\un}$, in which case we denote the resulting lax prestack simply by $\Ran_{\un, X^I}$;
\item Recursively, we can let $R = \Ran_{\un, R'}$ via $\tx{pr}_\Ran$ for any $R'$, so on and so forth.
\end{itemize}
\end{example}

\subsection{Factorization Categories and Algebras}
The functor $\tx{ShvCat}(-): \tx{PreStk}^\op \to 1\tx{-Cat}$ corresponds to a Cartesian fibration $\tx{PreStk}^\tx{ShvCat} \to \tx{PreStk}$. Similarly, the functor $\tx{ShvCatSect}(-): \tx{PreStk}^\op \to 1\tx{-Cat}$ which sends $S$ to the category of pairs $(\mc{F} \in \tx{ShvCat}(S), f \in \Gamma(S, \mc{F}))$ also gives rise to a Cartesian fibration $\tx{PreStk}^{\tx{ShvCatSect}} \to \tx{PreStk}$. The natural forgetful functor then gives a tower of maps
\[\tx{PreStk}^{\tx{ShvCatSect}} \to \tx{PreStk}^\tx{ShvCat} \to \tx{PreStk};\]
Note that the first map is not a (co)Cartesian fibration.
Box product then equips each level with a symmetric monoidal structure, such that all maps above are symmetric monoidal.
\begin{definition}
	Let $Z$ be a non-unital factorization space. Let $\tx{FactCat}_{\nonun, Z}$ denote the full subcategory of
	\[\tx{Fun}^{\tx{oplax}, \tensor}_\nonun(\fSetnu, \tx{PreStk}^\tx{ShvCat}) \times_{ \tx{Fun}^{\tx{oplax}, \tensor}_\nonun(\fSetnu, \tx{PreStk})} \{Z\}\]
	consisting of those maps $F$ satisfying the following conditions:
\begin{itemize}
\item[(Condition R)]For every $I \surj J \in \fSetnu$, the arrow $F(I) \to F(J)$ is Cartesian with respect to the Cartesian fibration $\tx{PreStk}^\tx{ShvCat} \to \tx{PreStk}$;
\item[(Condition L)]For every $I \surj J \in \fSetnu$, the arrow $F(I) \to \bigotimes_{j} F(I_j)$ is Cartesian with respect to the Cartesian fibration $\tx{PreStk}^\tx{ShvCat} \to \tx{PreStk}$.
\end{itemize}	
	An element of this category is called a non-unital factorization category over $Z$.
\end{definition}

\begin{convention}
	Our \emph{default} choice is $Z = \Ran_\dR$ (in the non-unital case) and $\Ran_{\un, \dR}$ (in the unital case); in such case, we omit the $Z$ subscript. When we work over other factorization spaces (e.g. $\Ran_\nonun$ itself), we will explicit warn our readers.
\end{convention}

\begin{definition}
	Let $Z$ be a non-unital factorization space and $C \in \tx{FactCat}_\nonun(Z)$. Let $\tx{FactAlg}_\nonun(C)$ denote the full subcategory of
	\[\tx{Fun}^{\tx{oplax}, \tensor}_\nonun(\fSetnu, \tx{PreStk}^\tx{ShvCatSect}) \times_{ \tx{Fun}^{\tx{oplax}, \tensor}_\nonun(\fSetnu, \tx{PreStk}^\tx{ShvCat})} \{C\}\]
	satisfying the corresponding conditions L and R with respect to the Cartesian fibration $\tx{PreStk}^\tx{ShvCatSect} \to \tx{PreStk}$. An element of this category is called a non-unital factorization algebra within $C$.
\end{definition}

\begin{variant}
	Let $\tx{PreStk}^{\op, \tx{ShvCat}, \vee} \to \tx{PreStk}^\op$ denote the dual \emph{coCartesian} fibration associated with $\tx{ShvCat}$. Let $\tx{FactCat}_{\nonun}^{\CoMark}(Z)$ denote the full subcategory of
	\[\tx{Fun}^{\tx{lax}, \tensor}_\nonun(\fSetnu^\op, \tx{PreStk}^{\op, \tx{ShvCat}, \vee}) \times_{\tx{Fun}^{\tx{lax}, \tensor}_\nonun(\fSetnu^op, \tx{PreStk}^\op)} \{Z^\op\}\]
	satisfying the following conditions:
\begin{itemize}
\item[(Condition R)]For every $I \surj J \in \fSetnu$, the arrow $F(J) \to F(I)$ is coCartesian with respect to the coCartesian fibration $\tx{PreStk}^{\op, \tx{ShvCat}, \vee} \to \tx{PreStk}^\op$;
\item[(Condition L)]For every $I \surj J \in \fSetnu$, the arrow $\bigotimes_{j} F(I_j) \to F(I)$ is coCartesian with respect to the coCartesian fibration $\tx{PreStk}^{\op, \tx{ShvCat}, \vee} \to \tx{PreStk}^\op$.
\end{itemize}
The category $ \tx{FactAlg}_{\nonun}^{\CoMark}(Z)$ can be defined analogously.
\end{variant}

\begin{remark}
Obviously we have $\tx{FactCat}_{\nonun}^{\CoMark}(Z) \simeq \tx{FactCat}_{\nonun}(Z)$.
\end{remark}

\subsubsection{Partition Description}

In alignment with the choices made above, we make the following definition:

\begin{definition}
We define $\tx{FactCat}_{\nonun, \WeakMark}$ as the full subcategory of
\[\tx{Fun}^{\tx{oplax}, \tensor}_\nonun(\fSetnu, \tx{PreStk}^{\tx{ShvCat}}) \times_{ \tx{Fun}^{\tx{oplax}, \tensor}_\nonun(\fSetnu, \tx{PreStk})} \{\Ran_\dR\}\]
such that
\begin{itemize}
\item[(Condition L)]For every $I \to J \in \fSetnu$, the arrow $F(I) \to \bigotimes_{j} F(I_j)$ is Cartesian with respect to the Cartesian fibration $\tx{PreStk}^{\tx{ShvCat}} \to \tx{PreStk}$.
\end{itemize}
Similar we define $\tx{FactAlg}_{\nonun, \WeakMark}(C)$ for any factorization category $C$.
\end{definition}

For technical purposes, however, we will exclusively work with the dual version:

\begin{variant}
We define $\tx{FactCat}^{\CoMark}_{\nonun, \LWeakMark}$ as the full subcategory of
\[\tx{Fun}^{\tx{lax}, \tensor}_\nonun(\fSetnu^\op, \tx{PreStk}^{\op, \tx{ShvCat}, \vee}) \times_{ \tx{Fun}^{\tx{lax}, \tensor}_\nonun(\fSetnu^\op, \tx{PreStk}^{\op})} \{\Ran_\dR^\op\}\]
such that
\begin{itemize}
\item[(Condition R)]For every $I \to J \in \fSetnu$, the arrow $F(J) \to F(I)$ is coCartesian with respect to the coCartesian fibration $\tx{PreStk}^{\op, \tx{ShvCat}, \vee} \to \tx{PreStk}^{\op}$.
\end{itemize}
We can similarly define $\tx{FactAlg}^{\CoMark}_{\nonun, \LWeakMark}(C)$ for any factorization category $C$.
\end{variant}

\begin{remark}
	We note that the notion definition of a (strong) factorization algebra / module only requires a \emph{weak} factorization category as background.
\end{remark}

\begin{remark}
There exists an equivalence of categories $\tx{FactCat}_{\nonun, \WeakMark} \simeq \tx{FactCat}^{\CoMark}_{\nonun, \LWeakMark}$ compatible with the embeddings from $\tx{FactCat}_{\nonun}(Z) \simeq \tx{FactCat}_{\nonun}^{\CoMark}(Z)$. The analogous statement holds for factorization algebras.
\end{remark}

\begin{remark}
\label{rmk:comparison-with-chiralcats}
Note that, in \cite{raskin2015chiral}, a weak factorization category (resp. algebra) is defined as a commutative algebra in a certain correspondence category of Grothendieck fibrations; such notion coincides with our $\tx{FactCat}^{\CoMark}_{\nonun, \LWeakMark}$. The comparison follows easily from \cite[2.2.4.9]{HA} and Lemma \ref{lemma:twisted-arrow-maps}.\footnote{Note that our condition L is hidden in \cite[5.14]{raskin2015chiral}, and that composability condition in Lemma \ref{lemma:twisted-arrow-maps} is vacuous in this case as remarked in \cite[8.12]{raskin2015chiral}.} The algebra case follows by comparing the category $(\tx{PreStk}^\tx{ShvCatSect})_\tx{corr}$ with the category $\mc{G}$ in \cite[5.26]{raskin2015chiral}.
\end{remark}

Now we introduce the partition description:

\begin{definition}
We let $\tx{FactCat}^{\CoMark}_{\nonun, \LWeakMark, \PartitionMark}$ denote the full subcategory of 
\[\tx{Fun}^{\tx{lax}, \tensor}_\nonun(\ECOp, \tx{PreStk}^{\op, \tx{ShvCat}, \vee}) \times_{ \tx{Fun}^{\tx{lax}, \tensor}_\nonun(\ECOp, \tx{PreStk}^{\op})} \{U_{\Ran_\dR}^\op\}\]
consisting of maps $F$ such that:
\begin{itemize}
\item[(Condition R')]For every $(I \to J) \to (I' \to J') \in \ECOp$, the arrow $F(I \to J) \to F(I' \to J')$ is coCartesian with respect to the coCartesian fibration $\tx{PreStk}^{\op, \tx{ShvCat}, \vee} \to \tx{PreStk}^\op$.
\end{itemize}
We let $\tx{FactCat}^{\CoMark}_{\nonun, \PartitionMark}$ denote the full subcategory of $\tx{FactCat}^{\CoMark}_{\nonun, \LWeakMark, \PartitionMark}$ where additionally,
\begin{itemize}[resume]
\item[(Condition L')]For every $I \to J \to K \in \fSetnu$, the arrow $\bigotimes_{k} F(I_k \to J_k) \to F(I \to J)$ is coCartesian with respect to the coCartesian fibration $\tx{PreStk}^{\op, \tx{ShvCat}, \vee} \to \tx{PreStk}^{\op}$.
\end{itemize}
The category $\tx{FactAlg}^{\CoMark}_{\nonun, \LWeakMark, \PartitionMark}(C)$ and $\tx{FactAlg}^{\CoMark}_{\nonun, \LWeakMark, \PartitionMark}(C)$ for a fixed factorization category $C$ are defined analogously.
\end{definition}

\begin{proposition}
\label{prop:fact-cat-alg-partition-description}
There is a natural equivalence $\tx{FactCat}^{\CoMark}_{\nonun, \LWeakMark, \PartitionMark} \simeq \tx{FactCat}^{\CoMark}_{\nonun, \LWeakMark}$ which restricts to an equivalence $\tx{FactCat}^{\CoMark}_{\nonun, \PartitionMark} \simeq \tx{FactCat}^{\CoMark}_{\nonun}$. The analogous statement holds for factorization algebras in a factorization category $C$.
\end{proposition}
\begin{proof}
We treat the category case by essentially copying the proof from \cite[8.4.1]{raskin2015chiral} with added details; the algebra case is completely analogous.
Let $\pi: \ETOp \to \ECOp^\op$ denote the map $(I \to J \to K) \mapsto (J \to K)$.
Passing to the monoidal envelope we obtain an adjunction
\[\tx{Fun}^{\tx{lax}, \tensor}_\nonun(\fSetnu^\op, \tx{PreStk}^{\op, \tx{ShvCat}, \vee}) \simeq \tx{Fun}^{\tensor}_\nonun(\ECOp^\op, \tx{PreStk}^{\op, \tx{ShvCat}, \vee})\]
\[\tx{Fun}^{\tx{lax}, \tensor}_\nonun(\ECOp, \tx{PreStk}^{\op, \tx{ShvCat}, \vee}) \simeq \tx{Fun}^{\tensor}_\nonun(\ETOp, \tx{PreStk}^{\op, \tx{ShvCat}, \vee})\]
We hope to use monoidal right Kan extension along $\pi$ to relate the RHS categories. In general, however, box product of sheaves of categories do not distribute over limits, so RKE along $\pi$ only yields a map between lax-monoidal functor categories. Nevertheless, we claim that for any $F, G \in \tx{Fun}^{\tx{lax}, \tensor}_\nonun(\ECOp, \tx{PreStk}^{\op, \tx{ShvCat}, \vee})$ satisfying condition (R'),
\[(\limit_{I \in \fSetnu^\op} (\pi_I)_*(F(I))) \boxtimes (\limit_{J \in \fSetnu^\op} (\pi_J)_*(G(J)) ) \to \limit_{I, J \in \fSetnu^\op} (\pi_I)_*(F(I)) \boxtimes (\pi_J)_*(G(J))\]
(where $\pi_I: X^I_\dR \to \Ran_\dR$ is the structure map) is an isomorphism. Indeed, condition (R') guarantees that for any affine $S$ and any $I_1 \surj I_2$, the section $(\pi_{I_2})_*(F(I_2))(S) \to (\pi_{I_1})_*(F(I_1))(S)$ admits a left adjoint, so we can turn the limits above to colimits and commute with box product. It follows that when restricting to this subcategory, RKE (and $\tx{Res}$) preserve the strictly symmetric monoidal functors. Moreover, for the same reason the resulting right Kan extension will satisfying condition (R). Thus monoidal right Kan extension induces an adjunction
\[\tx{Conv}_{\Ran \to \PartitionMark}: \tx{FactCat}^{\CoMark}_{\nonun, \LWeakMark} \adjoint \tx{FactCat}^{\CoMark}_{\nonun, \PartitionMark, \LWeakMark}: \tx{Conv}_{\PartitionMark \to \Ran};\] %
explicitly, the right adjoint sends $F: \fSetnu^\op \to \tx{PreStk}^{\op, \tx{ShvCat}, \vee}$ to
\[p = (I \surj J) \mapsto \pi_{p}^*(F(J)) \in \tx{ShvCat}(U_{\Ran_\dR}(p)) \hspace{1em} \pi_p: U_{\Ran_\dR}(p) \to \Ran_\dR(J);\]
and the left adjoint sends $T: \ECOp \to \tx{PreStk}^{\op, \tx{ShvCat}, \vee}$ to
\[J \mapsto \limit_{p = (J' \surj J)} (\pi_p)_*(T(p)) \in \tx{ShvCat}(\Ran_\dR(J)).\]

It is easy to check that the full subcategories of factorization objects are preserved by this adjunction, and it remains to check the unit and counit of this adjunction are isomorphisms.
Condition R and R' guarantee that forgetful functors $\tx{FactCat}^{\CoMark}_{\nonun, \LWeakMark} \to \tx{ShvCat}(\Ran_\dR)$ and
\[\tx{FactCat}^{\CoMark}_{\nonun, \LWeakMark, \PartitionMark} \to \tx{Fun}(\fSetnu, \tx{PreStk}^{\op, \tx{ShvCat}, \vee}) \times_{\tx{Fun}(\fSetnu, \tx{PreStk}^\op)} \{I \mapsto X^I_\dR\}\]
are both conservative. Unfolding definition, we see that we need to check
\[F \simeq \limit_{I \in \fSetnu^\op} (\pi_I)_* (\pi_I)^* F \hspace{1em} \hspace{1em} T(J) \simeq \pi_J^*(\limit_{I} (\pi_I)_* T(I)) \forall J \in \fSetnu^\op\]
for every $F$ which is the restriction of an element in $\tx{FactCat}^{\CoMark}_{\nonun, \LWeakMark}$ (resp. $T$ a restriction of an element in $\tx{FactCat}^{\CoMark}_{\nonun, \LWeakMark, \PartitionMark}$). The unit follows from the fact that $\tx{ShvCat}$ is a right Kan extension.

Set $T_J := \Gamma(X^J_\dR, T(J))$.
For $S$ affine mapping to $X^J_\dR$, the two sides of the counit evaluates to
\[\QCoh(S) \tensor_{\DMod(X^J)} T_J \to \limit_I \limit_{R \in \mc{D}_{I, J}^\op} \QCoh(R) \tensor_{\DMod(X^I)} T_I\]
where $\mc{D}_{I, J}$ is the category of diagrams (where $R$ is affine)
\[\xymatrix{
R \ar[dd] \ar[rrd] \ar@{-->}[rd] \\
 & X^I_\dR \times_{\Ran_\dR} X^J_\dR \ar[r] \ar[d] & X^I_\dR \ar[d] \\
S \ar[r] & X^J_\dR \ar[r] & \Ran_\dR
}\]
Note, however, for every $|I| > |J|$, the fiber product is simply $X^J_\dR$ and $\limit_{R \in \mc{D}_{I, J}^\op} \QCoh(R) \tensor_{\DMod(X^I)} T_I \simeq T_J$. We can restrict $I$ to the subcategory $\{I \in \fSetnu^\op \mid |I| > |J|\}$ since it is left cofinal in $\fSetnu^\op$; then RHS becomes a connected \emph{constant} diagram with value $\QCoh(S) \tensor_{\DMod(X^J)} T_J$.
\end{proof}

\subsubsection{The Unital Case}
First note that $\tx{ShvCat}$ can be defined for all lax prestacks (\cite[4.19]{raskin2015chiral}). Explicitly, if $S$ is affine and $\mc{Y}$ is a lax prestack, a sheaf of categories $\mc{F}$ over $\mc{Y}$ attaches to each $\varphi \Rightarrow \psi \in \mc{Y}(S)$ an $\QCoh(S)$-linear map $\mc{F}_\varphi \to \mc{F}_\psi$.
We now define unital factorization categories and algebras.

\begin{definition}
Let $Z$ be a unital factorization space.
We define $\tx{FactCat}^{\CoMark}_{\un, \LWeakMark, Z}$ as the full subcategory of
\[\tx{Fun}^{\tx{lax}, \tensor}(\Setli^\op, \tx{LaxPreStk}^{\op, \tx{ShvCat}, \vee}) \times_{ \tx{Fun}^{\tx{lax}, \tensor}(\Setli^\op, \tx{LaxPreStk}^{\op})} \{Z^\op\}\]
consisting of maps $F$ such that
\begin{itemize}
\item[(Condition R)]For every surjective $I \surj J \in \fSetnu$, the arrow $F(J) \to F(I)$ is coCartesian with respect to the coCartesian fibration $\tx{LaxPreStk}^{\op, \tx{ShvCat}, \vee} \to \tx{LaxPreStk}^{\op}$;
\item[(Condition Z)]For every map $[0] \to I \in \Setli$, the arrow $F(I) \to F([0])$ is coCartesian with respect to the coCartesian fibration $\tx{LaxPreStk}^{\op, \tx{ShvCat}, \vee} \to \tx{LaxPreStk}^{\op}$;
\item[(Condition C)] For every composable pair of morphisms $I \to J \to K$ in $\Setli$, the commutative square
\[\xymatrix{
	\bigotimes_{j = 1}^{J} F([1]) \ar[r] \ar[d] & \bigotimes_{k = 1}^{K} F(J_k) \ar[d] \\
	\bigotimes_{j = 1}^{J} F(I_j) \ar[r] & \bigotimes_{k = 1}^{K} F(I_k)
}\]
	is coCartesian within $\tx{LaxPreStk}^{\op, \tx{ShvCat}, \vee}$;
\end{itemize}
we define $\tx{FactCat}^{\CoMark}_{\un, Z}$ as the full subcategory consisting of objects that additionally satisfy
\begin{itemize}[resume]
\item[(Condition L)]For every surjective $I \surj J \in \fSetnu$, the arrow $\bigotimes_{j} F(I_j) \to F(I)$ is coCartesian with respect to the coCartesian fibration $\tx{LaxPreStk}^{\op, \tx{ShvCat}, \vee} \to \tx{LaxPreStk}^\op$.
\end{itemize}
\end{definition}

\begin{remark}
Again, in the non-unital case condition C was subsumed.
\end{remark}

\begin{variant}
In accordance with previous notations, we let $\tx{FactCat}_{\un, Z}$ denote the category obtained from $\tx{FactCat}^{\CoMark}_{\un, Z}$ by passing to dual fibrations; elements of $\tx{FactCat}_{\un, Z}$ will in particular be \emph{oplax} functors.
\end{variant}

\begin{variant}
For a given unital factorization category $C$, the categories
\[\tx{FactAlg}^{\CoMark}_{\un, \LWeakMark, Z}(C) \hspace{1em} \tx{FactAlg}^{\CoMark}_{\un, Z}(C) \hspace{1em} \tx{FactAlg}_{\un, Z}(C)\]
can be defined analogously.
\end{variant}

\begin{remark}
In definition of \emph{weak} factorization algebras there is no room for laxness for the unit (because the background category is $(\infty, 1)$), i.e. the value of $[0]$ in $C([0])$ (which obtains a symmetric monoidal category structure) must be the monoidal unit.
\end{remark}

Let $U_{\Ran_{\un, \dR}}: \ECOp^\op_{\un} \to \tx{PreStk}$ be the functor defined using the same formula as in Definition \ref{def:u-ran}. The following partition description of unital factorization categories and algebras is obtainable  following the same strategy as Proposition \ref{prop:fact-cat-alg-partition-description} above:

\begin{definition}
We let $\tx{FactCat}^{\CoMark}_{\un, \LWeakMark, \PartitionMark}$ denote the full subcategory of 
\[\tx{Fun}^{\tx{lax}, \tensor}(\ECOp_\un, \tx{PreStk}^{\op, \tx{ShvCat}, \vee}) \times_{ \tx{Fun}^{\tx{lax}, \tensor}(\ECOp_un, \tx{PreStk}^{\op})} \{U_{\Ran_{\un, \dR}}^\op\}\]
consisting of maps $F$ such that:
\begin{itemize}
\item[(Condition R')]For every $(I \to J) \to (I' \to J') \in \ECOp \subseteq \ECOp_\un$, the arrow $F(I \to J) \to F(I' \to J')$ is coCartesian with respect to the coCartesian fibration $\tx{PreStk}^{\op, \tx{ShvCat}, \vee} \to \tx{PreStk}^\op$;
\item[(Condition Z')]For every $([0] \to I) \to ([0] \to J)$, the arrow $F([0] \to I) \to F([0] \to J)$ is coCartesian with respect to the coCartesian fibration $\tx{PreStk}^{\op, \tx{ShvCat}, \vee} \to \tx{PreStk}^\op$.
\end{itemize}
We let $\tx{FactCat}^{\CoMark}_{\un, \PartitionMark}$ denote the full subcategory of $\tx{FactCat}^{\CoMark}_{\un, \LWeakMark, \PartitionMark}$ where additionally,
\begin{itemize}[resume]
\item[(Condition L')]For every $I \to J \to K$ maps of finite sets, the arrow $\bigotimes_{k} F(I_k \to J_k) \to F(I \to J)$ is coCartesian with respect to the coCartesian fibration $\tx{PreStk}^{\op, \tx{ShvCat}, \vee} \to \tx{PreStk}^{\op}$.
\end{itemize}
The category $\tx{FactAlg}^{\CoMark}_{\un, \LWeakMark, \PartitionMark}(C)$ and $\tx{FactAlg}^{\CoMark}_{\un, \LWeakMark, \PartitionMark}(C)$ for a fixed factorization category $C$ are defined analogously.
\end{definition}

\begin{proposition}[\cite{raskin2015chiral}]
\label{prop:fact-cat-alg-partition-description-unital}
There is a natural equivalence $\tx{FactCat}^{\CoMark}_{\un, \LWeakMark, \PartitionMark} \simeq \tx{FactCat}^{\CoMark}_{\un, \LWeakMark}$ which restricts to an equivalence $\tx{FactCat}^{\CoMark}_{\un, \PartitionMark} \simeq \tx{FactCat}^{\CoMark}_{\un}$. The analogous statement holds for factorization algebras in a factorization category $C$.
\end{proposition}

\medskip

Suppose $C$ is a \emph{unital} factorization category over $\Ran_{\un, \dR}$; in particular, this means that $C(0) \simeq \Vect$. The unital setting is of particular interest to us for the following reasons:

\begin{proposition}
\label{prop:unital-basic-properties}
Fix a map $\epsilon: R \to \Ran_\un$ for some laft prestack $R$, and let
\[\epsilon_\dR: R_\dR \to \Ran_{\un, \dR}\]
be the induced map, which defines a lax prestack $\Ran_{\un, \dR, \epsilon_\dR}$.
Then:
	\begin{enumerate}
	\item Let $C$ be a unital factorization category over $\Ran_{\un, \dR}$, then there exists a canonically defined \emph{unit factorization algebra} $\tx{unit}_C$ and, for every unital factorization algebra $A$ in $C$, there is an unique map of unital factorization algebras $\eta_A: \tx{unit}_C \to A$;
	\item Let $M$ be a unital $C$-factorization module category over $\Ran_{\un, \dR, \epsilon_\dR}$, and let $e: R_\dR \to \Ran_{\un, \dR, \epsilon_\dR}$ be the inclusion map, then we have an equivalence of categories
	\[\tx{unit}_C\tx{-FactMod}_\un(M) \simeq \Gamma(R_\dR, e^*(M));\]
	\item Let $F$ be a lax unital factorization functor between two pairs $(C, M) \to (D, N)$ of factorization categories and module categories over the pair $(\Ran_{\un, \dR}, \Ran_{\un, \dR, \epsilon_\dR})$, then we have a map of factorization algebras $\tx{unit}_D \to F(\tx{unit}_C)$ and a map
	\[\tx{unit}_C\tx{-FactMod}_\un(M) \to F(\tx{unit}_C)\tx{-FactMod}_\un(N).\]
	\end{enumerate}
\end{proposition}

\begin{proof}
	The proof in \cite[Section 6]{raskin2015chiral} applies \emph{mutatis mutandis}.
\end{proof}

\subsection{Restriction of Modules}

\label{sect:fact-restriction}

We shall adopt the following notations:
\begin{itemize}
\item Fix a map $\epsilon: R \to \Ran_\un$ for $R$ a laft prestack. Let $\Ran_{\dR, \un, \epsilon_\dR}$ denote the corresponding factorization module space given by $\epsilon_\dR$;
\item Let $C_M$ be a weak unital $\Ran_{\dR, \un, \epsilon_\dR}$-factorization module category, and let $C$ be the algebra part of $C_M$;
\item Let $f: A \to B$ be a map of (strong) \emph{unital} factorization algebras in $C$;
\item For any $X$ in $\Setli^\op$ ($\ECMOp$ in the module case), set $\Gamma_C(X) := \Gamma(\Ran_{\un, \dR}(X), C(X))$ and $\Gamma_M(X) := \Gamma(\Ran_{\un, \dR, \epsilon_\dR}(X), C_M(X))$;
\end{itemize}

The main fact we establish in this section is:

\begin{proposition}
\label{prop:restriction-functor}
	In the above setting, there exists an a well-defined continuous functor
	\[\tx{res}_{f}: B\tx{-FactMod}_\un(C_M) \to A\tx{-FactMod}_\un(C_M);\]
This functor has the following properties:
\begin{enumerate}
\item It is functorial in $f$;
\item $\epsilon_\dR^! \circ \tx{res}_f \simeq \epsilon_\dR$ for every $f$.
\end{enumerate}
An analogous non-unital version also exists.
\end{proposition}

\begin{example}
Let $B$ be a unital factorization algebra within $\DMod(\Ran_\un)$ and let $M_0$ denote a unital $B$-factorizable module within $\DMod(\Ran_{\un, X})$. Let $f: A \to B$ be a map of unital factorization algebras. Then $\tx{res}_f(M_0)|_X \simeq M_0|_X$, and $\tx{res}_f(M_0)|_{X^2}$ is the pullback
\[\xymatrix{
\tx{res}_f(M_0)|_{X^2} \ar[r] \ar[d] & M_0|_{X^2} \ar[d] \\
j_* j^*(A|_{X} \boxtimes M_0|_{X}) \ar[r] & j_* j^*(B|_{X} \boxtimes M_0|_{X})
}\]
where $j: X^2 \setminus \Delta_X \inj X^2$ is the inclusion of the complement of the diagonal. The unital structure at this level amounts to two maps
\[\tx{res}_f(M_0)|_X \boxtimes \omega_X \to \tx{res}_f(M_0)|_{X^2} \hspace{1em} \omega_X \boxtimes \tx{res}_f(M_0)|_X \to \tx{res}_f(M_0)|_{X^2}\]
which are induced by the original unital $B$-module structure and the unital structure of $A$.
\end{example}

\begin{remark}
This construction requires working in the stable setting; in particular, restriction of factorization module \emph{categories} do not exist in general. A \emph{factorization induction} functor
\[\tx{ind}_f: A\tx{-FactMod}_\un(C_M) \to B\tx{-FactMod}_\un(C_M)\]
also does not exist in general.
\end{remark}

\begin{remark}
In the \emph{non-unital} situation, following \cite{francis2012chiral}, one can equip $\Gamma_C([1])$ with the structure of a non-unital symmetric monoidal algebra within $\DGCat$ (known as the \emph{chiral tensor} structure) by
\[\Gamma_C([1]) \tensor \Gamma_C([1]) \xrightarrow{\tensor} \Gamma_C([2]) \xrightarrow{\varphi_{2 \to 1}^R} \Gamma_C([1])\]
where $\varphi_{2 \to 1}^R$ is the right adjoint to the structural functor. Similarly, $\Gamma_M([1]_\mf{m})$ obtains a structure of a (non-unital) commutative module over $\Gamma_C([1])$.
In such case, we have a chiral Koszul duality (\cite{francis2012chiral}, generalized in \cite[Section 6]{raskin2015chiral}) between non-unital factorization modules and non-unital \emph{chiral} modules.
For the sake of reader's convenience, we recall the statement in the next proposition.
One can show that the evident restriction functor of non-unital chiral algebras coincides with the non-unital version of the construction we give below under chiral Koszul duality, but this method does not extend easily to the unital setting.
\end{remark}

\begin{proposition}[Chiral Koszul Duality]
	\label{prop:chiral-koszul-duality}
	The following are true:
	\begin{enumerate}
	\item The category $\tx{FactAlg}_{\nonun}(C)$ identifies with the full subcategory of $\tx{CocommCoalg}_\nonun(\Gamma_C)$ consisting of objects satisfying the \emph{factorizability} condition specified in \cite{francis2012chiral};
	\item Let $\tx{ChiralAlg}_{\nonun}(C)$ denote the full subcategory of $\tx{LieAlg}_\nonun(\Gamma_C)$ consisting of objects lies in the full subcategory $\Gamma(X_\dR, C([1])|_{X_\dR}) \subseteq \Gamma_C$. (Objects in this category will be referred to as non-unital \emph{chiral} algebras in $C$.) Then we have an equivalence
	\[\tx{KD}: \tx{ChiralAlg}_{\nonun}(C) \simeq \tx{FactAlg}_{\nonun}(C): \tx{KD}^{-1}\]
	\item For every $A \in \tx{FactAlg}_{\nonun}(C)$ and $M$ as above, $A\tx{-FactMod}_{\nonun}(M)$ identifies with the full subcategory of $A\tx{-CocommComod}_\nonun(\Gamma_M)$ consisting of objects satisfying the factorizability condition;
	\item For $L \in \tx{ChiralAlg}_{\nonun}(C)$, let $L\tx{-ChiralMod}_\nonun(M)$ denote the full subcategory of $L\tx{-LieMod}_{\nonun}(\Gamma_M)$ whose underlying object lies in the full subcategory\footnote{Here we use the fact that $\iota_!: \DMod(R) \to \DMod((\Ran_{\un})_R)$ is fully faithful, as follows from Kashiwara's lemma.} $\Gamma(R_\dR, M|_{R_\dR}) \subseteq \Gamma_M$, then we have an equivalence
	\[\tx{KD}: L\tx{-ChiralMod}_\nonun(M) \simeq \tx{KD}(L)\tx{-FactMod}_{\nonun}(M): \tx{KD}^{-1}.\]
	\end{enumerate}
\end{proposition}

The idea of the construction is simple: we first forget the data $(B, M)$ to a \emph{weak} $B$-factorization module, which then allows us to restrict to a \emph{weak} $A$-factorization module $(A, M)$. We then apply the right adjoint to the forgetful functor from factorization modules to weak factorization modules to rebuild a factorization module. The work then lies in giving a homotopically coherent definition of this ``rebuilding'' procedure---for which we use the partition description of factorization algebras/modules and the theory of operadic Kan extension.
The rest of the proof is pure combinatorics.

We first treat the non-unital case.

\begin{convention}
The rest of the story is applicable to both the algebra and the module case, so we 
\begin{itemize}
\item let $\baseOp$ denote either $\fSetnu$ or $\CMnu$;
\item Let $\ECOp$ denote $\tx{Tw}(\baseOp)^\op$, which is again $\baseOp$-monoidal;
\item Let $\ETOp := \tx{Env}_{\baseOp}(\ECOp)$;
\item let $S$ be $\Ran_\dR$ for $\mc{O}^{\tensor} \simeq \fSetnu$ and $(\Ran_\dR, \Ran_{\dR, \epsilon_\dR})$ for $\mc{O}^{\tensor} \simeq \CMnu$; and
\item let $\mc{C}$ be a factorization category for $\mc{O}^{\tensor} \simeq\fSetnu$, and a factorization category with a module category over $Z$ for $\mc{O}^{\tensor} \simeq \CMnu$.
\end{itemize}
\end{convention}

We will often implicitly consider $\ECOp$ as a full subcategory of $\ETOp$ by
\[\iota: (I \surj J) \mapsto (I \surj J \surj 1);\] this embedding $\iota$ has a left adjoint
\[\tx{pr}: (I \surj J \surj K) \mapsto (I \surj J).\]

Fix a factorization category $\mc{C}$, considered as a $\baseOp$-monoidal functor $\ECOp \to \tx{PreStk}^{\op, \tx{ShvCat}, \vee}$ via Proposition \ref{prop:fact-cat-alg-partition-description}. Let us consider the functor
\[\ECOp \to \tx{DGCat}_{\tx{cont}}, p \mapsto \Gamma(U_\Ran(p), \mc{C})\]
and its Grothendieck fibration $\msc{C} \to \ECOp$, which is a map of $\baseOp$-monoidal categories. A section of this fibration is exactly a section to $\tx{PreStk}^{\op, \tx{ShvCatSect}, \vee}$ that lifts $\mc{C}$, so we see that
\begin{corollary}
\label{cor:algebra-partition-non-unital}
	The category of weak non-unital factorizable $\baseOp$-monoidal algebras in $\mc{C}$ is equivalent to the category of lax $\baseOp$-monoidal sections of $\mathscr{C} \to \ECOp$ sending all arrows to coCartesian arrows. The subcategory of factorization objects corresponds to the $\baseOp$-monoidal sections.
\end{corollary}

From now on, fix $F: \ECOp \to \mathscr{C}$ a weak factorization object, and $\widetilde{F}: \ETOp \to \mathscr{C}$ the corresponding symmetric monoidal functor. There is an evident forgetful functor from strong factorization objects to weak ones, which we denote by $\oblv_{s \to w}$. We claim that this functor has a \emph{right} adjoint, given by $\baseOp$-monoidally right Kan extending $\widetilde{F}$ along $\ms{pr}$. To apply monoidal right Kan extension, we first check the combinatorial condition in Lemma \ref{lemma:sm-right-kan-extension}:

\begin{lemma}
\label{lemma:initial-et}
	For each $p, q \in \ECOp$, The tensor product functor
	\[\ETOp_{p/} \times \ETOp_{q/} \to \ETOp_{p \sqcup q /}\]
	is left cofinal.
\end{lemma}

\begin{convention}
	For any $p: I_1 \surj J_1$, we define a full subcategory $\ETOp^\tx{top}_{p/} \subseteq \ETOp_{p/}$ consisting of $(I'_1, J'_1, K'_1)$ such that $I'_1 \to I_1$ is an isomorphism. Similarly, we define a full subcategory $\ECOp^\tx{top}_{p/} \subseteq \ECOp_{p/}$ consisting of $(I'_1, J'_1)$ such that $I'_1 \to I_1$ is an isomorphism.
\end{convention}

\begin{proof}
The $\CMnu$ case immediately reduces to the $\fSetnu$ case. Suppose $p = I_1 \surj J_1$ and $q = I_2 \surj J_2$.
It is clear (by Quillen's theorem A) that $\ETOp^\tx{top}_{p/}$ is left cofinal within $\ETOp_{p/}$, so it suffices to prove the stronger statement that
	\[\ETOp^\tx{top}_{p/} \times \ETOp^\tx{top}_{q/} \to \ETOp^\tx{top}_{p \sqcup q /}\]
	is left cofinal. An element $\ETOp^\tx{top}_{p \sqcup q /}$ is the following data:
		\[
\xymatrix{
	I_1 \sqcup I_2 \ar@{->>}[rd] \ar@{->>}[rr]  && J_1 \sqcup J_2 \\
	& J \ar@{->>}[ru] \ar@{->>}[r] & K
}	
	\]
	The essential image of the tensor product, on the other hand, consists of the following data:
	\[
\xymatrix{
	I_1 \sqcup I_2 \ar@{->>}[rd] \ar@{->>}[rr]  && J_1 \sqcup J_2 \ar@{->>}^s[r] & \tb{2} \\
	& J \ar@{->>}[ru] \ar@{->>}[r] & K \ar@{->>}[ru]
}	
	\]
	where $s$ is the fixed morphism that sends $J_1$ to 1 and $J_2$ to 2. Now we fix $(q: J \surj K) \in \ETOp^\tx{top}_{p \sqcup q /}$, let $J \surj J_1 \sqcup J_2 \surj^s \tb{2}$ be denoted $\pi$, and define $L$ to be the quotient by the following equivalence relation: for all $x, y \in J$, $x \sim y$ iff $q(x) = q(y), \pi(x) = \pi(y)$. It is then clear that $J \surj L$ maps to $q$ and is final within the fiber of $q$ in $\ETOp^\tx{top}_{p/} \times \ETOp^\tx{top}_{q/}$.
\end{proof}

Since we work in the stable setting, that the necessary $\ETOp_{p/}$-limits in Lemma \ref{lemma:sm-right-kan-extension} (for $\mc{K}$ the class of \emph{finite} simplicial sets) in the coCartesian fibration exist, and can be computed by first sending all entries to $\Gamma(\mc{C}, U_\Ran(p))$ via the right adjoints, then compute the limit there. The symmetric monoidal structure is the evident one. It then follows that the remaining conditions of Lemma~\ref{lemma:sm-right-kan-extension} are met, and we are allowed to monoidally Kan extend and produce an $\baseOp$-monoidal functor $\bar{F} := \tx{RKE}_{\baseOp}(\widetilde{F})$. Evidently $F \mapsto \bar{F}$ is the right adjoint to the forgetful functor $\tx{Fun}^{\tensor}(\ECOp, \msc{C}) \to \tx{Fun}^{\tensor, \tx{lax}}(\ECOp, \msc{C})$.

\begin{lemma}
\label{lemma:coCartesian-check}
	If $F$ sends all arrows to coCartesian arrows, then so does $\bar{F}$.
\end{lemma}

We demonstrate the case for $\fSetnu$; the $\CMnu$ case is completely analogous.

\begin{proof}
For $p_1 \to p_2 \in \ECOp$, let $\varphi_{p_1 \to p_2}: \Gamma(U_\Ran(p_1), \mc{C}) \to \Gamma(U_\Ran(p_2), \mc{C})$ be the structural morphism. Let $\pi: \ETOp \to \ECOp$ be the projection. Concretely, we want the statement that for any $p_1 \to p_2 \in \ECOp$ and $\mc{F}$ an $\baseOp$-weak factorization algebra, the natural morphism
\[\varphi_{p_1 \to p_2}(\limit_{q_1 \in \ETOp_{p_1/}} \varphi_{p_1 \to \pi(q_1)}^R(\mc{F}(q_1))) \to \limit_{q_2 \in \ETOp_{p_2/}} \varphi_{p_2 \to \pi(q_2)}^R(\mc{F}(q_2))\]
is an isomorphism. Using assumption on $\varphi$ we can rewrite the morphism as
\[\limit_{q_1 \in \ETOp_{p_1/}} \varphi_{p_1 \to p_2} \circ \varphi_{p_1 \to \pi(q_1)}^R(\mc{F}(q_1)) \to \limit_{q_2 \in \ETOp_{p_2/}} \varphi_{p_2 \to \pi(q_2)}^R(\mc{F}(q_2))\]
and the map is restriction to a subdiagram $Q$ of those $q$ such that $p_1 \to q$ factors as $p_1 \to p_2 \to q$, followed by the morphism obtained via
\[[\varphi_{p_1 \to p_2} \circ \varphi^R_{p_1 \to \pi(q)}] \simeq [\varphi_{p_1 \to p_2} \circ \varphi^R_{p_1 \to p_2} \circ \varphi^R_{p_2 \to \pi(q)}] \xrightarrow{\tx{counit}} [\varphi^R_{p_2 \to \pi(q)}]\]
Write $p_1 = I_1 \surj J_1$, $p_2 = I_2 \surj J_2$, and running induction on $|I_1 - I_2| + |J_2 - J_1|$ we reduce to the following two cases:

\paragraph{Case 1} $I_1 = I_2 =: I, |J_2| = |J_1| + 1$: this is the open embedding case. We can replace $\ETOp$ with $\ETOp^\tx{top}$ on both sides. Without loss of generality, we can assume the map $p_1 \to p_2$ is induced by the map $J_2 \surj J_1: \{1, 2\} \mapsto 1, \ldots, |J_1| \mapsto |J_1| - 1, |J_1| + 1 \mapsto |J_1|$. 
For any term $I \to J \to K$ included in the subdiagram $Q$, the map is given by (we omit the $I$-part)
\[[\varphi_{(J_1) \to (J_2)} \circ \varphi_{(J_1) \to (J_2)}^R \circ \varphi_{(J_2) \to (J)}^R] \to [\varphi_{(J_2) \to (J)}^R]\]
which is an isomorphism by assumption. It remains to check that restriction to the smaller diagram does nothing. Note that $\ECOp^\tx{top}_{p_1/}$ as an under-category admits a coproduct structure $\sqcap$, such that
\[U_\Ran(p \sqcap (I \to J_2)) \simeq U_\Ran(p) \cap U_\Ran(I \to J_2);\]
we define $\Psi(I \to J) = (I \to J) \sqcap (I \to J_2)$ and extend it to $\ETOp^\tx{top}_{p_1/}$ via identity. By base change, it suffices to show that $\ETOp^\tx{top}_{p_2/}$ is left cofinal within $\Psi(\ETOp^\tx{top}_{p_1/})$, but this is clear: after inverting isomorphisms, $\Psi(\ETOp^\tx{top}_{p_1/})$ indeed retracts onto $\ETOp^\tx{top}_{p_2/}$. %

\paragraph{Case 2} $|I_1| = |I_2| + 1, J_1 = J_2 =: J$: this is the diagonal embedding case. To avoid considering non-continuous functors we again replace $\ETOp$ with $\ETOp^\tx{top}$ on both sides. Without loss of generality we assume the mapping on $I$ is $\{1, 2\} \mapsto 1, \ldots, |I_2| + 1 \mapsto |I_2|$; in other words, the corresponding closed subscheme is $\{x_1 = x_2\}$; let $\Delta$ denote this closed embedding.
We claim that, upon taking $\varphi_{p_1 \to p_2}$, those nodes in $\mc{F}(\ETOp^\tx{top}_{p_1 / })$ that do not become zero form a diagram isomorphic to $\mc{F}(\ETOp^\tx{top}_{p_2 / })$. By (weak) factorization, it suffices to note that the full subcategory of arrows $I_1 \to J'$ in $\ECOp^\tx{top}_{p_1/}$ that admit a dotted arrow (which is by definition unique) like follows:
\[\xymatrix{
I_1 \ar[r] \ar[rd] \ar[d] & I_2 \ar[dl] \ar@{-->}[d] \\
J & J' \ar[l]
}\]
is isomorphic to $\ECOp^\tx{top}_{p_2/}$ (the map is by sending it to the dotted arrow).
\end{proof}

To summarize, the forgetful functor
\[\baseOp\tx{-FactAlg}^{\CoMark}_{\nonun}(\mc{C}) \to \baseOp\tx{-FactAlg}^{\CoMark}_{\nonun, \LWeakMark}(\mc{C})\]
admits a right adjoint $\ms{strc}_{\baseOp}$, which can be explicitly computed as
\[\ms{strc}_{\baseOp}(\mc{F})(p) = \limit_{q \in \ETOp^\tx{top}_{p/}} \varphi_{p \to \pi(q)}^R \left(\mc{F}(q)\right)\]
	for each $p \in \ECOp$, where $\pi: \ETOp \to \ECOp$ is the projection, $\varphi_{p \to \pi(q)}^R$ is right adjoint to the structure morphism $\mc{C}(p) \to \mc{C}(\pi(q))$.

\subsubsection{The Unital Case}
Now we carry the above construction to the unital situation. We will use subscript $\bullet_{\un}$ to denote the unital version of the various gadgets constructed above.

Using the partition description from Proposition \ref{prop:fact-cat-alg-partition-description-unital}, the statement of Corollary~\ref{cor:algebra-partition-non-unital} is replaced by the following:

\begin{corollary}
\label{cor:algebra-partition-unital}
	The category of weak unital $\baseOp_\un$-monoidal factorization algebra in $\mc{C}_\un$ is equivalent to the category of lax $\baseOp_\un$-monoidal sections of $\mathscr{C}_\un \to \ECOp_{\un}$ sending all arrows in $\ECOp \subseteq \ECOp_\un$ and $\{[0] \to \bullet\} \subseteq \ECOp_\un$ to coCartesian arrows\footnote{Note that the latter is necessarily of color $\mf{a}$ in the $\CM$ case.}. The subcategory of factorization objects corresponds to the $\baseOp$-monoidal sections.
\end{corollary}

In the unital case, $(\ETOp_\un)^\tx{top}_{p/}$ is still initial within $(\ETOp_\un)_{p/}$ for every $p \in \ECOp_\un$, but the cofinality statement no longer holds. Therefore we will manually check that, for every $a: I \to J \in \baseOp$, the map $\bar{F}^{J} \circ a_! \to a_! \circ \bar{F}^I$ is an isomorphism. The inert and injective cases are evident, so it suffices to prove that both
\begin{equation}
\label{eq:res-unital-eq1}
\underset{q \in (\ETOp_\un)^\tx{top}_{p_1 \sqcup p_2 /}}{\limit} \varphi_{p_1 \sqcup p_2 \to \pi(q)}^R(\mc{F}_\un(q)) \to
\underset{(q_1, q_2) \in (\ETOp_\un)^\tx{top}_{p_1 /} \times (\ETOp_\un)^\tx{top}_{p_2 /}}{\limit} \varphi_{p_1 \to \pi(q_1)}^R (\mc{F}_\un(q_1)) \tensor \varphi_{p_2 \to \pi(q_2)}^R(\mc{F}_\un(q_2))
\end{equation}
and
\begin{multline}
\label{eq:res-unital-eq2}
\underset{q_1 \in (\ETOp_\un)^\tx{top}_{p_1 /}}{\limit} \varphi_{p_1 \to \pi(q_1)}^R (\mc{F}_\un(q_1)) \bigotimes \underset{q_2 \in (\ETOp_\un)^\tx{top}_{p_2 /}}{\limit} \varphi_{p_2 \to \pi(q_2)}^R(\mc{F}_\un(q_2)) \to \\
\underset{(q_1, q_2) \in (\ETOp_\un)^\tx{top}_{p_1 /} \times (\ETOp_\un)^\tx{top}_{p_2 /}}{\limit} \varphi_{p_1 \to \pi(q_1)}^R (\mc{F}_\un(q_1)) \tensor \varphi_{p_2 \to \pi(q_2)}^R(\mc{F}_\un(q_2))
\end{multline}
are isomorphisms for $p_1, p_2 \in \ECOp_\un$.

Note that every $p \in \ECOp_\un$ canonically decomposes as $p \simeq p_s \sqcup p_e$ where $p_s$ is surjective and $p_e$ has $I$-component being $[0]$; similarly any $q \in \ETOp_\un$ canonically decomposes as $q \simeq q_s \sqcup q_e \sqcup q_v$ where $q_s$ is $I \surj J \surj K$, $q_e$ is $[0] \to J' \surj K'$ and $q_v$ is $[0] \to [0] \to K''$. These decompositions are compatible with disjoint unions.

\begin{remark}
$p_s$ and $q_s$ do not quite lie in $\ECOp$ and $\ETOp$ because they could be $[0] \to [0]$ and $[0] \to [0] \to [0]$ respectively. For ease of notation we will abuse notation and add in these monoidal units for $\ECOp$ and $\ETOp$ for the rest of this subsection. We extend the definition of $\mc{F}$ to $\mc{F}([0] \to [0] \to [0]) := \tx{unit}_{C([0])}$.

Note that $\mc{F}_\un(q) \simeq \mc{F}_\un(q_s) \simeq \mc{F}(q_s)$, because $\mc{F}_\un([0] \to \bullet \to \bullet)$ is the monoidal unit even for weak factorization objects.
\end{remark}

We now reduce to the non-unital case by proving that for any $p$ we have
\begin{equation}
\label{eq:res-unital-to-non-unital}
\underset{q \in (\ETOp_\un)^\tx{top}_{p /}}{\limit} \varphi_{p \to \pi(q)}^R(\mc{F}_\un(q)) \simeq \underset{r \in (\ETOp)^\tx{top}_{p_s /}}{\limit} \varphi_{p_s \to \pi(r)}^R(\mc{F}(r));
\end{equation}
Note that this would imply the relevant limit is actually finite so commutes with tensor product, which in turn implies
\[\underset{(q_1, q_2) \in (\ETOp_\un)^\tx{top}_{p_1 /} \times (\ETOp_\un)^\tx{top}_{p_2 /}}{\limit} \varphi_{p_1 \to \pi(q_1)}^R (\mc{F}_\un(q_1)) \tensor \varphi_{p_2 \to \pi(q_2)}^R(\mc{F}_\un(q_2))\]
\[ \simeq \underset{(r_1, r_2) \in (\ETOp)^\tx{top}_{(p_1)_s /} \times (\ETOp)^\tx{top}_{(p_2)_s /}}{\limit} \varphi_{(p_1)_s \to \pi(r_1)}^R (\mc{F}(r_1)) \tensor \varphi_{(p_2)_s \to \pi(r_2)}^R(\mc{F}(r_2))\]
and thus both Equation \eqref{eq:res-unital-eq1} and \eqref{eq:res-unital-eq2} above are equivalences. The case $p \simeq [0] \to \bullet$ is evident, so we can assume $p$ is a map of non-empty finite sets. Note that for any $p \to \pi(q)$ (which means $p_s \to \pi(q_s)$ within $\ECOp$), the commutative square
\[\xymatrix{
C_\un(p_s) \tensor C_\un(p_e) \ar^{\tensor}[r] \ar^{\varphi_{p_s \to \pi(q_s)} \tensor \tx{id}}[d] & C_\un(p) \ar^{\varphi_{p \to \pi(q)}}[d] \\
C_\un(\pi(q_s)) \tensor C_\un(\pi(q_e \sqcup q_v)) \ar^(0.65){\tensor}[r] & C_\un(\pi(q))
}\]
is right adjointable vertically since the horizontal maps are isomorphisms. Thus we can rewrite LHS as
\[\underset{q \in (\ETOp_\un)^\tx{top}_{p /}}{\limit} \varphi_{p \to \pi(q)}^R(\mc{F}_\un(q)) \simeq \underset{q \in (\ETOp_\un)^\tx{top}_{p /}}{\limit} \varphi_{p_s \to \pi(q_s)}^R(\mc{F}(q_s));\]
The discussion above means we may localize $(\ETOp_\un)^\tx{top}_{p /}$ with respect to maps that are injective on the $J$-component or the $K$-component without changing the limit, so it suffices to prove that for every $p$ we have a localization
\[(\ETOp_\un)^\tx{top}_{p/}[J, K\tx{ injective maps}^{-1}] \simeq \ETOp^{\tx{top}}_{p_s/}\]
that sends each $q$ to $q_s$. This now follows from the following zig-zag diagram
\[\xymatrix{
I \ar^{\simeq}[d] \ar^{f}[r] & J \ar^{g}[r] & K \ar^{\simeq}[d] \\
I \ar@{->>}[r] & \im(f) \ar@{^{(}->}[u] \ar^{\simeq}[d] \ar^{g|_{\im(f)}}[r] & K \\
I \ar^{\simeq}[u] \ar@{->>}[r] & \im(f) \ar@{->>}[r] & \im(g \circ f) \ar@{^{(}->}[u]
}\]
This finishes the proof of Equation \eqref{eq:res-unital-to-non-unital}. The condition of sending coCartesian arrows now follows Lemma \ref{lemma:coCartesian-check} above (the $0 \to \bullet$ case being evident). To summarize, the strictification functor also exists in the unital situation, and is computed on each $U_{\Ran_{\un}}(p)$ by the same formula as in the non-unital case.

Now we are ready to prove the original statement.

\begin{proof}[Proof of Proposition \ref{prop:restriction-functor}]
The following discussion applies to both the unital and the non-unital case. Consider the tautological inclusion
\[\oblv_{\tx{mod}}^{s \to w}: A\tx{-FactMod}^{\CoMark}_{\PartitionMark}(C_M) \to A\tx{-FactMod}^{\CoMark}_{\LWeakMark, \PartitionMark}(C_M)\]
By considerations above, this functor admits a right adjoint, which we denote by $\oblv_{\tx{mod}}^{s \to w, R}$ and is computed by $\ms{strc}_{\CM}$. Since weak modules are actual modules in the correspondence category (see Remark \ref{rmk:comparison-with-chiralcats}, or use the interpretation of lax monoidal functors as algebra objects under Day convolution), by \cite[4.2.3.1]{HA}, for any morphism $f: A \to B$ in $\tx{FactAlg}^{\CoMark}_{\LWeakMark, \PartitionMark}(C)$, we have a functor
\[\tx{res}^w_{f}: B\tx{-FactMod}^{\CoMark}_{\LWeakMark, \PartitionMark}(C_M) \to A\tx{-FactMod}^{\CoMark}_{\LWeakMark, \PartitionMark}(C_M)\]
which commutes with all small limits and colimits.
\begin{definition}
We \emph{define} the restriction functor of (strong) factorization modules as
\[\tx{res}_{f}:= B\tx{-FactMod}(C_M) \xrightarrow{\oblv_{\tx{mod}}^{s \to w}} B\tx{-FactMod}\tx{-FactMod}^{\CoMark}_{\LWeakMark, \PartitionMark}(C_M)\]
\[\xrightarrow{\tx{res}^w_{f}} A\tx{-FactMod}\tx{-FactMod}^{\CoMark}_{\LWeakMark, \PartitionMark}(C_M) \xrightarrow{\oblv_{\tx{mod}}^{s \to w, R}} A\tx{-FactMod}(C_M).\]
\end{definition}
Using the description obtained in previous section, we get an explicit expression
\[\tx{res}_{f}(M)(p) = \limit_{q \in \tx{Env}_{\CM_\nonun^{\tensor}}(\tx{Tw}(\CM_\nonun)^\op)^\tx{top}_{p/}} \varphi_{p \to \pi(q)}^R (A \boxtimes M)(q),\]
from which the second property is also evident.
\end{proof}

\subsection{Commutative Case}
We now consider the counterpart of $\mb{E}_\infty$-objects in the factorization world.

On the space level, we have the notion of \emph{commutative} factorization spaces. Let $\Ran_\nonun^\tx{comm}$ denote the oplax functor corresponding to $\tx{Ran}_\nonun^\tx{comm}$ described above. There is then a natural transformation $\eta: \Ran_\nonun \Rightarrow \Ran_\nonun^\tx{comm}$. It is easy to see that pullback along $\eta$ gives a functor
\[\eta^*: \tx{Fun}^{\tx{oplax}, \tensor}_\nonun(\fSetnu, \tx{PreStk})_{/\Ran_\nonun^\tx{comm}} \to \tx{Fun}^{\tx{oplax}, \tensor}_\nonun(\fSetnu, \tx{PreStk})_{/\Ran_\nonun}.\]

\begin{definition}
Let $\tx{CommFactSpace}_\nonun$ denote the full subcategory of
\[\tx{Fun}^{\tx{oplax}, \tensor}_\nonun(\fSetnu, \tx{PreStk})_{/\Ran_\nonun^\tx{comm}}\]
consisting of object $F$ satisfying condition L and
\begin{itemize}
\item[(Condition C)] For every composable pair of morphisms $I \to J \to K$ in $\fSetnu$, the commutative square
\[\xymatrix{
	\bigotimes_{k = 1}^{K} F(I_k) \ar[r] \ar[d] & \bigotimes_{k = 1}^{K} F(J_k) \ar[d] \\
	\bigotimes_{j = 1}^{J} F(I_j) \ar[r] & \bigotimes_{j = 1}^{J} F(j)
}\]
	is Cartesian;
\end{itemize}
but such that $\eta^*(F)$ satisfies both condition L and R. An element of this category is called a non-unital \emph{commutative} factorization space.
\end{definition}

\begin{variant}
The category $\tx{CommFactSpace}_{\un}$ of unital commutative factorization spaces can be analogously defined.
\end{variant}

\cite{raskin2015chiral} introduced a functor which attaches to every symmetric monoidal category a commutative factorization category. This can be written, in our notation, as
\[\tx{Loc}: \tx{CAlg}_\nonun(\DGCat) \to \tx{CommFactCat}_\nonun\]
\[\tx{Loc}: \tx{CAlg}(\DGCat) \to \tx{CommFactCat}_\tx{un}\]
A key feature is that the restriction of $\tx{Loc}(\ms{C})$ to $X_\dR$ is canonically isomorphic to $\ms{C} \tensor \DMod(X)$. The following is well-known (see \cite[8.10.1]{raskin2016chiral} \cite[5.6.4]{gaitsgory2019weil}):
\begin{proposition}
\label{prop:comm-fact-main-diagonal-same}
	For each $\ms{C} \in \tx{CAlg}(\DGCat)$, the diagonal restriction upgrades to an equivalence
	\[\tx{Fact}: \tx{CAlg}_\nonun(\ms{C} \tensor \DMod(X)) \simeq \tx{CommFactAlg}_\nonun(\tx{Loc}(\ms{C})): \tx{res}\]
	which further restricts to an equivalence of full subcategories
	\[\tx{Fact}: \tx{CAlg}(\ms{C} \tensor \DMod(X)) \simeq \tx{CommFactAlg}_{\un}(\tx{Loc}(\ms{C})): \tx{res}\]
	In particular, being a unital \emph{commutative} factorization algebra is not extra data.
\end{proposition}

\begin{example}
Taking $\ms{C} = \tx{Vect}$, this is the equivalence between commutative $D$-algebras and commutative factorization algebras from \cite{beilinson2004chiral}.
\end{example}

\begin{remark}
Discussion in Section \ref{sect:factvse2} can be seen as a generalization of this statement to the $\mb{E}_2$ case, albeit in the topological setting.
\end{remark}

\subsection{Configuration Space versus Affine Grassmannian}
\label{sect:conf-compare-gr}

Recall the configuration space $\tx{Conf}(X)$ and the the $\crho$-shifted affine Grassmannian $[\tx{Gr}_{\check{T}}^{\omega^{\crho}}]_{\Ran_\dR}$ introduced in Section \ref{sect:conf-space-in-minimal}. Consider the subfunctors $[\tx{Gr}_{\check{T}}^{\omega^{\crho}, \le 0}]_{\Ran_{\un, \dR}}$ and $[\tx{Gr}_{\check{T}}^{\omega^{\crho}, < 0}]_{\Ran_{\un, \dR}}$ of bundles satisfying \emph{regularity} and \emph{non-redundancy}, as defined in \cite[4.6.2]{gaitsgory2019metaplectic}.

As in \cite[4.6]{gaitsgory2019metaplectic}, we have a diagram as below
\[\xymatrix{
& [\tx{Gr}_{\check{T}}^{\omega^{\crho}, \le 0}]_{\Ran_{\un, \dR}} \ar@^{^(->}[r] & [\tx{Gr}_{\check{T}}^{\omega^{\crho}}]_{\Ran_{\un, \dR}} \\
[\tx{Gr}_{\check{T}}^{\omega^{\crho}, <0}]_{\Ran_\dR} \ar^{(\star)}[d] \ar@^{^(->}[r] & [\tx{Gr}_{\check{T}}^{\omega^{\crho}, \le 0}]_{\Ran_\dR} \ar@^{^(->}[r] \ar[u] & [\tx{Gr}^{\omega^{\crho}}_{\check{T}}]_{\Ran_\dR} \ar[u] \\
\tx{Conf}
}\]
where the arrow from second row to first is the inclusion of the non-unital version to the unital version. We note that $\tx{Gr}_{\check{T}}^{\omega^{\crho}, \le 0}$ admits a unital structure, but $\tx{Gr}_{\check{T}}^{\omega^{\crho}, < 0}$ does not.
The map denoted $(\star)$ is an isomorphism after sheafifying with respect to finite surjective morphisms; consequently it induces an isomorphism on both de Rham tame $\mb{G}_m$-gerbes and analytic $\mb{G}_m$-gerbes, and induces an isomorphism on twisted D-modules as well as constructible sheaves.

Let $\mathscr{G}$ be a $\mb{G}_m$-gerbe on $[\tx{Gr}_{\check{T}}^{\omega^{\crho}}]_{\Ran_{\un, \dR}}$, we will let the same symbol denote the pullback along the diagram above to $\tx{Conf}$. A similar discussion holds for $\tx{Conf}_{\infty \cdot \modulePoints}$ as well.

\begin{convention}
All discussions below in this subsection apply equally to the D-module and the constructible settings; for simplicity we will only write out the D-module case.
\end{convention}

As in \cite[5.5]{gaitsgory2019metaplectic}, pullback along $(\star)$ then pushforward will then give rise to an equivalence of categories
\[\tx{Conv}_{\tx{Conf}}^{\tx{Alg}}: \tx{FactAlg}_{\nonun}(\DMod_{\mathscr{G}}(\tx{Conf})) \simeq \tx{FactAlg}_\nonun(\DMod_{\mathscr{G}}([\tx{Gr}_{\check{T}}^{\omega^{\crho}, <0}]_{\Ran_\dR}))\]
and, for every factorization algebra $A$ on the configuration space, an equivalence of categories
\[\tx{Conv}_{\tx{Conf}}^{A}: A\tx{-FactMod}_{\nonun}(\DMod_{\mathscr{G}}(\tx{Conf}_{\infty \cdot \modulePoints})) \simeq \tx{Conv}_{\tx{Conf}}^{\tx{Alg}}(A)\tx{-FactMod}_\nonun(\DMod_{\mathscr{G}}([\tx{Gr}_{\check{T}}^{\omega^{\crho}, <0, \infty \cdot \modulePoints}]_{\Ran_{\dR, \modulePoints}}))\]
\[\simeq \tx{Conv}_{\tx{Conf}}^{\tx{Alg}}(A)\tx{-FactMod}_\nonun(\DMod_{\mathscr{G}}([\tx{Gr}_{\check{T}}^{\omega^{\crho}}]_{\Ran_{\dR, \modulePoints}})),\]
which immediately follows from factorization. We likewise have a fully faithful embedding
\[\tx{FactAlg}_\un(\DMod_{\mathscr{G}}([\tx{Gr}_{\check{T}}^{\omega^{\crho}, \le 0}]_{\Ran_{\un, \dR}})) \subseteq \tx{FactAlg}_\un(\DMod_{\mathscr{G}}([\tx{Gr}_{\check{T}}^{\omega^{\crho}}]_{\Ran_{\un, \dR}}))\]
whose essential image consisting of those unital factorization algebras whose underlying sheaf is suppported on the sublocus; for a fixed unital factorization algebra $A$ that lies in the former, we have an equivalence of categories
\[A\tx{-FactMod}_{\un}(\DMod_{\mathscr{G}}([\tx{Gr}_{\check{T}}^{\omega^{\crho}, \le 0, \infty \cdot \modulePoints}]_{\Ran_{\dR, \modulePoints}})) \simeq A\tx{-FactMod}_\un(\DMod_{\mathscr{G}}([\tx{Gr}_{\check{T}}^{\omega^{\crho}}]_{\Ran_{\dR, \modulePoints}})).\]
By pulling back along
\[[\tx{Gr}_{\check{T}}^{\omega^{\crho}, <0}]_{\Ran_\dR} \to [\tx{Gr}_{\check{T}}^{\omega^{\crho}, \le 0}]_{\Ran_\dR} \to [\tx{Gr}_{\check{T}}^{\omega^{\crho}, \le 0}]_{\Ran_{\un, \dR}}\]
we have a functor
\[\tx{Red}: \tx{FactAlg}_\un(\DMod_{\mathscr{G}}([\tx{Gr}_{\check{T}}^{\omega^{\crho}, \le 0}]_{\Ran_{\un, \dR}})) \to \tx{FactAlg}_\nonun(\DMod_{\mathscr{G}}([\tx{Gr}_{\check{T}}^{\omega^{\crho}, <0}]_{\Ran_\dR}))\]
and for every unital factorization algebra $A$ in the former, a functor
\[\tx{Red}_A: A\tx{-FactMod}_\un(\DMod_{\mathscr{G}}([\tx{Gr}_{\check{T}}^{\omega^{\crho}}]_{\Ran_{\dR, \modulePoints}})) \simeq A\tx{-FactMod}_{\un}(\DMod_{\mathscr{G}}([\tx{Gr}_{\check{T}}^{\omega^{\crho}, \le 0, \infty \cdot \modulePoints}]_{\Ran_{\dR, \modulePoints}}))\]
\[\xrightarrow{!\tx{-pull}} \tx{Red}(A)\tx{-FactMod}_\nonun(\DMod_{\mathscr{G}}([\tx{Gr}_{\check{T}}^{\omega^{\crho}, <0, \infty \cdot \modulePoints}]_{\Ran_{\dR, \modulePoints}})).\]

Let us set $A_{\tx{red}, \tx{Conf}} := (\tx{Conv}_{\tx{Conf}}^{\tx{Alg}})^{-1} \circ \tx{Red}(A)$, so that for every unital factorization algebra $A$ that lies in $\tx{FactAlg}_\un(\DMod_{\mathscr{G}}([\tx{Gr}_{\check{T}}^{\omega^{\crho}, \le 0}]_{\Ran_{\un, \dR}}))$, we have a well-defined functor\footnote{We note (without proving) that for $A \simeq \TheFactAlg$, this is actually an equivalence of categories.}
\[\tx{Red}_{\tx{Conf}}: A\tx{-FactMod}_\un(\DMod_{\mathscr{G}}([\tx{Gr}_{\check{T}}^{\omega^{\crho}}]_{\Ran_{\dR, \modulePoints}})) \to A_{\tx{red}, \tx{Conf}}\tx{-FactMod}_{\nonun}(\DMod_{\mathscr{G}}(\tx{Conf}_{\infty \cdot \modulePoints})).\]

\section{Appendix C: Ind-coherent Sheaves}

\label{sect:appendix-indcoh}

\begin{convention}
	In this chapter, all occurrences of $\tx{PreStk}$ and the term ``prestack'' should be understood as \emph{convergent prestack}. For a prestack $\mc{Y}$ to be convergent means that, for every affine $S$, we require the natural map
	\[\mc{Y}(S) \to \limit_{n} \mc{Y}(\tau^{\le n}(S))\]
	(where $\tau^{\le n}$ denotes the Postnikov truncation of $S$) to be an isomorphism.
	
	We note that being convergent is a condition that is satisfied by almost all ``geometric'' objects (e.g. all Artin stacks). In particular, all prestacks we consider in this paper are convergent.
\end{convention}

\subsection{Ind-coherent sheaves on qcqs schemes}

In this subsection, we study the $\IndCoh_*$-theory on eventually coconnective qcqs schemes.

\begin{definition} Following \cite[$\S$ 6.4]{raskin2020homological}, for $S\in \,^{>-\infty}\on{Sch}_{{\on{qcqs}}}$, define $\Coh(S)\subset \QCoh(S)$ to be the full subcategory of bounded \alert{almost perfect} object, i.e., those objects $\mCF\in \QCoh(S)^+$ such that for any integer $N$, $\mCF\in \QCoh(S)^{\ge -N}$ implies $\mCF$ is a compact object in $\QCoh(S)^{\ge -N}$. Define
$$ \IndCoh_*(S):= \on{Ind}(\Coh(S)).$$
Let
$$\Psi:=\Psi_S: \IndCoh_*(S) \to \QCoh(S)$$
be the continuous functor obtained by ind-extending from $\Coh(S) \to \QCoh(S)$.
\end{definition}

\begin{remark} In \cite{raskin2020homological}, $\IndCoh_*(S)$ is denoted by $\IndCoh^*(S)$. 
\end{remark}

\begin{remark} In the above definition, we can replace ``any integer $N$'' by ``any large enough integer $N$''.
\end{remark}

\begin{remark} It is easy to see $\Coh(S)$ is Karoubi complete, hence $\Coh(S)$ is the category of compact objects in $\IndCoh(S)$.
\end{remark}

\begin{remark} It is easy to see $\Coh(S)\subset \QCoh(S)$ is closed under finite colimits. Hence $\IndCoh_*(S)$ is a cocomplete DG-category.
\end{remark}

\begin{remark} If $S$ is almost finite type, i.e., $S \in \,^{>-\infty} \on{Sch}_{\on{aft}}$ (which is equivalent to $S\in\,^{\ge -m}\on{Sch}_{\on{ft}}$ for some $m$), the above definition of $\on{Coh}(S)$ coincides with that in \cite{gaitsgory2013ind}. Hence we write $\IndCoh(S)$ for $\IndCoh_*(S)$ in this case.

\end{remark}

\begin{construction} For $S\in \,^{>-\infty}\on{Sch}_{\on{qcqs}}$, there exists a compactly generated t-structure on $\IndCoh_*(S)$ such that $\IndCoh_*(S)^{\le 0}$ is generated under colimits and extensions by $\on{Coh(S)}\cap \QCoh(S)^{\le 0}$.
\end{construction}

\begin{remark} 
In particular, the above t-structure on $\IndCoh_*(S)$ is right separated and compatible with filtered colimits.
\end{remark}

\begin{remark} In the above definition, we can replace ``generated under colimits and extensions'' by ``generated under filtered colimits''.
\end{remark}

\begin{lemma} (\cite[Lemma 6.4.1]{raskin2020homological}) \label{lem-indcoh*-vs-qcoh-sch}
For $S\in \,^{>-\infty}\on{Sch}_{\on{qcqs}}$, the functor $\Psi: \IndCoh_*(S) \to \QCoh(S)$ is t-exact and induces an equivalence $\IndCoh(S)^+ \simeq \QCoh(S)^+$.
\end{lemma}

\begin{construction} For $f: S\to T$ in $ ^{>-\infty}\on{Sch}_{\on{qcqs}}$, the functor $f_*:\QCoh(S)\to \QCoh(T)$ is left t-exact. Consider the functor
$$ \Coh(S) \to \QCoh(S)^+ \xrightarrow{f_*}  \QCoh(T)^+ \xrightarrow{\Psi^{-1}} \IndCoh(T)^+ \to \IndCoh(T).$$
By ind-extending, we obtain a functor\footnote{The functor below was denoted by $f_*^{\IndCoh}$ in \cite{raskin2020homological}.}
$$ f_*: \IndCoh_*(S) \to \IndCoh_*(T)$$
and a commutative diagram
\begin{equation}
\label{eqn-*-push-Psi-square-sch}
\xymatrix{
	\IndCoh_*(S) \ar[r]^-{ f_* } \ar[d]^-{\Psi_S} &
	\IndCoh_*(T) \ar[d]^-{\Psi_T} \\
	\QCoh(S) \ar[r]^-{f_*} &
	\QCoh(T).
}
\end{equation}
\end{construction}

\begin{warning} The functor $ f_*: \IndCoh_*(S) \to \IndCoh_*(T)$ is \alert{not} necessarily left t-exact\footnote{This was falsely claimed in \cite{raskin2020homological}.}. However, we have the following lemma.
\end{warning}

\begin{lemma} \label{lem-f_*-left-t-exact}
For $f: S\to T$ in $ ^{>-\infty}\on{Sch}_{\on{qcqs}}$, suppose $\Coh(S) \subset \QCoh(S)$ is closed under truncations, then $f_*: \IndCoh_*(S) \to \IndCoh_*(T)$ is left t-exact.
\end{lemma}

\proof Let $\mCF = \colim \mCF_i$ be an object in $\IndCoh_*(S)^{\ge 0}$ with $\mCF_i\in \on{Coh}(S)$. Then $\mCF \simeq \colim  \tau^{\ge 0}(\mCF_i)$ and $\tau^{\ge 0}(\mCF_i)\in  \on{Coh}(S)\cap \IndCoh_*(S)^{\ge 0}$. It follows that $f_*( \mCF ) \simeq \colim f_*(  \tau^{\ge 0}(\mCF_i) )$ is contained in $\IndCoh_*(S)^{\ge 0}$.

\qed

\begin{lemma} \label{lemma-f_*-right-t-exact}
For $f: S\to T$ in $ ^{>-\infty}\on{Sch}_{\on{qcqs}}$, suppose $[n]\circ f_*: \QCoh(S)\to \QCoh(T)$ is \alert{right} t-exact, then $[n]\circ f_*: \IndCoh(S)\to \IndCoh(T)$ is also right t-exact.
\end{lemma}

\proof Follows from Lemma \ref{lem-indcoh*-vs-qcoh-sch}.

\qed

\begin{lemma} (C.f. \cite[Lemma 6.6.1]{raskin2020homological}) Suppose $f: S\to T$ in $ ^{>-\infty}\on{Sch}_{\on{qcqs}}$ is of bounded Tor dimension, then $f^*:\QCoh(T) \to \QCoh(S)$ sends $\Coh(T)$ into $\Coh(S)$. By ind-extending, we obtain a functor $f^{*}: \IndCoh_*(T) \to \IndCoh_*(S)$, which is left adjoint to $f_*: \IndCoh_*(S)\to \IndCoh_*(T)$. Also, the commutative square (\ref{eqn-*-push-Psi-square-sch}) is left adjointable along the horizontal direction, i.e., we have a commutative diagram
$$
\xymatrix{
	\IndCoh_*(T) \ar[r]^-{ f^* } \ar[d]^-{\Psi_S} &
	\IndCoh_*(S) \ar[d]^-{\Psi_T} \\
	\QCoh(T) \ar[r]^-{f^*} &
	\QCoh(S).
}
$$
\end{lemma}

\proof  \cite{raskin2020homological} proved the case when $f$ is flat. The proof can be generalized to our case.

\qed

\begin{warning} 
The functor $f^{*}: \IndCoh_*(T) \to \IndCoh_*(S)$ is \alert{not} necessarily left t-exact up to a shift\footnote{This was falsely claimed in \cite{raskin2020homological}.}. However, we have the following lemma.
\end{warning}

\begin{lemma} \label{lem-*-pull-left-t-exact-sch}
Suppose $f: S\to T$ in $ ^{>-\infty}\on{Sch}_{\on{qcqs}}$ is of Tor dimension $\le n$, and suppose $\Coh(T) \subset \QCoh(T)$ is closed under truncations, then $[-n]\circ f^{*}: \IndCoh_*(T) \to \IndCoh_*(S)$ is left t-exact.
\end{lemma}

\proof Similar to Lemma \ref{lem-f_*-left-t-exact}.

\qed

\begin{lemma} \label{lem-*-pull-right-t-exact-sch}
Suppose $f: S\to T$ in $ ^{>-\infty}\on{Sch}_{\on{qcqs}}$ is of bounded Tor dimension, then $f^{*}:  \IndCoh_*(T) \to \IndCoh_*(S)$ is right t-exact. In particular, $f_*:  \IndCoh_*(S) \to \IndCoh_*(T)$ is left t-exact.
\end{lemma}

\proof Follows from Lemma \ref{lem-indcoh*-vs-qcoh-sch}.

\qed

\begin{lemma} (C.f. \cite[Lemma 6.6.2, Remark 6.6.3]{raskin2020homological}) 
\label{lem-*-push-*-pull-base-change-sch}
Let $g:T'\to T$ be a map of bounded Tor dimension and $f:S \to T$ be any map. Consider the following Cartesian diagram
$$
\xymatrix{
	S' \ar[r]^-{\varphi} \ar[d]^-\phi & T'\ar[d]^-g \\
	S \ar[r]^-f & T
}
$$
in $^{>-\infty}\on{Sch}_{\on{qcqs}}$. Suppose
\begin{itemize}
	\item $\Coh(T)\subset \QCoh(T)$ is closed under truncations,
\end{itemize}
then the Beck-Chevalley natural transformation (for the $\IndCoh_*$-theory)
$$g^{*} \circ f_*\to \varphi_* \circ \phi^{*}$$
is an isomorphism.
\end{lemma}

\begin{warning} If we do not assume $\Coh(T)\subset \QCoh(T)$ to be closed under truncations, the above lemma might be false.
\end{warning}

\proof When $g$ is flat, under the additional assumption (ii), the proof in \cite{raskin2020homological} is correct. The proof can be generalized when $g$ is of bounded Tor dimension.

\qed

\begin{proposition} (C.f. \cite[Proposition 6.7.2]{raskin2020homological}) \label{prop-flat-descent-indcoh-sch}

Let $f:T\to S$ be a faithfully flat map in $^{>-\infty}\on{Sch}_{\on{qcqs}}$, then $*$-pullback induces a fully faithful functor
$$ \IndCoh_*(S) \to \on{Tot}_{\on{semi}} ( \IndCoh_*(T^{\times_S (\bullet+1)}) ).$$
If $\Coh(S)\subset \QCoh(S)$ is closed under truncations, then the above functor is an equivalence.
\end{proposition}

\begin{warning} If we do not assume $\Coh(S)\subset \QCoh(S)$ to be closed under truncations, the above proposition might be false. Under this additional assumption, the proof in \cite{raskin2020homological} is correct.
\end{warning}

\begin{lemma} \label{lem-*-push-right-t-exact-up-to-shift}
For any $f: S\to T$ in $ ^{>-\infty}\on{Sch}_{\on{qcqs}}$, the functor $f_*:\IndCoh_*(S)\to \IndCoh_*(T)$ is right t-exact up to a shift.
\end{lemma}

\proof Follows from Lemma \ref{lemma-f_*-right-t-exact} and Zariski descent.

\qed

\subsection{``Apaisant'' schemes: when ind-coherent sheaves are well-behaved}

Let $S$ be any eventually coconnective qcqs scheme. One can show that $\on{Coh}(S) \subset \on{QCoh}(S)$ is closed under truncations if and only if $S$ is \emph{coherent}, i.e., locally of the form $\on{Spec}(A)$ such that $A$ is (left) \emph{coherent} (see \cite[Definition 7.2.4.16]{HA}). By the previous discussions, we want to restrict our study to locally coherent schemes. However, the notion of locally coherence is not preserved by standard geometric constructions\footnote{For example, $A$ being coherent does not imply $A[t]$ being coherent.}. Hence we want to restrict to a smaller class of schemes, which we call \emph{apaisant schemes}\footnote{The relation between apaisant schemes and coherent schemes is similar to that between finite type schemes and Noetherian schemes.}. 

\begin{definition} $S\in \,^{>-\infty} \on{Sch}_{\on{qcqs}}$ is said to be \emph{apaisant}, or $S\in \on{Sch}_{\on{apai}}$ if it can be written as a (co)filtered limit $S=\lim_i S_i$ under faithfully flat affine structure maps with $S_i\in \,^{>-\infty} \on{Sch}_{\on{aft}}$.
\end{definition}

\begin{lemma} \label{lem-indcoh-on-placid} (C.f. \cite[Lemma 6.36.3]{raskin2020homological}) 
Suppose $S\in \,^{>-\infty} \on{Sch}_{\on{qcqs}}$ can be written as a filtered limit $S=\lim_i S_i$ under flat affine structure maps with $S_i\in \,^{>-\infty} \on{Sch}_{\on{aft}}$. Then $*$-pullback functors induce an equivalence
$$ \colim_{*\on{-pull}} \IndCoh(S_i) \xrightarrow{\simeq} \IndCoh_*(S).$$
In particular, $*$-pushforward functors induce an equivalence
$$ \IndCoh_*(S) \xrightarrow{\simeq} \lim_{*\on{-push}} \IndCoh(S_i) .$$
\end{lemma}

\begin{warning} If we do not have finiteness assumption on $S_i$, the above lemma might be false. Under this additional assumption, the proof in \cite{raskin2020homological} is correct.
\end{warning}

\begin{remark} In the definition of apaisant schemes, we required the connecting maps $S_i \to S_j$ to be surjective. It is possible to develop the theory without it, but this requirement simplifies many discussions (e.g. the statement and proof of Lemma \ref{lemma-pro-afp-map-is-apaisant} below). 
\end{remark}

\begin{remark} \label{rem-alternative-def-for-indcoh*-apai}
Let $\on{Pro}^{\on{aff,ff}}(^{>-\infty} \on{Sch}_{\on{aft}})$ be the full subcategory of $\on{Pro}(^{>-\infty} \on{Sch}_{\on{aft}})$ containing those filtered systems connected by faithfully flat affine maps. Using Noetherian approximation (see \cite[$\S$ 4.4]{lurie2018spectral}), one can show the functor
$$ \on{Pro}^{\on{aff,ff}}(^{>-\infty} \on{Sch}_{\on{aft}}) \to \,^{>-\infty}\on{Sch}_{\on{qcqs}} $$
obtained by right Kan extension along $^{>-\infty} \on{Sch}_{\on{aft}}\to  \,^{>-\infty}\on{Sch}_{\on{qcqs}}$ is fully faithful. By definition, the image of this functor is $\on{Sch}_{\on{apai}}$. Hence by Lemma \ref{lem-indcoh-on-placid}, the functor $\IndCoh_*: \on{Sch}_{\on{apai}} \to \on{DGCat}_{\on{cont}}$ that encodes the $*$-pushforward functors is the right Kan extension of $\IndCoh: \,^{>-\infty} \on{Sch}_{\on{aft}} \to \on{DGCat}_{\on{cont}} $.
\end{remark}

\begin{corollary} \label{cor-apai-coh-truncation}
Suppose $S\in \,^{>-\infty} \on{Sch}_{\on{qcqs}}$ is apaisant, then $\Coh(S)\subset \QCoh(S)$ is closed under truncations.
\end{corollary}

\proof Follows from Lemma \ref{lem-indcoh-on-placid} and the fact that $\Coh(S_i)$ are closed under truncations.

\qed

\begin{remark} \label{rem-t-structure-apai-sch}
Let $S=\lim_i S_i$ be an apaisant presentation. By Lemma \ref{lem-f_*-left-t-exact} and Lemma \ref{lemma-f_*-right-t-exact}, the equivalence $ \IndCoh_*(S) \xrightarrow{\simeq} \lim_{*\on{-push}} \IndCoh(S_i)$ induces equivalences
\begin{eqnarray*}
	 \IndCoh_*(S)^{\ge 0} \xrightarrow{\simeq} \lim_{*\on{-push}} \IndCoh(S_i)^{\ge 0},\;  \IndCoh_*(S)^{\le 0} \xrightarrow{\simeq} \lim_{*\on{-push}} \IndCoh(S_i)^{\le 0}
\end{eqnarray*}
\end{remark}

\begin{proposition} \label{prop-indcoh-zariski-sch}
Let $S_1 \to S$ and $S_2\to S$ be two open embeddings in $\on{Sch}_{\on{apai}}$ such that $S_1\sqcup S_2 \to S$ is surjective. Then $*$-pullbacks induce an equivalence
\begin{equation} \label{eqn-prop-indcoh-zariski-sch-1}
 \IndCoh_*(S) \to \IndCoh_*(S_1) \times_{ \IndCoh_*(S_1\cap S_2) } \IndCoh_*(S_2).
 \end{equation}
\end{proposition}

\proof This functor has a right adjoint defined as follows. Let $u_i: S_i \to S$, $v_{i}:S_1\cap S_2 \to S_i$ and $v: S_1\cap S_2 \to S$ be the schematic open embeddings. For an object 
$$(  \mathcal{F}_1, \mathcal{F}_2, \mathcal{F}_{12}, (v_i)^*(\mathcal{F}_i) \simeq  \mathcal{F}_{12}) \in  \IndCoh_*(S_1) \times_{ \IndCoh_*(S_1\cap S_2) } \IndCoh_*(S_2), $$
consider the object
$$  (u_1)_{*}(\mathcal{F}_1) \times_{  v_{*}(\mathcal{F}_{12})  } (u_2)_{*}(\mathcal{F}_2) \in \IndCoh_*(S) . $$
This defines a functor
$$ \IndCoh_*(S_1) \times_{ \IndCoh_*(S_1\cap S_2) } \IndCoh_*(S_2)  \to \IndCoh_*(S),$$
which is right adjoint to (\ref{eqn-prop-indcoh-zariski-sch-1}). It remains to show these adjoint functors are inverse to each other. We need to check the following:
\begin{itemize}
	\item[(i)] For any $\mathcal{F} \in  \IndCoh_*(S)$, the morphism
	$$ \mathcal{F} \to (u_1)_{*}\circ (u_1)^*(\mathcal{F}) \times_{  v_{*}\circ (v)^*(\mathcal{F})  } (u_2)_{*}\circ (u_2)^*(\mathcal{F}) $$
	is an equivalence.

	\item[(ii)] For any $(  \mathcal{F}_1, \mathcal{F}_2, \mathcal{F}_{12}, (v_i)^*(\mathcal{F}_i) \simeq  \mathcal{F}_{12})$ as before, the morphism
	$$ (u_i)^{*}[(u_1)_{*}(\mathcal{F}_1) \times_{  v_{*}(\mathcal{F}_{12})  } (u_2)_{*}(\mathcal{F}_2) ] \to \mathcal{F}_i $$
	is an equivalence.
\end{itemize}
(ii) follows from the base-change isomorphisms. To prove (i), note that both sides commute with taking filtered colimits for $\mathcal{F}$. Hence we can assume $\mathcal{F}\in \on{Coh}(S)$. This allows us to replace $\IndCoh_*(S)$ by $\QCoh^*(S)$, where the claim is well-known. 

\qed

Corollary \ref{cor-apai-coh-truncation} implies $\IndCoh_*$-theory on apaisant schemes is well-behaved. For example, we always have base-change isomorphisms. Note that the full subcategory $\on{Sch}_{apsi} \subset \on{Sch}$ is not stable under fiber products. However, we have the following results.

\begin{lemma} \label{lemma-afp-over-apaisant-is-apaisant}
Let $f:S\to T$ be an afp (i.e., almost of finite presentation) map in $^{>-\infty}\on{Sch}_{\on{qcqs}}$. If $T$ is apaisant, then $S$ is apaisant.
\end{lemma}

\proof Let $T= \lim_{i\in I^\op} T_i$ be an apaisant presentation, where $I$ is filtered. By assumption, there exists an integer $m$ such that $S$, $T$ and $T_i$ are contained in $^{\ge -m}\on{Sch}_{\on{qcqs}}$. Note that $f:S\to T$ is finitely $m$-presented (see \cite[Definition 4.2.3.1]{lurie2018spectral}). Hence by Noetherian approximation (see \cite[Proposition 4.4.4.1]{lurie2018spectral}), there exist an index $j$ and a finitely $m$-presented map $S_j\to T_j$ such that $S \simeq \tau^{>-m}(T \times_{T_j} S_j)$. Since $T\to T_j$ is flat, the truncation $\tau^{>-m}$ is superfluous, i.e., $S\simeq T \times_{T_j} S_j$. Now it is easy to see $S \simeq \lim_{ i \in (I_{j/ })^\op} (T_i \times_{T_j} S_j) $ is an apaisant presentation for $S$.

\qed

\begin{definition} A map $f:S\to T$ in $\on{PreStk}$ is called \emph{relatively apaisant} if it is qcqs, schematic, of Tor dimension $\le n$ for some integer $n$, and for any apaisant scheme $T'$ over $T$, the fiber product $S\times_T T'$ is apaisant.
\end{definition}

\begin{lemma} \label{lemma-afp-map-is-apaisant}
Suppose $f:S\to T$ in $\on{PreStk}$ is qcqs, schematic, afp and of Tor dimensinon $\le n$ for some integer $n$. Then $f$ is relatively apaisant.
\end{lemma}

\proof Follows from Lemma \ref{lemma-afp-over-apaisant-is-apaisant}.  

\qed

\begin{lemma} \label{lemma-pro-afp-map-is-apaisant}
Let $f:S \to T$ be a map in $\on{PreStk}$. Suppose $S$ can be written as a limit $S = \lim_{n\in \mathbb{N}} S_n$ under faithfully flat affine afp structure maps with $S_i\to T$ being relatively apaisant. Then $f$ is relatively apaisant.
\end{lemma}

\proof We only need to show $S$ is an apaisant scheme when each $S_i$ is assumed to be an apaisant sheme. Let $S_0 \simeq \lim_{\alpha\in J^\op} S_0^\alpha$ be an apaisant presentation. As in the proof of Lemma \ref{lemma-afp-over-apaisant-is-apaisant}, using induction, there exists a chain $\alpha_1 \to \alpha_2 \to \cdots$ in $J$ and apaisant presentations $S_n \simeq \lim_{\alpha\in (J_{\alpha_n/  })^\op} S_n^\alpha $ such that $S_{n+1}^\alpha \simeq S_n^\alpha \times_{  S_n^{\alpha_{n+1}} } S_{n+1}^{\alpha_{n+1}}$.  Now it is easy to see $S \simeq \lim_{ n\in \mathbb{N},  \alpha \in (J_{\alpha_n/  })^\op  } S_n^\alpha $ is an apaisant presentation for $S$.

\qed

\begin{example} \label{exam-classifying-stack-apaisant}
Let $H=\lim_{n} H_n$ be an apaisant presentation for an apaisant affine flat group scheme over a base $B\in \,^{>-\infty}\on{Sch}_{\on{aft}}$. Then $B \to \mathbb{B} H$ is relatively apaisant. Indeed, for any apaisant scheme $S$ over $\mathbb{B} H$, we only need to show the corresponding $H$-torsor $P_H$ over $S$ is an apaisant scheme. Let $P_{H_n}$ be the induced $H_n$-torsor, which is apaisant by Lemma \ref{lemma-afp-over-apaisant-is-apaisant}. Therefore $P_H \simeq \lim_n P_{H_n}$ is apaisant by Lemma \ref{lemma-pro-afp-map-is-apaisant}.
\end{example}

\begin{construction} \label{const-!-pull-afp-apsi-sch}
Let $f: S\to T$ be an afp map in $\on{Sch}_{\on{apai}}$. Choose an apaisant presentation $T= \lim_{i\in I^\op} T_i$. As in the proof of Lemma \ref{lemma-afp-over-apaisant-is-apaisant}, there exists an index $j$ and $S_j\in (^{>-\infty} \on{Sch}_{\on{aft}})_{/T_j}$ such that $S\simeq T\times_{T_j} S_j$, and $S \simeq \lim_{ i \in (I_{j/ })^\op} (T_i \times_{T_j} S_j) $ is an apaisant presentation for $S$. By Remark \ref{rem-alternative-def-for-indcoh*-apai}, we have the following functor
$$ \IndCoh_*(T) \simeq \lim_{i \in (I_{j/ })^\op} \IndCoh(T_j) \xrightarrow{!\on{-pull}}  \lim_{i \in(I_{j/ })^\op} \IndCoh(S_j) \simeq \IndCoh_*(S),$$
where the middle functor is defined using the base-change isomorphisms between $!$-pullback and $*$-pushforward functors (for the $\IndCoh$-theory on $^{-\infty}\on{Sch}_{\on{ft}}$ ). Alternatively, it can be defined as
$$ \IndCoh_*(T) \simeq \colim_{i \in I_{j/ }} \IndCoh(T_j) \xrightarrow{!\on{-pull}}  \colim_{i \in I_{j/ }} \IndCoh(S_j) \simeq \IndCoh_*(S),$$
where the middle functor is defined using the commutativity between $!$-pullback and $*$-pullback functors (see \cite[Proposition 7.1.6]{gaitsgory2013ind}). Note that the obatined functor 
$$f^!: \IndCoh_*(T)\to  \IndCoh_*(S)$$
does not depend on the choice of the apaisant presentation.
\end{construction}

\begin{construction} \label{constr-base-change-!*-apai-sch}
Let 
$$
\xymatrix{
	S' \ar[r]^-{\varphi} \ar[d]^-\phi & T'\ar[d]^-g \\
	S \ar[r]^-f & T
}
$$
be a Cartesian diagram contained in $\on{Sch}_{\on{apai}}$ such that $f$ is afp (hence so is $\varphi$). Then we have a canonical isomorphism (for $\IndCoh_*$-theory):
\begin{eqnarray*}
\phi_* \circ \varphi^! \simeq f^! \circ g_*
\end{eqnarray*}
induced via right Kan extensions from the similar isomorphisms for $\IndCoh$-theory on aft schemes.

If $g$ is further assumed to be of bounded Tor dimension (hence so is $\phi$), then the Beck-Chevalley natural transformation 
\begin{eqnarray*}
\phi^{*} \circ f^! \to \varphi^! \circ g^{*}
\end{eqnarray*}
induced by the above isomorphism is also an isomorphism.

\end{construction}

\begin{construction}
Let $f:S\to T$ be an afp open immersion in $\on{Sch}_{\on{apai}}$, by \cite[Corollary 4.6.2.2]{lurie2018spectral}, the map $S_j\to T_j$ in Contruction \ref{const-!-pull-afp-apsi-sch} can be chosen to be an open immersion. It follows that there is a canonical isomorphism $f^!\simeq f^*$.
\end{construction}

\begin{construction}
Let $f:S\to T$ be an afp proper map in $\on{Sch}_{\on{apai}}$. Then the map $S_j\to T_j$ in Contruction \ref{const-!-pull-afp-apsi-sch} is also proper (because $T\to T_j$ is \emph{faithfully} flat). It follows that there is a canonical adjoint pair $(f_*f^!)$.
\end{construction}

\begin{lemma} \label{lem-left-t-exact-!-pull-apai-sch}
Let $f: S\to T$ be an afp map in $\on{Sch}_{\on{apai}}$. The functor $f^!$ is left t-exact up to a shift.
\end{lemma}

\proof Follows from Remark \ref{rem-t-structure-apai-sch} and \cite[Lemma 7.1.7]{gaitsgory2013ind}.

\qed

We also need the following sharper result:

\begin{lemma} \label{lem-left-t-exact-!-pull-apai-sch-sharper}
Let $f: S\to T$ be an afp map in $\on{Sch}_{\on{apai}}$. Then there exists a natural number $n$ such that for any Cartesian square
$$
\xymatrix{
	S' \ar[r]^-{\varphi} \ar[d]^-\phi & T'\ar[d]^-g \\
	S \ar[r]^-f & T
}
$$
contained in $\on{Sch}_{\on{apai}}$, the functor $[-n]\circ \varphi^!$ is left t-exact.
\end{lemma}

\proof We can replace $S\to T$ by $S_j\to T_j$ in Contruction \ref{const-!-pull-afp-apsi-sch}. Hence we can assume $S,T\in \,^{>-\infty}\on{Sch}_{\on{aft}}$. Similar, using Remark \ref{rem-t-structure-apai-sch}, we can also assume $S', T' \in \,^{>-\infty}\on{Sch}_{\on{aft}}$. As in the proof of \cite[Lemma 7.1.7]{gaitsgory2013ind}\footnote{Warning: the derived Nagata compactification, which this proof depends on, currently has a gap.}, we can reduce to the case when $f$ is proper. Hence we only need to show $[n]\circ \varphi_*$ is right t-exact, where $n$ is an integer that depends only on $f$. 

The remaining argument is standard . If $f$ is affine, we can take $n:=0$. 

If $f$ is separated, we can find finite many open subschemes $\{S_i\}_{i\in I}$ of $S$ that covers $S$ such that each composition $S_i \to S\to T$ is affine. Note that for any $J\subset I$, the intersection $\cap_{j\in J} S_j$ is affine over $T$ (because $f$ is separated). Then we can take $n:=|I|$. 

For general map $f$, note that $f$ is quasi-separated. We can find finite many open subschemes $\{S_i\}_{i\in I}$ of $S$ that covers $S$ such that each composition $S_i \to S\to T$ is separated. Note that for any $J\subset I$, the intersection $\cap_{j\in J} S_j$ is separated over $T$ (because $f$ is quasi-separated). Let $n_J$ be the natural number constructed in the last paragraph for the map $\cap_{j\in J} S_j \to T$. Then we can take $n:= |I|+ \on{max}_J\{ n_J \}$.

\qed

\subsection{Ind-coherent sheaves on indschemes}

In this subsection, we study the $\IndCoh_*$-theory on (apaisant) indschemes. However, we postpone the study of $!$-pullback functors to the next subsection.

\begin{definition} A \alert{convergent} prestack $S$ is an \alert{apaisant indscheme}, or $S\in \on{IndSch}_{\on{apai}}$, if it can be written as a filtered colimit $S \simeq \colim_i S_i$ such that
\begin{itemize}
	\item $S_i \in \on{Sch}_{\on{apai}}$,
	\item $S_i \to S_j$ are afp closed embeddings.
\end{itemize}
\end{definition}

\begin{remark} Any apaisant indscheme is \emph{reasonable} in the sense of \cite[$\S 6.8$]{raskin2020homological}.

\end{remark}

\begin{construction} \label{constr-indcoh-indsch}
We define
$$\IndCoh_*: \on{IndSch}_{\on{apai}} \to \on{DGCat}_{\on{cont}}$$
by left Kan extension along $\on{Sch}_{\on{apai}} \subset \on{IndSch}_{\on{apai}}$. Explicitly, we have
$$ \IndCoh_*(S) \simeq \colim_{*\on{-push}} \IndCoh_*(S_i) \simeq \lim_{!\on{-pull}} \IndCoh_*(S_i)$$
\end{construction}

\begin{remark} Note that the $*$-pushforward functor $\IndCoh_*(S_i) \to \IndCoh_*(S_j)$ sends $\Coh(S_i)$ into $\Coh(S_j)$. In particular, $\IndCoh_*(S)$ is compactly generated by 
$$ \Coh(S):= \colim_{*\on{-push}} \Coh(S_i),$$
where the colimit is taken in $\on{Cat}$.
\end{remark}

\begin{remark} If $S$ is locally almost finite type, i.e., $S\in \on{IndSch}_{laft}$, the obtained $\IndCoh_*(S)$ coincides with $\IndCoh(S)$ defined in \cite{gaitsgory2013ind}. Hence we write $\IndCoh(S)$ for $\IndCoh_*(S)$ in this case.
\end{remark}

\begin{construction} \label{const-t-structure-indcoh-indsch}
By Lemma \ref{lemma-f_*-right-t-exact}, for $S \simeq \colim_i S_i \in \on{IndSch}_{\on{apai}}$, there exists a t-structure on $\IndCoh_*(S)$ such that 
$$\IndCoh_*(S)^{\ge 0} \simeq \lim_{!\on{-pull}} \IndCoh_*(S_i)^{\ge 0}.$$
This t-structure does not depend on the choice of the presentation $S \simeq \colim_i S_i$.

It is not hard to show the above t-structure is compactly generated and 
$$\IndCoh_*(S)^{\le 0} \simeq \colim_{*\on{-push}} \IndCoh_*(S_i)^{\le 0},$$
where the colimit is taken in $\on{Pr}^L$.
\end{construction}

\begin{corollary} \label{cor-coh-closed-truncation-apai-indsch}
For $S\in \on{IndSch}_{\on{apai}}$, $\on{Coh}(S) \subset \IndCoh_*(S)$ is closed under truncations.
\end{corollary}

\proof Follows from Corollary \ref{cor-apai-coh-truncation}.

\qed

\begin{lemma} (C.f. \cite[Lemma 6.11.2]{raskin2020homological}) 
\label{lem-desc-coh-indsch}
For $S\in \on{IndSch}_{\on{apai}}$, $\mCF\in \IndCoh_*(S)$ is coherent if and only if $\mCF$ is bounded and almost perfect, i.e., $\mCF\in \IndCoh_*(S)^+$ and for any integer $N$ with $\mCF\in \IndCoh_*(S)^{\ge -N}$, $\mCF$ is compact in $\IndCoh_*(S)^{\ge -N}$.
\end{lemma}

\begin{lemma} \label{lem-*-push-left-t-exact-indsch}
For $f:S\to T$ in $\on{IndSch}_{\on{apai}}$, the functor $f_*: \IndCoh_*(S) \to \IndCoh_*(T)$ is left t-exact. If $f$ is a closed immersion\footnote{This means for some apaisant presentations $S = \colim_{i\in \mathcal{I}} S_i$ and $T=\colim_{j\in \mathcal{J}} T_j$, and any index $i \in \mathcal{I}$, there exists $j\in \mathcal{J}$ such that $S_i \to S\to T$ factors through a closed immersion $S_i \to T_j$.}, then $f_*$ is t-exact.
\end{lemma}

\proof If $f$ is a closed immersion, $f_*$ is right t-exact by definition. Hence it remains to prove the first claim.

Let $T=\colim T_i$ be an apaisant presentation. We first prove in the case $S=T_j$ for some index $j$. We only need to show 
$$ \IndCoh_*(T_j)\xrightarrow{*\on{-push}}  \IndCoh_*(T) \xrightarrow{!\on{-pull}}  \IndCoh_*(T_i)$$
is left t-exact for any $i$. Since the $*$-pushforward and $!$-pullback functors along any $T_{i_1} \to T_{i_2}$ is left t-exact, we only need to prove in the case $i=j$. Then the above composition is the filtered colimit of the following functors
$$ \IndCoh_*(T_j) \xrightarrow{*\on{-push}}  \IndCoh_*(T_k)\xrightarrow{!\on{-pull}} \IndCoh_*(T_j),$$
where the colimit is taken for all arrows $j\to k$. Then we are done because the t-structure on $\IndCoh_*(T_j)$ is compatible with filtered colimits.

We then prove in the case when $S\in \on{Sch}_{\on{apai}}$. Note that $f:S\to T$ factors through some $T_i$. This reduces the problem to the previous case.

Finally, for general $S\in \on{IndSch}_{\on{apai}}$, let $S= \colim S_\alpha$ be an apaisant presentation. Note that the identity functor on $\IndCoh_*(S)$ is the filtered colimit of the following left t-exact functors
$$ \IndCoh_*(S) \xrightarrow{!\on{-pull}}  \IndCoh_*(S_\alpha) \xrightarrow{*\on{-push}}  \IndCoh_*(S).$$
Hence the previous case (applied to each $S_\alpha \to T$) implies $f_*$ is left t-exact.

\qed

\begin{lemma} (C.f. \cite[Lemma 6.16.1, Corollary 6.16.2]{raskin2020homological})
\label{lem-*-pull-indsch}
Suppose $f: S\to T$ in $ \on{IndSch}_{\on{apai}} $ is relatively apaisant. Then:
\begin{itemize}
	\item $f_*:\IndCoh_*(S) \to \IndCoh_*(T)$ has a right t-exact left adjoint $f^{*}$. 

	\item If $f$ is of Tor dimension $\le n$, then $[-n]\circ f^{*}$ is left t-exact.

	\item For a Cartesian diagram
$$
\xymatrix{
	S' \ar[r]^-{\varphi} \ar[d]^-\phi & T'\ar[d]^-g \\
	S \ar[r]^-f & T
}
$$
contained in $\on{IndSch}_{\on{apai}}$, if $f$ is relatively apaisant (hence so is $\varphi$), then the Beck-Chevalley natural transformation (for $\IndCoh_*$-theory)
$$ f^{*}\circ g_* \to \phi_* \circ \varphi^{*} $$
is an isomorphism.
\end{itemize}

\end{lemma}

\proof When $f$ is flat, the proof in \cite{raskin2020homological} is correct when restricted to apaisant indschemes. The proof can be generalized to maps of bounded Tor dimensions.

\qed

\begin{lemma} \label{lem-indcoh-on-apai-indsch} Suppose $S\in \on{Sch}_{qcqs}$ can be written as a filtered limit $S=\lim_{i\in I^\op} S_i$ under faithful flat affine structure maps with $S_i$ being finite type. Then $S$ is an apaisant \emph{indscheme} and $*$-pullback functors induce an equivalence
$$ \colim_{*\on{-pull}} \IndCoh(S_i) \xrightarrow{\simeq} \IndCoh_*(S).$$
In particular, $*$-pushforward functors induce an equivalence
$$ \IndCoh_*(S) \xrightarrow{\simeq} \lim_{*\on{-push}} \IndCoh(S_i) .$$
\end{lemma}

\proof Note that $\tau^{\ge -n} (S) \simeq \lim_{i\in I^\op} \tau^{\ge -n} (S_i)$. Therefore $\tau^{\ge -n} (S) $ is an apaisant scheme. Moreover, for $m\ge n$, the map $\tau^{\ge -n} (S) \to \tau^{\ge -m} (S)$ is an afp closed embedding because it is the base-change of a closed embedding $\tau^{\ge -n}(S_i)\to  \tau^{\ge -m}(S_i )$ contained in $^{>-\infty}\on{Sch}_{ft}$. Hence $S \simeq \colim_n \tau^{\ge -n} (S)$ is an apaisant indscheme. Now we have
\begin{eqnarray*}
 \colim_{i,*\on{-pull}} \IndCoh(S_i) &\simeq& \colim_{i,*\on{-pull}} \colim_{n,*\on{-push}}  \IndCoh(\tau^{\ge -n} (S_i)) \\
 &\simeq&  \colim_{n,*\on{-push}}  \colim_{i,*\on{-pull}} \IndCoh(\tau^{\ge -n} (S_i))  \\
 &\simeq& \colim_{n,*\on{-push}}  \IndCoh_*(\tau^{\ge -n} (S)) \\
 &\simeq& \IndCoh_*(S)
\end{eqnarray*}
where the first equivalence is due to \cite[Proposition 4.3.4]{gaitsgory2013ind}, the second is due to the base-change isomorphisms, the third is by Lemma \ref{lem-indcoh-on-placid}, and the last is by definition.

\qed

\subsection{Ind-coherent sheaves on indschemes: \texorpdfstring{$!$}{!}-pullbacks}

In this subsection, we study the $!$-pullback functors for $\IndCoh_*$-thoery on apaisant indschemes.

\begin{definition}  \label{def-ind-afp-map-indsch}
A map $f:S\to T$ in $\on{IndSch}_{\on{apai}}$ is called \emph{ind-afp} if for some apaisant presentations $S = \colim_{i\in \mathcal{I}} S_i$ and $T=\colim_{j\in \mathcal{J}} T_j$, and any index $i \in \mathcal{I}$, there exists $j\in \mathcal{J}$ such that $S_i \to S\to T$ factors through an afp map $S_i \to T_j$.

We say $f:S\to T$ is ind-afp and \emph{ind-proper} (resp. an \emph{ind-afp closed embedding}) if the above maps $S_i \to T_j$ are afp and proper (resp. afp closed embeddings).
\end{definition}

\begin{remark} By \cite[Proposition 6.12.1]{raskin2020homological} (whose proof is correct for apaisant indschemes), we see that in the above definition, we can replace ``some apaisant presentations'' by ``any apaisant presentations''.
\end{remark}

\begin{remark} \label{rem-ind-of-indsch-is-indsch}
A filtered colimit of apaisant indschemes connected by ind-afp closed embeddings is an apaisant indscheme.
\end{remark}

\begin{construction} Let $(\on{Sch}_{\on{apai}})_{\on{afp}}$ (resp. $(\on{IndSch}_{\on{apai}})_{\on{ind-afp}}$) be the 1-full subcategory of $\on{Sch}_{\on{apai}}$ (resp. $\on{IndSch}_{\on{apai}}$) containing afp (resp. ind-afp) morphisms. Let 
$$ \IndCoh_*: (\on{Sch}_{\on{apai}})_{\on{afp}}^\op \to \on{DGCat}_{\on{cont}}$$
be the functor that encodes $!$-pullback functors for the $\IndCoh_*$-theory on apaisant schemes. By right Kan extension, we obtain a functor
$$ \IndCoh_*: (\on{IndSch}_{\on{apai}})_{\on{ind-afp}}^\op \to \on{DGCat}_{\on{cont}}$$
which encodes $!$-pullback functors for the $\IndCoh_*$-theory on apaisant indschemes. Note that for any apaisant indscheme $S$, the obtained category $\IndCoh_*(S)$ is equivalent to that in Construction \ref{constr-indcoh-indsch}.

For any ind-proper map $f:S\to T$ in $\on{IndSch}_{\on{apai}}$, we have a canonical adjoint pair $(f_*,f^!)$.

For any schematic open embedding map $f:S\to T$ in $\on{IndSch}_{\on{apai}}$, we have a canonical adjoint pair $(f^!,f_*)$.
\end{construction}

We are going to discuss the base-change isomorphisms between $!$-pullback and $*$-pushforward (resp. $*$-pullback) functors. As a preparation, we have:

\begin{lemma} \label{lem-fiber-product-apai-indsch-is-apai-if-afp}
Let $f:S\to T$ be an ind-afp map in $\on{IndSch}_{\on{apai}}$ and $g:T'\to T$ be any map in $\on{IndSch}_{\on{apai}}$. Then the fiber product $S\times_T T'$ (in $\on{PreStk}_{\on{conv}}$) is an apaisant indscheme, and the map $S'\times_T T' \to T'$ is ind-afp.

If $f$ is further assumed to be ind-proper (resp. an ind-afp closed embedding), then so is $S'\times_T T' \to T'$.
\end{lemma}

\proof We first assume $S, T, T'$ are apaisant schemes. Let $S_j\to T_j$ be as in the proof of Lemma \ref{lemma-afp-over-apaisant-is-apaisant}. Let $T' \simeq \lim_{k\in K^\op} T'_k$ be an apaisant presentation of $T'$ defined over $T_j$. Then $S\times_T T' \simeq \lim_{k\in K^\op}  S_j\times_{T_j} T'_k$ is an apaisant indscheme by Lemma \ref{lem-indcoh-on-apai-indsch}. Moreover, each map $\tau^{\ge -n} (S\times_T T')\to T'$ is afp (resp. afp and proper, an afp closed embedding) because it is the base-change of $\tau^{\ge -n} ( S_j\times_{T_j} T'_k)\to T'_k$.

Then we only assume $S$ and $T'$ are apaisant schemes. Let $T \simeq \colim_\alpha T_\alpha$ be as in the definition of apaisant indschemes. We can assume this filtered diagram has an initial object $T_{\alpha_0}$ and the maps $S\to T$ and $S\to T'$ factor through $T_{\alpha_0}$. By the previous case and Remark \ref{rem-ind-of-indsch-is-indsch}, we only need to show $S\times_{T_\alpha} T' \to S\times_{T_\beta} T'$ is an afp closed embedding. However, this follows from the fact that it is the base-change of $T_\alpha \to T_\alpha\times_{T_\beta} T_\alpha$, which is an afp closed embedding because the maps from both sides to $T_\beta$ are afp.

Finally, the general case follows from Remark \ref{rem-ind-of-indsch-is-indsch}.

\qed

\begin{lemma} (C.f. \cite[Lemma 6.17.1, Lemma 6.17.2]{raskin2020homological})
\label{lem-base-change-indsch}
Let
$$
\xymatrix{
	S' \ar[r]^-{\varphi} \ar[d]^-\phi & T'\ar[d]^-g \\
	S \ar[r]^-f & T
}
$$
be a Cartesian diagram in $\on{IndSch}_{\on{apai}}$ such that $f$ is ind-afp and ind-proper (hence so is $\varphi$). Then the Beck-Chevalley natural transformation
\[
\phi_* \circ \varphi^!  \to f^! \circ g_*.
\]
is an isomorphism. If moreover $g$ is relatively apaisant (hence so is $\phi$), then the Beck-Chevalley natural transformation
\[
\phi^{*} \circ f^! \to \varphi^! \circ g^{*}
\]
obtained from the above isomorphism is also an isomorphism.
\end{lemma}

\proof We first assume $S, T, T'$ are apaisant schemes. Using Lemma \ref{lem-indcoh-on-apai-indsch}, as in the proof of Lemma \ref{lem-fiber-product-apai-indsch-is-apai-if-afp}, the desired isomorphisms can be induced from the base-change isomorphisms for $\IndCoh$-theory on $\on{Sch}_{ft}$.

Then we only assume $S$ and $T'$ are apaisant schemes. Let $T \simeq \colim_\alpha T_\alpha$ be as in the definition of apaisant indschemes. We can assume this filtered diagram has an initial object $T_{\alpha_0}$ and the maps $S\to T$ and $S\to T'$ factor through $T_{\alpha_0}$. Consider the Cartesian diagram
$$
\xymatrix{
	S'_\alpha \ar[r]^-{\varphi_\alpha} \ar[d]^-{\phi_\alpha} & T'\ar[d]^-{g_\alpha} \\
	S \ar[r]^-{f_\alpha} & T_\alpha
}
$$
and the ind-afp closed embeddings $\iota_\alpha: T_\alpha\to T$ and $\iota_\alpha': S_\alpha'\to S$. We have
\begin{equation} \label{eqn-proof-lem-base-change-indsch-1}
 \phi_* \circ \varphi^! \simeq \colim_\alpha \phi_* \circ (\iota_\alpha')_*\circ (\iota_\alpha')^! \circ \varphi^! \simeq  \colim_\alpha\phi_{\alpha,*} \circ \varphi_\alpha^! \simeq \colim_\alpha f_\alpha^! \circ g_{\alpha, *} \simeq f^!\circ g_*.
 \end{equation}

In the general case, let $T' \simeq \colim T'_\beta$ and $S \simeq \colim S_\gamma$ be as in the definition of apaisant indschemes. Then the base-change isomorphisms between $!$-pullback and $*$-pushforward functors for the Cartesian squares
$$
\xymatrix{
	S_\gamma\times_T T'_\beta  \ar[r] \ar[d] & T'_\beta\ar[d]\\
	S_\gamma \ar[r] & T
}
$$
imply the desired isomorphism $\phi_* \circ \varphi^!  \to f^! \circ g_*$.

To prove $\phi^{*} \circ f^! \simeq \varphi^! \circ g^{*}$ when $g$ is assumed to be relatively apaisant, we can first reduce to the case when $S$ is an apaisant scheme, then reduce to the case when $T$ is an apaisant scheme (hence so is $T'$), and finally using Noetherian approximation to reduce to the finite type case.

\qed

\begin{proposition} \label{prop-indcoh-zariski-indsch}
Let $S_1 \to S$ and $S_2\to S$ be two schematic open embeddings in $\on{IndSch}_{\on{apai}}$ such that $S_1\sqcup S_2 \to S$ is surjective. Then $*$-pullbacks induce an equivalence
$$
 \IndCoh_*(S) \to \IndCoh_*(S_1) \times_{ \IndCoh_*(S_1\cap S_2) } \IndCoh_*(S_2).
 $$
\end{proposition}

\proof Follows from Proposition \ref{prop-indcoh-zariski-sch} and Lemma \ref{lem-base-change-indsch}.

\qed

\begin{proposition}  \label{theorem-flat-descent-indcoh-indsch}
Let $f:T\to S$ be a faithfully flat and relatively apaisant map in $\on{IndSch}_{\on{apai}}$, then $*$-pullback induces an equivalence
$$ \IndCoh_*(S) \to \on{Tot}_{\on{semi}} ( \IndCoh_*(T^{\times_S (\bullet+1)}) ).$$
\end{proposition}

\proof Follows from Proposition \ref{prop-flat-descent-indcoh-sch} and Lemma \ref{lem-base-change-indsch}.

\qed

\begin{lemma} (C.f. \cite[Lemma 6.18.1]{raskin2020homological}) \label{lem-pullback-being-coh-implies-coh}
Suppose $f: S\to T$ in $ \on{IndSch}_{\on{apai}} $ is faithfully flat and relatively apaisant, then $\mCF\in \IndCoh_*(T)$ lies in $\Coh(T)$ if and only if $f^{*}(\mCF)$ lies in $\Coh(S) \subset \IndCoh_*(S)$.
\end{lemma}

\proof The proof in \cite{raskin2020homological} is correct when restricted to apaisant indschemes.

\qed

\begin{remark} \label{rem-why-not-all-base-change}
It is very possible that for any ind-afp map $f$, there are isomorphisms $\phi_* \circ \varphi^!  \simeq f^! \circ g_*$ and $\phi^{*} \circ f^! \simeq \varphi^! \circ g^{*}$ and higher compatibilities for them. In fact, the proof of Lemma \ref{lem-base-change-indsch} can be generalized to all ind-afp map $f$ as long as the isomorphisms $\phi_{\alpha,*} \circ \varphi_\alpha^!\simeq f_\alpha^! \circ g_{\alpha, *} $ in (\ref{eqn-proof-lem-base-change-indsch-1}) can be shown to be compatible with the colimits, which we do \emph{not} know how to prove homotopy-coherently. However, for the purpose of this paper, we only need the following special case.
\end{remark}

\begin{definition} A morphism $f:S \to T$ in $\on{IndSch}_{\on{apai}}$ is \emph{uniformly afp}, or \emph{uafp}, if there exists a presentation $T=\colim_i T_i$ of $T$ as an apaisant indscheme such that each $T_i \times_T S$ is an apaisant \emph{scheme} and each map $T_i\times_T S \to T_i$ is afp.
\end{definition}

\begin{remark} \label{rem-uafp-vs-sch&afp}
By Lemma \ref{lemma-afp-over-apaisant-is-apaisant}, $f$ is uafp iff
\begin{itemize}
	\item $f:S \to T$ is a schematic and afp morphism in $\on{IndSch}_{\on{afp}}$;
	\item there exists a presentation $T=\colim_i T_i$ of $T$ as an apaisant indscheme such that each $T_i \times_T S$ is eventually coconnective.
\end{itemize}
\end{remark}

\begin{remark} If $f$ is uafp and ind-proper, then $f$ is proper.
\end{remark}

\begin{remark} \label{rem-uafp-base-change-indsch}
Let $f:S \to T$ be a uafp morphism in $\on{IndSch}_{\on{afp}}$, and $g:T'\to T$ be a relatively apaisant morphism. Then $S\times_T T'$ is an apaisant indscheme and $S\times_T T' \to T'$ is uafp.
\end{remark}

\begin{warning} \label{warn-uafp-comp}
The composition of two uafp morphisms might \emph{not} be uafp.
\end{warning}

\begin{construction} 
\label{constr-base-change-!*-apai-indsch}
Let
$$
\xymatrix{
	S' \ar[r]^-{\varphi} \ar[d]^-\phi & T'\ar[d]^-g \\
	S \ar[r]^-f & T
}
$$
be a Cartesian diagram in $\on{IndSch}_{\on{apai}}$ such that both $f$ and $\varphi$ are uafp. Then we have a canonical isomorphism (for $\IndCoh_*$-theory):
\begin{eqnarray*}
\phi_* \circ \varphi^! \simeq f^! \circ g_*
\end{eqnarray*}
induced via left Kan extensions from the similar isomorphisms for apaisant schemes. In Remark \ref{rem-higher-comp-!-pull-indsch} we will explain how to construct the higher compatibilities for these base-change isomoprhisms.

It follows from construction that:
\begin{itemize}
	\item if $f$ is ind-proper (hence so is $\varphi$), then the above isomorphism can be identified with the Beck-Chevalley isomorphism in Lemma \ref{lem-base-change-indsch}

	\item if $g$ is relatively apaisant (hence so is $\phi$), then the Beck-Chevalley natural transformation 
	\[
	\phi^{*} \circ f^! \to \varphi^! \circ g^{*}
	\]
	obtained from the above isomorphism is also an isomorphism.
\end{itemize}
\end{construction}

\begin{lemma} \label{lem-left-t-exact-!-pull-apai-indsch-sharper}
Let $f: S\to T$ be an afp map in $\on{Sch}_{\on{apai}}$. Then there exists a natural number $n$ such that for any Cartesian square
$$
\xymatrix{
	S' \ar[r]^-{\varphi} \ar[d]^-\phi & T'\ar[d]^-g \\
	S \ar[r]^-f & T
}
$$
contained in $\on{IndSch}_{\on{apai}}$ with $\varphi$ being uafp, the functor $[-n]\circ \varphi^!$ is left t-exact.
\end{lemma}

\proof Follows from Lemma \ref{lem-left-t-exact-!-pull-apai-sch-sharper}.

\qed

\begin{construction} 
Now we sketch how to construct the higher compatibilities for the base-change isomorphisms in Construction \ref{constr-base-change-!*-apai-indsch}. Due to difficulties similar to those in Remark \ref{rem-why-not-all-base-change}, we do \emph{not} know how to construct a functor
$$ \on{Corr} ( \on{IndSch}_{\on{apai}} )_{\on{all,ind-afp}} \to \on{DGCat}_{\on{cont}} $$
that encodes the $!$-pullback, $*$-pushforward and the base-change isomorphisms. Instead, we are going to construct a similar functor
\begin{equation} \label{eqn-indcoh-out-ind-pro}
\mathcal{I}nd\mathcal{C}oh^*: \on{Corr} ( \on{Ind}( \on{Pro} (  \on{Sch}_{\on{aft}} ) ) )_{\on{all},\on{fafp}} \to \on{DGCat}_{\on{cont}} 
 \end{equation}
where
\begin{itemize}
	\item a morphism in $\on{Pro} (  \on{Sch}_{\on{aft}} )$ is called \emph{formally afp}, or \emph{fafp}, if it is a base-change\footnote{$ \on{Sch}_{\on{aft}}$ has finite limits, hence so is $\on{Pro} (  \on{Sch}_{\on{aft}} )$ and the functor $\on{Sch}_{\on{aft}}\to \on{Pro} (  \on{Sch}_{\on{aft}} )$ commutes with finite limits.} of a map contained in the full subcategory $\on{Sch}_{\on{aft}}\subset \on{Pro} (  \on{Sch}_{\on{aft}} )$;

	\item a morphism $\mathcal{S}\to \mathcal{T}$ in $\on{Ind}(\on{Pro} (  \on{Sch}_{\on{aft}} ))$ is called \emph{fafp} if for any $\mathcal{T}' \to \mathcal{T}$ with $\mathcal{T}'$ contained in $\on{Pro} (  \on{Sch}_{\on{aft}} )$, the base-change morphism\footnote{$ \on{Pro} (  \on{Sch}_{\on{aft}} )$ has finite limits, hence so is $\on{Ind} (  \on{Pro} (  \on{Sch}_{\on{aft}} ))$ and the functor $\on{Pro} (  \on{Sch}_{\on{aft}} )\to \on{Ind} ( \on{Pro} (  \on{Sch}_{\on{aft}} ) )$ commutes with finite limits.  } $\mathcal{S}\times_{\mathcal{T}} \mathcal{T}' \to \mathcal{T}'$ is a fafp morphism contained in $\on{Pro} (  \on{Sch}_{\on{aft}} )$.
\end{itemize}
Moreover, the functor (\ref{eqn-indcoh-out-ind-pro}) will have a canonical \emph{right-lax} symmetric monoidal structure, where we equip $\on{Corr} ( \on{Ind}( \on{Pro} (  \on{Sch}_{\on{aft}} ) ) )_{\on{all},\on{fafp}}$ with the symmetric monoidal structure induced by the Cartesian products on $ \on{Ind}( \on{Pro} (  \on{Sch}_{\on{aft}} ) ) $.

We will explain how to extract information about $\IndCoh_*$-theory on apaisant indschemes out of this functor in Remark \ref{rem-higher-comp-!-pull-indsch} below.

The desired functor (\ref{eqn-indcoh-out-ind-pro}) is constructed as follows.
\begin{itemize}
	\item We start from the symmetric monoidal functor 
	$$ \IndCoh: \on{Corr}(\on{Sch}_{\on{aft}}  )_{\on{all,all}} \to \on{DGCat}_{\on{cont}}  $$
	in \cite{GR-DAG1}. Using \emph{right} Kan extension along 
	$$ \on{Corr}(\on{Sch}_{\on{aft}}  )_{\on{all,all}}  \subset \on{Corr}( \on{Pro}(\on{Sch}_{\on{aft}} ) )_{\on{all,fafp}},$$
	we obtain a right-lax symmetric monoidal functor
	$$ \mathcal{I}nd\mathcal{C}oh^*: \on{Corr}( \on{Pro}(\on{Sch}_{\on{aft}} ) )_{\on{all,fafp}} \to \on{DGCat}_{\on{cont}}.$$
	It is pure formal to show that the composition 
	$$\on{Pro}(\on{Sch}_{\on{aft}} )\to \on{Corr}( \on{Pro}(\on{Sch}_{\on{aft}} ) )_{\on{all,fafp}} \to \on{DGCat}_{\on{cont}}$$
	is also the right Kan extension of the composition
	$$ \on{Sch}_{\on{aft}} \to \on{Corr}(\on{Sch}_{\on{aft}}  )_{\on{all,all}} \to \on{DGCat}_{\on{cont}}.$$

	\item Using \emph{left} Kan extension along 
	$$ \on{Corr}( \on{Pro}( \on{Sch}_{\on{aft}} ) )_{\on{all,fafp}}  \subset \on{Corr}( \on{Ind}( \on{Pro}(\on{Sch}_{\on{aft}} ) ) )_{\on{all,fafp}},$$
	we obtain the desired right-lax symmetric monoidal functor (\ref{eqn-indcoh-out-ind-pro}).
	It is pure formal to show that the composition 
	$$\on{Ind}( \on{Pro}(\on{Sch}_{\on{aft}} ) ) \to \on{Corr}( \on{Ind}( \on{Pro}(\on{Sch}_{\on{aft}} ) ) )_{\on{all,fafp}} \to \on{DGCat}_{\on{cont}}$$
	is also the left Kan extension of the composition
	$$  \on{Pro}(\on{Sch}_{\on{aft}} )  \to \on{Corr}( \on{Pro}(\on{Sch}_{\on{aft}} )  )_{\on{all,fafp}} \to \on{DGCat}_{\on{cont}}.$$
\end{itemize}
\end{construction}

\begin{remark} \label{rem-higher-comp-!-pull-indsch}

We have a fully faithful functor
$$ \iota: \on{IndSch}_{\on{apai}} \to \on{Ind}( \on{Sch}_{\on{apai}} ) \to  \on{Ind}( \on{Pro}( \on{Sch}_{\on{aft}} ) ).$$
Note that a uafp morphism $f:S\to T$ in $ \on{IndSch}_{\on{apai}} $ is sent to a fafp morphism by this functor. Note however that this functor does \emph{not} preserve fiber products\footnote{To see this, let $S,T,T'$ be apaisant schemes such that $f$ is afp and $S\times_T T'$ is not eventually coconnective, then the above functor sends $S\times_T T'$ to an object that is \emph{not} contained in $\on{Pro} (  \on{Sch}_{\on{aft}} )$, hence it can not be the fiber products of the images of $S,T,T'$, which are all contained in $\on{Pro} (  \on{Sch}_{\on{aft}} )$.}. Nevertheless, this functor preserves those Cartesian diagrams 
$$
\xymatrix{
	S' \ar[r]^-{\varphi} \ar[d]^-\phi & T'\ar[d]^-g \\
	S \ar[r]^-f & T
}
$$
in $\on{IndSch}_{\on{apai}} $ such that both $f$ and $\varphi$ are uafp. It will follow from construction that:
\begin{itemize}
	\item The composition
	$$ \on{IndSch}_{\on{apai}} \to  \on{Ind}( \on{Pro}( \on{Sch}_{\on{aft}} ) ) \to   \on{Corr} ( \on{Ind}( \on{Pro} (  \on{Sch}_{\on{aft}} ) ) )_{\on{all},\on{fafp}} \to \on{DGCat}_{\on{cont}}  $$
	can be identified with the functor
	$$ \IndCoh_*:  \on{IndSch}_{\on{apai}} \to \on{DGCat}_{\on{cont}} $$
	that encodes the $*$-pushforward functors for the $\IndCoh_*$-theory.

	\item Let $(\on{IndSch}_{\on{apai}})_{\on{uafp}} \subset (\on{IndSch}_{\on{apai}})_{\on{ind-afp}}$ be the \emph{simplicial subset} whose $1$-simplices are required to be uafp\footnote{By Warning \ref{warn-uafp-comp}, $(\on{IndSch}_{\on{apai}})_{\on{uafp}} $ is not a $\infty$-category}. Then the composition
	$$I^\op \to \on{Corr} ( \on{Ind}( \on{Pro} (  \on{Sch}_{\on{aft}} ) ) )_{\on{all},\on{fafp}} \to \on{DGCat}_{\on{cont}}  $$
	can be identified with
	$$ I^\op \to [(\on{IndSch}_{\on{apai}})_{\on{ind-afp}}]^\op \to \on{DGCat}_{\on{cont}} $$
	that encodes the $!$-pullback functors for the $\IndCoh_*$-theory\footnote{Warning: an ind-afp morphism $f:S\to T$ in $\on{IndSch}_{\on{apai}}$ that is not uafp might be sent to a fafp morphism in $\on{Ind}( \on{Pro}( \on{Sch}_{\on{aft}} ) )$. If this happens, we can obtain two functors $\IndCoh_*(T) \to \IndCoh_*(S)$. We do \emph{not} know if they coincide.}.	
\end{itemize}

Therefore the right-lax symmetric monoidal functor (\ref{eqn-indcoh-out-ind-pro}) provides the following. Let $\on{Cart}^{\bullet,\bullet}_{\on{all,uafp}}( \on{IndSch}_{\on{apai}} )$ be the bi-simplicial groupoid whose $(m,n)$-simplices is the full subgroupoid of $\on{Funct}( [m]\times[n], \on{IndSch}_{\on{apai}} )$
that consists of commutative diagrams
$$
\xymatrix{
	S_{0,0} \ar[d] \ar[r] & S_{1,0} \ar[r] \ar[d] & \cdots \ar[r] \ar[d] & S_{m,0} \ar[d] \\
	S_{0,1} \ar[d] \ar[r] & S_{1,1} \ar[r] \ar[d] & \cdots \ar[r] \ar[d] & S_{m,1} \ar[d] \\
 	\cdots	\ar[d] \ar[r] & \cdots \ar[r] \ar[d] & \cdots \ar[r] \ar[d] & \cdots \ar[d] \\
 	S_{0,n}  \ar[r] & S_{1,n} \ar[r]  & \cdots \ar[r]  & S_{m,n} 
}
$$
where all squares are Cartesian and all horizontal morphisms are uafp. The Cartesian products on $ \on{IndSch}_{\on{apai}}$ induce a symmetric monoidal structure on $\on{Cart}^{\bullet,\bullet}_{\on{all,uafp}}( \on{IndSch}_{\on{apai}} )$. Then we have a bi-simplicial right-lax\footnote{This right-lax symmetric monoidal structure is \emph{strict}. We will explain this in the next paragraph.} symmetric monoidal functor
$$ \IndCoh_*(-): \on{Cart}^{m,n}_{\on{all,uafp}}( \on{IndSch}_{\on{apai}} ) \to \on{Funct}( [m]^\op\times [n], \on{DGCat}_{\on{cont}})  $$
sending the above commutative diagram to
$$
\xymatrix{
	\IndCoh_*(S_{0,0}) \ar[d] & \IndCoh_*(S_{1,0}) \ar[l] \ar[d] & \cdots \ar[l] \ar[d] \ar[l] & \IndCoh_*(S_{m,0}) \ar[d] \ar[l] \\
	\IndCoh_*(S_{0,1}) \ar[d] & \IndCoh_*(S_{1,1}) \ar[l] \ar[d] & \cdots \ar[l] \ar[d] \ar[l] & \IndCoh_*(S_{m,1}) \ar[d] \ar[l] \\
 	\cdots \ar[d] & \cdots  \ar[l] \ar[d] & \cdots \ar[l] \ar[d] \ar[l] & \cdots  \ar[d] \ar[l] \\
 	\IndCoh_*(S_{0,n}) & \IndCoh_*(S_{1,n}) \ar[l]  & \cdots \ar[l] & \IndCoh_*(S_{m,n}) \ar[l] 
}
$$
where all horizontal functors are $!$-pullback functors and all vertical functors are $*$-pushforward functors. Moreover, if each $S_{i,j} \to S_{i,j'}$ is relatively apaisant, then all the squares in the above diagram are left adjointable along the vertical direction. 

The above right-lax symmetric monoidal structure means the following. For another similar diagram
$$
\xymatrix{
	T_{0,0} \ar[d] \ar[r] & T_{1,0} \ar[r] \ar[d] & \cdots \ar[r] \ar[d] & T_{m,0} \ar[d] \\
	T_{0,1} \ar[d] \ar[r] & T_{1,1} \ar[r] \ar[d] & \cdots \ar[r] \ar[d] & T_{m,1} \ar[d] \\
 	\cdots	\ar[d] \ar[r] & \cdots \ar[r] \ar[d] & \cdots \ar[r] \ar[d] & \cdots \ar[d] \\
 	T_{0,n}  \ar[r] & T_{1,n} \ar[r]  & \cdots \ar[r]  & T_{m,n} ,
}
$$
there is a canonical functor, i.e., the \emph{external tensor product functor}, from
$$
\xymatrix{
	\IndCoh_*(S_{0,0})\otimes \IndCoh_*(T_{0,0})  \ar[d] & \cdots \ar[l] \ar[d] \ar[l] & \IndCoh_*(S_{m,0}) \otimes \IndCoh_*(T_{m,0}) \ar[d] \ar[l] \\
 	\cdots \ar[d] & \cdots \ar[l] \ar[d] \ar[l] & \cdots  \ar[d] \ar[l] \\
 	\IndCoh_*(S_{0,n})  \otimes \IndCoh_*(T_{0,n})  & \cdots  \ar[l] & \IndCoh_*(S_{m,n}) \otimes \IndCoh_*(T_{m,n})  \ar[l] 
}
$$
to
$$
\xymatrix{
	\IndCoh_*(S_{0,0} \times T_{0,0})  \ar[d] & \cdots \ar[l] \ar[d] \ar[l] & \IndCoh_*(S_{m,0} \times T_{m,0}) \ar[d] \ar[l] \\
 	\cdots \ar[d] & \cdots \ar[l] \ar[d] \ar[l] & \cdots  \ar[d] \ar[l] \\
 	\IndCoh_*(S_{0,n} \times T_{0,n})  & \cdots  \ar[l] & \IndCoh_*(S_{m,n}\times T_{m,n}).  \ar[l] 
}
$$
By \cite[Proposition 6.35.2, Proposition 6.36.4]{raskin2020homological}, this right-lax symmetric monoidal structure is strict.
\end{remark}

\subsection{Ind-coherent sheaves on renormalizable prestacks: definitions}

In the following subsections, we study the $\IndCoh_*$-theory on nice enough prestacks such as the double quotient $\Arc G \backslash \Loop G / \Arc G$.

\begin{definition} (C.f. \cite[$\S$ 6.21-6.22]{raskin2020homological}) A convergent prestack $S\in \on{PreStk}_{\on{conv}}$ is \emph{weakly renormalizable}, or $S\in \on{PreStk}_{\on{w.ren}}$, if there exists a faithfully flat and relatively apaisant map $T\to S$ with $T$ being an apaisant indscheme. We refer such a map $T\to S$ as an \emph{flat apaisant cover} for $S$. 
\end{definition}

\begin{example} The (fppf) double quotient $\Arc G \backslash \Loop G / \Arc G$ is a weakly renormalizable prestack.

\end{example}

\begin{remark} \label{rem-cov-contractible-w.ren}
Consider the category $\on{Cov}(S)$ of flat apaisant covers of $S$, where a morphism from a cover $S_0\to S$ to $S_0' \to S$ is a faithfully flat and relatively apaisant map $S_0\to S_0'$ defined over $S$. Note that the fiber product of two flat apaisant covers is also a flat apaisant cover. Hence $\on{Cov}(S)$ is weakly contractible.
\end{remark}

\begin{definition} A map $f:S\to T$ in $\on{PreStk}_{\on{w.ren}}$ is \emph{apaisant indschematic} if for any relatively apaisant map $T' \to T$ with $T'\in \on{IndSch}_{\on{apai}}$, the fiber product $S\times_T T'$ is contained in $\on{IndSch}_{\on{apai}}$.
\end{definition}

\begin{remark} If $S\in \on{IndSch}_{\on{apai}}$, then $f$ is apaisant indschematic. 
\end{remark}

\begin{remark} If $f$ is relatively apaisant, then $f$ is apaisant indschematic.
\end{remark}

\begin{remark} Let $f:S\to T$ be an apaisant indschematic morphism and $g:T' \to T$ be a relatively apaisant morphism. Then $S\times_T T'$ is weakly renormalizable and $S\times_T T' \to T'$ is apaisant indschematic.
\end{remark}

\begin{warning} We do \emph{not} know if the class of apaisant indschematic maps is stable under general base-changes.
\end{warning}

\begin{construction} (C.f. \cite[Definition 6.24.1]{raskin2020homological}) \label{def-push-along-indshematic}

Recall we have a symmetric monoidal functor
$$ \IndCoh_*:   \on{Corr} (  \on{IndSch}_{\on{apai}} )_{ \on{all,rel.apai} }  \to \on{DGCat}_{\on{cont}} $$
that encodes the $*$-pullback and $*$-pushforward functors for the $\IndCoh_*$-theory on apaisant indschemes. Using right Kan extension along the fully faithful functor
$$ \on{Corr} (  \on{IndSch}_{\on{apai}} )_{ \on{all,rel.apai} } \to \on{Corr} (  \on{PreStk}_{\on{w.ren}} )_{ \on{apai.indsch,rel.apai} } ,$$
we obtain a right-lax symmetric monoidal functor
$$  \IndCoh_*: \on{Corr} (  \on{PreStk}_{\on{w.ren}} )_{ \on{apai.indsch,rel.apai} } \to \on{DGCat}_{\on{cont}}.$$
In other words, we extended the $\IndCoh_*$-theory to weakly renormalizable prestacks such that 
\begin{itemize}
	\item we have $*$-pullback functors along relatively apaisant morphisms,
	\item we have $*$-pushforward functors along apaisant indschematic morphisms,
	\item we have base-change isomorphisms between the above $*$-pullback and $*$-pushforward functors.
	\item we have external tensor product functors compatible with the $*$-pullback and $*$-pushforward functors.
\end{itemize}

By \cite[Chapter 8, Theorem 6.1.5]{GR-DAG1}, when restricted to horizontal morphisms, the obtained functor
$$  [(\on{PreStk}_{\on{w.ren}} )_{ \on{rel.apai} }]^\op \to \on{DGCat}_{\on{cont}} $$
is also the right Kan extension of the functor
$$  [(\on{IndSch}_{\on{apai}} )_{ \on{rel.apai} }]^\op \to \on{DGCat}_{\on{cont}} $$
that only encodes the $*$-pullback functors. In other words, for a weakly renormalizable prestack $S$, we have
$$ \IndCoh_*(S) \simeq \lim_{*\on{-pull}} \IndCoh_*(U),$$
where the limit is taken for all relatively apaisant morphisms $U\to S$ with $U \in (\on{IndSch}_{\on{apai}})_{\on{rel.apai}}^\op$.

By construction, for any relatively apaisant morphism $f:S\to T$ in $\on{PreStk}_{\on{w.ren}}$, the functor $f^*$ is left adjoint to $f_*$. Also, for a Cartesian diagram
$$
\xymatrix{
	S' \ar[r]^-{\varphi} \ar[d]^-\phi & T'\ar[d]^-g \\
	S \ar[r]^-f & T
}
$$
contained in $\on{PreStk}_{\on{w.ren}}$, if $f$ is relatively apaisant (hence so is $\varphi$) and $g$ is apaisant indschematic (hence so is $\phi$), then the Beck-Chevalley natural transformation
$$ f^{*}\circ g_* \to \phi_* \circ \varphi^{*} $$
is an isomorphism.
\end{construction}

\begin{remark} \label{rem-qca-indcoh*}
When $S$ is a QCA stack, by Theorem \ref{theorem-flat-descent-indcoh-w-ren} below, the category $\IndCoh_*(S)$ coincides with $\IndCoh(S)$ defined in \cite{drinfeld2013some}. Hence we write $\IndCoh(S)$ for $\IndCoh_*(S)$ in this case.
\end{remark}

\begin{warning} For general $S$, the category $\IndCoh_*(S)$ is not as well-behaved as in the case when $S$ is an apaisant indscheme. For instance, it is not compactly generated. Also, we do not have a \emph{continuous} $*$-pushforward functor along maps like $\mathbb{B} G(O) \to \on{pt}$. We will use the \emph{renormalization} precedure to obatin a more robust $\IndCoh_{*,\on{ren}}$-theory later.
\end{warning}

\begin{construction} \label{const-t-structure-w-ren-indcoh}
For $S\in \on{PreStk}_{\on{w.ren}}$, there exists a unique t-structure on $\IndCoh_*(S)$ such that for any flat and relatively apaisant morphism $f:U\to S$ with $U\in \on{IndSch}_{\on{apai}}$, the functor $f^{*}$ is t-exact. Equivalently, it is the unique t-exact such that for a given flat apaisant cover $f:U\to S$, the functor $f^{*}$ is t-exact.
\end{construction}

\begin{remark} It is easy to see the above t-structure is compatible with filtered colimits and right separated.
\end{remark}

\begin{lemma} \label{lem-*-push-*-pull-t-exact-w-ren}
Let $f:S\to T$ be an apaisant indschematic map in $\on{PreStk}_{\on{w.ren}}$, then $f_*: \IndCoh_*(S) \to \IndCoh_*(T)$ is left t-exact. If $f$ is further assumed to be relatively apaisant, then $f^*$ is right t-exact, and is left t-exact up to a cohomological shift.
\end{lemma}

\proof Follows from Lemma \ref{lem-*-push-left-t-exact-indsch} and Lemma \ref{lem-*-pull-indsch}.

\qed

\begin{lemma}  \label{lem-box-tensor-left-t-exact-w.ren}
Let $S$ and $T$ be in $\on{PreStk}_{\on{w.ren}}$. For $\mathcal{F}\in \IndCoh_*(S)^{\ge 0}$ and $\mathcal{G}\in \IndCoh_*(T)^{\ge 0}$, the external tensor product $\mathcal{F}\boxtimes \mathcal{G}$ is contained in $\IndCoh_*(S\times T)$.
\end{lemma}

\proof By Construction \ref{const-t-structure-w-ren-indcoh}, Construction \ref{const-t-structure-indcoh-indsch} and Remark \ref{rem-t-structure-apai-sch}, we can reduce to the case when $S,T\in\, ^{>-\infty}\on{Sch}_{\on{aft}}$, and then reduce to a similar obvious claim about $\on{QCoh}$.

\qed

\subsection{Ind-coherent sheaves on renormalizable prestacks: descent}

In this subsection, we show that the $\IndCoh_*$-theory on weakly renormalizable prestacks satisfies descents for flat apaisant covers.

\begin{theorem} (C.f. \cite[Theorem 6.25.1]{raskin2020homological}) \label{theorem-flat-descent-indcoh-w-ren}
Let $f:T\to S$ be a relatively apaisant and faithfully flat map in $\on{PreStk}_{\on{w.ren}}$, then $*$-pullback induces an equivalence
$$ \IndCoh_*(S) \to \on{Tot}_{\on{semi}} ( \IndCoh_*(T^{\times_S (\bullet+1)}) ).$$
\end{theorem}

\proof The proof in \cite{raskin2020homological} is correct as long as we replace \emph{Step 1} of that proof by Proposition \ref{theorem-flat-descent-indcoh-indsch}. In general, under mild assumptions, a functor satisfying certain descents formally implies its right Kan extension satisfies similar descents. See Proposition \ref{prop-general-descent-RKE} below for a precise statement.

\qed

\begin{proposition}\label{prop-general-descent-RKE}
 Let $\mathcal{D}$ be any category and $\mathcal{C}\to \mathcal{D}$ be a full subcategory. Let (P) be a property on morphisms in $\mathcal{D}$ that is stable under compositions and contains all isomorphisms. Let $\mathcal{D}_P$ (resp. $\mathcal{C}_P$) be the corresponding 1-full subcategory of $\mathcal{D}$ (resp. $\mathcal{C}$). Let (Q) be a property on morphisms in $\mathcal{D}$ such that
\begin{itemize}
	\item[(a)] it is stronger than (P), i.e., any (Q)-morphism is a (P)-morphism,
	\item[(b)] it is stable under compositions and contains all isomorphisms,
	\item[(c)] it is stable under base-changes, i.e., if $T\to S$ is a (Q)-morphism and $S'\to S$ is any morphism, then the fiber prodcut $T\times_S S'$ is well-defined in $\mathcal{D}$ and $T\times_S S' \to S'$ is a (Q)-morphism,
	\item[(d)] if $T\to S$ is a (Q)-morphism and $S$ is contained in $\mathcal{C}$, then so is $T$,
	\item[(e)] if $T\to S$ is a (P)-morphism with $T\in \mathcal{C}$ and $S'\to S$ is a (Q)-morphism in $\mathcal{D}$, then the fiber prodcut $T\times_S S'$ is well-defined in $\mathcal{D}$ and $T\times_S S' \to S'$ is a (P)-morphism,
	\item[(f)] for any $S\in \mathcal{D}$, there is a (Q)-morphism $S'\to S$ with $S'\in \mathcal{C}$.
\end{itemize}
Now suppose $I:\mathcal{C}_P^\op \to \mathcal{E}$ is any functor satisfying (Q)-descent, i.e., such that for any (Q)-morphism $T\to S$ in $\mathcal{C}$, the canonical morphism\footnote{The fiber product $T^{\times_S (\bullet+1)}$ is taken in $\mathcal{D}$. By assumption (d), it is actually contained in $\mathcal{C}$. However, it is \emph{not} the fiber product in $\mathcal{C}_P$ or $\mathcal{C}_Q$.}
$$ I(S) \to \on{Tot}_{\on{semi}} (  I(T^{\times_S (\bullet+1)}) ) $$
is an isomorphism. Suppose the right Kan extension $\on{RKE}(I): \mathcal{D}_P^\op \to \mathcal{E}$ of $I$ along $\mathcal{C}_P^\op \to \mathcal{D}_P^\op$ exists. Then $\on{RKE} (I)$ also satisfies (Q)-descent.
\end{proposition}

\proof Our argument is essentially borrowed from \emph{Step 2} of the proof of \cite[Theorem 6.25.1]{raskin2020homological}. 

We abuse notation by writing $I$ for $\on{RKE} (I)$. For any (Q)-morphism $T\to S$ in $\mathcal{D}$, we need to prove
$$F: I (S) \to \lim_{[m] \in \Delta_{inj}}  I  (T^{\times_S (m+1)}) $$
is an isomorphism.  By (b) and (c), the semi-cosimplicial diagram $T^{\times_S (\bullet+1)}$ is contained in $\mathcal{D}_Q$. Recall
$$ I(S) \simeq \lim J,$$
where $J: (\mathcal{C}_{P}^\op)_{S/} \to \mathcal{E} $ sends a (P)-morphisms $U\to S$ to $I(U)$.

We first prove in the case when $T \in \mathcal{C}$. By assumption (d), the semi-cosimplicial diagram $T^{\times_S (\bullet+1)}$ is contained in $\mathcal{C}_Q$. Hence by assumption (e), we have canonical functors
\begin{eqnarray*}
\mathbf{B}: (\mathcal{C}_{P}^\op)_{S/} \times \Delta_{inj} \to (\mathcal{C}_{P}^\op)_{S/}, & & (U,[m])\mapsto U\times_S T^{\times_S (m+1)}  ,\\
\mathbf{T}:  \Delta_{inj} \to (\mathcal{C}_{P}^\op)_{S/}, & & [m]\mapsto T^{\times_S (m+1)} ,
\end{eqnarray*}
and natural transformations $ \alpha: \mathbf{T} \circ \on{pr}_2 \to \mathbf{B}$, $\beta:\on{pr}_1 \to \mathbf{B}$. Note that the morphism $F$ can be identified with the canonical morphism
$$ \on{res}_{\mathbf{T}}: \lim J \to \lim J\circ \mathbf{T} $$
induced by restricting the limit diagrams along the functor $\mathbf{T}$. Since $I$ satisfies (Q)-descent on $\mathcal{C}$, for any $U\in (\mathcal{C}_{P}^\op)_{S/}$, we have 
\begin{equation} \label{eqn-proof-general-descent-RKE-0}
 J(U) \simeq \lim_{ [m]\in \Delta_{inj} }  J( \mathbf{B}(U,[m])).
 \end{equation}
Hence an isomorphism
\begin{equation} \label{eqn-proof-general-descent-RKE-1}
 \lim J \to \lim J \circ \mathbf{B}.  
 \end{equation}
Consider the following composition
\begin{eqnarray*}
G: \lim_{ [m]\in \Delta_{inj} } J\circ \mathbf{T}  & \xrightarrow{\on{res}_{ \on{pr}_2 }} & \lim J\circ \mathbf{T} \circ \on{pr}_2  \\
& \xrightarrow{ \lim(\mathbf{Id}_J \bigstar \alpha) } & \lim J\circ  \mathbf{B}  \\
& \to & \lim J,
\end{eqnarray*}
where 
\begin{itemize}
	\item the first morphism is restricting the limit diagrams along $\on{pr}_2:  (\mathcal{C}_{P}^\op)_{S/} \times \Delta_{inj}  \to  \Delta_{inj} $,
	\item the second morphism is induced by the natural transformation $\alpha:\mathbf{T} \circ \on{pr}_2 \to \mathbf{B}$,
	\item the third morphism is the inverse of (\ref{eqn-proof-general-descent-RKE-1}).
\end{itemize}
We claim $F$ and $G$ are inverse to each other, and in particular, they are isomorphism. Note that 
$$ \lim(\mathbf{Id}_J \bigstar \alpha) \circ \on{res}_{ \on{pr}_2 } \circ F \simeq \lim(\mathbf{Id}_J \bigstar \alpha) \circ \on{res}_{  \on{pr}_2 } \circ \on{res}_{ \mathbf{T} } \simeq  \lim(\mathbf{Id}_J \bigstar \alpha) \circ \on{res}_{ \mathbf{T} \circ  \on{pr}_2 } \simeq \on{res}_{ \mathbf{B} }.$$
Also note that (\ref{eqn-proof-general-descent-RKE-1}) is just $ \lim(\mathbf{Id}_J \bigstar \beta) \circ \on{res}_{\on{pr}_1}$, which is also isomorphic to $\on{res}_{ \mathbf{B} }$. Combining these, we get $G\circ F\simeq \on{Id}$. Now we show $F\circ G\simeq \on{Id}$. By above, we have
$$F\circ G \simeq   \on{res}_{\mathbf{T}}  \circ \res_{\mathbf{B}} ^{-1} \circ   \lim(\mathbf{Id}_J \bigstar \alpha) \circ \on{res}_{ \on{pr}_2 }.$$
Note that (\ref{eqn-proof-general-descent-RKE-0}) also induces an isomorphism
$$ J\circ \mathbf{T}([n]) \simeq  \lim_{ [m]\in \Delta_{inj} }  J( \mathbf{B}(  \mathbf{T}([n]),[m])),$$
hence an isomorphism
\begin{equation} \label{eqn-proof-general-descent-RKE-2}
\lim J\circ \mathbf{T} \to \lim J\circ \mathbf{T} \circ \on{join} 
\end{equation}
where 
\begin{eqnarray*}
\on{join}:  \Delta_{inj} \times \Delta_{inj} \to  \Delta_{inj}, & & ([n],[m])\mapsto [m+n+1]
\end{eqnarray*}
is the concatenation functor, and we are using the fact
$$ \mathbf{T}\circ \on{join} \simeq  \mathbf{B} \circ ( \mathbf{T}\times \mathbf{Id} ) . $$
As before, we can show (\ref{eqn-proof-general-descent-RKE-2}) is equivalent to $\on{res}_{\on{join}}$.

Hence
$$ \on{res}_{ \on{join}  } \circ \on{res}_{\mathbf{T}} \simeq \on{res}_{\mathbf{T}\times \mathbf{Id}} \circ \res_{\mathbf{B}} $$
and therefore
$$   \on{res}_{\mathbf{T}}  \circ \res_{\mathbf{B}} ^{-1} \simeq \on{res}_{ \on{join}  } ^{-1}\circ \on{res}_{\mathbf{T}\times \mathbf{Id}}.$$
Hence 
$$F\circ G \simeq    \on{res}_{ \on{join}  } ^{-1}\circ \on{res}_{\mathbf{T}\times \mathbf{Id}} \circ  \lim(\mathbf{Id}_J \bigstar \alpha) \circ \on{res}_{ \on{pr}_2 } \simeq \on{res}_{ \on{join}  } ^{-1}\circ \lim(\mathbf{Id}_J \bigstar \alpha \bigstar \mathbf{Id}_{ \mathbf{T}\times \mathbf{Id} } ) \circ   \on{res}_{\mathbf{T}\times \mathbf{Id}}  \circ \on{res}_{ \on{pr}_2 }. $$
Note that $\on{pr}_2 \circ ( \mathbf{T}\times\mathbf{Id} )$ is the second projection map
$$p_2: \Delta_{inj} \times \Delta_{inj} \to  \Delta_{inj},$$
and the natural transformation 
$$ \alpha \star \mathbf{Id}_{ \mathbf{T}\times \mathbf{Id} }: \mathbf{T} \circ p_2 \to \mathbf{T}\circ \on{join}$$ is induced by a natural transformation $\beta: p_2\to \on{join}$. Hence 
$$\lim(\mathbf{Id}_J \bigstar \alpha \bigstar \mathbf{Id}_{ \mathbf{T}\times \mathbf{Id} } ) \circ   \on{res}_{\mathbf{T}\times \mathbf{Id}}  \circ \on{res}_{ \on{pr}_2 }\simeq \lim( \mathbf{Id}_{J\circ \mathbf{T}} \bigstar \beta  ) \circ \on{res}_{p_2} \simeq \on{res}_{\on{join}}. $$
This shows $F\circ G \simeq \on{Id}$ and finishes the proof in the case $T \in \mathcal{C}$.

For the general case, let $S'\to S$ be as in assumption (f). We have already proved $I$ satisfies descent along $S' \to S$, $S' \times_S T^{\times_S (\bullet+1)} \to T^{\times_S (\bullet+1)} $ and $S'\times_S T\to S'$. These formally imply $I$ satisfies descent along $T\to S$.

\qed

\begin{proposition} \label{prop-indcoh-zariski-w-ren}
Let $S_1 \to S$ and $S_2\to S$ be two schematic open embeddings in $\on{PreStk}_{\on{w.ren}}$ such that $S_1\sqcup S_2 \to S$ is surjective. Then $*$-pullbacks induce an equivalence
$$
 \IndCoh_*(S) \to \IndCoh_*(S_1) \times_{ \IndCoh_*(S_1\cap S_2) } \IndCoh_*(S_2).
 $$
\end{proposition}

\proof Follows from Proposition \ref{prop-indcoh-zariski-indsch} and Theorem \ref{theorem-flat-descent-indcoh-w-ren}.

\qed

\subsection{Ind-coherent sheaves on renormalizable prestacks: more \texorpdfstring{$*$}{*}-pushforwards}

By Lemma \ref{lem-*-push-*-pull-t-exact-w-ren}, for any apaisant indschematic map $f:S\to T$, we have a functor
 $$ (f_*)^+: \IndCoh_*(S)^+ \to \IndCoh_*(T)^+.$$
If $f$ is relatively apaisant, then $(f_*)^+$ has a left adjoint $(f^*)^+$.

By Lemma \ref{lem-box-tensor-left-t-exact-w.ren}, the above $*$-pushforward functors can be assembled into a right-lax symmetric monoidal functor
$$ (\IndCoh_*)^+: (\on{PreStk}_{\on{w.ren}})_{\on{apai.indsch}} \to \on{Vect}^+\on{-mod}(\on{Cat}) ,\; S\mapsto \IndCoh_*(S)^+.$$

For the purpose of this paper, we need to extend the definition of $(f_*)^+$ to non-indschematic maps such as $\mathbb{B} B(O) \to \mathbb{B} T(O)$. The following construction is similar to that in \cite[$\S$ 3.7]{drinfeld2013some}.

\begin{proposition-construction}  \label{prop-general-push-w-ren}
\begin{itemize}
	\item[(1)] The right-lax symmetric monoidal functor 
$$ (\IndCoh_*)^+: (\on{PreStk}_{\on{w.ren}})_{\on{apai.indsch}} \to \on{Vect}^+\on{-mod}(\on{Cat}) ,\; S\mapsto \IndCoh_*(S)^+$$
	 can be uniquely extended to a right-lax symmetric monoidal functor
$$ (\IndCoh_*)^+: \on{PreStk}_{\on{w.ren}} \to \on{Vect}^+\on{-mod}(\on{Cat})  $$
satifying the following condition. Let $f:S\to T$ be a morphism in $\on{PreStk}_{\on{w.ren}}$ and $q_0:S_0\to S$ be any flat apaisant cover. Consider the corresponding Cech nerve $S_\bullet$ and projections $q_\bullet:S_\bullet \to S$. Then the canonical natural transformation
	$$ f_*^+ \to \on{Tot}_{\on{semi}}  [((  f\circ q_\bullet )_*)^+ \circ (q_\bullet^*)^+ ] $$
	is required to be an equivalence\footnote{Note that each $f\circ q_\bullet$ is apaisant indschematic. Hence the well-definedness of the RHS does not depend on the proposition.}.

	\item[(2)] For any $f:S\to T$, the functor $f_*^+$ defined above is left t-exact.

	\item[(3)] The retriction of $f_*^+$ on $\IndCoh_*(S)^{\ge 0}$ preserves filtered colimits.

\end{itemize}

\end{proposition-construction}

\proof By Remark \ref{rem-cov-contractible-w.ren}, $\on{Cov}(S)$ is weakly contractible. Hence to prove (1), we only need to check
\begin{itemize}
	\item[(i)] Let $p_0: S_0' \to S$ and $q_0:S_0\to S$ be two flat apaisant covers, and $S_0' \to S_0$ be a morphism between them. Then we need to show the canonical natural transformation
	$$  \on{Tot}_{\on{semi}} [ ((  f\circ q_\bullet )_*)^+ \circ (q_\bullet^*)^+] \to \on{Tot}_{\on{semi}}  [((  f\circ p_\bullet )_*)^+ \circ (p_\bullet^*)^+]  $$
	is an equivalence.

	\item[(ii)] Let $ S \xrightarrow{f}  T \xrightarrow{g}  U$ be a chain in $\on{PreStk}_{\on{w.ren}}$. Let $q_0:S_0 \to S$ and $r_0: T_0\to T$ be flat apaisant covers. Define $S_0:=S_0\times_{T} T_0$ and let $p_0: S_0'\to S$ be the projection. Consider the corresponding Cech nerves. Note that the square in the following diagram is Cartesian.
	\[
	\xymatrix{
		S_\bullet' \ar[d] \ar[dr]^-{p_\bullet} \ar[rr]
	&&	T_\bullet  \ar[d]^-{r_\bullet} \\
	  	S_\bullet \ar[r]^-{q_\bullet}
	& 	S \ar[r]^-{f}
	&	T \ar[r]^-{g}
	&	U
	}
	\]
	Then we need to show the following composition
	\begin{eqnarray*}
	 & &  \underset{m}{\on{Tot}_{\on{semi}}}  [((  g\circ r_m )_*)^+ \circ (r_m^*)^+ ]\circ 
	  \underset{n}{\on{Tot}_{\on{semi}}} [ ((  f\circ q_n)_*)^+ \circ (q_n^*)^+ ] \\
	  &\to& \underset{m}{\on{Tot}_{\on{semi}}}\underset{n}{\on{Tot}_{\on{semi}}}  [ ((  g\circ r_m )_*)^+ \circ (r_m^*)^+\circ ((  f\circ q_n )_*)^+ \circ (q_n^*)^+] \\
	  &\simeq& \on{Tot}_{\on{semi}}  [((  g\circ r_\bullet )_*)^+ \circ (r_\bullet^*)^+\circ ((  f\circ q_\bullet )_*)^+ \circ (q_\bullet^*)^+ ]\\
	  & \simeq & 
	  \on{Tot}_{\on{semi}}  [((  g\circ f\circ p_\bullet )_*)^+ \circ (p_\bullet^*)^+]
	  \end{eqnarray*}
	is an equivalence\footnote{In above, the second natural transformation is an equivalence because semi-totalization is a sifted colimit; the last equivalence is due to the base-change isomorphisms.}.
\end{itemize}

Let $\mathcal{F}\in \IndCoh_*(S)^{\ge 0}$. To prove (i), we only need to show the canonical morphism
$$ \lim_{h\in \on{Cov}(S)^\op}  [ ((  f\circ h )_*)^+ \circ (h^*)^+ (\mathcal{F})] \to \on{Tot}_{\on{semi}} [ ((  f\circ q_\bullet )_*)^+ \circ (q_\bullet^*)^+ (\mathcal{F})] $$
is an equivalence. We apply Proposition \ref{prop-general-descent-RKE} to the following case. Let $\mathcal{D}$ be the category of flat and relatively apaisant morphisms $h:R\to S$ with $R\in\on{PreStk}_{\on{w.ren}}$, where a morphism from $ R\to S$ to $R' \to S$ is \emph{any} morphism $R\to R'$ defined over $S$. A morphism $R\to R'$ is a (P)-morphism if it is relatively apaisant and faithfully flat. Any (P)-morphism is a (Q)-morphism. Let $\mathcal{C}\subset \mathcal{D}$ be the full subcategory of those $h:R\to S$ such that $R\in \on{IndSch}_{\on{apai}}$. Consider the functor
$$  \on{Cov}(S)^\op =  \mathcal{C}_P^{\op} \to \IndCoh_*(T)^{\ge 0},\; h\mapsto ((  f\circ h )_*)^+ \circ (h^*)^+ (\mathcal{F}).$$
To prove the desired equivalence, we only need to show the right Kan extension of the above functor along $\mathcal{C}_P^{\op} \subset \mathcal{D}_P^{\op}$ satisfies (Q)-descent. Hence by Proposition \ref{prop-general-descent-RKE}, we only need to show the above functor satisfies (Q)-descent. In other words, we have reduced (i) to the following statement:
\begin{itemize}
	\item[(iii)] Let $r_0: R_0\to R$ be a morphism in $\on{Cov}(S)$. Let $h_0: R_0 \to S$ to $h: R\to S$ be the corresponding covering maps. Consider the Cech nerve $r_\bullet: R_\bullet\to R$. Then we need to show the canonical natural transformation
	$$  ((f\circ h)_*)^+ \circ (h^*)^+ \to \on{Tot}_{\on{semi}}  [((f\circ h \circ r_\bullet)_*)^+ \circ ( ( h \circ r_\bullet )^*)^+] $$
	is an equivalence.
\end{itemize}
By Thoerem \ref{theorem-flat-descent-indcoh-w-ren}, we have
$$  (h^{*})^+(\mathcal{F})  \simeq  
 \tx{Tot}_{\on{semi}}   [((r_\bullet)_*)^+ \circ ((h\circ r_\bullet)^{*})^+ (\mathcal{F})] .$$
Hence it remains to show
\begin{equation} \label{eqn-prop-indcoh*-*push-naive-1} 
((f\circ h)_*)^+ \circ   \tx{Tot}_{\on{semi}}  [ ((r_\bullet)_*)^+ \circ ((h\circ r_\bullet)^{*})^+ ](\mathcal{F})  \to   \tx{Tot}_{\on{semi}} [ ((f\circ h)_*)^+\circ  ((r_\bullet)_*)^+ \circ ((h\circ r_\bullet)^{*})^+] (\mathcal{F})
\end{equation}
is an equivalence. Note that both sides are contained in $\IndCoh_*(T)^{\ge 0}$. Hence by right-separatedness, we only need to show it induces isomorphisms on cohomologies $H^m$. This allows us to replace $\tx{Tot}_{\on{semi}}$ in (\ref{eqn-prop-indcoh*-*push-naive-1}) by the $(m+2)$-truncated semi-totalizations. Now the desired equivalence is manifest because the stable functor $((f\circ h)_*)^+$ preserves any finite limit. This proves (iii) and therefore finishes the proof of (i).

The proof of (ii) is similar to that of (iii). Namely, we only need to show for any $\mathcal{F} \in \IndCoh_*(S)^{\ge 0}$, the canonical morphism
$$  ((  g\circ r_m )_*)^+ \circ (r_m^*)^+  \circ \underset{n}{\on{Tot}_{\on{semi}}}  [((  f\circ q_n)_*)^+ \circ (q_n^*)^+ ] (\mathcal{F}) \to \underset{n}{\on{Tot}_{\on{semi}}} [((  g\circ r_m )_*)^+ \circ (r_m^*)^+  \circ ((  f\circ q_n)_*)^+ \circ (q_n^*)^+] (\mathcal{F}) $$
is an equivalence. We can take adavantage of the t-structures to replace $\on{Tot}_{\on{semi}}$ by a finite limit, which makes the claim manifest.

This finishs the proof of (1).

The statement (2) follows from construction.

The proof of (3) is again similar to that of (iii). Namely, we want to show certain filtered colimits commute with certain semi-totalizations in $\IndCoh_*(T)$. We can take adavantage of the t-structures to replace $\on{Tot}_{\on{semi}}$ by a finite limit. Then we are done because filtered colimits commute with finite limits in any presentable stable categories.

\qed

\begin{definition} \label{def-*-pull-weak-ren}
A map $f:S\to T$ in $\on{PreStk}_{\on{w.ren}}$ is called \emph{$*$-pullable} if there exists a flat apaisant cover $S'\to S$ such that the composition $S'\to T$ is relatively apaisant.
\end{definition}

\begin{remark} Let $f:S\to T$ be an apaisant indschematic morphism and $g:T' \to T$ be a $*$-pullable morphism. Then $S\times_T T'$ is weakly renormalizable and $S\times_T T' \to S$ is $*$-pullable.
\end{remark}

\begin{lemma} 
\begin{itemize}
	\item[(1)] If $f:S\to T$ is $*$-pullable, then the functor $(f_*)^+: \IndCoh_*(S)^+ \to \IndCoh_*(T)^+$ has a left adjoint $(f^*)^+$. 

	\item[(2)] For a Cartesian diagram
$$
\xymatrix{
	S' \ar[r]^-{\varphi} \ar[d]^-\phi & T'\ar[d]^-g \\
	S \ar[r]^-f & T
}
$$
contained in $\on{PreStk}_{\on{w.ren}}$ such that $f$ and $\varphi$ are $*$-pullable, the Beck-Chevalley natural transformation (for $\IndCoh_*$-theory)
$$ f^{*}\circ g_* \to \phi_* \circ \varphi^{*} $$
is an isomorphism.
\end{itemize}
\end{lemma}

\proof Let $q_0: S_0\to S$ be a flat apaisant cover such that $f\circ q_0$ is relatively apaisant. Let $q_\bullet: S_\bullet\to S$ be the corresponding Cech nerve. Consider the functor
$$ \IndCoh_*(T)^+ \to \on{Tot}_{\on{semi}} \IndCoh_*(S_\bullet)^+,\; \mathcal{F} \mapsto (( f\circ q_\bullet )^*)^+ \mathcal{F}.$$
By Theorem \ref{theorem-flat-descent-indcoh-w-ren}, the target category is equivalent to $\IndCoh_*(S)^+$. It follows from definitions that the corresponding functor $\IndCoh_*(T) \to \IndCoh_*(S)$ is left adjoint to $(f_*)^+$. This proves (1).

To prove (2), consider the following pullback diagram:
$$
\xymatrix{
	S_\bullet \times_S S' \ar[r] \ar[d] & S' \ar[r]^-{\varphi} \ar[d]^-\phi & T'\ar[d]^-g \\
	S_\bullet \ar[r]^-{q_\bullet} &  S \ar[r]^-f & T.
}
$$
It follows that we only need to prove the base-change isomorphisms for the left square and the outer square. Note that both $q_\bullet$ and $f\circ q_\bullet$ are relatively apaisant. Hence we only need to prove (2) in the case when $f$ is relatively apaisant. Now let $p_0: T'_0\to T'$ be any flat apaisant cover and $p_\bullet: T'_\bullet \to T'$ be the corresponding Cech nerves. Consider the following pullback diagram:
$$
\xymatrix{
	S'\times_{T'}T'_\bullet \ar[r]^-{\varphi_\bullet} \ar[d]^-{r_\bullet} & T'_\bullet\ar[d]^-{p_\bullet} \\
	S' \ar[r]^-{\varphi} \ar[d]^-\phi & T'\ar[d]^-g \\
	S \ar[r]^-f & T.
}
$$
Note that $r_0: S'\times_{T'} T_0' \to S'$ is also a flat apaisant cover. Hence by definition, we only need to show the following composition is an equivalence:
\begin{eqnarray*}
	(f^*)^+ \circ (g_*)^+
	& \simeq & (f^*)^+ \circ \on{Tot}_{\on{semi}} [((g\circ p_\bullet )_* )^+ \circ( p_\bullet^*)^+   ] \\
	& \to & \on{Tot}_{\on{semi}} [(f^*)^+ \circ ((g\circ p_\bullet )_* )^+ \circ( p_\bullet^*)^+ ]   \\
	& \simeq & \on{Tot}_{\on{semi}} [((\phi \circ r_\bullet)_*)^+ \circ ((r_\bullet )^* )^+ \circ( \varphi^*)^+ ] \\
	& \simeq & (\phi_*)^+ \circ (\varphi^*)^*,
\end{eqnarray*}
where the equivalence between the second and the third lines are due to the base-change isomorphisms in Construction \ref{def-push-along-indshematic}. In other words, we only need to show
$$ (f^*)^+ \circ \on{Tot}_{\on{semi}} [((g\circ p_\bullet )_* )^+ \circ( p_\bullet^*)^+   ] \to  \on{Tot}_{\on{semi}} [(f^*)^+ \circ ((g\circ p_\bullet )_* )^+ \circ( p_\bullet^*)^+ ]  $$
is an equivalence. But this can be proved by using the same trick as in the proof of Proposition \ref{prop-general-push-w-ren}.

\qed

\subsection{Ind-coherent sheaves on renormalizable prestacks: \texorpdfstring{$!$}{!}-pullbacks}

In this subsection, we study the $!$-pullback functors for $\IndCoh_*$-theory on weakly renormalizable prestacks.

\begin{definition} An apaisant indschematic morphism $f:S \to T$ in $\on{PreStk}_{\on{w.ren}}$ is \emph{uafp} if for any relatively apaisant map $T'\to T$ with $T'\in \on{IndSch}_{\on{apai}}$, the map $S\times_T T' \to T'$ is a uafp morphism in $\on{IndSch}_{\on{apai}}$.

A uafp morphism $f:S \to T$ in $\on{PreStk}_{\on{w.ren}}$ is \emph{proper} (resp. a \emph{closed immersion}) if the uafp map $S\times_T T' \to T'$ as above is proper (resp. a closed immersion)\footnote{Recall that a uafp map in $\on{IndSch}_{\on{apai}}$ is schematic. Hence the notion of proper maps (resp. closed immersions) is well-defined.}.
\end{definition}

\begin{remark} Let $f:S \to T$ be a uafp morphism in $\on{PreStk}_{\on{w.ren}}$, and $g:T'\to T$ be a relatively apaisant morphism. Then by Remark \ref{rem-uafp-base-change-indsch}, $S\times_T T' \to T'$ is uafp.
\end{remark}

\begin{warning} The composition of two uafp morphisms might \emph{not} be uafp.
\end{warning}

\begin{construction} \label{const-!-pull-w-ren-uafp}
Let $f: S\to T$ be a uafp morphism in $\on{PreStk}_{\on{w.ren}}$. For any flat apaisant cover $T_0\in \on{Cov}(T)$, consider the flat apaisant cover $S_0:=S\times_T T_0 \to S$ of $S$. Let $T_\bullet$ and $S_\bullet$ be the corresponding Cech nerves. Note that $S_\bullet \to T_\bullet$ is a uafp morphism in $\on{IndSch}_{\on{apai}}$. By Construction \ref{constr-base-change-!*-apai-indsch}, the $!$-pullback functors induce a functor
$$ \on{Tot}_{\on{semi}} \IndCoh_*(T_\bullet) \to \on{Tot}_{\on{semi}}  \IndCoh_*(S\times_T T_\bullet).$$
Moreover, the obtained functor is functorial in $T_0$, i.e., we have a functor $\on{Cov}(T)^\op \to \on{Arrow}( \on{DGCat}_{\on{cont}} ) $ sending $T_0\in \on{Cov}(T)$ to the above arrow in $\on{DGCat}_{\on{cont}}$. By Theorem \ref{theorem-flat-descent-indcoh-w-ren} and Remark \ref{rem-cov-contractible-w.ren}, the functor $\on{Cov}(T)^\op \to \on{Arrow}( \on{DGCat}_{\on{cont}} ) $ is locally constant and the source category is weakly contractible. Hence it determines a \emph{canonical} functor
$$ f^!: \IndCoh_*(T) \to \IndCoh_*(S).$$

It follows from contruction that:
\begin{itemize}
	\item If $f$ is proper, we have a canonical adjoint pair $(f_*, f^!)$.

	\item If $f$ is an open immersion, we have a canonical adjoint pair $(f^!,f_*)$.

	\item The bi-simplicial symmetric monoidal  functor 
	$$ \IndCoh_*(-): \on{Cart}^{m,n}_{\on{all,uafp}}( \on{IndSch}_{\on{apai}} ) \to \on{Funct}( [m]^\op\times [n], \on{DGCat}_{\on{cont}})$$
	in Remark \ref{rem-higher-comp-!-pull-indsch} can be extended to a right-lax symmetric monoidal functor
	$$ \IndCoh_*(-): \on{Cart}^{m,n}_{\on{apai.indsch,uafp}}( \on{PreStk}_{\on{w.ren}} ) \to \on{Funct}( [m]^\op\times [n], \on{DGCat}_{\on{cont}})$$
	satisfying similar properties.
\end{itemize}

\end{construction}

Now we restrict to those $!$-pullback functors that send bounded-below objects to bounded-below ones.

\begin{lemma} \label{lem-left-t-exact-!-pull-w.ren-sharper}
Let $f: S\to T$ be an afp map in $\on{Sch}_{\on{apai}}$. Then there exists a natural number $n$ such that for any Cartesian square
$$
\xymatrix{
	S' \ar[r]^-{\varphi} \ar[d]^-\phi & T'\ar[d]^-g \\
	S \ar[r]^-f & T
}
$$
contained in $\on{PreStk}_{\on{w.ren}}$ with $\varphi$ being uafp, the functor $[-n]\circ \varphi^!$ is left t-exact. In particular, we obtain a functor
$$ (\varphi^!)^+: \IndCoh_*(T')^+  \to \IndCoh_*(S)^+.$$
\end{lemma}

\proof Follows from Lemma \ref{lem-left-t-exact-!-pull-apai-indsch-sharper}.

\qed

\begin{definition} A uafp morphism $\varphi:S'\to T'$ in $\on{PreStk}_{w.ren}$ is \emph{very afp}, or \emph{vafp}, if it is a base-change of an afp map $f: S\to T$ in $\on{Sch}_{\on{apai}}$.
\end{definition}

\begin{lemma} \label{lem-!-pull-base-change-w-ren-ult}
Let 
$$
\xymatrix{
	S' \ar[r]^-{\varphi} \ar[d]^-\phi & T'\ar[d]^-g \\
	S \ar[r]^-f & T
}
$$
be a Cartesian diagram in $\on{PreStk}_{\on{w.ren}}$ such that $f$ and $\varphi$ are vafp. For any flat apaisant cover $ q_0: T_0'\in T'$, let $p_0:S_0'\to S'$ be the flat apaisant cover of $S'$ induced by it. Consider the corresponding Cech nerves and Cartesian squares
$$
\xymatrix{
	S'_\bullet \ar[d]^-{p_\bullet} \ar[r]^-{\varphi_\bullet} & T'_\bullet \ar[d]^-{q_\bullet} \\
	S' \ar[r]^-{\varphi} \ar[d]^-\phi & T'\ar[d]^-g \\
	S \ar[r]^-f & T
}
$$
Then the following natural tranformation 
\begin{eqnarray*}
 (f^!)^+\circ (g_*)^+ & \simeq & (f^!)^+ \circ \on{Tot}_{\on{semi}} [   (( g\circ q_\bullet )_*)^+\circ  (q_\bullet^*)^+  ] \\
 & \to & \on{Tot}_{\on{semi}} [ (f^!)^+ \circ (( g\circ q_\bullet )_*)^+\circ  (q_\bullet^*)^+  ]  \\
  & \simeq & \on{Tot}_{\on{semi}} [  (( \phi\circ p_\bullet )_*)^+\circ (\varphi_\bullet^!)^+ \circ (q_\bullet^*)^+  ]   \\
  & \simeq & \on{Tot}_{\on{semi}} [  (( \phi\circ p_\bullet )_*)^+\circ (p_\bullet^*)^+ \circ (\varphi^!)^+   ]   \\
  & \simeq & (\phi_*)^+ \circ (\varphi^!)^+
\end{eqnarray*}
is an isomoprhism.

Moreover, if $g$ is $*$-pullable (hence so is $\phi$), the Bech-Chevalley natural transformation
\[
(\phi^{*})^+ \circ (f^!)^+ \to (\varphi^!)^+ \circ (g^{*})^+
\]
obtained from the above isomorphism is also an isomorphism.
\end{lemma}

\proof The first claim can be proved by using the same trick as in the proof of Proposition \ref{prop-general-push-w-ren}. To prove the second claim, we can use Theorem \ref{theorem-flat-descent-indcoh-w-ren} to reduce to the case when $g$ is relatively apaisant, which has already been treated in Construction \ref{const-!-pull-w-ren-uafp}.

\qed

\begin{construction} \label{constr-base-change-w-ren-ult}
By Lemma \ref{lem-!-pull-base-change-w-ren-ult} and Remark \ref{rem-higher-comp-!-pull-indsch}, we have a bi-simplicial right-lax symmetric monoidal functor 
	$$ \IndCoh_*(-)^+: \on{Cart}^{m,n}_{\on{all,vafp}}( \on{PreStk}_{\on{w.ren}} ) \to \on{Funct}( [m]^\op\times [n], \on{Vect}^+\on{-mod})$$
that encodes all the previous constructions. Namely, for a diagram 
	$$
\xymatrix{
	S_{0,0} \ar[d] \ar[r] & S_{1,0} \ar[r] \ar[d] & \cdots \ar[r] \ar[d] & S_{m,0} \ar[d] \\
	S_{0,1} \ar[d] \ar[r] & S_{1,1} \ar[r] \ar[d] & \cdots \ar[r] \ar[d] & S_{m,1} \ar[d] \\
 	\cdots	\ar[d] \ar[r] & \cdots \ar[r] \ar[d] & \cdots \ar[r] \ar[d] & \cdots \ar[d] \\
 	S_{0,n}  \ar[r] & S_{1,n} \ar[r]  & \cdots \ar[r]  & S_{m,n} 
}
$$
where all the squares are Cartesian and all the horizontal maps are vafp, it is sent to $$
\xymatrix{
	\IndCoh_*(S_{0,0})^+ \ar[d] & \IndCoh_*(S_{1,0})^+ \ar[l] \ar[d] & \cdots \ar[l] \ar[d] \ar[l] & \IndCoh_*(S_{m,0})^+ \ar[d] \ar[l] \\
	\IndCoh_*(S_{0,1})^+ \ar[d] & \IndCoh_*(S_{1,1})^+ \ar[l] \ar[d] & \cdots \ar[l] \ar[d] \ar[l] & \IndCoh_*(S_{m,1})^+ \ar[d] \ar[l] \\
 	\cdots \ar[d] & \cdots  \ar[l] \ar[d] & \cdots \ar[l] \ar[d] \ar[l] & \cdots  \ar[d] \ar[l] \\
 	\IndCoh_*(S_{0,n})^+ & \IndCoh_*(S_{1,n})^+ \ar[l]  & \cdots \ar[l]\ar[l] & \IndCoh_*(S_{m,n})^+ \ar[l] 
}
$$
where all horizontal functors are $!$-pullback functors and all vertical functors are $*$-pushforward functors. Moreover, if each $S_{i,j} \to S_{i,j'}$ is $*$-pullable, then all the squares in the above diagram are left adjointable along the vertical direction.
\end{construction}

\subsection{Renormalization}

\begin{construction}
For $S \in \on{PreStk}_{\on{w.ren}}$, let $\Coh(S) \subset \IndCoh_{*}(S)$ be the full subcategory of objects $\mathcal{F}$ such that its $*$-pullback along any flat apaisant cover $S'\to S$ is contained in $\Coh(S') \subset \IndCoh_{*}(S')$.
\end{construction}

\begin{remark} \label{rem-coh-bounded-w-ren}
By Lemma \ref{lem-desc-coh-indsch}, $\on{Coh}(S) \subset \IndCoh_{*}(S)^b$.
\end{remark}

\begin{lemma} \label{lem-coh-on-w-ren}
For $S \in \on{PreStk}_{\on{w.ren}}$ and $\mathcal{F}\in \IndCoh_{*}(S)$, the following are equivalent:
\begin{itemize}
	\item[(i)] $\mathcal{F} \in \on{Coh}(S) $;
	\item[(ii)] there exists a flat apaisant cover $S'\to S$ such that the $*$-pullback of $\mathcal{F}$ is contained in $\Coh(S') \subset \IndCoh_{*}(S')$;
	\item[(iii)] for any $*$-pullable morphism $T\to S$ in $ \on{PreStk}_{\on{w.ren}}$, the $*$-pullback of $\mathcal{F}$ is contained in $\Coh(T) \subset \IndCoh_{*}(T)$.
\end{itemize}
\end{lemma}

\proof Follows from Lemma \ref{lem-pullback-being-coh-implies-coh}.

\qed

\begin{lemma} \label{lem-coh-closed-under-truncation-w.ren}
For $S \in \on{PreStk}_{\on{w.ren}}$, $\on{Coh}(S)$ is contained in $\IndCoh_{*}(S)^+$ and is closed under truncations.
\end{lemma}

\proof Follows from Corollary \ref{cor-coh-closed-truncation-apai-indsch}.

\qed

\begin{construction} For $S \in \on{PreStk}_{\on{w.ren}}$, define
$$ \IndCoh_{*,\ren}(S):= \tx{Ind}( \tx{Coh}(S) ).$$
There is an obvious functor ${\ren}: \IndCoh_{*,\ren}(S)\to \tx{IndCoh}^*(S)$ defined by ind-extending.
\end{construction}

\begin{remark} When $S$ is an apaisant indscheme, recall that $\IndCoh_{*}(S)$ is compactly generated by $\on{Coh}(S)$. Hence ${\ren}: \IndCoh_{*,\ren}(S)\to \tx{IndCoh}^*(S)$ is an equivalence. Therefore we omit the subscript \emph{ren} in this case.
\end{remark}

\begin{remark} \label{rem-qca-indcoh*-ren}
When $S$ is a QCA stack, via the identification in Remark \ref{rem-qca-indcoh*}, the subcategory $\on{Coh}(S)$ coincides with the same-named category in \cite{drinfeld2013some}. Hence by \cite[Theorem 3.3.5]{drinfeld2013some}, the category $\IndCoh_{*,\ren}(S)$ coincides with $\IndCoh(S)$. Hence we write $\IndCoh(S)$ for $\IndCoh_{*,\ren}(S)$ in this case.
\end{remark}

\begin{remark} It is easy to see $\Coh(S) \subset \IndCoh_{*}(S)$ is Karoubi complete and closed under finite colimits, hence $\IndCoh(S)^c \simeq \Coh(S)$ and $\IndCoh_{*,\ren}(S)$ is a cocomplete DG-category.
\end{remark}

\begin{construction}
There exists a compactly generated t-structure on $\IndCoh_{*,\ren}(S)$ such that $\IndCoh_{*,\ren}(S)^{\le 0}$ is generated under colimits and extensions by $\on{Coh(S)}\cap \IndCoh_{*}(S)^{\le 0}$.
\end{construction}

\begin{proposition} (C.f. \cite[Proposition 6.28.5]{raskin2020homological}) For $S \in \on{PreStk}_{\on{w.ren}}$, the functor $\on{\ren}: \IndCoh_{*,\ren}(S) \to \IndCoh_{*}(S)$ is t-exact and induces an equivalence $\IndCoh_{*,\ren}(S)^+ \simeq \IndCoh_{*}(S)^+ $.

\end{proposition}

\proof If $S$ is \emph{renormalizable}, i.e., it has a flat apaisant cover $S'\to S$ that is laft, the proof in \cite{raskin2020homological} is correct (after adding the word ``apaisant'' at suitable places). However, the proof there only uses the \emph{laft} assumption on $S'$ to guarantee $\on{Coh}(S)\subset \IndCoh_*(S)$ is closed under truncations, which we have already proved in Lemma \ref{lem-coh-closed-under-truncation-w.ren}.

\qed

\begin{construction} For any morphism $f:S\to T$ in $\on{PreStk}_{{\on{w.ren}}}$, the left t-exact functor $(f_*)^+:\IndCoh_{*}(S)^+ \to \IndCoh_{*}(T)^+$ induces the following functor:
$$  \Coh(S) \to \IndCoh_{*}(S)^+ \to \IndCoh_{*}(T)^+ \simeq \IndCoh_{*,\ren}(T)^+ \to \IndCoh_{*,\ren}(T) $$
By ind-extending, we obtain the \emph{renormalized $*$-pushforward} functor
$$ f_{*,\ren}: \IndCoh_{*,\ren}(S) \to  \IndCoh_{*,\ren}(T).$$
Since $\Coh(S)\subset \IndCoh_{*,\ren}(S)$ is closed under truncations, the functor $f_{*,\ren}$ is left t-exact. Since the retriction of $f_*^+$ on $\IndCoh_{*}(S)^{\ge 0}$ preserves filtered colimits, we have a commutative diagram
$$
\xymatrix{
	\IndCoh_{*,\ren}(S)^+ \ar[rr]^-{ (f_{*,\ren})^+ }  \ar[d]^-{\simeq}
	& & \IndCoh_{*,\ren}(T)^+ \ar[d]^-{\simeq} \\
	\IndCoh_{*}(S)^+ \ar[rr]^-{ (f_*)^+ } 
	& & \IndCoh_{*}(T)^+.
}
$$
Also, if $f$ is apaisant indschematic, then the above commutative diagram can be extended to the unbounded categories.
\end{construction} 

\begin{construction} Similarly, we can use the above \emph{renormalization} procedure to construct the \emph{renormalized $*$-pullback functors} (along any $*$-pullback morphisms) and the \emph{renormalized $!$-pullback functors} (along any vafp morphisms). By Construction \ref{constr-base-change-w-ren-ult}, these functors can be encoded a bi-simplicial right-lax symmetric monoidal functor 
	$$ \IndCoh_{*,\ren}(-): \on{Cart}^{m,n}_{\on{all,vafp}}( \on{PreStk}_{\on{w.ren}} ) \to \on{Funct}( [m]^\op\times [n], \on{DGCat}_{\on{cont}})$$
satisfying similar properties as in Construction \ref{constr-base-change-w-ren-ult}.
\end{construction}

\begin{warning} The $\IndCoh_{*,\ren}$-theory does \emph{not} satisfy flat descent. For example, the $*$-pullback functor $\IndCoh_{*,\ren}(G(O)) \to \IndCoh_{*,\ren}(\on{pt})$ is not conservative. 
\end{warning}

\begin{proposition} \label{prop-indcoh-ren-zariski}
Let $S_1 \to S$ and $S_2\to S$ be two schematic open embeddings in $\on{PreStk}_{\on{w.ren}}$ such that $S_1\sqcup S_2 \to S$ is surjective. Then $*$-pullbacks induce an equivalence
$$
 \IndCoh_{*,\ren}(S) \to \IndCoh_{*,\ren}(S_1) \times_{ \IndCoh_{*,\ren}(S_1\cap S_2) } \IndCoh_{*,\ren}(S_2).
 $$
\end{proposition}

\proof The proof is similar to that of Proposition \ref{prop-indcoh-zariski-sch}. Using the notations of that proof, we only need to check
\begin{itemize}
	\item[(i)] For any $\mathcal{F} \in  \IndCoh_{*,\ren}(S)$, the morphism
	$$ \mathcal{F} \to (u_1)_{*,\ren}\circ (u_1)^*_\ren(\mathcal{F}) \times_{  v_{*,\ren}\circ (v)^*_\ren(\mathcal{F})  } (u_2)_{*,\ren}\circ (u_2)^*_\ren(\mathcal{F}) $$
	is an equivalence.

	\item[(ii)] For any $(  \mathcal{F}_1, \mathcal{F}_2, \mathcal{F}_{12}, (v_i)^*_\ren(\mathcal{F}_i) \simeq  \mathcal{F}_{12})$ as before, the morphism
	$$ (u_i)^{*}_\ren[(u_1)_{*,\ren}(\mathcal{F}_1) \times_{  v_{*,\ren}(\mathcal{F}_{12})  } (u_2)_{*,\ren}(\mathcal{F}_2) ] \to \mathcal{F}_i $$
	is an equivalence.
\end{itemize}
As before, (ii) follows from the base-change isomorphisms. To prove (i), we can reduce to the case when $\mathcal{F}\in \Coh(S)$. Then the claim follows from Proposition \ref{prop-indcoh-zariski-w-ren}.

\qed

\begin{proposition} \label{prop-prod-formula-indcoh-w-ren}
For any $S,T\in \on{PreStk}_{\on{w.ren}}$, the external product functor
$$-\boxtimes - : \IndCoh_{*,\ren}(S) \otimes \IndCoh_{*,\ren}(T) \to \IndCoh_{*,\ren}(S\times T) $$
is fully faithful and preserves compact objects. If $S$ and $T$ are apaisant indschemes, it is an equivalence.
\end{proposition}

\proof When $S$ and $T$ are apaisant indschemes, there is nothing new to prove (see the last paragraph of Remark \ref{rem-higher-comp-!-pull-indsch}).

In the general case, we first show the external product functor sends compact objects to compact objects. Recall $\IndCoh_{*,\ren}(S) \otimes \IndCoh_{*,\ren}(T)$ is compactly generated by the image of
\begin{equation} \label{eqn-prop-prod-formula-indcoh-w-ren-1}
 \on{Coh}(S) \times \on{Coh}(T) \to \IndCoh_{*,\ren}(S) \times \IndCoh_{*,\ren}(T) \to \IndCoh_{*,\ren}(S) \otimes \IndCoh_{*,\ren}(T).
 \end{equation}
Hence we only need to show
$$ \on{Coh}(S) \times \on{Coh}(T) \to \IndCoh_{*}(S) \times \IndCoh_{*}(T)  \to  \IndCoh_{*}(S\times T)   $$
factors through $\on{Coh}(S\times T)$. Using Lemma \ref{lem-coh-on-w-ren}, we can reduce to the case when $S$ and $T$ are apaisant indschemes. Then the claim is obvious because 
$$ \IndCoh_{*}(S) \otimes \IndCoh_{*}(T)  \to  \IndCoh_{*}(S\times T)  $$
is an equivalence hence sends compact objects to compact objects.

We then show the external product functor is fully faithful. Let $\mathcal{F}^i\in \on{Coh}(S)$ and $\mathcal{G}^i\in \on{Coh}(T)$. Let $\mathcal{M}^i$ be the image of $(\mathcal{F}^i, \mathcal{G}^i)$ under (\ref{eqn-prop-prod-formula-indcoh-w-ren-1}). Recall we have
$$ \underline{\on{Hom}}( \mathcal{M}^1,\mathcal{M}^2 ) \simeq \underline{\on{Hom}}( \mathcal{F}^1,\mathcal{F}^2 ) \otimes  \underline{\on{Hom}}( \mathcal{G}^1,\mathcal{G}^2 ),$$
where $\underline{\on{Hom}}(-,-)\in \on{Vect}$ is the \emph{mapping object}. Hence we only need to show 
$$\underline{\on{Hom}}( \mathcal{F}^1,\mathcal{F}^2 ) \otimes  \underline{\on{Hom}}( \mathcal{G}^1,\mathcal{G}^2 ) \to \underline{\on{Hom}}( \mathcal{F}^1\boxtimes \mathcal{G}^1,\mathcal{F}^2\boxtimes \mathcal{G}^2 ) $$
is an equivalence. Let $S_0\to S$ and $T_0\to T$ be flat apaisant covers. Consider the Cech nerves $S_\bullet\to S$ and $T_\bullet \to T$. Let $\mathcal{F}^i_\bullet \in \on{Coh}(S_\bullet)$ be the $*$-pullback of $\mathcal{F}^i$, and similarly define $\mathcal{G}^i_\bullet$. By Theorem \ref{theorem-flat-descent-indcoh-w-ren}, the above morphism is equivalent to the following composition:
\begin{eqnarray*}
\underline{\on{Hom}}( \mathcal{F}^1,\mathcal{F}^2 ) \otimes  \underline{\on{Hom}}( \mathcal{G}^1,\mathcal{G}^2 ) & \simeq & \underset{m}{ \on{Tot}_{\on{semi}}} [ \underline{\on{Hom}}( \mathcal{F}_m^1,\mathcal{F}_m^2 ) ] \otimes \underset{n}{ \on{Tot}_{\on{semi}}} [ \underline{\on{Hom}}( \mathcal{G}_n^1,\mathcal{G}_n^2 ) ]   \\
& \to &  \underset{m}{ \on{Tot}_{\on{semi}}} \underset{n}{ \on{Tot}_{\on{semi}}} [ \underline{\on{Hom}}( \mathcal{F}_m^1,\mathcal{F}_m^2 ) \otimes   \underline{\on{Hom}}( \mathcal{G}_n^1,\mathcal{G}_n^2 ) ] \\
& \simeq &  \on{Tot}_{\on{semi}} [ \underline{\on{Hom}}( \mathcal{F}_\bullet^1,\mathcal{F}_\bullet^2 ) \otimes   \underline{\on{Hom}}( \mathcal{G}_\bullet^1,\mathcal{G}_\bullet^2 ) ]\\
& \to &  \on{Tot}_{\on{semi}} [ \underline{\on{Hom}}( \mathcal{F}_\bullet^1 \boxtimes \mathcal{G}_\bullet^1,\mathcal{F}_\bullet^2 \boxtimes \mathcal{G}_\bullet^2) ] \\
& \simeq & \underline{\on{Hom}}( \mathcal{F}^1\boxtimes \mathcal{G}^1,\mathcal{F}^2\boxtimes \mathcal{G}^2 ).
 \end{eqnarray*}
Note that each 
$$ \underline{\on{Hom}}( \mathcal{F}_\bullet^1,\mathcal{F}_\bullet^2 ) \otimes   \underline{\on{Hom}}( \mathcal{G}_\bullet^1,\mathcal{G}_\bullet^2 ) \to \underline{\on{Hom}}( \mathcal{F}_\bullet^1 \boxtimes \mathcal{G}_\bullet^1,\mathcal{F}_\bullet^2 \boxtimes \mathcal{G}_\bullet^2) $$
is an isomorphism because $S_\bullet$ and $T_\bullet$ are apaisant indschemes. Hence we only need to show
$$  \underset{m}{ \on{Tot}_{\on{semi}}} [ \underline{\on{Hom}}( \mathcal{F}_m^1,\mathcal{F}_m^2 ) ] \otimes \underset{n}{ \on{Tot}_{\on{semi}}} [ \underline{\on{Hom}}( \mathcal{G}_n^1,\mathcal{G}_n^2 ) ]  \to  \underset{m}{ \on{Tot}_{\on{semi}}} \underset{n}{ \on{Tot}_{\on{semi}}} [ \underline{\on{Hom}}( \mathcal{F}_m^1,\mathcal{F}_m^2 ) \otimes   \underline{\on{Hom}}( \mathcal{G}_n^1,\mathcal{G}_n^2 ) ] $$
is an isomorphism. By Remark \ref{rem-coh-bounded-w-ren}, we can assume $F^i, G^i$ are concentrated in degrees $[0,d]$ for some integer $d$, and so are $F^i_\bullet, G^i_\bullet$. Hence each $\on{Hom}(-,-)$ as above is contained in $\on{Vect}^{\ge -d}$. Therefore we only need to show the above morphism induces isomorphisms on $H^j$ for $j\ge -2d$. For each fixed $j$, we can replace the semi-totalization by a truncated one. Then we are done because tensor products in $\on{Vect}$ commutes with finite limits.

\qed

So far we only allow $!$-pullback functors along vafp morphisms, which does not include \emph{ind}-proper morphisms. These functors can be supplied easily as follows.

\begin{definition} A morphism $f:S\to T$ be a morphism in $\on{PreStk}_{\on{w.ren}}$ is \emph{ind-afp and ind-proper} (resp. \emph{an ind-afp closed immersion}) if there exists a flat apaisant cover $T_0\in \on{Cov}(T)$ such that $S\times_T T_0 \to T_0$ is an ind-afp and ind-proper map (resp. an ind-afp closed immersion) in $\on{IndSch}_{\on{apai}}$. 

An ind-afp and ind-proper morphism $f:S\to T$ is \emph{$(?,\ren)$-pullable} if $f_*$ is left t-exact up to a shift (hence so is $f_{*,\ren}$). 
\end{definition}

\begin{remark} We do \emph{not} require ind-afp and ind-proper map to be apaisant indschematic.
\end{remark}

\begin{lemma} 
\label{lem-base-change-ren-weak}
\label{lem-base-change-ren}
\begin{itemize}
	\item[(1)] Let $f:S\to T$ be an ind-afp and ind-proper morphism in $\on{PreStk}_{\on{w.ren}}$, then the functor $f_{*,\ren}$ has a \emph{continuous} right adjoints.

	\item[(2)] Let $f:S\to T$ be a $(?,\ren)$-pullable map. Then the right adjoints $f^?$ and $f^?_\ren$ of $f_*$ and $f_{*,\ren}$ are left t-exact up to shifts, and $f^?_\ren$ can be obtained from $f^?$ via renormalization. I.e., it can be identified with the ind-extension of
	$$ \on{Coh}(T) \to \IndCoh_{*}(T)^+ \xrightarrow{f^?}  \IndCoh_{*}(S)^+ \simeq \IndCoh_{*,\ren}(S)^+  \to \IndCoh_{*}(S)$$
	such that we have a canonical commutative diagram
	$$
	\xymatrix{
		\IndCoh_{*,\ren}(T) \ar[r]^-{  f_\ren^? } \ar[d]^-{\ren} & 
		\IndCoh_{*,\ren}(S) \ar[d]^-{\ren} \\ 
		\IndCoh_{*}(T) \ar[r]^-{ f^? } & 
		\IndCoh_{*}(S).
	}
	$$

\item[(3)] Let 
$$
\xymatrix{
	S' \ar[r]^-{\varphi} \ar[d]^-\phi & T'\ar[d]^-g \\
	S \ar[r]^-f & T
}
$$
be a Cartesian diagram in $\on{PreStk}_{\on{w.ren}}$ such that $f^?_\ren$ and $\varphi^?_\ren$ are well-defined. Then the Beck-Chevalley natural transformation
\[
\phi_{*,\ren} \circ \varphi^?_\ren \to f^?_\ren \circ g_{*,\ren}.
\]
is an isomorphism. If $g$ is $*$-pullable (hence so is $\phi$), then the Beck-Chevalley natural transformation
\[
\phi^*_{\ren} \circ f^?_\ren \to
 \varphi^?_\ren   \circ g^*_{\ren}.
\]
is also an isomorphism. If $g$ and $\phi$ are vafp, then the Beck-Chevalley natural transformation
\[
 \phi^{!}_\ren \circ  f^?_\ren \to \varphi^?_\ren\circ g^{!}_\ren 
\]
is also an equivalence.
\end{itemize}
\end{lemma}

\proof To prove (1), we only need to show $(f_*)^+$ sends $\on{Coh}(S)$ to $\on{Coh}(T)$, which follows from Lemma \ref{lem-pullback-being-coh-implies-coh}.

Note that $S_0:S\times_T T_0$ is a flat apaisant cover of $S$. Let $T_\bullet \to T$ and $S_\bullet \to S$ be the Cech nerves. By Theorem \ref{theorem-flat-descent-indcoh-w-ren} and the base-change isomorphisms, $f_*$ can be identified with 
$$ \on{Tot}_{\on{semi}} \IndCoh_{*}(S_\bullet) \xrightarrow{*\on{-push}} \on{Tot}_{\on{semi}} \IndCoh_{*}(T_\bullet).$$
By Lemma \ref{lem-base-change-indsch}, this functor has a continuous right adjoint
$$\on{Tot}_{\on{semi}} \IndCoh_{*}(T_\bullet) \xrightarrow{!\on{-pull}} \on{Tot}_{\on{semi}} \IndCoh_{*}(S_\bullet) . $$
This proves the existence of $f^?$. The remaining claims of (2) follows from definitions.

(3) can be deduced from the non-renormalized version.

\qed

\begin{remark} Using Proposition \ref{prop-indcoh-on-ind-stack} below, Lemma \ref{lem-base-change-ren}(3) can be generalized to the case when $S$ can be written as a nice enough filtered colimit $S\simeq \colim S_\alpha$ such that $S_\alpha \to T$ and $S_\alpha\times_S T' \to T'$ are $(?,\ren)$-pullable. For example, this means for the map $f: \mathcal{L}^+ G\backslash \mathcal{L}G / \mathcal{L}^+G \to \on{pt} / \mathcal{L}^+G$, the right adjoint of $f_{*,\ren}$ is well-behaved although it can \emph{not} be obtained via a renormalization procedure. However, for the purpose of this paper, Lemma \ref{lem-base-change-ren}(3) is enough (because we use the \emph{forma} stack $\mathcal{L}^+ G\backslash \widehat{\mathcal{L}G} / \mathcal{L}^+G $ instead of the \emph{ind}-stack $\mathcal{L}^+ G\backslash \mathcal{L}G / \mathcal{L}^+G$).
\end{remark}

\begin{variant} For $S \in \on{PreStk}_{\on{w.ren}}$, the category $\IndCoh_{*,\ren}(S)$ is compactly generated hence dualizable. We define
$$ \IndCoh^{!,\ren}(S):= \IndCoh_{*,\ren}(S)^{\vee} \simeq \on{Ind}( \on{Coh}(S)^\op ).$$
For $f:S \to T$, we define the following functors:
\begin{itemize}
	\item The \emph{renormalized $!$-pullback functor for $\IndCoh^!$-theory}
	$$ f^{!,\ren}: \IndCoh^{!,\ren}(T) \to \IndCoh^{!,\ren}(S) $$
	is defined as the dual of
	$$ f_{*,\ren}: \IndCoh_{*,\ren}(S) \to \IndCoh_{*,\ren}(T). $$ 

	\item When $f$ is $*$-pullable, the \emph{renormalized $?$-pushforward functor for $\IndCoh^!$-theory}
	$$ f_?^{\ren}: \IndCoh^{!,\ren}(S) \to \IndCoh^{!,\ren}(T),$$
	is defined as the dual of
	$$ f_{\ren}^*: \IndCoh_{*,\ren}(T) \to \IndCoh_{*,\ren}(S), $$
	or equivalently as the \emph{continuous} right adjoint of $f^{!,\ren}$.

	We use \emph{$?$-pushable} as a synonym for \emph{$*$-pullable}.

	\item When $f$ is vafp, the \emph{renormalized $*$-pushforward functor for $\IndCoh^!$-theory}
	$$ f_*^{\ren}: \IndCoh^{!,\ren}(S) \to \IndCoh^{!,\ren}(T) $$
	is defined as the dual of
	$$ f_{\ren}^!: \IndCoh_{*,\ren}(T) \to \IndCoh_{*,\ren}(S). $$ 

	\item When $f$ is $(?,\ren)$-pullable, the \emph{renormalized $!$-pushforward functor for $\IndCoh^!$-theory}
	$$ f_!^\ren: \IndCoh^{!,\ren}(S) \to \IndCoh^{!,\ren}(T), $$
	is defined as the dual of 
	$$ f_{\ren}^?: \IndCoh_{*,\ren}(T) \to \IndCoh_{*,\ren}(S),$$
	or equivalently as the left adjoint of $f^{!,\ren}$.

	We use \emph{$(!,\ren)$-pushable} as a synonym for $(?,\ren)$-pullable.
\end{itemize}
For $S_1, S_2$ such that the external tensor product functor
$$ \IndCoh_{*,\ren}(S_1) \otimes \IndCoh_{*,\ren}(S_2) \to \IndCoh_{*,\ren}(S_1\times S_2)$$
is an equivalence, we also have an equivalence
$$ \IndCoh^{!,\ren}(S_1) \otimes \IndCoh^{!,\ren}(S_2) \to \IndCoh^{!,\ren}(S_1\times S_2)$$
obtained by taking conjugation.
\end{variant}

\begin{notation} For $S \in \on{PreStk}_{\on{w.ren}}$, let
$$\omega_S\in \IndCoh^{!,\ren}(S)$$
be the image of the unit object in $\IndCoh(\on{pt})\simeq \on{Vect}$ under the $!$-pullback functor. We often omit the subscript $S$ from $\omega_S$ if there is no danger of ambiguity.
\end{notation}

\subsection{Ind-coherent sheaves on quotient stacks}
In this subsection, we compare Raskin's definition of $\IndCoh_{*,\ren}(S)$ with the definition sketched in \cite{jerusalemday43} for quotient prestacks such as $\Arc G \backslash \Loop G / \Arc G$.

\begin{proposition} \label{prop-indcoh-on-ind-stack}
Suppose $S \in \on{PreStk}_{\on{w.ren}}$ can be written as a filtered colimit $S \simeq \colim_{\alpha\in I} S^\alpha$ such that: 
\begin{itemize}
	\item  The connecting maps $\iota_{\alpha\to \beta}:S^\alpha \to S^\beta$ and $\iota_\alpha:S^\alpha \to S$ are apaisant indschematic;
	\item There exists a flat apaisant cover $S_0\to S$ such that the map $S^\alpha\times_S S_0 \to S^\beta \times_S S_0$ is an ind-afp closed immersion in $\on{IndSch}_{\on{apai}}$.
\end{itemize}
Then the functors
\begin{eqnarray*}
 \colim_{*\on{-push}}  \IndCoh_{*}(S^\alpha) \to   \IndCoh_{*}(S)\\
 \colim_{*\on{-push}}  \IndCoh_{*,\ren}(S^\alpha) \to   \IndCoh_{*,\ren}(S)
 \end{eqnarray*}
are equivalences.
\end{proposition}

\proof  Let $S_0^\alpha:=S^\alpha\times_S S_0$. Note that $S_0^\alpha \to S^\alpha$ is also an apaisant cover. Let $S_\bullet\to S$ be the Cech nerve for $S_0\to S$ and $S_\bullet^\alpha \to S^\alpha$ be defined similarly. Note that we also have $S_\bullet \simeq \colim_\alpha S_\bullet^\alpha$. By assumption, each $S^\alpha_0 \to S_0$ is ind-afp and ind-proper. Hence $S^\alpha_0  \to S^\beta_0 $ and the base-changes $S^\alpha_\bullet \to S_\bullet$, $S^\alpha_\bullet  \to S^\beta_\bullet $ are also ind-afp and ind-proper (Lemma \ref{lem-fiber-product-apai-indsch-is-apai-if-afp}).

We first prove the non-renormalized version. Using the base-change isomorphisms and Theorem \ref{theorem-flat-descent-indcoh-w-ren}, we have the following commutative diagram
$$
\xymatrix{
	\underset{*\on{-push}}{\colim}  \IndCoh_{*}(S^\alpha) \ar[r]^-\simeq 
	\ar[d] &
	\underset{*\on{-push}}{\colim}  \underset{*\on{-pull}} {\on{Tot}_{\on{semi}}}  \IndCoh_{*}(S^\alpha_\bullet) \ar[r] &
	\underset{*\on{-pull}} {\on{Tot}_{\on{semi}}} \underset{*\on{-push}}{\colim}   \IndCoh_{*}(S^\alpha_\bullet) \ar[d]  \\
	\IndCoh_{*}(S) \ar[rr]^-\simeq & &
		\underset{*\on{-pull}} {\on{Tot}_{\on{semi}}} \IndCoh_{*}(S_\bullet)
}
$$
Since $S_\bullet \simeq \colim_\alpha S_\bullet^\alpha$ in $\on{IndSch}_{\on{apai}}$, the right vertical functor is an equivalence. Hence it remains to show 
$$
\underset{*\on{-push}}{\colim}  \underset{*\on{-pull}} {\on{Tot}_{\on{semi}}}  \IndCoh_{*}(S^\alpha_\bullet) \to
	\underset{*\on{-pull}} {\on{Tot}_{\on{semi}}} \underset{*\on{-push}}{\colim}   \IndCoh_{*}(S^\alpha_\bullet) 
$$
is an equivalence. By Lemma \ref{lem-base-change-indsch}, this functor is equivalent to
$$
\underset{!\on{-pull}}{\lim}  \underset{*\on{-pull}} {\on{Tot}_{\on{semi}}}  \IndCoh_{*}(S^\alpha_\bullet) \to
	\underset{*\on{-pull}} {\on{Tot}_{\on{semi}}} \underset{!\on{-pull}}{\lim}   \IndCoh_{*}(S^\alpha_\bullet),
$$
which is an equivalence because taking limits commute with taking limits.

Now we prove the renormalized version. We first show the functor is fully faithful. By Lemma \ref{lem-base-change-ren-weak}(1), $\colim_{*\on{-push}}  \IndCoh_{*,\ren}(S^\alpha)$ is compactly generated by the images of
$$\on{ins}_{\alpha,\ren}: \on{Coh}(S^\alpha) \to \colim_{*\on{-push}}  \IndCoh_{*,\ren}(S^\alpha) $$
 and the functor
$$\colim_{*\on{-push}}  \IndCoh_{*,\ren}(S^\alpha) \to   \IndCoh_{*,\ren}(S)$$
preserves compact objects. Hence we only need to show for fixed $\alpha$ and $\mathcal{F}, \mathcal{G}\in \on{Coh}(S^{\alpha})$, the morphism
$$ \on{Maps}( \on{ins}_{\alpha,\ren}( \mathcal{F} ),  \on{ins}_{\alpha,\ren}( \mathcal{G} )  ) \to \on{Maps} ( (\iota_{\alpha})_{*,\ren}(\mathcal{F}), (\iota_{\alpha})_{*,\ren}(\mathcal{G}) ) $$
is an isomorphism. Note that we have a commutative diagram
$$
\xymatrix{
	\underset{*\on{-push}}{\colim}  \IndCoh_{*,\ren}(S^\alpha) \ar[d] \ar[r]^-{\ren} &  
	\underset{*\on{-push}}{\colim}  \IndCoh_{*}(S^\alpha) \ar[d]^-\simeq \\
	\IndCoh_{*,\ren}(S) \ar[r]^-{\ren} &
	 \IndCoh_{*}(S)
}
$$
which induces
$$
\xymatrix{
	\on{Maps}( \on{ins}_{\alpha,\ren}( \mathcal{F} ),  \on{ins}_{\alpha,\ren}( \mathcal{G} )  ) \ar[d] \ar[r] &  
	\on{Maps}( \on{ins}_{\alpha}\circ \ren ( \mathcal{F} ),  \on{ins}_{\alpha}\circ \ren( \mathcal{G} )  ) \ar[d]^-\simeq \\
	\on{Maps} ( (\iota_{\alpha})_{*,\ren}(\mathcal{F}), (\iota_{\alpha})_{*,\ren}(\mathcal{G}) ) \ar[r] &
	\on{Maps} ( (\iota_{\alpha})_{*}\circ  \ren(\mathcal{F}), (\iota_{\alpha})_{*} \circ \ren (\mathcal{G}) ).
}
$$
By definition, the bottom morphism is an isomorphism. Hence we only need to show the top morphism is an isomorphism. By Lemma \ref{lem-*-push-left-t-exact-indsch}, the functor $(\iota_{\alpha \to \beta})_{*}$ is t-exact. Hence by Lemma \ref{lem-base-change-ren-weak}, the t-exact connecting functors $(\iota_{\alpha \to \beta})_{*,\ren}$ and $(\iota_{\alpha \to \beta})_{*}$ have continuous right adjoints. Moreover  $(\iota_{\alpha \to \beta})_\ren^?$ and $(\iota_{\alpha \to \beta})^?$ are compatible with the functor $\ren$. This implies the top morphism is an equivalence because formally it can be written as
$$ \on{Maps}(  \mathcal{F} ,  \colim_{\beta\in I_{\alpha/}} (\iota_{\alpha \to \beta})_\ren^?\circ (\iota_{\alpha \to \beta})_{*,\ren} ( \mathcal{F} )  ) \to \on{Maps}(   \ren(\mathcal{F}) ,  \colim_{\beta\in I_{\alpha/}} (\iota_{\alpha \to \beta})^?\circ (\iota_{\alpha \to \beta})_{*} \circ \ren( \mathcal{F} )  ).   $$

This finishes the proof that
$$ \colim_{*\on{-push}}  \IndCoh_{*,\ren}(S^\alpha) \to   \IndCoh_{*,\ren}(S)$$
is fully faithful. It remains to show it is surjective. We only need to show any $\mathcal{M}\in \on{Coh}(S)$,
$$  \colim_\alpha   (\iota_{\alpha})_{*,\ren} \circ (\iota_{\alpha})_\ren^?(  \mathcal{M}) \to \mathcal{M}$$
is an equivalence. Note that both sides are bounded below, hence we only need to show it is sent to an equivalence by $\ren$. Using Lemma \ref{lem-*-push-left-t-exact-indsch} again, $\ren$ sends it to
$$  \colim_\alpha   (\iota_{\alpha})_{*} \circ (\iota_{\alpha})^?\circ \ren(  \mathcal{M}) \to \ren(\mathcal{M}),$$
which is indeed an equivalence because
$$\colim_{*\on{-push}}  \IndCoh_{*}(S^\alpha) \to   \IndCoh_{*}(S).$$

\qed

\begin{proposition} \label{prop-indcoh-on-quot-stack}
Let $B\in \,^{>-\infty}\on{Sch}_{\on{aft}}$ be a base scheme, $T = \lim_{n\in \mathbb{N}} T^n$ be an apaisant presentation for an apaisant scheme over $B$, and $H=\lim_{n\in \mathbb{N}} H^n$ be an apaisant presentation for an apaisant affine flat group scheme over $B$. Suppose we have compatible actions of $H^n$ on $T^n$ relative to $B$. Then:
\begin{itemize}
	\item[(0)] The fppf quotient stack $T^n/H^n$ is a QCA stack.

	\item[(1)] The fppf quotient stack $T/H$ is weakly renormalizable and $T\to T/H$ is a flat apaisant cover.

	\item[(2)] The maps $ T/H \to T^n/H^n$ and $T^n/H^n \to T^m/H^m$ are $*$-pullable. 

	\item[(3)] The functor
$$ \colim_{*\on{-pull}} \IndCoh(T^n/H^n) \to \IndCoh_{*,\ren}(T/H) $$
is an equivalence. In particular, the functor
$$ \IndCoh_{*,\ren}(T/H) \to \lim_{*\on{-push}} \IndCoh(T^n/H^n) $$
is also an equivalence. 
	\end{itemize}
\end{proposition}

\proof (0) is well-known because $T^n, H^n \in \,^{>-\infty}\on{Sch}_{\on{aft}}$ and $H^n$ is affine. 

To prove (1), we only need to prove $T\to T/H$ is relatively apaisant. But this follows from Example \ref{exam-classifying-stack-apaisant}.

To prove (2), we only need to check $T\to T^n$ and $T^n\to T_j$ are relatively apaisant. But both claims are obvious.

We now prove (3). We first prove 
$$ \colim_{*\on{-pull}} \IndCoh(T^n/H^n) \to \IndCoh_{*,\ren}(T/H) $$
is fully faithful. Note that it preserves compact objects. Hence we only need to check full-faithfulness on compact objects. Consider the Cech nerve $T_\bullet \to T/H$ for the cover $T\to T/H$ and the similar Cech nerves $T_\bullet^n \to T^n/H^n$. For fixed $n$ and $\mathcal{F}^n,\mathcal{G}^n\in \on{Coh}(T^n/H^n)$, let $\mathcal{F}^m,\mathcal{G}^m\in \on{Coh}(T^m/H^m)$ be the $*$-pullbacks (for $m\ge n$) and similarly define $\mathcal{F}, \mathcal{G}$. Also consider the $*$-pullbacks $\mathcal{F}_\bullet^m,\mathcal{G}_\bullet^m \in \on{Coh}(T^m_\bullet)$ and $\mathcal{F}_\bullet,\mathcal{G}_\bullet \in \on{Coh}(T_\bullet)$. We only need to show 
$$ \colim_{m\ge n} \on{Maps}( \mathcal{F}^m,\mathcal{G}^m  ) \to \on{Maps} ( \mathcal{F},\mathcal{G} ) $$
is an equivalence. Note that $T_\bullet \simeq \lim_n T_\bullet^n$ is an apaisant presentation. Hence we have
$$  \colim_{*\on{-pull}} \IndCoh( T_\bullet^n ) \simeq \IndCoh_{*}( T_\bullet )$$
and therefore
$$ \colim_{m\ge n} \on{Maps}( \mathcal{F}^m_\bullet,\mathcal{G}^m_\bullet  ) \to \on{Maps} ( \mathcal{F}_\bullet,\mathcal{G}_\bullet ). $$
By Theorem \ref{theorem-flat-descent-indcoh-w-ren}, we have 
$$\on{Maps} ( \mathcal{F},\mathcal{G} ) \simeq \on{Tot}_{\on{semi}} \on{Maps} ( \mathcal{F}_\bullet,\mathcal{G}_\bullet ),\; \on{Maps} ( \mathcal{F}^m,\mathcal{G}^m ) \simeq \on{Tot}_{\on{semi}} \on{Maps} ( \mathcal{F}_\bullet^m,\mathcal{G}^m_\bullet ).$$
Hence we only need to show
$$ \colim_{m\ge n} \on{Tot}_{\on{semi}} \on{Maps} ( \mathcal{F}_\bullet^m,\mathcal{G}^m_\bullet ) \to \on{Tot}_{\on{semi}}  \colim_{m\ge n} \on{Maps} ( \mathcal{F}_\bullet^m,\mathcal{G}^m_\bullet )    $$
is an equivalence. But this can be proved by using the same trick as in the proof of Proposition \ref{prop-general-push-w-ren}.

It remains to show the image of
$$ \colim_{*\on{-pull}} \IndCoh(T^n/H^n) \to \IndCoh_{*,\ren}(T/H) $$
generates the target. Note that
$$\colim_{*\on{-pull}} \IndCoh(T^n/H^n)  \simeq \colim_{m,*\on{-pull}} \colim_{n\ge m,*\on{-pull}} \IndCoh_{*,\ren}(T^m/H^n).$$
We only need to show:
\begin{itemize}
	\item[(i)] The image of
	$$ \colim_{*\on{-pull}} \IndCoh_{*,\ren}(T^m/H) \to \IndCoh_{*,\ren}(T/H) $$
	generates the target.

	\item[(ii)] For each $m$, the image of
	$$\colim_{n\ge m,*\on{-pull}} \IndCoh(T^m/H^n) \to \IndCoh_{*,\ren}(T^m/H)$$
	generates the target.
\end{itemize}
Note that $p_m: T/H \to T^m/H$ is affine schematic. Hence for any $\mathcal{M}\in \on{Coh}(T/H)$, we have a well-defined morphism
$$ \colim_m  (p_m)^*_{\ren} \circ  (p_m)_{*,\ren}( \mathcal{M}) \to \mathcal{M}.$$
To prove (i), we only need to show this is an equivalence. Note that each term in the colimit is bounded. Hence we only need to prove 
$$ \colim_m  (p_m)^* \circ  (p_m)_{*} \circ \ren( \mathcal{M}) \to \ren( \mathcal{M})$$
is an equivalence in $\IndCoh_{*}(T/H).$ By Theorem \ref{theorem-flat-descent-indcoh-w-ren}, we only need to prove its $*$-pullback along $T\to T/H$ is an equivalence. But this follows from the base-change isomorphisms and the equivalence
$$ \colim_{*\on{-pull}} \IndCoh_{*}(T^m) \to \IndCoh_{*}(T).$$
This proves (i).

It remains to prove (ii). Since $\on{Coh}(T^m/H)$ is closed under truncations, we only need to show $\on{Coh}(T^m/H)^\heartsuit$ is contained in the image. Write $q: T^m \to T^m/H$ and $q_n: T^m \to T^m/H^n$ for the projections. Using Theorem \ref{theorem-flat-descent-indcoh-w-ren} and the base-change isomorphisms, it is easy to see the adjoint pairs $(q^*, q_{*})$, $((q_n)^*, (q_{n})_*)$ are comonadic. We claim the natural transformation
$$  \colim_{n\ge m}  (q_n)^* \circ  (q_{n})_*\to q^*\circ q_{*} $$
between comonads on $\IndCoh(T^m)$ is an equivalence. To prove the claim, we only need to show it is an equivalence between underlying functors. Then it follows from the base-change isomorphisms and the equivalence
$$  \colim_{*\on{-pull}} \IndCoh_{*}(T^m\times_B H^n) \to \IndCoh_{*}(T^m\times_B H).$$
It follows from definition that 
$$\on{Coh}(T^m/H)^\heartsuit \simeq  [q^*\circ q_{*}]\on{-comod}(\IndCoh(T^m)) \times_{\IndCoh(T^m)} \on{Coh}(T^m)^\heartsuit $$
and
$$\on{Coh}(T^m/H^n)^\heartsuit \simeq  [q_n^*\circ (q_n)_{*}]\on{-comod}(\IndCoh(T^m)) \times_{\IndCoh(T^m)} \on{Coh}(T^m)^\heartsuit.$$
Note that the comonads $q^*\circ q_{*}$ and $q_n^*\circ (q_n)_{*}$ are t-exact. Now consider a given $(q^*\circ q_{*})$-comodule $\mathcal{M}$ such that the underlying object $\underline{\mathcal{M}}$ is contained in $\on{Coh}(T^m)^\heartsuit$. Since $\mathcal{M}$ is compact in $\IndCoh(T^m)$, the coaction morphism
$$\underline{\mathcal{M}} \to  q^*\circ q_{*}(\underline{\mathcal{M}}) \simeq \colim_{n\ge m}  (q_n)^* \circ  (q_{n})_* (\underline{\mathcal{M}})$$ factors through $\underline{\mathcal{M}} \to  (q_n)^* \circ  (q_{n})_* (\underline{\mathcal{M}})$ for some $n$. By enlarging $n$, we can assume this map is a coaction morphism\footnote{There are no higher compatibilities for coactions in an abelian category!}. By definition, $\mathcal{M}$ is the co-restriction of the above $q_n^*\circ (q_n)_{*}$-comodule structure on $\underline{\mathcal{M}}$ along $q_n^*\circ (q_n)_{*}\to q^*\circ q_{*}$. In other words, any object in $\on{Coh}(T^m/H)^\heartsuit$ is an image of an object in $\on{Coh}(T^m/H^n)^\heartsuit$. This finishes the proof.

\qed

\begin{definition}  
Let $B\in \,^{>-\infty} \on{AffSch}_{\on{aft}}$ be an affine base scheme. A weakly renormalizable prestack $S$ over $B$ is \emph{of quotient type} if it can be written as $S \simeq \colim_{\alpha\in I} S^\alpha$ as in Proposition \ref{prop-indcoh-on-ind-stack} such that each $S^\alpha$ is of the form $S^\alpha \simeq T/H$ as in Proposition \ref{prop-indcoh-on-quot-stack} with $T\to B$ of bounded Tor dimension.
\end{definition}

\begin{lemma} \label{lemma-vafp-quotient-stack}
Let $S= \colim_{\alpha\in I} S^\alpha$ be of quotient type over $B\in \,^{>-\infty} \on{AffSch}_{\on{aft}}$. Then for any $B'\to B$ in $\,^{>-\infty} \on{AffSch}_{\on{aft}}$, the fiber product $S\times_B B'$ is of quotient type over $B'$, and the map $S\times_B B' \to S$ is vafp.
\end{lemma}

\proof The first claim is obvious. To prove the second claim, since $S\times_B B' \to S$ is a base-change of $B'\to B$, we only need to show $S\times_B B' \to S$ is uafp. Let $U\in\on{Cov}(S)$ be a flat apaisant cover, by Remark \ref{rem-uafp-vs-sch&afp}, we only need to show there exists a presentation $U=\colim_\alpha U^\alpha$ of $U$ as an apaisant indscheme such that each $U^\alpha \times_B B'$ is eventually coconnective. Now $U^\alpha:= S^\alpha\times_ S U$ satisfies this condition because it is an apaisant scheme of bounded Tor dimension over $B$.

\qed

\begin{definition} For a laft prestack $B \in \,^{>-\infty} \on{AffSch}_{\on{aft}}$. A convergent prestack $Y_{B}$ over $B$ is \emph{of quotient type} if for \emph{any} test scheme $s:S\to B$ with $S\in \,^{>-\infty}\on{AffSch}_{\on{aft}}$, the fiber product $Y_s:= Y_{B} \times_{ B } S$ is of quotient type over $S$. By Lemma \ref{lemma-vafp-quotient-stack}, when $B$ is affine, this definition is compatible with the previous one.
\end{definition}

\begin{example}\label{exam-hecke-apaisant-quot-stack}
	$[\Arc G \backslash \Loop G / \Arc G]_{\on{Ran}_\dR} \to \Ran_\dR$ is of quotient type over $\Ran_\dR$. To see this, note that any $s:S\to \Ran_{\on{dR}}$ factors (non-canonically) through $X^I \to \on{Ran}_\dR$ for some $I$. Hence we only need to check $[\Arc G \backslash \Loop G / \Arc G]_{X^I}$ is of quotient type over $X^I$. Recall that $[\on{Gr}_{G}]_{X^I} = [\Loop G / \Arc G]_{X^I}$ is a \emph{classical} ind-finite type indscheme, hence there exists a presentation $[\Loop G / \Arc G]_{X^I} \simeq \colim Y_\alpha$ such that each $Y_\alpha$ is stablized by the (left) $[\Arc G]_{X^I}$-action. Then we only need to check each $[\Arc G]_{X^I} \backslash Y_\alpha$ is of quotient type over $X^I$. Let $[\Arc G]_{X^I} \simeq \lim_m H_m$ be as in Section \ref{sect:unital-fact-spaces}. For fixed $\alpha$, there exists a large enough $m$ such that the action of $[\Arc G]_{X^I}$ on $Y_\alpha$ factors through $H_m$. Hence we obtain a desired cofiltered system of actions $\{H_n, Y_\alpha\}_{n\ge m}$ whose limit is $([\Arc G]_{X^I},Y_\alpha)$. It remains to check $Y_\alpha \to X^I$ is of bounded Tor dimension, but this is obvious because $Y_\alpha$ is classical and $X^I$ is smooth.
\end{example}

\begin{remark} If $Y_B$ is of quotient type over $B \in \,^{>-\infty} \on{Sch}_{\on{aft}}$, then $Y_B$ is weakly renormalizable. To see this, let $\{B_i \to B \}$ be a finite Zariski cover of $B$. By definition, each $Y_{B_i}:= Y_B \times_{B} B_i$ is weakly renormalizable. Let $S_i \to Y_{B_i}$ be a flat apaisant cover. Then the disjoint union $\sqcup S_i \to Y_B$ is a flat apaisant cover.
\end{remark}

\begin{corollary} \label{cor-box-strict-quotient-stack}
If $Y_B$ (resp. $Y'_{B'}$) is of quotient type over $B \in \,^{>-\infty} \on{Sch}_{\on{aft}}$ (resp. $B'\in \,^{>-\infty} \on{Sch}_{\on{aft}}$). Then the external tensor product functor
$$ \IndCoh_{*,\ren}(Y_B) \otimes \IndCoh_{*,\ren}(Y'_{B'}) \to \IndCoh_{*,\ren}(Y_B \times Y'_{B'}) $$
is an equivalence.
\end{corollary}

\proof Follows from Proposition \ref{prop-indcoh-ren-zariski}, Proposition \ref{prop-indcoh-on-ind-stack}, Proposition \ref{prop-indcoh-on-quot-stack} and \cite[Corollary 4.2.3]{drinfeld2013some}.

\qed

\subsection{IndCoh as a crystal of categories}
\label{sect:indcoh-crys}

In this subsection, we show that for a prestack $Y_{B_{\on{dR}}}$ over a laft de-Rham prestack $B_{\on{dR}}$ such that its fibers $Y_x$ at $x\in B_{\on{dR}}$ are nice enough, there is a canonical \emph{crystal of categories} on $B_{\on{dR}}$ such that its fibers at $x$ is $\IndCoh_{*,\ren}(Y_x)$ (resp. $\IndCoh^{!,\ren}(Y_x)$). In particular, given a (nice enough) factorization space, we can obtain two factorization categories that are dual to each other.

\begin{warning} If we replace $B_\dR$ by a general laft prestack $B$, then we can \emph{not} obtain a \emph{quasi-coherent} sheaf of of categories on $B$ such that its fibers at $x$ is $\IndCoh_{*,\ren}(Y_x)$ or $\IndCoh^{!,\ren}(Y_x)$.
\end{warning}

\begin{construction} \label{const-indcoh-base-action}
Let $Y_B$ be of quotient type over $B \in \,^{>-\infty} \on{Sch}_{\on{aft}}$. There is a canonical $(\IndCoh(B),\overset{!}\otimes)$-action on $\IndCoh_{*,\ren}(Y_B)$ whose action functor is
$$ \IndCoh(B) \otimes \IndCoh_{*,\ren}(Y_B) \simeq \IndCoh_{*,\ren}(B\times Y_B) \xrightarrow{!\on{-pull}} \IndCoh_{*,\ren}(Y_B).$$
Note that the $!$-pullback functor is well-defined because $Y_B\to B\times Y_B$ is vafp.

The above $\IndCoh(B)$-action is functorial in $Y_B$ in theo following sense. For $f: Y_B \to Y_B'$ defined over $B$, the base-change isomorphims provide $\IndCoh(B)$-linear structure on $f_{*,\on{ren}}, f_{\on{ren}}^!, f^*_{\on{ren}}$ as long as whenever functors are defined.
\end{construction}

\begin{construction} \label{const-indcoh-crystal}
Let $Y_B$ be of quotient type over a laft prestack $B$. For any laft prestack $T \to B$ over $B$, we define
$$ \mathbf{\Gamma}( t:T\to B, \mathcal{I}nd\mathcal{C}oh_{*,\ren}(Y) ) := \lim_{!\on{-pull}, S\to T} \IndCoh_{*,\ren}(Y_s),$$
where the limit is taken for any test scheme $S\to T$ with $S\in \,^{>-\infty}\on{AffSch}_{\on{aft}}$, and $s$ is the composition $S\to T \to B$. Note that the $!$-pullback functors are well-defined because each $Y_{s'} \to Y_s$ is vafp (Lemma \ref{lemma-vafp-quotient-stack}). When there is little danger ambiguity, we also write
$$\mathbf{\Gamma}( T, \mathcal{I}nd\mathcal{C}oh_{*,\ren}(Y) ) $$
for this category.

By Construction, $\mathbf{\Gamma}( T, \mathcal{I}nd\mathcal{C}oh_{*,\ren}(Y) )$ has a module structure for the symmetric monoidal category
$$ \IndCoh(T) \simeq \lim_{!\on{-pull}, S\to T_{\on{dR}}} \IndCoh(S),$$
and such module structure is functorial in $T$, i.e., for $T' \to T$ defined over $B_\dR$, the retriction functor
$$\mathbf{\Gamma}( T, \mathcal{I}nd\mathcal{C}oh_{*,\ren}(Y) )  \to \mathbf{\Gamma}( T', \mathcal{I}nd\mathcal{C}oh_{*,\ren}(Y) ) $$
intertwines the symmetric monoidal functor $\IndCoh(T') \to \IndCoh(T)$.
\end{construction}

\begin{remark}  \label{rem-section-on-non-affine}
Let $T\in \,^{>-\infty}\on{Sch}_{\on{aft}}$ be defined over $B$. Using Lemma \ref{prop-indcoh-ren-zariski}, one can show the functor
$$ \IndCoh_{*,\ren}(Y_t) \to \mathbf{\Gamma}( T, \mathcal{I}nd\mathcal{C}oh_{*,\ren}(Y) ) $$
is an equivalence.
\end{remark}

\begin{remark} \label{rem-dr-section-compactly-generated}
Let $T_\dR$ be defined over $B$ such that $T\in \,^{>-\infty}\on{Sch}_{\on{aft}}$. By \cite[Chapter 2, Corollary 4.3.4]{GR-DAG2} and 
, we have
\begin{eqnarray*}
\mathbf{\Gamma}( T_{\on{dR}}, \mathcal{I}nd\mathcal{C}oh_{*,\ren}(Y) ) &\simeq &
\lim_{!\on{-pull},S \in (>^{-\infty} \on{Sch}_{\on{aft}})_{\on{nil-isom\,to\,}T_\dR} } \mathbf{\Gamma}( S, \mathcal{I}nd\mathcal{C}oh_{*,\ren}(Y) ) \\
&\simeq &
\lim_{!\on{-pull},S \in (>^{-\infty} \on{Sch}_{\on{aft}})_{\on{nil-isom\,to\,}T_\dR}} \IndCoh_{*,\ren}(Y_s) \\
&\simeq &
\colim_{*\on{-push}, S \in (>^{-\infty} \on{Sch}_{\on{aft}})_{\on{nil-isom\,to\,}T_\dR}} \IndCoh_{*,\ren}(Y_s).
\end{eqnarray*}
Note that the last line is compactly generated hence dualizable. It folllows that 
\begin{eqnarray*}
\mathbf{\Gamma}( T_{\on{dR}}, \mathcal{I}nd\mathcal{C}oh_{*,\ren}(Y) )^\vee &\simeq & \lim_{!\on{-pull}, S \in (>^{-\infty} \on{Sch}_{\on{aft}})_{\on{nil-isom\,to\,}T_\dR}} \IndCoh^{!,\ren}(Y_s)
 \\
&\simeq & \lim_{!\on{-pull}, S \in (>^{-\infty} \on{Sch}_{\on{aft}})_{/T_\dR}} \IndCoh^{!,\ren}(Y_s).
\end{eqnarray*}
\end{remark}

\begin{proposition} \label{prop-indcoh-*-ren-crystal}
Let $Y_{B_\dR}$ be of quotient type over a laft de-Rham prestack $B_\dR$. Then 
$$ [T\in (^{>-\infty}\on{AffSch}_{\on{aft}})_{/B}] \mapsto [\mathbf{\Gamma}( T_{\on{dR}}, \mathcal{I}nd\mathcal{C}oh_{*,\ren}(Y) ) \in \DMod(T)\on{-mod}] $$
defines a crystal of categories on $B_\dR$. In other words, for any $T' \to T$, the functor
$$ \mathbf{\Gamma}( T_{\on{dR}}, \mathcal{I}nd\mathcal{C}oh_{*,\ren}(Y) ) \otimes_{\DMod(T)} \DMod(T') \to \mathbf{\Gamma}( T'_{\on{dR}}, \mathcal{I}nd\mathcal{C}oh_{*,\ren}(Y) ) $$
is an equivalence.

\end{proposition}

\begin{remark} The proof below actually works for any \emph{separated} scheme $T\in (^{>-\infty}\on{Sch}_{\on{aft,sep}})_{/B}$.
\end{remark}

\proof It follows from definition that the desired functor is equivalent to the following composition
\begin{eqnarray*}
	& &		 \mathbf{\Gamma}( T_{\on{dR}}, \mathcal{I}nd\mathcal{C}oh_{*,\ren}(Y) ) \otimes_{\DMod(T)} \DMod(T') \\
	&\simeq & [\lim_{ S\to T_\dR }  \mathbf{\Gamma}( S, \mathcal{I}nd\mathcal{C}oh_{*,\ren}(Y) )] \otimes_{\DMod(T)} [\lim_{V\to T_\dR }  \IndCoh(V) ]   \\
	&\to & \lim_{ S\to T_\dR, V\to T'_\dR } [\mathbf{\Gamma}( S, \mathcal{I}nd\mathcal{C}oh_{*,\ren}(Y) ) \otimes_{\DMod(T)} \IndCoh(V) ] \\
	&\simeq & \lim_{ S\to T_\dR, V\to T'_\dR } [\mathbf{\Gamma}( S, \mathcal{I}nd\mathcal{C}oh_{*,\ren}(Y) ) \otimes_{\IndCoh(S)} \IndCoh(S) \otimes_{\DMod(T)} \IndCoh(V)]  \\
	& \to & \lim_{ S\to T_\dR, V\to T'_\dR }[\mathbf{\Gamma}( S, \mathcal{I}nd\mathcal{C}oh_{*,\ren}(Y) ) \otimes_{\IndCoh(S)} \IndCoh(S\times_{T_\dR} V)] \\
	& \to & \lim_{ S\to T_\dR, V\to T'_\dR } [\mathbf{\Gamma}( S\times_{T_\dR} V, \mathcal{I}nd\mathcal{C}oh_{*,\ren}(Y) )] \\
	& \simeq &  \mathbf{\Gamma}( T'_{\on{dR}}, \mathcal{I}nd\mathcal{C}oh_{*,\ren}(Y) ),
\end{eqnarray*}
where the last equivalence is due to 
$$ T'_{\on{dR}} \simeq \colim_{ S\to T_\dR, V\to T'_\dR }  S\times_{T_\dR} V  $$
as a prestack over $B_\dR$.

Hence we only need to prove:
\begin{itemize}
	\item[(i)] The functors 
	$\mathcal{C} \mapsto \mathcal{C} \otimes_{\DMod(T)} \IndCoh(V)$ and $\mathcal{D} \mapsto  \mathbf{\Gamma}( S, \mathcal{I}nd\mathcal{C}oh_{*,\ren}(Y) )] \otimes_{\DMod(T)} \mathcal{D}$ commute with limits;
	\item[(ii)] The functor 
$$\mathbf{\Gamma}( S, \mathcal{I}nd\mathcal{C}oh_{*,\ren}(Y) ) \otimes_{\DMod(T)} \IndCoh(V) \to \mathbf{\Gamma}( S\times_{T_\dR} V, \mathcal{I}nd\mathcal{C}oh_{*,\ren}(Y) )$$
	is an equivalence.
\end{itemize}

Too prove claim (i), we only need to show $\IndCoh(V)$ and $\mathbf{\Gamma}( S, \mathcal{I}nd\mathcal{C}oh_{*,\ren}(Y) )$ are dualizable as $\DMod(T)$-module categories. Note that they are compactly generated hence dualizable as plain categories. Then \cite[Lemma B.3.3]{chen2021thesis} implies the desired claim.

We now prove (ii). By Corollary \ref{cor-box-strict-quotient-stack}, we have
$$ \mathbf{\Gamma}( S, \mathcal{I}nd\mathcal{C}oh_{*,\ren}(Y) )\otimes \IndCoh(V) \simeq  \mathbf{\Gamma}( S\times V, \mathcal{I}nd\mathcal{C}oh_{*,\ren}(Y) ),$$
where $S\times V$ is defined over $B_\dR$ via the composition $S\times V \to S\to B_\dR$. Hence the desired functor can be identified with the functor
$$b: \mathbf{\Gamma}( S\times V, \mathcal{I}nd\mathcal{C}oh_{*,\ren}(Y) ) \otimes_{\IndCoh(T_\dR\times T_\dR)} \IndCoh(T_\dR) \to \mathbf{\Gamma}( (S\times_{T_\dR} V, \mathcal{I}nd\mathcal{C}oh_{*,\ren}(Y) ).$$
Consider the functor
$$ \textbf{Id}\otimes \Delta^!:  \mathbf{\Gamma}( S\times V, \mathcal{I}nd\mathcal{C}oh_{*,\ren}(Y) ) \to\mathbf{\Gamma}( S\times V, \mathcal{I}nd\mathcal{C}oh_{*,\ren}(Y) ) \otimes_{\IndCoh(T_\dR\times T_\dR)} \IndCoh(T_\dR).$$
Note that it has a fully faithful left adjoint $\textbf{Id}\otimes \Delta_*$. Hence we only need to show $b\circ (\textbf{Id}\otimes \Delta^!)$ also has a fully faithful left adjoint, and the two co-localization functors on 
$$\mathbf{\Gamma}( S\times V, \mathcal{I}nd\mathcal{C}oh_{*,\ren}(Y) ) \simeq \IndCoh_{*,\ren}( Y_s\times V ). $$ can be identified. 

By definition, the endo-functor $\textbf{Id}\otimes (\Delta_*\circ \Delta^!)$ is given by the action of $\Delta_*(\omega) \in \IndCoh(T_\dR\times T_\dR)$ on 
$\IndCoh_{*,\ren}( Y_s\times V )$. Let $\delta: S\times_{T_\dR} V \to S\times V$ be the base-change of $\Delta$. Note that $S\times_{T_\dR} V$ is the formal completion of $S\times V$ along $S^{\on{cl,red}}\times_T V^{\on{cl,red}}$. Hence by \cite[Proposition 6.5.2]{gaitsgory2014dg}, $S\times_{T_\dR} V$ is a (laft) indscheme. The base-change isomorphism implies $\Delta_*(\omega)$ is sent to $\delta_*(\omega)$ by the (symmetric monoidal) $!$-pullback functor
$$  \IndCoh(T_\dR\times T_\dR)\to \IndCoh(S\times V).$$
Hence the endo-functor $\textbf{Id}\otimes (\Delta_*\circ \Delta^!)$ is also given by the action of $\delta_*(\omega) \in \IndCoh(S\times  V)$ on 
$\IndCoh_{*,\ren}( Y_s\times V )$. Now consider the following base-change of $\delta$:
$$ d: Y_s\times_{T_\dR} V\to Y_s\times V.$$
Note that $Y_s\times_{T_\dR} V$ is weakly renormalizable\footnote{Proof: Let $U$ be a flat apaisant cover of $Y_s$. By Lemma \ref{lem-fiber-product-apai-indsch-is-apai-if-afp}, $U\times_{S} ( S\times_{T_\dR} V )$ is an apaisant indscheme, which is a flat apaisant cover of $Y_s\times_S (S\times_{T_\dR} V)$.}. Hence by the base-change isomorphism (see Lemma \ref{lem-base-change-ren}), the above action can be identified with $d_{*,\ren}\circ d_\ren^?$.

On the other hand, note that the base-change of $d$ along itself is an isomophism. The base-change isomorphism implies $d_{*,\ren}$ is fully faithful. Hence we only need to identify the functor
$$ \mathbf{\Gamma}( S\times V, \mathcal{I}nd\mathcal{C}oh_{*,\ren}(Y) )\to \mathbf{\Gamma}( S\times_{T_\dR} V, \mathcal{I}nd\mathcal{C}oh_{*,\ren}(Y) ) $$
with $ d_\ren^?$. Let $S\times_{T_\dR} V \simeq \colim_i Z_i$ be a presentation of the laft indscheme such that $Z_i \in \,^{>-\infty}\on{Sch}_{\on{aft}}$. It follows from definition (and Remark \ref{rem-section-on-non-affine}) that the left adjoint of the above functor is given by the following composition:
\begin{eqnarray*}
\mathbf{\Gamma}( S\times_{T_\dR} V, \mathcal{I}nd\mathcal{C}oh_{*,\ren}(Y) ) 
& \simeq & \lim_{!\on{-pull}} \IndCoh_{*,\ren}(  Y_s \times_S Z_i ) \\
& \simeq & \colim_{*\on{-push}} \IndCoh_{*,\ren}(  Y_s \times_S Z_i ) \\
& \to &  \IndCoh_{*,\ren}(  Y_s \times_S \colim_i Z_i ) \\
& = & \IndCoh_{*,\ren}(Y_s\times_{T_\dR} V) \\
& \to & \IndCoh_{*,\ren}(Y_s\times V).
\end{eqnarray*}
Proposition \ref{prop-indcoh-on-ind-stack} implies the third functor is an equivalence, which finishes the proof.

\qed

\begin{definition} \label{def-crys-category-indcoh}
Let $Y_{B_\dR}$ be of quotient type over a laft de-Rham prestack $B_\dR$. We denote by $\mathcal{I}nd\mathcal{C}oh_{*,\ren}(Y)$ the crystal of categories over $B_\dR$ provided by Proposition \ref{prop-indcoh-*-ren-crystal}. For $t:T\to B$, we denote the corresponding section as
$$  [ \mathcal{I}nd\mathcal{C}oh_{*,\ren}(Y) ]_{t_\dR}:=\mathbf{\Gamma}( t_\dR, \mathcal{I}nd\mathcal{C}oh_{*,\ren}(Y)   ) . $$
When there is no danger of ambiguity, we also use $ [ \mathcal{I}nd\mathcal{C}oh_{*,\ren}(Y) ]_{T_\dR}$.

By Remark \ref{rem-dr-section-compactly-generated} and \cite[Lemma B.3.3]{chen2021thesis}, $\mathcal{I}nd\mathcal{C}oh_{*,\ren}(Y)$ is a dualizable object in the symmetric monoidal category $\on{CrsyCat}(B_\dR)$, and the sections of its dual $[\mathcal{I}nd\mathcal{C}oh_{*,\ren}(Y)]^\vee$ are given by
$$ \mathbf{\Gamma}( T_\dR, [\mathcal{I}nd\mathcal{C}oh_{*,\ren}(Y)]^\vee ) \simeq  \mathbf{\Gamma}( T_\dR, \mathcal{I}nd\mathcal{C}oh_{*,\ren}(Y))^\vee \simeq \lim_{!\on{-pull},S\to T_\dR} \IndCoh^{!,\ren}(Y_s). $$
It follows from the proof of Proposition \ref{prop-indcoh-*-ren-crystal} that the above identification is functorial in $T_\dR$. Hence we use the notation
$$\mathcal{I}nd\mathcal{C}oh^{!,\ren}(Y):= [\mathcal{I}nd\mathcal{C}oh_{*,\ren}(Y)]^\vee.$$
\end{definition}

\begin{construction} \label{const-crystal-pullback}
Let $Y_{B_\dR}$ be of quotient type over a laft de-Rham prestack $B_\dR$ and $f:B'_\dR \to B_\dR$ be any map between laft de-Rham prestacks. Then the fiber product $Y'_{B'_\dR}:= Y_{B_\dR} \times_{B_\dR} B'_\dR$ is of quotient type over $B'_\dR$ and Proposition \ref{prop-indcoh-*-ren-crystal} provides an identification
$$ \mathbf{cores}_f( \mathcal{I}nd\mathcal{C}oh_{*,\ren}(Y)) \simeq \mathcal{I}nd\mathcal{C}oh_{*,\ren}(Y')$$
as crystals of categories on $B'_\dR$. Therefore we also have
$$ \mathbf{cores}_f( \mathcal{I}nd\mathcal{C}oh^{!,\ren}(Y)) \simeq \mathcal{I}nd\mathcal{C}oh^{!,\ren}(Y').$$
\end{construction}

\begin{construction}  \label{const-crystal-boxtimes}
Let $Y_{B_\dR}$ (resp. $Y'_{B'_\dR}$) be of quotient type over a laft de-Rham prestack $B_\dR$ (resp. $B_\dR'$). By Corollary \ref{cor-box-strict-quotient-stack} and Lemma \ref{lemma-vafp-quotient-stack}, 
$$ (Y\boxtimes Y')_{ B_\dR\times B_\dR' }:= Y_{B_\dR}\times Y'_{B'_\dR}$$
is of quotient type over $B_\dR\times B_\dR'$. Moreover, for $T_\dR$ (resp. $T'_\dR$) over $B_\dR$ (resp. $B_\dR'$) such that $T, T'\in \,^{>-\infty}\on{Sch}_{\on{aft}}$, we have an equivalence
\begin{eqnarray*}
  & & \mathbf{\Gamma}( T_{\on{dR}}, \mathcal{I}nd\mathcal{C}oh_{*,\ren}(Y) )  \otimes \mathbf{\Gamma}( T'_{\on{dR}}, \mathcal{I}nd\mathcal{C}oh_{*,\ren}(Y') ) \\
  & \simeq &  [\lim_{!\on{-pull}, S\to T_{\on{dR}}} \IndCoh_{*,\ren}(Y_s)] \otimes [\lim_{!\on{-pull}, S'\to T'_{\on{dR}}} \IndCoh_{*,\ren}(Y'_{s'})] \\
   & \simeq &  \lim_{!\on{-pull}, S\to T_{\on{dR}}, S'\to T'_{\on{dR}}} [\IndCoh_{*,\ren}(Y_s ) \otimes  \IndCoh_{*,\ren}(Y'_{s'})] \\
   & \simeq &  \lim_{!\on{-pull}, S\to T_{\on{dR}}, S'\to T'_{\on{dR}}} \IndCoh_{*,\ren}(Y_s \times Y'_{s'}) \\
	& \simeq & \mathbf{\Gamma}( T_{\on{dR}}\times T'_\dR, \mathcal{I}nd\mathcal{C}oh_{*,\ren}(Y\boxtimes Y') ),
 \end{eqnarray*}
where the second line is equivalent to the third because $\mathbf{\Gamma}( T_{\on{dR}}, \mathcal{I}nd\mathcal{C}oh_{*,\ren}(Y) )$ and $\mathbf{\Gamma}( T'_{\on{dR}}, \mathcal{I}nd\mathcal{C}oh_{*,\ren}(Y') )$ are compactly generated. Hence we obtain an equivalence
$$ \mathcal{I}nd\mathcal{C}oh_{*,\ren}(Y)\boxtimes \mathcal{I}nd\mathcal{C}oh_{*,\ren}(Y) \simeq \mathcal{I}nd\mathcal{C}oh_{*,\ren}(Y\boxtimes Y') $$
as crystals of categories on $B_\dR \times B'_\dR$. Therefore we also have
$$ \mathcal{I}nd\mathcal{C}oh^{!,\ren}(Y)\boxtimes \mathcal{I}nd\mathcal{C}oh^{!,\ren}(Y) \simeq \mathcal{I}nd\mathcal{C}oh^{!,\ren}(Y\boxtimes Y').$$
\end{construction}

\begin{construction}  \label{const-crystal-tensor}
Let $Y_{B_\dR}$ and $Y'_{B_\dR}$ be of quotient type over a laft de-Rham prestack $B_\dR$. Combining Construction \ref{const-crystal-boxtimes} and Construction \ref{const-crystal-pullback}, we see 
$$ (Y\times Y')_{B_\dR}:= Y_{B_\dR}\times_{ B_\dR } Y'_{B_\dR}$$
is of quotient type over $B_\dR$, and we have
\[
\begin{aligned}
 \mathcal{I}nd\mathcal{C}oh_{*,\ren}(Y)\otimes \mathcal{I}nd\mathcal{C}oh_{*,\ren}(Y') \simeq  \mathcal{I}nd\mathcal{C}oh_{*,\ren}(Y\times Y') \\
 \mathcal{I}nd\mathcal{C}oh^{!,\ren}(Y)\otimes \mathcal{I}nd\mathcal{C}oh^{!,\ren}(Y) \simeq \mathcal{I}nd\mathcal{C}oh^{!,\ren}(Y\times Y').
\end{aligned}
\]
\end{construction}

\begin{construction}  \label{constr-fact-prestack-indcoh-fact-cat}
Let $Y_{\Ran_\dR}\to \Ran_\dR$ be a factorization space over $\Ran_\dR$ such that $Y_{\Ran_\dR}$ is of quotient type over $\Ran_\dR$. Combining Construction \ref{const-crystal-boxtimes} and Construction \ref{const-crystal-pullback}, we obtain a factorization structure on the crystals of categories $\mathcal{I}nd\mathcal{C}oh_{*,\ren}(Y)$ and $\mathcal{I}nd\mathcal{C}oh^{!,\ren}(Y)$ over $\Ran_\dR$, which are dual to each other in the symmetric monoidal category of factorization categories over $\Ran_\dR$. By definition, for any closed point $x\in X$, the fibers of these factorization categories at $x$ are $\IndCoh_{*,\ren}(Y_x)$ and $\IndCoh^{!,\ren}(Y_x)$.
\end{construction}

For the purpose of this paper, we also need to a \emph{unital} factorization version of the above construction.

\begin{construction} \label{constr-unit-fact-prestack-indcoh-fact-cat}
Let $Y_{ \Ran_{\on{un,dR}}}$ be a corr-unital space over $\Ran_{\on{un,dR}}$ such that:
\begin{itemize}
	\item[$(\spadesuit)$] For any affine test scheme $S\in \,^{>-\infty}\on{AffSch}_{\on{aft}}$ and a morphism $ \mathfrak{s}: \varphi \to \psi$ in $ \Ran_{\on{un,dR}}(S) $, the correspondence
	$$ Y_\psi\xleftarrow{\on{proj}_1}   Y_{\mathfrak{s}}   \xrightarrow{\on{proj}_0}   Y_\varphi $$
	is a correspondence in $(\on{PreStk}_{\on{w.ren}})_{/S}$ such that 
	\begin{itemize}
			\item $ Y_\psi,Y_{\mathfrak{s}},Y_\varphi$ are of quotient type over $S$;
			\item $\on{proj}_0$ is $*$-pullable;
			\item $\on{proj}_1$ is $(?,\ren)$-pullable.
	\end{itemize}
\end{itemize}
Then we obtain a crystal of categories $\mathcal{I}nd\mathcal{C}oh_{*,\ren}(Y)$ over $\Ran_{\on{un,dR}}$ such that:
\begin{itemize}
	\item Its section at $f_\dR:T_\dR \to \Ran_{\on{un,dR}}$ is given by
$$ \mathbf{\Gamma}( f_\dR,  \mathcal{I}nd\mathcal{C}oh_{*,\ren}(Y) ) := \lim_{!\on{-pull},u:S\to T_\dR} \IndCoh_{*,\ren} (Y_{ f_\dR\circ u }) \in \DMod(T)\on{-mod}.$$
Note that this is just the global section of $\mathcal{I}nd\mathcal{C}oh_{*,\ren}(Z)$ defined in Definition \ref{def-crys-category-indcoh}, where $Z_{T_\dR}:= Y_{ { \Ran_{\on{un,dR}}} }\times_{  { \Ran_{\on{un,dR}}} } T_\dR $ is of quotient type over $T_\dR$.

	\item For a morphism $\mathfrak{t}:f\to g$ in $\Ran_{\on{un}}(T)$, the corresponding $\DMod(T)$-linear functor (a.k.a. the \emph{unital} functor)
$$ \mathfrak{t}_!:  \mathbf{\Gamma}( f_\dR,  \mathcal{I}nd\mathcal{C}oh_{*,\ren}(Y) ) \to  \mathbf{\Gamma}( g_\dR,  \mathcal{I}nd\mathcal{C}oh_{*,\ren}(Y) )   $$
is given by
\begin{equation} \label{eqn-unit-structure-fact-indcoh}
\begin{aligned}
  \lim_{!\on{-pull},u:S\to T_\dR} \IndCoh_{*,\ren} (Y_{ f_\dR\circ u }) \xrightarrow{*\on{-pull}}  \lim_{!\on{-pull},u:S\to T_\dR} \IndCoh_{*,\ren} (Y_{ \mathfrak{t}_\dR \bigstar \mathbf{Id}_u })\to \\ \xrightarrow{*\on{-push}}  \lim_{!\on{-pull},u:S\to T_\dR} \IndCoh_{*,\ren} (Y_{   g_\dR\circ u   }).
\end{aligned}
\end{equation}
Note that the above $*$-pullback functors are well-defined because each $Y_{ \mathfrak{t}_\dR \bigstar \mathbf{Id}_u } \to Y_{ f_\dR\circ u }$ is assumed be to $*$-pullable. Also note tha the above $*$-pullback and $*$-pushforward functors are compatible with the limit diagrams because of the base-change isomorphisms.
\end{itemize}
Note that the fiber of $\mathcal{I}nd\mathcal{C}oh_{*,\ren}(Y)$ at $x\in X$ is $\IndCoh_{*,\ren}(Y_x)$.

As before, each section $\mathbf{\Gamma}( f_\dR,  \mathcal{I}nd\mathcal{C}oh_{*,\ren}(Y) )$ is dualizable and
$$
 \mathbf{\Gamma}( f_\dR,  \mathcal{I}nd\mathcal{C}oh_{*,\ren}(Y) )^\vee \simeq  \lim_{!\on{-pull},u:S\to T_\dR} \IndCoh^{!,\ren} (Y_{ f_\dR\circ u }) \in \DMod(T)\on{-mod}.
  $$
Via this identification, the composition (\ref{eqn-unit-structure-fact-indcoh}) is conjugate to\footnote{To see this, replace the limits in (\ref{eqn-unit-structure-fact-indcoh}) by certain colimits as in Remark \ref{rem-dr-section-compactly-generated}.
}
\begin{equation} \label{eqn-unit-structure-fact-indcoh-!}
\begin{aligned}
  \lim_{!\on{-pull},u:S\to T_\dR} \IndCoh^{!,\ren} (Y_{ f_\dR\circ u })  \xrightarrow{!\on{-pull}} \lim_{!\on{-pull},u:S\to T_\dR} \IndCoh^{!,\ren} (Y_{ \mathfrak{t}_\dR \bigstar \mathbf{Id}_u })\to \\ \xrightarrow{!\on{-push}}  \lim_{!\on{-pull},u:S\to T_\dR} \IndCoh^{!,\ren} (Y_{   g_\dR\circ u   }).
\end{aligned}
\end{equation}
Therefore we obtain a crystal of category $ \mathcal{I}nd\mathcal{C}oh^{!,\ren}(Y)$ over $\on{Ran}_{\on{un,dR}}$ whose sections are given by
$$
 \mathbf{\Gamma}( f_\dR,  \mathcal{I}nd\mathcal{C}oh^{!,\ren}(Y) ) \simeq  \lim_{!\on{-pull},u:S\to T_\dR} \IndCoh^{!,\ren} (Y_{ f_\dR\circ u }) \in \DMod(T)\on{-mod}.
  $$
 and whose unital functors are given by (\ref{eqn-unit-structure-fact-indcoh-!}). It follows formally that $\mathcal{I}nd\mathcal{C}oh_{*,\ren}(Y)$ and $ \mathcal{I}nd\mathcal{C}oh^{!,\ren}(Y) $ are dual to each other in the category of crystal of categories over $\Ran_{\on{un,dR}}$ and \emph{lax-unital} functors. See Remark \ref{rem-duality-unital-factorization-category} below for more details.

If $Y$ is further equipped with a factorization structure, then we obtain unital factorization categories $\mathcal{I}nd\mathcal{C}oh_{*,\ren}(Y)$ and $ \mathcal{I}nd\mathcal{C}oh^{!,\ren}(Y) $.
\end{construction}

\begin{remark} \label{rem-constr-unit-fact-prestack-indcoh-fact-cat}
Using the fact that $X$ is smooth, one can show that for a given $S\in \,^{>-\infty}\on{AffSch}_{\on{aft}}$ and a given morphism $ \mathfrak{s}: \varphi \to \psi$ in $ \Ran_{\on{un,dR}}(S) $, there exists an injective morphism $I\to J$ between finite set and a map $S\to X^J$ such thats $\mathfrak{s}$ can be identified with the horizontal composition of
$$
\xymatrix{
X^J \ar[rr] \ar[rd] & & \Ran_{\on{un,dR}} \\
& X^I \ar[ru] \ar@{=>}[u]
}
$$
with $S\to X^J$. Hence to verify $Y$ satisfies $(\spadesuit)$, we only need to check for any injective $\phi:I\to J$, the correspondence
$$  Y_{J} \gets Y_{\phi}  \to Y_{I}\times_{X^I} X^J$$
satisfies the required conditions in $(\spadesuit)$.
\end{remark}

\begin{construction} \label{constr-functorial-unital-fact-indcoh}
Construction \ref{constr-unit-fact-prestack-indcoh-fact-cat} is functorial in $Y$ in the following sense. For a morphism\footnote{We do not assume $f$ is strictly unital or strictly co-unital for now.} $f_{ \Ran_{\on{un,dR}}}: Y_{ \Ran_{\on{un,dR}}}\to Z_{ \Ran_{\on{un,dR}}}$ between corr-unital spaces over $\Ran_{\on{un,dR}}$ such that both $Y$ and $Z$ satisfy $(\spadesuit)$, then:
\begin{itemize}
	\item[(a)] We have a canonical lax unital functor
	$$ f_{*,\ren}: \mathcal{I}nd\mathcal{C}oh_{*,\ren}(Y)\to \mathcal{I}nd\mathcal{C}oh_{*,\ren}(Z) $$
  whose fiber at a closed point $x\in X$ is the $*$-pushforward functor
	$$ (f_x)_{*,\ren}: \IndCoh_{*,\ren}(Y_x) \to  \IndCoh_{*,\ren}(Z_x). $$
	Dually, we have a lax unital functor
	$$ f^{!,\ren}: \mathcal{I}nd\mathcal{C}oh^{!,\ren}(Z)\to \mathcal{I}nd\mathcal{C}oh^{!,\ren}(Y).$$
	When $f$ is strictly unital, $f_{*,\ren}$ is strictly unital. When $f$ is strictly co-unital, $f^{!,\ren}$ is strictly unital.

	\item[(b)] If $f$ is strictly co-unital and if for any $\phi:I\to J$ between finite sets, the map $f_\phi: Y_\phi \to Z_\phi$ is $*$-pullable, then we have a canonical \emph{strictly} unital functor
	$$ f^*_{\ren}: \mathcal{I}nd\mathcal{C}oh_{*,\ren}(Z)\to \mathcal{I}nd\mathcal{C}oh_{*,\ren}(Y) $$
	whose fiber at a closed point $x\in X$ is the $*$-pullback functor
	$$ (f_x)^*_{\ren}: \IndCoh_{*,\ren}(Z_x) \to  \IndCoh_{*,\ren}(Y_x). $$
	Dually, we have a \emph{lax} unital functor
	$$ f^{\ren}_?: \mathcal{I}nd\mathcal{C}oh^{!,\ren}(Y)\to \mathcal{I}nd\mathcal{C}oh^{!,\ren}(Z).$$
	When $f$ is also strictly unital, $f^{\ren}_?$ is strictly unital.

	\item[(c)] If $f$ is strictly unital and if for any $\phi:I\to J$ between finite sets, the map $f_\phi: Y_\phi \to Z_\phi$ is $(?,\ren)$-pullable, then we have a canonical \emph{strictly} unital functor
$$ f_!^\ren: \mathcal{I}nd\mathcal{C}oh_{*,\ren}(Y)\to \mathcal{I}nd\mathcal{C}oh_{*,\ren}(Z)$$
	whose fiber at a closed point $x\in X$ is the functor
	$$ (f_x)_!^\ren: \mathcal{I}nd\mathcal{C}oh_{*,\ren}(Y_x)\to \mathcal{I}nd\mathcal{C}oh_{*,\ren}(Z_x).$$
	Dually, we have a \emph{lax} unital functor
	$$ f_\ren^?: \mathcal{I}nd\mathcal{C}oh_{*,\ren}(Z)\to \mathcal{I}nd\mathcal{C}oh_{*,\ren}(Y).$$	
	When $f$ is also strictly co-unital, $f_\ren^?$ is strictly unital.
\end{itemize}

If $f_{ \Ran_{\on{un,dR}}}: Y_{ \Ran_{\on{un,dR}}}\to Z_{ \Ran_{\on{un,dR}}}$ is a morphism between corr-unital factorization spaces, then the above functors are factorization functors.
\end{construction}

\begin{remark} Construction \ref{constr-functorial-unital-fact-indcoh} breaks the symmetry between the $\IndCoh^{!,\ren}$ and $\IndCoh_{*,\ren}$ theory because in this paper we do not consider \emph{co-unital} factorization categories. Namely, as the dual of the unital factorization category $ \mathcal{I}nd\mathcal{C}oh_{*,\ren}(Y)$, the factorization category $ \mathcal{I}nd\mathcal{C}oh^{!,\ren}(Y)$ should be equipped with a co-unital structure. However, since we are more accustomed to unital structures\footnote{We prefer algebras than co-algebras.}, we produced a unital structure out of the co-unital one by replacing the \emph{dual functors} by \emph{conjugate functors}. It is this ad-hoc procedure that breaks the symmetry between $\IndCoh^{!,\ren}$ and $\IndCoh_{*,\ren}$.
\end{remark}

\begin{remark}  \label{rem-duality-unital-factorization-category}
Let $[\on{CrysCat}(\Ran_{\on{un,dR}})]_{\on{lax-unital}}$ be the category of crystal of categories over $\Ran_{\on{un,dR}}$ and lax-unital factorization functors. Here \emph{lax-unital} means for a morphism $\mathcal{C} \to \mathcal{D}$, we only have a \emph{lax} morphism between crystals of categories over the \emph{lax-prestack} $\on{Ran}_{\on{un,dR}}$. Explicitly, this means for a morphism $\mathfrak{t}:f\to g$ in $\Ran_{\on{un}}(T)$, we only have a $\DMod(T)$-linear natural transformation 
$$
\xymatrix{
		 \mathbf{\Gamma}( f_\dR,  \mathcal{C} ) \ar[r] \ar[d]_-{\mathfrak{t}_!(\mathcal{C})} &
		  \mathbf{\Gamma}( f_\dR, \mathcal{D} ) \ar[d]^-{\mathfrak{t}_!(\mathcal{D})} \ar@{=>}[dl] \\
		   \mathbf{\Gamma}( g_\dR,  \mathcal{C} ) \ar[r]   &
		  \mathbf{\Gamma}( g_\dR, \mathcal{D} ) 
}
$$
which is \emph{not} required to be invertible.

The category $[\on{CrysCat}(\Ran_{\on{un,dR}})]_{\on{lax-unital}}$ has an obvious symmetric monoidal structure induced by the tensor products of crystals of categories. Recall a duality beteen $\mathcal{M}$ and $\mathcal{N}$ inside this symmetric monoidal category consists of morphisms $u: \mathbf{1} \to \mathcal{M}\otimes \mathcal{N}$ and $c:\mathcal{N}\otimes \mathcal{M} \to \mathbf{1}$ satisfying the well-known conditions, where $\mathbf{1}$ is the monoidal unit in $[\on{CrysCat}(\Ran_{\on{un,dR}})]_{\on{lax-unital}}$ whose sections are $ \mathbf{\Gamma}( f_\dR,  \mathbf{1} ) \simeq \DMod(T)$. Unwinding the definitions, for any $\mathfrak{t}:f\to g$ as above, we have two $\DMod(T)$-linear natural transformation 
$$
\xymatrix{
		 \DMod(T) \ar[r]^-{u_{f_\dR}} \ar[d]^-= &
		 \mathbf{\Gamma}( f_\dR,  \mathcal{M} ) \otimes_{ \DMod(T) } \mathbf{\Gamma}( f_\dR,  \mathcal{N} ) \ar[d]^-{ \mathfrak{t}_!(\mathcal{M})\otimes \mathfrak{t}_!(\mathcal{N})  }  \ar@{=>}[dl] &
		  \mathbf{\Gamma}( f_\dR,  \mathcal{N} ) \otimes_{ \DMod(T) } \mathbf{\Gamma}( f_\dR,  \mathcal{M} ) \ar[r]^-{c_{f_\dR}} \ar[d]_-{\mathfrak{t}_!(\mathcal{N})\otimes \mathfrak{t}_!(\mathcal{M})} &
		 \DMod(T) 
		 \ar[d]^-=  \ar@{=>}[dl] \\
		    \DMod(T) \ar[r]_-{u_{g_\dR}}  &
		  \mathbf{\Gamma}( g_\dR,  \mathcal{M} ) \otimes_{ \DMod(T) } \mathbf{\Gamma}( g_\dR,  \mathcal{N} ) &
		   \mathbf{\Gamma}( g_\dR,  \mathcal{M} ) \otimes_{ \DMod(T) } \mathbf{\Gamma}( g_\dR,  \mathcal{N} ) \ar[r]_-{u_{g_\dR}}  &
		   \DMod(T)
} 
$$
satisfying certain conditions. Part of these conditions guarantee that $(u_{f_\dR},c_{g_\dR})$ exhibits a duality between $\mathbf{\Gamma}( f_\dR,  \mathcal{M} ) $ and $\mathbf{\Gamma}( f_\dR,  \mathcal{N} )$ as $\DMod(T)$-module categories. In other words, we have
\begin{equation}\label{eqn-rem-duality-unital-factorization-category-1}
\mathbb{D}_{f_\dR}: \mathbf{\Gamma}( f_\dR,  \mathcal{N} )^\vee \simeq \mathbf{\Gamma}( f_\dR,  \mathcal{M} ),\; \mathbb{D}_{g_\dR}: \mathbf{\Gamma}( g_\dR,  \mathcal{N} )^\vee \simeq \mathbf{\Gamma}( g_\dR,  \mathcal{M} ).
\end{equation}
Now knowing the above two natural transformations are equivalent to knowing
$$
\xymatrix{
		 \mathbf{\Gamma}( f_\dR,  \mathcal{N} )^\vee \ar[rr]^-{ \mathbb{D}_{f_\dR} }  & 
		 \ar@{=>}[d]&
		 \mathbf{\Gamma}( f_\dR,  \mathcal{M} )   \ar[d]^-{ \mathfrak{t}_!(\mathcal{M}) }   &
		  \mathbf{\Gamma}( f_\dR,  \mathcal{M} ) \ar[rr]^-{ \mathbb{D}_{f_\dR}^{-1} } \ar[d]_-{\mathfrak{t}_!(\mathcal{M})}  & \ar@{=>}[d] &
		 \mathbf{\Gamma}( f_\dR,  \mathcal{N} )^\vee 
		   \\
		    \mathbf{\Gamma}( g_\dR,  \mathcal{N} )^\vee \ar[u]^-{  \mathfrak{t}_!(\mathcal{N})^\vee } \ar[rr]_-{ \mathbb{D}_{g_\dR} }  & &
		  \mathbf{\Gamma}( g_\dR,  \mathcal{M} )  &
		   \mathbf{\Gamma}( g_\dR,  \mathcal{M} )  \ar[rr]_-{ \mathbb{D}_{g_\dR}^{-1} }  & &
		  \mathbf{\Gamma}( g_\dR,  \mathcal{N} )^\vee \ar[u]_-{\mathfrak{t}_!(\mathcal{N})^\vee}.
} 
$$
Now the duality conditions for $(u,c)$ can be translated exactly to the following condition: the above two natural transformations are the adjunction natural transformations for an adjoint pair 
$$  (  \mathfrak{t}_!(\mathcal{M}),   \mathbb{D}_{f_\dR} \circ \mathfrak{t}_!(\mathcal{N})^\vee \circ \mathbb{D}_{g_\dR}^{-1}).$$
In other words, they exhibit $ \mathfrak{t}_!(\mathcal{M})$ as the conjugate functor of $ \mathfrak{t}_!(\mathcal{N})$ via the dualities (\ref{eqn-rem-duality-unital-factorization-category-1}). This explains why $\mathcal{I}nd\mathcal{C}oh_{*,\ren}(Y)$ and $ \mathcal{I}nd\mathcal{C}oh^{!,\ren}(Y) $ are dual to each other in $[\on{CrysCat}(\Ran_{\on{un,dR}})]_{\on{lax-unital}}$.
\end{remark}

\section{Appendix D: \texorpdfstring{$\mathbb{E}_n$}{En}-Hochschild Cohomology and \texorpdfstring{$\mathbb{E}_n$}{En}-modules}

\label{sect:en}

\begin{notation}
In this section, we fix the following notations:
\begin{itemize}
	\item $n$ is a positive integer;
	\item $\mathcal{C}$ is a DG $\mathbb{E}_n$-category, i.e., an $\mathbb{E}_n$-algebra object in $\on{DGCat}$.
	\item $A\in \Enalg(\mathcal{C})$ is an $\mathbb{E}_n$-aglebra in $\mathcal{C}$.
\end{itemize}
\end{notation}

The goal of this section is to prove the following generalization of \cite[Proposition 4.36]{francistangent}:

\begin{theorem} \label{thm-HH-E_n}
Let $\tx{HH}_{\mathbb{E}_n,\,/\mathcal{C}} ( A\lmod(\mathcal{C}) ) $ be the $\mathbb{E}_n$-Hochschild cohomology of $A\lmod(\mathcal{C})$ relative to to $\mathcal{C}$. There is a canonical equivalence between DG $\mathbb{E}_n$-categories:
$$ \tx{HH}_{\mathbb{E}_n,\,/\mathcal{C}} ( A\lmod(\mathcal{C}) ) \simeq A\Enmod(\mathcal{C}) $$
compatible with the forgetful functors to $A\lmod(\mathcal{C})$.
\end{theorem}

We will explain the definitions of the relative $\mathbb{E}_n$-Hochschild cohomology later. However, let us point out when $n=2$, it is the \emph{relative Drinfeld center} in \cite[Section 4.4.1]{gaitsgory2021conjectural} (see also Definition \ref{def:rel-drinfeld-center}). Therefore we have:

\begin{corollary} Let $\on{Rep}_q(T)$ be the DG $\mathbb{E}_2$-category of representations of the quantum torus (see Definition \ref{def:quantum-torus}), and $U$ be a bialgebra object in $\on{Rep}_q(T)$ satisfying the condition in Lemma \ref{lemma:KD-vs-coBar}. Let $\Omega_U := \tx{coBar}(U)$ be the corresponding $\mathbb{E}_2$-algebra in $\on{Rep}_q(T)$. Then we have a canonical commutative diagram
$$
\xymatrix{
	U\lmod( \on{Rep}_q(T) )_{\on{loc.nilp.}} \ar[r]^-{\on{coBar}}_-\simeq &
	\Omega_U\lmod( \on{Rep}_q(T) ) \\
 	\on{Z}_{\on{Dr},/ \on{Rep}_q(T)} (U\lmod( \on{Rep}_q(T) )_{\on{loc.nilp.}}) \ar[u]^-{\on{oblv}} \ar[r]_-\simeq &
	\Omega_U{ -\textup{mod}^{\mathbb{E}_{2}} }( \on{Rep}_q(T) ). \ar[u]^-{\on{oblv}}
}
$$
\end{corollary}

\begin{remark}

To motivate the statement of the thoerem, let us consider the baby case where $n = 2$, $\mc{C} \simeq \tx{Vect}$, and suppose $A$ admits an $\mb{E}_\infty$-refinement $\tilde{A}$. Let us define $X := \tx{Spec}(\tilde{A})$, then the Drinfeld center becomes the category of $\QCoh(X \times X)$-linear endomorphisms of $\QCoh(X)$, which, by \cite{ben2010integral}, can be identified with $\QCoh(X \times_{X \times X} X)$, i.e. quasi-coherent sheaves on the derived loop space $\mc{L} X := \tx{Map}(S^1, X)$, by acting as kernel via the diagram
	\[\xymatrix{
		\mc{L} X \ar[r] \ar[d] & X \ar[d] \\
		X \ar[r] & X \times X
	}\]
Write $\QCoh^{\mb{E}_2}(X) := 
\tx{Mod}_{A}^{\mb{E}_2}(\tx{Vect})$, so our statement becomes $\QCoh^{\mb{E}_2}(X) \simeq \QCoh(\mc{L} X)$. More generally, the category $\tx{QCoh}^{\mb{E}_n}(X)$ can be identified with quasi-coherent sheaves on the $n$-fold based loop space. This provides a geometric interpretation of the notion of $\mb{E}_n$-modules.
\end{remark}

\subsection{Relative \texorpdfstring{$\mathbb{E}_n$}{En}-Hochschild Cohomology}

In this subsection, we define the relative \texorpdfstring{$\mathbb{E}_n$}{En}-Hochschild cohomology and state a precise version of Theorem \ref{thm-HH-E_n}.

\begin{construction} Consider the (very large) category $\mathcal{C}\lmod(\tx{DGCat})$ of left $\mathcal{C}$-modules in $\tx{DGCat}$. By \cite[5.1.4]{HA}, it has an $\mathbb{E}_{n-1}$-structure.
\end{construction}

\begin{example} When $n=1$, $\mathcal{C}\lmod(\tx{DGCat})$ is a $\mathbb{E}_{0}$-category with the distinguished object given by $\mathcal{C}$, considered as a left module of itself.

\end{example}

\begin{example} When $n=2$, $\mathcal{C}\lmod(\tx{DGCat})$ is a monoidal category with the multiplication given by
$$ \mathcal{C}\lmod(\tx{DGCat}) \times \mathcal{C}\lmod(\tx{DGCat}) \to \mathcal{C}\lmod(\tx{DGCat}),\; (\mathcal{M},\mathcal{N})\mapsto \mathcal{C}\otimes_{(\mathcal{C}\otimes \mathcal{C})} ( \mathcal{M}\otimes \mathcal{N}) ,$$
where the relative tensor product is formed by using the multiplication functor $m: \mathcal{C}\otimes \mathcal{C} \to \mathcal{C}$. Note that $m$ is a monoidal functor because $\mathcal{C}$ is assumed to have an $\mathbb{E}_2$-structure.

\end{example}

\begin{construction}
Consider the category $A\lmod(\mathcal{C})$ of left $A$-modules in $\mathcal{C}$. Note that it has a natural right $\mathcal{C}$-module structure. Similarly to \cite[4.8.5.20]{HA}, we will show
\begin{itemize}
	\item $A\lmod(\mathcal{C})\in \mathcal{C}\rmod(\tx{DGCat})$ has a canonical $\mathbb{E}_{n-1}$-algebra structure.
\end{itemize}

To define this structure, we consider the (very large) category $\DGCatMarked$ of \emph{DG categories with marked objects}. Its objects are pairs $(\mathcal{D},B)$ where $\mathcal{D}$ is a cocomplete DG-category and $B$ is an object in it; the morphisms $(\mathcal{D}_1,B_1)\to (\mathcal{D}_2,B_2)$ are given by continuous functors $F:\mathcal{D}_1\to \mathcal{D}_2$ together with morphisms $b:F(B_1)\to B_2$. A morphism as above in $\DGCatMarked$ is \emph{strict} if $b$ is an isomorphism, i.e., if it preserves the marked objects.

It is clear that $\DGCatMarked$ carries a natural symmetric monoidal structure given by 
$$(\mathcal{D}_1,B_1)\otimes (\mathcal{D}_2,B_2) := (\mathcal{D}_1\otimes \mathcal{D}_2,B_1\boxtimes B_2).$$
By definition, $(\mathcal{C},A)$ is an $\mathbb{E}_n$-algebra in $\DGCatMarked$. 

On the other hand, we consider the (very large) category $\alglmodpair(\tx{DGCat})$ of \emph{DG monoidal categories with right modules}. Its objects are pairs $(\mathcal{O},\mathcal{M})$ where $\mathcal{O}$ is an $\mathbb{E}_1$-algebra in $\tx{DGCat}$ and $\mathcal{M}$ is a right $\mathcal{O}$-module; the morphisms $( \mathcal{O}_1,\mathcal{M}_1 )\to ( \mathcal{O}_2,\mathcal{M}_2 ) $ are given by continuous $\mathbb{E}_1$-functors $F:\mathcal{O}_1\to \mathcal{O}_2$ and continuous functors $G: \mathcal{M}_1\to \mathcal{M}_2$ such that $F,G$ are compatible with the module structures in the usual sense. It is clear that $\alglmodpair(\tx{DGCat})$ inherits a symmetric monoidal structure from $\tx{DGCat}$.

Consider the functor
$$\begin{aligned}\Phi: \Eonealg( \DGCatMarked )  \to \alglmodpair(\tx{DGCat}),\; (\mathcal{D},B)\mapsto ( \mathcal{D},B\lmod(\mathcal{D}) ),\\
( (F,b): (\mathcal{D}_1,B_1)\to (\mathcal{D}_2,B_2)  ) \mapsto ( (F,\mathbf{ind}_b\circ F): ( \mathcal{D}_1,B_1\lmod(\mathcal{D}_1) )\to ( \mathcal{D}_2,B_2\lmod(\mathcal{D}_2) ) ),
\end{aligned}
$$
where the functor $ \mathbf{ind}_b\circ F $ is the composition
$$ B_1\lmod(\mathcal{D}_1) \overset{F}{\longrightarrow} F(B_1)\lmod(\mathcal{D}_2) \overset{\mathbf{ind}_b}{\longrightarrow} B_2\lmod(\mathcal{D}_2). $$
It is clear that $\Phi$ has a natural symmetric monoidal structure.

By enlarging the universe, we have a notion of the (very very large) category $\VVLCat$ of very large categories, such that $\DGCat$ is an object in it. Consider the (very very large) category $\VVLCatMarked$ defined similarly as $\DGCatMarked$. We have a functor
$$
\Psi: \alglmodpair(\tx{DGCat}) \to \VVLCatMarked ,\; (\mathcal{O,\mathcal{M}})\mapsto ( \mathcal{O}\rmod(\tx{DGCat}), \mathcal{M} ).
$$
It also has a natural \emph{right-lax} symmetric monoidal structure. Combining the above two paragraphs, we obtain a \emph{right-lax} symmetric monoidal functor $\Psi\circ  \Phi$. Hence it induces a functor $\Enalg( \DGCatMarked )  \to \Enminusonealg (\VVLCatMarked )$.
By construction, it sends $(\mathcal C,A)$ to 
$$(\mathcal{C}\rmod(\tx{DGCat}), A\lmod(\mathcal{C}) ) \in   \Enminusonealg (\VVLCatMarked ).$$
Unwinding definitions, this means $A\lmod(\mathcal{C})\in \mathcal{C}\rmod(\tx{DGCat})$ has a canonical $\mathbb{E}_{n-1}$-algebra structure as desired.
\end{construction}

\begin{example} When $n=1$, we obtain an $\mathbb{E}_0$-structure on $A\lmod(\mathcal{C})\in \mathcal{C}\rmod(\tx{DGCat})$, i.e., a right $\mathcal{C}$-linear functors $\mathcal{C}\to A\lmod(\mathcal{C})$. By construction, it is the functor that sends $X \in \mathcal{C}$ to $m(A, X)$ equipped with the obvious left $A$-module structure.
\end{example}

\begin{example} When $n=2$, we obtain an $\mathbb{E}_1$-structure on $A\lmod(\mathcal{C})\in \mathcal{C}\rmod(\tx{DGCat})$. By construction, the multiplication is given by
$$ (A\lmod(\mathcal{C}) \otimes A\lmod(\mathcal{C})) \otimes_{ (\mathcal{C}\otimes \mathcal{C})  }\mathcal{C} \simeq (A\boxtimes A)\lmod(\mathcal{C}\otimes \mathcal{C}  ) \otimes_{(\mathcal{C}\otimes \mathcal{C})} \mathcal{C} \simeq m(A, A)\lmod (\mathcal{C}) \overset{ \on{ind} }{\longrightarrow} A\lmod (\mathcal{C}),$$
where the last functor is induced by the $\mathbb{E}_1$-map $m(A, A) \to A$.
\end{example}

\begin{definition} \label{def:En-HH-definition}
Rougly speaking, for $\mathcal{D} \in \tx{Alg}^{\mathbb{E}_{n-1}}(\mathcal{C}\rmod(\tx{DGCat}))$, the \emph{$\mathbb{E}_n$-Hochschild cohomology of $\mathcal{D}$ relative to $\mathcal{C}$} is the universal DG $\mathbb{E}_n$-category, viewed as an algebra in $ \tx{Alg}^{\mathbb{E}_{n-1}}( \on{DGCat} )$, that acts on $\mathcal{D}\in \tx{Alg}^{\mathbb{E}_{n-1}}( \on{DGCat} )$ in a way compatible with the right action of $\mathcal{C}$ on $\mathcal{D}$.

More precisly,
\[\tx{HH}_{\mathbb{E}_n,\,/\mathcal{C}}( \mathcal{D} ) := \mathbf{Maps}_{ \mathcal{D}-\mathbb{E}_{n-1} } (  \mathcal{D},\mathcal{D} ),\]
where the RHS is the mapping DG category inside the $\on{DGCat}$-enriched category $\mathcal{D}\Enminusonemod(\mathcal{C}\rmod(\tx{DGCat}))$.
As in \cite[5.3.1.15]{HA}, $\tx{HH}_{\mathbb{E}_n,\,/\mathcal{C}}( \mathcal{D} )$ carries the structure of a DG $\mathbb{E}_n$-category.
\end{definition}

\begin{example} We have $\tx{HH}_{\mathbb{E}_n,\,/\mathcal{C}}( \mathcal{C} ) \simeq \mathcal{C}$ as DG $\mathbb{E}_n$-categories.

\end{example}

\begin{construction} It follows from definition that we have an action of $\tx{HH}_{\mathbb{E}_n,\,/\mathcal{C}}( \mathcal{D} )$ on $\mathcal{D}$ viewed as objects in $ \tx{Alg}^{\mathbb{E}_{n-1}}( \on{DGCat} )$. In particular, we have a forgetful functor
$$ \on{oblv}: \tx{HH}_{\mathbb{E}_n,\,/\mathcal{C}}( \mathcal{D} )\to \mathcal{D} $$
between DG $\mathbb{E}_{n-1}$-categories. Explicitly, it is defined as the composition
$$ \mathbf{Maps}_{\mathcal{D}-\mathbb{E}_{n-1} } (  \mathcal{D},\mathcal{D} ) \to \on{LFun}(  \mathcal{D},\mathcal{D} ) \to \mathcal{D},$$
where the last functor is taking value at the monoidal unit of $\mathcal{D}$.
\end{construction}

To describe the functoriality of $\tx{HH}_{\mathbb{E}_n,\,/\mathcal{C}}( \mathcal{D} )$ in $\mathcal{D}$, we need the following definition:

\begin{definition}

Let $F: \mathcal{D}_1\to \mathcal{D}_2$ be a morphism in $\tx{Alg}^{\mathbb{E}_{n-1}}(\mathcal{C}\rmod(\tx{DGCat}))$. We say $F$ \emph{satisfies the projection formula} if the right adjoint $F^R: \mathcal{D}_2\to \mathcal{D}_1$, which is a priori right-lax $\mathcal{D}_1$-linear, is strictly $\mathcal{D}_1$-linear. In particular, $F^R$ is also strictly $\mathcal{C}$-linear and continuous.
\end{definition}

\begin{remark} The above definition also makes sense when $n=\infty$. In this case, the pullback functor $f^*: \QCoh(Y)\to \QCoh(X)$ for quasi-coherent sheaves on schemes satisfies the projection formula in the above sense.
\end{remark}

\begin{example} One can verify that the induction functor $\on{ind}: \mathcal{C}\to A\lmod(\mathcal{C})$ satisfying the projection formula.
\end{example}

\begin{construction} Let $F: \mathcal{D}_1\to \mathcal{D}_2$ be a morphism in $\tx{Alg}^{\mathbb{E}_{n-1}}(\mathcal{C}\rmod(\tx{DGCat}))$ satisfying the projection formula. Then we have a functor
$$ \tx{HH}_{\mathbb{E}_n,\,/\mathcal{C}}( \mathcal{D}_2 ) \to \tx{HH}_{\mathbb{E}_n,\,/\mathcal{C}}( \mathcal{D}_1 ) $$
defined as
$$  \mathbf{Maps}_{ \mathcal{D}_2-\mathbb{E}_{n-1}} (  \mathcal{D}_2,\mathcal{D}_2 ) \to \mathbf{Maps}_{ \mathcal{D}_1-\mathbb{E}_{n-1} } (  \mathcal{D}_1,\mathcal{D}_1 ),\; T\mapsto F^R \circ T \circ F.$$
It follows formally that the above functor has a canonical \emph{right-lax} $\mathbb{E}_n$-structure\footnote{The details will be provided in the second author's thesis. Also, we do not need this claim to prove the main result of this paper. Namely, without it, we can still prove a weaker version of Theorem \ref{thm-HH-E_n} which says there is a canonical equivalence $\tx{HH}_{\mathbb{E}_n,\,/\mathcal{C}} ( A\lmod(\mathcal{C}) ) \simeq A\Enmod(\mathcal{C})$ between plain categories. This weaker result is enough for our purpose.
}.
\end{construction}

\begin{example} \label{exam-E_1-strucure-functorial-HH}
Explicitly, the right-lax $\mathbb{E}_1$-structure on $\tx{HH}_{\mathbb{E}_n,\,/\mathcal{C}}( \mathcal{D}_2 ) \to \tx{HH}_{\mathbb{E}_n,\,/\mathcal{C}}( \mathcal{D}_1 )$ is given by the natural transformation
$$ F^R \circ T \circ F \circ F^R \circ T' \circ F \to F^R \circ T\circ T' \circ F $$
induced by the counit adjunction $F\circ F^R \to \on{Id}$.
\end{example}

\begin{construction} We apply the above construction to the induction functor $\on{ind}: \mathcal{C}\to A\lmod(\mathcal{C})$ and obtain a right-lax $\mathbb{E}_n$-functor
$$  \tx{HH}_{\mathbb{E}_n,\,/\mathcal{C}}( A\lmod(\mathcal{C}) ) \to \tx{HH}_{\mathbb{E}_n,\,/\mathcal{C}}( \mathcal{C} )  \simeq \mathcal{C}. $$
It follows from constructinon that the unit of LHS is sent to the $\mathbb{E}_n$-algebra $A$ in RHS. Hence the above functor factors through a 
\end{construction}

\begin{lemma} The right-lax $\mathbb{E}_n$-functor $\tx{HH}_{\mathbb{E}_n,\,/\mathcal{C}}( A\lmod(\mathcal{C}) ) \to A\Enmod (\mathcal{C})$ is strict.
\end{lemma}

\proof We only need to prove the underlying right-lax $\mathbb{E}_1$-functor is strict. However, this follows from the explicit description in Example \ref{exam-E_1-strucure-functorial-HH}.

\qed

The rest of this section is to prove the following precise version of Theorem \ref{thm-HH-E_n}:

\begin{proposition}
\label{prop:HH-is-E_n-mod}
The $\mathbb{E}_n$-functor $\tx{HH}_{\mathbb{E}_n,\,/\mathcal{C}}( A\lmod(\mathcal{C}) ) \to A\Enmod (\mathcal{C})$ constructed above is an equivalence.
\end{proposition}

\begin{remark} When $\mathcal{C}$ can be promoted to a symmetric monoidal category, the proposition was proved in \cite[Proposition 4.36]{francis2013tangent}. The idea is to express both sides as $U_{\mathbb{E}_n}(A)\lmod(\mathcal{C})$, where $U_{\mathbb{E}_n}(A)$ is the enveloping algebra of the $\mathbb{E}_n$-aglebra $A$. Unfortunately, the construction of $U_{\mathbb{E}_n}(A)$ only makes sense when $\mathcal{C}$ is symmetric monoidal. Indeed, in the general case, $A\Enmod (\mathcal{C})$ does not have a $\mathcal{C}$-action.

Nevertheless, we will prove the proposition based on the following observation: there exists a \emph{relative} enveloping algrbra $U_{\mathbb{E}_n}^{\on{rel}}(A)$ of $A$ inside $U_{\mathbb{E}_n}(\mathcal{C})$ such that we have a canonical equivalence
$$ A\Enmod (\mathcal{C})\simeq U_{\mathbb{E}_n}^{\on{rel}}(A)\lmod( \mathcal{C} ), $$
where the definition of the RHS uses the left action of $U_{\mathbb{E}_n}(\mathcal{C})$ on $\mathcal{C}$ that corresponds to the $\mathbb{E}_n$-action on $\mathcal{C}$ on itself. Then one can generalize the proof of \cite[Proposition 4.36]{francis2013tangent} by carefully replace the enveloping algrbras by the relative ones.

We will construct and study $U_{\mathbb{E}_n}^{\on{rel}}(A)$ in $\S$ \ref{ssec-enveloping-algebras} and use it to prove the proposition in $\S$ \ref{ssec-proof-E_n-HH}.
\end{remark}

\subsection{Enveloping Algebras}
\label{ssec-enveloping-algebras}

\begin{construction} Recall the (very large) symmetric monoidal category $\on{DGCat}_*$ of DG categories with a marked object in it. Note that the tensor product in $\on{DGCat}_*$ distributes with sifted colimits (in fact with all small colimits). Hence for $(\mathcal{D}, B)\in \tx{Alg}^{\mathbb{E}_{n}}(\on{DGCat}_*)$, we have an adjoint pair
$$ \on{ind}_{(\mathcal{D}, B)} :\on{DGCat}_* \adjoint (\mathcal{D}, B)\Enmod(\on{DGCat}_*): \on{oblv}_{(\mathcal{D}, B)}$$
such that both functors are $\on{DGCat}_*$-linear. Hence the corresponding monad is given by tensoring with an associative algebra object $U_{\mathbb{E}_n}(\mathcal{D},B)$ in $ \on{DGCat}_*$, which is known as the \emph{enveloping algebra} of $(\mathcal{D},B)$. We have
\begin{equation} \label{eqn-universal-property-of-enveloping-marked}
  (\mathcal{D},B)\tx{-mod}^{\mb{E}_{n}}(\on{DGCat}_*)\simeq U_{\mathbb{E}_n}(\mathcal{D},B)\lmod( \on{DGCat}_* ) 
  \end{equation}
compatible with the forgetful functors to $\on{DGCat}_*$. Also, the functor
$$ U_{\mathbb{E}_n }: \tx{Alg}^{\mathbb{E}_{n}}(\on{DGCat}_* ) \to \tx{Alg}(\on{DGCat}_* ) $$
is symmetric monoidal.

Since $\on{oblv}_*: \on{DGCat}_* \to \on{DGCat}$ is symmetric monoidal, the algebra $\on{oblv}_* (U_{\mathbb{E}_n}(\mathcal{D},B))$ is just the enveloping algebra $U_{\mathbb{E}_n}(\mathcal{D})$ of $\mathcal{D}$ in $\on{DGCat}$. We have
\begin{equation} \label{eqn-universal-property-of-enveloping}
  \mathcal{D}\tx{-mod}^{\mb{E}_{n}}(\on{DGCat})\simeq U_{\mathbb{E}_n}(\mathcal{D})\lmod( \on{DGCat} ).
  \end{equation}

Recall that knowing $U_{\mathbb{E}_n}(\mathcal{D},B)$ is equivalent to knowing an associative algebra object inside $U_{\mathbb{E}_n}(\mathcal{D})$. We refer this algebra as the \emph{relative enveloping algebra} of $B$ and denote it by $U_{\mathbb{E}_n}^{\on{rel}}(B) \in U_{\mathbb{E}_n}(\mathcal{D})$.

Let $\mathcal{M}$ be any object in $ \mathcal{D}\tx{-mod}^{\mb{E}_{n}}(\on{DGCat})$. e can use (\ref{eqn-universal-property-of-enveloping}) to view it as a left module of $U_{\mathbb{E}_n}(\mathcal{D})$. By (\ref{eqn-universal-property-of-enveloping-marked}), we have
$$ 
(\mathcal{D},B)\tx{-mod}^{\mb{E}_{n}}(\on{DGCat}_*) \times_{  \mathcal{D}\tx{-mod}^{\mb{E}_{n}}(\on{DGCat}) } \{ \mathcal{M} \} \simeq U_{\mathbb{E}_n}(\mathcal{D},B)\lmod( \on{DGCat}_* ) \times_{  U_{\mathbb{E}_n}(\mathcal{D})\lmod( \on{DGCat} ) }  \{ \mathcal{M} \}.
$$
In other words,
$$  B\Enmod(\mathcal{M}) \simeq U_{\mathbb{E}_n}^{\on{rel}}(B)\lmod(\mathcal{M}),$$
where the LHS is defined via the $\mb{E}_{n}$-$\mathcal{D}$-module on $\mathcal{M}$ while the RHS is defined via the left $U_{\mathbb{E}_n}(\mathcal{D})$-module structure on $\mathcal{M}$.

\end{construction}

\begin{example} When $n=1$, we have $U_{\mb{E}_{1 }}(\mathcal{D}) \simeq \mathcal{D}\otimes \mathcal{D}^{\on{rev}}$ and $U_{\mb{E}_{1 }}^{\on{rel}}(B) \simeq B\boxtimes B^{\on{rev}}$. Note that $B^{\on{rev}}$ only makes sense inside the monoidal category $\mathcal{D}^{\on{rev}}$.
\end{example}

\begin{remark} Argument from \cite[3.16]{francis2013tangent} carries through and allows us to compute it as the factorization homology
\[U_{\mb{E}_{n }}(\mathcal{D},B) \simeq \int_{S^{ n - 1} \times \mb{R}} (\mathcal{D},B)\]
where on RHS we consider $(\mathcal{D},B)$ as a factorizable cosheaf valued in $\VVLCatMarked$ over $\mb{R}^{n }$. In partuclar, we have
$$
U_{\mb{E}_{n }}(\mathcal{D},B)  \simeq (\mathcal{D},B)^{\on{rev}}\otimes_{U_{\mb{E}_{n-1 }}(\mathcal{D},B)^{\on{rev}} } (\mathcal{D},B).
$$
Hence we have
$$ U_{\mb{E}_{n }}(\mathcal{D}) \simeq \mathcal{D}^{\on{rev}}\otimes_{ U_{\mb{E}_{n-1 }}(\mathcal{D})^{\on{rev}}} \mathcal{D},\; U_{\mathbb{E}_n}^{\on{rel}}(B)\simeq B^{\on{rev}}\boxtimes_{ U_{\mathbb{E}_{n-1}}^{\on{rel}}(B)^{\on{rev}} } B^{\on{rev}}. $$

\end{remark}

\begin{construction} The functor
$$\on{Mod}:\on{Alg}( \DGCatMarked )  \to \DGCat,\; (\mathcal{D},B)\mapsto B\lmod(\mathcal{D}) 
$$
is symmetric monoidal and commutes with sifted colimits. It follows that we have a commutative diagram
$$
\xymatrix{
	\on{Alg}^{\mathbb{E}_{n-1}}( \on{Alg}( \DGCatMarked )  ) \ar[r]^-{\on{Mod}}  \ar[d]^-{U_{\mathbb{E}_{n-1}}}&
	\on{Alg}^{\mathbb{E}_{n-1}}(\DGCat)  \ar[d]^-{U_{\mathbb{E}_{n-1}}} \\
	\on{Alg}( \on{Alg}( \DGCatMarked )  ) \ar[r]^-{\on{Mod}} &
	\on{Alg}(\DGCat).
}
$$
In other words, for $(\mathcal{C},A)\in \on{Alg}^{\mathbb{E}_{n}}( \DGCatMarked )$, we have
\begin{equation} \label{eqn-U-mod=mod-U}
 U_{\mathbb{E}_{n-1}}(  A\lmod(\mathcal{C} )) \simeq U_{\mathbb{E}_{n-1}}^{\on{rel}}(A)\lmod(  U_{\mathbb{E}_{n-1}}(\mathcal{C})) 
 \end{equation}
as DG monoidal categories.
\end{construction}

\begin{example} When $n=2$, the above equivalence reads as
$$ A\lmod(\mathcal{C} )\otimes A\lmod(\mathcal{C} )^{\on{rev}} \simeq (A\boxtimes A^{\on{rev}}) \lmod( \mathcal{C}\otimes \mathcal{C}^{\on{rev}} ). $$
Indeed, this follows from $A\lmod(\mathcal{C} )^{\on{rev}}\simeq A^{\on{rev}}\lmod(\mathcal{C} ^{\on{rev}})$.
\end{example}

\begin{construction} Similarly, we consider the (very large) symmetric monoidal category $\alglmodpair(\tx{DGCat})$ of DG monoidal categories with right modules. Note that the two forgerful functors
$$\alglmodpair(\tx{DGCat})  \to \on{DGCat} $$
are symmetric monoidal and commute with sifted colimits. Hence we can define
$$ U_{\mathbb{E}_{n-1}}: \on{Alg}^{\mathbb{E}_{n-1}}( \alglmodpair(\tx{DGCat})  ) \to \on{Alg}^{\mathbb{E}_{1}}( \alglmodpair(\tx{DGCat})  ) $$
such that
$$ U_{\mathbb{E}_{n-1}}( \mathcal{O},\mathcal{M} ) \lmod(  \alglmodpair(\tx{DGCat}) ) \simeq ( \mathcal{O},\mathcal{M} )\Enminusonemod (  \alglmodpair(\tx{DGCat}) ), $$
and
$$ U_{\mathbb{E}_{n-1}}( \mathcal{O},\mathcal{M} )\simeq (U_{\mathbb{E}_{n-1}}(\mathcal{O}) , U_{\mathbb{E}_{n-1}}(\mathcal{M}) ).$$

Note that we can view $\mathcal{O}$ as an object in
$$ \mathcal{O}\Enminusonemod ( \on{Alg}( \DGCat )  ) \simeq  U_{\mathbb{E}_{n-1}}( \mathcal{O}) \lmod ( \on{Alg}( \DGCat )  ).$$
Hence we have
\[
\begin{aligned}
	U_{\mathbb{E}_{n-1}}( \mathcal{O},\mathcal{M} ) \lmod(  \alglmodpair(\tx{DGCat}) )\times_{ U_{\mathbb{E}_{n-1}}( \mathcal{O} ) \lmod( \on{Alg}( \DGCat ))  } \{ \mathcal{O} \} \simeq \\
	\simeq ( \mathcal{O},\mathcal{M} )\Enminusonemod (  \alglmodpair(\tx{DGCat}) )\times_{  \mathcal{O}\Enminusonemod (  \on{Alg}(\tx{DGCat}) )  }\{ \mathcal{O} \}.
\end{aligned}
\]
Unwinding the definitions, this says
$$ U_{\mathbb{E}_{n-1}}( \mathcal{M} ) \lmod ( \mathcal{O}\rmod(\DGCat) ) \simeq \mathcal{M} \Enminusonemod (  \mathcal{O}\rmod(\DGCat) ),  $$
where $U_{\mathbb{E}_{n-1}}( \mathcal{M} ) \in \on{Alg}(U_{\mathbb{E}_{n-1}}( \mathcal{O} ) \rmod(\DGCat) )$ and the LHS is defined using the $U_{\mathbb{E}_{n-1}}( \mathcal{O} ) \rmod(\DGCat)$-action on $ \mathcal{O}\rmod(\DGCat)$.

\end{construction}

\begin{example} When $n=2$, the above equivalence reads as
$$ (\mathcal{M}\otimes \mathcal{M}^{\on{rev}})  \lmod ( \mathcal{O}\rmod(\DGCat) ) \simeq \mathcal{M} \on{-bimod} (  \mathcal{O}\rmod(\DGCat) ). $$
\end{example}

\begin{example}\label{exam-key-result-proof-E_n-HH}
 Applying the above construction to $( \mathcal{O},\mathcal{M} ):=(\mathcal{C}, A\mod(\mathcal{C}))$ and using the equivalence (\ref{eqn-U-mod=mod-U}), we obtain
$$ (U_{\mathbb{E}_{n-1}}^{\on{rel}}( A )\mod( U_{\mathbb{E}_{n-1}}(\mathcal{C} ) ) )\mod(  \mathcal{C}\rmod(\DGCat) ) \simeq  (A\lmod(\mathcal{C})) \Enminusonemod(  \mathcal{C}\rmod(\DGCat) ) .$$
We emphasize that $U_{\mathbb{E}_{n-1}}^{\on{rel}}( A )\mod( U_{\mathbb{E}_{n-1}}(\mathcal{C} ) )$ is an algebra in $U_{\mathbb{E}_{n-1}}( \mathcal{C} ) \rmod(\DGCat)$, which acts on $\mathcal{C}\rmod(\DGCat)$.
\end{example}

\subsection{Proof of Proposition \ref{prop:HH-is-E_n-mod}}
\label{ssec-proof-E_n-HH}

We only need to show the underlying functor $\tx{HH}_{\mathbb{E}_n,\,/\mathcal{C}}( A\lmod(\mathcal{C}) ) \to A\Enmod (\mathcal{C})$ is an equivalence. We will construct an equivalence
$$\tx{HH}_{\mathbb{E}_n,\,/\mathcal{C}}( A\lmod(\mathcal{C}) ) \to U_{\mathbb{E}_n}^{\on{rel}}(A)\lmod(\mathcal{C})$$
compatible with the forgetful functors to $\mathcal{C}$, and leave it to the reader to check the composition
$$\tx{HH}_{\mathbb{E}_n,\,/\mathcal{C}}( A\lmod(\mathcal{C}) ) \to  U_{\mathbb{E}_n}^{\on{rel}}(A)\lmod(\mathcal{C}) \simeq A\Enmod (\mathcal{C})$$
is indeed the functor constructed before\footnote{For the purpose of proving the main theorem of this paper, we only need \emph{an} equivalence $\tx{HH}_{\mathbb{E}_n,\,/\mathcal{C}}( A\lmod(\mathcal{C}) ) \to A\Enmod (\mathcal{C})$ between \emph{plain DG categories} compatible with the forgetful functors to $\mathcal{C}$. Hence the construction below is enough.}.

By Example \ref{exam-key-result-proof-E_n-HH}, we have
\begin{eqnarray*}
& & \tx{HH}_{\mathbb{E}_n,\,/\mathcal{C}}(  A\lmod(\mathcal{C}))  \\
&\simeq	&  \mathbf{Maps}_{ ( A\lmod(\mathcal{C}))-\mathbb{E}_{n-1} } (   A\lmod(\mathcal{C}), A\lmod(\mathcal{C}))  \\
&\simeq	&  \mathbf{Maps}_{   U_{\mathbb{E}_{n-1}}^{\on{rel}}( A )\mod( U_{\mathbb{E}_{n-1}}(\mathcal{C} ) )  } (   A\lmod(\mathcal{C}), A\lmod(\mathcal{C})) .
\end{eqnarray*}
Unwinding the definitions, the endo-functor on $\mathcal{C}\rmod(\DGCat)$ given by the algebra $ U_{\mathbb{E}_{n-1}}^{\on{rel}}( A )\mod( U_{\mathbb{E}_{n-1}}(\mathcal{C} ) )\in \on{Alg}(U_{\mathbb{E}_{n-1}}( \mathcal{C} ) \rmod(\DGCat) )$ is
$$\mathcal{M} \mapsto ( U_{\mathbb{E}_{n-1}}^{\on{rel}}( A )\mod( U_{\mathbb{E}_{n-1}}(\mathcal{C} ) ) \otimes \mathcal{M}) \otimes_{  U_{\mathbb{E}_{n-1}}(\mathcal{C} )\otimes \mathcal{C}  } \mathcal{C},$$
where $U_{\mathbb{E}_{n-1}}(\mathcal{C} )\otimes \mathcal{C} \to \mathcal{C}$ is the action of $U_{\mathbb{E}_{n-1}}(\mathcal{C} )$ on $\mathcal{C}$ viewed as objects in $\on{Alg}(\DGCat)$. It follows we have
\begin{eqnarray*}
& &  \mathbf{Maps}_{   U_{\mathbb{E}_{n-1}}^{\on{rel}}( A )\mod( U_{\mathbb{E}_{n-1}}(\mathcal{C} ) )  } (   A\lmod(\mathcal{C}), A\lmod(\mathcal{C})) \\
& \simeq &  \on{Tot} \mathbf{Maps}( (U_{\mathbb{E}_{n-1}}^{\on{rel}}( A )^{\boxtimes \bullet}\boxtimes A)\lmod(\mathcal{C}), A\lmod(\mathcal{C})),
\end{eqnarray*}
where $U_{\mathbb{E}_{n-1}}^{\on{rel}}( A )^{\boxtimes \bullet}\boxtimes A$ is an algebra in $U_{\mathbb{E}_{n-1}}(\mathcal{C} )^{\otimes \bullet} \otimes \mathcal{C}$ which acts on $\mathcal{C}$. Note that the above $\mathbf{Maps}(-,-)$ is the mapping DG category inside $\mathcal{C}\rmod(\DGCat)$. Hence we have
\begin{eqnarray*}
& &    \on{Tot} \mathbf{Maps}(( U_{\mathbb{E}_{n-1}}^{\on{rel}}( A )^{\boxtimes \bullet}\boxtimes A)\lmod(\mathcal{C}), A\lmod(\mathcal{C})) \\
&\simeq & \on{Tot} [(A\boxtimes  (U_{\mathbb{E}_{n-1}}^{\on{rel}}( A )^{\boxtimes \bullet}\boxtimes A)^{\on{rev}})\lmod(\mathcal{C})] \\
&\simeq & (\colim (A\boxtimes  (U_{\mathbb{E}_{n-1}}^{\on{rel}}( A )^{\boxtimes \bullet}\boxtimes A)^{\on{rev}})\lmod(\mathcal{C}) \\
&\simeq & A\boxtimes_{ U_{\mathbb{E}_{n-1}}^{\on{rel}}( A )^{\on{rev}} } A^{\on{rev}}\lmod(\mathcal{C}) \\
& \simeq & U_{\mathbb{E}_n}(A)\lmod(\mathcal{C})
\end{eqnarray*}
as desired.
\qed

\subsection{Proof of Lemma \ref{lemma:oblv-mxd-to-KD}}
The forgetful functor
\[U_{\mb{E}_2}(\Rep_q(T))\mod(\DGCat) \simeq \Rep_q(T)\mod^{\mb{E}_2}(\DGCat) \xrightarrow{\oblv} \Rep_q(T)\mod(\DGCat)\]
preserves the base point $\Rep_q(T)$, thus induces a monoidal functor $\tx{ins}: \Rep_q(T) \to U_{\mb{E}_2}(\Rep_q(T))$ (topologically, this corresponds to the inclusion $\mb{R} \setminus \{0\} \subseteq \mb{R}^2 \setminus \{0\}$).
Furthermore, we have an algebra
\[U_{\mb{E}_2}^\tx{rel}(\Omega_q) \in \tx{Alg}(U_{\mb{E}_2}(\Rep_q(T)))\]
such that $U_{\mb{E}_2}^\tx{rel}(\Omega_q)\mod(\Rep_q(T)) \simeq \Omega_q\mod^{\mb{E}_2}(\Rep_q(T))$. Looking again at $\mb{R} \setminus \{0\} \subseteq \mb{R}^2 \setminus \{0\}$ but this time considering the $\mb{E}_2$-algebra $(\Rep_q(T), \Omega_q)$ within $\mc{CAT}_*$, we see that there is a morphism $\tx{ins}(\Omega_q) \to U_{\mb{E}_2}^\tx{rel}(\Omega_q)$ of $\mb{E}_1$-algebras within $U_{\mb{E}_2}(\Rep_q(T))$.

The category $U_{\mb{E}_2}(\Rep_q(T))$ is a $(\Rep_q(T), U_{\mb{E}_2}(\Rep_q(T)))$-bimodule category, where $\Rep_q(T)$ acts via the functor $\tx{ins}$. We consider the object
\[\mf{K} := k \tensor_{\Omega_q} U_{\mb{E}_2}^\tx{rel}(\Omega_q) := |k \boxtimes \Omega_q^\bullet \boxtimes U_{\mb{E}_2}^\tx{rel}(\Omega_q)| \in \Rep_q(T) \tensor_{\Rep_q(T)} U_{\mb{E}_2}(\Rep_q(T)) \simeq U_{\mb{E}_2}(\Rep_q(T));\]
by above, this is a well-defined expression, and it is a bimodule over
\[\Hom_{\Omega_q\mod^r(\Rep_q(T))}(k, k) \simeq \Hom_{U_q^\tx{Lus}(N)\mod(\Rep_q(T))_{\tx{loc.nilp.}}}(\tx{Bar}_{\Omega_q}(k), \tx{Bar}_{\Omega_q}(k)) \simeq U_q^\tx{KD}(N^-)\]
as an algebra in $\Rep_q(T)$ acting from left, and $U_{\mb{E}_2}^\tx{rel}(\Omega_q)$ from the right. Thus, using it as the kernel, we obtain a functor
\[\Rep_q^\tx{mxd}(G) \simeq U_{\mb{E}_2}^\tx{rel}(\Omega_q)\mod(\Rep_q(T)) \xrightarrow{\mf{K} \tensor_{U_{\mb{E}_2}^\tx{rel}} (-)} U_q^\tx{KD}(N^-)\mod(\Rep_q(T));\]
it is straightforward to check that further forgetting to $\Rep_q(T)$ gives the expected answer.
\qed

\section{Appendix E: Riemann-Hilbert Correspondence}
\label{sect:rh}

For any prestack $\mc{Y}$:
\begin{itemize}
\item a de Rham $\mb{G}_m$-local system is a map from $\mc{Y}$ to $\mb{B}^1_{et}(\mb{G}_m)$, the \'{e}tale classifying stack for $\mb{G}_m$ local systems;
\item a de Rham $\mb{G}_m$-gerbe is a map from $\mc{Y}$ to $\mb{B}^2_{et}(\mb{G}_m)$;
\item an analytic gerbe is a map from $\mc{Y}^\tx{an}$, the analytification of $\mc{Y}$ (\ref{sect:rh}), to the topological space $K(\mb{C}^\times, 2)$;
\item a $\mb{G}_m$-twisting is an element of
\[\tx{Fib}[\tx{Maps}_{\tx{PreStk}}(\mc{Y}_\dR, \mb{B}^2_{et}(\mb{G}_m)) \to \tx{Maps}_{\tx{PreStk}}(\mc{Y}, \mb{B}^2_{et}(\mb{G}_m))].\]
\end{itemize}

For de Rham $\mb{G}_m$-gerbe, local system and twisting we also have the notion of tameness; we refer readers to \cite{zhao2020tame} for details. Following \emph{loc.cit.} we let
\begin{itemize}
\item $\tb{Ge}_\dR$ (resp.\ $\overset{\circ}{\tb{Ge}}$, $\tb{Ge}_\tx{an}$) denote the classical prestack of all de Rham $\mb{G}_m$-gerbes (resp.\ tame de Rham $\mb{G}_m$-gerbes, analytic $\mb{C}^\times$-gerbes). It is in fact an h-sheaf (resp. \'{e}h-sheaf, h-sheaf);
\item $\overset{\circ}{\tb{Loc}_1}$ denote the prestack of tame $\mb{G}_m$-local systems, and $\tb{Loc}_\tx{an}$ that of analytic $\mb{C}^\times$-local systems;
\item $\tb{Tw}$ (resp.\ $\overset{\circ}{\tb{Tw}}$) denote the prestack of all $\mb{G}_m$-twistings (resp. tame twistings); it is a derived h-sheaf (resp. derived \'{e}h-sheaf).
\end{itemize}
There are natural forgetful functors
\[\overset{\circ}{\tb{Tw}} \to \tb{Tw} \to \tb{Ge}_\dR \hspace{1em} \overset{\circ}{\tb{Tw}} \to \overset{\circ}{\tb{Ge}}\]
To avoid confusion, we let $\oblv(\deRhamGerbe)$ denote the de Rham $\mb{G}_m$-gerbe under the first forgetful functor.

One can further define, in a natural manner, the notion of a \emph{factorizable} tame twisting. The fundamental result here is (\cite[Theorem C]{zhao2020tame}):

\begin{proposition}
\label{prop:yifei-twisting-conversion}
For $G$ reductive, the $\infty$-groupoid (in fact a strict 2-groupoid) of factorizable tame twistings on $[\tx{Gr}_G]_{\Ran_\dR}$ is isomorphic to the $\mb{C}$-points of $\tx{Par}^\circ(G)_{\tx{tame}}$.
\end{proposition}

\subsection{Regular-Holonomic D-modules}
For $S$ any classical aft smooth scheme, we have (non-cocomplete) full subcategories
\[\DMod^{\tx{hol}, \tx{RS}}(S) \subseteq \DMod_\tx{coh}(S) \subseteq \DMod(S)\]
consists of objects whose cohomologies are regular holonomic D-modules (resp. coherent D-modules). Since $S$ is a aft scheme, we have $\DMod(S) \simeq \tx{Ind}(\DMod_\tx{coh}(S))$; it follows that
\[\DModRH(S) := \tx{Ind}(\DMod^{\tx{hol}, \tx{RS}}(S))\]
is a \emph{full subcategory} of $\DMod(S)$. We refer to it as the category of (ind)-regular holonomic D-modules.

More generally, for any aft scheme $S$, we choose a closed embedding $\iota: S^\tx{cl} \to T$ where $T$ is classical, aft and smooth. There is then a fully faithful functor
\[\DMod(S) \simeq \DMod(S^\tx{cl}) \xrightarrow{\iota_!} \DMod(T)\]
and we define
\[\DModRH(S) := \DMod(S) \times_{\DMod(T)} \DModRH(T) \subseteq \DMod(S).\]
Because for any two $S^\tx{cl} \to T_1, S^\tx{cl} \to T_2$ there exists another smooth classical $T_3$ such that
\[\xymatrix{
S^\tx{cl} \ar[d] \ar[r] & T_1 \ar[d] \\
T_2 \ar[r] & T_3
}\]
commutes, where all arrows are closed embeddings (take the pushout $T_1 \sqcup_{S^\tx{cl}} T_2$ then embed into another smooth classical scheme), by Kashiwara's lemma this definition does not depend on the choice of $\iota$.
\begin{lemma}
For arbitrary map $f: X \to Y$ between aft schemes, $f^!$ and $f_*$ both preserve the (ind-)regular holonomic subcategory.
\end{lemma}
\begin{proof}
	By nil-invariance we readily reduce to the case of $X, Y$ classical. It suffices to show that any $f$ can be completed to a diagram
	\[\xymatrix{
	X \ar^{g}[r] \ar^{f}[d] & \widetilde{X} \ar^{p}[d] \\
	Y \ar^{h}[r] & \widetilde{Y}	
	}\]
	where $g, h$ are closed immersions and $\widetilde{X}$, $\widetilde{Y}$ are smooth. Pick arbitrary $X \inj X', Y \inj Y'$ of closed immersion into smooth schemes, then set
	\[h := (Y \to Y' =: \widetilde{Y}) \hspace{1em} g := (X \xrightarrow{\Gamma_{h \circ f}} (X \times Y') \inj (X' \times Y') =: \widetilde{X}) \hspace{1em} p := \tx{pr}_2.\]
\end{proof}

For any laft prestack $\mc{X}$ and a tame twisting $\deRhamGerbe$ on $\mc{X}$, let
\[\tx{Split}(\deRhamGerbe, \mc{X}) \subseteq (\tx{Sch}^\tx{aff}_{\tx{aft}})_{/\mc{X}}\]
be the subcategory spanned by $f: S \to \mc{X}$ where $f^*(\deRhamGerbe)$ admits a trivialization on $S$. By definition of tame twisting, this is a basis for the derived \'{e}h-topology on the latter. Moreover, there exists a Grothendieck topology on $\tx{Split}(\deRhamGerbe, \mc{X})$ such that $U \to V$ is a covering iff it is in the \'{e}h-topology of $(\tx{Sch}^\tx{aff}_{\tx{aft}})_{/\mc{X}}$.

Let $\DMod_{\deRhamGerbe}^\tx{pre}(\mc{X}): \tx{Split}(\deRhamGerbe, \mc{X})^\op \to \tx{DGCat}$ denote the presheaf
\[f:(S \to \mc{X}) \mapsto \left(\vcenter{\xymatrix{\bullet \ar[r] \ar[d] & \tx{pt} \ar^{\tx{triv}}[d] \\
\tx{pt} \ar^(0.3){f^*(\deRhamGerbe)}[r] & \tx{Maps}(S, \overset{\circ}{\tb{Ge}})}}\right) \times^{\tx{Maps}(S, \overset{\circ}{\tb{Loc}_1})} \DMod(S)\]
where recall $\overset{\circ}{\tb{Loc}_1}$ is the prestack of \emph{tame} local systems. Note that this presheaf is in fact an \'{e}h-sheaf: upon choosing a trivialization for $\deRhamGerbe$, the sheaf condition reduces to the case of untwisted D-modules.

We then set
\[\tx{Twisted}(\DMod, \acute{e}h, \overset{\circ}{\tb{Ge}}, \overset{\circ}{\tb{Loc}_1}, \deRhamGerbe, \mc{X}) := \tx{RKE}_{\tx{Split}(\deRhamGerbe, \mc{X})^\op \to (\tx{Sch}^\tx{aff}_{\tx{aft}})_{/\mc{X}}^\op \xrightarrow{\tx{Yoneda}} (\tx{PreStk}_\tx{laft})_{/\mc{X}}^\op}\DMod_{\deRhamGerbe}^\tx{pre}(\mc{X})\]
i.e. the Kan extension of the \'{e}h-sheaf, and finally set
\[\DMod_{\deRhamGerbe}(\mc{X}) := \Gamma(\mc{X}, \tx{Twisted}(\DMod, \acute{e}h, \overset{\circ}{\tb{Ge}}, \overset{\circ}{\tb{Loc}_1}, \deRhamGerbe, \mc{X})).\]
Similarly we define
\[\DModRH_{\deRhamGerbe}(\mc{X}) := \Gamma(\mc{X}, \tx{Twisted}(\DModRH, \acute{e}h, \overset{\circ}{\tb{Ge}}, \overset{\circ}{\tb{Loc}_1}, \deRhamGerbe, \mc{X}))\]
\[\DModrhs_{\deRhamGerbe}(\mc{X}) := \Gamma(\mc{X}, \tx{Twisted}(\DModrhs, \acute{e}h, \overset{\circ}{\tb{Ge}}, \overset{\circ}{\tb{Loc}_1}, \deRhamGerbe, \mc{X}))\]

Derived h-descent of $\DMod$ guarantees that we recover the existing notion of $\DMod(\mc{X})$ when $\deRhamGerbe$ is in fact trivial. Moreover, the definition of $\oblv(\deRhamGerbe)$ and h-descent of $\DMod$ guarantees that
\[\DMod_\deRhamGerbe(\mc{X}) \simeq \Gamma(\mc{X}, \tx{Twisted}(\DMod, \acute{e}tale, \tb{Ge}_\dR, \tb{Loc}_1, \oblv(\deRhamGerbe), \mc{X}))\]
so our notation is unambiguous. The following fact justifies our definition in the ind-regular holonomic case as well:
\begin{proposition}
\label{prop:dmodrh-satisfies-h-descent}
	$\DModRH: \tx{Sch}_{\tx{aft}} \to \tx{DGCat}$ satisfies derived h-descent.
\end{proposition}
\begin{proof}
Let $f: U \to V$ be a derived h-cover of aft schemes. Since $\DMod$ itself satisfies derived h-descent, we have
	\[\DMod(V) \simeq \tx{Tot}(\DMod(U^\bullet / V))\]
	Now we have two full categories $\DModRH(V)$ and $\tx{Tot}(\DModRH(U^\bullet / V))$, where the latter is well-defined because $!$-pullback preserves the regular holonomic subcategories. Moreover, the first is a full subcategory of the second for the same reason. To prove the opposite inclusion, it suffices to show that, for every $\mc{F} \in \DMod(V)$, $f^!(\mc{F})$ being regular holonomic implies $\mc{F}$ is as well.
	
	Since $\DModRH$ is nil-invariant by definition, by replacing $f$ with $f^\tx{cl}$ we reduce to $U, V$ being classical schemes and $f$ a classical h-cover. By \cite[Proposition 3.4]{huber2013differential} one can assume that $f$ factors as $U \to U' \to V$ where $U \to U'$ is a Zariski covering and $U' \to V$ is proper surjective.
	Thus one reduce to each of the two cases.
	
	Suppose $f: U = \bigsqcup_i V_i \to V$ is a Zariski cover. By quasi-compactness assumption the index set $I$ is finite, so we can run induction on the size of $I$, and the claim now follows from the open-closed distinguished triangle and the previous lemma. Note that this argument also carries through in the non-ind setting.
	
	When $f$ is proper surjective, one can directly prove
	\[\DModRH(V) \to \tx{Tot}(\DModRH(U^\bullet / V))\]
	is an equivalence using the Barr-Beck-Lurie formalism (\cite[Corollary 4.7.5.3]{HA}): note that $f^!$ is continuous and conservative, and the diagrams in condition (2) of \emph{loc.cit.} are all left adjointable by proper base change of D-modules.
\end{proof}

\begin{remark}
	If we consider $\DModrhs$ i.e. the ``small'' category, then the argument above shows that it satisfies Zariski descent; nevertheless, it is unclear to us whether h-descent holds.
\end{remark}

For laft prestacks, the following constructions (carried out in \cite{gaitsgory2011crystals}) extend to the twisted D-module setting in an evident manner, which we will assume without further elaboration:
\begin{itemize}
\item The usual ``six functor'' formalism; and
\item Base change for arbitrary $!$-pull and schematic quasi-compact $*$-push.
\end{itemize}

\subsection{Analytic Sheaves}

To make contact with the analytic world it is more natural to switch to the locally ringed model of derived geometry. (The necessary comparison theorem is established in \cite{porta2017comparison}.) For any $S \in \tx{Sch}_{\tx{aft}}$, we let $S^\tx{an}$ denote the derived $\mb{C}$-analytic space as described in \cite{portathesis}, and let $S^\tx{topos}$ denote the underlying $\infty$-topos. We need some basic properties:
\begin{itemize}
\item By \cite[Lemma 2.3.2]{portathesis} $S^\tx{topos}$ is hypercomplete and (by definition) nil-invariant;
\item Recall what we call ``schemes'' are 0-geometric stacks, thus by \cite[Section 2.4]{porta2017comparison}, $S^\tx{topos}$ is 0-localic. Since it always (\cite[Lemma 3.1.24]{portathesis}) has enough points, by \cite[Proposition 2.5.15]{DAGV} there exists an underlying sober topological space $S^\tx{top}$ such that the underlying topos is equivalent to sheaves on the locale of $S^\tx{top}$. In fact, $S^\tx{top}$ can be identified with the set of $\mb{C}$-points of $S^\tx{cl}$ equipped with the analytic topology;
\item If $S_1 \to S_2$ is an closed immersion of (derived) schemes, then $S_1^\tx{topos} \to S_2^\tx{topos}$ is an closed immersion of topoi. It follows that the induced map $S_1^\tx{top} \to S_2^\tx{top}$ identifies $S_1^\tx{top}$ as a closed subspace of $S_2^\tx{top}$;
\item By \cite[XII 1.2]{SGA1} and the fact that classical truncation commutes with fiber product, the assignment $X \mapsto X^\tx{top}$ commutes with fiber products.
\end{itemize}

We set
\[\tx{Shv}^!(S) := \tx{Shv}(S^\tx{an}; \Vect) \simeq \tx{Shv}(S^\tx{an}) \tensor \Vect.\]
Using hypercompleteness of the topos, we know this category is equivalent to the (unbounded) derived category of the abelian category of $\Vect^\heartsuit$-valued sheaves on $S^\tx{an}$. It is thus equipped with a \emph{na\"ive} $t$-structure, where an object is coconnective iff it lands in $\Vect^{\ge 0}$ object-wise.

For any map $f: S_1 \to S_2$ between two such schemes, we have the usual adjunction $(f^*, f_*)$.
By the covariant Verdier duality\footnote{Our schemes $S$ are assumed to be separated, so $S^\tx{top}$ is Hausdorff.}, the category is also equipped with another pair of adjunctions $(f_!, f^!)$.

\begin{remark}
Since the formation $X \mapsto X^\tx{top}$ commutes with fiber products, it follows from \cite[7.3.1.19]{HTT} that we have ``$*$-pull, $*$-push'' and ``$!$-pull, $!$-push'' proper base change for $\tx{Shv}^!$. Because $!$-push and $*$-push are the same for proper morphisms, we also have ``$!$-pull, $*$-push'' base change.
\end{remark}

Now we specialize to the case of $S$ smooth and classical. Within $\tx{Shv}^!(S)$ we have the full subcategory $\tx{Shv}^!_{\tx{constr}}(S)$, the category of \emph{algebraically constructible} sheaves, which consists of objects whose cohomologies are constructible (by our convention, this in particular means all fibers are finite-dimensional) with respect to \emph{some} algebraic stratification. This category carries a \emph{perverse} $t$-structure, whose heart is the category of perverse sheaves on $S$.
We then define
\[\tx{Shv}_c^!(S) := \tx{Ind}(\tx{Shv}^!_\tx{constr}(S)).\]
It carries an unique $t$-structure such that its connective part is generated by $\tx{Shv}_\tx{constr}^{!}(S)^{\le 0}$.

\begin{convention}
Unless otherwise specified, we will only use the perverse $t$-structure on $\tx{Shv}^!_c(S)$.
\end{convention}

For every $f: X \to Y$ an closed embedding between aft schemes, $f^\tx{top}_*: \tx{Shv}_\tx{an}(X) \to \tx{Shv}_\tx{an}(Y)$ is fully faithful, as the underlying map of topological spaces is a closed embedding. Thus we take the same approach as before by specifying, for $S \in \tx{Sch}_{\tx{aft}}$:
\[\tx{Shv}_\tx{constr}^!(S) := \tx{Shv}_\tx{an}(S) \cap \tx{Shv}_\tx{constr}(T) \subseteq \tx{Shv}_\tx{an}(T) \hspace{1em} \tx{Shv}^!_c(S) := \tx{Ind}(\tx{Shv}_\tx{constr}^!(S))\]
for some (equiv.\ all) closed embedding $S^\tx{cl} \to T$ where $T$ is classical aft and smooth. These categories then inherit perverse $t$-structures. Same argument as above shows that arbitrary $!$-pull and $*$-push preserve the (ind-)constructible category.

\begin{remark}
The complex analytic setting differs from the D-module setting in many ways. For instance, the map $\tx{Shv}^!_c(S) \to \tx{Shv}^!(S)$, which ind-extends the inclusion $\tx{Shv}^!_\tx{constr}(S) \subseteq \tx{Shv}^!(S)$, is far from being fully faithful. (c.f. ``On the derived category of sheaves on a manifold'' by Neeman.) Also, the box product
\[\boxtimes: \tx{Shv}(S_1) \tensor \tx{Shv}(S_2) \to \tx{Shv}(S_1 \times S_2)\]
is in general not an equivalence.
\end{remark}

\begin{remark}
\label{remark:shvc-h-descent-problem}
	The proper base change noted above (and the conserativeness of $!$-fibers for proper surjection, which follows from Verdier duality) together imply (via same argument as in \cite{gaitsgory2013ind}) that the category $\tx{Shv}^!(S)$ satisfies h-descent. The same argument in Proposition~\ref{prop:dmodrh-satisfies-h-descent} then shows that $\tx{Shv}^!_c$ satisfies h-descent, and $\tx{Shv}^!_{\tx{constr}}$ satisfies at least Zariski descent. %
\end{remark}

We now define, for arbitrary laft prestack $\mc{X}$ and analytic gerbe $\deRhamGerbe_\tx{an}$ on $\mc{X}$,
\[\tx{Shv}^!_{\BettiGerbe}(\mc{X}):= \Gamma(\mc{X}, \tx{Twisted}(\tx{Shv}^!, \acute{e}h, \tb{Ge}_\tx{an}, \tb{Loc}_\tx{an}, \BettiGerbe, \mc{X}))\]
\[\tx{Shv}^!_{c, \BettiGerbe}(\mc{X}):= \Gamma(\mc{X}, \tx{Twisted}(\tx{Shv}_c^!, \acute{e}h, \tb{Ge}_\tx{an}, \tb{Loc}_\tx{an}, \BettiGerbe, \mc{X}))\]
\[\tx{Shv}^!_{\tx{constr}, \BettiGerbe}(\mc{X}):= \Gamma(\mc{X}, \tx{Twisted}(\tx{Shv}_\tx{constr}^!, \acute{e}h, \tb{Ge}_\tx{an}, \tb{Loc}_\tx{an}, \BettiGerbe, \mc{X}))\]

For laft prestacks, the following constructions extend to the (ind)-constructible setting in an evident manner, which we will assume without further elaboration:
\begin{itemize}
\item The usual ``six functor'' formalism; and
\item Proper base change \cite[7.3.1.19]{HTT}.
\end{itemize}

The following fact is useful:

\begin{lemma}
\label{lemma:shv-is-limit}
Let $\tb{S} \in \{\DMod, \DModrhs, \tx{Shv}^!, \tx{Shv}^!_{\tx{constr}}, \tx{Shv}^!_c\}$ and let $\tb{G}$ denote either $\overset{\circ}{\tb{Tw}}$ or $\tb{Ge}_\tx{an}$, depending on the context. Suppose we have $\mc{X} \simeq \colimit_{a \in A} \mc{X}_a$ a colimit in $\tx{PreStk}$, and $\mc{G}$ a $\tb{G}$-structure on $\mc{X}$. Then we have
	\[\tb{S}_\mc{G}(\mc{X}) \simeq \limit_{a \in A} \tb{S}_\mc{G}(\mc{X}_a).\]
\end{lemma}
\begin{proof}
	One can consider $\tb{S}_\mc{G}$ as a functor $(\tx{PreStk}_{/\mc{X}})^\tx{op} \to \tx{DGCat}$, which coincides with what we denoted above as $\tx{Twisted}(\tb{S}, \acute{e}h, \tb{G}, \tb{L}, \mc{G}, \mc{X})$ (where $\tb{L}$ denote either $\overset{\circ}{\tb{Loc}_1}$ or $\tb{Loc}_\tx{an}$ depending on the context), and it would suffice to prove that $\tb{S}_\mc{G}$ is the right Kan extension of its restriction to $(\tx{Sch}^\tx{aff}_{/\mc{X}})^\op$ along the Yoneda embedding. The same argument as \cite[Chapter 5, Lemma 3.2.2]{GR-DAG1} shows that it suffices to know that $\tb{S}_\mc{G}$ satisfies Zariski descent, which we have by construction.
\end{proof}

\subsection{Riemann-Hilbert}
The usual Riemann-Hilbert equivalence states that, when $X$ is a smooth aft scheme, there exists a $t$-exact, compact-preserving equivalence
\[\tx{RH}: \DModrhs(S) \simeq \tx{Shv}^!_{\tx{constr}}(S)\hspace{1em} \tx{RH}: \DModRH(S) \simeq \tx{Shv}^!_c(S)\]
that intertwines the right D-module $t$-structure one LHS and perverse $t$-structure on RHS, as well as all six functors.

For $S$ a general aft scheme, fix some closed embedding $f: S^\tx{cl} \to T$, $T$ classical smooth aft, and let $j: U \to T$ be the inclusion of the complement. Then $\DModrhs(S)$ can be characterized as the subcategory of $\DModrhs(T)$ consisting of objects $\mc{F}$ such that $j^!(\mc{F}) \simeq 0$; the same description works for $\tx{Shv}_\tx{constr}^!(S)$. It follows that the above equivalence can be extended to arbitrary aft scheme $S$.

Next, \cite[Lemma 3.7]{zhao2020tame} gives a the gerbe analytification map
\[\tx{RH}: \overset{\circ}{\tb{Ge}} \to \tb{Ge}_{\tx{an}}\]
which intertwines the action of (tame) local systems on two sides of Riemann-Hilbert. Thus, termwise application of the equivalence, we have
\begin{corollary}
\label{cor:riemann-hilbert}
For any laft prestack $\mc{X}$ and a tame de Rham gerbe $\deRhamGerbe$, there are equivalences
\[\tx{RH}: \DModrhs_\deRhamGerbe(\mc{X}) \simeq \tx{Shv}^!_{\tx{constr}, \tx{RH}(\deRhamGerbe)}(\mc{X})\]
\[\tx{RH}: \DMod^\tx{rh}_\deRhamGerbe(\mc{X}) \simeq \tx{Shv}^!_{c, \tx{RH}(\deRhamGerbe)}(\mc{X}).\]
\end{corollary}

We will now specialize to the case of $[\tx{Gr}_{\cT}]_{\Ran_\dR}$ for any smooth curve $X$, with a fixed $\mb{C}$ point $x \in X(\mb{C})$. We will further fix a factorizable tame twisting $\deRhamGerbe$.

\begin{convention}
For any space $\mc{Y}$ that maps to $[\tx{Gr}_{\cT}]_{\Ran_\dR}$ in an unambiguous manner, we will abuse notation and let $\deRhamGerbe$ denote the pullback to $\mc{Y}$ as well.
\end{convention}

For each $p \in \ECOp$, we have the corresponding indscheme $[\tx{Gr}_{\cT}]_{U_\Ran(p)} \subseteq [\tx{Gr}_{\cT}]_{\Ran}$; we set
\[\DMod_{\deRhamGerbe}^\tx{Ran-constr}([\tx{Gr}_{\cT}]_{U_\Ran(p)}) \subseteq \DMod_{\deRhamGerbe}^{\tx{hol}, \tx{RS}}([\tx{Gr}_{\cT}]_{U_\Ran(p)})\]
denote the full subcategories of objects that are locally constant along each strata of the cardinality stratification of the Ran space; we refer to such objects as \emph{Ran-constructible} twisted D-modules. We also set
\[\DMod_{\deRhamGerbe}^\tx{Ran-c}([\tx{Gr}_{\cT}]_{U_\Ran(p)}) := \tx{Ind}(\DMod_{\deRhamGerbe}^\tx{Ran-constr}([\tx{Gr}_{\cT}]_{U_\Ran(p)})).\]
Similarly, let
\[\tx{Shv}^!_{\tx{Ran-constr}, \BettiGerbe}([\tx{Gr}_{\cT}]_{U_\Ran(p)}) \subseteq \tx{Shv}^!_{\tx{constr}, \BettiGerbe}([\tx{Gr}_{\cT}]_{U_\Ran(p)})\]
denote the full subcategory consisting of objects locally constant along the same strata, and set
\[\tx{Shv}^!_{\tx{Ran-c}, \BettiGerbe}([\tx{Gr}_{\cT}]_{U_\Ran(p)}) := \tx{Ind}(\tx{Shv}^!_{\tx{Ran-constr}, \BettiGerbe}([\tx{Gr}_{\cT}]_{U_\Ran(p)}).\]
Clearly, Riemann-Hilbert preserves these subcategories.
\begin{definition}
We let
\[\tx{FactAlg}_{\nonun}(\DMod_{\deRhamGerbe}^\tx{Ran-c}([\tx{Gr}_{\cT}]_\Ran)) \subseteq \tx{FactAlg}_{\nonun}(\DMod_{\deRhamGerbe}([\tx{Gr}_{\cT}]_\Ran)\]
denote the full subcategory consisting of factorization algebras whose underlying sheaf belongs to the full subcategory $\DMod_{\deRhamGerbe}^\tx{Ran-c}([\tx{Gr}_{\cT}]_{U_\Ran(p)})$ for each $p$ under the equivalence from Proposition \ref{prop:fact-cat-alg-partition-description}. Similarly, for each $A$ we have
\[A\tx{-FactMod}_{\nonun, x}(\DMod_{\deRhamGerbe}^\tx{Ran-c}([\tx{Gr}_{\cT}]_{\Ran_x}))\]
analogously defined using the module version of $U_\Ran$.
\end{definition}

\begin{remark}
It follows from the partition description that a factorization module supported at a single point $x$ is ind-regular holonomic (resp. ind-Ran-constructible) if and only if the underlying algebra object is. In other words, we have
\[A\tx{-FactMod}_{\nonun, x}(\DMod_{\deRhamGerbe}^\tx{Ran-c}([\tx{Gr}_{\cT}]_{\Ran_x})) \simeq A\tx{-FactMod}_{\nonun, x}(\DMod_{\deRhamGerbe}([\tx{Gr}_{\cT}]_{\Ran_x}))\]
\end{remark}

\begin{definition}
\label{def:shv-ranc-weak-constr-factcat}
We let
\[\tx{Shv}^!_{\tx{Ran-c}, \BettiGerbe}([\tx{Gr}_{\cT}]_\Ran)\]
denote the weakly constructible factorization category corresponding to the system $\{\tx{Shv}^!_{\tx{Ran-c}, \BettiGerbe}([\tx{Gr}_{\cT}]_{U_\Ran(p)})\}$ as defined in Definition \ref{def:constructible-fact-objects}. We likewise have a constructible factorization module category
\[\tx{Shv}^!_{\tx{Ran-c}, \BettiGerbe}([\tx{Gr}_{\cT}]_{\Ran_x}).\]
\end{definition}

Combine discussions above, we arrive at

\begin{corollary}
\label{cor:RH-for-fact-mod}
	Let $X$ be a smooth $\mb{C}$-curve and $x \in X(\mb{C})$. Let $\deRhamGerbe$ be a tame twisting on $X$.
	The Riemann-Hilbert equivalence induces an equivalence
	\[\tx{RH}: \tx{FactAlg}_{\nonun}(\DMod_{\deRhamGerbe}^\tx{Ran-c}([\tx{Gr}_{\cT}]_\Ran)) \simeq \tx{FactAlg}^\CoMark_{\nonun, \ConstructibleMark, \PartitionMark}(\tx{Shv}^!_{\tx{Ran-c}, \BettiGerbe}([\tx{Gr}_{\cT}]_\Ran))\]
	and, for every $A$ on the LHS, an equivalence
	\[\tx{RH}: A\tx{-FactMod}_{\nonun, x}(\DMod_{\deRhamGerbe}([\tx{Gr}_{\cT}]_{\Ran_x})) \simeq \tx{RH}(A)\tx{-FactMod}_{\nonun, \ConstructibleMark, \PartitionMark, x}(\tx{Shv}^!_{\tx{Ran-c}, \BettiGerbe}([\tx{Gr}_{\cT}]_{\Ran_x})).\]
\end{corollary}

\newpage
\printbibliography[heading=bibliography]

\end{document}